\documentclass[12pt,reqno]{amsart}
\usepackage{etoolbox}
\usepackage{mathtools,amssymb,amsthm,thmtools,thm-restate,enumitem}
\usepackage[dvipsnames,svgnames]{xcolor}
\usepackage{tikz,tikz-cd,lastpage,indentfirst,calc}
\usepackage[letterpaper,margin=1in,footskip=0.5in]{geometry}
\usepackage[colorlinks=true,linkcolor=Purple,filecolor=black,urlcolor=blue,citecolor=blue,pdfpagemode=FullScreen,linktocpage=true,hyperfootnotes=true,plainpages=false,hypertexnames=false,pdfpagelabels]{hyperref}
\usepackage[capitalize,noabbrev]{cleveref}
\usetikzlibrary{patterns.meta,calc,decorations.markings,decorations.pathreplacing,arrows.meta,fadings}

\makeatletter

\sloppypar\allowdisplaybreaks
\let\ut@originalitem\item
\newcommand{\ut@customitem}[1][]{\ifx\relax#1\relax\ut@originalitem\else\ut@originalitem[#1]\protected@edef\@currentlabel{\textup{#1}}\fi}
\newcommand{\enableCustomItemLabels}{\let\item\ut@customitem}
\SetEnumitemKey{widestlabel}{labelwidth = \widthof{\ref{enum-\EnumitemId}},after = \label{enum-\EnumitemId}}
\setlist[itemize]{font=\upshape,listparindent=\parindent,leftmargin=*}
\newlist{lst}{enumerate}{4}
\setlist[lst]{
    before=\enableCustomItemLabels, 
    font=\upshape, 
    label=(\roman*), 
    ref=\textup{(\roman*)}, 
    listparindent=\parindent, 
    labelindent=0pt, 
    widestlabel, 
    leftmargin=!,
}

\@ifundefined{newcounteralias}{}{%
  \renewcommand\thmt@autorefsetup{\@xa\def\csname\thmt@envname autorefname\@xa\endcsname\@xa{\thmt@thmname}}%
}

\newcounter{alpha}
\newcommand{\enunciationparent}{subsection}% enunciations are numbered within this level

\newcommand{\defineEnunciation}[4]{%
  \ifcsname c@childof\enunciationparent\endcsname\else
    \newcounter{childof\enunciationparent}\counterwithin{childof\enunciationparent}{\enunciationparent}%
  \fi
  \declaretheoremstyle[spaceabove=\medskipamount,spacebelow=\medskipamount,headpunct={},headfont=\bfseries,notebraces={(}{)},headformat={\NAME~\NUMBER\NOTE.},postheadhook={\phantom{}{}},bodyfont=\upshape]{#1Style}%
  \declaretheorem[style=#1Style,name=#2,numberlike=childof\enunciationparent]{#1}%
  \declaretheorem[style=#1Style,name=#3,numberlike=childof\enunciationparent]{#1s}%
  \declaretheorem[style=#1Style,name=#2,sibling=alpha]{#1alpha}%
  \crefname{#1}{#2}{#3}\crefname{#1s}{#2}{#3}\crefname{#1alpha}{#2}{#3}%
}

\defineEnunciation{corollary}{Corollary}{Corollaries}{gray}
\defineEnunciation{theorem}{Theorem}{Theorems}{blue}
\defineEnunciation{conjecture}{Conjecture}{Conjectures}{blue}
\defineEnunciation{proposal}{Proposal}{Proposals}{blue}
\defineEnunciation{lemma}{Lemma}{Lemmas}{Dandelion}
\defineEnunciation{proposition}{Proposition}{Propositions}{purple}
\defineEnunciation{claim}{Claim}{Claims}{teal}
\defineEnunciation{construction}{Construction}{Constructions}{Purple}
\defineEnunciation{warning}{Warning}{Warnings}{magenta}
\defineEnunciation{example}{Example}{Examples}{ForestGreen}
\defineEnunciation{definition}{Definition}{Definitions}{olive}
\defineEnunciation{fact}{Fact}{Facts}{olive}
\defineEnunciation{convention}{Convention}{Conventions}{Red}
\defineEnunciation{remark}{Remark}{Remarks}{gray}
\defineEnunciation{notation}{Notation}{Notations}{blue}

\newcommand{\ut@sectionref}[1]{%
  \crefformat{#1}{\S##2##1##3}\Crefformat{#1}{\S##2##1##3}%
  \crefrangeformat{#1}{\S\S##3##1##4--##5##2##6}%
  \crefmultiformat{#1}{\S\S##2##1##3}{ and ##2##1##3}{, ##2##1##3}{ and ##2##1##3}%
}
\forcsvlist{\ut@sectionref}{section,subsection,subsubsection}

\newcommand{\ut@alias}[3]{%
  \edef\ut@tmp{\noexpand\@namedef{#1#3}%
    {{\expandafter\noexpand\csname #2\endcsname{#3}}}}\ut@tmp}
\newcommand{\ut@op}[2]{%
  \edef\ut@tmp{\noexpand\@namedef{#2}%
    {\expandafter\noexpand\csname #1\endcsname{#2}}}\ut@tmp}
\DeclareMathAlphabet{\matheulscr}{U}{eus}{m}{n}
\forcsvlist{\ut@alias{b}{mathbb}}{C,R,Z}
\forcsvlist{\ut@alias{c}{mathcal}}{A,B,G,H,L,M,P,Q,R,S,X,Y,Z}
\forcsvlist{\ut@alias{e}{matheulscr}}{A,C,D,F,G,H,L,M,N,P,R,S,T,V,X,Y,Z}
\forcsvlist{\ut@alias{f}{mathfrak}}{X,Y}
\forcsvlist{\ut@alias{r}{mathrm}}{A,B,E,N,O,T}
\forcsvlist{\ut@alias{sf}{mathsf}}{X}

\forcsvlist{\ut@op{mathrm}}{op,pt,id,ev,coev}
\forcsvlist{\ut@op{operatorname}}{colim,End,FPdim,Hom,Int,Irr,Obj,PL,pr,res,SO,tr,Tr}
\forcsvlist{\ut@alias{}{mathsf}}{Fun,Hilb,Mod,Rep,Set,Vec}
\newcommand{\src}{\mathrm{src}}
\newcommand{\targ}{\mathrm{targ}}
\newcommand{\glu}{\mathrm{gl}}
\newcommand{\cl}{\mathrm{cl}}
\newcommand{\Sk}{\mathrm{Sk}}
\newcommand{\Mfld}{\mathsf{Mfld}}
\newcommand{\Disk}{\mathsf{Disk}}
\newcommand{\eRep}{{\eR}\mathrm{ep}}
\newcommand{\eFun}{{\eH}\mathrm{om}}
\newcommand{\eHom}{{\eH}\mathrm{om}}
\newcommand{\eHilb}{{\eH}\mathrm{ilb}}
\newcommand{\eval}{\mathrm{eval}}
\newcommand{\inn}{\mathrm{in}}
\newcommand{\out}{\mathrm{out}}
\newcommand{\psn}{\Phi}
\newcommand{\A}{\eA}
\newcommand{\C}{\eC}
\newcommand{\D}{\eD}
\newcommand{\Z}{\eZ}

\newcommand{\e}{\varepsilon}
\newcommand{\bdy}{\partial}
\newcommand{\blt}{\bullet}

\newcommand{\tld}[1]{{\widetilde{#1}}}
\newcommand{\fcj}[1]{\overline{#1}}
\newcommand{\orev}{\overline}

\newcommand{\rad}{\sqrt}
\newcommand{\bra}[1]{\langle#1\rvert}
\newcommand{\ket}[1]{\lvert#1\rangle}
\newcommand{\bkt}[2]{\langle #1|#2 \rangle}
\newcommand{\Cstar}{\ensuremath{\mathrm{C}^{*}}}
\newcommand{\Hstar}{\ensuremath{\mathrm{H}^{*}}}
\newcommand{\defn}[1]{\emph{#1}}
\newcommand{\cent}{\textup{\textcent}}
\newcommand{\bigboxplus}{\mathop{\vcenter{\hbox{\tikz[very thick,transform shape,scale=1.2,inner sep=0pt,minimum width=0pt,outer sep=0pt,x=1em,y=1em,yshift=1em,baseline={(current bounding box.center)}]\draw (-.5,-.5)rectangle(.5,.5) (.5,0)--(-.5,0) (0,.5)--(0,-.5);}}}}
\newcommand{\nn}[1]{\bgroup\color{orange}[[#1]]\egroup}
\newcommand{\itemstep}[1]{\smallskip\noindent\underline{#1}}
\renewcommand{\setminus}{\smallsetminus}
\renewcommand{\amalg}{\sqcup}
\renewcommand{\star}{\dag}
\let\oldforall\forall\renewcommand{\forall}{\operatorname{\oldforall}}
\AtBeginDocument{\let\bar\overline}
\DeclareRobustCommand{\coloneq}{\mathrel{\vcenter{\baselineskip0.5ex\lineskiplimit0pt\hbox{\scriptsize{\upshape.}}\hbox{\scriptsize{\upshape.}}}}=}
\DeclareRobustCommand{\coloneqq}{\coloneq}
\DeclareRobustCommand{\eqcolon}{\mathrel{\scalebox{-1}[1]{$\coloneq$}}}
\NewDocumentCommand{\defeq}{o}{%
  \IfNoValueTF{#1}%
    {\mathrel{\smash{\overset{\mathclap{\text{def}}}{=}}}}%
    {\mathrel{\smash{\underset{#1}{\overset{\text{\tiny{}\mathclap{def}}}{=}}}}}%
}
\DeclareFontFamily{U}{dmjhira}{}\DeclareFontShape{U}{dmjhira}{m}{n}{ <-> dmjhira }{}
\DeclareRobustCommand{\yo}{\text{\usefont{U}{dmjhira}{m}{n}\symbol{"48}}}
\newcommand{\axiomtag}[2]{\phantomsection\protected@edef\@currentlabel{\protect\textup{(#1)}}\label{#2}{(#1)}}

\newcommand{\undCtemp}{\bgroup\smash[b]{\text{$\underset{\smash{\raisebox{.225em}{$\scriptstyle\kern.05ex\rightarrow$}}}{\smash{\C}}$}}\egroup}
\newcommand{\undC}{\vphantom{\C}{\undCtemp}}
\newcommand{\undD}{\D}
\NewDocumentCommand{\fldl}{o m m o}{\A_{#2}(#3\IfValueT{#4}{,#4})}
\NewDocumentCommand{\undfldl}{o m m o}{\A_{#2}(#3\IfValueT{#4}{,#4})}
\NewDocumentCommand{\fld}{o m m o}{#2\IfValueT{#1}{_{#1}}(#3\IfValueT{#4}{;#4})}
\NewDocumentCommand{\undfld}{o m m o}{#2\IfValueT{#1}{_{#1}}(#3\IfValueT{#4}{;#4})}
\newcommand{\lU}[1]{#1_U}
\NewDocumentCommand{\fldlU}{o m m o}{\lU{#2}(#3\IfValueT{#4}{,#4})}

\definecolor{violet}{RGB}{148,0,211}
\newcommand{\colllGr}{gray!10}
\newcommand{\colX}{red!80!black}   \newcommand{\colY}{blue!85!black}
\newcommand{\colZ}{green!55!black} \newcommand{\colW}{violet}
\newcommand{\colQ}{orange!95!red}  \newcommand{\colR}{cyan!60!black}
\newcommand{\colG}{Green!50!black} \newcommand{\cola}{Green!75!white}
\newcommand{\colb}{orange}         \newcommand{\colc}{violet!55!black}
\newcommand{\filf}{brown!6}        \newcommand{\filg}{blue!6}
\newcommand{\filh}{green!8}        \newcommand{\filk}{violet!7}
\newcommand{\filid}{black!6}

\NewDocumentEnvironment{tkz}{O{}}
  {\sbox0{$\relax$}%
   \csname begin\endcsname{tikzpicture}[baseline={([yshift=-\the\dimexpr\fontdimen22\textfont2\relax]current bounding box.center)},#1]}
  {\csname end\endcsname{tikzpicture}}

\pgfdeclarelayer{bottom}\pgfdeclarelayer{back}\pgfsetlayers{bottom,back,main}
\def\tkzframingwidth{0.1cm}

\tikzset{
  mid>/.style={decoration={markings, mark=at position #1 with {\arrow{>}}}, postaction={decorate}}, mid>/.default=0.5,
  mid</.style={decoration={markings, mark=at position #1 with {\arrow{<}}}, postaction={decorate}}, mid</.default=0.5,
  coL/.style={decoration={markings, mark=at position #1 with
      {\draw[-{Stealth[length=2.8pt,width=2.6pt]},line width=0.7pt] (0,-0.13)--(0,0.13);}}, postaction={decorate}}, coL/.default=0.5,
  coR/.style={decoration={markings, mark=at position #1 with
      {\draw[-{Stealth[length=2.8pt,width=2.6pt]},line width=0.7pt] (0,0.13)--(0,-0.13);}}, postaction={decorate}}, coR/.default=0.5,
  strand/.style={thick, blue!70!black},
  Fa/.style={red!10}, Fb/.style={blue!12}, Fc/.style={blue!46}, Fd/.style={blue!26}, Fe/.style={blue!38},
  bd/.style={thick, black!80}, eqt/.style={black!45, line width=0.5pt},
  dm/.style={line width=1.2pt, Orange}, dp/.style={line width=1.2pt, Green!45},
  fl/.style={font=\scriptsize}, tl/.style={font=\tiny},
  brace/.style={decorate, decoration={brace, amplitude=5pt}, thick},
  zshift/.style={shift={(0,0,#1)}},
  cpn/.style={rounded corners, draw=black, line width=1pt, font=\scriptsize\sffamily, align=center},
  pics/rbox/.style n args={4}{code={%
      \def\hh{#1}\def\ll{#2}\def\rr{#3}\def\txt{#4}%
      \draw[cpn, pic actions] (-\ll-\hh, -\hh) rectangle (\rr+\hh, \hh);
      \node[font=\scriptsize\sffamily] at (0, 0) {\txt};
      \coordinate (.north) at (0, \hh);      \coordinate (.south) at (0, -\hh);
      \coordinate (.east)  at (\rr+\hh, 0);  \coordinate (.west)  at (-\ll-\hh, 0);
      \coordinate (.center) at (0,0);}},
  pics/link/.style={code={%
      \def\lkr{0.3}%
      \begin{scope}[#1]
        \draw[thin,dashed] (0,0) circle (\lkr);
        \draw[thick,dashed] (-\lkr,0) -- (\lkr,0);
        \draw[->,thick,cyan!75!RoyalBlue] ({-0.75*\lkr},0) -- ++(90:{0.75*\lkr});
        \draw[->,thick,cyan!75!RoyalBlue] ({0.75*\lkr},0) -- ++(90:{0.7*\lkr});
      \end{scope}}},
  pics/tro/.style args={#1}{code={\draw[cyan!75!RoyalBlue,pic actions,->] (-#1,0) -- (#1,0);}},
  pics/tro/.default=0.15,
  frameL/.style={preaction={draw=gray, line width=\tkzframingwidth, opacity=0.5,
      every path/.style={}, line cap=butt, decorate,
      decoration={curveto, amplitude=0, raise=0.5*\tkzframingwidth}}},
  frameR/.style={preaction={draw=gray, line width=\tkzframingwidth, opacity=0.5,
      every path/.style={}, line cap=butt, decorate,
      decoration={curveto, amplitude=0, raise=-0.5*\tkzframingwidth}}},
}

\newcommand{\rbox}[6][]{\pic[fill=white,#1] at #2 {rbox={#3}{#4}{#5}{#6}};} % Usage: \rbox[opts]{coordinate}{half-height}{left dist}{right dist}{text}

\newcommand{\lk}[3][0]{%
  \begin{scope}[shift={#2},rotate=#1]
    \draw[thin, dash pattern=on 1.2pt off 1.0pt] (180:#3) arc (180:0:#3);
    \draw[thin, dash pattern=on 1.2pt off 1.0pt] (0:#3) arc (0:-180:#3);
    \draw[line width=0.8pt, dash pattern=on 1.9pt off 1.2pt] (-#3,0) -- (#3,0);
    \draw[-{Stealth[length=2.2pt,width=2.0pt]},line width=0.6pt] (-{0.62*#3},0) -- ++(90:{0.62*#3});
    \draw[-{Stealth[length=2.2pt,width=2.0pt]},line width=0.6pt] ({0.62*#3},0) -- ++(90:{0.62*#3});
  \end{scope}%
  \fill #2 circle (1.5pt);
}

\def\centerarc(#1)(#2:#3:#4){([shift=(#2:#4)]#1) arc (#2:#3:#4)}

\tikzset{even odd clip/.code={\pgfseteorule}}
\newcommand{\invclip}[1][]{%
  \def\invclip@parse##1;{%
    \begin{pgfinterruptboundingbox}%
      \path[clip, even odd clip, #1] ##1 (-2\paperwidth,-2\paperheight) rectangle (2\paperwidth,2\paperheight);%
    \end{pgfinterruptboundingbox}}%
  \invclip@parse
}

\tikzset{
  annulus/inner/.style={draw=none}, annulus/outer/.style={draw=none},
  inner/.style={annulus/inner/.append style={draw, #1}},
  outer/.style={annulus/outer/.append style={draw, #1}}
}
\newcommand{\annulus}[4][]{%
  \begin{scope}
    \tikzset{annulus/inner/.style={draw=none}, annulus/outer/.style={draw=none}}\tikzset{#1}%
    \invclip #2 circle (#3);
    \path[#1] #2 circle (#4);
    \draw[annulus/outer] #2 circle (#4);
  \end{scope}%
  \begin{scope}
    \tikzset{annulus/inner/.style={draw=none}, annulus/outer/.style={draw=none}}\tikzset{#1}%
    \draw[annulus/inner] #2 circle (#3);
  \end{scope}
}

\tikzdeclarepattern{name=primeddots, type=uncolored,
  bounding box={(-0.6pt,-0.6pt) and (0.6pt,0.6pt)}, tile size={(5pt,5pt)},
  tile transformation={rotate=60}, code={\fill (0pt,0pt) circle (.35pt);}}
\tikzdeclarepattern{name=primedbox, type=uncolored,
  bounding box={(-1pt,-1pt) and (1pt,1pt)}, tile size={(5pt,5pt)},
  tile transformation={rotate=60}, code={\draw[very thin] (-0.9pt,-0.9pt) rectangle (0.9pt,0.9pt);}}
\tikzdeclarepattern{name=primedstar, type=uncolored,
  bounding box={(-1.2pt,-1.2pt) and (1.2pt,1.2pt)}, tile size={(5pt,5pt)},
  tile transformation={rotate=30}, code={\draw[very thin] (-1.2pt,0pt) -- (1.2pt,0pt); \draw[very thin] (-0.6pt,-1.04pt) -- (0.6pt,1.04pt); \draw[very thin] (-0.6pt,1.04pt) -- (0.6pt,-1.04pt);}}
\tikzdeclarepattern{name=primedtarget, type=uncolored,
  bounding box={(-2pt,-2pt) and (2pt,2pt)}, tile size={(10pt,10pt)},
  tile transformation={rotate=15}, code={\fill (0,0) circle (0.4pt); \draw[very thin] (0,0) circle (1.8pt);}}
\tikzdeclarepattern{name=primedbubbles, type=uncolored,
  bounding box={(-3pt,-3pt) and (4pt,4pt)}, tile size={(11pt,11pt)},
  tile transformation={rotate=35}, code={\draw[very thin] (-2pt,1pt) circle (1.2pt); \draw[very thin] (2.5pt,-1.5pt) circle (2.2pt);}}
\tikzdeclarepattern{name=primedreticle, type=uncolored,
  bounding box={(-3pt,-3pt) and (3pt,3pt)}, tile size={(10pt,10pt)},
  tile transformation={rotate=20}, code={\draw[very thin] (0,0) circle (1.2pt); \draw[very thin] (-2.5pt,0) -- (-1.4pt,0); \draw[very thin] (1.4pt,0) -- (2.5pt,0); \draw[very thin] (0,-2.5pt) -- (0,-1.4pt); \draw[very thin] (0,1.4pt) -- (0,2.5pt);}}
\tikzdeclarepattern{name=primedhexnuts, type=uncolored,
  bounding box={(-2pt,-2pt) and (2pt,2pt)}, tile size={(8.5pt,8.5pt)},
  tile transformation={rotate=30}, code={\draw[very thin] (0:1.6pt) -- (60:1.6pt) -- (120:1.6pt) -- (180:1.6pt) -- (240:1.6pt) -- (300:1.6pt) -- cycle;}}
\tikzdeclarepattern{name=primedxmarks, type=uncolored,
  bounding box={(-1.5pt,-1.5pt) and (1.5pt,1.5pt)}, tile size={(8pt,8pt)},
  tile transformation={rotate=10}, code={\draw[very thin] (-1.2pt,-1.2pt) -- (1.2pt,1.2pt); \draw[very thin] (-1.2pt,1.2pt) -- (1.2pt,-1.2pt);}}

\tikzset{
  primedregion/.style={preaction={fill=#1}, pattern=primeddots, draw=#1},
  boxregion/.style={preaction={fill=#1}, pattern=primedbox, draw=#1},
  starregion/.style={preaction={fill=#1}, pattern=primedstar, draw=#1},
}

\tikzset{
  region opts/.cd,
    type/.initial=lines, color/.initial=black, angle/.initial=45, distance/.initial=5pt,
    .unknown/.code={\pgfkeyssetvalue{/tikz/region opts/type}{\pgfkeyscurrentname}},
  /tikz/.cd,
  region/.code={%
    \pgfkeys{/tikz/region opts/.cd, #1}%
    \edef\regtype{\pgfkeysvalueof{/tikz/region opts/type}}%
    \pgfkeysalso{/tikz/region apply/.expand once=\regtype}%
  },
  region apply/.is choice,
  region apply/lines/.style={/tikz/postaction/.expanded={pattern={Lines[angle=\pgfkeysvalueof{/tikz/region opts/angle}, distance=\pgfkeysvalueof{/tikz/region opts/distance}]}, pattern color=\pgfkeysvalueof{/tikz/region opts/color}}},
  region apply/grid/.style={
    /tikz/postaction/.expanded={pattern={Lines[angle=\pgfkeysvalueof{/tikz/region opts/angle}, distance=\pgfkeysvalueof{/tikz/region opts/distance}]}, pattern color=\pgfkeysvalueof{/tikz/region opts/color}},
    /tikz/postaction/.expanded={pattern={Lines[angle=\pgfkeysvalueof{/tikz/region opts/angle}+90, distance=\pgfkeysvalueof{/tikz/region opts/distance}]}, pattern color=\pgfkeysvalueof{/tikz/region opts/color}}},
  region apply/dots/.style={/tikz/postaction/.expanded={pattern=primeddots, pattern color=\pgfkeysvalueof{/tikz/region opts/color}}},
  region apply/box/.style={/tikz/postaction/.expanded={pattern=primedbox, pattern color=\pgfkeysvalueof{/tikz/region opts/color}}},
  region apply/hyperchaos/.style={
    /tikz/postaction/.expanded={pattern=primedtarget,  pattern color=red},
    /tikz/postaction/.expanded={pattern=primedbubbles, pattern color=orange},
    /tikz/postaction/.expanded={pattern=primedreticle, pattern color=green!80!black},
    /tikz/postaction/.expanded={pattern=primedhexnuts, pattern color=cyan},
    /tikz/postaction/.expanded={pattern=primedxmarks,  pattern color=magenta}},
}

\def\tkzangx{30}
\def\tkzangy{30+90}
\def\tkzangz{-70}
\def\tkzfout{3}
\newcommand{\hgt}{0.85}
\newcommand{\makediagL}[3][]{
\def\tkzsep{#2}
\def\tkzbdysep{#3}
\def\tkzglabel{\scriptsize$\fcj{g}$}
\begin{tkz}[
  x={({cos(\tkzangx)},{sin(\tkzangx)})},
  y={({cos(\tkzangy)},{sin(\tkzangy)})},
  z={({cos(\tkzangz)},{sin(\tkzangz)})},
  scale=0.5,
  cap=round,
  #1
]
\begin{scope}[zshift=-\tkzbdysep]
\begin{scope}[zshift={-2-\tkzsep}]
\fill[gray!7,fading angle=-30,path fading=west,opacity=0.75]
(0,0)
  .. controls +(0:0.8) and +(180:0.8) .. (2,-1)
  .. controls +(0:0.8) and +(180:0.8) .. (4,0)
  .. controls +(180:0.8) and +(0:0.8) .. (2,1)
  .. controls +(180:0.8) and +(0:0.8) .. (0,0);
\draw[red,very thick,cap=round]
(0,0)
.. controls +(0:1) and +(0:-1) ..
(2,-1)
.. controls +(0:1) and +(0:-1) ..
(4,0);
\node[violet!50!black,inner sep=1pt](vil) at (-1,0,-1) {\scriptsize$\fcj{v_i}$};
\draw[violet!50!black,->](vil)to[out=-10,in=180-(90-115)](2,0,-1.5);
\fill[white]
  (0,0,0)
  .. controls +(0:1) and +(0:-1) .. (2,1,0)
  .. controls +(0:1) and +(0:-1) .. (4,0,0)
  -- (4,0,-\tkzfout)
  .. controls +(180:1) and +(180:-1) .. (2,1,-\tkzfout)
  .. controls +(180:1) and +(180:-1) .. (0,0,-\tkzfout)
  -- cycle;
\fill[gray!15,fading angle=-30,path fading=west,opacity=1]
  (0,0,0)
  .. controls +(0:1) and +(0:-1) .. (2,1,0)
  .. controls +(0:1) and +(0:-1) .. (4,0,0)
  -- (4,0,-\tkzfout)
  .. controls +(180:1) and +(180:-1) .. (2,1,-\tkzfout)
  .. controls +(180:1) and +(180:-1) .. (0,0,-\tkzfout)
  -- cycle;
\node[violet,inner sep=1pt](vjl) at (0,2,-1) {\scriptsize$w_j$};
\draw[violet,->](vjl)to[out=45,in=115](2,0,-2);
\draw[blue,very thick,cap=round]
(0,0,0)
.. controls +(0:1) and +(0:-1) ..
(2,1,0)
.. controls +(0:1) and +(0:-1) ..
(4,0,0);
\draw[Green,very thick,fading angle=-30,path fading=west](0,0,0)--(0,0,-\tkzfout);
\draw[Green,very thick,fading angle=-30,path fading=west](4,0,0)--(4,0,-\tkzfout);
\end{scope}
\begin{scope}
\node[Brown!50!black,inner sep=1pt](vil) at (-1,0,-1) {\scriptsize$\fcj{p_i}$};
\draw[Brown!50!black,->](vil)to[out=-10,in=180-(90-115)](2,0,-1.5);
\end{scope}
\fill[gray!7,fading angle=-30,path fading=west,opacity=0.75]
(0,0)
  .. controls +(0:0.8) and +(180:0.8) .. (2,-1)
  .. controls +(0:0.8) and +(180:0.8) .. (4,0)
  .. controls +(180:0.8) and +(0:0.8) .. (2,1)
  .. controls +(180:0.8) and +(0:0.8) .. (0,0);
\fill[brown!20]
  (0,0,0)
  .. controls +(0:1) and +(0:-1) .. (2,1,0)
  .. controls +(0:1) and +(0:-1) .. (4,0,0)
  -- (4,0,-2)
  .. controls +(180:1) and +(180:-1) .. (2,1,-2)
  .. controls +(180:1) and +(180:-1) .. (0,0,-2)
  -- cycle;
\draw[blue,very thick,cap=round]
(0,0)
.. controls +(0:1) and +(0:-1) ..
(2,1)
.. controls +(0:1) and +(0:-1) ..
(4,0);
\draw[red,very thick,cap=round]
(0,0)
.. controls +(0:1) and +(0:-1) ..
(2,-1)
.. controls +(0:1) and +(0:-1) ..
(4,0);
\draw[purple,very thick](0,0,0)--(0,0,-2);
\draw[purple,very thick](4,0,0)--(4,0,-2);
\draw[blue,very thick,cap=round]
(0,0,-2)
.. controls +(0:1) and +(0:-1) ..
(2,1,-2)
.. controls +(0:1) and +(0:-1) ..
(4,0,-2);
\node[Gold!50!black,inner sep=1pt](vjl) at (0,2,-1) {\scriptsize$p_j$};
\draw[Gold!50!black,->](vjl)to[out=45,in=115](2,0,-2);
\end{scope}
\fill[violet!10]
(0,0)
  .. controls +(0:0.8) and +(180:0.8) .. (2,-1)
  .. controls +(0:0.8) and +(180:0.8) .. (4,0)
  .. controls +(180:0.8) and +(0:0.8) .. (2,1)
  .. controls +(180:0.8) and +(0:0.8) .. (0,0);
\node at (2,0) {\tkzglabel};
\draw[blue,very thick,cap=round]
(0,0,0)
.. controls +(0:1) and +(0:-1) ..
(2,1,0)
.. controls +(0:1) and +(0:-1) ..
(4,0,0);
\draw[red,very thick]
(0,0)
.. controls +(0:1) and +(0:-1) ..
(2,-1)
.. controls +(0:1) and +(0:-1) ..
(4,0);
\end{tkz}
}

\newcommand{\makediagR}[3][]{
\def\tkzsep{#2}
\def\tkzbdysep{#3}
\def\tkzglabel{\scriptsize$g$}
\begin{tkz}[
  x={({cos(\tkzangx)},{sin(\tkzangx)})},
  y={({cos(\tkzangy)},{sin(\tkzangy)})},
  z={({cos(\tkzangz)},{sin(\tkzangz)})},
  scale=0.5,
  xscale=-1,
  cap=round,
  #1
]
\begin{scope}[zshift=-\tkzbdysep]
\begin{scope}[zshift={-2-\tkzsep}]
\fill[gray!7,fading angle=30,path fading=east,opacity=0.75]
(0,0)
  .. controls +(0:0.8) and +(180:0.8) .. (2,-1)
  .. controls +(0:0.8) and +(180:0.8) .. (4,0)
  .. controls +(180:0.8) and +(0:0.8) .. (2,1)
  .. controls +(180:0.8) and +(0:0.8) .. (0,0);
\draw[red,very thick,cap=round]
(0,0)
.. controls +(0:1) and +(0:-1) ..
(2,-1)
.. controls +(0:1) and +(0:-1) ..
(4,0);
\node[violet!50!black,inner sep=1pt](vil) at (-1,0,-1) {\scriptsize$\fcj{v_i}$};
\draw[violet!50!black,->](vil)to[out=10,in=180-(90-115)](2,0,-1.5);
\fill[white]
  (0,0,0)
  .. controls +(0:1) and +(0:-1) .. (2,1,0)
  .. controls +(0:1) and +(0:-1) .. (4,0,0)
  -- (4,0,-\tkzfout)
  .. controls +(180:1) and +(180:-1) .. (2,1,-\tkzfout)
  .. controls +(180:1) and +(180:-1) .. (0,0,-\tkzfout)
  -- cycle;
\fill[gray!15,fading angle=30,path fading=east,opacity=1]
  (0,0,0)
  .. controls +(0:1) and +(0:-1) .. (2,1,0)
  .. controls +(0:1) and +(0:-1) .. (4,0,0)
  -- (4,0,-\tkzfout)
  .. controls +(180:1) and +(180:-1) .. (2,1,-\tkzfout)
  .. controls +(180:1) and +(180:-1) .. (0,0,-\tkzfout)
  -- cycle;
\node[violet,inner sep=1pt](vjl) at (0,2,-1) {\scriptsize$w_j$};
\draw[violet,->](vjl)to[out=30,in=115](2,0,-2);
\draw[blue,very thick,cap=round]
(0,0,0)
.. controls +(0:1) and +(0:-1) ..
(2,1,0)
.. controls +(0:1) and +(0:-1) ..
(4,0,0);
\draw[Green,very thick,fading angle=30,path fading=east](0,0,0)--(0,0,-\tkzfout);
\draw[Green,very thick,fading angle=30,path fading=east](4,0,0)--(4,0,-\tkzfout);
\end{scope}
\begin{scope}
\node[Brown!50!black,inner sep=1pt](vil) at (-1,0,-1) {\scriptsize$\fcj{p_i}$};
\draw[Brown!50!black,->](vil)to[out=10,in=180-(90-115)](2,0,-1.5);
\end{scope}
\fill[gray!7,fading angle=30,path fading=east,opacity=0.75]
(0,0)
  .. controls +(0:0.8) and +(180:0.8) .. (2,-1)
  .. controls +(0:0.8) and +(180:0.8) .. (4,0)
  .. controls +(180:0.8) and +(0:0.8) .. (2,1)
  .. controls +(180:0.8) and +(0:0.8) .. (0,0);
\fill[brown!20]
  (0,0,0)
  .. controls +(0:1) and +(0:-1) .. (2,1,0)
  .. controls +(0:1) and +(0:-1) .. (4,0,0)
  -- (4,0,-2)
  .. controls +(180:1) and +(180:-1) .. (2,1,-2)
  .. controls +(180:1) and +(180:-1) .. (0,0,-2)
  -- cycle;
\draw[blue,very thick,cap=round]
(0,0)
.. controls +(0:1) and +(0:-1) ..
(2,1)
.. controls +(0:1) and +(0:-1) ..
(4,0);
\draw[red,very thick,cap=round]
(0,0)
.. controls +(0:1) and +(0:-1) ..
(2,-1)
.. controls +(0:1) and +(0:-1) ..
(4,0);
\draw[purple,very thick](0,0,0)--(0,0,-2);
\draw[purple,very thick](4,0,0)--(4,0,-2);
\draw[blue,very thick,cap=round]
(0,0,-2)
.. controls +(0:1) and +(0:-1) ..
(2,1,-2)
.. controls +(0:1) and +(0:-1) ..
(4,0,-2);
\node[Gold!50!black,inner sep=1pt](vjl) at (0,2,-1) {\scriptsize$p_j$};
\draw[Gold!50!black,->](vjl)to[out=30,in=115](2,0,-2);
\end{scope}
\fill[violet!10]
(0,0)
  .. controls +(0:0.8) and +(180:0.8) .. (2,-1)
  .. controls +(0:0.8) and +(180:0.8) .. (4,0)
  .. controls +(180:0.8) and +(0:0.8) .. (2,1)
  .. controls +(180:0.8) and +(0:0.8) .. (0,0);
\node at (2,0) {\tkzglabel};
\draw[blue,very thick,cap=round]
(0,0,0)
.. controls +(0:1) and +(0:-1) ..
(2,1,0)
.. controls +(0:1) and +(0:-1) ..
(4,0,0);
\draw[red,very thick]
(0,0)
.. controls +(0:1) and +(0:-1) ..
(2,-1)
.. controls +(0:1) and +(0:-1) ..
(4,0);
\end{tkz}}

\makeatletter

\renewcommand{\l@section}{\@tocline{1}{\medskipamount}{0pt}{1.75pc}{}}
\renewcommand{\l@subsection}{\@tocline{2}{0pt}{0pt}{2.5pc}{}}
\renewcommand{\l@subsubsection}{\@tocline{3}{0pt}{0pt}{3.5pc}{}}

\makeatother

\title{Unitary TQFTs, unitary disk-like {\large \lowercase{$n$}}-categories, and higher Hilbert spaces}

\author{Greyson Wesley}

\date{\today}

\begin{document}

\begin{abstract}
We introduce the notion of a unitary disk-like $n$-category, which is a disk-like $n$-category equipped with a reflection structure and a sphere trace inducing positive-definite pairings. Since a finite unitary disk-like $n$-category is defined to be the local field data of a fully extended $(n+1)$D unitary TQFT, we propose a complete finite unitary disk-like $n$-category as the definition of a finite $(n+1)$-Hilbert space for all $n$. For $n=1$ and $n=2$ we verify this proposal, proving that complete finite unitary disk-like 1- and 2-categories are isometrically equivalent to finite 2- and 3-Hilbert spaces respectively, and that these equivalences are functorial. For $n=1$ we recover a unitary refinement of Schommer-Pries' classification of oriented $(1+1)$D TQFTs in terms of $\Hstar$-Morita equivalence classes of $\Hstar$-algebras, and for $n=2$ we categorify this to classify oriented $(2+1)$D unitary TQFTs by $\Hstar$-Morita equivalence classes of $\Hstar$-multifusion categories. Along the way, we investigate fully incomplete 3-Hilbert spaces, prove a strictification result for pivotal dagger 2-categories, and define functors and higher transformations between disk-like $n$-categories.
\end{abstract}

\begingroup
\maketitle
\endgroup

\begingroup\small
\tableofcontents
\endgroup

\section{Introduction}

There are several ways to axiomatize a (fully extended) topological quantum field theory (TQFT). One way is to axiomatize its manifold invariants directly, as in the fully extended version of the Atiyah--Segal axioms, where a theory is a symmetric monoidal functor out of a suitable higher category of bordisms \cite{Ati88,Seg04,BD95,Lurie09}. Another way is to axiomatize the fields and the path-integral boundary conditions from which these invariants are assembled, which is the viewpoint of Walker \cite{W06,W21} and Morrison--Walker \cite{MW12}, and is the one we take in this article. A \defn{disk-like $n$-category} in the sense of \cite{MW12} is an axiomatization of the local structure of these fields: it assigns to each $0\leq k\leq n$ and each $k$-ball $W$ a set of \defn{fields} on $W$, with boundary maps, gluing maps, and pullback maps recording how fields restrict, compose, and thicken, together with a subspace of local relations in the top spatial dimension for linearization. Throughout, the two examples of fields on manifolds to keep in mind are string diagrams and maps to a fixed topological space. This article develops the unitary theory of these structures and proposes a definition of finite $(n{+}1)$-Hilbert space for all $n$. More specifically, we propose the \emph{complete} finite unitary disk-like $n$-categories---those equivalent to disk-like $n$-categories of representations---as a definition of finite $(n{+}1)$-Hilbert space, and we verify this proposal against the existing algebraic definitions for $n=1$ and $n=2$. Here a \defn{unitary disk-like $n$-category} is a disk-like $n$-category equipped with a reflection structure $\fcj{\,\cdot\,}$, which can roughly be thought of as a disk-like analog of a dagger structure, together with a \defn{sphere trace} $\psn$ on the fields on $S^n$ whose induced sesquilinear pairings on $n$-balls are positive-definite. The sphere trace is the local manifestation of the path integral on $D^{n+1}$, and a finite unitary disk-like $n$-category is by design the local data of a fully extended $(n+1)$D unitary TQFT.

Given a disk-like $n$-category $\C$, \defn{skeletonization} extracts an ordinary (weak) $n$-category $\eX_\C$ whose $k$-morphisms are the fields on the standard $k$-ball $D^k$. In the other direction, the \defn{string diagram construction} produces from finite $(n+1)$-Hilbert spaces a disk-like $n$-category $\C^\eX$ whose fields are $\eX$-labeled string diagrams. For $n=1$ these are mutually inverse isometric equivalences between finite unitary disk-like 1-categories and finite pre-2-Hilbert spaces (\cref{TWODUCH}), and for $n=2$ between finite unitary disk-like 2-categories and finite proto-3-Hilbert spaces (\cref{THREEDUCH}). In both cases this equivalence restricts to one between the complete objects on one side and the (complete) higher Hilbert spaces on the other. These equivalences are moreover functorial in that disk-like functors and their higher analogs---which we call \defn{transfors}---assemble into a disk-like $n$-category $\eFun(\C{\to}\D)$ that recovers (skeletonizes to) the higher Hilbert spaces $\Hom(\eX\to\eY)$ of morphisms (\cref{HOMEQUIV}).

These results have the following applications. For $n=1$ we obtain a unitary refinement of Schommer-Pries' classification of oriented $(1+1)$D TQFTs, with $\Hstar$-Morita equivalence classes of $\Hstar$-algebras in place of Morita classes of standard separable Frobenius algebras (\cref{TWODCLASS}); this refinement is strict, as there are Morita equivalent algebras that are not $\Hstar$-Morita equivalent (\cref{rmk:n1failure}). For $n=2$ we categorify this to classify (fully extended) oriented $(2+1)$D unitary TQFTs by $\Hstar$-Morita equivalence classes of $\Hstar$-multifusion categories (\cref{THREEDCLASS}), where the refinement turns out not to be strict due to the existence of a higher categorical version of polar decomposition (\cref{n2nonfailure}). Along the way, we equip the disk-like functor categories with a canonical sphere trace making them unitary (\cref{subsec:eFun-sphere-weight}) and we prove that every pivotal dagger 2-category is equivalent to a strictly pivotal dagger 2-category with strict underlying 2-category in a way that preserves the unitary pivotal structure (\cref{strictification-sec}).

We suggest that the reader more interested in the concrete story skip the general-$n$ section \cref{sec:n+1D} on a first reading in favor of the $(1+1)$D and $(2+1)$D sections \cref{sec:1+1D} and \cref{sec:2+1D}; the former is entirely self-contained and the latter is mostly self-contained as well.

\subsection{Disk-like \texorpdfstring{$n$}{n}-categories and \texorpdfstring{$(n+\e)$}{(n+e)}D TQFTs}

Among other reasons, fields can be motivated by axiomatizing the expected behavior of boundary conditions for path integrals. Roughly, the idea is that an $(n+1)$D TQFT assigns to each $k$-manifold $W$ for $0\leq k\leq n$ a collection $\fld[k]{\C}{W}$ of \defn{fields} on $W$ together with rules for how these fields combine when their underlying manifolds are glued together. The higher category encodes the \emph{local} data of these assignments, which in the topological setting just means the behavior of fields on $k$-balls for $0\leq k\leq n$, together with their behavior when balls are glued to form another ball. In dimensions $n$ and lower, no finiteness is needed: a colimit construction \cite[\S6.3]{MW12} uniquely extends the local (ball-level) field data to all $k$-manifolds for $0\leq k\leq n$ (see \cref{sec:dlncats} below). In dimension $n+1$, however, a notion of \emph{finiteness} of the field data is required to allow it to be patched together under gluing. This finiteness allows an $(n+1)$-manifold to be cut into pieces, so the path integral on the pieces can be patched together by summing over a finite basis of intermediate states. It is this finite summation under gluing that reconstructs the path integral from the local data, as formalized in the Walker Extension Theorem \cite{W06,W21} (see \cref{WalkerTheorem} below). 

In the pivotal setting one can freely rotate the boundary of a morphism, thereby exchanging which of its parts are ``incoming'' and ``outgoing.'' This suggests looking for a model of a pivotal $n$-category in which source and target are not singled out at all. In \cite{MW12}, Morrison--Walker define a \defn{disk-like $n$-category} as a model in which this pivotal symmetry is built in from the start: for $0\leq k\leq n$, its $k$-morphisms---which we will call \defn{$k$-fields} in this article---may take \emph{any} combinatorial shape, so long as they are homeomorphic to the standard $k$-ball. For instance, a $k$-field could take the shape of a bigon/globe, a simplex, a square, an opetope, and so on. So, if we want to extract an ordinary/traditional $n$-category from a disk-like $n$-category $\C$, we have to \emph{choose} a shape for our morphisms, which we refer to as \defn{skeletonizing}.

Briefly, a disk-like $n$-category $\C$ assigns to each $k$-ball $W$ a set $\fld[k]{\C}{W}$ of \defn{$k$-fields} of shape $W$ (the morphisms), together with boundary maps (combined source and target), gluing maps along $(k-1)$-balls (composition), pullbacks along pinched product maps (identities), and, in dimension $n$, a subspace of \defn{local relations} $\fldlU{\C}{X}[c]\subset\bC\{\fld[n]{\C}{X}[c]\}$ for boundary conditions $c$ of $n$-balls $X$. The vector space $\A_\C(X,c)$ is then given by quotienting these local relations $\fldlU{\C}{X}[c]$ out of the linearized pure $n$-fields $\bC\{\fld[n]{\C}{X}[c]\}$. All of this data is required to be compatible in essentially every way that can be written down; see \cref{sec:full-axiomatic-definition} below. We will usually require $\C$ to be finite (see \cref{psnF} below). The two main classes of examples of disk-like $n$-categories are string diagrams labeled by some ``pivotal $n$-category''---for example, pivotal monoidal categories and more generally pivotal 2-categories, premodular categories \cite{W14premodular}, and spherical fusion 2-categories \cite{DR18} giving rise respectively to Turaev--Viro, Crane--Yetter, and Douglas--Reutter theories---and maps to a fixed target space $T$, giving Dijkgraaf--Witten theories (take $T\coloneq\rB G$ for a finite group $G$) and more generally $\sigma$-model TQFTs. Examples outside of these two main classes include bordism categories \cite[Exm. 6.2.6]{MW12} and twisted maps to a space \cite[Exm. 6.2.3]{MW12}. There are also more exotic examples; see \cite{walkeruscslides,W09fields}. 
While the axioms are long, they simply record properties visibly shared by string diagrams and by maps to a space; we recommend keeping these two classes in mind throughout.

Morrison--Walker also introduce disk-like $(\infty,n)$-categories, which they call \defn{$\rA_\infty$ disk-like $n$-categories}, and several examples are provided in \cite[\S6.2]{MW12}. While their axioms differ only slightly from those of ordinary disk-like $n$-categories (spelled out in \cref{sec:dlncats} below), we do not explore them in this paper.

By the \defn{colimit construction} \cite[\S6.3]{MW12} (see also \cite{morrison2011higher}), any disk-like $n$-category $\C$ uniquely extends to a gadget $\undC$ that gives spaces of fields on general $k$-manifolds for each $0\leq k\leq n$. While passing from an $(n+1)$D TQFT $\Z$ to its corresponding disk-like $n$-category $\C$ forgets the behavior of the data of the path integral $\Z$ in the top dimension (see \cref{defUnitaryTQFT} below), $\C$ also contains the data of the action of orientation-preserving homeomorphisms of $n$-balls $\varphi\colon X\to X'$ inducing linear maps $\varphi_*\colon\undfldl[n]{\C}{X}[c]\to\undfldl[n]{\C}{X'}[\varphi_* c]$, which in turn determine the value $\Z(M_\varphi)$ of the path integral on the mapping cylinder $M_\varphi$ of $\varphi$. This motivates calling $\undC$ an \defn{$(n+\e)$D TQFT} (also called a \defn{once-categorified $n$D TQFT}). 

One of the key advantages of Kevin Walker's approach to TQFT is that the higher categorical structure is automatically present down to points: for every $(n-\ell)$-manifold $Q$, the theory produces a disk-like $\ell$-category $\A_\C(Q)$ called the \defn{skein disk-like $\ell$-category} on $Q$, which may also be called the \defn{cylinder $\ell$-category} over $Q$, \defn{dimensional reduction} along $Q$, or \defn{compactification} by $Q$, whose objects (0-fields) are $(n-\ell)$-fields of $\C$ on $Q$, whose 1-fields on 1-balls $J$ are $(n-\ell+1)$-fields of $\C$ on $Q\times J$, whose 2-fields on 2-balls $X$ are $(n-\ell+2)$-fields of $\C$ on $Q\times X$, and so on. When $Q$ is a ball, $\A_\C(Q)$ is the disk-like analog of a hom $\ell$-category in an ordinary/traditional $n$-category. 

We warn the reader familiar with the Atiyah--Segal axioms \cite{Ati88,Seg04} that this axiomatization of TQFT differs in that we remember the \emph{predual} $\A$ to the usual functorial TQFT $\Z$-invariants, which are representations of (modules for) the $\A$-invariant. As the usual functorial TQFT invariant $\Z=\Z_\C$ is defined by $\Z_\C(W)\coloneq\eRep(\A_\C(W))$, it is precisely representations of skein disk-like $n$-categories for unitary TQFTs that our main results identify with finite $(n+1)$-Hilbert spaces for $n=1,2$. 

\subsection{Unitarity}
Absent from this picture has been any systematic treatment of \emph{unitarity}, a gap that persists largely because the structure of unitarity in higher categories has historically been poorly understood.
Motivated by this gap, there has been recent progress in constructing finite $n$-Hilbert spaces inductively up to $n = 3$, as developed in \cite{3Hilb}, which recover finite 2-Hilbert spaces in the sense of Baez \cite{Baez97} and propose a definition and theory of 3-Hilbert spaces. In parallel, the article \cite{FHJF24} defines \defn{dagger $(\infty, n)$-categories with unitary duality} (also called \defn{$\PL(n)$-dagger categories}) and proposes that $n$-dimensional unitary TQFTs should be formalized as symmetric monoidal functors valued in such structures defined on suitable higher bordism categories. For $n=1$, the coherent notion of a dagger structure coincides with the existing algebraic definition of dagger category \cite{SS23}, and for $n=2$ it is expected to coincide with that of a dagger $2$-category \cite{CHPJP22}, which we use in this article. See \cite{FH21} for the fully extended invertible case, and see \cite{M26} and \cite{MS26} for reflection positivity in defect and in once-extended theories respectively. Nonetheless, a complete understanding of the ``correct'' target category for $(n+1)$D unitary TQFTs in full generality remains an open problem.

We use the prefixes ``proto-'' and ``pre-'' for the versions of these notions before completion: a ``pre-'' higher Hilbert space is complete in every degree below the top, and a ``proto-'' one is subject to no completeness condition at all. A \defn{pre-2-Hilbert space} is a unitary category equipped with a faithful positive trace, and a \defn{$2$-Hilbert space} is one that is moreover complete, meaning that it has orthogonal direct sums and its projections split. Likewise a \defn{proto-$3$-Hilbert space} is a pivotal $\Cstar$-$2$-category equipped with a spherical weight, and a \defn{pre-$3$-Hilbert space} is a locally complete one.
We warn the reader that our conventions here are slightly non-standard: finiteness is often built into the definitions of these notions in the literature, whereas for us it is an extra condition (see \cref{def:finiteunitarycategory,def:pre-2-Hilbert space,def:3hilb-spaces}). Accordingly, and for the reader's convenience, we always write ``finite'' explicitly in the hypotheses and conclusions of our results.

As in \cite{3Hilb}, we will use the term \defn{unitary} generically to mean not only the presence of dagger structures at multiple levels together with positivity conditions, but also the presence of a trace providing positive-definite inner products. A guiding philosophy of this article, shared with \cite{3Hilb}, is that unitarity should be \emph{manifestly built in from the start}, rather than imposed as extra conditions after the fact. We adopt this approach by insisting that the local data of a unitary TQFT be unitary, and then prove that unitarity is preserved when extending the local data to a theory on all manifolds. In other words, we axiomatize the local data---the data visible on balls of dimension at most $n$---in such a way that positivity is already present at the local level, and then reconstruct the full path integral from this local data.

Thus, since ``suitably finite pivotal $n$-categories'' correspond to $(n+1)$D TQFTs, one would expect the local part of an $(n+1)$D \emph{unitary} TQFT to model the correct notion of a suitably finite pivotal $n$-category equipped with a positive trace, where this trace provides ``unitary duality in dimension $n+1$,'' and can be thought of as the path integral on the $(n+1)$-ball $D^{n+1}$. The purpose of this article is to make this expectation precise: we propose complete finite unitary disk-like $n$-categories as a definition of finite $(n+1)$-Hilbert spaces for \emph{all} $n$, and for $n = 1$ and $n = 2$ we verify the proposal by constructing an explicit correspondence between unitary disk-like $n$-categories and the existing notions of higher Hilbert spaces.

To carry out this program, we will build unitarity into the notion of \defn{disk-like $n$-category} in the sense of Morrison--Walker \cite{MW12}. A disk-like $n$-category is defined to capture the local data of an $(n+\e)$D TQFT, that is, of an $(n+1)$D TQFT \emph{without} the top-dimensional (path integral) data, but that retains the behavior of the theory on mapping cylinders. Since unitarity inherently requires positive-definite inner products compatible with all the data of the theory, if we are to adapt the disk-like framework to the unitary case then we need to bring the path integral into the fold. 
An $(n+1)$D unitary TQFT (\cref{defUnitaryTQFT} below) is an $(n+\e)$D TQFT $\undC$ equipped with a \defn{path integral} (also called a \defn{functional integral} or \defn{partition function}) $\Z_\C(M)\colon \undfldl[n]{\C}{\partial M}\to\bC$ for each $(n+1)$-manifold $M$ satisfying reflection positivity, compatibility with gluing, and homeomorphism invariance.

\subsection{Unitary disk-like \texorpdfstring{$n$}{n}-categories}
With these conventions and motivation in hand, we now describe the additional data making a disk-like $n$-category ``unitary.'' The first piece of additional structure is a \defn{reflection structure} $\fcj{\,\cdot\,}$ on $\C$, which consists of strictly involutive bijections $\fcj{\,\cdot\,}\colon\fld[k]{\C}{W}\to\fld[k]{\C}{\orev{W}}$ for $0\leq k<n$ and conjugate-linear isomorphisms $\fcj{\,\cdot\,}\colon\fldl[n]{\C}{X}[c]\to\fldl[n]{\C}{\orev{X}}[\fcj{c}]$ that are compatible with homeomorphisms, gluing, boundaries, and pinched product maps.

The second piece of additional structure is a \defn{sphere trace} (or \defn{trace}) $\psn\colon\undfldl{\C}{S^n}\to\bC$, which is a linear functional on the space of fields on the $n$-sphere satisfying \ref{psnP} below. The sphere trace is the disk-like analog of a positive trace (for $n=1$) or a spherical weight (for $n=2$) on an ordinary $n$-category; it should be thought of as the ``local manifestation of the path integral'' in the sense that it is the path integral on $D^{n+1}$. In skein-theoretic contexts this functional is sometimes called an \defn{evaluation map}. The sphere trace induces sesquilinear pairings $\bkt{f}{g}_{X,c}\coloneq\psn(\fcj{f}\blt_c g)$ on the vector spaces $\fldl{\C}{X}[c]$ of fields on $n$-balls, where $\blt_c$ denotes gluing along $\partial X$ (see \ref{DGk} below) so that $\fcj f\blt_c g$ lives on $\orev X\cup_{\partial X}X\cong S^n$, and we demand that these pairings be positive-definite.
Our sphere trace plays the same conceptual role as the so-called $S^n$-functional of \cite{Z24,LZ26}, and the importance of such a functional was first demonstrated in \cite{W06} and again in \cite{W21}. Our work develops the unitary aspects directly within Walker's general framework.

\begin{restatable}[Unitary disk-like $n$-category]{definitionalpha}{defUnitaryDiskLikeNCategory}
    \label{def:intro-UDL}
    A \defn{unitary disk-like $n$-category} $(\C,\fcj{\,\cdot\,},\psn)$ is a disk-like $n$-category $\C$ equipped with a reflection structure $\fcj{\,\cdot\,}$ and a linear functional $\psn\colon\undfldl{\C}{S^n}\to\bC$ satisfying the following condition.
    \begin{lst}[font=\upshape,leftmargin=0.425in]
        \item[($\psn$P)]
        (\emph{Positivity})
        For each $n$-ball $X$ and $c\in\undfld[n-1]{\C}{\partial X}$, the induced pairing
        \begin{equation}\label{intro-disk-pairings}
            \bkt{f}{g}_{X,c}\coloneq\psn(\fcj f\blt_c g)
        \end{equation}
        on $\fldl[n]{\C}{X}[c]$ is positive-definite.
    \end{lst}
\end{restatable}

Positivity implies $\psn(\fcj\alpha)=\overline{\psn(\alpha)}$ for all $\alpha\in\undfldl[n]{\C}{S^n}$, where $\fcj\alpha$ is computed by transporting from $S^n$ back to $S^n$ by the unique-up-to-isotopy orientation-reversing self-homeomorphism of $S^n$; see \cref{reflection and functionals} and \cref{prop:psnP implies}.

We note that a significant portion of this article could be done similarly in the nondegenerate case. In doing so, we could define a (faithfully) \defn{traced} or \defn{tracial} (or less descriptively, \defn{Calabi--Yau}) \defn{disk-like $n$-category} to be the nondegenerate analog of a unitary disk-like $n$-category, but we only pursue the unitary case in this paper.

\subsection{Outline and main results}

We now state the main results of this article. The definition of a (finite) unitary disk-like $n$-category makes sense for all $n$, and we propose the complete ones as a definition of (finite) $(n+1)$-Hilbert spaces; the main results verify this proposal for $n=1$ and $n=2$ by identifying (finite) unitary disk-like $n$-categories with the (finite) higher Hilbert spaces already in the literature, and the complete ones with the complete ones. We verify the proposal only for $n=1$ and $n=2$ because finite $4$-Hilbert spaces have not yet been defined, not because the constructions below break down in higher dimensions.

A unitary disk-like $n$-category is \defn{complete} when there is an isometric weak equivalence between it and its disk-like $n$-category of \defn{representations}, defined in \cref{def:completeDL1cat} for $n=1$ and \cref{def:completeDL2cat} for $n=2$. 

Given a disk-like 1-category $\C$, we write $\cX_\eC$ for the \defn{skeletonization} of $\C$ (\cref{skein 1-categories} below), which is the (ordinary/traditional) dagger 1-category whose objects are 0-fields on $D^0=\pt$ and whose morphisms are 1-fields on the standard 1-ball $D^1$. When $\C$ is unitary with sphere trace $\psn$, we can construct a pre-2-Hilbert space trace (see \cref{2hilbtrace} below) $\Tr^{\cX_\eC}$ on $\cX_\eC$ given roughly by closing up morphisms on $D^1$ to $S^1$ and applying $\psn$.

For the other direction, if $\cX$ is a dagger 1-category, then $\C^\cX$ is the disk-like 1-category of $\cX$-string diagrams. When $\cX$ is equipped with a pre-2-Hilbert space trace, we obtain a sphere trace $\psn^{\C^\cX}$ roughly by applying $\Tr^\cX$ to string diagrams on the circle $S^1$.
 
\begin{restatable}{theoremalpha}{TWODUCH}
    \label{TWODUCH}
Skeletonization $\C\mapsto\cX_\C$ and the string diagram construction $\cX\mapsto\C^\cX$ are mutually inverse isometric equivalences between finite unitary disk-like $1$-categories and finite pre-$2$-Hilbert spaces. Moreover, these restrict to isometric equivalences between complete finite unitary disk-like $1$-categories and finite $2$-Hilbert spaces.
\end{restatable}

For a unitary disk-like 2-category $\C$, its \defn{skeletonization} $\fX_\C$ (\cref{sec:skein 2-category}) is defined similarly to its $n=1$ analog---objects are fields on $D^0$, 1-morphisms are fields on $D^1$, and 2-morphisms are fields on $D^2$---and admits a canonical pivotal dagger structure (see \cref{def:pivotal 2-category}). Its 1-morphisms have unitary adjoints given by (spatially) reflecting 1-fields about the midpoint $D^0\hookrightarrow D^1$, while the dagger on 2-morphisms is given by (spatially) reflecting about the equator $D^1\hookrightarrow D^2$. If $\C$ is finite unitary, then $\fX_\C$ is a finite proto-$3$-Hilbert space whose spherical weight is obtained by closing up $2$-fields on $D^2$ to $S^2$ and applying $\psn$.

To go the other direction, we write $\C^\fX$ for the disk-like 2-category of $\fX$-string diagrams.
\cref{THREEDUCH} categorifies \cref{TWODUCH}.

\begin{restatable}{theoremalpha}{THREEDUCH}
    \label{THREEDUCH}
    Skeletonization $\C\mapsto\fX_\C$ and the string diagram construction $\fX\mapsto\C^\fX$ are mutually inverse isometric equivalences between finite unitary disk-like 2-categories and finite proto-3-Hilbert spaces. Moreover, these restrict to isometric equivalences between complete finite unitary disk-like $2$-categories and finite $3$-Hilbert spaces.
\end{restatable}

A \defn{disk-like functor} $\C\to\D$ is a collection of maps on fields commuting with the boundary, gluing, and pullback structure, and an \defn{$(n,\ell)$-transfor} is the $\ell$-fold higher analog of a natural transformation or modification. Unfortunately, this strict notion is too rigid to construct functors between the main classes of examples, such as from a disk-like $n$-category of string diagrams to one given by maps to a fixed space. To resolve this, we introduce \defn{unrestricted disk-like functors}, which allow for the choice of parameterized regular neighborhoods as extra data. This framework is discussed for $n=1$ in \cref{sec:1+1D} and the full details are in \cref{sec:unnrestricted} for general $n$. Nevertheless, this is not necessary for the case of string diagrams, and we can always pass to the string diagram construction. 

In proving the completeness assertions of the above theorems, we prove that the above equivalences are functorial in the following sense. Write $\eFun(\C{\to}\D)$ for the disk-like $n$-category whose $0$-fields are disk-like functors $\C\to\D$ and whose $\ell$-fields are disk-like $(n,\ell)$-transfors, and $\Hom(\cX\to\cY)$ for the higher Hilbert space of morphisms of $(n+1)$-Hilbert spaces (dagger functors for $n=1$ and UAF-preserving dagger functors for $n=2$) and their higher transformations.

\begin{restatable}[Functor category equivalence]{theoremalpha}{HOMEQUIV}
    \label{HOMEQUIV}
For finite pre-2-Hilbert spaces $\cX$ and $\cY$ and finite proto-3-Hilbert spaces $\fX$ and $\fY$, there are canonical isometric equivalences $\cX_{\eFun(\eC^\cX\to\eC^\cY)}\cong^\dag\Hom(\cX\to\cY)$ and
$\fX_{\eFun(\C^\fX\to\C^\fY)}\cong^\dag\Hom(\fX\to\fY)$ and canonical isometric weak equivalences
$\eC^{\Hom(\cX\to\cY)}\cong^\dag\eFun(\eC^\cX\to\eC^\cY)$ and
$\C^{\Hom(\fX\to\fY)}\cong^\dag\eFun(\C^\fX\to\C^\fY)$.
\end{restatable}

The string diagram equivalences follow from the observation that a string diagram labeled by a functor category is the same thing as an overlay pattern in the sense of the overlay graphical calculus of \cite[\S2.4]{CP22}. Briefly, such a pattern $\xi$ on an $\ell$-ball $Q$ acts on a string diagram $\zeta$ on a $k$-ball $W$ to give the product diagram $\zeta\times\xi$ on $W\times Q$: its strata are the products of a stratum of $\zeta$ with a stratum of $\xi$, with each being labeled by the component of its $\xi$-label at the label of the stratum of $\zeta$ it lies over. Combining these string diagram equivalences with \cref{TWODUCH} and \cref{THREEDUCH} implies the assertions regarding the skeletonizations. In particular, $\eRep(\eC)\coloneq\eFun(\eC^{\cX_\eC^\op}{\to}\eHilb)$ skeletonizes to $\cX_{\eRep(\eC)}\cong^\dag\Rep(\cX_\eC^\op)$ for $n=1$ and $\eRep(\C)\coloneq\eFun(\C^{\fX_\C^{1\op}}{\to}2\eHilb)$ skeletonizes to $\fX_{\eRep(\C)}\cong^\dag\Rep(\fX_\C^{1\op})$ for $n=2$, which implies the completeness assertions in \cref{TWODUCH} and \cref{THREEDUCH} above. By \cref{complete 2-Hilb equiv to its reps} and \cref{3hilb equiv to its reps} respectively, this means that for $n=1,2$, a unitary disk-like $n$-category is complete precisely when its skeletonization is complete. 

Recall that an \defn{$\Hstar$-algebra} is a finite-dimensional $\Cstar$-algebra $A$ equipped with a faithful positive trace $\Tr_A$. Thus an $\Hstar$-algebra $A$ is canonically a Hilbert space $L^2A$ with inner product $\bkt{a}{b}_{L^2A}\coloneq\Tr_A(a^* b)$. Two $\Hstar$-algebras $A$ and $B$ are \defn{$\Hstar$-Morita equivalent} (or \defn{isometrically Morita equivalent}) when $A$ is isometrically isomorphic (\cref{H*algisomeq}) to the commutant (\cref{H*AlgComm}) of $B$ with respect to the action of $B$ on some faithful module \cite[Exm. 3.7.9]{UQSL}; this is \cref{def:HstarMorita} below.

\begin{restatable}[Classification of oriented $(1+1)$D unitary TQFTs]{theoremalpha}{TWODCLASS}
    \label{TWODCLASS}
    Complete finite unitary disk-like $1$-categories (equivalently, fully extended oriented $(1+1)$D unitary TQFTs) are classified up to isometric weak equivalence by $\Hstar$-Morita equivalence classes of $\Hstar$-algebras.
\end{restatable}

As $\Hstar$-algebras are equivalently standard separable unitary Frobenius algebras in $\Hilb$ \cite[§I.3.6]{UQSL}, \cref{TWODCLASS} is the unitary refinement of the classification of fully extended oriented 2D TQFTs in Schommer-Pries' thesis \cite{SP09}. And this is genuinely a refinement: there are $\Hstar$-algebras that are Morita equivalent but not $\Hstar$-Morita equivalent; see \cref{rmk:n1failure} below.

In \cref{sec:2+1D}, we categorify \cref{TWODCLASS}. The categorification of an $\Hstar$-algebra is an \defn{$\Hstar$-multifusion category}, that is, a unitary multifusion category $\cA$ equipped with a unitary dual functor (UDF) \cite{Pen20} and a spherical weight $\psi^\cA$ \cite[Defn. 3.1]{3Hilb}. An $\Hstar$-multifusion category is canonically a 2-Hilbert space $L^2\cA$ with trace $\Tr^{L^2\cA}$ induced by the spherical weight (see \cite[Rmk. 3.2]{3Hilb}). In \cref{classificationHstarmFC}, we define \defn{isometric equivalence} of $\Hstar$-multifusion categories (\cref{def:HmFC-isometric-equivalence}), the \defn{commutant} of an $\Hstar$-multifusion category with respect to a module (\cref{def:commutant}), and \defn{$\Hstar$-Morita equivalence} of $\Hstar$-multifusion categories (\cref{def:hstar-morita}). The aforementioned notion of isometric equivalence of $\Hstar$-multifusion categories is the ``correct'' categorification of isometric isomorphism in the following sense. At $n=1$, a $*$-isomorphism of $\Hstar$-algebras is isometric if it preserves the trace, or equivalently if it is an isometry of the underlying Hilbert spaces. At $n=2$, there are a priori \emph{three} different ways to ask that a dagger monoidal equivalence of $\Hstar$-multifusion categories $\alpha\colon\cA\to\cB$ preserve the unitary structure: preservation of the spherical weight, the induced trace $\Tr^{L^2\cA}$, and the \defn{$\Hilb$-valued trace} $\mathrm{TR}^\cA\coloneq\cA(1_\cA \to -)$ of \cite[Rmk. 3.30]{bases}. We show in \cref{lem:isometric-tfae} that these three conditions coincide, and thus we say $\alpha$ is an \defn{isometric equivalence} of $\Hstar$-multifusion categories if any (hence all) of the above conditions hold.

\begin{restatable}[Classification of oriented $(2+1)$D unitary TQFTs]{theoremalpha}{THREEDCLASS}
\label{THREEDCLASS}
Complete finite unitary disk-like $2$-categories (equivalently, fully extended oriented $(2+1)$D unitary TQFTs) are classified up to isometric weak equivalence by $\Hstar$-Morita equivalence classes of $\Hstar$-multifusion categories, or equivalently ordinary Morita equivalence classes.
\end{restatable}

By \cite[\S6]{3Hilb}, $\Hstar$-multifusion categories are precisely the $\Hstar$-algebras---the standard separable unitary Frobenius algebras---internal to the symmetric monoidal $\Cstar$-$2$-category $2\Hilb$ of finite 2-Hilbert spaces. (By ``internal to'' we mean in the sense of \cite[Defn. 6.8]{3Hilb}; for a non-unitary counterpart, see \cite{Xu26}.) Thus \cref{THREEDCLASS} is a direct categorification of \cref{TWODCLASS}.

The $n=1$ and $n=2$ classifications differ in one remarkable feature. For $\Hstar$-algebras the refinement is strict, as just mentioned, whereas for $\Hstar$-multifusion categories every Morita equivalence can be augmented to an $\Hstar$-Morita equivalence (see \cref{n2nonfailure} below).
 
The $n=2$ story appeals to a strictification result that may be of independent interest. We define a \defn{pivotal dagger 2-category} $\fX$ to be a rigid dagger 2-category (\cref{dagger 2-category}) equipped with choices of adjoint data for each 1-morphism for which the induced adjoint functor $\vee\colon\fX\to\fX^{1\op,2\op}$ is a $\dag$-functor; such a choice of adjoint data is called a \defn{unitary adjoint functor} (UAF) (see \cref{def:uaf} below). Such a choice induces a canonical pivotal structure $\phi\colon\id\Rightarrow\vee\circ\vee$ (see \cref{canpivstr}). A pivotal dagger structure on $\fX$ is called \defn{strict} when $\phi_X=\id_X$ for every 1-morphism $X$ in $\fX$.
 
\begin{restatable}{theoremalpha}{strictification}
\label{strictification}
Every pivotal dagger 2-category is (UAF-preservingly) equivalent to a strict dagger 2-category equipped with a UAF whose induced pivotal structure is strict.
\end{restatable}

We expect that if one were to instead define the evaluation map in \cref{eval-facts-2} below by choosing a fixed parenthesization scheme (e.g., left-nested parenthesizations), then the constructions in \cref{sec:2+1D} would imply this strictification directly. 
 
In \cref{sec:n+1D} we outline the main definitions and tools to develop the theory for general $n$. This means introducing disk-like $n$-categories, reflection structures, the dimensional reduction/skein construction, disk-like functors and transfors, unitary disk-like $n$-categories, and the path integral. The technical work behind the extension from local data to a theory on all manifolds---including gluing lemmas, finiteness and orthogonality results, and Walker's handle construction of the path integral---is collected in \cref{sec:gluing-appendix}.

Extending the local data of a finite unitary disk-like $n$-category to an $(n+1)$D unitary TQFT follows from a sequence of gluing lemmas, a finiteness-and-orthogonality analysis of the resulting skein modules, and Walker's construction of the path integral by handle attachment. The argument is similar to that of Walker \cite[\S6]{W06}, \cite[\S4]{W21}, and it is not hard to adapt it to the unitary case; this is what we do in \cref{sec:gluing-appendix}, where we work out details that have not appeared in the literature and recast the results in the disk-like framework. Aside from providing the details, our contribution is that Walker's handle construction of the path integral propagates reflection positivity; the upshot is the Unitary Walker Extension Theorem (\cref{UnitaryWalkerTheorem}), which asserts that the local data of a unitary disk-like $n$-category suffices to reconstruct the full $(n+1)$D unitary TQFT; the inductive proof in \cref{cstr:path-integral-via-handle-decomp} shows that unitarity is preserved under gluing thanks to the finiteness of the local data.

\begin{restatable}[Unitary Walker Extension Theorem]{corollaryalpha}{UnitaryWalkerTheorem}
\label{UnitaryWalkerTheorem}
A finite unitary disk-like $n$-category $\C$ extends to a unique $(n+1)$D unitary TQFT $(\C,\Z_\C)$ with $\Z_\C(D^{n+1})=\psn^\C$.
\end{restatable}

The arguments of \cite{W06,W21} assume (pre)semisimplicity, which in our setting is supplied by unitarity and finiteness. In upcoming work, Reutter and Walker construct the path integral by a rather different and more abstract argument that drops the (pre)semisimplicity assumption \cite{RW,R20}. We do not use that argument here; the older argument works naturally in the unitary framework where the constants/renormalization factors matter (see \cref{subsec:eFun-sphere-weight} and \cref{sec:gluing-appendix}).

In \cref{sec:1+1D} we specialize to $n=1$. We construct the skeletonization functor $\eC\mapsto\cX_\eC$ and the string diagram construction $\cX\mapsto\eC^\cX$, prove they are mutually inverse isometric equivalences (\cref{TWODUCH}), establish the functor category equivalence (\cref{HOMEQUIV} for $n=1$), and classify $(1+1)$D unitary TQFTs by $\Hstar$-Morita equivalence classes of $\Hstar$-algebras (\cref{TWODCLASS}).
 
In \cref{sec:2+1D} we specialize to $n=2$. We construct the skeletonization $\C\mapsto\fX_\C$ and the string diagram construction $\fX\mapsto\C^\fX$, prove the $(2+1)$D Unitary Cobordism Hypothesis (\cref{THREEDUCH}), establish the functor category equivalence (\cref{HOMEQUIV} for $n=2$), and classify $(2+1)$D unitary TQFTs by $\Hstar$-Morita classes of $\Hstar$-multifusion categories (\cref{THREEDCLASS}).

We have tried to make \cref{sec:1+1D} and \cref{sec:2+1D} readable independently of \cref{sec:n+1D}. We have done so completely for $n=1$ and mostly for $n=2$, where some details are deferred to the appendix or \cref{sec:n+1D}. The reader who wants the concrete before the general may therefore start with \cref{sec:1+1D}, or read it in parallel with \cref{sec:n+1D}.

We include appendices with the details on gluing and the construction of the path integral (\cref{sec:gluing-appendix}), unrestricted disk-like functors (\cref{sec:unnrestricted}), disk-like $n$-subcategories and liftable $n$-fields (\cref{sec:liftable}), and the proof of \cref{strictification} (\cref{strictification-sec}).

\subsection*{Acknowledgments}
The author is incredibly grateful to Kevin Walker for his valuable feedback and for many helpful and insightful conversations. The author also thanks David Penneys, Giovanni Ferrer, and Brett Hungar for helpful discussions. The author was supported by NSF grants 2244045 and 2154389.

\section{Overview: \texorpdfstring{$(n+1)$}{(n+1)}D}

\label{sec:n+1D}

Everything below is specialized and unpacked for $n=1$ in \cref{sec:1+1D} and most is spelled out for $n=2$ in \cref{sec:2+1D}; we will sprinkle in pointers to these specializations throughout. The reader who prefers the concrete low-dimensional cases over the general could instead start in those sections and return here when appropriate.

\subsection{Conventions}
\label{sec:conventions}
We collect here the conventions on manifolds, gluing, and notation that we use throughout this article.

\subsubsection{Manifolds}
\label{PL-justification}
We work throughout with compact piecewise linear (PL) manifolds. As we only work with compact PL manifolds in this paper, we will henceforth suppress these adjectives and moreover will suppress the fact that whenever constructions obtain possibly noncompact manifolds, we take their closures. There are several reasons for this choice, but we point out three as follows. First, $\PL(n)$ is thought to be the largest group acting on rigid $(\infty, n)$-categories \cite[Statement 6.2]{FHJF24}. Second, smoothing theory \cite[Essay IV]{KS77} guarantees that framed PL manifolds admit unique smoothings compatible with the PL structure. Finally, the graphical calculus for higher categories naturally lives in the PL setting \cite[§6]{FHJF24}. 

\subsubsection{Oriented manifolds}
Throughout this article, we fix a spatial dimension $n\in\bZ_{\geq 0}$ so that the corresponding spacetime dimension is $n+1$.\footnote{We call attention to this convention because it differs from a portion of the literature wherein $n$ denotes the spacetime dimension, so that our $(n+1)$D TQFT would be called an $n$D TQFT in those references. If we took that approach, we would have to say ``disk-like $(n-1)$-category'' instead of ``disk-like $n$-category.''}
By ``$k$-manifold'' we will mean a $k$-manifold $W$ equipped with the germ of a thickening to an oriented $(n+1)$-manifold $\tld W$ together with a framing $\tau_W$ of its normal bundle $\nu_{W\hookrightarrow\tld{W}}$ in that thickening. By the standard short exact sequence
\[
0 \longrightarrow \rT_{W} \longrightarrow \rT_{\tld W}\vert{}_{W} \longrightarrow \nu_{W\hookrightarrow\tld W} \longrightarrow 0,
\]
this data uniquely determines an orientation on $W$.

Henceforth in this article, unless indicated otherwise, by a ``homeomorphism'' of $k$-manifolds $\varphi\colon W\to W'$ we will mean an \defn{orientation-preserving homeomorphism}, i.e., a homeomorphism $\varphi\colon W\to W'$ of the underlying manifolds that extends to an orientation-preserving homeomorphism of their germs $\tld\varphi\colon\tld W\to\tld {W'}$ such that the induced isomorphism of normal bundles $\tld\varphi_* \colon \nu_{W\hookrightarrow\tld W} \to \nu_{W'\hookrightarrow\tld {W'}}$ preserves the framings in that $\tld\varphi_*\tau_W=\tau_{W'}$.

\subsubsection{Gluing}
Let $W$ be a $k$-manifold and let $Q_1\hookrightarrow \partial W$ and $Q_2\hookrightarrow \partial W$ be embedded codimension-1 submanifolds. To \defn{glue} $W$ along $Q_1$ and $Q_2$ means to identify their restricted germs in $\tld W$, i.e., to provide a homeomorphism $\phi \colon Q_1 \to Q_2$ that extends to an orientation-preserving homeomorphism of the restricted germs $\tld\phi \colon \tld W\vert{}_{Q_1}\to \tld W\vert{}_{Q_2}$ such that the induced isomorphism of normal bundles $\widetilde{\phi}_* \colon \nu_{W \hookrightarrow \widetilde{W}}\vert{}_{Q_1} \to \nu_{W \hookrightarrow \widetilde{W}}\vert{}_{Q_2}$ preserves the restricted framings in that $\widetilde{\phi}_*(\tau_{W}\vert{}_{Q_1})=\tau_{W}\vert{}_{Q_2}$.

The canonical isomorphism $\rT_{W}\vert{}_{Q_i} \cong \nu_{Q_i\hookrightarrow W} \oplus \rT_{Q_i}$ dictates the induced boundary orientation on $Q_i$. Because gluing identifies the outward-pointing normal of $Q_1$ with the inward-pointing normal of $Q_2$, the extension $\widetilde{\phi}$ negates the normal summand $\nu_{Q_i\hookrightarrow W}$. To simultaneously preserve the total orientation of $\rT_{W}\vert_{Q_i}$, the underlying map $\phi \colon Q_1 \to Q_2$ must then negate the tangent summand, so $\phi$ is orientation-\emph{reversing}.

Because this orientation-reversal is uniquely determined by the preserved auxiliary data, we henceforth suppress the ambient germs and framings from our notation. We denote a gluing simply by identifying $Q_1$ with its orientation-reversed copy $Q_2$, writing these as $Q$ and $\orev{Q}$ respectively.

\subsubsection{\texorpdfstring{$k$}{k}-balls and the Alexander Trick}
By ``$k$-ball'' we mean a $k$-manifold homeomorphic to the standard $k$-ball. For $0\leq k\leq n$, we let $\Disk_k$ and $\Mfld_k$ denote the groupoids of homeomorphisms of $k$-balls and $k$-manifolds respectively.

Finally, we will use the broad term ``the Alexander Trick'' to mean that (i) any (orientation-preserving!) homeomorphism of an $n$-ball is isotopic to the identity, (ii) any homeomorphism of an $n$-ball restricting to the identity on the boundary is isotopic-rel-boundary to the identity, and (iii) any homeomorphism of $S^n$ extends to a homeomorphism of $D^{n+1}$ (and is thus isotopic to the identity).

\begin{convention}
\label{sec:fd-convention}
All linear categories and linear 2-categories in this article are assumed to have finite-dimensional hom spaces (for 2-categories, finite-dimensional 2-morphism spaces).
\end{convention}

\subsection{Disk-like \texorpdfstring{$n$}{n}-categories}
\label{sec:dlncats}
Before we introduce the unitary version of a disk-like $n$-category, we quickly review the usual notion of disk-like $n$-category; see \cite[\S6.1]{MW12} for the details. As in \cite{MW12}, we find it best to think in terms of examples; several are presented in \cite[\S6.2]{MW12}, and we focus throughout on the two main classes described in the introduction, namely string diagrams labeled by some ``pivotal $n$-category'' \cite{W06,W14premodular,DR18} and maps to some fixed target space $T$. Examples of the latter type include Dijkgraaf--Witten theories \cite{DW90} (take $T\coloneq \rB G$ for some finite group $G$), the finite homotopy TQFTs of Quinn \cite[Lec. 4]{Q95}, and more generally $\sigma$-model TQFTs. The definitions below are specialized and unpacked for $n=1$ in \cref{sec:def-dl1cat} and for $n=2$ in \cref{sec:def-dl2cat}.

Recall that a disk-like $n$-category captures the local data of an $(n+\e)$D TQFT, that is, of an $(n+1)$D TQFT stripped of its top-dimensional path integral data but retaining the action of homeomorphisms of $n$-balls, and hence the path integral on mapping cylinders. An equivalent way to specify the local data of an $(n+\e)$D TQFT---an $(n+1)$D TQFT without its top-dimensional path integral data---is as a \emph{presentation} (``generators and relations''). One specifies a collection of ``pure fields'' $\eF(X)$ on each $n$-manifold $X$ and, for each $n$-ball $D\hookrightarrow X$, a subspace $U_D\subset\bC\{\eF(D)\}$ of \defn{local relations}, such that two formal linear combinations of pure $n$-fields on $X$ are identified when they differ by an element of the subspace of $\bC\{\eF(X)\}$ generated by the elements $u\blt_{\partial D}\eta$ ($u$ glued to $\eta$ along $\partial D$) with $u\in U_D$ and $\eta$ a field on $X\setminus D$.
This package $(\eF,U)$ is called a \defn{system of fields and local relations}.
On the other hand, instead of carrying around subspaces of local relations, we can impose the local relations from the start, directly enriching our top-dimensional fields in vector spaces.
This package $\C$ is what a disk-like $n$-category is in the sense of \cite{MW12}. These two ways of packaging an $(n+\e)$D TQFT are manifestly equivalent \cite[\S6.1]{MW12}, so there is no loss of generality by using only disk-like $n$-categories. In this article, we will take the former approach: we will define a disk-like $n$-category $\C$ to be the restriction of a system of fields and local relations to balls, while retaining the subspaces $\lU{\C}$ of local relations.

We warn the reader that while the definition below seems long, its content amounts to a list of properties that hold for both classes of the two running examples above; the axioms will seem less complicated if the two examples are kept in mind throughout. We now spell this out, though we point the reader to \cite[\S6.1--3]{MW12} for the details. The definition proceeds by recursion on $k=0,\ldots,n$: for each $1\leq k\leq n$, round $k$ of the recursion takes $\undC_{k-1}$ as input and produces $\undC_k$ via the colimit construction below; all data and conditions are supplied at every round unless marked by ``at $k=n$.''
\phantomsection\label{sec:full-axiomatic-definition}
Following the labeling conventions of \cite{DR18}, axioms whose names begin with ``D'' record \emph{data} and axioms whose names begin with ``C'' record \emph{conditions}. A \defn{disk-like $n$-category} consists of the following data.
\begin{lst}

    \item[(D$\C$)]
    \label{DCk}
    \emph{Fields.}
    Functors $\C_k\colon\Disk_k\to\Set$, $0\leq k\leq n$, from oriented $k$-balls and homeomorphisms to sets; we call $\fld[k]{\C}{W}$ the set of (pure) \defn{$k$-fields} on $W$.
    Functoriality packages a homeomorphism action \axiomtag{D$\varphi$}{Dvarphi} $\varphi_*\colon\fld[k]{\C}{W}\to\fld[k]{\C}{W'}$ for each $\varphi\colon W\xrightarrow{\sim}W'$.
    In the examples, a $k$-field on $W$ is a string diagram drawn on $W$ (which we also refer to as a \defn{$W$-shaped string diagram}) and $\varphi_*$ is the pushforward diagram, while for $\sigma$-models a $k$-field is a map $W\to T$ and $\varphi_*f\coloneq f\circ\varphi^{-1}$.

    \item[(D$\partial$)]
    \label{Dbdk}
    \emph{Boundary maps.}
    Natural transformations $\partial_k\colon \C_k\Rightarrow\undC_{k-1}\circ\partial$. Naturality unpacks as maps $\partial\colon\fld[k]{\C}{W}\to\undfld[k-1]{\undC}{\partial W}$ that are natural with respect to homeomorphisms \axiomtag{C$\partial\varphi$}{Cbdvarphi}; write $\fld[k]{\C}{W}[r]\coloneq\partial^{-1}(r)$. This is vacuous at $k=0$, where we set $\undfld[-1]{\undC}{\varnothing}\coloneq\{*\}$.
    This is restriction of a diagram to $\partial W$ for string diagrams, and $f\mapsto f|_{\partial W}$ for $\sigma$-models.
    (We usually drop the subscript on $\partial_k$.)
    \item[(DG)]
    \label{DGk}
    \emph{Gluing.}
    Maps $-\blt_Q-\colon\fld[k]{\C}{W_1}\times_{\undfld[k-1]{\undC}{Q}}\fld[k]{\C}{W_2}\to\fld[k]{\C}{W_1\cup_Q W_2}$ for each splitting (defined below) along a $(k-1)$-manifold $Q$.
    For string diagrams (resp. $\sigma$-models), this is concatenation of two string diagrams (resp. maps) agreeing along the gluing locus into a new string diagram (resp. map) on $W_1\cup_Q W_2$.

    \item[(D$\pi$)]\label{Dpik}
    \emph{Pullbacks.}
    Maps $\pi^*\colon\fld[k]{\C}{W}_{\pitchfork\theta_\pi}\to\fld[k+m]{\C}{E}$, defined on the fields splittable along the stratification $\theta_\pi$ of $\partial W$ by fiber dimension, for \defn{pinched product maps} $\pi\colon E\to W$, which are PL maps locally modeled on simplex degeneracy maps $\Delta^{k+m}\to\Delta^k$ with $m\geq 1$; we write $\xi\times E\coloneq\pi^*\xi$ when $\pi$ is understood. This is thickening of the diagram over the fiber for string diagrams and precomposition $\pi^*f\coloneq f\circ\pi$ for $\sigma$-models. 
    \item[(D$U$)]\label{DCU}
    \emph{Local relations (at $k=n$).} For each $n$-ball $X$ and $c\in\undfld[n-1]{\undC}{\partial X}$, a subspace $\fldlU{\C}{X}[c]\subset\bC\{\fld[n]{\C}{X}[c]\}$.
\end{lst}

A \defn{splitting} of a $k$-manifold $W$ is a collection of finitely many sub-$k$-balls $\alpha=\{W_i\}$ in $W$ with disjoint interiors and $W=\bigcup_i W_i$. We say a field is \defn{splittable along} (or \defn{transverse to}) a splitting $\alpha=\{W_i\}$ of a $k$-ball $W$ if it lies in the image $\fld[k]{\C}{W}_{\pitchfork\alpha}$ of the corresponding iterated gluing map; since gluing is injective \ref{CGi} and associative \ref{CGa}, such a field has well-defined \defn{restriction maps} to the subballs $W_i$.

\itemstep{Compatibility.} The above data are subject to the condition that they be compatible in essentially every possible way that can be written down. Explicitly, we impose the following conditions.
\begin{lst}
    \item[(C$\varphi$)]\label{Cvarphi} $\C_k$ is a functor $\Disk_k \to \Set$, i.e., $\varphi_*$ acts associatively and unitally.
    \item[(CGi)]\label{CGi} \emph{Gluing is injective.}
    \item[(CGa)]\label{CGa} \emph{Gluing is strictly associative.}
    \item[(CG$\varphi$)]\label{CGvarphi} \emph{Gluing is natural with respect to homeomorphisms.}
    \item[(CG$\partial$)]\label{CGbk} \emph{Gluing commutes with boundary and restriction maps.}
    \item[(CG$\pi$)]\label{Cpigk} \emph{Gluing is compatible with pullbacks.}
    \item[(C$\pi$a)]\label{Cpiak} \emph{Pullbacks are associative.}
    \item[(C$\pi\varphi$)]\label{Cpivarphi} \emph{Pullbacks are natural with respect to homeomorphisms.}
    \item[(C$\pi\partial$)]\label{Cpirk} \emph{Pullbacks are compatible with boundary and restriction maps.}
    \item[(C$U\varphi$)]\label{CUvarphi} \emph{Local relations are natural}: $\varphi_*(\fldlU{\C}{X}[c])=\fldlU{\C}{X'}[\varphi_*c]$ for every homeomorphism $\varphi\colon X\to X'$.
    \item[(C$U$G)]\label{CUi} \emph{Local relations form an ideal under gluing}: if $u\in\fldlU{\C}{X_1}[c_1]$ and $\eta\in\fld[n]{\C}{X_2}[c_2]$ are compatible along $Y$, then $u\blt_Y\eta\in\fldlU{\C}{X}[c]$.
    \item[(CS)]\label{CS} \emph{Splitting}: for every pure $k$-field $\xi$ on a $k$-ball $W$ there is an embedded cell decomposition (splitting) $S_\xi$ of $W$ such that $\xi$ is splittable along every splitting of $W$ transverse to $S_\xi$. Such an $S_\xi$ is called a \defn{string locus} for $\xi$.
\end{lst}
Every one of these conditions is easily seen to be satisfied by our running examples. For instance, the examples of local relations to keep in mind are $\lU{\C}$ being the kernel of an ``evaluation map'' that sends string diagrams to the morphisms they represent, and $\lU{\C}$ being generated by the homotopy relation for maps to a fixed space $T$; both \ref{CUvarphi} and \ref{CUi} hold in these cases, since the evaluation map and the homotopy relation are each equivariant under homeomorphisms, and gluing an element of $\ker(\eval)$---resp.\ a nullhomotopic difference---against an arbitrary field on the rest of $X$ yields another such element. As for \ref{CS}, for a string diagram $\xi$ we can take $S_\xi$ to be the underlying string diagram stratification itself, while the empty set is always a string locus for maps to a space $T$.

\itemstep{Colimit construction \cite[\S6.3]{MW12}.} Each $\C_k$ uniquely extends from balls to any $k$-manifold $W$ by considering \defn{permissible} (non-pathological\footnote{By ``non-pathological'' we mean \defn{permissible ball decomposition} in the sense of \cite[\S6.3]{MW12}, which is defined as a partitioning of $W$ into balls with disjoint interiors such that there exists \emph{some} sequence of pairwise gluings ending with $W$ in which every intermediate space is a valid manifold. This is necessary because one could otherwise partition a manifold into balls that meet along non-manifold boundaries; see \cite[Exm. 3.1.2]{MW12}.}) splittings of $W$ whose constituent balls are equipped with \defn{compatible} fields of $\C_k$---meaning they agree on shared boundaries---and then taking the colimit over all such decompositions $\cP$ equipped with such labelings $\beta$:
\[
    \undfld[k]{\undC}{W} \coloneq \colim_{(\cP,\beta)}\prod_{W_i\in\cP}\fld[k]{\C}{W_i}.
\]
Here the colimit is taken in the poset $\Disk_\C(W)$ of such things ordered by \emph{antirefinement} via iterated gluing set maps. 

If $W$ is itself a $k$-ball, then $\fld[k]{\undC}{W}=\fld[k]{\C}{W}$ since the trivial decomposition is terminal in $\Disk_\C(W)$. As this extension is moreover unique, we will henceforth drop the $\undC$ notation in favor of $\C$.
As for balls, we write $\fld[k]{\C}{W}[c]$ to mean pure $k$-fields on $W$ restricting on $\partial W$ to $c\in \undfld[k-1]{\C}{\partial W}$. 

For string diagrams, $\undfld[k]{\C}{W}$ consists of just the diagrams drawn on $W$, thought of as assembled from compatible diagrams on ball pieces; for $\sigma$-models, $\undfld[k]{\C}{W}$ consists of maps $W\to T$ reconstructed from maps on the pieces that agree on overlaps. 

\begin{proposition}[{\cite[Prop. 6.3.4]{MW12}}]
\phantomsection\label{lem:inj}
For a given $(\cP=\{W_i\};\beta)\in\Disk_\C(W;c)$, the map $\prod_{W_i\in\cP}\fld[k]{\C}{W_i}[\beta|_{\partial W_i}]\hookrightarrow\undfld[k]{\C}{W}[c]$ is injective.
\end{proposition}

\itemstep{Linearization at $k=n$.}
\phantomsection\label{def:linearization}
We extend $\lU{\C}$ to a general $n$-manifold $X$ by declaring $\fldlU{\C}{X}[c]\subset\bC\{\undfld[n]{\C}{X}[c]\}$ to be the subspace generated by every element obtained from a permissible ball decomposition $\{X_i\}$ of $X$ by inserting a local relation $u\in\fldlU{\C}{X_j}[-]$ in one piece $X_j$ and gluing in arbitrary pure fields on the rest. For an $n$-ball $X$, quotienting by the local relations $\fldlU{\C}{X}[c]$ turns the pure $n$-fields $\fld[n]{\C}{X}[c]$ into a vector space $\fldl[n]{\C}{X}[c]$ of \defn{$n$-fields}. A field is \defn{pure} when it is an element of one of the sets $\fld[k]{\C}{W}$ rather than of the linearization $\fldl[n]{\C}{X}[c]$; since this only distinguishes anything in dimension $n$, we simply write ``$k$-field'' for ``pure $k$-field'' when $k<n$, and we use the same convention for transfors (see \cref{def: DL k-transfor} below).
We set $\undfldl[n]{\C}{X}[c]\coloneq\bC\{\undfld[n]{\C}{X}[c]\}/\fldlU{\C}{X}[c]$.
This again agrees on balls with the definition above, written $\fldl[n]{\C}{X}[c]$.

It is a standard exercise in category theory to show that for an $n$-manifold $X$ and a boundary condition $c\in\fld[n-1]{\C}{\partial X}$, the vector space $\undfldl[n]{\C}{X}[c]$ coincides with the vector space obtained by first linearizing the fields on balls of splittings of $X$, and then running a similar colimit construction with a similar poset $\Disk_\C(X,c)$ in $\Vec$. In other words, the result is the same regardless of whether we linearize before or after taking the colimit.

\begin{remark}[Failure of injectivity at $k=n$]
\phantomsection\label{rem:inj-at-n}
Gluing is injective at the set level \ref{CGi}, but the induced linear gluing map on $\fldl[n]{\C}{-}[-]$ need not be injective after quotienting by local relations. (This is why we retain pure $n$-fields rather than enriching in $\Vec$ from the start.\footnote{If we did not retain pure $n$-fields and local relations, then injectivity of gluing required for \cref{lem:eFun disk-like} would fail. To enrich the top dimension in a general distributive symmetric monoidal category $\cS$ in the sense of \cite[Defn. 5.3]{ST12} rather than in $\Vec$, one must instead demand that $\cS$ be equipped with a forgetful functor $\cS\to\Set$ sending the gluing morphisms in $\cS$ to injective set maps.})
\end{remark}

Lastly, we require each $n$-field to be fixed under \defn{extended isotopy}, that is, under any composite of the following two moves. Here $X$ is an $n$-ball, $c\in\undfld[n-1]{\C}{\partial X}$, and $\xi\in\fld[n]{\C}{X}[c]$.
\begin{lst}
    \item[(C$U$i)]\label{CTi} (\emph{Boundary-fixing isotopies}) If a homeomorphism $\varphi\colon X\to X$ fixes $c$ and is isotopic to $\id_X$ through homeomorphisms fixing $c$, then $\varphi_*\xi-\xi\in\fldlU{\C}{X}[c]$.
    \item[(C$U$c)]\label{CTc} (\emph{Collaring homeomorphisms}) $(\psi_{Y,J})_*(\xi\blt_Y(\res_Y\xi\times J))-\xi\in\fldlU{\C}{X}[c]$ for every collaring homeomorphism $\psi_{Y,J}$.
\end{lst}
Here, for an $(n-1)$-ball $Y$ embedded in $\partial X$ and a 1-ball $J$, a \defn{collaring homeomorphism} $\psi_{Y,J}\colon X\cup_Y(Y\times J)\xrightarrow{\sim}X$ is any homeomorphism obtained by gluing on the pinched product $Y\times J$ (pinched along $\partial Y\times J$) along $Y$ and then collapsing the result back into $X$, and in \ref{CTc} we carry the product field along. As usual, both \ref{CTi} and \ref{CTc} extend to all $n$-manifolds by the colimit construction. For string diagrams, extended isotopies act as ordinary isotopies rel boundary.

In this paper, we will always demand the following condition, which is the disk-like analog of requiring finite-dimensional vector spaces of $n$-morphisms.
\begin{lst}
    \item[(CFD)]\label{FDF}
    For all \emph{$n$-balls} $X$ and $c\in\undfld[n-1]{\C}{\partial X}$, $\dim\undfldl[n]{\C}{X}[c]<\infty$.
\end{lst}

\begin{definition}
We will call a disk-like $n$-category $\C$ \defn{finite} if
\begin{lst}
    \item[(F)]
    \label{psnF}
    (\emph{Finiteness})
    $\dim\undfldl[n]{\C}{X}[c]<\infty$ for all $n$-manifolds $X$ and $c\in\undfld[n-1]{\C}{\partial X}$.
\end{lst}
\end{definition}

\begin{remark}
We expect the above notion of ``finite'' to be equivalent to the inductive notion proposed in \cite[Rmk. 2.28]{3Hilb}; see \cref{cor:psnF-weaker} and the finiteness arguments in this article below.
\end{remark}

\begin{remark}
  \label{weakerpsnF}
Note that \ref{psnF} is a condition regarding general $n$-manifolds, not just $n$-balls. However, when $\C$ is unitary, it will follow from \cref{UnitaryWalkerTheorem} that \ref{psnF} is equivalent to $\dim\undfldl[n]{\C}{S^k\times D^{n-k}}[c]<\infty$ for all $k\in\{0,\dots,n\}$ and all boundary conditions $c$; see \cref{cor:psnF-weaker}.
\end{remark}

\subsection{Reflection structures}
\label{dlreflsec}
For each $0\leq k\leq n$, the groupoids $\Disk_k$ and $\Mfld_k$ admit a strictly involutive functor $\orev{\,\cdot\,}$ given on objects by orientation-reversal $W\mapsto \orev{W}$ and on morphisms $\varphi\colon W\to W'$ by $\varphi\mapsto \orev\varphi\colon\orev W\to\orev{W^\prime}$ where $\orev\varphi$ is given by $\varphi$ on the underlying manifolds but viewed as a map between the orientation-reversed manifolds. (Thus in particular $\orev\varphi$ is orientation-\emph{preserving}.) In the same way, we get maps $\orev\pi\colon \orev E\to \orev W$ for each pinched product map $\pi\colon E\to W$. 

As our manifolds are oriented in this article, we will equip all of our disk-like $n$-categories with a structure intertwining their field data with orientation-reversal. As \cite{MW12} gives the detailed definition of a disk-like $n$-category for unoriented manifolds, its authors did not need such a notion. The article \cite{W21} implicitly uses a similar notion, while \cite{W17aasen} gives the details for $n=2$ when considering spin manifolds. In addition to oriented manifolds, \cref{def:reflection structure} should work for manifolds with an $H$-structure and reversal in the sense of \cite{FH21}. For $n=1$ and $n=2$ in this article, we will give the reflection structure axioms alongside the usual disk-like 1- and 2-category axioms in \cref{sec:def-dl1cat} and \cref{sec:def-dl2cat} respectively.

\begin{definition}[Reflection structure]
\label{def:reflection structure}
Let $\C$ be a disk-like $n$-category. A \defn{reflection structure} (or \defn{bar structure}) $\fcj{\,\cdot\,}$ on $\C$ consists of the following data for each $0\leq k\leq n$.
\begin{lst}
\item[(D$\fcj{\,\cdot\vphantom\eta\,}$)]\label{Drefl} For each $k$-ball $W$, an involutive bijection $\fcj{\,\cdot\,}\colon\fld[k]{\C}{W}\to\fld[k]{\C}{\orev{W}}$.
\end{lst}
The above data are subject to the following conditions for each $0\leq k\leq n$.
\begin{lst}

\item[(C$\fcj{\,\cdot\vphantom\eta\,}\varphi$)]
\label{axiom:refl-homeo}
$\orev{\varphi}_*\fcj{\xi} = \fcj{\varphi_*\xi}$ for all homeomorphisms $\varphi\colon W \to W'$ and $k$-fields $\xi\in\fld[k]{\C}{W}$.

\item[(C$\fcj{\,\cdot\vphantom\eta\,}\partial$)]
\label{axiom:refl-bdy} $\partial\fcj{\xi} = \fcj{\partial \xi}$ for all $k$-fields $\xi$.

\item[(C$\fcj{\,\cdot\vphantom\eta\,}$G)]
\label{axiom:refl-gluing}
$\fcj{\xi\blt_E\eta}=\fcj\xi\blt_{\orev{E}}\fcj\eta$ for all $k$-fields $\xi$ and $\eta$ compatible along $E$.

\item[(C$\fcj{\,\cdot\vphantom\eta\,}\pi$)]
\label{axiom:refl-products} $\fcj{\pi^*\xi}=\orev{\pi}^*\fcj{\xi}$ for all pinched product maps $\pi\colon E\to W$ and fields $\xi\in\fld[k]{\C}{W}$.

\item[(C$\fcj{\,\cdot\vphantom\eta\,}U$)]
\label{axiom:refl-local-relations}
(\emph{{Only for $k=n$}}.) 
$\fcj{\fldlU{\C}{X}[c]}=\fldlU{\C}{\orev X}[\fcj c]$, i.e., the \emph{conjugate-linear} extension $\bC\{\fld[n]{\C}{X}\}\to\bC\{\fld[n]{\C}{\orev{X}}\}$ given by $\sum_i\alpha_i\xi_i\mapsto\sum_i\fcj{\alpha_i}\fcj{\xi_i}$ sends $\fldlU{\C}{X}[c]$ into $\fldlU{\C}{\orev X}[\fcj c]$ for all $n$-balls $X$ and $c\in\undfld[n-1]{\C}{\partial X}$. Thus $\fcj{\,\cdot\,}$ descends to a conjugate-linear involutive vector space isomorphism $\fldl[n]{\C}{X}[ c] \to \fldl[n]{\C}{\orev{X}}[ \fcj{c}]$.
\end{lst}
\end{definition}

\begin{convention}
A disk-like $n$-category equipped with a reflection structure is the disk-like analog of a dagger $n$-category; see \cite{FHJF24}.
\end{convention}

We will see that by the same argument used in \cref{prop:functor-extension} to extend disk-like functors, a reflection structure extends uniquely to all manifolds.

\begin{convention}
For the rest of this article, we assume all disk-like $n$-categories are (on oriented balls and) equipped with reflection structures.
\end{convention}

\begin{definition}
  \label{def:closure}
  Let $\C$ be a disk-like $n$-category. For an $n$-ball $X$ and a $0$-field $a\in\fld[0]{\C}{\pt}$, define the \defn{sphere closure} map by
  \[
    \begin{aligned}
      \cl_a\colon \fld{\C}{X}[a\times\partial X] & \longrightarrow \undfld{\C}{S^n},                          \\
      f                   & \longmapsto\cl_a(f)\coloneq \fcj{a\times X}\blt_{a\times\partial X} f.
    \end{aligned}
  \]
  By the gluing axioms, $\cl_a$ descends to a linear map $\cl_a\colon\undfldl[n]{\C}{X}[a\times\partial X]\to\undfldl[n]{\C}{S^n}$.
\end{definition}

\subsection{Disk-like skein \texorpdfstring{$\ell$}{l}-categories}

The following example/construction of a disk-like $n$-category is critically important, since it will allow us to deduce facts about the invariants produced by any unitary TQFT from the results on unitary disk-like higher categories and higher Hilbert spaces, and vice versa. 

The following definition is the disk-like analog of hom $\ell$-categories between $(n-\ell-1)$-morphisms in an ordinary/traditional $n$-category. Briefly, it takes a fixed $(n-\ell)$-manifold $Q$ and a boundary condition on it in $\C$, say $r\in\undfld[n-\ell-1]{\C}{\partial Q}$, and then views the fields of $\C$ on the $(n-\ell+k)$-balls $Q\times W$ as the $k$-fields of shape $W$ in a new disk-like $\ell$-category $\A_\C(Q,r)$. This disk-like $\ell$-category has many names: $\A_\C(Q,r)$ is often called the \defn{dimensional reduction along $(Q,r)$}, but it is also called the \defn{cylinder $\ell$-category over $(Q,r)$} or the \defn{compactification by $(Q,r)$}, whose objects (0-fields) are fields on $Q$, whose 1-fields on a 1-ball $J$ are fields in $\C$ on $Q\times J$, whose 2-fields on a 2-ball $X$ are fields in $\C$ on $Q\times X$. We will call $\A_\C(Q,r)$ the  \defn{skein disk-like $k$-category} on $(Q,r)$, since skeletonization produces the more familiar skein modules $(\ell=0)$, skein (1-)categories $(\ell=1)$, and so on.  

\begin{definition}[{\cite[Exm. 6.2.4]{MW12}}]
  \label{disk-like skein k-categories}
  Let $\C$ be a disk-like $n$-category, fix $0\leq\ell\leq n$, and fix an $(n-\ell)$-manifold $Q\in\Mfld_{n-\ell}$ together with a boundary condition $r\in\undfld[n-\ell-1]{\C}{\partial Q}$.
  Then there is a disk-like $\ell$-category $\A_\C(Q,r)$, called the \defn{disk-like skein $\ell$-category} of $\C$ on $(Q,r)$, defined as follows.
  \begin{lst}
    \item[\ref{DCk}] $\fld[k]{\A_\C(Q,r)}{W}\coloneq \undfld[n-\ell+k]{\C}{Q\times W}[r]$ for $0\leq k\leq \ell$ and $k$-balls $W$, where $Q\times W$ is pinched along $\partial Q\times W$ so that $\partial(Q\times W)=Q\times\partial W$.
    \item[\ref{DCU}] $\fldlU{\A_\C(Q,r)}{X}[c]$ for $\ell$-balls $X$ and $c\in\undfld[\ell-1]{\A_\C(Q,r)}{\bdy X}$ is the subspace of $\bC\{\fld[\ell]{\A_\C(Q,r)}{X}[ c]\}$ consisting of formal linear combinations $\sum_i\lambda_i\xi_i$ such that $\sum_i\lambda_i\xi_i\in\fldlU{\C}{Q\times X}[c]$.
    \item[\ref{Dvarphi}] For a homeomorphism $\varphi \colon W \to W'$, set $\varphi_*\xi \coloneq (\id_Q \times \varphi)_*\xi$.
    \item[\ref{Dbdk}] For $\xi\in\fld[k]{\A_\C(Q,r)}{W}$, the boundary $\partial\xi\in\undfld[k-1]{\A_\C(Q,r)}{\partial W}$ is obtained by restricting $\xi$ to $Q\times\partial W=\partial(Q\times W)$.
    \item[\ref{DGk}] For a splitting $W = W_1 \cup_E W_2$, set $\xi_1 \blt_E^{\A_\C(Q,r)} \xi_2 \coloneq \xi_1 \blt_{Q \times E}^\C \xi_2$.
    \item[\ref{Dpik}] For a pinched product map $\pi \colon E \to W$, set $\pi^*\xi \coloneq (\id_Q \times \pi)^*\xi$.
    \item[\ref{Drefl}] For $\xi\in\fld[k]{\A_\C(Q,r)}{W}$, define $\fcj\xi \in \fld[k]{\A_\C(Q,r)}{\orev{W}}$ by $\fcj{\xi}^{\A_\C(Q,r)} \coloneq \fcj{\xi}^{\C} \in \undfld[n-\ell+k]{\C}{\orev{Q \times W}}[r] = \undfld[n-\ell+k]{\C}{Q \times \orev{W}}[r]$, where we identify $\orev{Q\times W}$ with $Q\times\orev{W}$.
  \end{lst}
The fact that $\A_\C(Q,r)$ is a disk-like $\ell$-category is easily seen to be inherited from $\C$, as is the fact that finiteness of $\C$ implies finiteness of $\A_\C(Q,r)$.
\end{definition}

For the case $n=2$ and $\ell=1$ worked out concretely, see \cref{disk-like skein 1-category}.

\subsection{(Co)isometries and unitary equivalences of fields}
\label{sec:coisometries-general-n}

The definitions in this subsection are the disk-like analogs of isometries, coisometries, and unitaries in a dagger category. The $n=1$ and $n=2$ cases are unpacked in \cref{sec:coisometries-1} and \cref{sec:coisometries-2} respectively, where they are shown to recover the usual familiar notions for dagger 1- and 2-categories.

\begin{definition}\label{def:isometry-1field}
Let $\C$ be a disk-like $n$-category and let $\iota\colon D^1\to\orev{D^1}$ be the orientation-\emph{preserving} homeomorphism $x\mapsto-x$. For $a,b\in\fld[0]{\C}{\pt}$ and a 1-field $\xi\in\fld[1]{\C}{D^1}[\fcj{a}\amalg b]$, set $\xi^\star\coloneq\fcj{\iota_*\xi}$ and define the ``bubble'' $n$-field
\[
  \beta_\xi\coloneq(\xi\times S^{n-1})\cup(b\times D^n)
\]
which is the $n$-field on $D^n$ obtained by thickening $\xi$ over $S^{n-1}$ and capping off with $b\times D^n$.
\begin{itemize}
  \item We call $\xi$ an \defn{isometry} if $\beta_\xi=a\times D^n$ as (quotient-level) $n$-fields.
  \item We call $\xi$ a \defn{coisometry} if $\xi^\star$ is an isometry.
\end{itemize}
We call $\xi$ \defn{unitary} or a \defn{unitary equivalence} if $\xi$ is both an isometry and a coisometry. We write $a\cong^\star b$ if there is a unitary equivalence $\xi$ between $a$ and $b$.

We call a general (pure or quotient-level) 1-field $\xi\in\fld[1]{\C}{X}$ an isometry (resp. coisometry, unitary) if there is a homeomorphism $\varphi\colon X\to D^1$ such that $\varphi_*\xi$ is an isometry (resp. coisometry, unitary) as just defined.
\end{definition}

More generally, we have the following definition, which recovers \cref{def:isometry-1field} by taking $k=1$ as $\A_\C(\pt)=\C$.

\begin{definition}\label{def:isometry-general}
Let $\C$ be a disk-like $n$-category and let $1\leq k\leq n$. Suppose we have a $(k-1)$-manifold $Q$, a $(k-2)$-field $r\in\undfld[k-2]{\C}{\partial Q}$, and a 1-ball $J$. For $(k-1)$-fields $\xi,\eta\in\undfld[k-1]{\C}{Q}[r]$, we call a $k$-field $\alpha\in\undfld[k]{\C}{Q\times J}[\,\fcj{\xi}\cup\eta]$ an \defn{isometry} (resp. \defn{coisometry}, \defn{unitary equivalence}) in $\C$ if $\alpha$ is an isometry (resp. coisometry, unitary equivalence) in the sense of \cref{def:isometry-1field} when viewed as a 1-field on $J$ in the disk-like $(n-k+1)$-category $\A_\C(Q,r)$ of \cref{disk-like skein k-categories}. We write $\xi\cong^\star\eta$ (in $\C$) if there is a unitary equivalence $\alpha$ between $\xi$ and $\eta$ in $\C$.
\end{definition}

\subsection{Disk-like functors and transfors}
In this subsection we introduce disk-like functors---the disk-like analog of traditional/ordinary functors---which are just collections of set maps between fields of two disk-like $n$-categories that commute with boundary maps, gluing maps, pullbacks, and the reflection structure. We also define disk-like $(n,\ell)$-transfors, which are the higher analogs of natural transformations and modifications. 

\subsubsection{Disk-like functors}

\

\begin{definition}[Disk-like functor]
\label{def: DL functor}
Let $\C$ and $\D$ be disk-like $n$-categories.
A \defn{disk-like functor} $\eF\colon\C\to\D$ consists of the following data for each $0\leq k\leq n$.
\begin{lst}
\item[(D$\eF$)]\label{axiom:eF-data}
To each $k$-ball $W$ and each $k$-field $\xi\in\fld[k]{\C}{W}$, a pure $k$-field $\eF_W(\xi)\in\fld[k]{\D}{W}$.
\end{lst}
The data \ref{axiom:eF-data} are subject to the following conditions.
\begin{lst}
\item[($\eF\varphi$)]\label{axiom:eF-nat}
$\eF_{W'}(\varphi_*\xi)=\varphi_*\eF_W(\xi)$ for all homeomorphisms $\varphi\colon W\xrightarrow{\sim} W'$ and $k$-fields $\xi\in\fld[k]{\C}{W}$.
\item[($\eF\partial$)]\label{axiom:eF-boundary} 
$\partial(\eF_W(\xi)) = \eF_{\partial W}(\partial \xi)$ for all $k$-fields $\xi\in\fld[k]{\C}{W}$. More generally, by applying $\eF$ on the component balls of a permissible ball decomposition, $\eF$ extends uniquely to all manifolds; see \cref{prop:functor-extension} for the details.
\item[($\eF$G)]\label{axiom:eF-gluing} 
$\eF_{W_1\cup_E W_2}(\xi\blt_E \eta) = \eF_{W_1}(\xi)\blt_E \eF_{W_2}(\eta)$ for all $k$-fields $\xi\in\fld[k]{\C}{W_1}$ and $\eta\in\fld[k]{\C}{W_2}$ compatible along a $(k-1)$-ball $E$.
\item[($\eF\pi$)]\label{axiom:eF-products} 
$\eF_E(\pi^* \xi) = \pi^* \eF_W(\xi)$ for all pinched product maps $\pi \colon E \to W$ and fields $\xi \in \fld[k]{\C}{W}$.
\item[($\eF\kern.2ex\fcj{\,\cdot\,}$)]\label{axiom:eF-refl}
$\eF_{\orev{W}}(\fcj{\xi}) = \fcj{\eF_W(\xi)}$ for all $k$-fields $\xi\in\fld[k]{\C}{W}$.
\item[($\eF$U)]\label{axiom:eF-local-relations} 
When $k=n$, the linear extension $\bC\{\fld[n]{\C}{X}[ c]\} \to \bC\{\fld[n]{\D}{X}[ \eF_{\partial X}(c)]\}$ of $\eF$ sends $\fldlU{\C}{X}[c]$ into $\fldlU{\D}{X}[ \eF_{\partial X}(c)]$ for all $n$-balls $X$ and $c\in\undfld[n-1]{\C}{\partial X}$.
\end{lst}
We will usually omit the subscript on $\eF$.

Observe that \ref{axiom:eF-local-relations} gives rise to well-defined linear maps $\fldl[n]{\C}{X}[c]\to\fldl[n]{\D}{X}[\eF(c)]$, which we denote by $\widehat\eF_{X,c}$ (or sometimes also by $\eF$).
\end{definition}

Notice that \ref{axiom:eF-data} and \ref{axiom:eF-nat} are together equivalent to asserting that a disk-like functor is a tuple $\eF=(\eF^0,\eF^1,\dots,\eF^n)$ of ordinary natural transformations $\eF^k\colon\C_k\Rightarrow\D_k$ satisfying \ref{axiom:eF-boundary}--\ref{axiom:eF-local-relations}, as naturality is precisely \ref{axiom:eF-nat}.

\begin{remark}[Extension to all manifolds]
\label{prop:functor-extension}
For a $k$-manifold $W$ with $c\in\undfld[k-1]{\C}{\partial W}$ and a permissible ball decomposition $(\cP=\{W_i\};\beta)\in\Disk_\C(W;c)$, the assignment $\{f_i\}\mapsto\{\eF(f_i)\}$ is antirefinement-compatible by \ref{axiom:eF-gluing}, and thus descends to a well-defined set map $\underrightarrow{\eF}_k\colon\undfld[k]{\C}{W}[c]\to\undfld[k]{\D}{W}[\eF(c)]$. Compatibility with homeomorphisms, $\partial$, $\glu$, $\pi^*$, and $\fcj{\,\cdot\,}$, and (when $k=n$) local relations on manifolds follows componentwise on decompositions into balls by \ref{axiom:eF-nat}, \ref{axiom:eF-boundary}, \ref{axiom:eF-gluing}, \ref{axiom:eF-products}, \ref{axiom:eF-refl}, and \ref{axiom:eF-local-relations} respectively.

This extension is unique: if $\underrightarrow{\eG}\colon\undC_k\to\undD_k$ is any other extension of $\eF$ to $k$-manifolds that preserves gluing, then for any $c=\glu\{ f_i\}_i$ in $\undfld[k]{\C}{W}$ we have $\underrightarrow{\eG}(c)=\glu\{\eF(f_i)\}_i=\underrightarrow{\eF}_k(c)$. Thus we will usually write $\eF$ in place of $\underrightarrow{\eF}$. 
\end{remark}

\begin{example}
\label{example:functors-preserve-closure}
  A disk-like functor $\eF \colon \C \to \D$ of disk-like $n$-categories satisfies $\eF(\cl_a(f)) = \cl_{\eF(a)}(\eF(f))$ for all $f \in \fldl{\C}{X}[ a\times S^{n-1}]$ and $a \in \undfldl[0]{\C}{\pt}$.
  Indeed, $\eF(f)$ is defined on $\undfldl{\C}{S^n}$ as the image of $\{ \eF(f|_{X_i})\}_i$ in the
  colimit $\undfldl{\D}{S^n}$ for some (any) permissible ball decomposition $\{X_i\}_i$ of $S^n$.
\end{example}

\subsubsection{Disk-like transfors}

\

\begin{definition}[Disk-like \texorpdfstring{$(n,\ell)$}{(n,l)}-transfor]
  \label{def: DL k-transfor}
  Let $\C$ and $\D$ be disk-like $n$-categories and let $Q\in\Disk_\ell$.
  A \defn{$Q$-shaped $(n,\ell)$-transformation}, also called an \defn{$(n,\ell)$-transfor} or (if $n$ is understood) \defn{$\ell$-transfor}, consists of the following data for each $0\leq k\leq n-\ell$.
\begin{lst}
\item[(D$\eT$)]\label{axiom:eT-data}
To each $k$-ball $W$ and each pure $k$-field $\xi\in\fld[k]{\C}{W}$, a pure $(k+\ell)$-field $\eT_W(\xi)\in\fld[k+\ell]{\D}{W\times Q}$.
\end{lst}
The data \ref{axiom:eT-data} is subject to the following conditions.
  \begin{lst}
    \item[($\eT\varphi$)]\label{axiom:eT-nat}
    $\eT_{W'}(\varphi_*\xi) = (\varphi\times\id_Q)_*\eT_W(\xi)$ for all homeomorphisms $\varphi\colon W\xrightarrow{\sim}W'$ and $k$-fields $\xi\in\fld[k]{\C}{W}$.
    \item[($\eT\partial$)]\label{axiom:eT-boundary} 
    $\partial(\eT_W(\xi)) = \eT_{\partial W}(\partial \xi)\blt_{\partial W \times \partial Q}(\partial \eT)_W(\xi)$ for all $k$-fields $\xi\in\fld[k]{\C}{W}$, where $\partial \eT$ denotes the $\partial Q$-shaped $(\ell-1)$-transfor obtained by restricting to $W\times\partial Q$.

    \item[($\eT$G)]\label{axiom:eT-gluing} $\eT_{W_1\cup_E W_2}(\xi\blt_E \eta) = \eT_{W_1}(\xi)\blt_{E \times Q} \eT_{W_2}(\eta)$ for all splittings $W=W_1\cup_E W_2$ and compatible $k$-fields $\xi,\eta$.

    \item[($\eT\pi$)]\label{axiom:eT-products} $\eT_E(\pi^* \xi)=(\pi\times \id_Q)^* \eT_W(\xi)$ for all pinched product maps $\pi \colon E \to W$ and fields $\xi \in \fld[k]{\C}{W}$.

    \item[($\eT\kern0.2ex\Box$)]\label{axiom:eT-naturality}
    (\emph{Naturality}.) Let $W$ be an $(n{-}\ell{+}1)$-ball and let $\xi\in\fld[n-\ell+1]{\C}{W}$, so that by \cref{rmk:transfor boundary axiom} $\eT_{\partial W}(\partial \xi)\blt_{\partial W \times \partial Q}(\partial \eT)_W(\xi)$ is a pure $n$-field on the $n$-sphere $(\partial W\times Q)\cup(W\times\partial Q)$. Then for every pair of two-ball splittings $\partial W=W_1\cup W_2$ and $\partial Q=Q_1\cup Q_2$ into hemispheres along which this field splits, its restrictions $\alpha_i$ to the $n$-balls $X_i\coloneq(W_i\times Q)\cup(W\times Q_i)$ for $i\in\{1,2\}$ satisfy $\varphi_*[\fcj{\alpha_1}]=[\alpha_2]$ for some homeomorphism $\varphi\colon\orev{X_1}\to X_2$ with $\varphi|_{X_1\cap X_2}=\id$.
    \item[($\eT\kern.2ex\fcj{\,\cdot\,}$)]\label{axiom:eT-refl} 
    $\eT_{\orev{W}}(\fcj{\xi}) = \fcj{\eT_W(\xi)}$ for all $k$-fields $\xi\in\fld[k]{\C}{W}$.

    \item[($\eT$S)]\label{axiom:eT-tameness}
    $\eT$ is \defn{tame}, i.e., there is an embedded cell decomposition $S_\eT \hookrightarrow Q$ such that for every $0\leq k \leq n-\ell$, every $k$-ball $W$, and every pure $k$-field $\xi\in\fld[k]{\C}{W}$, there is a string locus $S_\xi$ such that $S_{\eT_W(\xi)} \subset (S_\xi \times Q) \cup (W \times S_\eT)$.
\end{lst}
\end{definition}

This definition may look unfamiliar at first, and the reader is encouraged to keep the small cases in mind: for $n=1$ an $(n,1)$-transfor is exactly a natural transformation of dagger functors (\cref{dag-thing-i-2}), and for $n=2$ an $(n,1)$- and an $(n,2)$-transfor are a natural transformation and a modification of $\dag$-functors of pivotal $\dag$-2-categories respectively; see \cref{dag-thing-ii-2,dag-thing-ii-3}. 

Notice that \ref{axiom:eT-data} and \ref{axiom:eT-nat} are together equivalent to asserting that $\eT$ is a tuple $(\eT^0,\dots, \eT^{n-\ell})$ of ordinary natural transformations $\eT^k\colon\C_k\Rightarrow\fld[k+\ell]{\D}{-\times Q}$. 

Observe that a disk-like $(n,0)$-transfor is just a disk-like functor. Indeed, if in \cref{def: DL k-transfor} we set $\ell=0$, then axioms \ref{axiom:eT-boundary}, \ref{axiom:eT-gluing}, \ref{axiom:eT-products}, and \ref{axiom:eT-refl} become \ref{axiom:eF-boundary}, \ref{axiom:eF-gluing}, \ref{axiom:eF-products}, \ref{axiom:eF-refl} respectively, while \ref{axiom:eT-tameness} becomes trivial and \ref{axiom:eT-naturality} is vacuous; on the other hand, \ref{axiom:eF-local-relations} is the requirement that $\eF$ preserves local relations, which is vacuous for $(n,\ell)$-transfors unless $\ell=0$. Moreover, like disk-like functors, disk-like $(n,\ell)$-transfors uniquely extend to all manifolds by the same argument as in \cref{prop:functor-extension}.

Tameness \ref{axiom:eT-tameness} is a technical condition to ensure the splittability axiom \ref{CS} holds in the disk-like $n$-category $\eFun(\C{\to}\D)$ of $\ell$-transfors, which we construct below in \cref{def: functor disk-like n-category}.\footnote{It was pointed out by Kevin Walker that if we instead adopted the dense open set variant of the splitting axiom proposed in \cite{MW12}, then tameness \ref{axiom:eT-tameness} would not be necessary. As all examples of disk-like transfors in this article satisfy this property, we deem it reasonable to include it in the definition.}

\begin{remark}
  \label{rmk:transfor boundary axiom}
  Let $\C$ and $\D$ be disk-like $n$-categories. Recall that $\partial(W \times Q) = (\partial W \times Q) \cup_{\partial W \times \partial Q} 
  (-1)^{\dim(W)}(W \times \partial Q)$ where $(-1)^{\dim(W)}(W\times\partial Q)$ denotes $W\times\partial Q$ if $\dim(W)$ is even and $\orev{W\times \partial Q}$ if $\dim(W)$ is odd. Then \ref{axiom:eT-boundary} requires that for all $W \in \Disk_k$ and $\xi \in \fld[k]{\C}{W}$,
  \[
    \underbrace{\partial(\eT_W(\xi))}_{\in \,\fld[k+\ell-1]{\D}{\partial(W \times Q)}}
    =
    \underbrace{\eT_{\partial W}(\partial \xi)}_{\in\, \fld[k+\ell-1]{\D}{\partial W \times Q}}
    \blt_{\partial W\times\partial Q}
    \underbrace{(-1)^{k}(\partial \eT)_W(\xi)}_{\in\,\fld[k+\ell-1]{\D}{(-1)^{k}(W \times \partial Q)}}\kern-.75ex.
  \]
  The following special cases of \ref{axiom:eT-boundary} are useful to know.
  \begin{itemize}
    \item For $\ell=0$ (disk-like functors), $\partial Q = \varnothing$ and $\partial(W\times Q)= \partial W \times Q$, so $\partial (\eT_W(\xi)) = \eT_{\partial W}(\partial \xi)$.
    \item For $k=0$, $\partial W = \varnothing$ and $\partial (W \times Q)= W \times \partial Q$, so we have $\partial( \eT_W(\xi)) = (\partial \eT)_W(\xi)$.
  \end{itemize}
\end{remark}

\subsubsection{Disk-like \texorpdfstring{$n$}{n}-category of \texorpdfstring{$(n,\ell)$}{(n,l)}-transfors}

\

\begin{construction}[Disk-like \texorpdfstring{$n$}{n}-category of disk-like transfors]
\label{def: functor disk-like n-category}
Let $\C$ and $\D$ be disk-like $n$-categories. Then there is a disk-like $n$-category $\eFun(\C{\to}\D)$ with
  \begin{lst}
    \item[\ref{DCk}] $\fld[\ell]{\eFun(\C{\to}\D)}{Q}\coloneq\{Q\text{-shaped }\ell\text{-transfors }\C\to\D\}$ for $0\leq \ell\leq n$ and $\ell$-balls $Q$.
    \item[\ref{DCU}] $\fldlU{\eFun(\C{\to}\D)}{X}[\eT_0]$ for $n$-balls $X$ and $\eT_0\in\undfld[n-1]{\eFun(\C{\to}\D)}{\partial X}$ is the subspace of $\bC\{\fld[n]{\eFun(\C{\to}\D)}{X}[ \eT_0]\}$ consisting of formal linear combinations $\sum_i\lambda_i\eS_i$ such that $\sum_i\lambda_i(\eS_i)_P(\xi)\in\fldlU{\D}{P\times X}[ (\eT_0)_P(\xi)]$ for all $0$-balls $P$ and $\xi\in\fld[0]{\C}{P}$.
  \end{lst}
  The remaining data of $\eFun(\C{\to}\D)$ is defined as follows.
  \begin{lst}
    \item[\ref{Dvarphi}] For a homeomorphism $\varphi \colon Q \to Q'$, set $(\varphi_*\eT)_W(\xi) \coloneq (\id_W \times \varphi)_*(\eT_W(\xi))$.
    \item[\ref{Dbdk}] For $\eT\in\fld[\ell]{\eFun(\C{\to}\D)}{Q}$, the boundary $\partial\eT\in\fld[\ell-1]{\eFun(\C{\to}\D)}{\partial Q}$ is obtained by restricting each $\eT_W(\xi)$ to $W\times \partial Q$.
    \item[\ref{DGk}] For a splitting $Q = Q_1 \cup_E Q_2$, set $(\eT_1 \blt_E \eT_2)_W(\xi) \coloneq (\eT_1)_W(\xi) \blt_{W \times E} (\eT_2)_W(\xi)$.
    \item[\ref{Dpik}] For a pinched product map $\pi \colon E \to Q$, set $(\pi^*\eT)_W(\xi) \coloneq (\id_W \times \pi)^*(\eT_W(\xi))$.
    \item[\ref{Drefl}] For $\eT\in\fld[\ell]{\eFun(\C{\to}\D)}{Q}$, define $\fcj\eT \in \fld[\ell]{\eFun(\C{\to}\D)}{\orev{Q}}$ at a $k$-field $\xi\in\fld[k]{\C}{W}$ by $\fcj{\eT}_W(\xi) \coloneq \fcj{\eT_W(\xi)}  \in \fld[k+\ell]{\D}{W \times \orev{Q}}$ (identifying $\orev{W\times Q}$ with $W\times\orev{Q}$).
  \end{lst}
\end{construction}

To avoid confusion, we will say \defn{pure $(n,n)$-transfor} to mean a pure $n$-field in {$\eFun(\C{\to}\D)$}, i.e., an $(n,n)$-transfor in the sense of \cref{def: DL k-transfor}.

\begin{proposition}
  \label{lem:eFun disk-like}
  If $\C$ and $\D$ are disk-like $n$-categories, then $\eFun(\C{\to}\D)$ is a disk-like $n$-category, except possibly not satisfying \ref{FDF}.
\end{proposition}

\begin{pf}
The assertion is a straightforward check; most axioms follow from checking their counterparts pointwise in $\D$.
\end{pf}

To obtain \ref{FDF}, see \cref{thm:eFun-sphere-weight}.

\subsubsection{Equivalences of disk-like \texorpdfstring{$n$}{n}-categories}
\label{sec:equivalences-general-n}
Let $\C$ and $\D$ be disk-like $n$-categories. 

\begin{definition} 
\label{def:wkequiv}
A disk-like functor $\eF\colon\C\to\D$ is a \defn{weak equivalence} if
\begin{lst}
  \item[(WE$k$)]\label{WEgenk} $\eF$ is \defn{unitarily essentially surjective on $k$-fields} for each $0\leq k<n$: for every $k$-ball $W$ and every $\eta\in\fld[k]{\D}{W}$ whose boundary lies in the image of $\eF$, there is a $\xi\in\fld[k]{\C}{W}$ with $\eF(\xi)\cong^\star\eta$ in $\D$; and
  \item[(WE$n$)]\label{WEgenn} $\eF$ is a linear isomorphism on $n$-fields: for every $n$-ball $X$ and every $c\in\undfld[n-1]{\C}{\partial X}$, the induced linear map $\widehat{\eF}_{X,c}\colon\fldl[n]{\C}{X}[c]\to\fldl[n]{\D}{X}[\eF(c)]$ is an isomorphism.
\end{lst}
We write $\C\cong\D$ if there is a zig-zag of weak equivalences between $\C$ and $\D$.
\end{definition}

\subsection{Unitary disk-like \texorpdfstring{$n$}{n}-categories}
\label{sec:unitary disk-like n-categories}

We can now define unitary disk-like $n$-categories. The $n=1$ and $n=2$ cases are treated in \cref{sec:unitary disk-like 1-categories} and \cref{sec:unitary disk-like 2-categories} respectively.

\begin{remark}
\label{lem:reparam-gen}
If $\Sigma$ is an $n$-sphere, meaning there is some homeomorphism $\varphi\colon \Sigma\to S^n$, then there is a canonical linear isomorphism
$\undfldl[n]{\C}{\Sigma}\to\undfldl[n]{\C}{S^n}$ given by sending an $n$-field $\xi\in\undfldl[n]{\C}{\Sigma}$ to $\varphi_*\xi\in\undfldl[n]{\C}{S^n}$. Indeed, by the Alexander Trick, any two orientation-preserving homeomorphisms $\varphi,\psi\colon \Sigma\to S^n$ are isotopic, and hence by \ref{CTi} induce equal maps $\varphi_*=\psi_*\colon\undfldl[n]{\C}{\Sigma}\to\undfldl[n]{\C}{S^n}$. This remark is also true if we replace $n$-spheres with $n$-balls, provided the homeomorphisms agree on the boundary (as is automatic for $n=1$), and we will use it very frequently.
\end{remark}

\begin{notation}
\label{rmk:sphere-functional-evaluation}
If $\Sigma$ is an $n$-sphere, meaning there is some homeomorphism $\Sigma\to S^n$, then by \cref{lem:reparam-gen} any linear functional $\psn \colon \undfldl[n]{\C}{S^n} \to \bC$ uniquely induces a functional on $\undfldl[n]{\C}{\Sigma}$, which justifies the following notation: for any linear functional $\psn\colon\undfldl[n]{\C}{S^n}\to\bC$ and any $\xi\in\undfldl[n]{\C}{\Sigma}$, we define
\[
\psn(\xi) \coloneq \psn(\varphi_*\xi)
\]
for some (any) orientation-preserving homeomorphism $\varphi\colon\Sigma\to S^n$. In particular, for any $n$-ball $X$ and any $c\in\undfld[n-1]{\C}{\partial X}$ and any $f,g \in \fldl[n]{\C}{X}[c]$, the scalar $\psn(\fcj{f}\blt_{\partial X} g)$ is well-defined because $\orev{X} \cup_{\partial X} X$ is an $n$-sphere.
\end{notation}

\begin{notation}
\label{reflection and functionals}
Up to isotopy, there is a unique orientation-\emph{reversing} homeomorphism $\iota\colon S^n\to S^n$. Postcomposing this with the orientation-reversal $\orev{\,\cdot\,}$ gives an orientation-\emph{preserving} homeomorphism $\orev{\,\cdot\,}\circ\iota\colon S^n\to \orev{S^n}$, which gives a distinguished linear isomorphism $(\fcj{\,\cdot\,}\circ \iota)_*\colon\undfldl[n]{\C}{S^n}\to\undfldl[n]{\C}{\orev{S^n}}$. As $\undfldl[n]{\C}{S^n}$ is invariant under the actions of isotopies of the sphere by \ref{CTi}, this justifies the following notation: for any linear functional $\psn\colon\undfldl[n]{\C}{S^n}\to\bC$ and any $\xi\in \undfldl[n]{\C}{S^n}$, we define
\[
\Phi(\fcj\xi)\coloneq \Phi(\iota_*^{-1}\fcj\xi).
\]
\end{notation}

\begin{definition}[Unitary disk-like \texorpdfstring{$n$}{n}-category]
  \label{def:unitary disk-like n-category}
  A \defn{unitary disk-like $n$-category} $(\C,\psn)$ is a disk-like $n$-category $\C$ equipped with a \defn{sphere trace} (or \defn{trace}), that is, a linear functional $\psn\colon\undfldl[n]{\C}{S^n}\to\bC$ satisfying the following condition.
  \begin{lst}[font=\upshape]
    \item[($\psn$P)]\label{psnP}
    (\emph{Positivity})
    For each $n$-ball $X$ and each $c\in \undfld[n-1]{\C}{\partial X}$, the sesquilinear pairing
    \begin{equation}\label{disk pairings}
      \begin{aligned}
        \orev{\fldl[n]{\C}{X}[c]}\otimes_{\bC}\fldl[n]{\C}{X}[c] & \longrightarrow \bC,
        \\
        f\otimes g                      & \longmapsto \bkt{f}{g}_{X,c}\coloneq\psn(\fcj f\blt_{\partial X} g)
      \end{aligned}
    \end{equation}
    is positive-definite.
  \end{lst}
\end{definition}

\begin{proposition}
  \label{prop:psnP implies}
  Let $\C$ be a unitary disk-like $n$-category with sphere trace $\psn$.
  \begin{lst}[font=\upshape]
    \item[($\psn$R)]\label{psnR}
    $\psn(\fcj\xi)=\overline{\psn(\xi)}$ for all $\xi\in\undfldl[n]{\C}{S^n}$.
    \item[($\psn$I)]\label{psnI}
    For any homeomorphism $\varphi\colon D^{n+1}\to D^{n+1}$ and any $\xi\in\undfldl[n]{\C}{S^n}$, we have $\psn((\varphi  |_{S^n})_*\xi)=\psn(\xi)$.
  \end{lst}
\end{proposition}

\begin{pf}
  \itemstep{\ref{psnR}.} Since we can always isotope $\xi$ to be transverse to a string locus, we can write $\xi$ as $\xi = \fcj f \blt_c g$ so that $\psn(\xi) = \bkt{f}{g}_{X,c}$. Note that $\fcj \xi = f \blt_{\fcj c} \fcj g \in \undfldl{\C}{\orev{S^n}}$. Applying $\iota^{-1}$ swaps the hemispheres of $S^n$, and thus (up to isotopy) corresponds to swapping the position of the fields on $S^n$, thus mapping $f \blt \fcj g$ to $\fcj g \blt f$. It follows that
  \[
    \psn(\fcj \xi) \underset{\eqref{reflection and functionals}}=\psn(\iota_*^{-1}(\fcj \xi)) = \psn(\fcj g \blt_c f) = \bkt{g}{f}_{X,c}= \bar{\bkt{f}{g}_{X,c}}=\bar{\psn(\xi)}
  \]
  where we used that the pairing is conjugate-symmetric (Hermitian) by \ref{psnP}. 
  
  \itemstep{\ref{psnI}.} Elements of $\undfldl{\C}{S^n}$ enjoy isotopy invariance, and by the Alexander Trick every homeomorphism of $S^n$ is isotopic to the identity.
\end{pf}

\begin{definition} 
For unitary disk-like $n$-categories $\C$ and $\D$, we call a weak equivalence $\eF\colon\C\to\D$ \defn{isometric} if $\psn^\C=\psn^\D\circ\widehat\eF_{S^n}$. We write $\C\cong^\dag\D$ if there is a zig-zag of isometric weak equivalences between $\C$ and $\D$.
\end{definition}

As with \cref{sec:coisometries-general-n}, the notion of weak equivalence in \cref{def:wkequiv} requires only the reflection structures on $\C$ and $\D$; only \emph{isometric} weak equivalence refers to the sphere traces.

\subsection{\texorpdfstring{$(n+1)$D}{(n+1)D} unitary TQFTs}
Recall that for a finite-dimensional complex vector space $V$ equipped with a positive-definite sesquilinear form $\bkt--$, the \defn{canonical element} of $V\otimes V^\vee$ is the element $\Omega_V$ corresponding to the identity $\id_V\in\End(V)$ under the isomorphism $\End(V)\overset\sim\to V\otimes V^\vee$ with inverse $\ket\xi\otimes\bra\eta\mapsto \ket\xi\bra\eta$. Thus $\Omega_V$ is manifestly basis-independent. 

\begin{remark}
\label{rmk:omega-basis}
Observe that if $\{e_i\}$ is a $\bkt--$-orthogonal basis of $V$, then
\[
    \Omega_V=\sum_i\frac{\ket{e_i}\otimes\bra{e_i}}{\bkt{e_i}{e_i}}.
\]
If $W$ is another finite-dimensional vector space and $\psi\in V^\vee\otimes V\otimes W$, we will write $\Omega_V\triangleright \psi$ to mean the element of $W$ obtained by viewing $\psi$ as a linear map $V\otimes V^\vee\to W$ and applying it to the canonical element $\Omega_V$:
\begin{equation}\label{116i}
    \Omega_V\triangleright \psi \coloneq \psi(\Omega_V) = \sum_i \frac{\psi(\ket{e_i}\otimes\bra{e_i})}{\bkt{e_i}{e_i}}\in W.
\end{equation}
\end{remark}

The following facts follow immediately from \cref{rmk:omega-basis}.

\begin{facts}\label{lem:Omega-scale}
For a finite-dimensional complex vector space $V$ with inner product $\bkt--$, the following hold.
\begin{lst}
  \item\label{it:Omega-rescale} $\Omega_{(V,\lambda\bkt--)}=\lambda^{-1}\Omega_V$ for $\lambda\in\bR_{>0}$.
  \item\label{it:Omega-additive} If $V$ is an orthogonal direct sum $V=V_1\oplus V_2$, then $\Omega_V=\Omega_{V_1}+\Omega_{V_2}$.
\end{lst}
\end{facts}

\begin{definition}
\label{defUnitaryTQFT}
An \defn{$(n+1)$D unitary TQFT} is an $(n+\e)$D TQFT $\C$ equipped with a \defn{path integral}
(also called a \defn{functional integral} or \defn{partition function}), 
that is, a collection of linear functionals
\[
\Z_\C(M)\colon \undfldl[n]{\C}{\partial M}\to\bC
\]
indexed by $(n+1)$-manifolds $M$ that satisfies the following conditions. (The sphere trace is recovered from the path integral as $\psn^\C=\Z_\C(D^{n+1})$; see \cref{UnitaryWalkerTheorem}.)
\begin{lst}[leftmargin=1.15cm]
    \item[(Rfl)]\label{Rfl} For each $(n+1)$-manifold $M$ and each $\psi\in\undfldl[n]{\C}{\partial M}$, $\Z_\C(\orev{M})\fcj{\psi}=\overline{\Z_\C(M)\psi}$.
    \item[(Pos)]\label{Pos} For each $n$-manifold $X$ and each $c\in \undfld[n-1]{\C}{\partial X}$, the sesquilinear pairing
    \begin{equation}
    \label{pair}
      \begin{aligned}
        \overline{\undfldl[n]{\C}{X}[c]}\otimes_{\bC}\undfldl[n]{\C}{X}[c] & \longrightarrow \bC,
        \\
        f\otimes g                      & \longmapsto \bkt{f}{g}_{X,c}\coloneq\Z_\C(X\times I)(\fcj f\blt_c g)
      \end{aligned}
    \end{equation}
    is positive-definite (where $X\times I$ is pinched so that $\partial (X\times I)=\orev{X}\cup_{\partial X} X$).
    \item[(Glu)]\label{Glu}
    If an $(n+1)$-manifold $M_\glu$ is obtained from $M$ by gluing along $n$-manifolds $\orev{X}$ and $X$ in $\bdy M$, then we define
      \begin{equation}\label{eq:gluing-formula}
          \Z_\C(M_\glu)(\phi_\glu)
          \;\coloneq\;
          \Omega_{\undfldl[n]{\C}{X}[c]}\triangleright\Z_\C(M)\!\left((-)\blt\phi\right)
      \end{equation}
      for $\phi\in\undfldl{\C}{R}[\fcj c\amalg c]$ where $R\coloneq\partial M\setminus(\orev X\amalg X)$ and $c\in\undfld[n-1]{\C}{\bdy X}$.
    \item[(Inv)]\label{Inv}
    For $(n+1)$-manifolds $M$ and $M'$ and a homeomorphism $\varphi\colon M\to M'$, we have $\Z_\C(M')({\varphi|_{\partial M*}}\psi)=\Z_\C(M)\psi$ for all $\psi\in\undfldl[n]{\C}{\partial M}$.
\end{lst}
\end{definition}

Axiom \ref{Rfl} is equivalent to requiring that orientation reversal of bordisms corresponds to taking adjoint operators with respect to the pairings \eqref{pair}; see \cref{RvsDagR} below.

\begin{remark}
 By using the reflection structure $\fcj{\,\cdot\,}$ and the positive-definite pairing on $\undfldl[n]{\C}{X}[c]$ to identify $\undfldl[n]{\C}{X}[c]^\vee\cong\overline{\undfldl{\C}{X}[c]}\cong \undfldl[n]{\C}{\orev X}[\fcj c]$ via $\bra{e_i}\mapsto\fcj{e_i}$, we have
  \[
      \Z_\C(M_\glu)(\phi_\glu)
      \underset{\eqref{eq:gluing-formula}}{=}
      \Omega_{\text{$\undfldl{\C}{X}[c]$}}\triangleright\Z_\C(M)\!\left((-)\blt\phi\right)
      \underset{\eqref{116i}}{=}
      \sum_i\frac{\Z_\C(M)
      ((\fcj{e_i}\amalg e_i)\blt\phi)}{\bkt{e_i}{e_i}_{X,c}}
  \]
  for any orthogonal basis $\{e_i\}$ of $\undfldl{\C}{X}[c]$. 
\end{remark}

Once a sphere trace $\psn$ is fixed on a finite disk-like $n$-category $\C$, by the Unitary Walker Extension Theorem (\cref{UnitaryWalkerTheorem}) there is a unique path integral $\Z_\C$ extending $\C$ to an $(n+1)$D unitary TQFT such that $\Z_\C(D^{n+1})=\psn^\C$. The proof uses the gluing lemmas and Walker's handle construction of the path integral, which can be found in \cref{sec:gluing-appendix}.

\begin{example}
\label{inducspwt}
Let $\C$ be a finite unitary disk-like $n$-category, let $Q$ be an $(n-\ell)$-manifold for some $0\leq \ell\leq n$, and let $r\in\undfld[n-\ell-1]{\C}{\partial Q}$. Then the disk-like skein $\ell$-category $\A_\C(Q,r)$ from \cref{disk-like skein k-categories} is unitary when equipped with the sphere trace given by the path integral on $Q\times D^{\ell+1}$ (constructed in \cref{cstr:path-integral-via-handle-decomp} below):
\[
\psn^{\A_\C(Q,r)}(\xi)
\coloneq
\Z_\C(Q\times D^{\ell+1})(\xi).
\]
That $\psn^{\A_\C(Q,r)}$ satisfies \ref{psnP} follows from the Unitary Walker Extension Theorem (\cref{UnitaryWalkerTheorem}).
\end{example}

\begin{proposition}
\label{bordism gives a linear map}
\label{RvsDagR}
For all $(n+1)$-manifolds $M$, $n$-manifolds $X_{\inn},X_{\out},W\in\Mfld_{n}$ with $\bdy M=(\orev{X_{\inn}}\amalg X_{\out})\cup W$, and fields $c_\inn\in\undfld[n-1]{\C}{\partial X_\inn}$, $c_\out\in\undfld[n-1]{\C}{\partial X_\out}$, and $v\in \undfldl[n]{\C}{W}[\fcj{c_{\inn}}\amalg c_{\out}]$, let $F_{M}^v\colon \undfldl[n]{\C}{X_{\inn}}[c_{\inn}]\to \undfldl[n]{\C}{X_{\out}}[c_{\out}]$ be the linear map determined by 
\begin{equation}
\label{R}
\bkt{F_{M}^v\eta}{\xi}_{X_{\out},c_{\out}}=\Z_\C(M)(\fcj{\eta}\bullet
v\bullet
\xi)\qquad\forall \eta\in \undfldl[n]{\C}{X_{\inn}}[c_{\inn}],\,\xi\in \undfldl[n]{\C}{X_{\out}}[c_{\out}].
\end{equation}
Then the following are equivalent.
\begin{lst} 
\item\label{RvsDagRi}$\Z_\C(\orev{M})\fcj{\psi}=\bar{\Z_\C(M)\psi}$ for all $\psi\in \undfldl[n]{\C}{\bdy M}$ such that $\psi|_{W}=v$.
\item\label{RvsDagRii} The operator $F_{\orev M}^{\fcj{v}}$ characterized by Condition \ref{R} is the adjoint of $F_M^v$ with respect to the TQFT inner products:
\[
    \bkt{b}{F_M^v a}_{X_{\out},c_{\out}}=\bkt{F_{\orev M}^{\fcj v}b}{a}_{X_{\inn},c_{\inn}}\qquad\forall a\in \undfldl[n]{\C}{X_{\inn}}[c_{\inn}],\;b\in \undfldl[n]{\C}{X_{\out}}[c_{\out}].
\]
\end{lst}
\end{proposition}

\begin{pf}
Indeed, to see the forward implication, write
\[
\bkt{F_{\orev M}^{\fcj v}b}{a}_{X_{\inn},c_{\inn}}
\overset{\eqref{R}}{=}
\Z_\C(\orev M)(\fcj{b}\bullet \fcj v\bullet a)
\overset{\text{\ref{RvsDagRi}}}{=}
\overline{\Z_\C(M)(\fcj{a}\bullet v\bullet b)}
\overset{\eqref{R}}{=}
\overline{\bkt{F_{M}^{v} a}{b}_{X_{\out},c_{\out}}}
=
\bkt{b}{F_{M}^{v} a}_{X_{\out},c_{\out}},
\]
and for the reverse implication, note that if $\psi=\fcj a\bullet v\bullet b$, then
\[ 
\Z_\C(M)\psi
\overset{\eqref{R}}{=}
\bkt{F_M^{v} a}{b}_{X_{\out},c_\out}
=\overline{\bkt{b}{F_M^{v} a}_{X_{\out},c_\out}}
\overset{\text{\labelcref{RvsDagRii}}}{=}
\overline{\bkt{F_{\orev{M}}^{\fcj{v}} b}{a}_{X_\inn,c_\inn}}
\overset{\eqref{R}}{=}
\overline{\Z_\C(\orev{M})(\fcj{\psi})}.\qedhere
\]
\end{pf}

\subsection{Sphere traces on functor categories}
\label{subsec:eFun-sphere-weight}

In this subsection we equip the disk-like functor category $\eHom(\C{\to}\D)$ with a sphere trace making it a unitary disk-like $n$-category. 

The goal of this subsection is to prove \cref{thm:eFun-sphere-weight} below, which exhibits a canonical sphere trace on a disk-like hom category with which it becomes unitary. This sphere trace will be given by a sum over the components (see \cref{rmk:DLcomponents} below) of $\C$ weighted by certain scalars $w_\C(b)$ defined below. These scalars first appeared in \cite{W21}. 

Throughout, $\C$ and $\D$ are finite unitary disk-like $n$-categories. For a 0-field $a\in\fld[0]{\C}{\pt}$ we define the (quotient-level) $n$-fields $\Omega_a$ and $\varnothing_a$ by
\[
    \Omega_a
    \coloneq 
    [a\times D^n]
    \in
    \fldl[n]{\C}{D^n}[r_a]
    \qquad \text{and} \qquad
    \varnothing_a
    \coloneq
    [a\times S^n]\in\undfldl[n]{\C}{S^n}
\]
where $r_a \coloneq a\times S^{n-1}$ and $[x]$ denotes the quotient-level $n$-field with $x$ as a pure $n$-field representative.

The following is a special case of the well-known action of the tube algebra on an $n$D annulus.

\begin{definition}
Consider an $n$-field $e\in\undfldl{\C}{S^{n-1}\times D^1}[\fcj{r_b}\amalg r_a]$. For $\alpha\in\fldl[n]{\C}{D^n}[r_a]$ and $\beta\in\fldl[n]{\C}{D^n}[r_b]$, we define 
\[
	e\triangleright\beta\coloneq e\blt_{S^{n-1}}\beta\in\fldl[n]{\C}{D^n}[r_a]
    \qquad\text{and}\qquad
    \alpha\triangleleft e\coloneq e^\star\triangleright\alpha\in\fldl[n]{\C}{D^n}[r_b]
\]
where $e^\star\coloneq \fcj{(\id_{S^{n-1}}\times \iota^{(1)})_* e} \in\undfldl{\C}{S^{n-1}\times D^1}[\fcj{r_a}\amalg r_b]$ and $\iota^{(1)}$ is the orientation-\emph{preserving} homeomorphism $\iota^{(1)}\colon D^1\to \orev{D^1}$ given by $t\mapsto -t$. For instance, for $n=2$ we can illustrate these actions as
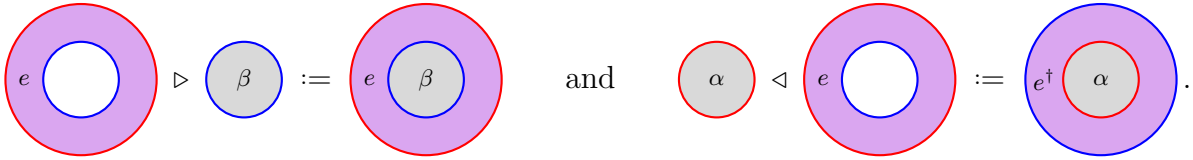
\[
\begin{tkz}[scale=0.5]
\begin{scope}
\invclip(0,0)circle(1);
\fill[violet!35](0,0)circle(2);
\end{scope}
\draw[blue,thick](0,0)circle(1);
\draw[red,thick](0,0)circle(2);
\node at (-1.5,0){\scriptsize$e$};
\end{tkz}
\;\triangleright\;
\begin{tkz}[scale=0.5]
\fill[gray!30] (0,0) circle (1);
\draw[blue,thick](0,0)circle(1);
\node at (0,0){\scriptsize$\beta$};
\end{tkz}
\;\coloneq\;
\begin{tkz}[scale=0.5]
\begin{scope}
\invclip(0,0)circle(1);
\fill[violet!35](0,0)circle(2);
\end{scope}
\fill[gray!30] (0,0) circle (1);
\node at (0,0){\scriptsize$\beta$};
\node at (-1.5,0){\scriptsize$e$};
\draw[blue,thick](0,0)circle(1);
\draw[red,thick](0,0)circle(2);
\end{tkz}
\qquad\text{and}\qquad
\begin{tkz}[scale=0.5]
\fill[gray!30] (0,0) circle (1);
\draw[red,thick](0,0)circle(1);
\node at (0,0){\scriptsize$\alpha$};
\end{tkz}
\;\triangleleft\;
\begin{tkz}[scale=0.5]
\begin{scope}
\invclip(0,0)circle(1);
\fill[violet!35](0,0)circle(2);
\end{scope}
\draw[blue,thick](0,0)circle(1);
\draw[red,thick](0,0)circle(2);
\node at (-1.5,0){\scriptsize$e$};
\end{tkz}
\;\coloneq\;
\begin{tkz}[scale=0.5]
\begin{scope}
\invclip(0,0)circle(1);
\fill[violet!35](0,0)circle(2);
\end{scope}
\fill[gray!30] (0,0) circle (1);
\node at (0,0){\scriptsize$\alpha$};
\node at (-1.5,0){\scriptsize$e^\star$};
\draw[red,thick](0,0)circle(1);
\draw[blue,thick](0,0)circle(2);
\end{tkz}\,.
\]
\end{definition}

We next recall the notion of a minimal 0-field from \cite{W21}, which we adapt to the disk-like framework as follows.
\begin{definition}[{\cite{W21}}]
  A 0-field $b$ in a disk-like $n$-category $\C$ is called \defn{minimal} if $\fldl[n]{\C}{D^n}[r_b]\cong\bC$.
\end{definition}

\begin{remark}
As we work over $\bC$ and only deal with equivalence classes of 0-fields with respect to a certain connectedness equivalence relation (see \cref{cnd0fld} below), requiring $\fldl[n]{\C}{D^n}[r_b]\cong\bC$ is the same as requiring that this algebra be simple. If we were to work over a ground field other than $\mathbb{C}$ or with manifolds such as spin manifolds, then the algebra of a minimal 0-field may be a more general Clifford algebra, in which case minimal 0-fields must instead be defined as those for which this algebra is simple over the ground field.
\end{remark}

\begin{lemma} 
\label{bnz}
If $b$ is a minimal 0-field, then $\Omega_b\neq0$.
\end{lemma} 

\begin{pf} 
We prove the contrapositive. Suppose $b$ is any 0-field with $\Omega_b=0$, so that $\tld\Omega_b\in\fldlU{\C}{D^n}[r_b]$ where $\tld\Omega_b\coloneq b\times D^n$ (a pure $n$-field). For any $\beta\in\fld[n]{\C}{D^n}[r_b]$, by \ref{CTc} we can write $\beta=\beta\blt_Y\tld\Omega_b$ for some $(n-1)$-ball $Y\hookrightarrow \partial D^n$. Since $\tld\Omega_b\in\lU{\C}$, by \ref{CUi} $\beta=\beta\blt_{Y}\tld\Omega_b\in\fldlU{\C}{D^n}[r_b]$. Thus $\beta=0$.
\end{pf}

\begin{corollary}
A 0-field $b$ in a disk-like $n$-category $\C$ is minimal if and only if $\fldl[n]{\C}{D^n}[r_b]=\bC\Omega_b$.
\end{corollary}

\begin{definition}[Sphere trace and Walker norm {\cite{W21}}]
\label{def:sphere-trace-walker-norm}
Let $\C$ be a unitary disk-like $n$-category and let $b$ be a 0-field in $\C$. In \cite{W21}, Walker defines scalars
\[
\tr_s^\C(b)\coloneq\psn^\C(\varnothing_b)
\qquad\text{and}\qquad
\rN^\C(b)\coloneq\bkt{\varnothing_b}{\varnothing_b}_{S^n}
\]
which we will respectively call the \defn{sphere trace} of $b$ and the \defn{Walker norm} of $b$. We will define the \defn{weight} of a 0-field $b$ in $\C$ by
\[
  w_\C(b)\coloneq\frac{\tr_s^\C(b)}{\rN^\C(b)}.
\]
\end{definition}

\begin{remark}
\label{posthings}
For a minimal 0-field $b$, since $\Omega_b\neq 0$, $\psn^\C(\varnothing_b)=\bkt{\Omega_b}{\Omega_b}_{D^n,r_b}>0$ by \ref{psnP}. Hence $\varnothing_b\neq0$, and so $\rN^\C(b)>0$ by positive-definiteness of $\bkt--_{S^n}$ (see \cref{UnitaryWalkerTheorem}). It follows that $w_\C(b)>0$ as well.
\end{remark}

\begin{notation}
\label{notation-lambda-rho}
Observe that when $b$ is a minimal 0-field and $e\in\undfldl{\C}{S^{n-1}\times D^1}[\fcj{r_b}\amalg r_a]$, then $\Omega_a\triangleleft e\in\fldl[n]{\C}{D^n}[r_b]=\bC\Omega_b$. Similarly, when $a$ is minimal, $e\triangleright\Omega_b\in\fldl[n]{\C}{D^n}[r_a]=\bC\Omega_a$. Thus if both $a$ and $b$ are minimal, then
\[
	e\triangleright\Omega_b=\lambda_e\Omega_a 
    \qquad \text{and}\qquad
    \Omega_a\triangleleft e=\rho_e\Omega_b
\]
for scalars $\lambda_e,\rho_e\in\bC$. Diagrammatically,
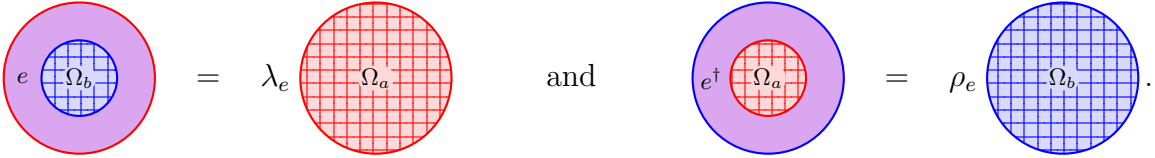
\[
\begin{tkz}[scale=0.5]
\begin{scope}
\invclip(0,0)circle(1);
\fill[violet!35](0,0)circle(2);
\end{scope}
\fill[blue!15,region={grid,color=blue,angle=0}] (0,0) circle (1);
\node[fill=blue!15,inner sep=1pt,rounded corners] at (0,0){\scriptsize$\Omega_b$};
\node at (-1.5,0){\scriptsize$e$};
\draw[blue,thick](0,0)circle(1);
\draw[red,thick](0,0)circle(2);
\end{tkz}
\quad
=
\quad
\lambda_e
\;
\begin{tkz}
\fill[red!15,region={grid,color=red,angle=0}] (0,0) circle (1);
\draw[red,thick](0,0)circle(1);
\node[fill=red!15,inner sep=1pt,rounded corners] at (0,0){\scriptsize$\Omega_a$};
\end{tkz}
\qquad\quad\text{and}\quad\qquad
\begin{tkz}[scale=0.5]
\begin{scope}
\invclip(0,0)circle(1);
\fill[violet!35](0,0)circle(2);
\end{scope}
\fill[red!15,region={grid,color=red,angle=0}] (0,0) circle (1);
\node[fill=red!15,inner sep=1pt,rounded corners] at (0,0){\scriptsize$\Omega_a$};
\node at (-1.5,0){\scriptsize$e^\star$};
\draw[red,thick](0,0)circle(1);
\draw[blue,thick](0,0)circle(2);
\end{tkz}
\quad
=
\quad
\rho_e
\;
\begin{tkz}
\fill[blue!15,region={grid,color=blue,angle=0}] (0,0) circle (1);
\draw[blue,thick](0,0)circle(1);
\node[fill=blue!15,inner sep=1pt,rounded corners] at (0,0){\scriptsize$\Omega_b$};
\end{tkz}\,.
\]
By construction, $\rho_e=\lambda_{e^\star}$.
\end{notation}

\begin{lemma}
	\label{lem:lambda-rho}
	Let $a$, $b$, and $c$ be minimal 0-fields in $\C$ and let $e\in\undfldl{\C}{S^{n-1}\times D^1}[\fcj{r_b}\amalg r_a]$ and $e'\in\undfldl{\C}{S^{n-1}\times D^1}[\fcj{r_c}\amalg r_b]$.
	\begin{lst}
		\item\label{lr1} $\lambda_{r_a\times D^1}=1$ and $\lambda_{e\blt e'}=\lambda_{e}\lambda_{e'}$.
		\item\label{lr2} $\lambda_e\varnothing_a=\overline{\rho_e}\varnothing_b$ in $\undfldl[n]{\C}{S^n}$. 
        \item\label{lr15} $\lambda_e\psn^\C(\varnothing_a)=\overline{\rho_e}\psn^\C(\varnothing_b)$.
        \item\label{lr4} $\lambda_e=0$ if and only if $\rho_e=0$.
        \end{lst} 
        Now suppose moreover that $\lambda_e\neq0$ and set
        \[
        t_{a/b}\coloneq\psn^\C(\varnothing_a)/\psn^\C(\varnothing_b)>0.
        \]
        Then
        \begin{lst}[resume]
        \item\label{lr8} $\rho_e=t_{a/b}\overline{\lambda_e}$, 
        \item\label{lr5} $\varnothing_a=t_{a/b}\varnothing_b$ in $\undfldl[n]{\C}{S^n}$,
        \item\label{lr6} $\rN^\C(a)=t_{a/b}^2\rN^\C(b)$, and
        \item\label{lr7} $w_\C(a)=t_{a/b}^{-1}w_\C(b)$.
	\end{lst}
\end{lemma}

\begin{pf}
\itemstep{\ref{lr1}.} By \ref{CTc}, gluing $r_a\times D^1$ onto $\partial D^n$ gives $\Omega_a$. By associativity of gluing \ref{CGa}, $(e\blt e')\triangleright\Omega_c=e\triangleright(e'\triangleright\Omega_c)=\lambda_{e'}\lambda_e\Omega_a$.

\itemstep{\ref{lr2}.} Indeed, 
\begin{equation}\label{lr2eq}
\lambda_e\varnothing_a
=\Omega_a\blt(e\triangleright\Omega_b)
=(\Omega_a\blt e)\blt\Omega_b
=\overline{\rho_e}(\Omega_b\blt\Omega_b)
=\overline{\rho_e}\varnothing_b.
\end{equation}

\itemstep{\ref{lr15}.} Apply $\psn^\C$ to \ref{lr2}.

\itemstep{\ref{lr4}.} Since $\psn^\C(\varnothing_a),\psn^\C(\varnothing_b)>0$, \ref{lr15} gives that $\lambda_e=0$ if and only if $\rho_e=0$.

\itemstep{\ref{lr8}.} Dividing through \ref{lr15} by $\psn^\C(\varnothing_b)\neq0$ gives $\displaystyle\overline{\rho_e}=\frac{\psn^\C(\varnothing_a)}{\psn^\C(\varnothing_b)}\lambda_e=t_{a/b}\lambda_e$.

\itemstep{\ref{lr5}.} Substituting $\overline{\rho_e}=t_{a/b}\lambda_e$ into \eqref{lr2eq} gives $\lambda_e\varnothing_a=\overline{\rho_e}\varnothing_b=t_{a/b}\lambda_e\varnothing_b$. Now divide through by $\lambda_e$.

\itemstep{\ref{lr6}.}
By \ref{lr5}, $\rN^\C(a)=\bkt{\varnothing_a}{\varnothing_a}_{S^n}=\Vert{}t_{a/b}\varnothing_b\Vert{}^2=t_{a/b}^2\Vert{}\varnothing_b\Vert{}^2=t_{a/b}^2\rN^\C(b)$.

\itemstep{\ref{lr7}.} By \ref{lr6}, $\displaystyle w_\C(a)=\frac{\psn^\C(\varnothing_a)}{\rN^\C(a)}=\frac{t_{a/b}\psn^\C(\varnothing_b)}{t_{a/b}^2\rN^\C(b)}=t_{a/b}^{-1}\frac{\psn^\C(\varnothing_b)}{\rN^\C(b)}=t_{a/b}^{-1}w_\C(b)$.
\end{pf}

In what follows, we will write $\Sk_\C(S^{n-1})$ to mean the \defn{skein 1-category} of $\C$ on the $(n-1)$-sphere $S^{n-1}$, which is the dagger 1-category whose objects are the $(n-1)$-fields on $S^{n-1}$ and whose morphisms are the $n$-fields on $S^{n-1}\times D^1$, and $\Sk_\C(S^{n-1})^\cent$ denotes the completion of $\Sk_\C(S^{n-1})$, whose objects are pairs $(\bigoplus_i r_i,p)$ consisting of formal direct sums of objects $r_i\in\Sk_\C(S^{n-1})$ and a projection $p\in\End(\bigoplus_i r_i)$; see \cref{skein 1-category} and \cref{completeunitarycategories} for the details.

\begin{definition}\label{cnd0fld}
We will call two minimal 0-fields $a,b\in\fld[0]{\C}{\pt}$ \defn{connected}, written $a\sim b$, if there is a morphism $e\in\Sk_\C(S^{n-1})(r_b\to r_a)$ with $\lambda_e\neq0$.
\end{definition}

\begin{remark}
\label{rmk:DLcomponents}
Connectedness is an equivalence relation on minimal 0-fields: reflexivity holds since $\lambda_{r_a\times D^1}=1$ by \ref{lr1}; symmetry holds since $\lambda_{e^\dag}=\rho_e\neq0$ whenever $\lambda_e\neq0$, by \ref{lr4}; and transitivity holds since $\lambda_{e_1\blt e_2}=\lambda_{e_1}\lambda_{e_2}\neq0$ by \ref{lr1}. 

We will write $\pi_0\C$ to mean any set consisting of exactly one representative from each $\sim$-class of minimal 0-fields in $\C$; elements of $\pi_0\C$ will be called (connected) \defn{components} of $\C$.
\end{remark}

\begin{remark}
\label{rem:minimal-0-field-simple-object}
  Here we show that every minimal 0-field $a$ in $\C$ gives a canonical simple object $(r_a, p_a)$ in $\Sk_\C(S^{n-1})^\cent$. Observe that $\fldl[n]{\C}{D^n}[ r_a] = \bC \Omega_a$ as $\Omega_a \neq 0$ by \cref{bnz}, so $\bC \Omega_a$ is a simple $\End_{\Sk_\C(S^{n-1})}(r_a)$-module under the gluing action, which by \cref{notation-lambda-rho} satisfies $e \triangleright \Omega_a = \lambda_e \Omega_a$. Hence there is a minimal projection $p_a \in \End_{\Sk_\C(S^{n-1})}(r_a)$ with $p_a \triangleright \Omega_a = \Omega_a$. Thus $(r_a, p_a)$ is simple in $\Sk_\C(S^{n-1})^\cent$ since $\End_{\Sk_\C(S^{n-1})^\cent}((r_a,p_a))\defeq p_a \End_{\Sk_\C(S^{n-1})}(r_a) p_a \cong \bC$ by minimality of $p_a$.
\end{remark}

\begin{lemma}\label{lem:connected-iff-isomorphic}
  For minimal 0-fields $a, b \in \C$, $a \sim b$ if and only if $(r_a, p_a) \cong (r_b, p_b)$ in $\Sk_\C(S^{n-1})^\cent$.
\end{lemma}

\begin{pf}
\itemstep{($\Rightarrow$)} Suppose $a\sim b$ and choose $e\in\Sk_\C(S^{n-1})(r_b\to r_a)$ with $\lambda_{e} \neq 0$. Then $p_a e p_b \in \Sk_\C(S^{n-1})^\cent((r_b, p_b) \to (r_a, p_a))$ satisfies
\[
(p_a e p_b) \triangleright \Omega_b = p_a \triangleright(e \triangleright (p_b\triangleright\Omega_b)) = p_a\triangleright(\lambda_e\Omega_a)=\lambda_e\Omega_a \neq 0,
\]
so $p_a e p_b \neq 0$ and thus $(r_a, p_a) \cong (r_b, p_b)$ by Schur's Lemma \cite[Lem. 7.10.5]{UQSL}.

\itemstep{($\Leftarrow$)} Suppose $(r_a,p_a)\cong(r_b,p_b)$, say witnessed by the isomorphism $u\in\Sk_\C(S^{n-1})^\cent((r_a,p_a)\to (r_b,p_b))$. Then $u\blt u^{-1}=p_b\in\Sk_\C(S^{n-1})(r_b\to r_b)$, so
\[
\lambda_u\lambda_{u^{-1}}\Omega_b = \lambda_{u\blt u^{-1}}\Omega_b = (u\blt u^{-1})\triangleright\Omega_b = p_b\triangleright\Omega_b = \Omega_b.
\]
As $\Omega_b\neq 0$ forces $\lambda_{u^{-1}}\neq 0$, we conclude $u^{-1}\in\Sk_\C(S^{n-1})(r_b\to r_a)$ has $\lambda_{u^{-1}}\neq0$. Thus $a\sim b$ in $\C$.
\end{pf}

\begin{proposition}[Finiteness of components]
\label{thm:finiteness-components}
The set of components $\pi_0\C$ of a finite unitary disk-like $n$-category $\C$ satisfies
\[
|\pi_0\C|\leq\dim \undfldl[n]{\C}{S^n}.
\]
\end{proposition}

\begin{pf}
Consider the gluing data $\cG$ given by gluing $D^n_+\amalg D^n_-$ to $S^n$ along $S^{n-1}$, so that for each $r\in\undfld[n-1]{\C}{S^{n-1}}$, $V_r = \undfldl[n]{\C}{D^n_+\amalg D^n_-}[\fcj r\amalg r] \cong \undfldl[n]{\C}{D^n_-}[\fcj r]\otimes\undfldl[n]{\C}{D^n_+}[r]$. By \cref{UnitaryWalkerTheorem}, $\psn^\C$ gives finite unitary path integral data $\psi^{\cG}$ for $\cG$, so we can write $\{(r_i,p_i)\}_{i=1}^N=\Irr(\Sk_\C(S^{n-1})^\cent)$. For each $b\in\pi_0\C$, set $v_b\coloneq\fcj{\Omega_b}\amalg\Omega_b\in V_{r_b}$, so that $\varnothing_b=\Gamma(v_b)$. 

Since $p_b\triangleright\Omega_b=\Omega_b$, we have $v_b\in p_bV_{r_b}p_b$. As $(r_b,p_b)$ is simple in $\Sk_\C(S^{n-1})^\cent$ by \cref{rem:minimal-0-field-simple-object}, $(r_b,p_b)\cong(r_{i_b},p_{i_b})$ for some unique $i_b\in\{1,\dots,N\}$. Thus $\varnothing_b\in\Gamma(p_{i_b}V_{r_{i_b}}p_{i_b})$ by \cref{glu-corners}.

For distinct $b,b'\in\pi_0\C$ we have $(r_b,p_b)\not\cong(r_{b'},p_{b'})$ by \cref{lem:connected-iff-isomorphic}, so $i_b\neq i_{b'}$. Writing $\varnothing_b=\Gamma(w_b)$ and $\varnothing_{b'}=\Gamma(w_{b'})$ for $w_b\in p_{i_b}V_{r_{i_b}}p_{i_b}$ and $w_{b'}\in p_{i_{b'}}V_{r_{i_{b'}}}p_{i_{b'}}$, \cref{glu-ip} gives
\[
\bkt{\varnothing_b}{\varnothing_{b'}}_{S^n}=\frac{\delta_{i_b=i_{b'}}}{\bkt{p_{i_b}}{p_{i_b}}_{S^{n-1}\times I,\,\fcj{r_{i_b}}\amalg r_{i_b}}}\bkt{w_b}{w_{b'}}_{D^n_+\amalg D^n_-,\fcj{r_{i_b}}\amalg r_{i_b}}=0.
\]
As each $\varnothing_b$ is nonzero by \cref{posthings}, $\bkt{\varnothing_b}{\varnothing_b}_{S^n}=\rN^\C(b)>0$ by \cref{UnitaryWalkerTheorem}. It follows that $\{\varnothing_b\}_{b\in\pi_0\C}$ is an orthogonal set of nonzero vectors in the finite-dimensional space $\undfldl[n]{\C}{S^n}$, so $|\pi_0\C|\leq\dim\undfldl[n]{\C}{S^n}$.
\end{pf}

\begin{definition}
\label{def:nondegenerate}
A pure $k$-field $\xi$ in $\C$ is called \defn{nondegenerate} if the
product $n$-field $[\xi \times D^{n-k}]$ is nonzero.
\end{definition}

\begin{lemma}\label{lem:lambda-nonneg}
Let $a$ and $b$ be minimal $0$-fields of $\C$. Then for any
$\xi \in \fld[1]{\C}{D^1}[ \fcj b \amalg a]$, $\lambda_{e_\xi} \geq 0$, and $\lambda_{e_\xi} \neq 0$ if and only if $\xi$ is nondegenerate.
\end{lemma}

\begin{pf}
Indeed, by \ref{CS} we can split $e_\xi \triangleright \Omega_b$ into equivalent halves: where $\Gamma$ denotes the gluing map back to $D^n$, we have
\[
\lambda_{e_\xi}
\;
\begin{tkz}
\fill[fill=red!15, draw=red,region={grid, color=red!50, angle=0, distance=4pt}] (0,0) circle (1);
\node[fill=red!15,inner sep=0.5pt,rounded corners] at (0,0){\scriptsize$\Omega_a$};
\end{tkz}
\;\;=\;\;
\begin{tkz}
\def\tkzlen{0.5}
\def\tkzsep{0}
\def\tkzfibstep{0.05}
\fill[fill=blue!15, draw=blue,region={grid, color=blue!50, angle=0, distance=4pt}] (0,0) circle (\tkzlen);
\annulus[fill=violet!25]{(0,0)}{1*\tkzlen+1*\tkzsep}{2*\tkzlen+1*\tkzsep}
\pgfmathsetmacro{\rstart}{1*\tkzlen+1*\tkzsep}
\pgfmathsetmacro{\rnext}{\rstart+\tkzfibstep}
\pgfmathsetmacro{\rend}{2*\tkzlen+1*\tkzsep}
\foreach \tkztemp in {\rstart, \rnext, ..., \rend} {\draw[thin,violet] (0,0) circle (\tkztemp cm);}
\annulus[inner={blue}, outer={red}]{(0,0)}{1*\tkzlen+1*\tkzsep}{2*\tkzlen+1*\tkzsep}
\node[fill=blue!15,inner sep=0pt,rounded corners] at (0,0){\scriptsize$\Omega_b$};
\node[fill=violet!25,inner sep=0.5pt,rounded corners] at (-1.5*\tkzlen,0){\scriptsize$e_\xi$};
\end{tkz}
\;\;=\;\;
\Gamma\!\left(
\begin{tkz}
\def\tkzlen{0.5}
\def\tkzsep{0}
\def\tkzfibstep{0.05}
\clip(-2*\tkzlen,-2*\tkzlen)rectangle(0,2*\tkzlen);
\fill[fill=blue!15, draw=blue,region={grid, color=blue!50, angle=0, distance=4pt}] (0,0) circle (\tkzlen);
\annulus[fill=violet!25]{(0,0)}{1*\tkzlen+1*\tkzsep}{2*\tkzlen+1*\tkzsep}
\pgfmathsetmacro{\rstart}{1*\tkzlen+1*\tkzsep}
\pgfmathsetmacro{\rnext}{\rstart+\tkzfibstep}
\pgfmathsetmacro{\rend}{2*\tkzlen+1*\tkzsep}
\foreach \tkztemp in {\rstart, \rnext, ..., \rend} {\draw[thin,violet] (0,0) circle (\tkztemp cm);}
\annulus[inner={blue}, outer={red}]{(0,0)}{1*\tkzlen+1*\tkzsep}{2*\tkzlen+1*\tkzsep}
\end{tkz}
\quad
\begin{tkz}[xscale=-1]
\def\tkzlen{0.5}
\def\tkzsep{0}
\def\tkzfibstep{0.05}
\clip(-2*\tkzlen,-2*\tkzlen)rectangle(0,2*\tkzlen);
\fill[fill=blue!15, draw=blue,region={grid, color=blue!50, angle=0, distance=4pt}] (0,0) circle (\tkzlen);
\annulus[fill=violet!25]{(0,0)}{1*\tkzlen+1*\tkzsep}{2*\tkzlen+1*\tkzsep}
\pgfmathsetmacro{\rstart}{1*\tkzlen+1*\tkzsep}
\pgfmathsetmacro{\rnext}{\rstart+\tkzfibstep}
\pgfmathsetmacro{\rend}{2*\tkzlen+1*\tkzsep}
\foreach \tkztemp in {\rstart, \rnext, ..., \rend} {\draw[thin,violet] (0,0) circle (\tkztemp cm);}
\annulus[inner={blue}, outer={red}]{(0,0)}{1*\tkzlen+1*\tkzsep}{2*\tkzlen+1*\tkzsep}
\end{tkz}
\right)
\;\underset{\ref{CTc}}{=}\;
\Gamma(\fcj{w_\xi} \blt_{D^{n-1}} w_\xi)
\]
where $w_\xi \coloneq \xi \times D^{n-1}$. Applying $\psn^\C\circ\cl_a$ gives $\lambda_{e_\xi}\, \psn^\C(\varnothing_a)
= \bkt{w_\xi}{w_\xi}_{D^n,\partial w_\xi}$. Since $\bkt{w_\xi}{w_\xi}_{D^n,\partial w_\xi}\geq0$ by \ref{psnP} and $\psn^\C(\varnothing_a) > 0$ by \cref{posthings}, the claim follows.
\end{pf}

\begin{proposition}\label{lemcon}
Let $a$ and $b$ be minimal $0$-fields of $\C$. The following are equivalent.
\begin{lst}
\item[(a)]\label{lemcona}
$a$ and $b$ are connected.
\item[(b)]\label{lemconb}
There is a pure 1-field $\xi \in \fld[1]{\C}{D^1}[ \fcj b \amalg a]$ such that
$e_\xi \coloneq \xi \times S^{n-1}$ satisfies $e_\xi \triangleright \Omega_b = \lambda_{e_\xi}\Omega_a$ with $\lambda_{e_\xi} \neq 0$.
\item[(c)]\label{lemconc}
There is some nondegenerate $\xi \in \fld[1]{\C}{D^1}[ \fcj b \amalg a]$.
\end{lst}
\end{proposition}

\begin{pf}
\itemstep{\ref{lemconb}$\Rightarrow$\ref{lemcona}.} This is immediate.

\itemstep{\ref{lemconc}$\Rightarrow$\ref{lemconb}.} This is \cref{lem:lambda-nonneg}.   

\itemstep{\ref{lemcona}$\Rightarrow$\ref{lemconc}.} Suppose $a \sim b$, say witnessed by
$e \in \undfldl[n]{\C}{A}[ \fcj{r_b}\amalg r_a]$ where $A\coloneq S^{n-1}\times D^1$. Since pure $n$-fields span, we may assume without loss of generality that $e$ is a pure $n$-field with $[e] \neq 0$. Identifying the outer boundary $S^{n-1}\times\{1\}\subset A$ with $S^{n-1}$ sitting inside $\bR^n$, choose an ordered basis $(e_1,\dots,e_n)$ of $\bR^n$ transverse to some string locus $S_e$ of $e$ in the sense that each cell of $S_e$ intersects $(S^{n-1}\cap\bR\{e_1,\dots,e_k\})\times D^1$ transversely for $1\leq k<n$. For $0 \leq k \leq n-1$, set
\[
E_k \coloneq S^{n-1} \cap \bR\{e_1, \dots, e_{k+1}\} \cap \{x\in\bR^n\mid x_{k+1} \geq 0\},   
\]
so that $E_k$ is a $k$-ball. Define the $(k+1)$-field
\[
\xi_k \coloneq e|_{E_k \times D^1}\in\fld[k+1]{\C}{E_k \times D^1}
\]
and define the (quotient-level) $n$-field $\theta_k \coloneq \xi_k \times D^{n-1-k}$. Observe that $\theta_k\neq 0$ if and only if $\xi_k$ is nondegenerate. Thus to obtain a nondegenerate 1-field $\xi$, it suffices to show $\theta_0\neq 0$. We show this by descending coinductively from $k=n-1$ to $k=0$.

\itemstep{Base case.} We claim $\theta_{n-1}\neq 0$. Since $E_{n-1}$ is a hemisphere of $S^{n-1}$, the decomposition $A = (E_{n-1}\times D^1)\cup_{\partial E_{n-1}\times D^1}(\overline{(S^{n-1}\setminus E_{n-1})}\times D^1)$ writes $e$ as a gluing of $\xi_{n-1}$ with its complement. If $\xi_{n-1}\in\lU{\C}$, then \ref{CUi} would give $e\in\lU{\C}$, contradicting $[e]\neq 0$; thus $\theta_{n-1}=\xi_{n-1}\neq 0$.

\itemstep{Inductive step.}
We claim that for $1\leq k\leq n-1$, $\theta_k\neq 0$ implies $\theta_{k-1}\neq 0$. The inclusion $E_{k-1}\hookrightarrow\partial E_k$ embeds $Y\coloneqq E_{k-1}\times D^1\times D^{n-k-1}$ as an $(n-1)$-ball into $\partial(E_k\times D^1\times D^{n-k-1})$, along which the boundary field of $\theta_k$ restricts to $\xi_{k-1}\times D^{n-k-1}$ by \ref{Cpirk}. Applying \cref{lem:nondeg-criterion} to $\theta_k\neq 0$ forces $\xi_{k-1}\times D^{n-k-1}$ to be nondegenerate, whence $\theta_{k-1} = \xi_{k-1}\times D^{n-k}\neq 0$, as desired. It follows that $\theta_0\neq 0$, so $\xi\coloneq\xi_0\in\fld[1]{\C}{D^1}[\fcj{b}\amalg a]$ is nondegenerate.
\end{pf}

\begin{lemma}
\label{lem:nondeg-criterion}
For an $n$-ball $X$, a field $c \in \undfld[n-1]{\C}{\partial X}$, and an $(n-1)$-ball $Y \hookrightarrow \partial X$ along which $c$ splits, if an
$n$-field $\alpha \in \fldl[n]{\C}{X}[c]$ is nonzero then $c|_Y$ is nondegenerate.
\end{lemma}

\begin{pf}
The contrapositive holds: by \ref{CTc}, $\alpha = \psi_{Y,D^1*}(\alpha \blt_Y [c|_Y \times D^1])$ for any collaring homeomorphism $\psi_{Y,D^1}$, so
\begin{align*}
c|_Y \text{ is degenerate}
&\Longrightarrow
[c|_Y \times D^1] = 0
% \\&
\underset{\ref{CUi}}{\Longrightarrow}
\alpha \blt_Y [c|_Y \times D^1] = 0
\\
&\underset{\ref{CUvarphi}}{\Longrightarrow}
\alpha
=
\psi_{Y,D^1*}(\alpha \blt_Y [c|_Y \times D^1])
=
0.\qedhere
\end{align*}
\end{pf}

The proof of the following lemma is the general-$n$ disk-like analog of the naturality assertions of \cref{dag-thing-i-2-pf} and \cref{dag-thing-ii-3-pf}.

\begin{lemma}[Naturality for \texorpdfstring{$(n,n)$}{(n,n)}-transfors]
	\label{lem:transfor-naturality}
    Let $\C$ and $\D$ be any disk-like $n$-categories, let $X$ be an $n$-ball, and consider a pure $(n,n)$-transfor $\eS\in\fld[n]{\eHom(\C{\to}\D)}{X}$
	with boundary $(n,n-1)$-transfor $\eN\coloneq\partial\eS\in\undfld[n-1]{\eHom(\C{\to}\D)}{\partial X}$. Then for all pure 1-fields $\xi\in\fld[1]{\C}{D^1}[\fcj{b}\amalg a]$ and all splittings
	$\partial X=Y\cup_{S^{n-2}}Y'$ into $(n-1)$-balls $Y$ and $Y'$ along which $\eN(\xi)$ splits,
	\begin{equation}
		\label{eq:transfor-naturality}
		\eS(b)\blt_Y\eN(\xi)|_{D^1\times Y}
		=
		\eS(a)\blt_{Y'}\eN(\xi)|_{D^1\times Y'}
	\end{equation}
	in $\undfldl[n]{\D}{X'}[\eN(a)|_Y\cup\eN(b)|_{Y'}]$ where $X'\cong X$ denotes $X$ glued to the corresponding collar.
\end{lemma}

\begin{pf}
This is an instance of the naturality axiom \ref{axiom:eT-naturality}. Applying it to $\eS$ with $W=D^1$, the associated $n$-sphere field is $(\fcj{\eS(b)}\amalg\eS(a))\blt\eN(\xi)$ on $(\partial D^1\times X)\cup(D^1\times\partial X)\cong S^n$. The unique splitting of $\partial D^1$ together with the given splitting $\partial X=Y\cup_{S^{n-2}}Y'$ cuts this sphere into two halves, and \ref{axiom:eT-naturality} identifies the fields on them, which gives \eqref{eq:transfor-naturality}.
\end{pf}

\begin{remark}
\cref{thm:intertwining} below is the general-$n$ disk-like analog of the identity
$d_b^{-1}\Psi^{\fY}_{F(b)}(\mu_b)=d_c^{-1}\Psi^{\fY}_{F(c)}(\mu_c)$ when $b$ and $c$ are connected simple objects in a $3$-Hilbert space $\fX$ (so $n=2$), which we use in \cite[Step 1 of Lem. 5.2]{bases} to show independence of the ONB in defining the spherical weight on
$\Hom(\fX\to\fY)$. And indeed, in \cref{cor:weight-balanced} below we will use \cref{thm:intertwining} to prove that \eqref{functorsphereweight} is independent of the choice of representatives of $\sim$-classes $\pi_0\C$ of minimal 0-fields of $\C$ for general $n$ as part of proving \cref{thm:eFun-sphere-weight}.
\end{remark}

\begin{lemma}
\label{thm:intertwining}
Let $\C$ and $\D$ be unitary disk-like $n$-categories, and let $a$ and $b$ be connected minimal 0-fields in $\C$.
Then for any (quotient-level) $(n,n)$-transfor $\mu\in\undfldl[n]{\eHom(\C{\to}\D)}{S^n}$,
\[
\mu_a=t_{a/b}\mu_b\qquad\text{in}\qquad \undfldl[n]{\D}{S^n}.
\]
\end{lemma}

\begin{pf}
Both sides are linear in $\mu$, so we may assume $\mu$ is pure. Consider a three-piece decomposition of $S^n$ consisting of two $n$-balls $X_-\cong D^n$ and $X_+\cong \orev{D^n}$ and an annulus $A\cong S^{n-1}\times D^1$. By extended-isotopy invariance \ref{CTc} in $\eFun(\C{\to}\D)$, we may further assume $\mu$ restricts on $X_-$ and $X_+$ to product fields $\mu|_{X_-}=\eG\times X_-$ and $\mu|_{X_+}=\eF\times X_+$ for disk-like functors $\eF,\eG\colon\C\to\D$; in particular $\mu_b|_{X_-}=\eG(\Omega_b)$ and $\mu_a|_{X_+}=\eF(\fcj{\Omega_a})$ by \ref{axiom:eT-products}.

Now set $Y_-\coloneq\partial X_-$, $E_-\coloneq S^n\setminus\Int X_-$, $Y_+\coloneq\partial X_+$, $E_+\coloneq S^n\setminus\Int X_+$, and fix collaring homeomorphisms $S^n\cong X_-\cup_{Y_-}C_-\cup_{Y_-}E_-$ and $S^n\cong E_+\cup_{Y_+}C_+\cup_{Y_+}X_+$ where $C_\pm\coloneq D^1\times Y_\pm$. By \cref{lemcon}, there is a pure 1-field $\xi$ in $\C$ with $\partial\xi=\fcj{b}\amalg a$ such that $\lambda_{e_\xi}\neq 0$ for $e_\xi\coloneq\xi\times S^{n-1}$. 

In the illustrations below, a shading indicates a field in $\C$; $a$ is red, $b$ is blue, and $\xi$ is violet:
\[
\Omega_a
=
\begin{tkz}[scale=0.75]\fill[red!15,rounded corners,region={grid, color=red!50,angle=0}](0,0)rectangle(1,1);\end{tkz}
\,,
\qquad
\Omega_b
=
\begin{tkz}[scale=0.75]\fill[blue!15,rounded corners,region={grid, color=blue!50,angle=0}](0,0)rectangle(1,1);\end{tkz}
\,,
\qquad
\text{and}
\,
\qquad
\xi\times D^{n-1}
=
\begin{tkz}[scale=0.75]\fill[violet!25,rounded corners,region={lines,angle=90,color=Purple,distance=2pt}](0,0)rectangle(1,1);\end{tkz}\,.
\]
As a disk-like $(n,n)$-transfor is the $n$-fold iterate of the sequence ``natural transformation, modification, perturbation, \ldots'', we represent disk-like functors ($(n,0)$-transfors) below as regular black overlay patterns and the $(n,n)$-transfor $\mu$ as a more ``chaotic'' rainbow pattern, so that a black or rainbow pattern overlaid on a field indicates the functor or $\mu$ respectively applied to that field: 
\[
\eF=\begin{tkz}[scale=0.75]\fill[rounded corners,primedregion=white,draw=black,dotted](0,0)rectangle(1,1);\end{tkz}
\,,\qquad
\eG=\begin{tkz}[scale=0.75]\fill[rounded corners,boxregion=white,draw=black,dotted](0,0)rectangle(1,1);\end{tkz}\,,\qquad\text{and}\,\qquad
\mu=\begin{tkz}[scale=0.75]\fill[white,rounded corners,region={hyperchaos, color=white},draw=black,dotted](0,0)rectangle(1,1);\end{tkz}\,.
\]

Since $Y_\pm\subset X_\pm$, the restrictions are the product $(n,n-1)$-transfors $\mu|_{Y_-}=\eG\times Y_-$ and $\mu|_{Y_+}=\eF\times Y_+$. Thus as $n$-fields in $\D$ on $D^1\times Y_\pm$, $\mu|_{Y_-}(\xi)=\eG(\xi)\times Y_-=\eG(e_\xi)$ and similarly $\mu|_{Y_+}(\xi)=\eF(e_\xi)$. By \ref{axiom:eT-boundary} we have $\partial(\mu|_{Y_\pm}(\xi))=\fcj{\mu_b|_{Y_\pm}}\amalg\mu_a|_{Y_\pm}$.
Set
\[
\Lambda_-\coloneq \eG(\Omega_b)\blt_{Y_-}\eG(e_\xi)\blt_{Y_-}(\mu_a|_{E_-})
\qquad\text{and}\qquad
\Lambda_+\coloneq (\mu_b|_{E_+})\blt_{Y_+}\eF(e_\xi)\blt_{Y_+}\eF(\fcj{\Omega_a}).
\]
Applying \cref{lem:transfor-naturality} to each ball in a ball decomposition of $\mu$ shifts the collar $C_-$ across the support of $\mu$ onto $C_+$, so $\Lambda_-=\Lambda_+$ in $\undfldl[n]{\D}{S^n}$. 
Now
\begin{align*}
\Lambda_-
=\begin{tkz}
\def\tkzlen{0.5}
\def\tkzsep{0.25}
\def\tkzfibstep{0.05}
\fill[fill=blue!15, draw=blue,region={grid, color=blue!50, angle=0, distance=4pt},region={box, color=black}] (0,0) circle (\tkzlen);
\end{tkz}
\;\blt\;
\begin{tkz}
\def\tkzlen{0.5}
\def\tkzsep{0.25}
\def\tkzfibstep{0.05}
\annulus[fill=violet!25]{(0,0)}{1*\tkzlen+1*\tkzsep}{2*\tkzlen+1*\tkzsep}
\pgfmathsetmacro{\rstart}{1*\tkzlen+1*\tkzsep}
\pgfmathsetmacro{\rnext}{\rstart+\tkzfibstep}
\pgfmathsetmacro{\rend}{2*\tkzlen+1*\tkzsep}
\foreach \tkztemp in {\rstart, \rnext, ..., \rend} {\draw[thin,violet] (0,0) circle (\tkztemp cm);}
\annulus[region={box},inner={blue}, outer={red}]{(0,0)}{1*\tkzlen+1*\tkzsep}{2*\tkzlen+1*\tkzsep}
\end{tkz}
\;\blt\;
\begin{tkz}
\def\tkzlen{0.5}
\def\tkzsep{0.25}
\def\tkzfibstep{0.05}
\fill[fill=red!15, draw=red,region={hyperchaos, color=red!15}] (0,0) circle (\tkzlen);
\end{tkz}
% \\&
&=
\begin{tkz}
\def\tkzlen{0.5}
\def\tkzsep{0}
\def\tkzfibstep{0.05}
\fill[fill=blue!15, draw=blue,region={grid, color=blue!50, angle=0, distance=4pt}] (0,0) circle (\tkzlen);
\annulus[fill=violet!25]{(0,0)}{1*\tkzlen+1*\tkzsep}{2*\tkzlen+1*\tkzsep}
\pgfmathsetmacro{\rstart}{1*\tkzlen+1*\tkzsep}
\pgfmathsetmacro{\rnext}{\rstart+\tkzfibstep}
\pgfmathsetmacro{\rend}{2*\tkzlen+1*\tkzsep}
\foreach \tkztemp in {\rstart, \rnext, ..., \rend} {\draw[thin,violet] (0,0) circle (\tkztemp cm);}
\annulus[inner={blue}, outer={red}]{(0,0)}{1*\tkzlen+1*\tkzsep}{2*\tkzlen+1*\tkzsep}
\path[region={box, color=black}](0,0) circle (2*\tkzlen+1*\tkzsep);
\end{tkz}
\;\blt\;
\begin{tkz}
\def\tkzlen{0.5}
\def\tkzsep{0.25}
\def\tkzfibstep{0.05}
\fill[fill=red!15, draw=red,region={hyperchaos, color=red!15}] (0,0) circle (\tkzlen);
\end{tkz}
\\&
=
\begin{tkz}
\def\tkzlen{0.5}
\def\tkzsep{0}
\def\tkzfibstep{0.05}
\fill[fill=red!15, draw=red,region={grid, color=red!50, angle=0, distance=4pt}] (0,0) circle (\tkzlen+\tkzsep+\tkzlen);
\end{tkz}
\;\blt\;
\lambda_{e_\xi}\;
\begin{tkz}
\def\tkzlen{0.5}
\def\tkzsep{0.25}
\def\tkzfibstep{0.05}
\fill[fill=red!15, draw=red,region={hyperchaos, color=red!15}] (0,0) circle (\tkzlen);
\end{tkz}
=
\lambda_{e_\xi}\,
\begin{tkz}[scale=0.5]
\def\tkzlen{0.5}
\fill[shading=ball,ball color=red!15] (0,0) circle (2*\tkzlen);
\path[fill=red!15,region={hyperchaos},opacity=0.5] (0,0) circle (2*\tkzlen);
\node at (0,0){$\mu_a$};
\end{tkz}\,.
\end{align*}
Similarly,
\begin{align*}
\Lambda_+
=
\begin{tkz}
\def\tkzlen{0.5}
\def\tkzsep{0.25}
\def\tkzfibstep{0.05}
\fill[fill=blue!15, draw=blue,region={hyperchaos, color=blue!15}] (0,0) circle (\tkzlen);
\end{tkz}
\;\blt\;
\begin{tkz}
\def\tkzlen{0.5}
\def\tkzsep{0.25}
\def\tkzfibstep{0.05}
\annulus[fill=violet!25]{(0,0)}{1*\tkzlen+1*\tkzsep}{2*\tkzlen+1*\tkzsep}
\pgfmathsetmacro{\rstart}{1*\tkzlen+1*\tkzsep}
\pgfmathsetmacro{\rnext}{\rstart+\tkzfibstep}
\pgfmathsetmacro{\rend}{2*\tkzlen+1*\tkzsep}
\foreach \tkztemp in {\rstart, \rnext, ..., \rend} {\draw[thin,violet] (0,0) circle (\tkztemp cm);}
\annulus[region={dots, color=black},inner={blue}, outer={red}]{(0,0)}{1*\tkzlen+1*\tkzsep}{2*\tkzlen+1*\tkzsep}
\end{tkz}
\;\blt\;
\begin{tkz}
\def\tkzlen{0.5}
\def\tkzsep{0.25}
\def\tkzfibstep{0.05}
\fill[fill=red!15, draw=red,region={grid, color=red!50, angle=0, distance=4pt},region={dots, color=black}] (0,0) circle (\tkzlen);
\end{tkz}
&=
\begin{tkz}
\def\tkzlen{0.5}
\def\tkzsep{0.25}
\def\tkzfibstep{0.05}
\fill[fill=blue!15, draw=blue,region={hyperchaos, color=blue!15}] (0,0) circle (\tkzlen);
\end{tkz}
\;\blt\;
\begin{tkz}
\def\tkzlen{0.5}
\def\tkzsep{0}
\def\tkzfibstep{0.05}
\fill[fill=red!15, draw=red,region={grid, color=red!50, angle=0, distance=4pt}] (0,0) circle (\tkzlen);
\annulus[fill=violet!25]{(0,0)}{1*\tkzlen+1*\tkzsep}{2*\tkzlen+1*\tkzsep}
\pgfmathsetmacro{\rstart}{1*\tkzlen+1*\tkzsep}
\pgfmathsetmacro{\rnext}{\rstart+\tkzfibstep}
\pgfmathsetmacro{\rend}{2*\tkzlen+1*\tkzsep}
\foreach \tkztemp in {\rstart, \rnext, ..., \rend} {\draw[thin,violet] (0,0) circle (\tkztemp cm);}
\annulus[inner={red}, outer={blue}]{(0,0)}{1*\tkzlen+1*\tkzsep}{2*\tkzlen+1*\tkzsep}
\path[region={dots, color=black}](0,0) circle (2*\tkzlen+1*\tkzsep);
\end{tkz}
\\&
=
\begin{tkz}
\def\tkzlen{0.5}
\def\tkzsep{0.25}
\def\tkzfibstep{0.05}
\fill[fill=blue!15, draw=blue,region={hyperchaos, color=blue!15}] (0,0) circle (\tkzlen);
\end{tkz}
\;\blt\;
\rho_{e_\xi}\;
\begin{tkz}
\def\tkzlen{0.5}
\def\tkzsep{0}
\def\tkzfibstep{0.05}
\fill[fill=blue!15, draw=blue,region={grid, color=blue!50, angle=0, distance=4pt},region={dots, color=black}] (0,0) circle (\tkzlen+\tkzsep+\tkzlen);
\end{tkz}
=
\rho_{e_\xi}\;
\begin{tkz}[scale=0.5]
\def\tkzlen{0.5}
\fill[shading=ball,ball color=blue!15] (0,0) circle (2*\tkzlen);
\path[fill=blue!15,region={hyperchaos},opacity=0.5] (0,0) circle (2*\tkzlen);
\node at (0,0){$\mu_b$};
\end{tkz}\,.
\end{align*}
Thus $\lambda_{e_\xi}\mu_a=\rho_{e_\xi}\mu_b$, so if $\lambda_{e_\xi}\neq0$ then $\mu_a=\lambda_{e_\xi}^{-1}\rho_{e_\xi}\mu_b=t_{a/b}\mu_b$ by \cref{lem:lambda-rho}\ref{lr8} and \cref{lem:lambda-nonneg}.
\end{pf}

In the above proof, observe that minimality of $a$ was used only in the computation of $\Lambda_-$. It follows that for arbitrary 0-fields $x$ in $\C$, we have
\begin{equation}\label{eq:half-intertwining}
\eG(e_\xi\triangleright\Omega_b)\blt_{Y_-}(\mu_x|_{E_-})=\Lambda_-=\Lambda_+=\rho_{e_\xi}\mu_b.
\end{equation}
Diagrammatically,
\[
\begin{tkz}
\def\tkzlen{0.5}
\def\tkzsep{0}
\def\tkzfibstep{0.05}
\fill[fill=blue!15, draw=blue,region={grid, color=blue!50, angle=0, distance=4pt}] (0,0) circle (\tkzlen);
\annulus[fill=violet!25]{(0,0)}{1*\tkzlen+1*\tkzsep}{2*\tkzlen+1*\tkzsep}
\pgfmathsetmacro{\rstart}{1*\tkzlen+1*\tkzsep}
\pgfmathsetmacro{\rnext}{\rstart+\tkzfibstep}
\pgfmathsetmacro{\rend}{2*\tkzlen+1*\tkzsep}
\foreach \tkztemp in {\rstart, \rnext, ..., \rend} {\draw[thin,violet] (0,0) circle (\tkztemp cm);}
\annulus[inner={blue}, outer={gray,thick}]{(0,0)}{1*\tkzlen+1*\tkzsep}{2*\tkzlen+1*\tkzsep}
\path[region={box, color=black}](0,0) circle (2*\tkzlen+1*\tkzsep);
\end{tkz}
\,\blt\,
\begin{tkz}
\def\tkzlen{0.5}
\def\tkzsep{0.25}
\def\tkzfibstep{0.05}
\fill[fill=gray!15, draw=gray,thick,region={hyperchaos, color=gray!15}] (0,0) circle (\tkzlen);
\end{tkz}
% \\&
\quad=\quad
\rho_{e_\xi}\;
\begin{tkz}[scale=0.5]
\def\tkzlen{0.5}
\fill[shading=ball,ball color=blue!15] (0,0) circle (2*\tkzlen);
\path[fill=blue!15,region={hyperchaos},opacity=0.5] (0,0) circle (2*\tkzlen);
\node at (0,0) {$\mu_b$};
\end{tkz}\;.
\]

\begin{proposition}
\label{cor:weight-balanced}
Let $\C$ and $\D$ be finite unitary disk-like $n$-categories and let $a$ and $b$ be minimal 0-fields in $\C$ with $a\sim b$. Then for any $(n,n)$-transfor $\mu\in\undfldl[n]{\eHom(\C{\to}\D)}{S^n}$,
\begin{equation}\label{eq:balance}
w_\C(a)\psn^\D(\mu_a)=w_\C(b)\psn^\D(\mu_b).
\end{equation}
\end{proposition}

\begin{pf}
Since $a\sim b$, 
\[
w_\C(a)\psn^\D(\mu_a) \underset{\eqref{thm:intertwining}\&\eqref{lem:lambda-rho}}{=} t_{a/b}^{-1}w_\C(b)t_{a/b}\psn^\D(\mu_b) = w_\C(b)\psn^\D(\mu_b).\qedhere
\]
\end{pf}

\begin{definition}[Weak 0-completeness]
\label{def:0-complete}
A disk-like $n$-category $\C$ is \defn{weakly 0-complete} if for every 0-field $x\in\fld[0]{\C}{\pt}$ there are minimal 0-fields $b_1,\dots,b_k$, 1-fields $\xi_i\in\fld[1]{\C}{D^1}[ \fcj{b_i}\amalg x]$, and scalars $c_i\in\bC$ such that $\Omega_x=\sum_{i=1}^k c_i\, e_i\triangleright\Omega_{b_i}$ in $\fldl[n]{\C}{D^n}[r_x]$, where $e_i\coloneq\xi_i\times S^{n-1}$. Diagrammatically,
\begin{equation}\label{eq:0complete}
\begin{tkz}[scale=0.5]
\fill[gray!15,region={grid, color=gray, angle=0}] (0,0) circle (1);
\draw[gray,thick](0,0)circle(1);
\node[fill=gray!15,inner sep=0pt,rounded corners] at (0,0){\scriptsize$\Omega_x$};
\end{tkz}
\quad = \quad
\sum_{i=1}^k c_i\;
\begin{tkz}
\def\tkzlen{0.5}
\def\tkzsep{0}
\def\tkzfibstep{0.05}
\fill[fill=blue!15, draw=blue,region={grid, color=blue!50, angle=0, distance=4pt}] (0,0) circle (\tkzlen);
\annulus[fill=violet!25,outer={red}]{(0,0)}{1*\tkzlen+1*\tkzsep}{2*\tkzlen+1*\tkzsep}
\pgfmathsetmacro{\rstart}{1*\tkzlen+1*\tkzsep}
\pgfmathsetmacro{\rnext}{\rstart+\tkzfibstep}
\pgfmathsetmacro{\rend}{2*\tkzlen+1*\tkzsep}
\foreach \tkztemp in {\rstart, \rnext, ..., \rend} {\draw[thin,violet] (0,0) circle (\tkztemp cm);}
\annulus[inner={blue}, outer={gray,thick}]{(0,0)}{1*\tkzlen+1*\tkzsep}{2*\tkzlen+1*\tkzsep}
\node[fill=blue!15,inner sep=0pt,rounded corners] at (0,0){\scriptsize$\Omega_{b_i}$};
\node[fill=violet!25,inner sep=1pt,rounded corners] at (-1.5*\tkzlen,0){\scriptsize$e_i$};
\end{tkz}\;.
\end{equation}
\end{definition}

Note that weak 0-completeness is a weaker condition than the notion of \defn{weak completeness} from \cite{W21}.

\begin{lemma}
\label{lem:0-complete-transport}
Let $\C$ and $\D$ be disk-like $n$-categories with $\C$ weakly 0-complete, let $x\in\fld[0]{\C}{\pt}$, and let $\mu\in\undfldl[n]{\eHom(\C{\to}\D)}{S^n}$. Then in the notation of \cref{def:0-complete}, 
\[
    \mu_x = \sum_{i=1}^k c_i\,\rho_{e_i}\mu_{b_i} \quad \text{in} \quad \undfldl[n]{\D}{S^n}.
\]
Diagrammatically,
\[
\begin{tkz}[scale=0.5]
\def\tkzlen{0.5}
\fill[shading=ball,ball color=gray!15] (0,0) circle (2*\tkzlen);
\path[fill=gray!15,region={hyperchaos},opacity=0.5] (0,0) circle (2*\tkzlen);
\node at (0,0) {$\mu_x$};
\end{tkz}
=
\sum_{i=1}^k c_i \rho_{e_i}\;
\begin{tkz}[scale=0.5]
\def\tkzlen{0.5}
\fill[shading=ball,ball color=blue!15] (0,0) circle (2*\tkzlen);
\path[fill=blue!15,region={hyperchaos},opacity=0.5] (0,0) circle (2*\tkzlen);
\node at (0,0) {$\mu_{b_i}$};
\end{tkz}
\]
\end{lemma}

\begin{pf}
Both sides are linear in $\mu$, so we may assume $\mu$ is pure. Then in the notation of \cref{def:0-complete}, we have
\[
\def\tkzlen{0.5}
\def\tkzsep{0.25}
\begin{tkz}
\def\tkzlen{0.25}
\fill[shading=ball,ball color=gray!15] (0,0) circle (2*\tkzlen);
\path[fill=gray!15,region={hyperchaos},opacity=0.5] (0,0) circle (2*\tkzlen);
\node at (0,0){$\mu_x$};
\end{tkz}
\quad
=
\quad
\begin{tkz}
\def\tkzlen{0.5}
\def\tkzsep{0.25}
\def\tkzfibstep{0.05}
\fill[fill=gray!15, draw=gray,thick,region={grid, color=gray!50, angle=0, distance=4pt},region={box, color=black}] (0,0) circle (\tkzlen);
\end{tkz}
\,\blt\,
\begin{tkz}
\def\tkzlen{0.5}
\def\tkzsep{0.25}
\def\tkzfibstep{0.05}
\fill[fill=gray!15, draw=gray,thick,region={hyperchaos, color=gray!15}] (0,0) circle (\tkzlen);
\end{tkz}
\quad
\underset{\eqref{eq:0complete}}{=}
\quad
\sum_{i=1}^k c_i\;
\begin{tkz}
\def\tkzlen{0.5}
\def\tkzsep{0}
\def\tkzfibstep{0.05}
\fill[fill=blue!15, draw=blue,region={grid, color=blue!50, angle=0, distance=4pt}] (0,0) circle (\tkzlen);
\annulus[fill=violet!25]{(0,0)}{1*\tkzlen+1*\tkzsep}{2*\tkzlen+1*\tkzsep}
\pgfmathsetmacro{\rstart}{1*\tkzlen+1*\tkzsep}
\pgfmathsetmacro{\rnext}{\rstart+\tkzfibstep}
\pgfmathsetmacro{\rend}{2*\tkzlen+1*\tkzsep}
\foreach \tkztemp in {\rstart, \rnext, ..., \rend} {\draw[thin,violet] (0,0) circle (\tkztemp cm);}
\annulus[inner={blue}, outer={gray,thick}]{(0,0)}{1*\tkzlen+1*\tkzsep}{2*\tkzlen+1*\tkzsep}
\path[region={box, color=black}](0,0) circle (2*\tkzlen+1*\tkzsep);
\end{tkz}
\,\blt\,
\begin{tkz}
\def\tkzlen{0.5}
\def\tkzsep{0.25}
\def\tkzfibstep{0.05}
\fill[fill=gray!15, draw=gray,thick,region={hyperchaos, color=gray!15}] (0,0) circle (\tkzlen);
\end{tkz}
% \\&
\quad
\underset{\eqref{eq:half-intertwining}}{=}
\quad
\sum_{i=1}^k c_i\,\rho_{e_i}\;
\begin{tkz}[scale=0.5]
\def\tkzlen{0.5}
\fill[shading=ball,ball color=blue!15] (0,0) circle (2*\tkzlen);
\path[fill=blue!15,region={hyperchaos},opacity=0.5] (0,0) circle (2*\tkzlen);
\node at (0,0){$\mu_{b_i}$};
\end{tkz}\;.\qedhere
\]
\end{pf}

\begin{theorem}[Sphere trace on {\texorpdfstring{$\eHom(\C{\to}\D)$}{Hom(CtoD)}}]
	\label{thm:eFun-sphere-weight}
	If $\C$ and $\D$ are unitary disk-like $n$-categories and $\C$ is finite and weakly 0-complete, then
	\begin{equation}
		\label{functorsphereweight}
		\psn^{\eHom}(\mu)\coloneq\sum_{a\in\pi_0\C}w_\C(a)\psn^\D(\mu_a)
	\end{equation}
	is a well-defined linear functional $\undfldl[n]{\eHom(\C{\to}\D)}{S^n}\to\bC$ satisfying \ref{psnP} that makes $(\eHom(\C{\to}\D),\psn^{\eHom})$ a unitary disk-like $n$-category.
\end{theorem}

\begin{pf}
\itemstep{Well-definedness.}
The sum in \eqref{functorsphereweight} is independent of the choice of representative minimal 0-fields by \cref{cor:weight-balanced}. Its index set $\pi_0\C$ is finite by \cref{thm:finiteness-components}. 

\itemstep{\ref{psnP}.}
For $m\in\undfldl[n]{\eHom(\C{\to}\D)}{D^n}[c]$, set $\mu\coloneq \fcj{m}\blt m\in\undfldl[n]{\eHom(\C{\to}\D)}{S^n}$. Then
\[
    \bkt mm_{D^n,c}=\psn^{\eHom}(\mu)=\sum_{b\in\pi_0\C}w_\C(b)\bkt{m_b}{m_b}_{D^n,c_b}\geq0,
\]
since each $w_\C(b)>0$ by \cref{posthings} and each $\bkt{m_b}{m_b}_{D^n}\geq0$ by \ref{psnP} for $\D$. If $\bkt mm_{D^n,c}=0$, then $m_a=0$ for all $a\in\pi_0\C$, whence $\mu_a=0$. For any minimal 0-field $b$ with $b\sim a\in\pi_0\C$, \cref{thm:intertwining} gives $\mu_b=t_{b/a}\mu_a=0$. For an arbitrary 0-field $x$ in $\C$, \cref{lem:0-complete-transport} gives $\mu_x=\sum_i c_i\rho_{e_i}\mu_{b_i}=0$, so $m_x=0$ by \ref{psnP} for $\D$. Since $x$ was arbitrary, $m=0$. Thus $\psn^{\eHom}$ satisfies \ref{psnP} for $\eHom(\C{\to}\D)$.

\itemstep{\ref{FDF}.}
For an $n$-ball $X$ and boundary condition $c\in\undfld[n-1]{\eHom(\C{\to}\D)}{\partial X}$,
the map
\[
    \mathrm{ev}\colon\undfldl[n]{\eHom(\C{\to}\D)}{X}[c]\to\bigoplus_{b\in\pi_0\C}\undfldl[n]{\D}{X}[c_b],\qquad
    m\mapsto(m_b)_{b\in\pi_0\C}
\]
is injective. Indeed, if $m_b=0$ for all $b\in\pi_0\C$, then $(\fcj{m}\blt m)_b=\fcj{m_b}\blt m_b=0$, so $\bkt{m}{m}_{X,c}=\psn^{\eHom}(\fcj{m}\blt m)=\sum_{b\in\pi_0\C}w_\C(b)\bkt{m_b}{m_b}_{X,c_b}=0$
and $m=0$ by \ref{psnP} for $\eHom(\C{\to}\D)$. And its target is finite-dimensional:
$|\pi_0\C|<\infty$ by \cref{thm:finiteness-components} and $\dim\undfldl[n]{\D}{X}[c_b]<\infty$
by \ref{FDF} for $\D$.
\end{pf}

\begin{remark}
Observe that $\psn^{\eHom}$ is the sphere trace transferred to $\eHom(\C{\to}\D)$ from the skeletonization of $\Hom$ by \cref{tsf} and \cref{tsf2} for $n=1$ and $n=2$ respectively. The purpose of \cref{tsf} and \cref{tsf2} therefore is to show that for $n=1,2$, disk-like weak equivalences induce isomorphisms on the vector spaces of $n$-fields for general $n$-manifolds (not just balls). The same holds for general $n$ by carrying out the obvious inductive construction/procedure evident from the proof of \cref{tsf2}. 

In addition, these two results show the corresponding hom disk-like 1- and 2-categories are finite. Nonetheless, we expect that if $\C$ is moreover complete, then this should follow by walking through the construction of the path integral and arguing similarly to the proof of \ref{FDF} above at each step (index) of the inductive handle construction in \cref{cstr:path-integral-via-handle-decomp}. 
\end{remark}

\section{\texorpdfstring{$(1+1)$D}{(1+1)D}}
\label{sec:1+1D}

\subsection{Disk-like 1-categories and functors}

\subsubsection{The definition of a disk-like 1-category}
\label{sec:def-dl1cat}
In this section, all 0-manifolds $P$ (resp. 1-manifolds $X$) are assumed to be PL, compact, and equipped with the germ of a thickening of $P$ (resp. $X$) to an oriented 2-manifold $M$. Gluing manifolds along embedded submanifolds of their boundaries means identifying these germs. We attempt to give more exposition in this section than in the next, opting to spell out the specializations of the general definitions from \cref{sec:n+1D}.

Here we give a relatively self-contained definition of a disk-like 1-category. For the full details, see \cite[\S6.1]{MW12} or \cref{sec:full-axiomatic-definition}. As in the rest of this article, we include the notion of reflection structure in the definition of disk-like $n$-category as extra data.

\itemstep{0-fields.}
To each oriented point $P$, which we call a \defn{0-ball}, we assign a set $\fld[0]{\C}{P}$ of ($P$-shaped, pure) \defn{0-fields}. To each (orientation-preserving) homeomorphism $\varphi\colon P\to P'$, we assign a bijection $\varphi_*\colon\fld[0]{\C}{P}\to\fld[0]{\C}{P'}$ such that $\C_0$ is a functor $\Disk_0\to\Set$. For each 0-ball $P$, we assign an involutive bijection $\fcj{\,\cdot\,}\colon\fld[0]{\C}{P}\to\fld[0]{\C}{\orev P}$ satisfying $\orev{\varphi}_*\fcj a=\fcj{\varphi_*a}$ for homeomorphisms $\varphi\colon P\to P'$. For any 0-manifold (a finite set of 0-balls) $S$, define $\undfld[0]{\undC}{S}\coloneq\prod_{Q\in S}\fld[0]{\C}{Q}$.

\itemstep{1-fields.}
For each oriented interval $X$, which we call a \defn{1-ball}, we assign the following data. 
\begin{lst}
\item[(D$\C_1$)]\label{DC1} (\emph{Pure 1-fields}) A set $\fld[1]{\C}{X}$ of \defn{pure 1-fields}.
\item[(D$\varphi_1$)]\label{Dvarphi1} (\emph{Homeomorphism actions}) A bijection $\varphi_*\colon\fld[1]{\C}{X}\to\fld[1]{\C}{X'}$ for each homeomorphism $\varphi\colon X\to X'$.
\item[(D$\partial_1$)]\label{Dbdy1} (\emph{Boundary maps}) A map $\partial_1\colon\fld[1]{\C}{X}\to\undfld[0]{\undC}{\partial X}$ sending a 1-field on $X$ to its boundary pair of 0-fields. Write $\fld[1]{\C}{X}[c]$ for the preimage $\partial_1^{-1}(c)\subset\fld[1]{\C}{X}$.
\item[(DG$_1$)]\label{DG1} (\emph{Gluing}) For each splitting $X=X_1\cup_P X_2$ along a point $P$, a gluing map $\glu_P\colon\fld[1]{\C}{X_1}\times_{\fld[0]{\C}{P}}\fld[1]{\C}{X_2}\to\fld[1]{\C}{X}$.
\item[(D$\pi_1$)]\label{Dpi1} (\emph{Pullbacks}) For each 1-ball $X$ and 0-ball $P$, a set map $\pi^*\colon\fld[0]{\C}{P}\to\fld[1]{\C}{X}$ (the pullback of the unique pinched product map $\pi\colon X\to P$). The (pure) 1-fields in its image are called \defn{product fields}.
\item[(D$\fcj{\,\cdot\,}_1$)]\label{Dfcj1} (\emph{Reflection}) An involutive bijection $\fcj{\,\cdot\,}\colon\fld[1]{\C}{X}\to\fld[1]{\C}{\orev X}$.
\item[(D$U_1$)]\label{DU1} (\emph{Local relations}) For each 1-ball $X$ and $c\in\undfld[0]{\undC}{\partial X}$, a subspace $\fldlU{\C}{X}[c]\subset\bC\{\fld[1]{\C}{X}[c]\}$ of so-called \defn{local relations}.
\end{lst}
The above data are subject to the following conditions.
\begin{lst}
\item[(C$\varphi$)]\label{Cvarphi1} $\C_1$ is a functor $\Disk_1\to\Set$.
\item[(C$\partial$)]\label{Cbdy1} $\partial_1$ is a natural transformation $\C_1\Rightarrow\undC_0\circ\partial$.
\item[(CGi)]\label{CGi1} Every gluing map $\glu_P$ is injective.
\item[(CGa)]\label{CGa1} Gluing is associative in that any two ways of starting with a collection of 1-balls and gluing them up sequentially to form another 1-ball induce equal gluing maps.
\item[(CG)]\label{CGvarphi1} A homeomorphism respecting a splitting sends the glued field to the gluing of its image. The boundary of a glued field restricts correctly to each piece's outer boundary.
\item[(C$\pi$)]\label{Cpi1} Product fields are compatible with gluing, homeomorphisms, boundary maps, and restriction to sub-1-balls.
\item[(C$\fcj{\,\cdot\,}$)]\label{Cfcj1} Reflection commutes with homeomorphisms, boundary maps, gluing, and products.
\item[(CS)]\label{CS1} \emph{Splittability.} Every pure 1-field $\xi$ on a 1-ball $X$ can be cut at all but finitely many (interior) points and still recovered by gluing the two halves back together. Call any such finite set a \defn{string locus} of $\xi$.
\item[(C$U$)]\label{CU1} Local relations are preserved under homeomorphisms and reflection, and form an ideal in the sense that gluing a local relation to anything results in a local relation.
\item[(C$U$i)]\label{CTi1} \emph{Isotopy invariance.} If all components of an isotopy $h$ on a 1-ball fix the boundary of a 1-field $\xi$, then $\xi$ and its image under $h$ differ only by a local relation.
\item[(C$U$c)]\label{CTc1} \emph{Collaring invariance.} For a 1-ball $X$, attaching a collar $X\mapsto X\cup_P J$ at an endpoint $P\in\partial X$ and collapsing it back via a homeomorphism $\psi_{P,J}\colon X\cup_P J\to X$ that is the identity outside a neighborhood of $P$ changes a pure 1-field only by a local relation: $(\psi_{P,J})_*(\xi\blt_P\pi^*(\partial_1(\xi)|_P))-\xi\in\fldlU{\C}{X}[c]$ for any $\xi\in\fld[1]{\C}{X}[c]$.
\end{lst}

\itemstep{Extending to all 1-manifolds.}
Let $\Disk(X,c;\C)$ be the poset of pairs $(\cP; \beta_\cP)$, where $\cP$ is a splitting of $X$ into 1-balls and $\beta_\cP$ consists of compatible boundary labels for each constituent 1-ball in $\cP$. The ordering is by \defn{antirefinement}: $(\cP; \beta_\cP) \leq (\cQ; \beta_\cQ)$ whenever $(\cP;\beta_\cP)$ admits a sequence of gluings of its constituent 1-balls that results in $(\cQ;\beta_\cQ)$. Define a functor $\Pi_{X;c}^{\C} \colon \Disk(X,c;\C) \to \Set$ via $\Pi_{X;c}^{\C}(\cP; \beta_\cP) \coloneq \prod_{X_i \in \cP} \fld[1]{\C}{X_i}[ \beta_\cP|_{\partial X_i}]$. We define
\[
\undfld[1]{\undC}{X}[c] \coloneq \colim \Pi_{X;c}^{\C}
\]
and on morphisms by the application of the corresponding (iterated) gluing map. 

For a $1$-ball $X$, the trivial one-ball splitting of $X$ consisting of $X$ itself is the terminal object in $\Disk(X,c;\C)$, so $\undfld[1]{\undC}{X}[c] = \fld[1]{\C}{X}[c]$. Thus we will henceforth drop the arrow decoration on $\undC$ in favor of $\C$.

\itemstep{Linearization.}
Set $\undfldl[1]{\C}{X}[c]\coloneq\bC\{\undfld[1]{\C}{X}[c]\}\big/\fldlU{\C}{X}[c]$.
If $\partial X=\varnothing$, we write $\undfldl[1]{\C}{X}$ to mean the vector space $\fldl[1]{\C}{X}[\varnothing]$. This defines a $\Vec$-valued functor on the category of pairs $(X,c)$ of 1-manifolds $X$ and boundary conditions $c$ whose morphisms $(X,c)\to (X',c')$ are homeomorphisms $\varphi\colon X\to X'$ with $\varphi_*c=c'$; by the above axioms, gluing, boundary maps, and reflection descend to this quotient.

Recall that in this paper we will always require \ref{FDF}, i.e., that for a disk-like $n$-category $\C$, the local relations $\fldlU{\C}{X}[c]$ make the vector space of $n$-fields $\fldl[n]{\C}{X}[c]$ finite-dimensional for all $n$-balls $X$ and $c\in\undfld[n-1]{\C}{\partial X}$. 

Observe that the result is the same regardless of whether we linearize before or after taking the colimit. That is,
\[
\undfldl[1]{\C}{X}[c] = \bC\{\colim\Pi_{X;c}^\C\}\Big/{\fldlU{\C}{X}[c]}=\colim\tld\Pi_{X;c}^{\C}
\]
where $\tld\Pi_{X;c}^{\C}\colon\Disk(X,c;\C)\to\Vec$ is given by $\tld\Pi_{X;c}^{\C}(\cP;\beta_\cP) \coloneq \bigotimes_{X_i\in\cP}\fldl[1]{\C}{X_i}[\,\beta_\cP|_{\partial X_i}]$ and on morphisms by the maps induced by the linearized gluing maps on the tensor products of the vector spaces of 1-fields on the constituent balls of the decomposition. Here $\fldlU{\C}{X}[c] \subset \bC\{\undfld[1]{\C}{X}[c]\}$ is the collection of all finite sums of the form $\sum_\gamma \lambda_\gamma [\eta_1, \dots, \eta_{j-1}, \xi_j^\gamma, \eta_{j+1}, \dots, \eta_N]$ indexed over ball splittings $(\cP = \{X_i\}_{i=1}^N;\beta_\cP)\in\Disk(X,c;\C)$, $1 \leq j \leq N$, local relations $\sum_\gamma \lambda_\gamma \xi_j^\gamma \in \fldlU{\C}{X_j}[ \beta|_{\partial X_j}]$, and compatible pure 1-fields $\eta_i \in \fld[1]{\C}{X_i}[ \beta|_{\partial X_i}]$ for $i \neq j$. The square bracket $[-]$ denotes the image in the (set-valued) colimit $\undfld[1]{\C}{X}[c]$ (not to be confused with the square bracket $[-]$ used to denote the (quotient-level) 1-field $[\alpha]$ represented by a pure 1-field $\alpha$).

We will call a disk-like 1-category $\C$ \defn{finite} if
\begin{lst}
    \item[(F)] 
    \label{1psnF}
    (\emph{Finiteness})
    $\dim\undfldl[1]{\C}{X}[c]<\infty$ for all 1-manifolds $X$ and $c\in\undfld[0]{\C}{\partial X}$.
\end{lst}

Observe that \ref{1psnF} regards general 1-manifolds, not just 1-balls.

\subsubsection{Unitary disk-like 1-categories}
\label{sec:unitary disk-like 1-categories}

\

\begin{definition}[Unitary disk-like 1-category]
\label{def:unitary disk-like 1-category}
A \defn{unitary disk-like 1-category} $(\C,\psn)$ is a disk-like 1-category $\C$ equipped with a \defn{sphere trace} (or \defn{trace}), that is, a linear functional $\psn\colon\undfldl[1]{\C}{S^1}\to\bC$ satisfying the following condition.
\begin{lst}[font=\upshape]
\item[($\psn$P)]\label{1psnP}
(\emph{Positivity})
For each 1-ball $X$ and each $c\in \undfld[0]{\C}{\partial X}$, the sesquilinear pairing
\begin{equation}
\begin{aligned}
\orev{\fldl[1]{\C}{X}[c]}\otimes_{\bC}\fldl[1]{\C}{X}[c] & \longrightarrow \bC,
\\
f\otimes g                      & \longmapsto \bkt{f}{g}_{X,c}\coloneq\psn(\fcj f\blt_{\partial X} g)
\end{aligned}
\end{equation}
is positive-definite.
\end{lst}
\end{definition}

\subsection{Dagger categories, unitary categories, and 2-Hilbert spaces}
Recall that a \defn{dagger category} is a linear category $\cX$ equipped with a \defn{dagger structure}, that is, a collection of conjugate-linear maps $\dag\colon\cX(x\to y)\to \cX(y\to x)$ for all $x,y\in\cX$ such that for all $f\in\cX(x\to y)$ and $g\in \cX(y\to z)$, $(g\circ f)^\dag=f^\dag\circ g^\dag$ and $f^{\dag\dag}=f$. 

Let $\cX$ be a linear category. For objects $a_1,\dots,a_\ell$, the \defn{linking algebra} $L(a_1,\dots,a_\ell)\coloneq\bigoplus_{i,j=1}^\ell\cX(a_j \to a_i)$ has elements formal matrices, with multiplication given by matrix multiplication. When $\cX$ has a dagger structure, $L(a_1,\dots,a_\ell)$ is a $*$-algebra with involution the dagger-transpose.

\begin{definition}
\label{def:unitarycategory}
A \defn{unitary category} is a dagger category whose linking algebras are \defn{unitary}, i.e., finite-dimensional $\Cstar$-algebras.
\end{definition}

\begin{definition}
\label{def:finiteunitarycategory}
A unitary category $\cX$ is called \defn{finite} if there is a global bound on the dimensions of the centers of all linking algebras, or equivalently if its completion $\cX^\cent$ has finitely many isomorphism classes of simple objects.
\end{definition}

A \defn{complete set of orthogonal minimal projections} in a finite-dimensional $\Cstar$-algebra $A$ is a finite set $\{p_i\}\subset A$ of projections (self-adjoint idempotents) such that $p_i A p_i\cong\bC$ (minimality) and $\sum_i p_i=1$ (completeness); completeness implies $p_i p_j=\delta_{ij}p_i$ (orthogonality). 

A linear category $\cX$ is called \defn{presemisimple} if all linking algebras are finite-dimensional semisimple algebras, and \defn{finite presemisimple} if moreover there is a global bound on the dimensions of their centers, or equivalently if its completion has finitely many simple objects (see \cref{completeunitarycategories} below). Since finite-dimensional $\Cstar$-algebras are semisimple, every unitary category is presemisimple. In particular, every linking algebra of a unitary category admits a complete set of orthogonal minimal projections.

\subsubsection{Unitary morphisms and unitary equivalences of fields}
\label{sec:coisometries-1}
Let $\cX$ be a dagger category. An \defn{isometry} is a morphism $f\in\cX(a\to b)$ with $f^\dag f=\id_a$. A \defn{coisometry} is a morphism $f\in\cX(a\to b)$ with $ff^\dag=\id_b$. A \defn{unitary isomorphism} (or a \defn{unitary}) is a morphism that is both an isometry and a coisometry.

The following are the analogous notions in the disk-like framework, which is just the unpacking of \cref{def:isometry-1field} for the case $n=1$.

\begin{definition} \label{def:cosiom1}
Consider a (pure or quotient-level) 1-field $\xi\in\fld[1]{\C}{D^1}[\fcj{a}\amalg b]$, and set $\xi^\star\coloneq\fcj{\iota_*\xi}$ where $\iota\colon D^1\to\orev{D^1}$ is $x\mapsto -x$.
\begin{itemize}
    \item We call $\xi$ an \defn{isometry} if $\beta_\xi\coloneq\xi\times S^0\cup b\times D^1$ equals $a\times D^1$ as (quotient-level) 1-fields: 
    \[
        \begin{tkz}
            \draw[mid<,violet,thick](0,0)--node[left]{\scriptsize$\xi^\star$}(0,.6);
            \draw[mid<,violet,thick](1.2,.6)--node[right]{\scriptsize$\xi$}(1.2,0);
            \draw[mid<,blue,thick] (0,.6) arc(180:0:.6)node[pos=.5,above]{\scriptsize$b\times D^1$};
        \end{tkz}
        =
        \quad
        \begin{tkz}
            \draw[red,thick](0,0)--(0,.6);
            \draw[red,thick](1.2,.6)--(1.2,0);
            \draw[mid<,red,thick] (0,.6) arc(180:0:.6)node[pos=.5,above]{\scriptsize$a\times D^1$};
        \end{tkz}
    \]
    \item We call $\xi$ a \defn{coisometry} if $\beta_{\xi^\star}\coloneq\xi^\star{\times} S^0\cup a{\times} D^1$ equals $b\times D^1$ as (quotient-level) 1-fields: 
    \[
        \begin{tkz}[yscale=-1]
            \draw[mid>,violet,thick](0,0)--node[left]{\scriptsize$\xi^\star$}(0,.6);
            \draw[mid>,violet,thick](1.2,.6)--node[right]{\scriptsize$\xi$}(1.2,0);
            \draw[mid>,red,thick] (0,.6) arc(180:0:.6)node[pos=.5,below]{\scriptsize$a\times D^1$};
        \end{tkz}
        =
        \quad
        \begin{tkz}[yscale=-1]
            \draw[blue,thick](0,0)--(0,.6);
            \draw[blue,thick](1.2,.6)--(1.2,0);
            \draw[mid>,blue,thick] (0,.6) arc(180:0:.6)node[pos=.5,below]{\scriptsize$b\times D^1$};
        \end{tkz}
    \]
\end{itemize} 
We call $\xi$ \defn{unitary} or a \defn{unitary equivalence} if $\xi$ is both an isometry and a coisometry. We write $a\cong^\star b$ if there is a unitary equivalence $\xi$ between $a$ and $b$.

We call a general (pure or quotient-level) 1-field $\xi\in\fld[1]{\C}{X}$ an isometry (resp. coisometry, unitary) if there is a homeomorphism $\varphi\colon X\to D^1$ such that $\varphi_*\xi$ is an isometry (resp. coisometry, unitary).
\end{definition}

\subsubsection{Complete unitary categories}
\label{completeunitarycategories}
A unitary category is called \defn{complete} if it has all finite orthogonal direct sums and splittings of projections. The \defn{completion} $\cX^\cent$ of a unitary category $\cX$ is the projection completion of its orthogonal additive envelope. That is, objects of $\cX^\cent$ are pairs $(\bigoplus_ix_i,p)$ consisting of a formal orthogonal direct sum and an orthogonal projection $p$ on $\bigoplus_ix_i$, and morphisms $(\bigoplus_i x_i, p) \to (\bigoplus_j y_j, q)$ are matrices $(f_{ji})$ with entries $f_{ji} \in \cX(x_i\to y_j)$ satisfying $q \circ (f_{ji}) \circ p = (f_{ji})$ where composition is given by matrix multiplication. The linking algebras of a unitary category $\cX$ are the endomorphism algebras of objects in its completion $\cX^\cent$.

\begin{proposition}[Universal property of completion {\cite[Ch. 8]{UQSL}}]
  \label{3.3.12}
  If $\cX$ and $\cY$ are unitary categories and $\cY$ is complete, then restriction along the canonical inclusion $\iota\colon\cX\hookrightarrow\cX^{\cent}$ gives an equivalence $\iota^*\colon\Hom(\cX^{\cent}\to\cY)\overset\sim\to\Hom(\cX\to\cY)$.
\end{proposition}

\subsubsection{Traces and pre-2-Hilbert spaces}

\

\begin{definition}[{{\cite[Defn. 2.9]{3Hilb}, notation adapted}}]
\label{2hilbtrace}
  A \defn{faithful positive trace} on a unitary category $\cX$ is a collection of linear maps $\Tr_a^\cX\colon \End_\cX(a)\to \bC$ for each $a\in\cX$ satisfying the following properties.
  \begin{lst}
    \item[\upshape($\dag$Tr1)]\label{dagTr2}The sesquilinear form $\bkt{f}{g}_{a\to b}\coloneq \Tr^\cX_a(f^\dag\circ  g)$ on $\cX(a\to b)$ is positive-definite.
    \item[\upshape($\dag$Tr2)]\label{dagTr1} $\Tr^\cX_a(g\circ f)=\Tr^\cX_b(f\circ g)$ for all $f\in\cX(a\to b)$ and $g\in \cX(b\to a)$.
  \end{lst}
\end{definition}

\begin{definition}
  \label{def:pre-2-Hilbert space}
  A \defn{pre-2-Hilbert space} is a unitary category equipped with a faithful positive trace. A \defn{2-Hilbert space} is a complete (equivalently, semisimple) pre-2-Hilbert space. A (pre-)2-Hilbert space is called \defn{finite} if it is finite as a unitary category.
\end{definition}

\begin{remark}
\label{rem:2Hilb-nonstandard}
As warned in the introduction, allowing non-finite 2-Hilbert spaces is non-standard: in \cite{Baez97,3Hilb,bases} a 2-Hilbert space is finite by definition. Our non-finite 2-Hilbert spaces are only the semisimple ones, and \cref{def:pre-2-Hilbert space} is not intended as a general notion of infinite-dimensional 2-Hilbert space. For instance, the category of Hilbert space modules over a von Neumann algebra should also be considered a 2-Hilbert space, but it is not semisimple and so is not of this form; see \cite{HPT24}.
\end{remark}

The opposite category $\cX^\op$ of a (finite) (pre-)2-Hilbert space $\cX$ is again one when equipped with the same dagger and trace. Moreover, $(\cX^\op)^\cent=(\cX^\cent)^\op$.

\begin{remark}
\label{trcomp}
The completion $\cX^\cent$ of a pre-2-Hilbert space $\cX$ is a 2-Hilbert space whose faithful positive trace $\Tr^{\cX^\cent}$ is the extension to $\cX^\cent$ of $\Tr^\cX$ given by $\Tr^{\cX^\cent}_{(\bigoplus a_{i},p)}(f)\coloneq\sum_{i} \Tr^\cX_{a_{i}}(p_{i}fp_{i})$ for all $f = pfp \in \End(\bigoplus a_{i})$, where $\bigoplus a_{i}$ is an orthogonal direct sum and $p \in \End(\bigoplus a_{i})$ is a projection.
\end{remark}

\subsection{From disk-like to traditional: the pre-2-Hilbert space \texorpdfstring{$\cX_\eC$}{X(C)}}

\subsubsection{\texorpdfstring{$\cX_\eC$}{X(C)} as a unitary category}
Except for the definition of the dagger and the notion of reflection structure, the following construction can mostly be found in \cite[Appendix C.1]{MW12}, but we include the details here for completeness.

\begin{construction}
  \label{skein 1-categories}
  \label{daggers for skeins}
  Let $\C$ be a disk-like 1-category. Define $\Obj(\cX_\C)\coloneq\fld[0]{\C}{\pt}$, and for each $a,b\in\cX_\C$ define $\cX_\C(a \to b) \coloneq \fldl[1]{\C}{D^1}[ \fcj{a} \amalg b]$, which is finite-dimensional by \ref{FDF}.
  Now fix a homeomorphism $\rho\colon I\cup_\pt I\to I$. (Since the choice of $\rho$ is unique up to isotopy-rel-boundary, $\rho$ is not extra data by \ref{CTi1}.) For $f \in \cX_\C(a \to b)$ and $g \in \cX_\C(b \to c)$, define $g \circ f \coloneq \rho_*(f \blt g)\in\cX_\C(a\to c)$. Composition is associative: the two reparameterization maps $I\cup_\pt I\cup_\pt I\to I$ built from $\rho$, namely $\rho^{(12)3}\coloneq \rho\circ(\rho\cup_\pt\id_I)$ and $\rho^{1(23)}\coloneq \rho\circ(\id_I\cup_\pt\rho)$, are isotopic-rel-boundary; thus, for composable morphisms $f,g,h$ in $\cX_\C$, we have
  \begin{align*}
    h\circ (g\circ f)
     & =
    \rho_*(\rho_*(f\blt g)\blt h)
    =
    \rho_*((\rho\cup_\pt\id_I)_*(f\blt g\blt h))
    =
    \rho^{(12)3}_*(f\blt g\blt h)
    \\&
    =
    \rho^{1(23)}_*(f\blt g\blt h)
    =
    \rho_*((\id_I\cup_\pt\rho)_*(f\blt g\blt h))
    =
    \rho_*(f\blt \rho_*(g\blt h))
    =
    (h\circ g)\circ f,
  \end{align*}
  where the second and sixth equalities are by \ref{CGvarphi1}. Finally, for $a\in\cX_\C$, the identity morphism is the product field $\id_a\coloneq a \times I$. Indeed, for any $f\in\cX_\C(a\to b)$, we have $f \circ \id_a=\rho_*((a\times I) \blt f)=f$ where the second equality follows from \ref{CTc1}. Showing $\id_b \circ f = f$ is similar.

  Lastly, we equip $\cX_\C$ with the following canonical dagger structure, which is the dagger structure implicit in \cite[Ch. 4]{W06}. Since the homeomorphism $\iota\colon D^1\to \orev{D^1}$ given by $x\mapsto -x$ is orientation-\emph{preserving}, postcomposition with orientation-reversal $\orev{D^1}\to D^1$ induces a conjugate-linear map
  \[
    \begin{aligned}
      \star\colon\cX_\C(a\to b) & \to \cX_\C(b\to a),                \\
      f                        & \mapsto f^\dag\coloneq \fcj{\iota_*f}
    \end{aligned}
  \]
  for all $a,b\in\cX_\C$.
\end{construction}

\begin{lemma}
  \cref{daggers for skeins} indeed gives a dagger structure on $\cX_\C$.
\end{lemma}

\begin{pf}
  Since $\iota$ and $\fcj{\,\cdot\,}$ are commuting involutions, $f^{\dag\dag}=f$ for all morphisms $f$ in $\cX_\C$ and $(g\circ f)^{\dag}=\fcj{\iota_*(g\circ f)}=\fcj{\iota_*f}\circ \fcj{\iota_*g}=f^{\dag}\circ g^{\dag}$ for composable morphisms $f$ and $g$ in $\cX_\C$.
\end{pf}

\begin{example}
\label{skein 1-category}
Let $\C$ be a disk-like $n$-category and fix an $(n-1)$-manifold $Y$ and a boundary condition $\eta\in\undfld[n-2]{\C}{\partial Y}$. Then there is an (ordinary/traditional) dagger 1-category $\Sk_\C(Y,\eta)$, called the \defn{skein 1-category} of $\C$ on $(Y,\eta)$, defined by the skeletonization of the disk-like skein 1-category $\A_\C(Y,\eta)$ from \cref{disk-like skein k-categories}:
\[
\Sk_\C(Y,\eta)\coloneq\cX_{\A_\C(Y,\eta)}.
\]
Thus
\begin{itemize}
    \item objects are $(n-1)$-fields $\xi\in\undfld[n-1]{\C}{Y}[\eta]$, and
    \item morphisms $f\in\Sk_\C(Y,\eta)(\xi\to\zeta)$ are (quotient-level) $n$-fields $f\in\fldl{\C}{Y\times D^1}[\fcj{\xi}\cup\zeta]$, where the product $Y\times D^1$ is pinched so that $\partial (Y\times D^1)=\orev{Y}\cup_{\partial Y}Y$.
\end{itemize}
\end{example}

\begin{lemma}[Unitarity of $\cX_\eC$]
  \label{unitarity of cX}
  If $\eC$ is a unitary disk-like 1-category, then $\cX_\eC$ is a unitary category.
\end{lemma}

\begin{pf}
  Fix objects $a_1,\dots,a_\ell\in \cX_\eC$ and consider the linking algebra $L\coloneq L(a_1,\dots,a_\ell)$. It is finite-dimensional since its summands are finite-dimensional by \ref{FDF}. Define a pairing $\bkt{-}{-}_{L}\colon L\times L\to\bC$ by $\bkt{(f_{ij})}{(g_{ij})}_{L}\coloneq \sum_{i,j=1}^\ell\bkt{f_{ij}}{g_{ij}}_{a_j\to a_i}$. This inherits positive-definiteness from the pairings $\bkt{-}{-}_{a_j\to a_i}$ from \ref{1psnP},
  and $L$ is $*$-definite in that $x^*x=0$ implies $x=0$ (with $*$ the dagger-transpose),
  so by the Characterization of Unitary Algebras \cite[Thm. I.2.4.4]{UQSL}, $L$ is unitary.
\end{pf}

\begin{proposition}[Finiteness of skein 1-categories]
\label{lem:closing-up}
Let $\C$ be a finite unitary disk-like $n$-category, let $Y$ be an $(n-1)$-manifold, and let $\eta\in\undfld[n-2]{\C}{\partial Y}$.
Then $\Sk_\C(Y,\eta)$ is finite unitary and
\begin{equation}
\label{dimct}
|\Irr(\Sk_\C(Y,\eta)^\cent)|=\dim\undfldl[n]{\C}{Y\times S^1}[\eta\times S^1].
\end{equation}
\end{proposition}

\begin{pf}
Since $\A_\C(Y,\eta)$ is a unitary disk-like 1-category by \cref{inducspwt}, $\Sk_\C(Y,\eta)=\cX_{\A_\C(Y,\eta)}$ is unitary by \cref{unitarity of cX}, say with $\Irr(\Sk_\C(Y,\eta)^\cent)=\{(r_i,p_i)\}_{i\in I}$. 
The subspaces $p_i\End_{\Sk_\C(Y,\eta)}(r_i)p_i$ are 1-dimensional by minimality of the $p_i$, and $\bkt--_{Y\times D^1,(\eta\times D^1)\cup(\fcj{r_i}\amalg r_i)}$ is positive-definite on them by \ref{Pos}. As $\C$ is finite unitary, the path integral $\Z_\C$ from \cref{UnitaryWalkerTheorem} gives finite unitary path integral data, so \cref{glu-basis} applies and \cref{glu-basis}\ref{dimformula} gives \eqref{dimct}, whose right side is finite by \ref{psnF}. Thus $\Sk_\C(Y,\eta)$ is finite.
\end{pf}

\begin{corollary}
  \label{finite cylinders}
  If $\C$ is a finite unitary disk-like 1-category, then $\cX_\eC$ is finite unitary and
  $|\Irr(\cX_\eC^\cent)|=\dim\undfldl[1]{\C}{S^1}$.
\end{corollary}
\begin{pf}
As $\Sk_\C(\pt)=\cX_\eC$, this follows from \cref{lem:closing-up} with $Y\coloneq\pt$ and $\eta\coloneq\varnothing$.
\end{pf}

\subsubsection{The faithful positive trace on \texorpdfstring{$\cX_\eC$}{X(C)}}

\

\begin{construction}
  \label{cstr:pre-2-Hilbert spaces from unitary disk-like 1-categories}
  For a unitary disk-like 1-category $(\eC,\psn^\eC)$, we equip the dagger 1-category $\cX_\eC$ constructed in \cref{skein 1-categories} with the faithful positive trace $\Tr^{\cX_\eC}$ defined for $a\in\cX_\C$ by the functional $\Tr_a^{\cX_\eC}\colon\End_{\cX_\eC}(a)\to\bC$ given by 
  \[
    \Tr_a^{\cX_\eC}(f)\coloneq \psn^\eC(\cl_af)=\bkt{\id_a}{f}_{a\to a}
    =
    \psn^\eC\Big(
    \begin{tkz}[scale=0.5]
      \draw[thick,Green,mid>] (0,0) arc(90:270:.65) node[pos=.5,left=.3ex,font=\scriptsize]{$\fcj{a{\times}I}$};
      \draw[thick,Orange,mid>] (0,-1.3) arc(-90:90:.65) node[pos=.5,right=.3ex,font=\scriptsize]{$f$};
    \end{tkz}
\Big).
  \]
\end{construction}

\begin{lemma}
  \label{pairing-composition compatibility}
  Let $\eC$ be a unitary disk-like 1-category, and let $f\in\cX_\eC(a\to b)$, $g\in\cX_\eC(b\to c)$, and $h\in\cX_\eC(a\to c)$. Then
  
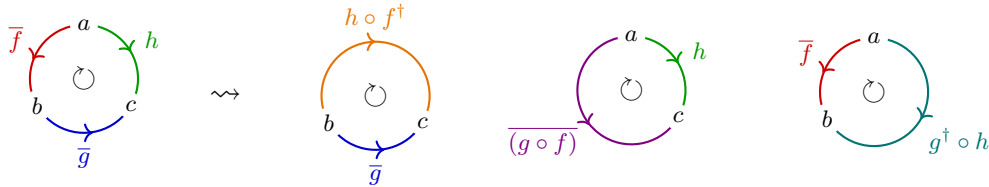
\[
    \begin{gathered}
    \bkt{g}{h\circ f^\dag}_{b\to c}=\bkt{g\circ f}{h}_{a\to c}=\bkt{f}{g^\dag \circ h}_{a\to b}.
    \\
    \def\rad{1.1}%
    \def\lblrad{1.6}%
    \begin{tkz}[scale=0.65,thick]
      \draw[mid>, green!60!black] (90:\rad) arc (90:-30:\rad);
      \node[green!60!black, font=\scriptsize] at (30:\lblrad) {$h$};
      \draw[mid<, blue!80!black] (-30:\rad) arc (-30:-150:\rad);
      \node[blue!80!black, font=\scriptsize] at (-90:\lblrad) {$\fcj{g}$};
      \draw[mid<, red!80!black] (-150:\rad) arc (-150:-270:\rad);
      \node[red!80!black, font=\scriptsize] at (150:\lblrad) {$\fcj{f}$};
      \node[circle, fill=white, inner sep=1.5pt, font=\scriptsize] at (90:\rad)  {$a$};
      \node[circle, fill=white, inner sep=1.5pt, font=\scriptsize] at (-30:\rad) {$c$};
      \node[circle, fill=white, inner sep=1.5pt, font=\scriptsize] at (210:\rad) {$b$};
      \node at (0,0) {$\circlearrowright$};
    \end{tkz}%
    \quad
    \rightsquigarrow
    \qquad
    \begin{tkz}[scale=0.65,thick]
      \draw[mid>, orange!90!black] (-150:\rad) arc (-150:-270:\rad) arc (-270:-390:\rad);
      \node[orange!90!black, font=\scriptsize] at (90:\lblrad) {$h \circ f^\dag$};
      \draw[mid<, blue!80!black] (-30:\rad) arc (-30:-150:\rad);
      \node[blue!80!black, font=\scriptsize] at (-90:\lblrad) {$\fcj{g}$};
      \node[circle, fill=white, inner sep=1.5pt, font=\scriptsize] at (-30:\rad) {$c$};
      \node[circle, fill=white, inner sep=1.5pt, font=\scriptsize] at (210:\rad) {$b$};
      \node at (0,0) {$\circlearrowright$};
    \end{tkz}
    \qquad
    \begin{tkz}[scale=0.65,thick]
      \draw[mid<, Purple] (-30:\rad) arc (-30:-150:\rad) arc (-150:-270:\rad);
      \node[Purple, font=\scriptsize] at (210:1.3*\lblrad) {$\fcj{(g \circ f)}$};
      \draw[mid>, green!60!black] (90:\rad) arc (90:-30:\rad);
      \node[green!60!black, font=\scriptsize] at (30:\lblrad) {$h$};
      \node[circle, fill=white, inner sep=1.5pt, font=\scriptsize] at (90:\rad) {$a$};
      \node[circle, fill=white, inner sep=1.5pt, font=\scriptsize] at (-30:\rad) {$c$};
      \node at (0,0) {$\circlearrowright$};
    \end{tkz}
    \;
    \qquad
    \begin{tkz}[scale=0.65,thick]
      \draw[mid>, teal!90!black] (90:\rad) arc (90:-30:\rad) arc (-30:-150:\rad);
      \node[teal!90!black, font=\scriptsize] at (330:1.25*\lblrad) {$g^\dag \circ h$};
      \draw[mid<, red!80!black] (-150:\rad) arc (-150:-270:\rad);
      \node[red!80!black, font=\scriptsize] at (150:\lblrad) {$\fcj{f}$};
      \node[circle, fill=white, inner sep=1.5pt, font=\scriptsize] at (90:\rad) {$a$};
      \node[circle, fill=white, inner sep=1.5pt, font=\scriptsize] at (210:\rad) {$b$};
      \node at (0,0) {$\circlearrowright$};
    \end{tkz}
    \end{gathered}
  \]
\end{lemma}

\begin{pf}
  Recall that $\undfldl[1]{\C}{S^1}$ is constructed as the colimit over splittings of $S^1$ into 1-balls. Observe that $\fcj{g}\blt_c(h\circ f^\dag)$, $\fcj{(g\circ f)}\blt_c h$, and $\fcj{f}\blt_b(g^\dag\circ h)$ represent the same element (colimit class) of $\undfldl[1]{\C}{S^1}$: they are all coarsenings of the splitting of $S^1$ consisting of three 1-balls whose labels in cyclic order are $h$, $\fcj{g}$, and $\fcj{f}$. By applying the sphere trace $\psn^\C$ (and recalling \cref{lem:reparam-gen}) we obtain the claimed equalities.
\end{pf}

\begin{lemma}
If $\C$ is (finite) unitary, then $\cX_\C$ equipped with the trace from \cref{cstr:pre-2-Hilbert spaces from unitary disk-like 1-categories} is a (finite) pre-2-Hilbert space. 
\end{lemma}

\begin{pf}
 First $\cX_\eC$ is unitary by \cref{unitarity of cX}, and is finite unitary when $\eC$ is finite by \cref{finite cylinders}. To see we get a pre-2-Hilbert space $(\cX_{\eC},\Tr^{\cX_\eC})$, observe that \ref{dagTr2} is just \ref{1psnP} and \ref{dagTr1} follows from \cref{pairing-composition compatibility}:
  \[
    \Tr_b^{\cX_\eC}(f\circ g)=\bkt{\id_b}{f\circ g}_{b\to b}=\bkt{g^\dag}{f}_{a\to b}=\bkt{\id_a}{g\circ f}_{a\to a}=\Tr_a^{\cX_\eC}(g\circ f).\qedhere
  \]
\end{pf}

\subsection{From traditional to disk-like: the unitary disk-like 1-category \texorpdfstring{$\eC^\cX$}{C(X)}}

\subsubsection{String diagrams and the disk-like 1-category \texorpdfstring{$\eC^\cX$}{C(X)}}
Let $\cX$ be a (traditional) 1-category. An \defn{$\cX$-string diagram} on a 1-ball $X \in \Disk_1$ is a pair $\xi=(\Gamma,\lambda)$, where
\begin{lst} 
\item $\Gamma$ is a \defn{1D string diagram stratification}, that is, a tuple $\Gamma=((E,V),\e)$ consisting of a set $V=\{v_1<\cdots<v_m\}$ of finitely many points in the interior of $X$ (called \defn{vertices} or \defn{0-cells}) ordered by the underlying orientation of $X$, together with a transverse orientation $\varepsilon_v \in \{\pm 1\}$ at each vertex $v\in V$; and
\item $\lambda$ consists of object labels $\lambda_{e}\in\cX$ on the resulting \defn{edges} (also called \defn{1-cells} $e\in E$), i.e., the connected components of $X\setminus V$, and morphism labels $\lambda_{v}$ at the vertices $v\in V$, where $\lambda_{v_i} \in\cX(\lambda_{e_{i-1}} \to \lambda_{e_i})$ when $\varepsilon_{v_i} = +1$ and $\lambda_{v_i} \in\cX(\lambda_{e_i} \to \lambda_{e_{i-1}})$ when $\varepsilon_{v_i} = -1$.
\end{lst}

\begin{construction}
  \label{cstr:disk-like 1-categories from ordinary 1-categories}
  A dagger 1-category $\cX$ gives a disk-like $1$-category $\eC^\cX$ as follows.
  \begin{lst}[leftmargin=.5in]
    \item[\ref{DC1}] A field on a 0-ball $P$ is a labeling of $P$ by an object of $\cX$. A pure $1$-field on a 1-ball $X$ is an $\cX$-labeled string diagram $\xi=(\Gamma,\lambda)$ on $X$.
    \item[\ref{Dvarphi1}] For a string diagram $\xi=(\Gamma,\lambda)$ on $X$ and a homeomorphism $\varphi\colon X\to X'$, we define $\varphi_*\xi\coloneq (\varphi(\Gamma),\lambda\circ \varphi^{-1})$.
      \item[\ref{Dbdy1}] For a string diagram $\xi\in\fld[1]{\eC^\cX}{X}$, define $\partial \xi\coloneq c_\inn\amalg c_\out\in\undfld[0]{\C^\cX}{P\amalg Q}$ where $\partial X=P\amalg Q$ and where $c_\inn\in\cX$ (resp. $c_\out\in\cX$) labels the incoming (resp. outgoing) half-1-cell in $\xi$.
      \item[\ref{DG1}] For compatible fields $\xi_1=(\Gamma_1,\lambda_1)$ and $\xi_2=(\Gamma_2,\lambda_2)$, define $\xi_1\blt\xi_2\coloneq(\Gamma_1\cup\Gamma_2,\lambda_1\cup\lambda_2)$, i.e., take the disjoint union of the two string diagram stratifications, then merge the shared edge.
      \item[\ref{Dpi1}] For a pinched product map $\pi\colon X\to P$ from a 1-ball $X$ to a 0-ball $P$ and a 0-field $a\in\eC^\cX_0(P)$, define $\pi^*a\coloneq(\varnothing,a)$, that is, $\pi^*a$ is the vertex-free string diagram on $X$ with a single edge labeled by $a$. 

    \item[\ref{Dfcj1}] For fields $a$ on 0-balls, set $\fcj a\coloneq a$ (as labels). For a 1-ball $X$ and a string diagram (pure 1-field) $\xi=(\Gamma,\lambda)\in\fld[1]{\eC^\cX}{X}[\fcj{a}\amalg b]$, define $\fcj{\xi}\in\fld[1]{\eC^\cX}{\orev X}[\fcj{b}\amalg a]$ by $\fcj{\xi}\coloneq(\orev\Gamma,\lambda)$ where $\orev\Gamma\coloneq((E,V),\orev\e)$ has reversed transverse orientations $\orev\e_v\coloneq-\e_v$ and unchanged labels. For example, for morphisms $f\in\cX(a\to b)$ and $g\in\cX(c\to b)$, 
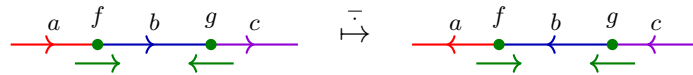
  \[
      \begin{tkz}
        \draw[thick, red, mid>={0.5}] (0,0) -- (1.15,0)
        node[midway, above=1pt, font=\scriptsize, text=black] {$a$};
        \draw[thick, blue!70!black, mid>={0.5}] (1.15,0) -- (2.65,0)
        node[midway, above=1pt, font=\scriptsize, text=black] {$b$};
        \draw[thick, violet, mid>={0.5}] (2.65,0) -- (3.8,0)
        node[midway, above=1pt, font=\scriptsize, text=black] {$c$};
        \node[circle, fill=Green, inner sep=1.5pt,
        label={[font=\scriptsize, label distance=0mm]above:{$f$}}] at (1.15,0) {};
        \draw[->, thick, Green] (0.85,-0.25) -- (1.45,-0.25);
        \node[circle, fill=Green, inner sep=1.5pt,
        label={[font=\scriptsize, label distance=0mm]above:{$g$}}] at (2.65,0) {};
        \draw[<-, thick, Green] (2.35,-0.25) -- (2.95,-0.25);
      \end{tkz}
    \quad
    \overset{\fcj{\,\cdot\,}}{\mapsto}
    \quad
      \begin{tkz}
        \draw[thick, red, mid<={0.5}] (0,0) -- (1.15,0)
        node[midway, above=1pt, font=\scriptsize, text=black] {$a$};
        \draw[thick, blue!70!black, mid<={0.5}] (1.15,0) -- (2.65,0)
        node[midway, above=1pt, font=\scriptsize, text=black] {$b$};
        \draw[thick, violet, mid<={0.5}] (2.65,0) -- (3.8,0)
        node[midway, above=1pt, font=\scriptsize, text=black] {$c$};
        \node[circle, fill=Green, inner sep=1.5pt,
        label={[font=\scriptsize, label distance=0mm]above:{$f$}}] at (1.15,0) {};
        \draw[->, thick, Green] (0.85,-0.25) -- (1.45,-0.25);
        \node[circle, fill=Green, inner sep=1.5pt,
        label={[font=\scriptsize, label distance=0mm]above:{$g$}}] at (2.65,0) {};
        \draw[<-, thick, Green] (2.35,-0.25) -- (2.95,-0.25);
      \end{tkz}
  \]
  where the green arrows depict the transverse orientations $\varepsilon$: $\varepsilon_v=+1$ if they agree with the underlying 1-ball orientation, and $\e_v=-1$ otherwise.
    \item[\ref{DU1}] For a 1-ball $X$ and $c_0,c_1\in\cX$, the subspace $\fldlU{\eC^\cX}{X}[c_0\amalg c_1]\subset\bC\{\fld[1]{\eC^\cX}{X}[c_0\amalg c_1]\}$ of local relations is generated by differences $\xi-\xi'$ where $\xi,\xi'$ are related in any of the following ways.
    \begin{lst}
      \item (\emph{Isotopies rel boundary}) $\xi'=H(-,1)_*\xi$ for an isotopy-rel-boundary $H\colon X\times I\to X$ of $X$ (consisting of orientation-preserving homeomorphisms).
      \item (\emph{Identity deletion}) $\xi$ has a vertex $v$ labeled $\id_x$, and $\xi'$ is obtained from $\xi$ by deleting $v$ and merging the adjacent edges into a single $x$-labeled edge.
      \item (\emph{Adjacent vertex composition}) $\xi$ has adjacent vertices $v_i<v_{i+1}$ with $\varepsilon_i=\varepsilon_{i+1}$, and $\xi'$ is obtained from $\xi$ by replacing both vertices with a single vertex (located anywhere in the closed interval with endpoints at $v_i$ and $v_{i+1}$) labeled by $\lambda_{v_{i+1}} \circ \lambda_{v_i}$ when $\varepsilon_i = +1$ and by $\lambda_{v_i} \circ \lambda_{v_{i+1}}$ when $\varepsilon_i = -1$.
      \item (\emph{$\e$-$\lambda$-flip}) $\xi'$ is obtained from $\xi$ by flipping both $\e_v\mapsto -\e_v$ and $\lambda_v\mapsto \lambda_v^\dag$ for some vertex $v$. For example, for morphisms $f\in\cX(a\to b)$ and $g\in\cX(c\to b)$, 
      \[
       \begin{tkz}
        \draw[thick, red, mid<={0.5}] (0,0) -- (1.15,0)
        node[midway, above=1pt, font=\scriptsize, text=black] {$a$};
        \draw[thick, blue!70!black, mid<={0.5}] (1.15,0) -- (2.65,0)
        node[midway, above=1pt, font=\scriptsize, text=black] {$b$};
        \draw[thick, violet, mid<={0.5}] (2.65,0) -- (3.8,0)
        node[midway, above=1pt, font=\scriptsize, text=black] {$c$};
        \node[circle, fill=Green, inner sep=1.5pt,
        label={[font=\scriptsize, label distance=0mm]above:{$f$}}] at (1.15,0) {};
        \draw[->, thick, Green] (0.85,-0.25) -- (1.45,-0.25);
        \node[circle, fill=Green, inner sep=1.5pt,
        label={[font=\scriptsize, label distance=0mm]above:{$g$}}] at (2.65,0) {};
        \draw[<-, thick, Green] (2.35,-0.25) -- (2.95,-0.25);
      \end{tkz}
    \;
    \overset{\substack{\varepsilon\text{-}\lambda\text{-flip}}}{\sim}\;\begin{tkz}[baseline=0.05ex]\draw[thick, red, mid<={0.5}] (0,0) -- (1.15,0)
      node[midway, above=1pt, font=\scriptsize, text=black] {$a$};\draw[thick, blue!70!black, mid<={0.5}] (1.15,0) -- (2.65,0)
      node[midway, above=1pt, font=\scriptsize, text=black] {$b$};\draw[thick, violet, mid<={0.5}] (2.65,0) -- (3.8,0)node[midway, above=1pt, font=\scriptsize, text=black] {$c$};\node[circle, fill=Green, inner sep=1.5pt,label={[font=\scriptsize, label distance=0mm]above:{$f^\dag$}}] at (1.15,0) {};\draw[<-, thick, Green] (0.85,-0.25) -- (1.45,-0.25);\node[circle, fill=Green, inner sep=1.5pt,label={[font=\scriptsize, label distance=0mm]above:{$g^\dag$}}] at (2.65,0) {};\draw[->, thick, Green] (2.35,-0.25) -- (2.95,-0.25);
    \end{tkz}.
  \]
\end{lst}
\end{lst}
\end{construction}

\subsubsection{Evaluation maps}
\begin{lemma}
    \label{lem:eval-quotient}
    If $\cX$ is a dagger $1$-category, then for a 1-ball $X$ and objects $a,b\in\cX$ the \defn{evaluation map}
    $\eval\colon\bC\{\fld[1]{\C^\cX}{X}[\fcj{a}\amalg b]\} \to\cX(a\to b)$ determined on string diagrams $\xi\in\fld[1]{\C^\cX}{X}[\fcj{a}\amalg b]$ by
    \begin{equation}
    \label{evalthingeqn}
        \eval(\xi)\coloneq
        \lambda_{v_m}'\circ\cdots\circ\lambda_{v_1}'
        \in\cX(a\to b)
        \qquad\text{where}\qquad
        \lambda_v'\coloneq
        \begin{cases}
            \lambda_v & \text{if }\e_v=+1,\\
            \lambda_v^\dag & \text{if }\e_v=-1,
        \end{cases}
    \end{equation}
    descends to an isomorphism $\bar\eval\colon\fldl{\C^\cX}{X}[\fcj{a}\amalg b]\overset\cong\to\cX(a\to b)$.
\end{lemma}

\begin{pf}
None of the four generators of $\lU{\C^\cX}$ change the evaluation, so $\lU{\C^\cX}\subset\ker(\eval)$. Conversely, for any $\cX$-string diagram $\xi$ on a 1-ball, we have $\eval(\xi)=\eval(\eta_{\eval(\xi)})$ where $\eta_{\eval(\xi)}$ is any single-vertex string diagram with vertex label $\eval(\xi)$ and transverse orientation $\varepsilon=+1$. Indeed, first apply $\varepsilon$-$\lambda$-flips to align the transverse orientations of all vertices, and then apply the adjacent vertex composition local relation until exactly one vertex remains. Thus $\xi - \eta_{\eval(\xi)} \in \lU{\C^\cX}$, so $\ker(\eval)\subset\lU{\C^\cX}$.
\end{pf}

In light of \cref{lem:reparam-gen}, we will often omit any decoration of $\bar\eval$ or $\langle-\rangle$ indicating the 1-ball $X$.

\begin{facts}
\label{factdag}
The assignment $f \mapsto \langle f \rangle\coloneq(\bar\eval)^{-1}(f)$ defines an identity-on-objects dagger functor $\cX \to \cX_{\C^\cX}$. More specifically, for a morphism $f \in \cX(a \to b)$, the 1-field $\langle f \rangle \in \fldl{\C^\cX}{I}[ \fcj{a}\amalg b]$ is in $\cX_{\C^\cX}(a\to b)$, and the following identities hold.
\begin{lst}[label=\text{($\langle$\arabic*$\rangle$)}]
\item\label{factdag-1} $\langle g \circ f \rangle = \langle g \rangle \circ \langle f \rangle$ for composable morphisms $f$ and $g$ in $\cX$;
\item\label{factdag-2}  $\langle \id_a \rangle = \id_a$ for all objects $a$ in $\cX$;
\item\label{factdag-3}  $\langle f^\dag \rangle=\langle f \rangle^\dag$ for all morphisms $f$ in $\cX$.
\end{lst}
\end{facts}

\begin{remark}
Observe that if $\cX$ is a dagger 1-category, then for a 1-ball $J$, a pure 1-field $\xi\in\fld[1]{\eC^\cX}{J}[\fcj{a}\amalg b]$ is an isometry (resp. coisometry, unitary equivalence) in $\eC^\cX$ if and only if $\eval(\xi)\in\cX(a\to b)$ is an isometry (resp. coisometry, unitary equivalence) in $\cX$. 
\end{remark}

\subsubsection{The sphere trace on \texorpdfstring{$\eC^\cX$}{C(X)}}

\

\begin{construction}[The sphere trace $\psn^{\C^\cX}$]
\label{cstr:unitary disk-like 1-categories from pre-2-Hilbert spaces}
Let $\cX$ be a pre-2-Hilbert space. We define a sphere trace on $\C^\cX$ as follows. Consider a string diagram (pure 1-field) $\xi\in\undfld[1]{\C^\cX}{S^1}$. Choose a point $p\in S^1\setminus V_\xi$ and let $e_p\in E_\xi$ be the edge of $\xi$ containing $p$. Now set $a_p\coloneq\lambda_{e_p}$. We define
\[
  \tld{\psn}^{\C^\cX}(\xi)\coloneq\Tr^\cX_{a_p}(\eval(\xi|_{S^1\setminus p}))
\]
where $\eval$ is the evaluation map $\fld[1]{\C^\cX}{S^1\setminus p}[\fcj{a_p}\amalg a_p]\to\cX(a_p\to a_p)$. For example, 
\[
\tld\psn^{\C^\cX}\Bigg(\!
\begin{tkz}[scale=0.6,thick,baseline=-0.5ex]
    \draw[mid>, red] (90:1) arc (90:210:1);
    \node[red, font=\scriptsize] at (150:1.35) {$a$};
    \draw[mid>, blue!70!black] (210:1) arc (210:330:1);
    \node[blue!70!black, font=\scriptsize,overlay] at (270:1.35) {$b$};
    \draw[mid>, orange!90!black] (330:1) arc (330:450:1);
    \node[orange!90!black, font=\scriptsize] at (30:1.35) {$c$};
    \node[circle, fill=Green, inner sep=1pt, label={[font=\scriptsize, Green]left:{$f$}}] at (210:1) {};
    \draw[->, Green] (210:1.3) ++(120:0.25) -- ++(-60:0.5);
    \node[circle, fill=Green, inner sep=1pt, label={[font=\scriptsize, Green]right:{$g$}}] at (327:1.1) {};
    \draw[<-, Green] (330:1.3) ++(240:0.25) -- ++(60:0.5);
    \node[above] at (90:1) {~};
    \node[circle, fill=Green, inner sep=1pt, label={[font=\scriptsize, Green,overlay]above:{$h$}},overlay] at (90:1.1) {};
    \draw[->, Green,overlay] (90:1.3) ++(0:0.25) -- ++(180:0.5);
    \draw[dashed, black, thin] (180:0.75) -- (180:1.25);
    \node[font=\scriptsize] at (180:0.55) {$p$};
\end{tkz}
\Bigg)
=
\Tr_a^\cX(\eval(
\begin{tkz}[scale=0.6,thick,baseline=-0.5ex,x=0.75cm]
    \draw[mid>={0.5}, blue!70!black] (0.8,0) -- (3.2,0) node[midway,above=1pt,font=\scriptsize,text=blue!70!black] {$b$};
    \draw[mid>={0.5}, orange!90!black] (3.2,0) -- (5.6,0) node[midway,above=1pt,font=\scriptsize,text=orange!90!black] {$c$};
    \draw[mid>={0.75}, red] (5.6,0) -- (7.4,0) node[midway,above=1pt,font=\scriptsize,text=red] {$a$};
    \node[circle, fill=Green, inner sep=1pt, label={[font=\scriptsize, Green, label distance=0mm]above:{$f$}}] at (0.8,0) {};
    \draw[->, Green] (0.5,-0.3) -- (1.1,-0.3);
    \node[circle, fill=Green, inner sep=1pt, label={[font=\scriptsize, Green, label distance=0mm]above:{$g$}}] at (3.2,0) {};
    \draw[<-, Green] (2.9,-0.3) -- (3.5,-0.3);
    \node[circle, fill=Green, inner sep=1pt, label={[font=\scriptsize, Green, label distance=0mm]above:{$h$}}] at (5.6,0) {};
    \draw[->, Green] (5.3,-0.3) -- (5.9,-0.3);
\end{tkz}
))
=
\Tr_a^\cX(h \circ g^\dag \circ f).
\]
This assignment is independent of the choice of $p$, and so extends to a linear map $\tld{\psn}^{\C^\cX}\colon\bC\{\undfld[1]{\C^\cX}{S^1}\}\to\bC$, which itself descends to a well-defined linear map $\psn^{\C^\cX}\colon\undfldl[1]{\C^\cX}{S^1}\to\bC$ satisfying \ref{1psnP}.
\end{construction}

\begin{proposition}
\label{prop:unitary disk-like 1-categories from pre-2-Hilbert spaces}
If $\cX$ is a pre-2-Hilbert space, then $\C^\cX$ equipped with the sphere trace $\psn^{\C^\cX}$ of \cref{cstr:unitary disk-like 1-categories from pre-2-Hilbert spaces} is a unitary disk-like 1-category. If $\cX$ is moreover finite, then $\C^\cX$ is finite. 
\end{proposition}

\begin{pf}
We first show $\psn^{\eC^\cX}$ is well-defined.

\itemstep{$\tld\psn^{\C^\cX}$ is independent of the cut point $p$.}
Let $Y\hookrightarrow S^1$ denote the arc from $q$ (incoming) to $p$ (outgoing) on $S^1$. Then for $Y'\coloneq S^1\setminus Y$, we have
\begin{multline*}
  \tld{\psn}^{\C^\cX}(\xi)
  \defeq{}\Tr^\cX_{a_q}(\eval_{S^1\setminus q}(\xi|_{S^1\setminus q}))
  =\Tr^\cX_{a_q}(\eval_{Y'}(\xi|_{Y'})\circ\eval_{Y}(\xi|_{Y}))
  \\
  \underset{\ref{dagTr1}}{=}\Tr^\cX_{a_p}(\eval_{Y}(\xi|_{Y})\circ\eval_{Y'}(\xi|_{Y'}))
  =\Tr^\cX_{a_p}(\eval_{S^1\setminus p}(\xi|_{S^1\setminus p}))
  \defeq{}\tld{\psn}^{\C^\cX}(\xi).
\end{multline*}
Thus we get a linear map $\tld{\psn}^{\C^\cX}\colon\bC\{\undfld[1]{\C^\cX}{S^1}\}\to\bC$.

\itemstep{$\tld{\psn}^{\C^\cX}$ descends to a linear map $\psn^{\C^\cX}\colon\undfldl[1]{\C^\cX}{S^1}\to\bC$.}
We need to show $\fldlU{\C^\cX}{S^1}\subset\ker\tld{\psn}^{\C^\cX}$. Recall that $\fldlU{\C^\cX}{S^1}$ is spanned by elements $\eta$ of $\bC\{\undfld[1]{\C^\cX}{S^1}\}$
represented in the colimit $\undfld[1]{\C^\cX}{S^1}$ by a splitting $\{\eta_1,\dots,\eta_N\}$
of $\eta$ corresponding to some permissible ball decomposition $\{X_j\}$ of $S^1$ such that
$\eta_j\in\fldlU{\C^\cX}{X_j}[\fcj{a_p}\amalg a_q]$ for some $j$ and some objects $a_p,a_q\in\cX$,
where $X_j\hookrightarrow S^1$ has incoming and outgoing boundary points $p$ and $q$
respectively in $S^1$. If we set $\zeta\coloneq\glu\{\textstyle\bigcup_{i\neq j}\eta_i\}$, then
\[
  \tld{\psn}^{\C^\cX}(\eta)
  \defeq\Tr^\cX_{a_p}(\eval_{S^1\setminus p}(\eta_j\blt\zeta))
  =\Tr_{a_p}^\cX(\eval_{S^1\setminus X_j}(\zeta)\circ\,\underbrace{\eval_{X_j}(\eta_j)}_{=0})
  =0,
\]
so $\fldlU{\C^\cX}{S^1}\subset\ker\tld{\psn}^{\C^\cX}$. Thus $\tld{\psn}^{\C^\cX}$ descends to a well-defined linear map $\psn^{\C^\cX}\colon\undfldl[1]{\C^\cX}{S^1}\to\bC$.

\itemstep{\ref{1psnP}.}
The induced pairings $\bkt{ - }{ - }_{X,c}$ satisfy \ref{1psnP} by \ref{dagTr2}.

\itemstep{\ref{1psnF}.} 
Now suppose that $\cX$ is finite.
  It suffices to show $\dim\undfldl[1]{\C^\cX}{S^1}<\infty$ and that $\dim\undfldl[1]{\C^\cX}{D^1}[\fcj{a}\amalg b]<\infty$ for all $a,b\in\cX$. 
  The latter is immediate because $\undfldl[1]{\C^\cX}{D^1}[\fcj{a}\amalg b]=\cX(a\to b)$ is finite-dimensional for each $a,b\in\cX$. 
  For the former, $\Sk_{\C^\cX}(\pt)\defeq\cX_{\C^\cX}\cong^\dag\cX$ by
  \cref{lem:eval-quotient} and \cref{factdag}, so $\Sk_{\C^\cX}(\pt)$ is unitary; as $\cX$ is finite, we can write $\Irr(\Sk_{\C^\cX}(\pt)^\cent)=\{(r_i,p_i)\}_{i=1}^N$ with $N<\infty$. Applying \cref{glu-span}
  to the gluing data of $D^1$ along $\partial D^1$ to form $S^1$ gives that $\undfldl[1]{\C^\cX}{S^1}=\sum_{i=1}^N\Gamma(p_iV_{r_i}p_i)$ with $p_iV_{r_i}p_i=p_i\cX(r_i\to r_i)p_i=\bC p_i$, so $\dim\undfldl[1]{\C^\cX}{S^1}\leq N<\infty$. 
\end{pf}

\subsection{The round-trip equivalences}

\subsubsection{The first round-trip: \texorpdfstring{$\cX_{\eC^\cX}\cong^\dag\cX$}{X(C(X)) = X}}
An \defn{equivalence} of dagger categories is a fully faithful dagger functor between them that is unitarily essentially surjective in that the isomorphisms witnessing essential surjectivity can be chosen to be unitary. It is well-known that all isomorphisms in a \emph{unitary} category can be augmented to unitary isomorphisms, so mere essential surjectivity suffices.

A \defn{unitary weak equivalence} of disk-like 1-categories is a weak equivalence of disk-like 1-categories that is \defn{unitarily essentially surjective}, i.e., for each 0-field $b$ in $\D$, there is a 0-field $a$ in $\C$ and a unitary equivalence $\eF(a)\cong^\star b$ in $\D$.

An \defn{isometric weak equivalence} of unitary disk-like 1-categories $\C$ and $\D$ is a weak equivalence $\eF\colon\C\to\D$ such that $\psn^\C=\psn^\D\circ\widehat\eF$.

\begin{lemma}
\label{lem:functor-commutes-closure}
A disk-like functor $\eF\colon\eC^\cX\to\D$ satisfies $\eF(\cl_a \xi)=\cl_{\eF(a)}\eF(\xi)$ in $\undfldl{\D}{S^1}$. Thus if $\psn^\D(\cl_{\eF(a)}\eF\langle f\rangle)=\psn^{\eC^\cX}(\cl_a\langle f\rangle)$ for all $f\in\cX(a\to a)$, then $\eF$ is isometric.
\end{lemma}

\begin{pf}
Since $\eF$ is a disk-like functor, it preserves gluing and product fields, so $\eF(\cl_a \xi) = \cl_{\eF(a)}\eF(\xi)$ as 1-fields. To see $\eF$ is isometric, note that by local relations any string diagram $\xi$ on $S^1$ is represented by a single-vertex string diagram on $S^1$, so $[\xi] = \cl_a\langle f\rangle$ for some $a\in\cX$ and some $f\in\cX(a\to a)$. Thus
\[
    \psn^\D(\eF([\xi])) = \psn^\D(\eF(\cl_a \langle f\rangle)) = \psn^\D(\cl_{\eF(a)}\eF\langle f\rangle) = \psn^{\eC^\cX}(\cl_a \langle f\rangle ) = \psn^{\eC^\cX}([\xi]).
\]
As any 1-field $\xi\in\undfldl{\eC^\cX}{S^1}$ is a linear combination of such 1-fields $[\xi]$, by linearity this completes the proof.
\end{pf}

\begin{proposition}
\label{prop:XCXcongX}
For a pre-2-Hilbert space $\cX$, there is a canonical isometric equivalence $\langle-\rangle\colon\cX\to\cX_{\eC^\cX}$.
\end{proposition}

\begin{pf}
By \cref{lem:eval-quotient} and \cref{factdag}, $\langle-\rangle$ is a dagger equivalence. For $f \in \cX(a\to a)$, we have $\Tr_a^{\cX_{\eC^\cX}}(\langle f\rangle) \defeq{} \psn^{\eC^\cX}(\cl_a \langle f\rangle) \defeq{} \Tr_a^\cX(f)$, so $\langle-\rangle$ is isometric.
\end{pf}

\begin{corollary}
\label{cor:XCXcongX implies homs}
For finite pre-2-Hilbert spaces $\cX$ and $\cY$, we have isometric equivalences of finite pre-2-Hilbert spaces
    \[
        \Hom(\cX\to\cY)
        \cong^\dag
        \Hom(\cX_{\eC^\cX}\to\cY)
        \cong^\dag
        \Hom(\cX\to\cX_{\eC^\cY})
        \cong^\dag
        \Hom(\cX_{\eC^\cX}\to\cX_{\eC^\cY}).
    \]
\end{corollary}

\subsubsection{Unrestricted disk-like functors}

\

The following unpacks \cref{def: DL functor} in the case $n=1$.
\begin{definition}
\label{def: DL functor1}
Let $\eC$ and $\D$ be disk-like 1-categories. A \defn{disk-like functor} $\eF\colon\eC\to\D$ is the data of an assignment to each $k\in\{0,1\}$, $k$-ball $W$, and pure $k$-field $\sigma\in\fld[k]{\C}{W}$ of a pure $k$-field $\eF_W(\sigma)\in\fld[k]{\D}{W}[\eF(\partial \sigma)]$. (Here $\eF(\partial \sigma)$ is defined by applying $\eF$ componentwise to the components of a splitting of $\partial W$ in the colimit; see \cref{prop:functor-extension}.) This data is subject to the following conditions for $k\in\{0,1\}$.
\begin{lst}
\item[($\eF\varphi$)]\label{1axiom:eF-nat}
$\eF_{W'}(\varphi_*\xi)=\varphi_*\eF_W(\xi)$ for all $k$-ball homeomorphisms $\varphi\colon W\overset\sim\to W'$ and $k$-fields $\xi\in\fld[k]{\C}{W}$.
\item[($\eF\partial$)]\label{1axiom:eF-boundary} 
$\partial(\eF_W(\xi)) = \eF_{\partial W}(\partial \xi)$ for all $k$-fields $\xi\in\fld[k]{\C}{W}$.
\item[($\eF$G)]\label{1axiom:eF-gluing} 
$\eF_{X_1\cup_P X_2}(\xi\blt_P \eta) = \eF_{X_1}(\xi)\blt_P \eF_{X_2}(\eta)$ for all pure 1-fields $\xi\in\fld[1]{\C}{X_1}$ and $\eta\in\fld[1]{\C}{X_2}$ compatible along a 0-ball $P$.
\item[($\eF\pi$)]\label{1axiom:eF-products} 
$\eF_X(\pi^*a) = \pi^* \eF_P(a)$ for all pinched product maps $\pi \colon X \to P$ and all 0-fields $a\in\fld[0]{\C}{P}$.
\item[($\eF\kern.2ex\fcj{\,\cdot\,}$)]\label{1axiom:eF-refl}
$\eF_{\orev{W}}(\fcj{\xi}) = \fcj{\eF_W(\xi)}$ for all $k$-fields $\xi\in\fld[k]{\C}{W}$.
\item[($\eF$U)]\label{1axiom:eF-local-relations} 
When $k=1$, the linear extension $\bC\{\fld[1]{\C}{X}[ c]\} \to \bC\{\fld[1]{\D}{X}[ \eF_{\partial X}(c)]\}$ of $\eF$ sends $\fldlU{\C}{X}[c]$ into $\fldlU{\D}{X}[ \eF_{\partial X}(c)]$ for all 1-balls $X$ and $c\in\undfld[0]{\C}{\partial X}$.
\end{lst}
Observe that \ref{1axiom:eF-local-relations} gives rise to well-defined linear maps $\fldl[1]{\C}{X}[c]\to\fldl[1]{\D}{X}[\eF(c)]$, which we denote by $\widehat\eF_{X,c}$ (or sometimes also by $\eF$).
\end{definition}

One may hope that the above notion of disk-like functor is flexible enough to map between the main examples of disk-like 1-categories, but unfortunately this is not the case. To see why, try to construct a valid disk-like functor $\C\to\C^{\cX_\C}$ or $\C^{\cX_\C}\to\C$, or more concretely, try to define a functor between a string diagram disk-like 1-category and a disk-like 1-category given by maps into some fixed space. The situation is the same when working with disk-like $n$-categories for general $n$. Thus we must weaken the definition to allow for the choice of regular neighborhoods of strata. We make this formal with the following framework of ``unrestricted disk-like functors'' between disk-like 1-categories. See \cref{def:parameterized-splitting-filtered} for the general-$n$ case; the following is the specialization of it to $n=1$.

If we were to define disk-like functors in the top dimension $n$ directly as linear maps between vector spaces of $n$-fields instead of first defining an underlying set map of pure $n$-fields that linearizes to one, then the unrestricted disk-like functor framework would not be necessary for $n=1$. But given that when $n=2$ we will need parameterized splittings in the sense defined below to define functors on (pure!) 1-fields, we deem it best to introduce the unrestricted disk-like functor framework for $n=1$ too. See \cref{sec:unnrestricted} for the definitions for general $n$.

The problem outlined above suggests a general distinction between disk-like $n$-categories of string diagrams and disk-like $n$-categories given by maps to some fixed space. Another reason we include the unrestricted functor framework is that we want to apply it to be able to discuss TQFTs of the latter type. 

\begin{definition}[Parameterized splitting]
Let $\C$ be a disk-like 1-category and let $X$ be a 1-manifold. A \defn{parameterized splitting} $\Theta$ of a pure 1-field $\sigma\in\fld[1]{\C}{X}$ with respect to a finite subset $V\subset \Int(X)$ is a collection of homeomorphisms
\[
  \Theta=\{\theta_v\colon D^1\overset\sim\to N_v \mid v\in V\}
\]
such that the $N_v$ have pairwise disjoint interiors, $\sigma$ splits along $\bigcup_{v\in V}\partial N_v$, $N_v\subset\Int(X)$ for each $v\in V$, and $\theta_v(0)=v$ for all $v\in V$.

We write $\mathrm{PS}(\sigma,V)$ for the collection of all parameterized splittings $\Theta$ with respect to $V$ and $\mathrm{PS}(\sigma)$ for the collection of $\Theta$ that are parameterized splittings with respect to some choice of $V$. For $\Theta\in\mathrm{PS}(\sigma)$, we denote by $V_\Theta$ the corresponding subset $V_\Theta\subset X$ with respect to which $\Theta$ is a parameterized splitting.
\end{definition}

The data of $\C$ acts on parameterized splittings as follows.
\begin{lst}  \item[($\Theta\varphi$)]\label{axiom:Theta-phi1} For all $\Theta \in \mathrm{PS}(\sigma)$ and all homeomorphisms $\varphi\colon X\to X'$, we define $\varphi_*\Theta\coloneq\{\varphi\circ\theta_v\colon D^1\to\varphi(N_v)\mid v\in V_\Theta\}\in\mathrm{PS}(\varphi_*\sigma,\varphi(V_\Theta))$.
  \item[($\Theta$G)]\label{1axiom:Theta-G} For all $\Theta_1\in\mathrm{PS}(\sigma_1)$ and $\Theta_2\in\mathrm{PS}(\sigma_2)$ such that $\sigma_1$ and $\sigma_2$ are compatible along a 0-ball $P$, we define $\Theta_1\blt_P\Theta_2\coloneq\Theta_1\amalg\Theta_2=\{\theta_v\in\Theta_i\mid i=1,2\}\in\mathrm{PS}(\sigma_1\blt_P\sigma_2,V_{\Theta_1}\amalg V_{\Theta_2})$.
  \item[($\Theta\orev{\,\cdot\,}$)]\label{axiom:Theta-refl1} For all $\Theta\in\mathrm{PS}(\sigma)$, $\orev\Theta\coloneq\{\orev{\theta_v}\circ\iota\colon D^1\to\orev{N_v}\mid v\in V_\Theta\}$ where $\orev{\theta_v}\coloneq\orev{\,\cdot\,}\circ\theta_v\colon\orev{D^1}\to\orev{N_v}$ and $\orev{V_\Theta}\coloneq\{\orev{v}\mid v\in V_\Theta\}\subset \orev{X}$.
\end{lst}

\begin{remark}
\label{1.33}
When $\C$ is a disk-like 1-category of string diagrams and $\xi\in\fld[1]{\C}{X}$, we can take parameterized splittings with respect to $V$ given by the vertices of $\xi$ itself. In this case $\Theta$ is just a choice of regular neighborhoods of the string diagram strata in $\xi$.
\end{remark}

\begin{definition}
A \defn{system of stratifications} for $\C$ is an assignment to each 1-manifold $X$ and pure 1-field $\sigma\in\undfld[1]{\C}{X}$ of a finite subset $V_\sigma\subset\Int(X)$ subject to the following conditions.
\begin{lst}
  \item[($V\varphi$)]
  \label{axiom:Gamma-nat1} $V_{\varphi_*\sigma}=\varphi(V_\sigma)$ for all homeomorphisms $\varphi\colon X\to X'$.
  \item[($V$G)]
  $V_{\sigma_1\blt_P\sigma_2}=V_{\sigma_1}\amalg V_{\sigma_2}$ for all $\sigma_1$ and $\sigma_2$ compatible along $P$.
  \item[($V\orev{\,\cdot\,}$)]\label{axiom:Gamma-refl1}
  \label{axiom:Gamma-refl} 
  $V_{\fcj\sigma}=\orev{V_\sigma}$.
  \item[($V\pi$)]\label{axiom:Gamma-products1} $V_{\pi^*a}=\varnothing$ for all $a\in\fld[0]{\C}{P}$ and all pinched product maps $\pi\colon X\to P$.
\end{lst}
\end{definition}

By \ref{axiom:Gamma-nat1}--\ref{axiom:Gamma-products1}, the actions \ref{axiom:Theta-phi1}--\ref{axiom:Theta-refl1} define maps $\varphi_*\colon\eF_{\mathrm{PS}}(\sigma)\to\eF_{\mathrm{PS}}(\varphi_*\sigma)$, $\partial\colon\eF_{\mathrm{PS}}(\sigma)\to\eF_{\mathrm{PS}}(\partial\sigma)$, 
$-\blt_E-\colon\eF_{\mathrm{PS}}(\sigma_1)\times\eF_{\mathrm{PS}}(\sigma_2)\to\eF_{\mathrm{PS}}(\sigma_1\blt_E\sigma_2)$, and $\orev{\,\cdot\,}\colon\eF_{\mathrm{PS}}(\sigma)\to\eF_{\mathrm{PS}}(\fcj\sigma)$.
\begin{definition}[Unrestricted disk-like functor]
\label{def:1unrestricted-functor}
For disk-like 1-categories $\C$ and $\D$, an \defn{unrestricted disk-like functor} $\eF\colon\C\to\D$ consists of the following data.
\begin{lst}
  \item[(D$\eF\Theta\Gamma$)]\label{1axiom:FGamma-data}
  A system of string diagram stratifications $\sigma\mapsto V_\sigma$ for $\C$.
  \item[(1$\eF\Theta_0$)]
  To each 0-manifold $P$ and each 0-field $a\in\fld[0]{\C}{P}$, a 0-field $\eF(a)\in\fld[0]{\D}{P}$.
  \item[(1$\eF\Theta_1$)] 
  For each 1-manifold $X$ and each pure 1-field $\sigma\in\fld[1]{\C}{X}$, a set map $\eF(\sigma,-)\colon\eF_{\mathrm{PS}}(\sigma)\to\undfld[1]{\D}{X}$ where $\eF_{\mathrm{PS}}(\sigma)\coloneq\mathrm{PS}(\sigma,V_\sigma)$.
\end{lst}
The above data for 0-fields are subject to conditions \ref{1axiom:eF-nat}, \ref{1axiom:eF-boundary} and \ref{1axiom:eF-refl} in the case $k=0$, together with the following conditions regarding 1-fields.
\begin{lst}
  \item[(F$\Theta\varphi$)]\label{1axiom:FTheta-nat} $\eF(\varphi_*\sigma,\varphi_*\Theta)=\varphi_*\eF(\sigma,\Theta)$.
  \item[(F$\Theta\partial$)]\label{1axiom:FTheta-bdy} $\eF(\partial\sigma)=\partial\eF(\sigma,\Theta)$ for all $\Theta\in\eF_{\mathrm{PS}}(\sigma)$.
  \item[(F$\Theta$G)]\label{1axiom:FTheta-gluing} $\eF(\sigma_1\blt_E\sigma_2,\Theta_1\blt_E\Theta_2)=\eF(\sigma_1,\Theta_1)\blt_E\eF(\sigma_2,\Theta_2)$ for gluings along any 0-manifold $E$, whenever each of (i) $\sigma_1$ and $\sigma_2$, (ii) $\Theta_1$ and $\Theta_2$, and (iii) $\eF(\sigma_1,\Theta_1)$ and $\eF(\sigma_2,\Theta_2)$ are compatible along $E$.
  \item[(F$\Theta\orev{\,\cdot\,}$)]\label{1axiom:FTheta-refl} $\eF(\fcj\sigma,\orev\Theta)=\fcj{\eF(\sigma,\Theta)}$.
  \item[(F$\Theta\pi$)]
  \label{1axiom:FTheta-products} $\eF(\pi^*a,\Theta)=\pi^*\eF(a)$ for all pinched product maps $\pi\colon X\to P$ and $a\in\fld[0]{\C}{P}$. (Since $V_{\pi^*a}=\varnothing$, $\eF_{\mathrm{PS}}(\pi^*a)=\{\varnothing\}$, so the $\Theta$ here is empty. In other words, this condition is simply \ref{1axiom:eF-products}.)
  \item[(F$\Theta$U)]\label{1axiom:FTheta-local-relations} If $X$ is a 1-manifold, $c\in\undfld[0]{\C}{\partial X}$, and $\sum_i\lambda_i\sigma_i\in\fldlU{\C}{X}[c]$, then $\sum_i\lambda_i\eF(\sigma_i,\Theta_i)\in\fldlU{\D}{X}[\eF(c)]$ for any $\Theta_i\in\eF_{\mathrm{PS}}(\sigma_i)$.
\end{lst}
\end{definition}

\begin{remark}
\label{1dldldl}
By applying \ref{1axiom:FTheta-local-relations} to the case $\sum_i\lambda_i\sigma_i=\sigma-\sigma\,(=0)$, we get that  $[\eF(\sigma,\Theta)]=[\eF(\sigma,\Theta^\prime)]$ (that is, $\eF(\sigma,\Theta)-\eF(\sigma,\Theta^\prime)\in\fldlU{\D}{X}[\eF(c)]$) for all $\Theta,\Theta^\prime\in\eF_{\mathrm{PS}}(\sigma)$. Thus $\eF$ descends to a linear map $\widehat\eF_{X,c}\colon\fldl[1]{\C}{X}[c]\to\fldl[1]{\D}{X}[\eF(c)]$.
\end{remark}

Disk-like functors are examples of unrestricted disk-like functors; see \cref{dlurl}.

\subsubsection{The second round-trip: \texorpdfstring{$\eC^{\cX_\eC}\cong^\dag\eC$}{C(X(C)) = C}}
Let $\C$ be a disk-like $1$-category. For each liftable morphism $f$ in $\cX_\C$ (see \cref{sec:liftable}), choose a pure 1-field $\sigma_f\in\fld[1]{\C}{D^1}$ with $[\sigma_f]=f$. For the orientation-preserving homeomorphism $\iota\colon D^1\to\orev{D^1}$ given by $x\mapsto -x$, recall that $\sigma^\dag\coloneq\fcj{\iota_*\sigma}\in\fld[1]{\C}{D^1}$ satisfies $[\sigma^\dag]=[\sigma]^\dag$ and $\iota_*\sigma^\dag=\fcj\sigma$. For each transversally oriented vertex $(v,\varepsilon_v)$ in $\xi$, set $\sigma_v\coloneq\sigma_{\lambda_v}^{(\varepsilon_v)}$ where $\sigma_{\lambda_v}^{(+1)}\coloneq\sigma_{\lambda_v}$ and $\sigma_{\lambda_v}^{(-1)}\coloneq\sigma_{\lambda_v}^\dag$, so that $[\sigma_v]=\lambda_v'$ (see \cref{lem:eval-quotient}) and $\sigma_{\orev v}=\sigma_v^\dag$. 

\begin{construction}
\label{eLn1construction}
\label{ex:eL-n1}
We define an unrestricted disk-like functor $\eL\colon\tld\C^{\cX_\C}\to\C$ on 0-fields by $\eL(a)\coloneq a$ and on pure 1-fields $\xi\in\fld[1]{\tld\C^{\cX_\C}}{X}$ by the assignment $\mathrm{PS}(\xi,V_\xi)\to\fld[1]{\C}{X}$ given by
\[
  \eL(\xi,\Theta)
  \coloneq
  \glu
  \{\,
  \{
  (\theta_v)_* \sigma_{v} \}_{v\in V_\xi}
  \cup
  \{
  \pi_e^*\lambda_e
  \}_{e\in E_\xi}
  \,\}
\]
where $\hat e\coloneq e\setminus\bigcup_{v\in V_\xi}\Int(N_v)$ and $\pi_e$ is the unique pinched product map $\pi_e\colon \hat e\to \pt$ for $e\in E_\xi$. For example,
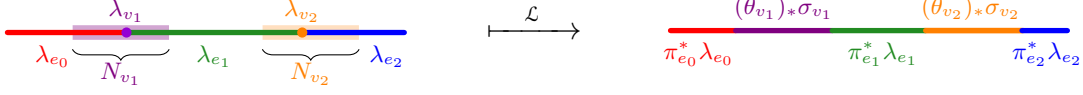
\[
\begin{tkz}[scale=0.75,baseline=-.475ex,cap=round]
    \draw[line width=5.5pt,draw=Purple, opacity=.35,cap=butt] (1.15,0)--(2.85,0);
    \draw[line width=5.5pt, draw=orange, opacity=.25,cap=butt] (4.5,0)--(6.2,0);
    \draw[line width=2.0pt, draw=red] (0,0)--(2.1,0);
    \draw[line width=2.0pt, draw=ForestGreen] (2.1,0)--(5.2,0);
    \draw[line width=2.0pt, draw=blue] (5.2,0)--(7,0);
    \node[circle, fill=violet, inner sep=1.35pt] at (2.1,0) {};
    \node[circle, fill=orange, inner sep=1.35pt] at (5.2,0) {};
    \node[below,text=red] at (.8,-.03) {\scriptsize$\lambda_{e_0}$};
    \node[below,text=ForestGreen] at (3.65,-.03) {\scriptsize$\lambda_{e_1}$};
    \node[below,text=blue] at (6.7,-.03) {\scriptsize$\lambda_{e_2}$};
    \node[above,text=Purple] at (2.1,0) {\scriptsize$\lambda_{v_1}$};
    \node[above,text=orange] at (5.2,0) {\scriptsize$\lambda_{v_2}$};
    \draw[decorate,decoration={brace,mirror,amplitude=4pt}] (1.15,-.3)--(2.85,-.3) node[Purple,midway,below=1pt] {\scriptsize$N_{v_1}$};
    \draw[decorate,decoration={brace,mirror,amplitude=4pt}] (4.5,-.3)--(6.2,-.3) node[orange,midway,below=1pt] {\scriptsize$N_{v_2}$};
\end{tkz}
\qquad\xmapsto{\quad\eL\quad}\qquad
\begin{tkz}[scale=0.75,cap=round]
    \draw[line width=2.0pt, draw=red] (0,0)--(1.15,0);
    \draw[Purple,line width=2.0pt] (1.15,0)--(2.85,0);
    \draw[line width=2.0pt, draw=ForestGreen] (2.85,0)--(4.5,0);
    \begin{pgfonlayer}{back}
        \draw[line width=2.6pt, draw=orange, opacity=.18,cap=round] (4.5,0)--(6.2,0);
    \end{pgfonlayer}
    \draw[line width=2.0pt, draw=orange] (4.5,0)--(6.2,0);
    \draw[line width=2.0pt, draw=blue] (6.2,0)--(7,0);
    \node[below,text=red] at (.5,0) {\scriptsize$\pi_{e_0}^*\lambda_{e_0}$};
    \node[above,text=Purple] at (2,0) {\scriptsize$(\theta_{v_1})_*\sigma_{v_1}$};
    \node[below,text=ForestGreen] at (3.75,0) {\scriptsize$\pi_{e_1}^*\lambda_{e_1}$};
    \node[above,text=orange] at (5.3,0) {\scriptsize$(\theta_{v_2})_*\sigma_{v_2}$};
    \node[below,text=blue] at (6.65,0) {\scriptsize$\pi_{e_2}^*\lambda_{e_2}$};
\end{tkz}
\]
where we are suppressing the transverse orientations on the vertices.
\end{construction}

\begin{pf}
The proof consists of straightforward computations; we give one as an example and leave the rest to the reader.
\itemstep{\ref{1axiom:FTheta-refl}.} Let $\xi\in\fld[1]{\tld\C^{\cX_\C}}{X}$ and let $v\in V_\xi$. Since in $\fld[1]{\C}{\orev{N_v}}$ we have
\begin{equation}\label{vertid}
(\orev{\theta_v}\circ\iota)_*\sigma_v^\dag \underset{\ref{Cvarphi1}}{=} \orev{\theta_v}_*\,\iota_*\sigma_v^\dag 
=
\orev{\theta_v}_*\fcj{\sigma_v} 
\underset{\ref{Cfcj1}}{=} 
\fcj{(\theta_v)_*\sigma_v},
\end{equation}
for an $\cX_\C$-string diagram $\xi\in\fld[1]{\tld\C^{\cX_\C}}{X}$ we compute that
\begin{align*}
&\fcj{\eL(\xi,\Theta)}
=
\fcj{\glu\{\{(\theta_v)_*\sigma_v\}_{v\in V_\xi}
\cup
\{\pi_e^*\lambda_e\}_{e\in E_\xi}\}}
% \\& 
\underset{\substack{\ref{Cfcj1}\&\ref{CGa1}}}{=}
\glu\{
\{\fcj{(\theta_v)_*\sigma_v}
\}_{v\in V_\xi}
\cup
\{\fcj{\pi_e^*\lambda_e}\}_{e\in E_\xi}
\}
\\
&\!\!\!\!\!\!\underset{\substack{\eqref{vertid}\&\ref{Cfcj1}}}{=}
\glu\{\{(\orev{\theta_v}\circ\iota)_*\sigma_v^\dag\}_{v\in V_\xi}
\cup
\{(\orev{\pi_e})^*\fcj{\lambda_e}\}_{e\in E_\xi}
\}
% \\&
=
\glu\{\{(\orev{\theta_v}\circ\iota)_*\sigma_{\orev v}\}_{v\in V_\xi}\cup\{(\orev{\pi_e})^*\lambda_{\orev e}\}_{e\in E_\xi}\}
% \\& 
=
\eL(\fcj\xi,\orev\Theta).\qedhere
\end{align*}
\end{pf}

Let $\widetilde{\eC}^{\cX_\eC}$ denote the disk-like 1-category obtained from $\eC^{\cX_\eC}$ by only considering string diagrams whose vertices are labeled by morphisms in $\cX_\C$ represented by pure fields in $\C$. We call such morphisms \defn{liftable} (see \cref{sec:liftable}). 

\begin{proposition}
\label{prop:CXCcongC}
    For a unitary disk-like 1-category $\eC$, there is an isometric weak equivalence $\eC\cong^\dag\eC^{\cX_\eC}$ witnessed by the zig-zag $\eC\xleftarrow[\eqref{eLn1construction}]{\eL}\widetilde{\eC}^{\cX_\eC}\xrightarrow[\eqref{lem:tilde-inclusion}]{\iota}\eC^{\cX_\eC}$.
\end{proposition}

\begin{pf}
By \cref{lem:tilde-inclusion}, the inclusion $\iota\colon\widetilde{\eC}^{\cX_\eC}\to\eC^{\cX_\eC}$ is an isometric weak equivalence. It remains to show $\eL$ is an isometric weak equivalence (see \cref{def:unrestricted-weak-equivalence}).
Observe that $\widehat{\eL}_{X,c}=\bar\eval\circ\widehat\iota$, so that because $\fldlU{\tld{\eC}^{\cX_\eC}}{X}[c]=\ker(\eval\circ\iota)$, the linear map $\widehat{\eL}_{X,c}\colon\fldl[1]{\tld\eC^{\cX_\eC}}{X}[c]\to\fldl[1]{\C}{X}[c]$ is injective. To see $\widehat{\eL}_{X,c}$ is surjective, recall that liftable 1-fields span $\fldl[1]{\C}{X}[c]$ and observe that if $f\in\fldl[1]{\C}{X}[c]$ is liftable, then any choice of pure 1-field $\sigma_f\in\fld[1]{\C}{X}$ with $[\sigma_f]=f$ has $\widehat{\eL}_{X,c}([\xi_f])=f$ for the single-vertex string diagram $\xi_f$ labeled by $f$; surjectivity follows by linearity. To see $\eL$ is isometric, let $\xi\in\fld[1]{\tld{\eC}^{\cX_\eC}}{S^1}$ and choose any point $p\in S^1\setminus V_\xi$. If $e_p$ and $a_p$ denote the corresponding edge and its label, then
\begin{align*}
\psn^\C([\eL(\xi,\Theta)]) 
&= \psn^\C\left(\cl_{a_p}\left(\glu\{\{(\theta_v)_*\sigma_v\}_{v}\cup\{\pi_e^*\lambda_e\}_{e}\}\Big|_{S^1\setminus p}\right)\right) 
\\&
= \psn^\C(\cl_{a_p}(\eval(\xi|_{S^1\setminus p}))) 
% \\&
= \Tr^{\cX_\C}_{a_p}(\eval(\xi|_{S^1\setminus p})) 
= \psn^{\eC^{\cX_\C}}([\xi]). \qedhere
\end{align*}
\end{pf}

\subsection{Representations, completion, and Hom-spaces}

\begin{example}[Faithful positive trace on functor categories {{\cite[Example 2.15]{3Hilb}}}]
  \label{faithful positive trace on functor categories}
  For a finite 2-Hilbert space $(\cX,\Tr^\cX)$ and a finite pre-2-Hilbert space $(\cY,\Tr^\cY)$, the category $\Hom(\cX\to\cY)$ of dagger functors is a finite pre-2-Hilbert space with faithful positive trace
  \[
    \Tr_F^{\Hom}(\rho)\coloneq \sum_{s\in\Irr(\cX)}d_s\Tr_{F(s)}^\cY(\rho_s).
  \]
More generally, if $(\cX,\Tr^\cX)$ and $(\cY,\Tr^\cY)$ are finite pre-2-Hilbert spaces, then $\Hom(\cX\to\cY)$ is again a finite pre-2-Hilbert space with trace
  \[
    \Tr_F^{\Hom}(\rho)= \sum_{s\in\Irr(\cX^\cent)}d_s\Tr_{F^\cent(s)}^{\cY^\cent}(\rho^\cent_s),
  \]
  where $F^\cent\colon\cX^\cent\to\cY^\cent$ (resp. $\rho^\cent$) is the unique extension of $F$ (resp. $\rho$) from \cref{3.3.12} and $\Tr^{\cY^\cent}$ is the extension of $\Tr^\cY$ from \cref{trcomp}.
\end{example}

\begin{remark}[{\cite[§2.1]{HPT24}, \cite[Eqns. (4)--(5)]{3Hilb}}: Unitary/Isometric Yoneda embedding] 
  \label{rmk:Unitary Yoneda}
  If $\cX$ is a finite $2$-Hilbert space, then the unitary/isometric Yoneda embedding $\!\yo \colon \cX \to \Hom(\cX^\op \to \Hilb)$ given by $a \mapsto \cX(- \to a)$ is an isometric equivalence of $2$-Hilbert spaces. Thus for any $f \in \cX(a \to a)$,
  \[
    \Tr^\cX_a(f)=\Tr^{\Hom}_{\!\yo(a)}(f \circ -) = \sum_{s \in \Irr(\cX)} d_s \Tr^{\Hilb}_{\cX(s \to a)}(f \circ -)
  \]
  where $d_s\coloneq\Tr^\cX_s(\id_s)$.
\end{remark}

\begin{lemma}
  \label{complete 2-Hilb equiv to its reps} 
  If $\cX$ is a finite pre-2-Hilbert space, then $\Rep(\cX)$ is a finite 2-Hilbert space. Moreover, for any finite 2-Hilbert space $\cX$, there is an isometric equivalence $\cX\cong^\dag\Rep(\cX^\op)$.
\end{lemma}

\begin{pf}
  By \cref{faithful positive trace on functor categories} and \cref{3.3.12}, $\Rep(\cX) \coloneq \Hom(\cX \to \Hilb)\cong^\dag\Hom(\cX^\cent \to \Hilb)$ is a finite 2-Hilbert space. Thus by \cref{rmk:Unitary Yoneda} the unitary/isometric Yoneda embedding $\yo$ gives an isometric equivalence $\cX \cong^\dag \Rep(\cX^\op)$.
\end{pf}

\subsection{Functoriality and the functor-category equivalence}

\subsubsection{Functors and transfors from the disk-like side}

\

\begin{construction}[Functors from unrestricted disk-like functors, \texorpdfstring{$n=1$}{n is 1}]
  \label{dag-thing-i-1}
  Let $\eC$ and $\D$ be disk-like 1-categories. An unrestricted disk-like functor $\eF\colon\eC\to\D$ induces a canonical dagger functor $\cX_\eF\colon\cX_\eC\to\cX_\D$ given on objects $a\in\Obj(\cX_\eC)=\fld[0]{\eC}{\pt}$ by $\cX_\eF(a)\coloneq \eF(a)$ and on morphisms $f\in\cX_\eC(a\to b)=\fldl[1]{\eC}{D^1}[\fcj{a}\amalg b]$ by $\cX_\eF(f)\coloneq \widehat\eF_{D^1,\fcj a\amalg b}(f)$. 
\end{construction}

\begin{lemma}
  \label{dag-thing-i-1-pf}
  \cref{dag-thing-i-1} indeed gives a dagger functor.
\end{lemma}

\begin{pf}
For $f\in\cX_\eC(a\to b)$ and $g\in\cX_\eC(b\to c)$, choose pure representatives $\sigma$ and $\tau$ of $f$ and $g$ respectively and choose parameterized splittings $\Theta\in\eF_{\mathrm{PS}}(\sigma)$ and $\Xi\in\eF_{\mathrm{PS}}(\tau)$. Then
\[
\cX_\eF(g\circ f)=[\eF(\rho_*(\sigma\blt\tau),\rho_*(\Theta\blt\Xi))]
\underset{\ref{1axiom:FTheta-nat}\&\ref{1axiom:FTheta-gluing}}{=}\rho_*([\eF(\sigma,\Theta)]\blt[\eF(\tau,\Xi)])=\cX_\eF(g)\circ\cX_\eF(f),
\]
$\cX_\eF(\id_a)=[\eF(\pi^*a,\varnothing)]\underset{\ref{1axiom:FTheta-products}}{=}\pi^*\eF(a)=\id_{\cX_\eF(a)}$, and
\[
\cX_\eF(f^\dag)=[\eF(\fcj{\iota_*\sigma},\orev{\iota_*\Theta})]\underset{\ref{1axiom:FTheta-nat}\&\ref{1axiom:FTheta-refl}}{=}[\fcj{\iota_*\eF(\sigma,\Theta)}]=\cX_\eF(f)^\dag.\qedhere
\]
\end{pf}

\begin{lemma}
\label{CcongD implies XCcongXD}
An unrestricted isometric weak equivalence of unitary disk-like 1-categories $\eF\colon\eC\to\D$ induces a canonical isometric equivalence of pre-2-Hilbert spaces $\cX_\eF\colon\cX_\eC\to\cX_\D$. Thus an unrestricted isometric weak equivalence $\eC\cong^\dag\D$ induces an isometric equivalence $\cX_\eC\cong^\dag\cX_\D$.
\end{lemma}

\begin{pf}
By \cref{dag-thing-i-1-pf} $\cX_\eF$ is a dagger functor. By \ref{WEurk} at $k=0$, $\cX_\eF$ is essentially surjective, while by \ref{WEurn} it is fully faithful as $\widehat\eF_{D^1,\fcj a\amalg b}\colon\cX_\eC(a\to b)\to\cX_\D(\eF(a)\to\eF(b))$ is a linear isomorphism.

To see $\cX_\eF$ is isometric, let $f\in\cX_\eC(a\to a)$. Choose a pure representative $\sigma_f$ and a parameterized splitting $\Theta\in\eF_{\mathrm{PS}}(\sigma_f)$. As $\fcj{a\times D^1}$ is a product field, $\eF_{\mathrm{PS}}(\fcj{a\times D^1})=\{\varnothing\}$ by \ref{axiom:Gamma-products1}, so \ref{1axiom:FTheta-gluing} and \ref{1axiom:FTheta-products} give $\eF(\cl_a\sigma_f,\varnothing\blt\Theta)=\cl_{\eF(a)}\eF(\sigma_f,\Theta)$, whence
\[
    \Tr^{\cX_\D}_{\eF(a)}(\cX_\eF(f)) 
    \defeq{}
    \psn^\D([\cl_{\eF(a)}\eF(\sigma_f,\Theta)])
    =
    \psn^\D([\eF(\cl_a\sigma_f,\varnothing\blt\Theta)])
    = 
    \psn^\eC([\cl_a\sigma_f])
    \defeq{}
    \Tr^{\cX_\eC}_a(f)
\]
where for the third equality we used that $\eF$ is isometric.
\end{pf}

\begin{definition}[Disk-like \texorpdfstring{$(1,1)$}{(1,1)}-transfor]
\label{def: DL 1-transfor}
Let $\C$ and $\D$ be disk-like $1$-categories and let $Q\in\Disk_1$.
A \defn{$Q$-shaped $(1,1)$-transformation}, also called a \defn{$(1,1)$-transfor} or simply \defn{$1$-transfor}, consists of the following data.
\begin{lst}
\item[(D$\eT$)]\label{1axiom:eT-data}
To each 0-ball $P$ and each 0-field $a\in\fld[0]{\C}{P}$, a pure $1$-field $\eT_P(a)\in\fld[1]{\D}{P\times Q}$.
\end{lst}
The data \ref{1axiom:eT-data} is subject to the following conditions.
\begin{lst}
\item[($\eT\varphi$)]\label{1axiom:eT-nat}
$\eT_{P'}(\varphi_*a) = (\varphi\times\id_Q)_*\eT_P(a)$ for all homeomorphisms $\varphi\colon P\overset\sim\to P'$ and 0-fields $a\in\fld[0]{\C}{P}$.

\item[($\eT\partial$)]\label{1axiom:eT-boundary} 
$\partial(\eT_P(a)) = (\partial\eT)_P(a)$ for all 0-fields $a\in\fld[0]{\C}{P}$, where $(\partial\eT)_P(a)\in\undfld[0]{\D}{P\times\partial Q}$ denotes the pair of values at $a$ of the boundary functors (see \cref{rmk:transfor boundary axiom}).

\item[($\eT\kern.2ex\fcj{\,\cdot\,}$)]\label{1axiom:eT-refl} 
$\eT_{\orev{P}}(\fcj{a}) = \fcj{\eT_P(a)}$ for all 0-fields $a\in\fld[0]{\C}{P}$.

\item[($\eT\kern0.2ex\Box$)]\label{1axiom:eT-naturality}
(\emph{Naturality}.) Let $W$ be a 1-ball and let $\xi\in\fld[1]{\C}{W}$ be a pure 1-field, so that gluing the 1-fields $\eT_{\partial W}(\partial\xi)$ on $\partial W\times Q$ to the values $(\partial\eT)_W(\xi)$ of the boundary functors on $W\times\partial Q$ gives a pure 1-field on the circle $(\partial W\times Q)\cup(W\times\partial Q)$. Then for each labeling $\partial W=\{p_1\}\amalg\{p_2\}$ and $\partial Q=\{q_1\}\amalg\{q_2\}$ of the boundary points, the restrictions $\alpha_i$ of this field to the two 1-balls $X_i\coloneq(\{p_i\}\times Q)\cup(W\times\{q_i\})$ satisfy $\varphi_*[\fcj{\alpha_1}]=[\alpha_2]$ for some homeomorphism $\varphi\colon X_1\to X_2$ fixing $X_1\cap X_2$.
\end{lst}
\end{definition}

Notice that \ref{1axiom:eT-data} and \ref{1axiom:eT-nat} are together equivalent to asserting that $\eT^0\colon\C_0\Rightarrow\fld[1]{\D}{-\times Q}$ is an ordinary natural transformation, while \ref{1axiom:eT-naturality} recovers the usual naturality square (see \cref{dag-thing-i-2-pf}). In the latter case, when $\C=\C^\cX$ is string diagrams for some dagger 1-category $\cX$, one labeling recovers the ordinary naturality square of $\eta$ at $f\coloneq\eval(\xi)$, while the other gives the $\dag$ of the naturality square at $f^\dag$ (see \cref{dag-thing-i-2-pf}); since $f\mapsto f^\dag$ is a bijection, the two conditions are equivalent, so together they impose a single naturality condition, as expected.

\begin{construction}[(1,1)-transfors from \texorpdfstring{$\Hom(\cX\to\cY)$}{Hom}-string diagrams]
  \label{DL 1-functor from 1-functor-2}
  For a $\Hom(\cX{\to}\cY)$-labeled string diagram $\xi$ on a 1-ball $J$, we define a $J$-shaped disk-like (1,1)-transfor $\eC^\xi\in\eFun(\eC^\cX{\to}\eC^\cY)$ as follows.
  \begin{itemize}
    \item For a 0-ball $P$ and a 0-field $a\in\eC^\cX(P)=\Obj(\cX)$, we define $\eC^\xi(a)\coloneq\xi(a)$ where $\xi(a)$ denotes the $\cY$-labeled string diagram obtained from $\xi$ by replacing edge labels $F_e\coloneq\lambda_e$ with $F_e(a)$ and vertex labels $\eta_v\coloneq\lambda_v$ with $(\eta_v)_a$. For example,
          \[
            \xi =
            \begin{tkz}
              \draw[thick, red] (0,0) -- (1,0) node[midway, above=-.5mm, font=\scriptsize] {$F$};
              \draw[thick, blue] (1,0) -- (2,0) node[midway, above=-.5mm, font=\scriptsize] {$G$};
              \draw[thick, purple] (2,0) -- (3,0) node[midway, above=-.5mm, font=\scriptsize] {$H$};
              \node[circle, fill=Green, inner sep=1.5pt, label={[font=\scriptsize, label distance=-.5mm]above:$\eta$}] (v1) at (1,0) {};
              \node[circle, fill=Green, inner sep=1.5pt, label={[font=\scriptsize, label distance=-.5mm]above:$\mu$}] (v2) at (2,0) {};
              \draw[->, thick, Green] (0.7, -0.25) -- (1.3, -0.25);
              \draw[<-, thick, Green] (1.7, -0.25) -- (2.3, -0.25);
            \end{tkz}
            \quad\text{and}\quad a\in\cX
            \qquad\rightsquigarrow\qquad
            \eC^\xi(a)\coloneq \xi(a)=
            \begin{tkz}
              \draw[thick, red] (0,0) -- (1,0) node[midway, above=-.5mm, font=\scriptsize] {$F(a)$};
              \draw[thick, blue] (1,0) -- (2,0) node[midway, above=-.5mm, font=\scriptsize] {$G(a)$};
              \draw[thick, purple] (2,0) -- (3,0) node[midway, above=-.5mm, font=\scriptsize] {$H(a)$};
              \node[circle, fill=Green, inner sep=1.5pt, label={[font=\scriptsize, label distance=-.5mm]above:$\eta_a$}] (v1) at (1,0) {};
              \node[circle, fill=Green, inner sep=1.5pt, label={[font=\scriptsize, label distance=-.5mm]above:$\mu_a$}] (v2) at (2,0) {};
              \draw[->, thick, Green] (0.7, -0.25) -- (1.3, -0.25);
              \draw[<-, thick, Green] (1.7, -0.25) -- (2.3, -0.25);
            \end{tkz}
          \]
\end{itemize}
  The axioms \ref{1axiom:eT-nat}, \ref{1axiom:eT-boundary}, and \ref{1axiom:eT-refl} are immediate. To see naturality \ref{1axiom:eT-naturality}, consider a pure 1-field $\zeta\in\fld[1]{\eC^\cX}{W}[\fcj{a}\amalg b]$ and set $\eta\coloneq\eval(\xi)\colon F_{e_-}\Rightarrow F_{e_+}$ and $f\coloneq\eval(\zeta)$. For one labeling, the two pieces carry the fields $\xi(a)\blt\,\eC^{F_{e_+}}(\zeta)$ and $\eC^{F_{e_-}}(\zeta)\blt\,\xi(b)$, whose evaluations $F_{e_+}(f)\circ\eta_a$ and $\eta_b\circ F_{e_-}(f)$ agree by naturality of $\eta$; thus the two fields are identified by \cref{lem:eval-quotient}. The other labeling follows by applying $\dag$ to the naturality square of $\eta$ at $f^\dag$.
\end{construction}

\subsubsection{Functors and transfors from the traditional side}

\

\begin{construction}[Disk-like functors from functors, \texorpdfstring{$n=1$}{n is 1}]
  \label{DL 1-functor from 1-functor-1}
  A dagger functor $F\colon \cX\to\cY$ induces a canonical disk-like functor $\eC^F\colon\eC^\cX\to\eC^\cY$ sending $\cX$-labeled string diagrams to the same string diagram with $F$ applied to the labels.
\end{construction}

\begin{lemma}
  \label{DL 1-functor from 1-functor-1-pf}
  \cref{DL 1-functor from 1-functor-1} indeed gives a disk-like functor $\eC^F\colon\eC^\cX\to\eC^\cY$.
\end{lemma}

\begin{pf}
  Let $\xi =(\Gamma,\lambda)$ be an $\cX$-labeled string diagram on a 1-ball $X$ with boundary $\fcj{a}\amalg b$. We define $\eC^F(\xi) \coloneq (\Gamma,F\circ\lambda)$ and extend linearly. We verify the disk-like functor axioms:
  \begin{lst}
    \item[\ref{1axiom:eF-nat}] $\eC^F(\varphi_*\xi) = (\varphi(\Gamma), F\circ\lambda\circ \varphi^{-1}) = \varphi_*(\Gamma, F\circ\lambda) = \varphi_*\eC^F(\xi)$.
    \item[\ref{1axiom:eF-boundary}] $\partial(\eC^F(\xi))=\fcj{F(\lambda_{e_0})}\amalg F(\lambda_{e_m})=\eC^F(\partial\xi)$.
    \item[\ref{1axiom:eF-gluing}] $\eC^F(\xi_1\blt_a\xi_2)=(\Gamma_1\cup\Gamma_2,F\circ(\lambda_1\cup\lambda_2))=\eC^F(\xi_1)\blt_{F(a)}\eC^F(\xi_2)$.
    \item[\ref{1axiom:eF-products}] $\eC^F(\pi^* a) = (\varnothing, F(a)) = \pi^* F(a) = \pi^* \eC^F(a)$.
    \item[\ref{1axiom:eF-refl}] $\eC^F(\fcj{\xi}) = (\orev{\Gamma}, F \circ \lambda) = \fcj{(\Gamma, F \circ \lambda)} = \fcj{\eC^F(\xi)}$.
    \item[\ref{1axiom:eF-local-relations}] We show that the linear extension of $\eC^F$ sends each generator of $\lU{\eC^\cX}$ to a generator of $\lU{\eC^\cY}$. Isotopy generators are preserved since $\eC^F$ does not affect $\Gamma$. Adjacent-vertex generators are preserved since $F(\lambda_{v_{i+1}}\circ\lambda_{v_i})=F(\lambda_{v_{i+1}})\circ F(\lambda_{v_i})$. Identity-deletion generators are preserved since $F(\id_x)=\id_{F(x)}$. The $\varepsilon$-$\lambda$-flip generators $(\varepsilon_v,\lambda_v)\sim(-\varepsilon_v,\lambda_v^\dag)$ map to $(\varepsilon_v,F(\lambda_v))\sim(-\varepsilon_v,F(\lambda_v)^\dag)$, since $F(\lambda_v^\dag)=F(\lambda_v)^\dag$. Thus $\eC^F(\fldlU{\eC^\cX}{X}[c])\subset\fldlU{\eC^\cY}{X}[\eC^F(c)]$. \qedhere
  \end{lst}
\end{pf}

\begin{lemma}
  \label{lem:C-preserves-equivs}
If $F\colon(\cX,\Tr^\cX)\to(\cY,\Tr^\cY)$ is an isometric equivalence of pre-2-Hilbert spaces, then $\eC^F\colon\eC^\cX\to\eC^\cY$ is an isometric weak equivalence of unitary disk-like 1-categories.
\end{lemma}

\begin{pf}
\cref{DL 1-functor from 1-functor-1-pf} shows that $\eC^F$ is a disk-like functor.

\itemstep{$\eC^F$ is essentially surjective.}
Let $b\in\Obj(\cY)=\eC^\cY_0(\pt)$. Since $F$ is an equivalence, there is an object $a\in\Obj(\cX)$ and an isomorphism $\eta\in\cY(F(a)\to b)$. The single-vertex diagram $\xi_\eta\in\fldl[1]{\eC^\cY}{I}[\fcj{F(a)}\amalg b]$ and its inverse $\xi_{\eta^{-1}}\in\fldl[1]{\eC^\cY}{I}[\fcj{b}\amalg F(a)]$ witness that $F(a)$ and $b$ are equivalent $0$-fields in $\eC^\cY$.

\itemstep{$\eC^F$ is fully faithful.} 
Fix $a,b\in\Obj(\cX)$. Since $\eC^F$ applies $F$ to every label, $\bar\eval\circ\eC^F\circ\langle-\rangle=F$. As $\bar\eval$ and $\langle-\rangle$ are isomorphisms by \cref{lem:eval-quotient} and $F$ is an isomorphism on hom spaces, $\eC^F$ is an isomorphism on $1$-fields.

\itemstep{$\eC^F$ is isometric.} Let $\xi$ be an $\cX$-labeled string diagram on $S^1$. As in the proof of \cref{lem:functor-commutes-closure}, $[\xi] = \cl_a\langle f\rangle$ for some $a\in\cX$ and $f\in\cX(a\to a)$.
Then $\eC^F(\xi) = \cl_{F(a)}\langle F(f)\rangle$, so
\[
\psn^{\eC^\cY}(\eC^F(\xi))
= \psn^{\eC^\cY}(\cl_{F(a)}\langle F(f)\rangle)
\defeq{} \Tr^\cY_{F(a)}(F(f))
= \Tr^\cX_a(f)
\defeq{} \psn^{\eC^\cX}(\cl_a\langle f\rangle)
= \psn^{\eC^\cX}(\xi)
\]
where the third equality uses that $F$ is isometric.
\end{pf}

\begin{construction}[Natural transformations from (1,1)-transfors]
  \label{dag-thing-i-2}
  An $I$-shaped disk-like 1-transfor $\eN$ with boundary $\fcj{\eF}\amalg\eG$ for disk-like functors $\eF,\eG\colon\eC\to\D$ induces a canonical natural transformation $\cX_\eN\colon\cX_\eF\Rightarrow\cX_\eG$ whose component at an object $a\in\cX_\eC$ is given by $(\cX_\eN)_a \coloneq \eN(a)\in \cX_\D(\cX_\eF(a)\to\cX_\eG(a))$.
\end{construction}

\begin{lemma}
  \label{dag-thing-i-2-pf}
  \cref{dag-thing-i-2} indeed gives a natural transformation.
\end{lemma}

\begin{pf}
  Let $f\in \cX_{\eC}(a\to b)=\fldl[1]{\eC}{I}[\fcj a\amalg b]$ be represented by a pure 1-field $\sigma$.
  The circle field of \ref{1axiom:eT-naturality} at $\sigma$ is obtained by gluing $\eN(a)$, $\eG(\sigma)$, $\eN(b)$, and $\eF(\sigma)$ in cyclic order.
  For the labeling whose two 1-ball components are $X_1$ and $X_2$ with fields $\xi\coloneq \eN(a)\blt \eG(\sigma)\in \fld[1]{\D}{X_1}$ and $\zeta\coloneq \eF(\sigma)\blt \eN(b)\in \fld[1]{\D}{X_2}$, the axiom \ref{1axiom:eT-naturality} gives $\varphi_*[\xi]=[\zeta]$ for some homeomorphism $\varphi\colon X_1\to X_2$, so $\eG(f) \circ \eN(a) = \eN(b) \circ \eF(f)$ in $\cX_\D$.
\end{pf}

\begin{construction}
  \label{cstr:F_eF}
  Given dagger categories $\cX$ and $\cY$ and a disk-like functor $\eF\colon\eC^\cX\to\eC^\cY$, define a dagger functor $F_\eF\colon\cX\to\cY$ on objects $a\in\cX$ by $F_\eF(a)\coloneq\eF(a)$ and on morphisms $f\in\cX(a\to b)$ by $F_\eF(f)\coloneq\eval_I^\cY(\eF\langle f\rangle)$.
\end{construction}

\begin{lemma}
\label{cstr:F_eF-pf}
\cref{cstr:F_eF} indeed gives a well-defined dagger functor $F_\eF$. Moreover, $\eF\cong^\star\eC^{F_\eF}$ in $\eFun(\eC^\cX{\to}\eC^\cY)$.
\end{lemma}

\begin{pf}
First observe that $F_\eF$ is a dagger functor as a composite of dagger functors, namely $F_\eF=\eval\circ\cX_\eF\circ\langle-\rangle$.

\itemstep{$\eF\cong^\star\eC^{F_\eF}$ in $\eFun(\eC^\cX{\to}\eC^\cY)$.}
On $0$-fields, $\eC^{F_\eF}(a)=F_\eF(a)=\eF(a)$ by definition. For a pure $1$-field $\xi\in\fld[1]{\eC^\cX}{W}[\fcj{a}\amalg b]$, set $f\coloneq\eval_W^\cX(\xi)\in\cX(a\to b)$. Then
\(
\eval_W^\cY(\eC^{F_\eF}(\xi))=F_\eF(f)=\eval_W^\cY(\eF\langle f\rangle)\underset{\eqref{lem:eval-quotient}}{=}\eval_W^\cY(\eF(\xi))
\),
where the third equality holds because $[\xi]=\langle f\rangle$ in $\fldl[1]{\eC^\cX}{W}[\fcj{a}\amalg b]$ and $\eF$ preserves local relations. By \cref{lem:eval-quotient}, $[\eC^{F_\eF}(\xi)]=[\eF(\xi)]$ in $\fldl[1]{\eC^\cY}{W}[-]$. Now the claim follows by observing that the $I$-shaped product 1-transfor $\alpha\colon \eC^{F_\eF}\to \eF$ is a unitary equivalence.
\end{pf}

\begin{construction} 
\label{etaNconstruction}
For dagger 1-categories $\cX$ and $\cY$, functors $F,G\colon\cX\to\cY$, and a pure $(1,1)$-transfor $\eN\in\fld[1]{\eFun(\eC^\cX\to\eC^\cY)}{J}[\fcj{\eC^F}\amalg\eC^G]$, define a natural transformation $\eta_\eN\colon F\Rightarrow G$ whose component at $a\in\cX$ is given by
\[
(\eta_\eN)_a\coloneq\eval^\cY_{J,F(a)\to G(a)}(\eN_P(a))\in\cY(F(a)\to G(a)).
\]
\end{construction}

\begin{lemma}
\label{etaNconstruction-pf}
The $\eta_\eN$ constructed in \cref{etaNconstruction} is a natural transformation and $\eC^{\langle\eta_\eN\rangle}=\eN$ as 1-fields in $\fldl{\eFun(\eC^\cX\to\eC^\cY)}{J}[\fcj{\C^F}\amalg \C^G]$.
\end{lemma}

\begin{pf} 
As a whiskering $\eta_\eN\coloneq\eval\circ\cX_\eN\circ\langle-\rangle$ of the natural transformation $\cX_\eN$ from \cref{dag-thing-i-2}, $\eta_\eN$ is a natural transformation. For the second point, we have that $[\eC^{\langle\eta_\eN\rangle}(a)]$ is the quotient class of the single-vertex string diagram with vertex label $\eval(\eN(a))$, so $[\eC^{\langle\eta_\eN\rangle}(a)]=[\eN(a)]$ by \cref{lem:eval-quotient}, whence $[\eC^{\langle\eta_\eN\rangle}]=[\eN]$.
\end{pf}

\subsubsection{The functor category equivalence}
For disk-like 1-categories $\C$ and $\D$, disk-like functors $\C\to\D$ may also be referred to as (1,0)-transfors. As described in \cref{def: functor disk-like n-category}, (1,0)- and (1,1)-transfors assemble into a disk-like 1-category $\eHom(\C{\to}\D)$. 

We next prove a result that will allow us to transfer sphere traces along weak equivalences.
\begin{proposition}
\label{tsf}
If $\eF \colon \eC \to \eD$ is a weak equivalence of disk-like 1-categories and $\psn^\C\colon \undfldl{\C}{S^1} \to \bC$ is a linear functional satisfying \ref{1psnP}, then the following hold.
\begin{lst} 
\item\label{tsfa} $\widehat\eF_{S^1}$ is an isomorphism.
\item\label{tsfb} $\psn^{\D} \coloneq \psn^{\C} \circ \widehat{\eF}_{S^1}^{-1} \colon \undfldl{\D}{S^1} \to \bC$ satisfies \ref{1psnP}.
\item\label{tsfc} $(\eD, \psn^{\eD})$ is a unitary disk-like 1-category, and is finite when $\eC$ is.
\end{lst}
\end{proposition}

\begin{pf} 
To see surjectivity for \ref{tsfa}, first choose for each 0-field $d$ in $\D$ a 0-field $c$ in $\C$ and a unitary equivalence $u\in\fldl{\D}{I}[\fcj{\eF(c)}\amalg d]$. Then for any decomposition $\{\alpha_i\}$ representing a 1-field $\alpha$ in the colimit $\undfldl[1]{\D}{S^1}$, say with $d_i\coloneq\partial\alpha_i$, by extended-isotopy invariance \ref{CTc1} we can also represent $\alpha$ by the decomposition $\{\alpha_i\}\cup\{d_i\times I\}$. Since $u$ is unitary, this is just $\{\alpha_i\}\cup \{{u_i^\star\blt(\eF(c_i)\times I)\blt u_i}\}$. Assuming that we have ordered our indices $i$ so that $\alpha_i$ is glued to $\alpha_{i+1}$ for each $i$ (possibly counting mod $m$ for some integer $m$), by gluing the $u_i$ to the $\alpha_i$ we get another representation of $\alpha$, namely $\{u_{i+1}\bullet\alpha_i\bullet u_{i}^\star\}\cup\{\eF(c_i)\times I\}$, which again by \ref{CTc1} is $\{u_{i+1}\bullet\alpha_i\bullet u_{i}^\star\}$. Since $u_{i+1}\bullet\alpha_i\bullet u_{i}^\star\in\fldl[1]{\D}{I}[\fcj{\eF(c_{i+1})}\amalg \eF(c_i)]$ and $\widehat\eF$ is an isomorphism on 1-balls, we may set $\gamma_i\coloneq \widehat\eF^{-1}(u_{i+1}\bullet\alpha_i\bullet u_{i}^\star)$. Then $\eF([\{\gamma_i\}])=[\{\eF(\gamma_i)\}]=\alpha$, so $\widehat\eF_{S^1}$ is surjective. 

To see $\eF$ is injective, suppose $\xi,\eta\in\undfldl{\C}{S^1}$ have $\widehat\eF_{S^1}(\xi)=\widehat\eF_{S^1}(\eta)$, i.e., $[\{\widehat\eF_{X_i}(\xi_i)\}_{X_i \in \cP}] = [\{\widehat\eF_{X_j}(\eta_j)\}_{X_j \in\cQ}]$
for some $(\cP,\beta_\cP),(\cQ,\beta_\cQ)\in\Disk_\C(S^1)$.
Then there is a refinement $(\cR, \beta_\cR) \in\Disk_\eD(S^1)$ that glues up to both $\{\widehat{\eF}_{X_i}(\xi_i)\}_{X_i \in\cP}$ and $\{\widehat{\eF}_{X_j}(\eta_j)\}_{X_j \in\cQ}$, say via gluing sequences $A$ and $B$ respectively.

By applying the same strategy as in the above proof of surjectivity to the decomposition $\{\widehat{\eF}_{X_k}(\zeta_k)\}_{X_k \in\cR}$, we can take the inverse images of the ball-component fields: set $\gamma_k\coloneq \widehat{\eF}_{X_k}^{-1}(u_{k+1}\blt \alpha_k\blt u_k^\star) \in \fldl{\C}{X_k}[ c_k]$ as above. After doing so, by applying gluing sequences $A$ and $B$ to $\{\gamma_k\}_{X_k \in \cR}$, we have now shown $\{\xi_i\}_{X_i \in \cP}$ and $\{\eta_j\}_{X_j \in \cQ}$ in $\Disk_\C(S^1)$ represent the same element of $\C$ in the colimit, i.e., are equal in $\undfldl[1]{\C}{S^1}$. Hence $\widehat{\eF}_{S^1}$ is injective.

For \ref{tsfc}, observe that \ref{1psnP} holds by \ref{tsfb} and \ref{1psnF} holds by \cref{cor:psnF-weaker} and \ref{tsfa}. Thus it only remains to show \ref{tsfb}.
Let $d=d_1\amalg d_2\in \fld[0]{\D}{\partial I}$ and $\alpha\in\fldl[1]{\D}{I}[d]$. It suffices to show $\bkt\alpha\alpha_{I,d}\geq 0$, with equality if and only if $\alpha=0$. Let $u\in\fldl[1]{\D}{I}[\fcj{\eF(c)}\amalg d]$ be a unitary equivalence in $\D$. Then as quotient-level 1-fields in $\D$ we have
\begin{equation}
\label{star1}
\begin{multlined}
\alpha^\star\blt_d\alpha
=
\begin{tkz}
\draw[orange,thick,mid>](0,0)arc(-90:90:.5)node[pos=.75,right]{\scriptsize$\alpha$};
\draw[orange,thick,mid<](0,0)arc(270:90:.5)node[pos=.75,left]{\scriptsize$\alpha^\star$};
\fill[blue] (0,0) circle (1pt) node[below]{\scriptsize$d_2$};
\fill[blue] (0,1) circle (1pt) node[above]{\scriptsize$d_1$};
\end{tkz}
\underset{\ref{CTc1}}{=}
\begin{tkz}
\draw[orange,thick,mid>](0,0)arc(-90:90:.5)node[pos=.75,right]{\scriptsize$\alpha$};
\draw[orange,thick,mid<](-2.4,0)arc(270:90:.5)node[pos=.75,left]{\scriptsize$\alpha^\star$};
\draw[blue,thick] (0,0) --node[below]{\scriptsize$d_2\times I$}(-2.4,0);
\draw[blue,thick] (0,1)--node[above]{\scriptsize$d_1\times I$}(-2.4,1);
\end{tkz}
=
\begin{tkz}
\draw[orange,thick,mid>](0,0)arc(-90:90:.5)node[pos=.75,right]{\scriptsize$\alpha$};
\draw[orange,thick,mid<](-2.4,0)arc(270:90:.5)node[pos=.75,left]{\scriptsize$\alpha^\star$};
\draw[violet,thick] (0,0) --node[below]{\scriptsize$u_2$}(-.6,0);
\draw[red,thick] (-.6,0) --node[below]{\scriptsize$\eF(c_2){\times} I$}(-1.8,0);
\draw[violet,thick] (-1.8,0) --node[below]{\scriptsize$u_2^\star$}(-2.4,0);
\draw[violet,thick] (0,1) --node[above]{\scriptsize$u_1^\star$}(-.6,1);
\draw[red,thick] (-.6,1) --node[above]{\scriptsize$\eF(c_1){\times} I$}(-1.8,1);
\draw[violet,thick] (-1.8,1) --node[above]{\scriptsize$u_1$}(-2.4,1);
\end{tkz}
\\
\underset{\ref{CTc1}}{=}
\begin{tkz}
\draw[orange,thick,mid>](0,0)arc(-90:90:.5)node[pos=.75,right]{\scriptsize$\alpha$};
\draw[orange,thick,mid<](-2,0)arc(270:90:.5)node[pos=.75,left]{\scriptsize$\alpha^\star$};
\draw[violet,thick] (0,0) --node[below,xshift=1mm]{\scriptsize$u_2$}(-1,0);
\draw[violet,thick] (-1,0) --node[xshift=-1mm,below]{\scriptsize$u_2^\star$}(-2,0);
\fill[red] (-1,0) circle (1pt) node[below]{\scriptsize$\eF(c_2)$};
\draw[violet,thick] (0,1) --node[xshift=1mm,above]{\scriptsize$u_1^\star$}(-1,1);
\draw[violet,thick] (-1,1) --node[xshift=-1mm,above]{\scriptsize$u_1$}(-2,1);
\fill[red] (-1,1) circle (1pt) node[above]{\scriptsize$\eF(c_1)$};
\end{tkz}
=
(u_2\blt\alpha\blt u_1^\star)^\star\blt_{\eF(c)}(u_2\blt\alpha\blt u_1^\star)
\end{multlined}
\end{equation}
where the third equality uses that $u$ is a unitary equivalence.
Set $\beta\coloneq u_2\blt\alpha\blt u_1^\star\in\fldl[1]{\D}{I}[\eF(c)]$ and $\gamma\coloneq\widehat\eF_I^{-1}(\beta)$. Then $\widehat\eF_{S^1}(\gamma^\star\blt_c\gamma)=\beta^\star\blt_{\eF(c)}\beta
\underset{\eqref{star1}}{=}
\alpha^\star\blt_d\alpha$,
so 
\[
\bkt\alpha\alpha_{I,d}^\D\defeq{}\psn^\D(\alpha^\star\blt_d\alpha)\defeq{}\psn^\C(\widehat\eF_{S^1}^{-1}(\alpha^\star\blt_d\alpha))=\psn^\C(\gamma^\star\blt_c\gamma)\defeq{}\bkt\gamma\gamma_{I,c}^\C.
\]
Thus, since $\psn^\C$ satisfies \ref{1psnP}, so does $\psn^\D$.
\end{pf}

\begin{theorem}
\label{DL-thing-i}
For a finite pre-2-Hilbert space $\cX$ and a pre-$2$-Hilbert space $\cY$, the constructions $F\smash{\overset{\eqref{DL 1-functor from 1-functor-1}}\longmapsto}\eC^F$ and $\xi\overset{\eqref{DL 1-functor from 1-functor-2}}\longmapsto\eC^\xi$ assemble into a weak equivalence
\[
  \eC^{(-)}\colon\eC^{\Hom(\cX\to\cY)}\to\eFun(\eC^\cX\to\eC^\cY).
\]
Thus, when $\cY$ is moreover finite, so that $\Hom(\cX\to\cY)$ is a finite pre-2-Hilbert space by \cref{faithful positive trace on functor categories}, \cref{tsf} makes $\eFun(\eC^\cX{\to}\eC^\cY)$ finite unitary with its induced sphere trace.
\end{theorem}

\begin{pf}
Aside from showing \ref{1axiom:eF-local-relations}, we leave the straightforward verification that $\C^{(-)}$ is a disk-like functor to the reader.

\itemstep{\ref{1axiom:eF-local-relations}.}
For $u\in\bC\{\fld[1]{\eC^{\Hom(\cX\to\cY)}}{J}[\fcj F\amalg G]\}$ with $1$-ball $J$ and objects $F,G$ of $\Hom(\cX\to\cY)$,
\begin{align*}
  u\in\fldlU{\eC^{\Hom(\cX\to\cY)}}{J}[\fcj F\sqcup G]
  &\iff\eval^{\Hom(\cX\to\cY)}_{J,F\to G}(u)=0 \\
  &\iff(\eval^{\Hom(\cX\to\cY)}_{J,F\to G}(u))_a=0\quad\forall a\in\cX \\
  &\iff\eval^\cY_{J,F(a)\to G(a)}(u(a))=0\quad\forall a\in\cX \\
  &\iff u(a)\in\fldlU{\eC^\cY}{J}[\fcj{F(a)}\amalg G(a)]\quad\forall a\in\cX \\
  &\overset{\mathrm{def}}{\iff}\eC^u\in\fldlU{\eFun(\eC^\cX\to\eC^\cY)}{J}[\overline{\eC^F}\sqcup\eC^G].
\end{align*}

\itemstep{$\eC^{(-)}$ is unitarily essentially surjective on $0$-fields.}
Given a disk-like functor $\eF\colon\eC^\cX\to\eC^\cY$, the functor $F_\eF\colon\cX\to\cY$ from \cref{cstr:F_eF} satisfies $\eC^{F_\eF}\cong^\star \eF$ by \cref{cstr:F_eF-pf}.

\itemstep{$\eC^{(-)}$ is a linear isomorphism on $1$-fields.}
Faithfulness follows by following the chain of equivalences in \ref{1axiom:eF-local-relations} in reverse, which gives $\eC^u=0\implies u=0$. To see fullness, recall that by \cref{etaNconstruction-pf}, any $\eN\in\fld[1]{\eFun(\eC^\cX\to\eC^\cY)}{J}$ equals $\eC^{\langle\eta_{\eN}\rangle}$ as (quotient-level) 1-fields in $\eFun(\eC^\cX\to\eC^\cY)$, where $\eta_{\eN}$ is from \cref{etaNconstruction}.

\itemstep{$\eC^{(-)}$ is isometric.}
This holds by definition of the sphere trace on $\eFun$: $\psn^{\eFun}\coloneq\psn^{\eC^{\Hom(\cX\to\cY)}}$ via \cref{tsf}\ref{tsfb}. 
\end{pf}

\begin{corollary}
\label{dag-thing-i-equiv}
\label{prop:XFun}
For finite pre-2-Hilbert spaces $\cX$ and $\cY$, there is a canonical isometric equivalence of finite pre-2-Hilbert spaces
\[
    \cX_{\eFun(\eC^\cX\to\eC^\cY)} \overset\sim\longrightarrow \Hom(\cX\to\cY).
\]
\end{corollary}

\begin{pf}
Composing the isometric equivalence $\cX_{\eC^{(-)}}$ with the inverse of $\eval$ gives the claim:
    \[
    \cX_{\eFun(\eC^\cX\to\eC^\cY)}
    \xleftarrow[\substack{\cX_{\eC^{(-)}}\\\eqref{DL-thing-i}\,\&\,\eqref{CcongD implies XCcongXD}}]{\cong^\dag}\cX_{\eC^{\Hom(\cX\to\cY)}}
    \xrightarrow[\substack{\eval\\\eqref{prop:XCXcongX}}]{\cong^\dag}
    \Hom(\cX\to\cY).\qedhere
    \]
\end{pf}

This proves the $n=1$ half of \cref{HOMEQUIV}.

\subsection{The \texorpdfstring{$(1+1)$D}{(1+1)D} Unitary Cobordism Hypothesis}
Let $(\eC,\psn)$ be a finite unitary disk-like 1-category. A \defn{unitary representation} of $\eC$ is a disk-like functor $\eF\colon\eC^{\cX_\C^\op}\to\eHilb$, where $\eHilb\coloneq\eC^\Hilb$.

\begin{remark}
\label{rem:eRep-strings}
Representations are taken of $\eC^{\cX_\eC^\op}$ rather than of $\eC$ because the functor category $\eFun(\eC{\to}\eHilb)$ is not invariant under isometric weak equivalences whose zig-zag legs are unrestricted. Indeed, let $\eC_T$ be the disk-like 1-category of maps into a fixed space $T$ and let $\cX_T\coloneq\cX_{\C_T}$. Every pure 1-field in $\C_T$ splits at every point in the interior, so by \ref{1axiom:eF-gluing} and injectivity \ref{CGi1} each disk-like functor $\eF\colon\eC_T\to\eHilb$ must send all pure 1-fields to a vertex-free $\Hilb$-string diagram. But by \ref{1axiom:eF-boundary} such an $\eF$ is constant, and $\cX_{\eFun(\eC_T{\to}\eHilb)}\cong^\dag\Hilb$. On the other hand $\eC_T\cong^\dag\eC^{\cX_T}$ by \cref{prop:CXCcongC}, while $\cX_{\eFun(\eC^{\cX_T^\op}{\to}\eHilb)}\cong^\dag\Rep(\cX_T^\op)$ by \cref{prop:XFun}, which when $T=\rB G$ for a finite group $G$ is just $\Rep(G)$. 

For a string diagram disk-like category $\C$, however, the two candidates for $\eRep(\C)$ agree: for any finite pre-2-Hilbert space $\cX$, there is a canonical zig-zag of isometric weak equivalences
\[
\eRep(\eC^\cX)
\underset{\eqref{DL-thing-i}}{\cong^\dag}
\eC^{\Hom(\cX_{\eC^\cX}^\op\to\Hilb)}
\underset{\eqref{lem:C-preserves-equivs}\&\eqref{prop:XCXcongX}}{\cong^\dag}
\eC^{\Hom(\cX^\op\to\Hilb)}
\underset{\eqref{DL-thing-i}}{\cong^\dag}
\eFun(\eC^{\cX^\op}{\to}\eHilb).
\]
\end{remark}

Unitary representations of $\eC$ form a disk-like 1-category $\eRep(\eC)\coloneq\eFun(\eC^{\cX_\eC^\op}{\to}\eHilb)$. For any finite pre-2-Hilbert space $\cX$, by \cref{DL-thing-i} and \cref{tsf} the sphere trace $\psn^{\C^{\Rep(\cX^\op)}}$ obtained from the faithful positive trace on $\Rep(\cX^\op)=\Hom(\cX^\op\to\Hilb)$ defined in \cref{faithful positive trace on functor categories} equips $\eFun(\eC^{\cX^\op}{\to}\eHilb)\cong^\dag\eRep(\eC^\cX)$ (see \cref{rem:eRep-strings})
with the structure of a unitary disk-like 1-category. By \cref{prop:CXCcongC} and \cref{finite cylinders}, this equips $\eRep(\eC)$ with a sphere trace for any finite unitary disk-like 1-category $\eC$. 

\begin{definition}
\label{def:completeDL1cat}
    A finite unitary disk-like 1-category $\eC$ is called \defn{complete} if there is an isometric weak equivalence $\eC\cong^\dag\eRep(\eC)$.
\end{definition}

\TWODUCH*

\begin{pf}
The isometric equivalence $\cX_{\eC^\cX}\cong^\dag\cX$ and the isometric weak equivalence $\eC^{\cX_\eC}\cong^\dag\eC$ are \cref{prop:XCXcongX} and \cref{prop:CXCcongC} respectively. Thus it only remains to show these equivalences preserve completeness in that (i) $\cX_{\eRep(\eC)}$ is complete for all finite unitary disk-like 1-categories $\eC$ and (ii) $\eC^\cX$ is complete for all finite 2-Hilbert spaces $\cX$. 

For (i), observe that 
\[
\cX_{\eRep(\eC)}\defeq{}\cX_{\eFun(\eC^{\cX_\eC^\op}\to\eC^{\Hilb})} \;\underset{\mathclap{\eqref{prop:XFun}}}{\,\cong^\dag}\; \Hom(\cX_\eC^\op\to\Hilb)\defeq{}\Rep(\cX_\eC^\op).
\]
As $\cX_\eC$ is finite unitary by \cref{finite cylinders}, it follows from \cref{complete 2-Hilb equiv to its reps} that $\cX_{\eRep(\eC)}$ is complete.

For (ii), we have
\[
\eRep(\eC^\cX)\underset{\eqref{rem:eRep-strings}}{\,\cong^\dag}\eFun(\eC^{\cX^\op}{\to}\eHilb)\underset{\eqref{DL-thing-i}}{\,\cong^\dag}\eC^{\Rep(\cX^\op)}\underset{\eqref{complete 2-Hilb equiv to its reps}\,\&\,\eqref{lem:C-preserves-equivs}}{\,\cong^\dag}\eC^\cX,
\]
where the third equivalence applies $\eC^{(-)}$ to the isometric equivalence $\!\yo^{\!-1}\colon\Rep(\cX^\op)\overset\sim\to\cX$ from \cref{complete 2-Hilb equiv to its reps}. Thus $\eC^\cX$ is complete. 
\end{pf}

\subsection{Application: Classification of (oriented) \texorpdfstring{$(1+1)$D}{(1+1)D} unitary TQFTs as \texorpdfstring{$\Hstar$}{H*}-Morita classes of \texorpdfstring{$\Hstar$}{H*}-algebras}
\label{classificationHstaralg}

In this subsection, we show how the above results recover the unitary version of the classification in Schommer-Pries' thesis \cite{SP09} of fully extended oriented 2D TQFTs as standard separable Frobenius algebras.

\begin{definition}[$\Hstar$-algebra, {\cite{categoricalquantum}}, {\cite[§I.3.6]{UQSL}}]
\label{def:Hstaralg}
  An \defn{$\Hstar$-algebra} is a finite-dimensional $\Cstar$-algebra $A$ equipped with a faithful positive trace $\Tr_A \colon A \to \mathbb{C}$.
\end{definition}

By \cite[§I.3.6]{UQSL}, standard separable unitary Frobenius algebras are precisely $\Hstar$-algebras, and by \cite{3Hilb} every $\Hstar$-algebra $A$ gives a pre-2-Hilbert space $\rB A$, called the \defn{delooping} of $A$, which has exactly one object $*$ and $\End_{\rB A}(*)\coloneq A$.

Denote by $\Mod(A)$ the category of right Hilbert space modules over an $\Hstar$-algebra $A$.

\begin{definition}[Commutant trace]
\label{H*AlgComm}
Let $(B, \Tr_B)$ be an $\Hstar$-algebra, and let $H_B$ be a right Hilbert space $B$-module. By \cite[Defn. 3.4.10]{UQSL}, the space of endomorphisms $\End(H_B)$ has a canonical \defn{commutant trace} $\Tr_B'$, which is the unique faithful positive trace satisfying
\[
\Tr_B'(\ket{\xi}_B\bra{\eta}) = \bkt{\eta}{\xi}_H \qquad \forall \eta, \xi \in H_B.
\]
\end{definition}

\begin{definition}
\label{H*algisomeq}
A $*$-algebra homomorphism $\alpha\colon A\to B$ between $\Hstar$-algebras is called \defn{isometric} if $\Tr_B(\alpha(a))=\Tr_A(a)$ for all $a\in A$.
\end{definition}

\begin{definition}
\label{def:HstarMorita}
We say that two $\Hstar$-algebras $(A, \Tr_A)$ and $(B, \Tr_B)$ are \defn{$\Hstar$-Morita equivalent} \cite[Exm. 3.7.9]{UQSL}, or \defn{isometrically Morita equivalent}, if there is a faithful right Hilbert space $B$-module $H_B$ and a unitary algebra isomorphism $\alpha \colon A \overset\sim\to \End(H_B)$ such that $\Tr_A = \Tr_B' \circ \alpha$.
\end{definition}

\begin{remark} 
For $\Hstar$-algebras $(A, \Tr_A)$ and $(B, \Tr_B)$, an $\Hstar$-Morita equivalence is precisely an isometric equivalence between these objects in the 3-Hilbert space $\mathsf{H^*Alg}$ in the sense of \cref{def:isometric-equivalence}. 
\end{remark}

\begin{proposition}\label{lem:hstar-morita-rep}
\label{precprop}
For $\Hstar$-algebras $A$ and $B$,
\[
\eRep(\C^{\rB A}) \cong^\dag \eRep(\C^{\rB B})
\quad\iff\quad
A \text{ and } B \text{ are } \Hstar\text{-Morita equivalent}.
\]
\end{proposition}

\begin{pf}
By \cref{rem:eRep-strings} and \cref{DL-thing-i},
\[
\eRep(\C^{\rB A})
\cong^{\dag}
\C^{\Rep((\rB A)^\op)}
\underset{\eqref{lem:C-preserves-equivs}}{\,\cong^\dag}
\C^{\Mod(A)},
\]
and similarly for $B$, so by
\cref{CcongD implies XCcongXD,lem:C-preserves-equivs} it suffices to show
\[
\Mod(A)\cong^{\dag} \Mod(B)
\quad
\iff
\quad
A \text{ and } B \text{ are } \Hstar\text{-Morita equivalent.}
\]

\itemstep{$(\Rightarrow)$} Suppose $F\colon \Mod(A) \to \Mod(B)$ is an isometric equivalence. Because $A$ is semisimple, by \cite[\S I.2.2]{UQSL} $L^2A_A$ has each simple $A$-module as a summand, so $L^2A_A$ is a generator of $\Mod(A)$. Thus $H_B \coloneq F(L^2A_A)$ is a generator of $\Mod(B)$, and so is faithful. As $F$ is isometric, it restricts to a trace-preserving $\Hstar$-algebra isomorphism $\alpha'_F \colon \End(L^2A_A) \overset\sim\to\End(H_B)$. Precomposing with the left-multiplication map $L\colon A \overset\sim\to \End(L^2A_A)$ then gives the desired trace-preserving $\Hstar$-algebra isomorphism $\alpha_F \coloneq \alpha'_F \circ L \colon A \to \End(H_B)$. Indeed, for all $a \in A$, since $L_a = \ket{a}_A\bra{1_A} \in \End(L^2A_A)$, we have
\[
\Tr_{\End(L^2A)}(L_a)
= \Tr_{\End(L^2A)}(\ket{a}_A\bra{1_A})
= \Tr_A(\bkt{1_A}{a}_A)
= \Tr_A(a).
\]
Thus $A$ and $B$ are $\Hstar$-Morita equivalent.

\itemstep{$(\Leftarrow)$} Suppose $A$ and $B$ are $\Hstar$-Morita equivalent, say via the faithful right Hilbert space $B$-module $H_B$ and the trace-preserving $\Hstar$-algebra isomorphism $\alpha\colon A \to \End(H_B)$. We seek an isometric equivalence $F\colon \Mod(A) \to \Mod(B)$. Define $F \colon \Mod(A) \to \Mod(B)$ by $F(K_A)\coloneq K \boxtimes_A H_B$, where the left $A$-action on $H$ is by $a \triangleright \xi \coloneq \alpha(a)(\xi)$, and on morphisms $f\colon K_A \to P_A$ by $F(f) \coloneq f \boxtimes_A \mathrm{id}_{H_B}$. Then $F$ is dagger since
\[
F(f^\dag) = f^\dag \boxtimes_A \id_{H_B}
\underset{\text{\cite[Exr. I.3.2.34]{UQSL}}}{=}
(f \boxtimes_A \mathrm{id}_{H_B})^\dag
= F(f)^\dag.
\]
Since $H_B$ is faithful, it now follows from \cite[Cor. I.3.3.12]{UQSL} that $F$ is an equivalence. It only remains to show $\Tr^{\Mod(B)}_{F(K_A)}(F(f)) = \Tr^{\Mod(A)}_{K_A}(f)$ for all $K_A \in \Mod(A)$ and $f \in \End(K_A)$.
As $\End(K_A)$ is spanned by rank-1 operators, it suffices to show this for $f = \ket{\xi}_A\bra{\eta}$. Let $u \colon F(L^2A_A) = L^2A \boxtimes_A H_B \to H_B$ denote the unitary isomorphism
$u(a \boxtimes \xi) \coloneq a \triangleright \xi = \alpha(a)(\xi)$. Observe that $u F(L_a) u^\dag = \alpha(a)$. Using this fact together with functoriality of $F$ and cyclicity of $\Tr^{\Mod(B)}$ and $\Tr^{\Mod(A)}$, we compute that
\begin{multline*}
\Tr_{F(K_A)}^{\Mod(B)} (F(\ket{\xi}_A \bra{\eta}))
= \Tr_{F(L^2A_A)}^{\Mod(B)} (F(\bra{\eta} \circ \ket{\xi}_A))
= \Tr_{F(L^2A_A)}^{\Mod(B)} (F(L_{\bkt{\eta}{\xi}_A})) \\
= \Tr_{H_B}^{\Mod(B)} (u F(L_{\bkt{\eta}{\xi}_A}) u^\dag)
= \Tr_{H_B}^{\Mod(B)} (\alpha(\bkt{\eta}{\xi}_A))
= \Tr_{L^2A_A}^{\Mod(A)} (L_{\bkt{\eta}{\xi}_A})
= \Tr_{K_A}^{\Mod(A)} ( \ket{\xi}_A\bra\eta).
\end{multline*}
Thus $F$ is isometric.
\end{pf}

\TWODCLASS*

\begin{pf}
  By \cref{TWODUCH}, classifying complete finite unitary disk-like $1$-categories up to isometric weak equivalence reduces to classifying finite $2$-Hilbert spaces up to isometric equivalence. Let $\cX$ be a finite $2$-Hilbert space. As $\cX$ is finite semisimple, it admits a generator $a=\bigoplus_{s\in\Irr(\cX)}s \in \cX$. By \cite[Exm. 2.14]{3Hilb}, the unitary algebra $A \coloneq \End_{\cX}(a)$ equipped with the faithful positive trace $\Tr_A \coloneq \Tr^{\cX}_a$ is an $\Hstar$-algebra, and $\cX \cong^\dag \Mod(A)$ isometrically as $2$-Hilbert spaces. Thus, by \cref{precprop}, we are done.
\end{pf}

\begin{remark}
\label{rmk:n1failure}
  If two $\Hstar$-algebras are merely von Neumann Morita equivalent in the sense of \cite[Defn. I.3.3.9]{UQSL}, then they need not give the same unitary TQFT, since their $\rO(2)$-homotopy fixed point structures need not be isomorphic. Indeed, the two $\Hstar$-algebras $(\bC,1)$ and $(\bC,\lambda)$ with $1\neq\lambda\in\bR_{>0}$ are trivially von Neumann Morita equivalent, but are not isometrically Morita equivalent by \cref{precprop}: if $F\colon(\Hilb,\lambda\Tr)\cong^\dag\Mod(\bC,\lambda)\to\Mod(\bC,1)\cong^\dag(\Hilb,\Tr)$ is an isometric equivalence of 2-Hilbert spaces, then
  \[
    \lambda=\lambda\Tr(\id_\bC)=\Tr(\lambda\id_\bC)=\Tr(F(\id_\bC))=\Tr(\id_{F(\bC)})=\dim F(\bC)=\dim\bC=1.
  \]
\end{remark}

\section{\texorpdfstring{$(2+1)$D}{(2+1)D}}
\label{sec:2+1D}

\subsection{Disk-like 2-categories}
In this section, all manifolds are assumed to be PL, compact, and equipped with the germ of a thickening to an oriented 3-manifold.

\subsubsection{The definition of a disk-like 2-category}
\label{sec:def-dl2cat}
Here we give a relatively self-contained definition of a disk-like 2-category. For the full details, see \cite[\S6.1]{MW12}. As we did for $n=1$, we will include the notion of reflection structure in the definition as extra data.

\itemstep{Dimensions $k=0$ and $k=1$.}
The data and conditions for 0-fields and 1-fields are exactly as in the definition of a disk-like 1-category in \cref{sec:1+1D}, with the following modifications. First, we will typically denote 1-manifolds by $Y$ and variants (e.g., $Y'$ and $Y_i$), so that, as in the $n=1$ case, we reserve $X$ and its variants for top-dimensional manifolds. We use $W$ and its variants when the dimension is irrelevant or arbitrary.
Second, the datum \ref{DU1} (local relations) and conditions \ref{CU1}, \ref{CTi1}, and \ref{CTc1} are omitted, as these are top-dimensional ($k=n$) conditions only. The colimit construction from round $k=1$ results in a functor $\undC_1\colon\Mfld_1\to\Set$ as described in \cref{sec:1+1D}.

Given a splitting $\partial X=Y_1\cup_E Y_2$ of the boundary of a 2-ball $X$ along a 0-manifold $E$, the colimit construction for 1-manifolds gives an injective gluing map $\undfld[1]{\undC}{Y_1}\times_{\undfld[0]{\undC}{E}}\undfld[1]{\undC}{Y_2}\hookrightarrow\undfld[1]{\undC}{\partial X}$; we write $\undfld[1]{\undC}{\partial X}_{\pitchfork E}$ for its image and set $\fld[2]{\C}{X}_{\pitchfork E}\coloneq\partial_2^{-1}(\undfld[1]{\undC}{\partial X}_{\pitchfork E})$ where $\partial_2$ is defined in \ref{2Dbdy2} below. Elements of $\fld[2]{\C}{X}_{\pitchfork E}$ are called \defn{splittable along $E$}. This condition was vacuous in dimension $k=1$, where the gluing locus was a point (which has empty boundary).

\itemstep{2-fields.}
We will refer to an oriented 2-manifold admitting some orientation-preserving homeomorphism to the standard 2-ball (bigon) $D^2$ as a \defn{2-ball}. Let $X$ be a 2-ball.
\begin{lst}
\item[(D$\C_2$)]\label{2DC2} (\emph{Pure 2-fields}) A set $\fld[2]{\C}{X}$ of \defn{pure 2-fields}.
\item[(D$\varphi_2$)]\label{2Dvarphi2} (\emph{Homeomorphism actions}) A bijection $\varphi_*\colon\fld[2]{\C}{X}\to\fld[2]{\C}{X'}$ for each homeomorphism $\varphi\colon X\to X'$.
\item[(D$\partial_2$)]\label{2Dbdy2} (\emph{Boundary maps}) A map $\partial_2\colon\fld[2]{\C}{X}\to\undfld[1]{\undC}{\partial X}$. Write $\fld[2]{\C}{X}[c]$ for the preimage $\partial_2^{-1}(c)\subset\fld[2]{\C}{X}$.
\item[(DG$_2$)]\label{2DG2} (\emph{Gluing}) For each splitting $X=X_1\cup_Y X_2$ along a 1-ball $Y$, a gluing map $\glu_Y\colon\fld[2]{\C}{X_1}_{\pitchfork(\partial Y)}\times_{\fld[1]{\C}{Y}}\fld[2]{\C}{X_2}_{\pitchfork(\partial Y)}\to\fld[2]{\C}{X}_{\pitchfork(\partial Y)}$, where the fiber product is over the common restriction to $\fld[1]{\C}{Y}$.
\item[(D$\pi_2$)]\label{2Dpi2} (\emph{Pullbacks}) For each pinched product map $\pi\colon X\to W$ with $\dim X=2$ and $\dim W<2$, a set map $\pi^*\colon\fld[\dim W]{\C}{W}_{\pitchfork\theta_\pi}\to\fld[2]{\C}{X}$, where $\theta_\pi$ denotes the stratification of $\partial W$ by fiber dimension. The pure 2-fields in its image are called \defn{product 2-fields}. 
\item[(D$\fcj{\,\cdot\,}_2$)]\label{2Dfcj2} (\emph{Reflection}) An involutive bijection $\fcj{\,\cdot\,}\colon\fld[2]{\C}{X}\to\fld[2]{\C}{\orev X}$.
\item[(D$U$)]\label{2DU2} (\emph{Local relations}) For each 2-ball $X$ and $c\in\undfld[1]{\undC}{\partial X}$, a subspace $\fldlU{\C}{X}[c]\subset\bC\{\fld[2]{\C}{X}[c]\}$ of so-called \defn{local relations}.
\end{lst}
The above data are subject to the following conditions.
\begin{lst}
\item[(C$\varphi$)]\label{2Cvarphi2} $\C_2$ is a functor $\Disk_2\to\Set$.
\item[(C$\partial$)]\label{2Cbdy2} $\partial_2$ is a natural transformation $\C_2\Rightarrow\undC_1\circ\partial$.
\item[(CGi)]\label{2CGi2} Every gluing map $\glu_Y$ is injective.
\item[(CGa)]\label{2CGa2} Gluing is associative in that any two ways of starting with a collection of 2-balls and gluing them up sequentially to form another 2-ball induce equal gluing maps.
\item[(CG)]\label{2CG2} A homeomorphism respecting a splitting sends the glued field to the gluing of its image. The boundary of a glued field restricts correctly to each piece's outer boundary.
\item[(C$\pi$)]\label{2Cpi2} Product fields are compatible with gluing, homeomorphisms, boundary maps, restriction to sub-balls, and composition of pinched product maps.
\item[(C$\fcj{\,\cdot\,}$)]\label{2Cfcj2} Reflection commutes with homeomorphisms, boundary maps, gluing, and products.
\item[(CS)]\label{2CS2} \emph{Splittability.} Every pure 2-field $\xi$ on a 2-ball $X$ has a string diagram stratification $S_\xi\subset X$ such that $\xi$ is splittable along every splitting of $X$ transverse to $S_\xi$. Call any such $S_\xi$ a \defn{string locus} of $\xi$.
\item[(C$U$)]\label{2CU2} Local relations are preserved under homeomorphisms and under the conjugate-linear extension of reflection, and form an ideal in the sense that gluing a local relation to anything results in a local relation.
\item[(C$U$i)]\label{2CTi2} \emph{Isotopy invariance.} If all components of an isotopy $h$ on a 2-ball fix the boundary of a pure 2-field $\xi$, then $\xi$ and its image under $h$ differ by a local relation.
\item[(C$U$c)]\label{CTc2} \emph{Collaring invariance.} For a 2-ball $X$, attaching a collar $X\mapsto X\cup_Y(Y\times J)$ along a 1-ball $Y\subset\partial X$ and collapsing it back via a collaring homeomorphism $\psi_{Y,J}\colon X\cup_Y(Y\times J)\to X$ changes a pure 2-field only by a local relation: $(\psi_{Y,J})_*(\xi\blt_Y(\res_Y(\xi)\times J))-\xi\in\fldlU{\C}{X}[c]$ for any $\xi\in\fld[2]{\C}{X}[c]$.
\end{lst}
As before, \ref{2CTi2} and \ref{CTc2} are together known as extended-isotopy invariance.

\itemstep{Extending to all 2-manifolds.} Let $X$ be a 2-manifold and let $\Disk(X,c;\C)$ be the poset of pairs $(\cP; \beta_\cP)$, where $\cP$ is a splitting of $X$ into 2-balls and $\beta_\cP$ consists of compatible boundary labels for each constituent 2-ball in $\cP$ that restrict to $c$ on $\partial X$. As it was for $n=1$, the poset ordering is by antirefinement. We define
\[
\undfld[2]{\undC}{X}[c] \coloneq \colim \Pi_{X;c}^{\C}
\]
where $\Pi_{X;c}^{\C} \colon \Disk(X,c;\C) \to \Set$ is the functor given by $\Pi_{X;c}^{\C}(\cP; \beta_\cP) \coloneq \prod_{X_i \in \cP} \fld[2]{\C}{X_i}[ \beta_\cP|_{\partial X_i}]$ and on morphisms by the corresponding iterated gluing map. For a 2-ball $X$, the trivial decomposition $\{X\}$ is terminal, so $\undfld[2]{\undC}{X}[c]=\fld[2]{\C}{X}[c]$. Thus we will drop the arrow decoration on $\undC$ in favor of $\C$.

\itemstep{Linearization.}
Set $\undfldl[2]{\C}{X}[c]\coloneq\bC\{\undfld[2]{\C}{X}[c]\}\big/\fldlU{\C}{X}[c]$.
This defines a $\Vec$-valued functor on the category of pairs $(X,c)$ of 2-manifolds $X$ and boundary conditions $c$ whose morphisms $(X,c)\to(X',c')$ are homeomorphisms $\varphi\colon X\to X'$ with $\varphi_*c=c'$; by the above axioms, gluing, boundary maps, and reflection descend to this quotient.

As for $n=1$, the above colimit is the same regardless of whether we linearize before or after taking the colimit. That is,
\[
\undfldl[2]{\C}{X}[c] = \bC\{\colim\Pi_{X;c}^\C\}\Big/{\fldlU{\C}{X}[c]}=\colim\tld\Pi_{X;c}^{\C}
\]
where $\tld\Pi_{X;c}^{\C}\colon\Disk(X,c;\C)\to\Vec$ is given by $\tld\Pi_{X;c}^{\C}(\cP;\beta_\cP) \coloneq \bigotimes_{X_i\in\cP}\fldl[2]{\C}{X_i}[\beta_\cP|_{\partial X_i}]$ and on morphisms by the maps induced by the linearized gluing maps on the tensor products of the vector spaces of 2-fields on the constituent balls of the decomposition. Exactly as it was defined for $n=1$, here $\fldlU{\C}{X}[c] \subset \bC\{\undfld[2]{\C}{X}[c]\}$ is the collection of all finite sums of the form $\sum_\gamma \lambda_\gamma [\eta_1, \dots, \eta_{j-1}, \xi_j^\gamma, \eta_{j+1}, \dots, \eta_N]$ indexed over ball splittings $(\cP = \{X_i\}_{i=1}^N,\beta)\in\Disk(X,c;\C)$, $1 \leq j \leq N$, local relations $\sum_\gamma \lambda_\gamma \xi_j^\gamma \in \fldlU{\C}{X_j}[ \beta|_{\partial X_j}]$, and compatible pure 2-fields $\eta_i \in \fld[2]{\C}{X_i}[ \beta|_{\partial X_i}]$ for $i \neq j$. The square bracket $[-]$ denotes the image in the (set-valued) colimit $\undfld[2]{\C}{X}[c]$, not to be confused with the square bracket $[-]$ used to denote the (quotient-level) 2-field $[\alpha]$ represented by a pure 2-field $\alpha$.

\begin{example}
\label{disk-like skein 1-category}
Let $\C$ be a disk-like 2-category. For any 1-manifold $Y$ and any $\eta\in\undfld[0]{\C}{\partial Y}$, there is a disk-like 1-category $\A_\C(Y,\eta)$ whose 0-fields on a 0-ball $P$ are 1-fields $\xi\in\fld{\C}{Y\times P}[\eta]$ and whose pure 1-fields on a 1-ball $J$ are pure 2-fields $\alpha\in\undfld[2]{\C}{Y\times J}[\eta\times J]$. By the boundary condition $\eta\times J$ we mean that pure 2-fields $\xi\in\fld[1]{\A_\C(Y,\eta)}{J}$ can have any boundary 0-field so long as $\xi|_{\partial Y\times J}=\eta$. (See \cref{disk-like skein k-categories} for the general case.) One could reasonably call this the \defn{disk-like skein 1-category} of $\C$ on $Y$ with boundary condition $\eta$, or alternatively the \defn{dimensional reduction} of $\C$ along $(Y,\eta)$. 
\end{example}

\subsubsection{Unitary disk-like 2-categories}
\label{sec:unitary disk-like 2-categories}

\

\begin{definition}
We will call a disk-like 2-category $\C$ \defn{finite} if
\begin{lst}
    \item[(F)] 
    \label{2psnF}
    (\emph{Finiteness})
    $\dim\undfldl[2]{\C}{X}[c]<\infty$ for all 2-manifolds $X$ and $c\in\undfld[1]{\C}{\partial X}$.
\end{lst}
\end{definition}

Note that \ref{2psnF} is a condition regarding general 2-manifolds, not just 2-balls.

\begin{definition}[Unitary disk-like 2-category]
  \label{def:unitary disk-like 2-category}
  A \defn{unitary disk-like 2-category} $(\C,\psn)$ is a disk-like 2-category $\C$ equipped with a \defn{sphere trace} (or \defn{trace}), that is, a linear functional $\psn\colon\undfldl[2]{\C}{S^2}\to\bC$ satisfying the following condition.
  \begin{lst}[font=\upshape]
    \item[($\psn$P)]\label{2psnP}
    (\emph{Positivity})
    For each 2-ball $X$ and each $c\in \undfld[1]{\C}{\partial X}$, the sesquilinear pairing
    \begin{equation}
      \begin{aligned}
        \orev{\fldl[2]{\C}{X}[c]}\otimes_{\bC}\fldl[2]{\C}{X}[c] & \longrightarrow \bC,
        \\
        f\otimes g                      & \longmapsto \bkt{f}{g}_{X,c}\coloneq\psn(\fcj f\blt_c g)
      \end{aligned}
    \end{equation}
    is positive-definite.
  \end{lst}
\end{definition}

\begin{example}
\label{disk-like skein 1-categories are unitary}
For a finite unitary disk-like 2-category $\C$, a 1-manifold $Y$, and $\eta\in\undfld[0]{\C}{\partial Y}$, the disk-like skein 1-category $\A_\C(Y,\eta)$ defined in \cref{disk-like skein 1-category} is finite unitary when equipped with the canonical sphere trace 
\[
\psn^{\A_\C(Y,\eta)}(\xi)
\coloneq
\Z_\C(Y\times D^{2})(\xi).
\]
Indeed, for a 1-ball $J$ the pairing induced by $\psn^{\A_\C(Y,\eta)}$ on $\fldl[1]{\A_\C(Y,\eta)}{J}[c]=\undfldl[2]{\C}{Y\times J}[c]$ is the path-integral pairing $\bkt--_{Y\times J,\,c}$, which is positive-definite by \cref{UnitaryWalkerTheorem}. Finiteness of $\A_\C(Y,\eta)$ follows from finiteness of $\C$.
\end{example}

\subsubsection{(Co)isometries and unitary equivalences of fields}
\label{sec:coisometries-2}

For a 1-manifold $Q$ and a 1-ball $Y$, consider a 2-field $\alpha\in\fld[2]{\C}{Q\times Y}[\fcj\xi\cup\eta]$ $=\fldl{\A_\C(Q,\fcj{a}\amalg b)}{Y}[\fcj{\xi}\amalg \eta]$, where $\xi,\eta\in\fld[1]{\C}{Q}[\fcj{a}\amalg b]$. We call $\alpha$ an \defn{isometry} (resp. \defn{coisometry}) in $\C$ if $\alpha$ is an isometry (resp. coisometry) in the sense of \cref{def:cosiom1} when viewed as a 1-field on $Y$ in $\A_\C(Q,\fcj{a}\amalg b)$. We call $\alpha$ \defn{unitary} or a \defn{unitary equivalence} if $\alpha$ is both an isometry and a coisometry. We write $\xi\cong^\star \eta$ (in $\C$) if there is a unitary equivalence $\alpha$ between $\xi$ and $\eta$ in $\C$.

\begin{definition} \label{2coisometryunitarydl}
Let $\C$ be a disk-like 2-category, let $a,b\in\fld[0]{\C}{\pt}$, let $\xi\in\fld[1]{\C}{D^1}[\fcj{a}\amalg b]$ be a 1-field in $\C$, and set $\xi^\star\coloneq\fcj{\iota_*\xi}$ where $\iota\colon D^1\to\orev{D^1}$ is $x\mapsto -x$.
\begin{itemize}
    \item We call $\xi$ an \defn{isometry} if $\beta_\xi\coloneq\xi\times S^1\cup b\times D^2$ equals $a\times D^2$ as (quotient-level) 2-fields:
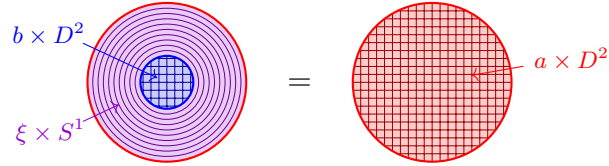
\[
\begin{tkz}[scale=0.7]
    \fill[violet!20] (0,0) circle (1.5);
    \foreach \r in {0.6,0.7,...,1.4} {
        \draw[violet!70!black, very thin] (0,0) circle (\r);
    }
    \fill[blue!20] (0,0) circle (0.5);
    \begin{scope}
        \clip (0,0) circle (0.5);
        \draw[blue!60!black, step=0.15, very thin] (-0.5,-0.5) grid (0.5,0.5);
    \end{scope}
    \draw[red,thick] (0,0) circle (1.5);
    \draw[blue,thick] (0,0) circle (0.5);
    \draw[<-,blue] (155:0.2) -- (155:1.75) node[above left=-1.75mm] {\scriptsize$b \times D^2$};
    \draw[<-,violet] (205:1.0) -- (205:1.75) node[below left=-1.75mm] {\scriptsize$\xi \times S^1$};
\end{tkz}
\quad
=
\quad
\begin{tkz}[scale=0.7]
    \fill[red!20] (0,0) circle (1.5);
    \begin{scope}
        \clip (0,0) circle (1.5);
        \draw[red!60!black, step=0.15, very thin] (-1.5,-1.5) grid (1.5,1.5);
    \end{scope}
    \draw[red,thick] (0,0) circle (1.5);
    \draw[<-,red] (10:0.75) -- (10:1.75) node[right,yshift=.5ex] {\scriptsize$a \times D^2$};
\end{tkz}
\]
\item We call $\xi$ a \defn{coisometry} if $\beta_{\xi^\star}\coloneq\xi^\star\times S^1\cup a\times D^2$ equals $b\times D^2$ as (quotient-level) 2-fields:
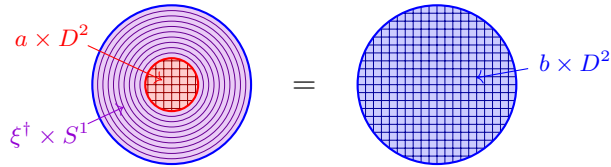
\[
\begin{tkz}[scale=0.7]
    \fill[violet!20] (0,0) circle (1.5);
    \foreach \r in {0.6,0.7,...,1.4} {
        \draw[violet!70!black, very thin] (0,0) circle (\r);
    }
    \fill[red!20] (0,0) circle (0.5);
    \begin{scope}
        \clip (0,0) circle (0.5);
        \draw[red!60!black, step=0.15, very thin] (-0.5,-0.5) grid (0.5,0.5);
    \end{scope}
    \draw[blue,thick] (0,0) circle (1.5);
    \draw[red,thick] (0,0) circle (0.5);
    \draw[<-,red] (155:0.2) -- (155:1.75) node[above left=-1.75mm] {\scriptsize$a \times D^2$};
    \draw[<-,violet] (205:1.0) -- (205:1.75) node[below left=-1.75mm] {\scriptsize$\xi^\star \times S^1$};
\end{tkz}
\quad
=
\quad
\begin{tkz}[scale=0.7]
    \fill[blue!20] (0,0) circle (1.5);
    \begin{scope}
        \clip (0,0) circle (1.5);
        \draw[blue!60!black, step=0.15, very thin] (-1.5,-1.5) grid (1.5,1.5);
    \end{scope}
    \draw[blue,thick] (0,0) circle (1.5);
    \draw[<-,blue] (10:0.75) -- (10:1.75) node[right,yshift=.5ex] {\scriptsize$b \times D^2$};
\end{tkz}
\]
\item We call $\xi$ \defn{unitary} or a \defn{unitary equivalence} if $\xi$ is both an isometry and a coisometry. We write $a\cong^\star b$ if there is a unitary equivalence $\xi$ between $a$ and $b$.
\end{itemize}
For a general 1-ball $Y$, we call a 1-field $\xi\in\fld[1]{\C}{Y}$ an isometry (resp. coisometry, unitary) if there is a homeomorphism $\varphi\colon Y\to D^1$ such that $\varphi_*\xi$ is an isometry (resp. coisometry, unitary).
\end{definition}

We leave the proof of the following lemma to the reader.
\begin{lemma}\label{equivunitaryequiv}
For a homeomorphism $\varphi\colon W\to W'$ of 1-balls, a 1-field $\xi$ is an isometry (resp. coisometry, unitary equivalence) if and only if $\varphi_*\xi$ is an isometry (resp. coisometry, unitary equivalence). 
\end{lemma}

\subsection{Pivotal \texorpdfstring{$\Cstar$}{C*}-2-categories and proto-3-Hilbert spaces}

\subsubsection{\texorpdfstring{$\dag$}{Dagger}- and \texorpdfstring{$\Cstar$}{C*}-2-categories}

\

\begin{definition}[{{\cite[Defn. 2.2]{CHPJP22}}}]
  \label{dagger 2-category}
  A \defn{dagger 2-category} is a linear 2-category $\fX$ equipped with a \defn{dagger structure}, that is, a collection of conjugate-linear maps $\dag\colon\fX({}_aX_b\Rightarrow{}_aY_b) \to \fX({}_aY_b\Rightarrow{}_aX_b)$ indexed over 1-morphisms ${}_aX_b, {}_aY_b \in \fX(a \to b)$ for all objects $a, b\in \fX$ such that
  \begin{lst}
    \item[(2$\dag$a)]\label{2daga} $f^{\dag\dag}=f$ for all $f\in \fX\left({ }_aX_b \Rightarrow{ }_aY_b\right)$,
    \item[(2$\dag$b)]\label{2dagb} $(g \circ f)^{\dag}=f^{\dag} \circ g^{\dag}$ for all $f \in \fX\left({ }_aX_b \Rightarrow{ }_aY_b\right)$ and $g \in \fX\left({ }_aY_b \Rightarrow{ }_aZ_b\right)$,
    \item[(2$\dag$c)]\label{2dagc} $\left(f \otimes g\right)^{\dag}=f^{\dag} \otimes g^{\dag}$ for all $f \in \fX\left({ }_aW_b \Rightarrow{ }_aX_b\right)$ and $g \in \fX\left({ }_bY_c \Rightarrow{ }_bZ_c\right)$, and
    \item[(2$\dag$d)]\label{2dagd} all unitors and associators in $\fX$ are unitary with respect to $\dag$.
  \end{lst}
  A (dagger) \defn{functor} between dagger 2-categories $\fX$ and $\fY$ is an ordinary functor of the underlying 2-categories $F\colon \fX\to\fY$ such that $F(f^\dag)=F(f)^\dag$ for all 2-morphisms $f$ in $\fX$ and whose unitors and tensorators are unitary.
\end{definition}

\begin{definition}[{{\cite[Defn. 2.25]{3Hilb}}}]
  \label{2.25}
  For a 2-category $\fX$, the \defn{linking ($\rE_1$-)algebra} of objects $a,b\in\fX$ is the monoidal category
  \[
    \cL(a, b)\coloneq
    \begin{bmatrix}\fX(a\to a)&\fX(b\to a) \\ \fX(a\to b) & \fX(b\to b)\end{bmatrix}
  \]
  whose monoidal product is defined by matrix multiplication and whose morphism composition is componentwise. The linking algebra $\cL(a_1,\dots,a_{\ell})$ for $a_1,\dots,a_\ell\in\fX$ is defined similarly.
\end{definition}

\begin{definition}[{{\cite[Defn. 2.26]{3Hilb}}}]\label{def:presemisimple}
  A linear 2-category is called \defn{presemisimple} if all linking algebras are semisimple multitensor categories. A presemisimple 2-category is called \defn{finite} if (i) all linking algebras are multi\emph{fusion} and (ii) there is a $K>0$ such that for any linking algebra $\cL$, $\dim\End_{\cZ(\cL)}(1)\leq K$, where $\cZ(\cL)$ denotes the 1-center (Drinfeld center) of $\cL$.
  A linear 2-category is (\defn{finite}) \defn{protosemisimple} if its local (Cauchy) completion is (finite) presemisimple (see \cref{localcomp} for more details on local completion). 
\end{definition}

A \defn{$\Cstar$-monoidal category} is a unitary category equipped with a monoidal structure whose tensor product is a dagger functor and whose associators and unitors are unitary; it is \defn{rigid} if every object has a dual.

\begin{definition}\label{def:Cstar-2-cat}
  A dagger 2-category $\fX$ is called a (\defn{rigid}) \defn{$\Cstar$-2-category} if every $\ell$-fold linking algebra is a (rigid) $\Cstar$-monoidal category. A rigid $\Cstar$-2-category $\fX$ is \defn{finite} if it is finite protosemisimple (\cref{def:presemisimple}), i.e., if its local completion $\fX^{1\cent}$ is finite presemisimple.

Equivalently, a rigid $\Cstar$-2-category $\fX$ is finite if (i) the completion $\cL^\cent$ of every linking algebra $\cL$ of $\fX$ is multi\emph{fusion} and (ii) $\fX$ has finitely many components, i.e., $|\pi_0\fX|<\infty$ (see \cref{def:connected} below and \cref{rem:components-blocks}). Here we remind the reader that all hom spaces in this article are finite-dimensional, so that the local completion of a rigid $\Cstar$-2-category is automatically presemisimple: each of its linking algebras is a complete rigid $\Cstar$-multitensor category with finite-dimensional hom spaces, and thus is semisimple.
\end{definition}

Notice that if $\fX$ is a $\Cstar$-2-category, then so is the delooping of every linking algebra $\cL=\cL(a_1,\dots,a_\ell)$ when equipped with the dagger-transpose.

\subsubsection{Pivotal dagger 2-categories}
\label{sec:pivotal 2-category}

We note that related material on $\Cstar$- and pivotal dagger 2-categories appears in Longo--Roberts \cite{LR97}, Giorgetti--Longo \cite{GL19}, and Verdon \cite{V20}; see also \cite{Pen20,3Hilb}.

A \defn{pivotal 2-category} is a rigid 2-category $\fX$ equipped with a \defn{pivotal structure}, i.e., an iconic natural equivalence $\phi\colon\id_\fX\overset\cong\Rightarrow \vee\circ\vee$. This definition unpacks as follows. A pivotal structure is the data of a 2-isomorphism $\phi_X\in\fX(X\Rightarrow X^{\vee\vee})$ for all 1-morphisms ${}_aX_b\in\fX$ satisfying the following conditions.
\begin{lst} 
\item[(C$\phi$N)]\label{CphiN} For any 2-morphism $f\in\fX({}_aX_b\Rightarrow{}_aY_b)$, the naturality condition $f^{\vee\vee}\circ\phi_X=\phi_Y\circ f$ holds.
\item[(C$\phi\otimes$)] For all composable 1-morphisms ${}_aX_b$ and ${}_bY_c$ in $\fX$, we have $(\vee^2_{X,Y})^\vee\circ \phi_{X\otimes Y}=\vee^2_{Y^\vee,X^\vee}\circ(\phi_X\otimes\phi_Y)$.
\item[(C$\phi u$)]\label{Cphiu} For all objects $a\in\fX$, we have $(\vee^0_a)^\vee\circ\phi_{1_a}=\vee^0_a$.
\end{lst}

By applying \cite[Lem. 4.11]{Sel11} to the linking algebras of $\fX$, we find that pivotal 2-categories $(\fX,\phi)$ enjoy the following property.
\begin{lst}
\item[(C$\phi\vee$)]\label{CphiV} For all 1-morphisms $X$ in $\fX$, $\phi^\vee_X=\phi^{-1}_{X^\vee}$.
\end{lst}

\begin{definition}
\label{def:uaf}
For a rigid dagger 2-category $\fX$, a \defn{unitary adjoint functor} (UAF) is a choice of adjoint data $(X^\vee,\coev_X,\ev_X)$ for each 1-morphism $X$ in $\fX$ such that $\vee$ together with its canonical unitors and tensorators assembles into a dagger 2-functor $\vee\colon\fX\to\fX^{1\op,2\op}$.
\end{definition} 

\begin{fact}[\cite{Pen20,3Hilb}]
\label{canpivstr}
A rigid dagger 2-category $\fX$ equipped with a choice of UAF $\vee$ carries a canonical pivotal structure given for a 1-morphism $X\in\fX(a\to b)$ by
\[
\phi_X\coloneq
\begin{tkz}
\draw(0,0)--node[left=-1mm,pos=0.3]{$\scriptstyle X$}(0,1.35);
\draw(0.7,1.35)--node[right=-1.2mm,pos=0.8]{$\scriptstyle X^{\vee}$}(0.7,0.75)arc(180:360:.35);
\draw(1.4,0.75)--node[right=-1mm,pos=0.7]{$\scriptstyle X^{\vee\vee}$}(1.4,2);
\rbox[fill=white]{(0.35,1.35)}{0.3}{0.3}{0.3}{$\coev_X^\dag$};
\end{tkz}
=
(\coev_{X}^\dag\otimes \id_{X^{\vee\vee}})\circ(\id_{X}\otimes\coev_{X^\vee})
\]
(where we have suppressed unitors, associators, and the shadings for the objects $a$ and $b$). 
\end{fact}

\cref{canpivstr} motivates the following definition.

\begin{definition} 
\label{def:pivotal 2-category}
A \defn{pivotal dagger 2-category} is a rigid dagger 2-category $\fX$ equipped with a choice of UAF $\vee$. 
\end{definition}

A pivotal structure $\phi$ on $\fX$ is called \defn{strict} if $\phi=\id_\fX$, i.e., if $\phi_X=\id_X$ for all 1-morphisms $X$ in $\fX$. See \cref{strictification-sec} for the proof of the following strictification result.

\strictification*

This justifies the following convention.

\begin{convention} 
Henceforth we will always assume that all pivotal dagger 2-categories $\fX$ are strictly pivotal and have strict underlying 2-categories.
\end{convention}

\subsubsection{Proto-3-Hilbert spaces}

\begin{definition}\label{def:connected}
  Let $\fX$ be a pivotal $\Cstar$-2-category. For each $a\in\fX$, write $P_a$ for the set of minimal projections of the (commutative) unitary algebra $\End_\fX(1_a)$. Define a \defn{point} of $\fX$ to be a pair $(a,p)$ consisting of an object $a\in \fX$ and a minimal projection $p\in P_a$. We will call points $(a,p)$ and $(b,q)$ \defn{connected} if there is some 1-morphism ${}_aX_b\in\fX(a\to b)$ such that $p\otimes\id_X\otimes q\neq 0$. Since 1-morphisms are adjointable and the horizontal composites $(p\otimes\id_X\otimes q)\otimes(q\otimes\id_Y\otimes s)$ of nonzero such 2-morphisms are nonzero, this is an equivalence relation on points. We will call its equivalence classes the \defn{components} of $\fX$, and we denote the set of components by $\pi_0\fX$.
\end{definition}

\begin{remark}
  \label{rem:points-completion}
When an object $a$ in $\fX$ is \defn{simple}, meaning $\End_\fX(1_a)\cong\bC$, we have $P_a=\{\id_{1_a}\}$, so \cref{def:connected} recovers the usual notion of connectedness \cite[Defn. 1.2.22]{DR18}: simple objects $a,b\in\fX$ are \defn{connected} if $\fX(a\to b)\neq0$. Conversely, each point $(a,p)$ determines a simple object $a_p$ of the completion $\fX^\cent$ (defined below), every component of $\fX^\cent$ contains such a point, and $\fX^\cent(a_p\to b_q)\neq0$ if and only if $p\otimes\id_X\otimes q\neq0$ for some $X\in\fX(a\to b)$. Thus the induced map $\pi_0\fX\to\pi_0\fX^\cent$ is a bijection, so a proto-3-Hilbert space $\fX$ is finite if and only if its completion $\fX^\cent$ is finite. 
\end{remark}

\begin{remark}
  \label{rem:components-blocks}
  Suppose $\fX$ is a pivotal $\Cstar$-2-category such that the linking algebra $\cL=\cL(a_1,\dots,a_\ell)$ is complete for all $a_1,\dots,a_\ell\in\fX$. The simple summands of $1_\cL$ are exactly the points of the $a_i$, and by \cite[Rmk. 2.27]{3Hilb} the simple summands of $1_{\cZ(\cL)}$ correspond to the indecomposable fusion blocks of $\cL$, i.e., to the equivalence classes of these points under the connectedness relation of \cref{def:connected}. Any collection consisting of representatives from every component realizes every class, so $\sup_\cL\dim\End_{\cZ(\cL)}(1_{\cZ(\cL)})=|\pi_0\fX|$, i.e., $|\pi_0\fX|$ is the minimal $K$ in condition (ii) of \cref{def:presemisimple}. Thus conditions (ii) of \cref{def:presemisimple} and of \cref{def:Cstar-2-cat} agree for pre-3-Hilbert spaces, so the notions of \emph{finite} in \cite[Defns. 2.26 and 4.3]{3Hilb} coincide.
\end{remark}

\begin{definition}
  A \defn{spherical weight} $\Psi$ for a rigid $\Cstar$-2-category $\fX$ equipped with a UAF $\vee$ consists of linear functionals $\Psi_{a}\colon \End_\fX(1_{a})\to\bC$ for all $a\in \fX$ satisfying
  \begin{lst}
    \item[(i)] (\emph{Faithful weight}) $\Psi_a(f^\dag f)\geq 0$ with equality if and only if $f=0$ for all $a\in\fX$ and $f\in\End_\fX(1_a)$, and
    \item[(ii)] (\emph{Sphericality}) $\Psi_b(\tr^\vee_R(f)) = \Psi_a(\tr^\vee_L(f))$ for all $a,b\in\fX$, ${}_aX_b\in\fX(a\to b)$, and $f\in\End_\fX(X)$, where $\tr^\vee_R(f)\in\End_\fX(1_b)$ and $\tr^\vee_L(f)\in\End_\fX(1_a)$ are respectively the right and left traces formed by closing $f$ with the cups and caps of $\vee$; see \cite[Rmk. 2.19]{3Hilb}.
  \end{lst}
\end{definition}

\begin{definition}
  \label{def:3hilb-spaces}
  A \defn{proto-3-Hilbert space} $(\fX,\vee,\Psi)$ is a pivotal $\Cstar$-2-category $(\fX,\vee)$ equipped with a spherical weight $\Psi$. 
  A \defn{pre-3-Hilbert space} $(\fX,\vee,\Psi)$ is a proto-3-Hilbert space whose linking algebras are unitary multitensor categories. When $\fX$ is finite, this is equivalent to \cite[Defn. 4.3]{3Hilb} by \cref{rem:components-blocks}.
\end{definition}

\begin{example}[{{\cite[Exm. 4.4]{3Hilb}}}]\label{ex:delooping-pre-3-hilb}
An example of a pre-3-Hilbert space is the delooping $\rB\cA$ for an $\Hstar$-multifusion category $\cA$ in the sense of {{\cite[Defn. 3.1]{3Hilb}}} (see \cref{classificationHstarmFC} below). An example of a proto-3-Hilbert space is $\rB^2 A$ for a commutative $\Hstar$-algebra $A$ (see \cref{def:Hstaralg}).  
\end{example}

\subsubsection{(Co)isometries and isometric equivalences}

\begin{definition}[UAF-preserving functors and isometries, {{\cite[Defn. 4.14]{3Hilb}}}]
  \label{UAF-preserving}
  For pivotal dagger 2-categories $(\fX,\vee_\fX)$ and $(\fY,\vee_\fY)$, a dagger functor $F\colon\fX\to\fY$ is called \defn{UAF-preserving} if for every 1-morphism ${}_aX_b\in\fX(a\to b)$, the canonical isomorphism $\delta_X\in\fY(F(X^\vee)\Rightarrow F(X)^\vee)$ given by
  \[
    \delta_X\coloneq(F(\ev_X)\otimes F(X)^\vee)\circ (F_{X^\vee,X}^2\otimes F(X)^\vee)\circ(F(X^\vee)\otimes \coev_{F(X)})
  \]
  (suppressing associators and unitors) is unitary. We write $\Hom(\fX\to\fY)$ or $\Fun^{\dag,\vee}(\fX\to\fY)$ for the full sub-2-category of $\Fun^\dag(\fX\to\fY)$ consisting of UAF-preserving dagger functors $\fX\to\fY$. For proto-3-Hilbert spaces $\fX$ and $\fY$, we call $F\in\Hom(\fX\to\fY)$ an \defn{isometry} if $\Psi^\fY_{F(a)}(F(f))=\Psi^\fX_a(f)$ for all $f\in\End_\fX(1_a)$ and $a\in\fX$.
\end{definition}

\begin{definition}[Isometric equivalence {{\cite[Defn. 4.14]{3Hilb}}}]\label{def:isometric-equivalence}
  Two proto-3-Hilbert spaces $\fX,\fY$ are \defn{isometrically equivalent} if there is a fully faithful isometry $F\colon\fX\to\fY$ in the sense of \cref{UAF-preserving} that is \defn{isometrically essentially surjective}, i.e., for every $b\in\fY$ there exist $a\in\fX$ and an isometric equivalence---that is, an adjoint equivalence in $\fY$ with respect to $\vee$---between $F(a)$ and $b$.
\end{definition}

\begin{definition}[{\cite[Defn. 4.10, Rmk. 4.12]{3Hilb}}]
A 1-morphism ${}_aX_b$ in a proto-3-Hilbert space $\fX$ is called an \defn{isometry} (resp. \defn{coisometry}) if $\coev_X\colon 1_a\Rightarrow {}_aX\otimes_b X_a^\vee$ (resp. $\ev_X\colon {}_bX^\vee\otimes_a X_b\Rightarrow 1_b$) is unitary. This definition makes sense for pivotal $\Cstar$-2-categories too.
\end{definition}

\begin{lemma}
\label{isomequiviff}
Let $\C$ be a unitary disk-like 2-category. A 1-field in $\C$ is a unitary equivalence if and only if it skeletonizes to an isometric equivalence in the sense of \cite[Defn. 4.10]{3Hilb}. In particular, 0-fields $a$ and $b$ in $\C$ satisfy $a\cong^\star b$ in $\C$ if and only if $a$ and $b$ are isometrically equivalent objects in $\fX_\C$.
\end{lemma}

\begin{pf} 
Observe that for $x\coloneq{}_a\xi_b\in\fX_\C(a\to b)$, as 2-morphisms of $\fX_\C$ we have $\beta_\xi = \coev_x^\dag \circ \coev_x$ and $\beta_{\xi^\star} = \ev_x \circ \ev_x^\dag$, so $\xi$ is isometric (resp. coisometric) if and only if $\beta_\xi = \id_{1_a}$ (resp. $\beta_{\xi^\star} = \id_{1_b}$), i.e., if and only if $\coev_x$ (resp. $\ev_x$) is an isometry (resp. coisometry). By the argument in \cite[Rmk. 4.12]{3Hilb}, this is equivalent to being an isometry (resp. coisometry) in $\fX_\C$.
\end{pf}

\subsubsection{Local completion of proto-3-Hilbert spaces}
\label{localcomp}
For a dagger 2-category $\fX$, its \defn{local completion} (or \defn{1-Cauchy completion}) $\fX^{1\cent}$ is the dagger 2-category with the same objects as $\fX$ and hom categories $\fX^{1\cent}(a\to b)\coloneq\fX(a\to b)^\cent$ the completions of those of $\fX$ (\cref{completeunitarycategories}). The composition, dagger, and coherence data are uniquely extended componentwise via the universal property of completion (\cref{3.3.12}). 

We next show that the local completion of a proto-3-Hilbert space has a canonical pre-3-Hilbert space structure.

\begin{lemma}
  \label{cstr:local-cauchy-3hilb}
  \begin{lst}
    \item[(a)]\label{lcc1} If $\cX$ is a 
    rigid 
    $\Cstar$-monoidal category equipped with a UDF $\vee$ and a spherical weight $\psi$ on $\End_\cX(1)$, and $\cX^\cent$ has finitely many isomorphism classes of simple objects, then its completion $\cX^\cent$ is canonically an $\Hstar$-multifusion category and the inclusion
    $\iota\colon\cX\hookrightarrow\cX^\cent$ is isometric.
    \item[(b)]\label{lcc2} If $(\fX,\vee,\Psi)$ is a finite proto-3-Hilbert space, then its local completion $\fX^{1\cent}$ admits a canonical finite pre-3-Hilbert space structure $(\vee^{1\cent},\Psi^{1\cent})$ and the inclusion $\iota_{1\cent}\colon\fX\hookrightarrow\fX^{1\cent}$ is isometric.
  \end{lst}
\end{lemma}

\begin{pf}
  For \ref{lcc1}, start by noting that $\cX \xrightarrow{\vee} \cX^{1\op,2\op} \hookrightarrow (\cX^{\cent})^{1\op,2\op}$ has a complete target, and thus by \cref{3.3.12} extends to a functor $\vee^{\cent} \colon \cX^{\cent} \to (\cX^{\cent})^{1\op,2\op}$. Define $\psi^{\cent}$ by $\psi^{\cent}(\iota_{\cent}(f)) \coloneq \psi(f)$, so that $\iota_\cent$ is automatically isometric; this is well-defined since $\iota_{\cent}$ restricts to a $*$-isomorphism $\End_\cX(1_\cX)\to\End_{\cX^{\cent}}(1_{\cX^{\cent}})$.
  
  It remains to show $\psi^\cent$ is spherical. For $f = pfp \in \End((Y,p))$ with $Y = \bigoplus_i X_i$, defining $\ev_{(Y,p)} \coloneq \ev_Y \circ (p^\vee \otimes p)$ and $\coev_{(Y,p)} \coloneq (p \otimes p^\vee) \circ \coev_Y$ gives $\tr^\vee_R(f)_{(Y,p)} = \tr^\vee_R(pfp)_Y$, and similarly for $\tr^\vee_L$. Writing $pfp = ((pfp)_{ij})$ as a matrix with $(pfp)_{ij} \in \Hom(X_i \to X_j)$, the fact that $\ev_Y$ and $\coev_Y$ act componentwise on the direct sum gives
  \[
    \psi^{\cent}(\tr^\vee_R(f)_{(Y,p)})
    = \sum_i \psi(\tr^\vee_R((pfp)_{ii})_{X_i})
    = \sum_i \psi(\tr^\vee_L((pfp)_{ii})_{X_i})
    = \psi^{\cent}(\tr^\vee_L(f)_{(Y,p)}),
  \]
  where the second equality uses sphericality of $\psi$.

    Finally, \ref{lcc2} follows from applying \ref{lcc1} to each linking algebra $\cL(a,b)$, whose completion $\cL(a,b)^\cent$ is multifusion since $\fX$ is finite and hence is organically an $\Hstar$-multifusion category (cf. \cite[Exm. 4.5]{3Hilb}); its UAF and spherical weight on $\fX^{1\cent}(a\to b)$ are the restrictions of those on $\cL(a,b)^{\cent}$ to its off-diagonal block. Finiteness of $\fX^{1\cent}$ follows by applying \cref{rem:points-completion} to $\fX$ and to $\fX^{1\cent}$, whose completions agree.
\end{pf}

\subsubsection{3-Hilbert spaces}
\begin{definition}[{\cite[Defn. 3.16]{bases}}]
  A \defn{Hilbert direct sum} of $a,b\in\fX$ is an object $a\boxplus b\in\fX$ equipped with isometries $I\colon a\hookrightarrow a\boxplus b$ and $J\colon b\hookrightarrow a\boxplus b$ satisfying $\ev_I\ev_I^\dag+\ev_J\ev_J^\dag=\id_{1_{a\boxplus b}}$.
\end{definition}

\begin{definition}[{\cite[Defn. 3.18]{bases}}]
An \defn{$\Hstar$-monad} on an object $a\in\fX$ is an algebra $({}_aA_a,\mu,\iota)$ internal to $\End_\fX(a)$ whose multiplication $\mu\colon {}_aA\otimes_a A_a\Rightarrow {}_aA_a$ is such that $\mu^\dag$ is an $A$--$A$-bimodule map, $\mu\mu^\dag$ is invertible, and $({}_aA_a,\mu,\iota)$ is \defn{standard} in the sense of \cite[Defn. 3.18, ($\Hstar$3)]{3Hilb}. A \defn{splitting} of an $\Hstar$-monad $({}_aA_a,\mu,\iota)$ is a 1-morphism ${}_aX_b$ in $\fX$ together with a unitary monad isomorphism $\gamma\colon {}_aA_a\Rightarrow {}_aX\otimes_b X_a^\vee$ such that $\ev_X\ev_X^\dag\in\End_\fX(1_b)$ is invertible.
\end{definition}

\begin{definition}[{\cite[Defns. 4.45 and 4.47]{3Hilb}}]
  \label{def:completion-3hilb}
  A \defn{3-Hilbert space} is a pre-3-Hilbert space admitting Hilbert direct sums and whose $\Hstar$-monads split. For a pre-3-Hilbert space $\fX$, its \defn{completion} $\fX^\cent\coloneq\mathsf{H^*Alg}(\Hilb_\boxplus(\fX))$ is the 3-Hilbert space obtained from $\fX$ by formally adjoining Hilbert direct sums of objects and splitting $\Hstar$-monads; see \cite{3Hilb} for the details of these two constructions. For a proto-3-Hilbert space $\fX$, we set $\fX^\cent\coloneq(\fX^{1\cent})^\cent$, where $\fX^{1\cent}$ is the local completion of \cref{localcomp}.
\end{definition}

\subsection{From disk-like to traditional: the proto-3-Hilbert space \texorpdfstring{$\fX_\C$}{X(C)}}
\label{skein 2-categories}
Let $\C$ be a disk-like 2-category. In this subsection, we expand on the sketch in \cite[Appendix C.2]{MW12} of the construction of an ordinary (traditional, weak, algebraic) pivotal dagger 2-category $\fX_\C$ from $\C$. 

For each $k$, let $D^k$ denote the standard $k$-dimensional globe, let $D^{k-1}$ denote the standard equatorial disk cutting $D^k$ into two halves, and let $S^{k-2}$ denote the equator $\partial D^{k-1}$. With the pinching convention of \cref{disk-like skein k-categories}, the globes satisfy $D^k=D^{k-1}\times D^1$, and $\iota^{(2)}=\id_{D^1}\times\iota^{(1)}$. Fix homeomorphisms $\rho^1 \colon D^1 \cup_{D^0} D^1 \to D^1$ and $\rho^2\coloneq\id_{D^1}\times\rho^1 \colon D^2 \cup_{D^1} D^2 \to D^2$. For $j=-,+$, let $\pr^1_j$ denote the projection of $\undfld[0]{\C}{\partial D^1}=\undfld[0]{\C}{\{-1\}}\times\undfld[0]{\C}{\{+1\}}$ onto its $j$th factor.

\subsubsection{\texorpdfstring{$\fX_\C$}{X(C)} as a linear 2-category}

\

\label{sec:skein 2-category}
\begin{construction}
    \label{construction:skein 2-category}
    The 2-category $\fX_\C$ is defined as follows.
    \begin{lst}[widestlabel]
        \item
        $\fX^0_\C\coloneq\Obj(\fX_\C)\coloneq\fld[0]{\C}{\pt}$.

        \item
        $\fX^1_\C\coloneq\fld[1]{\C}{D^1}$. For a 1-field $\alpha$ on $D^1$, define $\src(\alpha)\coloneq\fcj{\pr^1_-(\partial\alpha)}$ and $\targ(\alpha)\coloneq\pr^1_+(\bdy\alpha)$. The hom 1-category $\fX_\C(a \to b)$ has objects $\Obj(\fX_{\C}(a \to b))\coloneq\src^{-1}(a) \cap \targ^{-1}(b)$.

        \item
        $\fX_\C^2\coloneq\bigoplus_{c\in\undfld[1]{\C}{S^1}_{\pitchfork S^0}}\fldl{\C}{D^2}[c]$, and, for ${}_aX_b,{}_aY_b\in\fX_\C$, hom vector spaces $\fX_{\C}({}_aX_b \Rightarrow {}_aY_b) \coloneq \fldl{\C}{D^2}[\fcj{X}\cup Y]$, which are finite-dimensional by \ref{FDF}.

        \item[($\circ$)]\label{2skelvcomp}
        The 2-composition maps $\circ \colon \fX_{\C}({}_aY_b \Rightarrow {}_aZ_b) \otimes_\bC \fX_{\C}({}_aX_b \Rightarrow {}_aY_b)\to\fX_{\C}({}_aX_b \Rightarrow {}_aZ_b)$ are defined for $f\in\fX_\C({}_aX_b\Rightarrow{}_aY_b)$ and $g\in\fX_\C({}_aY_b\Rightarrow{}_aZ_b)$ by
        \[
            \begin{tkz}
                \fill[\filg](0,0) .. controls (0.5,\hgt) and (1.5,\hgt) .. (2,0).. controls (1.5,-\hgt) and (0.5,-\hgt) .. (0,0);
                \draw[\colZ,thick] (0,0) .. controls (0.5,\hgt) and (1.5,\hgt) .. node[above]{\scriptsize${}_aZ_b$}(2,0);
                \draw[\colY,thick] (0,0) .. controls (0.5,-\hgt) and (1.5,-\hgt) .. node[below]{\scriptsize$\fcj{{}_aY_b}$} (2,0);
                \filldraw[\cola] (0,0) circle (1.5pt) node[left]{\scriptsize$\fcj{a}$};
                \filldraw[\colb] (2,0) circle (1.5pt) node[right]{\scriptsize$b$};
                \node at (current bounding box.center){$g$};
            \end{tkz}
            \circ
            \begin{tkz}
                \fill[\filf](0,0) .. controls (0.5,\hgt) and (1.5,\hgt) .. (2,0).. controls (1.5,-\hgt) and (0.5,-\hgt) .. (0,0);
                \draw[\colY,thick] (0,0) .. controls (0.5,\hgt) and (1.5,\hgt) .. node[above]{\scriptsize${}_aY_b$} (2,0);
                \draw[\colX,thick] (0,0) .. controls (0.5,-\hgt) and (1.5,-\hgt) ..  node[below]{\scriptsize$\fcj{{}_aX_b}$} (2,0);
                \filldraw[\cola] (0,0) circle (1.5pt) node[left]{\scriptsize$\fcj{a}$};
                \filldraw[\colb] (2,0) circle (1.5pt) node[right]{\scriptsize$b$};
                \node at (current bounding box.center){$f$};
            \end{tkz}
            \;
            =
            \;
            g \circ f
            \coloneq
            \rho^2_*(f \blt g)
            \;
            =
            \;
            \begin{tkz}
                \begin{scope}
                    \clip(0,0) .. controls (0.5,\hgt) and (1.5,\hgt) .. (2,0).. controls (1.5,-\hgt) and (0.5,-\hgt) .. (0,0);
                    \fill[\filg](-1,0) rectangle (3,2);
                    \fill[\filf](-1,-2) rectangle (3,0);
                \end{scope}
                \draw[\colZ,thick] (0,0) .. controls (0.5,\hgt) and (1.5,\hgt) ..node[above]{\scriptsize${}_aZ_b$} (2,0);
                \draw[\colX,thick] (0,0) .. controls (0.5,-\hgt) and (1.5,-\hgt) ..node[below]{\scriptsize$\fcj{{}_aX_b}$} (2,0);
                \filldraw[\cola,thick] (0,0) circle (1.5pt);
                \filldraw[\colb,thick] (2,0) circle (1.5pt);
                \node at (1,.375*\hgt){\scriptsize$g$};
                \node at (1,-.4*\hgt){\scriptsize$f$};
                \draw[dashed](0,0)--(2,0);
            \end{tkz}\,.
        \]

        \item[(id)]\label{2skelvid}
        For ${}_aX_b\in\fX_\C(a\to b)$, the identity 2-morphism $\id_{{}_aX_b}$ is defined by
        \begin{equation*}
            \id_{X} \coloneq {}_aX_b \times I = \begin{tkz}
                \fill[\filid](0,0) .. controls (0.5,\hgt) and (1.5,\hgt) .. (2,0).. controls (1.5,-\hgt) and (0.5,-\hgt) .. (0,0);
                \draw[\colX,thick] (0,0) .. controls (0.5,\hgt) and (1.5,\hgt) .. node[above,\colX]{\scriptsize${}_aX_b$}(2,0);
                \draw[\colX,thick] (0,0) .. controls (0.5,-\hgt) and (1.5,-\hgt) .. node[below,\colX]{\scriptsize$\fcj{{}_aX_b}$} (2,0);
                \begin{scope}
                    \clip(0,0) .. controls (0.5,\hgt) and (1.5,\hgt) .. (2,0).. controls (1.5,-\hgt) and (0.5,-\hgt) .. (0,0);
                    \foreach \x in {0.2,0.4,...,1.8} {\draw[\colX,thick](\x,2)--(\x,-2);}
                \end{scope}
                \filldraw[\cola] (0,0) circle (1.5pt) node[left] {\scriptsize$\fcj{a}$};
                \filldraw[\colb] (2,0) circle (1.5pt) node[right] {\scriptsize$b$};
            \end{tkz}.
        \end{equation*}

        \item[($\otimes$)]\label{2skelhcomp}
        The 1-composition functors $\otimes \colon \fX_{\C}(a \to b) \times \fX_{\C}(b \to c) \to \fX_{\C}(a \to c)$ are defined on 1-morphisms ${}_aX_b$, ${}_bY_c$ in $\fX_\C$ by
        \[
            \begin{tkz}
                \draw[mid>={.55},thick,\colX](0,0)--node[above,overlay]{\scriptsize${}_aX_b$}(2,0);
                \foreach \x/\d/\l/\cc in {0/left/a/\cola,2/right/b/\colb}{
                        \filldraw[\cc] (\x,0) circle (1pt) node[\d] {\scriptsize$\l$};}
            \end{tkz}
            \otimes
            \begin{tkz}
                \draw[mid>={.55},thick,\colY](0,0)--node[above,overlay]{\scriptsize${}_bY_c$}(2,0);
                \foreach \x/\d/\l/\cc in {0/left/b/\colb,2/right/c/\colc}{
                        \filldraw[\cc] (\x,0) circle (1pt) node[\d] {\scriptsize$\l$};}
            \end{tkz}
            % \;
            =
            % \;
            {}_aX_b \otimes {}_bY_c
            % \;
            \coloneq
            % \;
            \rho^1_*({}_aX\blt_b Y_c)
            % \;
            =
            % \;
            \begin{tkz}
                \draw[mid>={.65},thick,\colX](0,0)--(1,0)node[above,overlay]{\scriptsize\color{violet}${}_aX\blt_b Y_c$};
                \draw[mid>={.65},thick,\colY](1,0)--(2,0);
                \foreach \x/\d/\l/\cc in {0/left/a/\cola,2/right/c/\colc}{\filldraw[\cc] (\x,0) circle (1pt) node[\d] {\scriptsize$\l$}; }
                \filldraw[\colb](1,0)circle(1pt)node[overlay,below=0mm]{\scriptsize$b$};
            \end{tkz}
        \]
        and on 2-morphisms $f\in\fX_\C({}_aX_b\Rightarrow {}_aW_b)$, $g\in\fX_\C({}_bY_c\Rightarrow {}_bZ_c)$ by
        \[
        f\otimes g\coloneq
        \begin{tkz}[scale=0.7]
            \coordinate (A) at (-3, 0);
            \coordinate (B) at (0, 0);
            \coordinate (C) at (3, 0);
            \coordinate (B') at (0, -2);
            \fill[\filf] (A) to[bend left=45] (B) to[bend left=45] (A);
            \fill[\filg] (B) to[bend left=45] (C) to[bend left=45] (B);
            \begin{scope}
                \clip (A) to[bend right=45] (B) -- (B') to[out=180, in=-60] (A) -- cycle;
                \foreach \x in {-3, -2.7, ..., 0} { \draw[\colX, thick] (\x, 0) -- (\x, -2.5); }
            \end{scope}
            \draw[\colX, thick] (A) to[out=-60, in=180] node[pos=.65,below] {\scriptsize$\fcj{{}_aX_b}$} (B');
            \draw[\colX,dotted,thick] (A) to[bend right=45] (B);
            \draw[\colW, thick] (A) to[bend left=45] node[pos=0.5,above, text=\colW] {\scriptsize${}_aW_b$} (B);
            \node at (-1.5, 0) {\scriptsize$f$};
            \begin{scope}
                \clip (B) to[bend right=45] (C) -- (C) to[out=-120, in=0] (B') -- (B) -- cycle;
                \foreach \x in {0.3, 0.6, ..., 3} { \draw[\colY,thick] (\x, 0) -- (\x, -2.5); }
            \end{scope}
            \draw[\colY, thick] (B') to[out=0, in=-120] node[pos=.65,below, text=\colY] {\scriptsize$\fcj{{}_bY_c}$} (C);
            \draw[\colY,dotted,thick] (B) to[bend right=45] (C);
            \draw[\colZ, thick] (B) to[bend left=45] node[above, text=\colZ] {\scriptsize${}_bZ_c$} (C);
            \node at (1.5, 0) {\scriptsize$g$};
            \draw[\colb, thick] (B) -- (B');
            \fill[\cola] (A) circle (1.5pt) node[left] {\scriptsize$\fcj{a}$};
            \fill[\colb] (B) circle (1.5pt);
            \fill[\colc] (C) circle (1.5pt) node[right] {\scriptsize$c$};
            \fill[\colb] (B') circle (1.5pt);
            \draw[<-, thick] (0.75,-1.15) to[in=90,out=-45] (1.5,-1.95) node[below] {\scriptsize$({}_aX \otimes_b Y_c) \times I$};
        \end{tkz}
        \]
        (which we reparameterize to the standard disk $D^2$).

        \item
        For $a\in\fX_\C$, the 1-morphism $1_a\in\fX_\C(a\to a)$ is the product field $1_a\coloneq a\times I$.

        \item
        In \cite[Appendix C.2]{MW12}, the unitors are constructed by reparameterizing half-pinched products, while the associators for $\otimes$ are defined using product fields by a ``shift'' between the two obvious parameterizations $D^1\cup_{D^0} D^1\cup_{D^0} D^1\to D^1$.
    \end{lst}
\end{construction}

\begin{proposition}
    \label{skein 2-categories are linear}
    For a disk-like 2-category $\C$, $\fX_\C$ is a linear 2-category.
\end{proposition}

\begin{pf}
    \itemstep{\ref{2skelvcomp}.}
    Vertical composition is strictly associative by the same argument from the $n=1$ case in \cref{skein 1-categories}.

    \itemstep{\ref{2skelvid}.}
    This indeed defines identities for $\circ$, again by the same reasoning from the $n=1$ case.

    \itemstep{\ref{2skelhcomp}.}
    Identity 2-morphisms are strictly preserved ($\id_{X\otimes Y}=\id_X\otimes\id_Y$) by extended-isotopy invariance, while composition is preserved since for all 2-morphisms $f\colon {}_aX_b\Rightarrow{}_aQ_b$, $g\colon{}_aQ_b\Rightarrow{}_aW_b$, $h\colon{}_bY_c\Rightarrow{}_bR_c$, and $k\colon{}_bR_c\Rightarrow{}_bZ_c$ in $\fX_\C$,
        \begin{align*}
            (g &\circ f) \otimes (k \circ h)
            =
            \begin{tkz}[scale=0.7]
                \coordinate (A) at (-3, 0);
                \coordinate (B) at (0, 0);
                \coordinate (C) at (3, 0);
                \coordinate (B') at (0, -2);
                \begin{scope}
                    \clip (A) to[bend left=45] (B) to[bend left=45] (A);
                    \fill[\filg](-4,0) rectangle (1,2);
                    \fill[\filf](-4,-2) rectangle (1,0);
                \end{scope}
                \begin{scope}
                    \clip (B) to[bend left=45] (C) to[bend left=45] (B);
                    \fill[\filk](-1,0) rectangle (4,2);
                    \fill[\filh](-1,-2) rectangle (4,0);
                \end{scope}
                \begin{scope}
                    \clip (A) to[bend right=45] (B) -- (B') to[out=180, in=-60] (A) -- cycle;
                    \foreach \x in {-3, -2.7, -2.4, -2.1, -1.8, -1.5, -1.2, -0.9, -0.6, -0.3, 0} {
                            \draw[\colX, thick] (\x, 0) -- (\x, -2.5);
                        }
                \end{scope}
                \draw[\colX, thick] (A) to[out=-60, in=180] node[below] {\scriptsize$\fcj{{}_aX_b}$} (B');
                \draw[\colX,thick,dotted] (A) to[bend right=45] (B);
                \draw[black, dashed] (A) -- (B);
                \draw[\colW, thick] (A) to[bend left=45] node[above] {\scriptsize${}_aW_b$} (B);
                \node at (-1.5, 0.3) {\scriptsize$g$};
                \node at (-1.5, -0.3) {\scriptsize$f$};
                \draw[<-, thick, shorten <=2pt] (-1, .25) to[out=90, in=180] (-0.65, 0.8) node[right] {\scriptsize$g \circ f$};
                \begin{scope}
                    \clip (B) to[bend right=45] (C) -- (C) to[out=-120, in=0] (B') -- (B) -- cycle;
                    \foreach \x in {0.3, 0.6, 0.9, 1.2, 1.5, 1.8, 2.1, 2.4, 2.7, 3} {
                            \draw[\colY, thick] (\x, 0) -- (\x, -2.5);
                        }
                \end{scope}
                \draw[\colY, thick] (B') to[out=0, in=-120] node[pos=.65,right] {\scriptsize$\fcj{{}_bY_c}$} (C);
                \draw[\colY,thick,dotted] (B) to[bend right=45] (C);
                \draw[black, dashed] (B) -- (C);
                \draw[\colZ, thick] (B) to[bend left=45] node[above] {\scriptsize${}_bZ_c$} (C);
                \node at (1.5, 0.3) {\scriptsize$k$};
                \node at (1.5, -0.3) {\scriptsize$h$};
                \draw[<-, thick, shorten <=2pt] (2, 0.25) to[out=90, in=180] (2.35, 0.8) node[right] {\scriptsize$k \circ h$};
                \draw[\colb, thick] (B) -- (B');
                \fill[\cola] (A) circle (1.5pt) node[left] {\scriptsize$\fcj{a}$};
                \fill[\colb] (B) circle (1.5pt);
                \fill[\colc] (C) circle (1.5pt) node[right] {\scriptsize$c$};
                \fill[\colb] (B') circle (1.5pt);
                \draw[<-, thick] (0.75,-1.15) to[in=90,out=-45] (1.5,-1.95) node[below] {\scriptsize$({}_aX \otimes_b Y_c) \times I$};
            \end{tkz}
            =
            \begin{tkz}[scale=0.7]
                \coordinate (A) at (-3, 0);
                \coordinate (B) at (0, 0);
                \coordinate (B2) at (0, 1.5);
                \coordinate (C) at (3, 0);
                \coordinate (B') at (0, -2);
                \fill[\filf] (A) to[bend right=45,looseness=1.5] (B) -- cycle;
                \fill[\filg] (A) to[bend left=50] (B2) to[bend left=20] (A);
                \fill[\filh] (B) to[bend right=45,looseness=1.5] (C) -- cycle;
                \fill[\filk] (B2) to[bend left=50] (C) to[bend left=20] (B2);
                \begin{scope}
                    \clip (A) to[bend right=45,looseness=1.5] (B) -- (B') to[out=180, in=-60] (A) -- cycle;
                    \foreach \x in {-3, -2.7, ..., 0} { \draw[\colX, thick] (\x, 0) -- (\x, -2.5); }
                \end{scope}
                \draw[\colX, thick] (A) to[out=-60, in=180] (B');
                \draw[\colX,thick,dotted] (A) to[bend right=45,looseness=1.5] (B);
                \draw[\colQ, dotted, thick] (A) -- (B);
                \node at (-1.5, -0.4) {\scriptsize$f$};
                \begin{scope}
                    \clip (A) to[bend right=20] (B2) -- (B) -- (A) -- cycle;
                    \foreach \x in {-3, -2.7, ..., 0} { \draw[\colQ, thick] (\x, -1.5) -- (\x, 1.5); }
                \end{scope}
                \draw[\colQ, dotted, thick] (A) to[bend right=20]  (B2);
                \draw[\colW, thick] (A) to[bend left=50] node[above=1mm, text=\colW] {\scriptsize${}_aW_b$}  (B2);
                \node at (-1.5, 0.95) {\scriptsize$g$};
                \begin{scope}
                    \clip (B) to[bend right=45,looseness=1.5] (C) -- (C) to[out=-120, in=0] (B') -- (B) -- cycle;
                    \foreach \x in {0.3, 0.6, ..., 3} { \draw[\colY, thick] (\x, 0) -- (\x, -2.5); }
                \end{scope}
                \draw[\colY, thick] (B') to[out=0, in=-120] (C);
                \draw[\colY,thick,dotted] (B) to[bend right=45,looseness=1.5] (C);
                \draw[\colR, dotted, thick] (B) -- (C);
                \node at (1.5, -0.4) {\scriptsize$h$};
                \begin{scope}
                    \clip (B2) to[bend right=20] (C) -- (B) -- cycle;
                    \foreach \x in {0.3, 0.6, ..., 3} { \draw[\colR, thick] (\x, -1.5) -- (\x, 1.5); }
                \end{scope}
                \draw[\colR, dotted, thick] (B2) to[bend right=20] (C);
                \draw[\colZ, thick] (B2) to[bend left=50] node[pos=.6,right=1mm, text=\colZ] {\scriptsize${}_bZ_c$} (C);
                \node at (1.5, 0.95) {\scriptsize$k$};
                \fill[\cola] (A) circle (1.5pt) node[left] {\scriptsize$\fcj{a}$};
                \fill[\colb] (B) circle (1.5pt);
                \fill[\colc] (C) circle (1.5pt) node[right] {\scriptsize$c$};
                \fill[\colb] (B2) circle (1.5pt);
                \fill[\colb] (B') circle (1.5pt);
                \draw[\colb,thick] (B2) -- (B) -- (B');
                \draw[->, thick] (1,1.75) node[above] {\scriptsize$({}_aQ \otimes_b R_c) \times I$} to[in=60,out=-90] (0.4,0.5);
            \end{tkz}
            \\
            &=
            \begin{tkz}[scale=0.7]
                \coordinate (A) at (-3, 0);
                \coordinate (B) at (0, 0);
                \coordinate (C) at (3, 0);
                \coordinate (B') at (0, -2);
                \fill[\filg] (A) to[bend left=45] (B) to[bend left=45] (A);
                \fill[\filk] (B) to[bend left=45] (C) to[bend left=45] (B);
                \begin{scope}
                    \clip (A) to[bend right=45] (B) -- (B') to[out=180, in=-60] (A) -- cycle;
                    \foreach \x in {-3, -2.7, ..., 0} { \draw[\colQ, thick] (\x, 0) -- (\x, -2.5); }
                \end{scope}
                \draw[\colQ, thick] (A) to[out=-60, in=180] node[below] {\scriptsize$\fcj{{}_aQ_b}$} (B');
                \draw[\colQ, dotted, thick] (A) to[bend right=45] (B);
                \draw[\colW, thick] (A) to[bend left=45] node[above] {\scriptsize${}_aW_b$} (B);
                \node at (-1.5, 0) {\scriptsize$g$};
                \begin{scope}
                    \clip (B) to[bend right=45] (C) -- (C) to[out=-120, in=0] (B') -- (B) -- cycle;
                    \foreach \x in {0.3, 0.6, ..., 3} { \draw[\colR, thick] (\x, 0) -- (\x, -2.5); }
                \end{scope}
                \draw[\colR,thick] (B') to[out=0, in=-120] node[pos=.65,below] {\scriptsize$\fcj{{}_bR_c}$} (C);
                \draw[\colR,dotted, thick] (B) to[bend right=45] (C);
                \draw[\colZ, thick] (B) to[bend left=45] node[above] {\scriptsize${}_bZ_c$} (C);
                \node at (1.5, 0) {\scriptsize$k$};
                \draw[\colb, thick] (B) -- (B');
                \fill[\cola] (A) circle (1.5pt) node[left] {\scriptsize$\fcj{a}$};
                \fill[\colb] (B) circle (1.5pt);
                \fill[\colc] (C) circle (1.5pt) node[right] {\scriptsize$c$};
                \fill[\colb] (B') circle (1.5pt);
                \draw[<-, thick] (0.75,-1.15) to[in=90,out=-45] (1.5,-1.95) node[below] {\scriptsize$({}_aQ \otimes_b R_c) \times I$};
            \end{tkz}
            \!\circ\!
            \begin{tkz}[scale=0.7]
                \coordinate (A) at (-3, 0);
                \coordinate (B) at (0, 0);
                \coordinate (C) at (3, 0);
                \coordinate (B') at (0, -2);
                \fill[\filf] (A) to[bend left=45] (B) to[bend left=45] (A);
                \fill[\filh] (B) to[bend left=45] (C) to[bend left=45] (B);
                \begin{scope}
                    \clip (A) to[bend right=45] (B) -- (B') to[out=180, in=-60] (A) -- cycle;
                    \foreach \x in {-3, -2.7, ..., 0} { \draw[\colX, thick] (\x, 0) -- (\x, -2.5); }
                \end{scope}
                \draw[\colX, thick] (A) to[out=-60, in=180] node[below] {\scriptsize$\fcj{{}_aX_b}$} (B');
                \draw[\colX,dotted,thick] (A) to[bend right=45] (B);
                \draw[\colQ, thick] (A) to[bend left=45] node[above, text=\colQ] {\scriptsize${}_aQ_b$} (B);
                \node at (-1.5, 0) {\scriptsize$f$};
                \begin{scope}
                    \clip (B) to[bend right=45] (C) -- (C) to[out=-120, in=0] (B') -- (B) -- cycle;
                    \foreach \x in {0.3, 0.6, ..., 3} { \draw[\colY, thick] (\x, 0) -- (\x, -2.5); }
                \end{scope}
                \draw[\colY, thick] (B') to[out=0, in=-120] node[below, text=\colY] {\scriptsize$\fcj{{}_bY_c}$} (C);
                \draw[\colY,dotted,thick] (B) to[bend right=45] (C);
                \draw[\colR, thick] (B) to[bend left=45] node[above, text=\colR] {\scriptsize${}_bR_c$} (C);
                \node at (1.5, 0) {\scriptsize$h$};
                \draw[\colb, thick] (B) -- (B');
                \fill[\cola] (A) circle (1.5pt) node[left] {\scriptsize$\fcj{a}$};
                \fill[\colb] (B) circle (1.5pt);
                \fill[\colc] (C) circle (1.5pt) node[right] {\scriptsize$c$};
                \fill[\colb] (B') circle (1.5pt);
            \end{tkz}
            = (g \otimes k) \circ (f \otimes h).
            \qedhere
        \end{align*}
\end{pf}

\subsubsection{Rigidity of \texorpdfstring{$\fX_\C$}{XC}}
Next we show that $\fX_\C$ is rigid and admits an adjoint functor that induces a strict pivotal structure.

\begin{construction}
    \label{construction:skein adjoint functor}
    On 1-morphisms ${}_aX_b\in\fX_\C(a\to b)=\fld[1]{\C}{D^1}[\fcj{a}\amalg b]$, define ${}_bX_a^\vee\in\fX_\C(b\to a)$ by $X^\vee\coloneq\fcj{\iota^{(1)}_*X}$ where $\iota^{(1)}\colon D^1\to \orev{D^1}$ is the orientation-\emph{preserving} homeomorphism $t\mapsto -t$. Because $\iota^{(1)}$ and $\orev{\,\cdot\,}$ are commuting involutions, $X^{\vee\vee}=X$. Define the associated evaluation and coevaluation 2-morphisms respectively by
    \begin{align*}
            \begin{tkz}[scale=0.7]
                \def\tkzWidth{2.5}
                \def\tkzHeight{1.0}
                \def\tkzAngle{30}
                \def\tkzVertLooseness{1.0}
                \def\tkzPointLooseness{.5}
                \begin{scope}
                    \draw[\colX, thick] (0,0) to[out=\tkzAngle,in=180, out looseness=\tkzPointLooseness, in looseness=\tkzVertLooseness] node[pos=0.4, above left=-3pt,\colX] {\scriptsize${}_bX_a^\vee$} (\tkzWidth,\tkzHeight);
                    \draw[\colX, thick] (0,0) to[out=-\tkzAngle,in=180, out looseness=\tkzPointLooseness, in looseness=\tkzVertLooseness] node[pos=0.4, below left=-2pt,\colX] {\scriptsize${}_aX_b$} (\tkzWidth,-\tkzHeight);
                    \draw[\colb, thick] (\tkzWidth,-\tkzHeight)--node[right]{\scriptsize$1_b$} (\tkzWidth,\tkzHeight);
                    \begin{scope}
                        \clip (0,0) to[out=\tkzAngle,in=180, out looseness=\tkzPointLooseness, in looseness=\tkzVertLooseness] (\tkzWidth,\tkzHeight) -- (\tkzWidth,-\tkzHeight) to[out=180,in=-\tkzAngle, out looseness=\tkzVertLooseness, in looseness=\tkzPointLooseness] cycle;
                        \foreach \x in {0.2,0.4,...,2.8} {
                                \draw[\colX, thick] (\x,-1.5) -- (\x,1.5);
                            }
                    \end{scope}
                \end{scope}
                \fill[\cola] (0,0) circle (2pt) node[left, \cola] {\scriptsize $a$};
                \fill[\colb] (\tkzWidth,\tkzHeight) circle (2pt) node[right, \colb] {\scriptsize $b$};
                \fill[\colb] (\tkzWidth,-\tkzHeight) circle (2pt) node[right, \colb] {\scriptsize $b$};
            \end{tkz}
             & \xrightarrow{\text{\scriptsize{}reparameterize}}
            \begin{tkz}[yscale=-1,xscale=-1,scale=0.7]
                \def\tkzWidth{4}
                \def\tkzMidX{2}
                \def\tkzHeight{1.1}
                \def\tkzAngleOut{45}
                \def\tkzAngleIn{135}
                \begin{scope}
                    \clip (0,0) to[out=-\tkzAngleOut,in=-\tkzAngleIn,looseness=1.25] (\tkzWidth,0) to[out=\tkzAngleIn,in=0] (\tkzMidX,\tkzHeight) to[out=180,in=\tkzAngleOut] cycle;
                    \foreach \t in {0.1,0.2,...,0.9} {
                            \coordinate (L) at ({\t*\tkzMidX - 0.2}, {\t*\tkzHeight + 0.4});
                            \coordinate (R) at ({\tkzWidth-\t*\tkzMidX + 0.2}, {\t*\tkzHeight + 0.4});
                            \draw[\colX, thick] (L) to[bend right=60,looseness=1.25] (R);
                        }
                \end{scope}
                \draw[\colb, thick] (0,0) to[bend right=45,looseness=1.25] node[above] {\scriptsize$1_b$} (\tkzWidth,0);
                \draw[\colX, thick] (0,0) to[out=\tkzAngleOut,in=180] node[below right] {\scriptsize ${}_aX_b$} (\tkzMidX,\tkzHeight);
                \draw[\colX, thick] (\tkzMidX,\tkzHeight) to[out=0,in=\tkzAngleIn] node[below left] {\scriptsize ${}_bX_a^\vee$} (\tkzWidth,0);
                \fill[\colb] (0,0) circle (2pt) node[right, \colb] {\scriptsize $b$};
                \fill[\cola] (\tkzMidX,\tkzHeight) circle (2pt) node[below, \cola] {\scriptsize $a$};
                \fill[\colb] (\tkzWidth,0) circle (2pt) node[left, \colb] {\scriptsize $b$};
            \end{tkz}
            \eqcolon \ev_X \in \fX_\C({}_bX^\vee \otimes_a X_b \Rightarrow 1_b),
            \\
            \begin{tkz}[scale=0.7]
                \def\tkzWidth{2.5}
                \def\tkzHeight{1.0}
                \def\tkzAngle{150}
                \def\tkzVertLooseness{1.0}
                \def\tkzPointLooseness{.5}
                \begin{scope}
                    \draw[\colX, thick] (0,\tkzHeight) to[out=0,in=\tkzAngle, out looseness=\tkzVertLooseness, in looseness=\tkzPointLooseness] node[pos=0.6, above right=-2pt] {\scriptsize${}_aX_b$} (\tkzWidth,0);
                    \draw[\colX, thick] (0,-\tkzHeight) to[out=0,in=-\tkzAngle, out looseness=\tkzVertLooseness, in looseness=\tkzPointLooseness] node[pos=0.6, below right=-3pt] {\scriptsize${}_bX_a^\vee$} (\tkzWidth,0);
                    \draw[\cola,thick] (0,-\tkzHeight) --node[left]{\scriptsize$1_a$} (0,\tkzHeight);
                    \begin{scope}
                        \clip (\tkzWidth,0) to[out=\tkzAngle,in=0, out looseness=\tkzPointLooseness, in looseness=\tkzVertLooseness] (0,\tkzHeight) -- (0,-\tkzHeight) to[out=0,in=-\tkzAngle, out looseness=\tkzVertLooseness, in looseness=\tkzPointLooseness] cycle;
                        \foreach \x in {0.2, 0.4, ..., 2.8} {
                                \draw[\colX, thick] (\x,-1.5) -- (\x,1.5);
                            }
                    \end{scope}
                \end{scope}
                \fill[\colb] (\tkzWidth,0) circle (2pt) node[right,\colb] {\scriptsize$b$};
                \fill[\cola] (0,\tkzHeight) circle (2pt) node[left,\cola] {\scriptsize$a$};
                \fill[\cola] (0,-\tkzHeight) circle (2pt) node[left,\cola] {\scriptsize$a$};
            \end{tkz}
             & \xrightarrow{\text{\scriptsize{}reparameterize}}
            \begin{tkz}[scale=0.7]
                \def\tkzWidth{4}
                \def\tkzMidX{2}
                \def\tkzHeight{1.1}
                \def\tkzAngleOut{45}
                \def\tkzAngleIn{135}
                \begin{scope}
                    \clip (0,0) to[out=-\tkzAngleOut,in=-\tkzAngleIn,looseness=1.25] (\tkzWidth,0) to[out=\tkzAngleIn,in=0] (\tkzMidX,\tkzHeight) to[out=180,in=\tkzAngleOut] cycle;
                    \foreach \t in {0.1,0.2,...,0.9} {
                            \coordinate (L) at ({\t*\tkzMidX - 0.2}, {\t*\tkzHeight + 0.4});
                            \coordinate (R) at ({\tkzWidth-\t*\tkzMidX + 0.2}, {\t*\tkzHeight + 0.4});
                            \draw[\colX, thick] (L) to[bend right=60,looseness=1.25] (R);
                        }
                \end{scope}
                \draw[\cola, thick] (0,0) to[bend right=45,looseness=1.25] node[below] {\scriptsize$1_a$} (\tkzWidth,0);
                \draw[\colX, thick] (0,0) to[out=\tkzAngleOut,in=180] node[above left] {\scriptsize ${}_aX_b$} (\tkzMidX,\tkzHeight);
                \draw[\colX, thick] (\tkzMidX,\tkzHeight) to[out=0,in=\tkzAngleIn] node[above right] {\scriptsize ${}_bX_a^\vee$} (\tkzWidth,0);
                \fill[\cola] (0,0) circle (2pt) node[left, \cola] {\scriptsize $a$};
                \fill[\colb] (\tkzMidX,\tkzHeight) circle (2pt) node[above, \colb] {\scriptsize $b$};
                \fill[\cola] (\tkzWidth,0) circle (2pt) node[right, \cola] {\scriptsize $a$};
            \end{tkz}
            \eqcolon \coev_X \in \fX_\C(1_a \Rightarrow {}_aX\otimes_b X_a^\vee).
        \end{align*}
\end{construction}

\begin{proposition}
    \label{skein 2-categories are rigid}
    For a disk-like 2-category $\C$, the linear 2-category $\fX_\C$ is rigid.
\end{proposition}

\begin{pf}
    The caps and cups of \cref{construction:skein adjoint functor} are by extended-isotopy invariance easily seen to satisfy the zig-zag/snake equations. Thus the induced adjoint functor $\vee\colon\fX_\C\to\fX_\C^{1\op,2\op}$ is given on 2-morphisms $f\in\fX_\C({}_aX_b\Rightarrow {}_aY_b)$ by
    \[
        \begin{tkz}[scale=0.75]
            \def\tkzxdst{5}
            \def\tkzbigth{65}
            \fill[\filf](0,0)to[out=\tkzbigth,in=180-\tkzbigth](\tkzxdst,0)to[out=180+\tkzbigth,in=-\tkzbigth](0,0);
            \draw[thick,\colX](0,0)to[out=-\tkzbigth,in=180+\tkzbigth](\tkzxdst,0);
            \draw[thick,\colY](0,0)to[out=\tkzbigth,in=180-\tkzbigth](\tkzxdst,0);

            \fill[\cola] (0,0) circle (1.5pt) node [left]{\scriptsize$\fcj{a}$};
            \fill[\colb] (\tkzxdst,0) circle (1.5pt) node [right]{\scriptsize$b$};

            \node at (0.5*\tkzxdst,0) {$f$};
        \end{tkz}
        \mapsto
        \begin{tkz}[scale=0.8,xscale=-1,yscale=-1]
            \def\tkzxdst{5}
            \def\tkzbigth{65}
            \def\tkzmd{.395}
            \def\tkzsepr{2.5}
            \draw[thick,\colX](0,0)to[out=-\tkzbigth,in=180+\tkzbigth](\tkzxdst,0);
            \draw[thick,\colY](0,0)to[out=\tkzbigth,in=180-\tkzbigth](\tkzxdst,0);
            \fill[\cola] (0,0) circle (1.5pt) node [right]{\scriptsize$a$};
            \fill[\colb] (\tkzxdst,0) circle (1.5pt) node [left]{\scriptsize$\fcj{b}$};
            \clip(0,0)to[out=\tkzbigth,in=180-\tkzbigth](\tkzxdst,0)to[out=180+\tkzbigth,in=-\tkzbigth](0,0);
            \fill[\filf](\tkzmd*\tkzxdst,0)
            to[out=0.8*\tkzbigth,in=180-0.8*\tkzbigth] ({\tkzxdst-\tkzmd*\tkzxdst},0)
            to[out=180+0.8*\tkzbigth,in=-0.8*\tkzbigth] (\tkzmd*\tkzxdst,0);
            \draw[thick,dotted,\colY](\tkzmd*\tkzxdst,0)
            to[out=-0.8*\tkzbigth,in=180+0.8*\tkzbigth]
            coordinate[pos=0] (bottom0)
            coordinate[pos=0.1] (bottom1) coordinate[pos=0.2] (bottom2)
            coordinate[pos=0.3] (bottom3) coordinate[pos=0.4] (bottom4)
            coordinate[pos=0.5] (bottom5) coordinate[pos=0.6] (bottom6)
            coordinate[pos=0.7] (bottom7) coordinate[pos=0.8] (bottom8)
            coordinate[pos=0.9] (bottom9) coordinate[pos=1.0] (bottom10)
            ({\tkzxdst-\tkzmd*\tkzxdst},0);
            \draw[thick,dotted,\colX](\tkzmd*\tkzxdst,0)
            to[out=0.8*\tkzbigth,in=180-0.8*\tkzbigth]
            coordinate[pos=0] (top0)
            coordinate[pos=0.1] (top1) coordinate[pos=0.2] (top2)
            coordinate[pos=0.3] (top3) coordinate[pos=0.4] (top4)
            coordinate[pos=0.5] (top5) coordinate[pos=0.6] (top6)
            coordinate[pos=0.7] (top7) coordinate[pos=0.8] (top8)
            coordinate[pos=0.9] (top9) coordinate[pos=1.0] (top10)
            ({\tkzxdst-\tkzmd*\tkzxdst},0);
            \foreach\x in {1,...,9} {
                    \draw[thick,\colY] (bottom\x)
                    to[out=-90,in=-90,looseness=1.4] ++({\tkzsepr*(1-0.1*\x)*((1-\tkzmd)*\tkzxdst-\tkzmd*\tkzxdst)},0)
                    to[out=90,in=-90] ++({1-.175*\x}, 1.6);
                    \draw[thick,\colX] (top\x)
                    to[out=90,in=90,looseness=1.4] ++({-\tkzsepr*(0.1*\x)*((1-\tkzmd)*\tkzxdst-\tkzmd*\tkzxdst)},0)
                    to[out=-90,in=90] ++({0.75-.175*\x}, -1.6);
                }
            \node at (0.5*\tkzxdst,0) {\scriptsize$f$};
            \draw[thick,\colb] (top0)
            to[out=-90,in=-90+-35,out looseness=1.1,in looseness=1.7] ++(0.55*\tkzxdst,-.05*\tkzxdst) to[in=200,out=-90+-35] (\tkzxdst,0);
            \draw[thick,\cola] (bottom10)
            to[out=90,in=90-35,out looseness=1.1,in looseness=1.7] ++(-0.55*\tkzxdst,.05*\tkzxdst) to[in=20,out=90-35] (0,0);
            \fill[\colb] (\tkzmd*\tkzxdst,0) circle (1pt);
            \fill[\cola] ({\tkzxdst-\tkzmd*\tkzxdst},0) circle (1pt);
        \end{tkz}
        \in\fX_\C({}_bY_a^\vee\Rightarrow {}_bX_a^\vee).
        \qedhere
    \]
\end{pf}

\subsubsection{\texorpdfstring{$\fX_\C$}{XC} as a pivotal \texorpdfstring{dagger}{dagger} 2-category}
Next we define the dagger structure on $\fX_\C$ and show it is compatible with the adjoint functor associated to the adjoint data constructed in \cref{construction:skein adjoint functor}.

\begin{proposition}
    \label{daggers for skeins ii}
    For a disk-like 2-category $\C$, the 2-category $\fX_\C$ has a canonical dagger structure.
\end{proposition}

\begin{pf}
    Define the dagger structure $\dag\colon\fX_\C({}_aX_b\Rightarrow{}_aY_b)\to\fX_\C({}_aY_b\Rightarrow{}_aX_b)$ by $\alpha^\dag\coloneq\fcj{\iota^{(2)}_*\alpha}$ where $\iota^{(2)}\colon D^2\to\orev{D^2}$ is the orientation-preserving homeomorphism whose underlying map is the reflection about the horizontal axis. Observe that $\dag$ is conjugate-linear since $\iota^{(2)}_*$ is linear and $\fcj{\,\cdot\,}$ is conjugate-linear. Then \ref{2daga} and \ref{2dagb} are immediate. Moreover, \ref{2dagd} holds since the coheretors are unitary with respect to $\dag$, which is immediate by the diagrammatic description of the associator and unitors in \cite[Appendix C.2]{MW12} and invariance under isotopy-rel-boundary. Finally, to see \ref{2dagc}, observe that for $f\in\fX_\C({}_aX_b\Rightarrow{}_aW_b)$ and $g\in\fX_\C({}_bY_c\Rightarrow{}_bZ_c)$,
        \begin{align*}
            (f&\otimes g)^\dag=
            \left(
            \begin{tkz}[scale=0.7]
                \coordinate (A) at (-3, 0);
                \coordinate (B) at (0, 0);
                \coordinate (C) at (3, 0);
                \coordinate (B') at (0, -2);
                \fill[\filf] (A) to[bend left=45] (B) to[bend left=45] (A);
                \fill[\filg] (B) to[bend left=45] (C) to[bend left=45] (B);
                \begin{scope}
                    \clip (A) to[bend right=45] (B) -- (B') to[out=180, in=-60] (A) -- cycle;
                    \foreach \x in {-3, -2.7, ..., 0} { \draw[\colX, thick] (\x, 0) -- (\x, -2.5); }
                \end{scope}
                \draw[\colX, thick] (A) to[out=-60, in=180] node[pos=.65,below] {\scriptsize$\fcj{{}_aX_b}$} (B');
                \draw[\colX,dotted,thick] (A) to[bend right=45] (B);
                \draw[\colW, thick] (A) to[bend left=45] node[pos=0.5,above, text=\colW] {\scriptsize${}_aW_b$} (B);
                \node at (-1.5, 0) {\scriptsize$f$};
                \begin{scope}
                    \clip (B) to[bend right=45] (C) -- (C) to[out=-120, in=0] (B') -- (B) -- cycle;
                    \foreach \x in {0.3, 0.6, ..., 3} { \draw[\colY, thick] (\x, 0) -- (\x, -2.5); }
                \end{scope}
                \draw[\colY, thick] (B') to[out=0, in=-120] node[pos=.65,below, text=\colY] {\scriptsize$\fcj{{}_bY_c}$} (C);
                \draw[\colY,dotted,thick] (B) to[bend right=45] (C);
                \draw[\colZ, thick] (B) to[bend left=45] node[above, text=\colZ] {\scriptsize${}_bZ_c$} (C);
                \node at (1.5, 0) {\scriptsize$g$};
                \draw[\colb, thick] (B) -- (B');
                \fill[\cola] (A) circle (1.5pt) node[left] {\scriptsize$\fcj{a}$};
                \fill[\colb] (B) circle (1.5pt);
                \fill[\colc] (C) circle (1.5pt) node[right] {\scriptsize$c$};
                \fill[\colb] (B') circle (1.5pt);
                \draw[<-, thick] (0.75,-1.15) to[in=90,out=-45] (1.5,-1.95) node[below] {\scriptsize$({}_aX \otimes_b Y_c) \times I$};
            \end{tkz}
            \right)^\dag=
            \begin{tkz}[scale=0.8,yscale=-1]
                \coordinate (A) at (-3, 0);
                \coordinate (B) at (0, 0);
                \coordinate (C) at (3, 0);
                \coordinate (B') at (0, -2);
                \fill[\filf] (A) to[bend left=45] (B) to[bend left=45] (A);
                \fill[\filg] (B) to[bend left=45] (C) to[bend left=45] (B);
                \begin{scope}
                    \clip (A) to[bend right=45] (B) -- (B') to[out=180, in=-60] (A) -- cycle;
                    \foreach \x in {-3, -2.7, ..., 0} { \draw[\colX, thick] (\x, 0) -- (\x, -2.5); }
                \end{scope}
                \draw[\colX, thick] (A) to[out=-60, in=180] node[pos=.65,above left] {\scriptsize$\fcj{{}_aX_b}$} (B');
                \draw[\colX,dotted,thick] (A) to[bend right=45] (B);
                \draw[\colW, thick] (A) to[bend left=45] node[below, text=\colW] {\scriptsize${}_aW_b$} (B);
                \node at (-1.5, 0) {\scriptsize$f^\dag$};
                \begin{scope}
                    \clip (B) to[bend right=45] (C) -- (C) to[out=-120, in=0] (B') -- (B) -- cycle;
                    \foreach \x in {0.3, 0.6, ..., 3} { \draw[\colY, thick] (\x, 0) -- (\x, -2.5); }
                \end{scope}
                \draw[\colY, thick] (B') to[out=0, in=-120] node[pos=.65,above right, text=\colY] {\scriptsize$\fcj{{}_bY_c}$} (C);
                \draw[\colY,dotted,thick] (B) to[bend right=45] (C);
                \draw[\colZ, thick] (B) to[bend left=45] node[below, text=\colZ] {\scriptsize${}_bZ_c$} (C);
                \node at (1.5, 0) {\scriptsize$g^\dag$};
                \draw[\colb, thick] (B) -- (B');
                \fill[\cola] (A) circle (1.5pt) node[left] {\scriptsize$\fcj{a}$};
                \fill[\colb] (B) circle (1.5pt);
                \fill[\colc] (C) circle (1.5pt) node[right] {\scriptsize$c$};
                \fill[\colb] (B') circle (1.5pt);
                \draw[<-, thick] (0.75,-1.15) to[in=90,out=-45] (1.5,-1.95) node[above] {\scriptsize$({}_aX \otimes_b Y_c) \times I$};
            \end{tkz}
            \\
            &=
            \begin{tkz}[scale=0.8,yscale=-1]
                \coordinate (A) at (-3, 0);
                \coordinate (B) at (0, 0);
                \coordinate (C) at (3, 0);
                \coordinate (B') at (0, -2);
                \coordinate (B'') at (0, 2);
                \fill[\filf] (A) to[bend left=45] (B) to[bend left=45] (A);
                \fill[\filg] (B) to[bend left=45] (C) to[bend left=45] (B);
                \begin{scope}
                    \clip (A) to[bend right=45] (B) -- (B') to[out=180, in=-60] (A) -- cycle;
                    \foreach \x in {-3, -2.7, ..., 0} { \draw[\colX, thick] (\x, 0) -- (\x, -2.5); }
                \end{scope}
                \draw[\colX, thick] (A) to[out=-60, in=180] node[pos=.65,above left] {\scriptsize$\fcj{{}_aX_b}$} (B');
                \draw[\colX,dotted,thick] (A) to[bend right=45] (B);
                \begin{scope}
                    \clip (A) to[bend left=45] (B) -- (B'') to[out=180, in=60] (A) -- cycle;
                    \foreach \x in {-3, -2.7, ..., 0} { \draw[\colW, thick] (\x, 0) -- (\x, 2.5); }
                \end{scope}
                \draw[\colW, thick] (A) to[out=60, in=180] node[pos=.5,below left, text=\colW] {\scriptsize$\fcj{{}_aW_b}$} (B'');
                \draw[\colW, thick, dotted] (A) to[bend left=45]
                (B);
                \node at (-1.5, 0) {\scriptsize$f^\dag$};
                \begin{scope}
                    \clip (B) to[bend right=45] (C) -- (C) to[out=-120, in=0] (B') -- (B) -- cycle;
                    \foreach \x in {0.3, 0.6, ..., 3} { \draw[\colY, thick] (\x, 0) -- (\x, -2.5); }
                \end{scope}
                \draw[\colY, thick] (B') to[out=0, in=-120] node[pos=.65,above right, text=\colY] {\scriptsize$\fcj{{}_bY_c}$} (C);
                \draw[\colY,dotted,thick] (B) to[bend right=45] (C);
                \begin{scope}
                    \clip (B) to[bend left=45] (C) -- (C) to[out=120, in=0] (B'') -- (B) -- cycle;
                    \foreach \x in {0.3, 0.6, ..., 3} { \draw[\colZ, thick] (\x, 0) -- (\x, 2.5); }
                \end{scope}
                \draw[\colZ, thick] (B'') to[out=0, in=120] node[pos=.65,below right, text=\colZ] {\scriptsize$\fcj{{}_bZ_c}$} (C);
                \draw[\colZ, thick, dotted] (B) to[bend left=45]
                (C);
                \node at (1.5, 0) {\scriptsize$g^\dag$};
                \draw[\colb, thick] (B) -- (B');
                \draw[\colb, thick] (B) -- (B'');
                \fill[\cola] (A) circle (1.5pt) node[left] {\scriptsize$\fcj{a}$};
                \fill[\colb] (B) circle (1.5pt);
                \fill[\colc] (C) circle (1.5pt) node[right] {\scriptsize$c$};
                \fill[\colb] (B') circle (1.5pt);
                \fill[\colb] (B'') circle (1.5pt);
                \draw[<-, thick] (0.75,-1.15) to[in=90,out=-45] (1.5,-1.95) node[above] {\scriptsize$({}_aX \otimes_b Y_c) \times I$};
                \draw[<-, thick] (0.75,1.15) to[in=-90,out=45] (1.5,1.95) node[below] {\scriptsize$({}_aW \otimes_b Z_c) \times I$};
            \end{tkz}
            =
            \begin{tkz}[scale=0.8,yscale=-1]
                \coordinate (A) at (-3, 0);
                \coordinate (B) at (0, 0);
                \coordinate (C) at (3, 0);
                \coordinate (B'') at (0, 2);
                \fill[\filf] (A) to[bend left=45] (B) to[bend left=45] (A);
                \fill[\filg] (B) to[bend left=45] (C) to[bend left=45] (B);
                \draw[\colX, thick] (A) to[bend right=45] node[above] {\scriptsize$\fcj{{}_aX_b}$} (B);
                \begin{scope}
                    \clip (A) to[bend left=45] (B) -- (B'') to[out=180, in=60] (A) -- cycle;
                    \foreach \x in {-3, -2.7, ..., 0} { \draw[\colW, thick] (\x, 0) -- (\x, 2.5); }
                \end{scope}
                \draw[\colW, thick] (A) to[out=60, in=180] node[pos=.5,below left, text=\colW] {\scriptsize$\fcj{{}_aW_b}$} (B'');
                \draw[\colW, thick, dotted] (A) to[bend left=45]
                (B);
                \node at (-1.5, 0) {\scriptsize$f^\dag$};
                \draw[\colY, thick] (B) to[bend right=45] node[above, text=\colY] {\scriptsize$\fcj{{}_bY_c}$} (C);
                \begin{scope}
                    \clip (B) to[bend left=45] (C) -- (C) to[out=120, in=0] (B'') -- (B) -- cycle;
                    \foreach \x in {0.3, 0.6, ..., 3} { \draw[\colZ, thick] (\x, 0) -- (\x, 2.5); }
                \end{scope}
                \draw[\colZ, thick] (B'') to[out=0, in=120] node[pos=.65,below right, text=\colZ] {\scriptsize$\fcj{{}_bZ_c}$} (C);
                \draw[\colZ, thick, dotted] (B) to[bend left=45]
                (C);
                \node at (1.5, 0) {\scriptsize$g^\dag$};
                \draw[\colb, thick] (B) -- (B'');
                \fill[\cola] (A) circle (1.5pt) node[left] {\scriptsize$\fcj{a}$};
                \fill[\colb] (B) circle (1.5pt);
                \fill[\colc] (C) circle (1.5pt) node[right] {\scriptsize$c$};
                \fill[\colb] (B'') circle (1.5pt);
                \draw[<-, thick] (0.75,1.15) to[in=-90,out=45] (1.5,1.95) node[below] {\scriptsize$({}_aW \otimes_b Z_c) \times I$};
            \end{tkz}
            =
            f^\dag\otimes g^\dag.
            \qedhere
        \end{align*}
\end{pf}

\begin{proposition}
    \label{skein 2-categories are unitary planar pivotal}
    For a disk-like 2-category $\C$, the rigid 2-category $\fX_\C$ has a canonical UAF $\vee$, and $\vee$ induces a strict pivotal structure on $\fX_\C$. In particular, $\fX_\C$ is a pivotal dagger 2-category whose pivotal structure is strict.
\end{proposition}

\begin{pf}
    \itemstep{$\vee\dag=\dag\vee$.}
    The fact that $\vee\dag=\dag\vee$ on 2-morphisms in $\fX_\C$ follows from the invariance of 2-fields under isotopy-rel-boundary:
    \begin{equation}
        \label{eq:2skel vee dag}
        \begin{multlined}
            f^{\vee\dag}
            =
            \begin{tkz}[scale=1,xscale=-1]
                \def\tkzxdst{5}
                \def\tkzbigth{65}
                \def\tkzmd{.395}
                \def\tkzsepr{2.5}
                \fill[\filf](0,0)to[out=\tkzbigth,in=180-\tkzbigth](\tkzxdst,0)to[out=180+\tkzbigth,in=-\tkzbigth](0,0);
                \draw[thick,\colX](0,0)to[out=-\tkzbigth,in=180+\tkzbigth](\tkzxdst,0);
                \draw[thick,\colY](0,0)to[out=\tkzbigth,in=180-\tkzbigth](\tkzxdst,0);
                \fill[\cola] (0,0) circle (1.5pt) node [right]{\scriptsize$a$};
                \fill[\colb] (\tkzxdst,0) circle (1.5pt) node [left]{\scriptsize$\fcj{b}$};
                \clip(0,0)to[out=\tkzbigth,in=180-\tkzbigth](\tkzxdst,0)to[out=180+\tkzbigth,in=-\tkzbigth](0,0);
                \draw[thick,dotted,\colY](\tkzmd*\tkzxdst,0)
                to[out=-0.8*\tkzbigth,in=180+0.8*\tkzbigth]
                coordinate[pos=0] (bottom0)
                coordinate[pos=0.1] (bottom1) coordinate[pos=0.2] (bottom2)
                coordinate[pos=0.3] (bottom3) coordinate[pos=0.4] (bottom4)
                coordinate[pos=0.5] (bottom5) coordinate[pos=0.6] (bottom6)
                coordinate[pos=0.7] (bottom7) coordinate[pos=0.8] (bottom8)
                coordinate[pos=0.9] (bottom9) coordinate[pos=1.0] (bottom10)
                ({\tkzxdst-\tkzmd*\tkzxdst},0);
                \draw[thick,dotted,\colX](\tkzmd*\tkzxdst,0)
                to[out=0.8*\tkzbigth,in=180-0.8*\tkzbigth]
                coordinate[pos=0] (top0)
                coordinate[pos=0.1] (top1) coordinate[pos=0.2] (top2)
                coordinate[pos=0.3] (top3) coordinate[pos=0.4] (top4)
                coordinate[pos=0.5] (top5) coordinate[pos=0.6] (top6)
                coordinate[pos=0.7] (top7) coordinate[pos=0.8] (top8)
                coordinate[pos=0.9] (top9) coordinate[pos=1.0] (top10)
                ({\tkzxdst-\tkzmd*\tkzxdst},0);
                \foreach\x in {1,...,9} {
                        \draw[thick,\colY] (bottom\x)
                        to[out=-90,in=-90,looseness=1.4] ++({\tkzsepr*(1-0.1*\x)*((1-\tkzmd)*\tkzxdst-\tkzmd*\tkzxdst)},0)
                        to[out=90,in=-90] ++({1-.175*\x}, 1.6);
                        \draw[thick,\colX] (top\x)
                        to[out=90,in=90,looseness=1.4] ++({-\tkzsepr*(0.1*\x)*((1-\tkzmd)*\tkzxdst-\tkzmd*\tkzxdst)},0)
                        to[out=-90,in=90] ++({0.75-.175*\x}, -1.6);
                    }
                \node at (0.5*\tkzxdst,0) {\scriptsize$f^\dag$};
                \draw[thick,\colb] (top0)
                to[out=-90,in=-90+-35,out looseness=1.1,in looseness=1.7] ++(0.55*\tkzxdst,-.05*\tkzxdst) to[in=200,out=-90+-35] (\tkzxdst,0);
                \draw[thick,\cola] (bottom10)
                to[out=90,in=90-35,out looseness=1.1,in looseness=1.7] ++(-0.55*\tkzxdst,.05*\tkzxdst) to[in=20,out=90-35] (0,0);
                \fill[\colb] (\tkzmd*\tkzxdst,0) circle (1pt);
                \fill[\cola] ({\tkzxdst-\tkzmd*\tkzxdst},0) circle (1pt);
            \end{tkz}
            =
            \begin{tkz}[scale=1]
                \def\tkzxdst{5}
                \def\tkzbigth{65}
                \def\tkzmd{.395}
                \def\tkzsepr{2.5}
                \fill[\filf](0,0)to[out=\tkzbigth,in=180-\tkzbigth](\tkzxdst,0)to[out=180+\tkzbigth,in=-\tkzbigth](0,0);
                \draw[thick,\colX](0,0)to[out=-\tkzbigth,in=180+\tkzbigth](\tkzxdst,0);
                \draw[thick,\colY](0,0)to[out=\tkzbigth,in=180-\tkzbigth](\tkzxdst,0);
                \fill[\colb] (0,0) circle (1.5pt) node [left]{\scriptsize$\fcj{b}$};
                \fill[\cola] (\tkzxdst,0) circle (1.5pt) node [right]{\scriptsize$a$};
                \clip(0,0)to[out=\tkzbigth,in=180-\tkzbigth](\tkzxdst,0)to[out=180+\tkzbigth,in=-\tkzbigth](0,0);
                \draw[thick,dotted,\colY](\tkzmd*\tkzxdst,0)
                to[out=-0.8*\tkzbigth,in=180+0.8*\tkzbigth]
                coordinate[pos=0] (bottom0)
                coordinate[pos=0.1] (bottom1) coordinate[pos=0.2] (bottom2)
                coordinate[pos=0.3] (bottom3) coordinate[pos=0.4] (bottom4)
                coordinate[pos=0.5] (bottom5) coordinate[pos=0.6] (bottom6)
                coordinate[pos=0.7] (bottom7) coordinate[pos=0.8] (bottom8)
                coordinate[pos=0.9] (bottom9) coordinate[pos=1.0] (bottom10)
                ({\tkzxdst-\tkzmd*\tkzxdst},0);
                \draw[thick,dotted,\colX](\tkzmd*\tkzxdst,0)
                to[out=0.8*\tkzbigth,in=180-0.8*\tkzbigth]
                coordinate[pos=0] (top0)
                coordinate[pos=0.1] (top1) coordinate[pos=0.2] (top2)
                coordinate[pos=0.3] (top3) coordinate[pos=0.4] (top4)
                coordinate[pos=0.5] (top5) coordinate[pos=0.6] (top6)
                coordinate[pos=0.7] (top7) coordinate[pos=0.8] (top8)
                coordinate[pos=0.9] (top9) coordinate[pos=1.0] (top10)
                ({\tkzxdst-\tkzmd*\tkzxdst},0);
                \foreach\x in {1,...,9} {
                        \draw[thick,\colY] (bottom\x)
                        to[out=-90,in=-90,looseness=1.4] ++({\tkzsepr*(1-0.1*\x)*((1-\tkzmd)*\tkzxdst-\tkzmd*\tkzxdst)},0)
                        to[out=90,in=-90] ++({1-.175*\x}, 1.6);
                        \draw[thick,\colX] (top\x)
                        to[out=90,in=90,looseness=1.4] ++({-\tkzsepr*(0.1*\x)*((1-\tkzmd)*\tkzxdst-\tkzmd*\tkzxdst)},0)
                        to[out=-90,in=90] ++({0.75-.175*\x}, -1.6);
                    }
                \node at (0.5*\tkzxdst,0) {\scriptsize$f^\dag$};
                \draw[thick,\colb] (top0)
                to[out=-90,in=-90+-35,out looseness=1.1,in looseness=1.7] ++(0.55*\tkzxdst,-.05*\tkzxdst) to[in=200,out=-90+-35] (\tkzxdst,0);
                \draw[thick,\cola] (bottom10)
                to[out=90,in=90-35,out looseness=1.1,in looseness=1.7] ++(-0.55*\tkzxdst,.05*\tkzxdst) to[in=20,out=90-35] (0,0);
                \fill[\colb] (\tkzmd*\tkzxdst,0) circle (1pt);
                \fill[\cola] ({\tkzxdst-\tkzmd*\tkzxdst},0) circle (1pt);
            \end{tkz}
            = f^{\dag\vee}.
        \end{multlined}
    \end{equation}

    \itemstep{Coheretors.}
    The fact that the canonical tensorators and unitors of $\vee$ are unitary also follows from invariance under isotopy-rel-boundary, and is left to the reader.

    \itemstep{Pivotality.}
    The induced pivotal structure $\phi=\{\phi_X\in\fX_\C(X\Rightarrow X^{\vee\vee})\}_X$ is strict by extended-isotopy invariance:
    \[
        \phi_X
        \;=\;
        \begin{tkz}[scale=0.7]
            \coordinate (L)  at (-3, 0);
            \coordinate (M)  at (-1, 0);
            \coordinate (P)  at (1, 0);
            \coordinate (R)  at (3, 0);
            \coordinate (M') at (-1, -2);
            \coordinate (Q)  at (1, 2);
            \begin{scope}
                \clip (L) to[out=-60,in=180] (M') to[out=0,in=-100] (R)to[out=120,in=0] (Q) to[out=180,in=60] cycle;
                \foreach \s in {0.2,0.4,...,1.6,1.8} {
                    \draw[\colX,thick]
                        (-1-\s,-2.5) -- (-1-\s,0)
                        to[out=110,in=90,looseness=1.4,xslant=0.5] (-1+\s,0)
                        to[out=-110,in=-90,looseness=1.4] (3-\s,0)
                        -- (3-\s,2.5);
                }
            \end{scope}
            \draw[\colb,thick] (M) -- (M');
            \draw[\cola,thick] (P) -- (Q);
            \draw[\colX,thick] (L) to[out=-60,in=180]node[pos=.6,below left=-2pt,\colX]{\scriptsize$\fcj{{}_aX_b}$} (M');
            \draw[\colb,thick] (M') to[out=0,in=-100]node[pos=.45,below right=-1pt]{\scriptsize$1_b$} (R);
            \draw[\cola,thick] (L) to[out=60,in=180]node[pos=.45,above left=-1pt]{\scriptsize$1_a$} (Q);
            \draw[\colX,thick] (Q) to[out=0,in=120]node[pos=.4,above right=-2pt,\colX]{\scriptsize${}_aX_b^{\vee\vee}={}_aX_b$} (R);
            \fill[\cola] (L) circle (1.5pt) node[left]{\scriptsize$\fcj{a}$};
            \fill[\colb] (R) circle (1.5pt);
            \fill[\colb] (M) circle (1pt);
            \fill[\cola] (P) circle (1pt);
            \fill[\colb] (M') circle (1.5pt);
            \fill[\cola] (Q) circle (1.5pt);
        \end{tkz}
        =
        \id_{X}.
        \qedhere
    \]
\end{pf}

\begin{example}
\label{skein 2-category}
Let $n\geq 2$, let $\C$ be a disk-like $n$-category, and fix an $(n-2)$-manifold $P$ and a boundary condition $\delta\in\undfld[n-3]{\C}{\bdy P}$. Then there is an (ordinary/traditional) pivotal dagger 2-category $\Sk_\C(P,\delta)$, called the \defn{skein 2-category} of $\C$ on $(P,\delta)$, defined by the skeletonization of the disk-like skein 2-category $\A_\C(P,\delta)$ from \cref{disk-like skein k-categories}:
\[
\Sk_\C(P,\delta)\coloneq\fX_{\A_\C(P,\delta)}.
\]
Thus
\begin{itemize}
    \item objects are $(n-2)$-fields $a\in\undfld[n-2]{\C}{P}[\delta]$, and
    \item 1-morphisms ${}_aX_b\in\Sk_\C(P,\delta)(a\to b)$ are $(n-1)$-fields $X\in\undfld[n-1]{\C}{P\times D^1}[\fcj a\cup b]$
    where the product $P\times D^1$ is pinched so that $\partial (P\times D^1)=\orev{P}\cup_{\partial P}P$, and
    \item 2-morphisms $f\in\Sk_\C(P,\delta)({}_aX_b\Rightarrow {}_aY_b)$ are (quotient-level) $n$-fields $f\in\undfldl[n]{\C}{P\times D^2}[(\fcj X\cup Y)\cup(\delta\times D^2)]$ 
    where the product $P\times D^2$ is pinched so that $\partial (P\times D^2)=\orev{P\times D^1}\cup_{\partial (P\times D^1)}(P\times D^1)$, where the two copies of $P\times D^1$ are themselves pinched as above.
\end{itemize}
\end{example}

\subsubsection{Unitarity and finiteness of \texorpdfstring{$\fX_\C$}{X(C)}}

\begin{lemma}
\label{lemXCunitary}
For a unitary disk-like 2-category $\C$, $\fX_\C$ is a pivotal $\Cstar$-2-category.
\end{lemma}

\begin{pf}
That $\fX_\C$ is a pivotal dagger 2-category is \cref{skein 2-categories are unitary planar pivotal}. 
Fix $a_1,\dots,a_\ell\in\fX_\C$ and let $\cL\coloneq\cL(a_1,\dots,a_\ell)$ as in \cref{2.25}. For $X^1,\dots,X^m\in\cL$, the linking algebra $L\coloneq\bigoplus_{p,q=1}^m\cL(X^q\to X^p)$ is finite-dimensional by \ref{FDF} with pairing $\bkt{(f^{pq})}{(g^{pq})}_L\coloneq\sum_{p,q=1}^m\sum_{i,j=1}^\ell\bkt{f^{pq}_{ij}}{g^{pq}_{ij}}_{D^2,\,\fcj{X^q_{ij}}\cup X^p_{ij}}$, which inherits positive-definiteness from \ref{2psnP}. Since $L$ is $*$-definite, it is unitary by the Characterization of Unitary Algebras \cite[Thm. I.2.4.4]{UQSL}. The monoidal product $-\otimes-$ on $\cL$ is a dagger functor by the computation in the proof of \cref{daggers for skeins ii} and has unitary associators and unitors, since these are computed entrywise from those of $\fX_\C$. Finally, for an object $X=(X_{ij})\in\cL$, the ``adjoint-transpose'' $(X^\vee)_{ij}\coloneq X_{ji}^\vee$ together with $(\ev_X)_{ii}\coloneq\bigoplus_{k=1}^\ell\ev_{X_{ki}}$ and $(\coev_X)_{ii}\coloneq\bigoplus_{k=1}^\ell\coev_{X_{ik}}$ (and zeros off the diagonal) satisfies the zig-zag equations entrywise, since $(\id_X\otimes\ev_X)\circ(\coev_X\otimes\id_X)$ has $(i,j)$-entry $(\id_{X_{ij}}\otimes\ev_{X_{ij}})\circ(\coev_{X_{ij}}\otimes\id_{X_{ij}})=\id_{X_{ij}}$. Thus $\cL$ is rigid $\Cstar$-monoidal, so $\fX_\C$ is $\Cstar$. 
\end{pf}

\begin{lemma}
  \label{hom 1-categories are finite unitary}
  If $\C$ is a finite unitary disk-like 2-category, then each hom 1-category $\fX_\C(a\to b)$ is finite unitary and
  \begin{equation}
    \label{hom count}
    |\Irr(\fX_\C(a\to b)^\cent)|
    =
    \dim\undfldl[2]{\C}{S^1\times D^1}[(\fcj a\amalg b)\times S^1].
  \end{equation}
  In particular $\fX_\C(a\to b)=0$ if and only if $\undfldl[2]{\C}{S^1\times D^1}[(\fcj a\amalg b)\times S^1]=0$.
\end{lemma}

\begin{pf}
  Since $\fX_\C(a\to b)=\Sk_\C(D^1,\fcj a\amalg b)$ as dagger 1-categories, applying \cref{lem:closing-up} with $Y\coloneq D^1$ and $\eta\coloneq\fcj a\amalg b$ gives both finite unitarity and \eqref{hom count}.
\end{pf}

\begin{remark} 
\label{pointsconnect}
Let $\Gamma$ denote the gluing map for the gluing of $D^1\times D^1$ to $S^1\times D^1$. Then for $v\in\End_{\fX_\C}(r)$ with $r\in\fX_\C(a\to b)$, $p,p'\in\End_{\fX_\C}(1_a)$, and $q\in\End_{\fX_\C}(1_b)$
\begin{equation} 
\label{eq:annular-insertion}
\Gamma(q)\circ\Gamma(v)\circ\Gamma(p)=\Gamma\big((p\otimes\id_r\otimes q)\circ v\big)\qquad\text{and}\qquad \Gamma(p')\circ\Gamma(p)=\Gamma(p'\circ p) 
\end{equation}
by \ref{CGa} and the interchange law (\cref{skein 2-categories are linear}). By \cref{glu-surj}, $\Sk_\C(S^1)(r_a\to r_b)=\sum_{X\in\fX_\C(a\to b)}\Gamma(\End_{\fX_\C}(X))$. Now let $(a,p)$ be a point of $\fX_\C$. Then $\Gamma(p)$ is a nonzero projection in the finite-dimensional $\Cstar$-algebra $\End_{\Sk_\C(S^1)}(r_a)$:
\begin{itemize} 
\item\emph{(Idempotent)} $\Gamma(p)\circ\Gamma(p)=\Gamma(p\circ p)=\Gamma(p)$ by \eqref{eq:annular-insertion}.
\item \emph{(Self-adjoint)} $\Gamma(p)^\dag=\Gamma(p^\dag)=\Gamma(p)$ by \ref{axiom:refl-gluing} and extended-isotopy invariance.
\item \emph{(Nonzero)} Observe $\Gamma(p)\triangleright p = p\circ p = p\neq 0$, where $\triangleright$ is the usual gluing action of annular fields onto disk fields.
\end{itemize}
Now choose any minimal projection $s_p$ with $s_{p}\leq\Gamma(p)$. This gives a simple object $(r_a,s_p)\in\Irr(\Sk_\C(S^1)^\cent)$.
\end{remark}

\begin{lemma}
\label{lem:XC-simple-connected}
Let $\C$ be a finite unitary disk-like 2-category and write $r_x\coloneq x\times S^1$ for 0-fields $x$ in $\C$.
\begin{lst}
\item[(a)]\label{XCsca} An object $a\in\fX_\C$ is simple if and only
if $a$ is a minimal 0-field of $\C$, and two simple objects $a,b\in\fX_\C$ are connected in the
sense of \cref{def:connected} if and only if $a\sim b$. Hence $a\mapsto[a]$ is an injection $\pi_0\C\hookrightarrow\pi_0\fX_\C$.
\item[(b)]\label{XCscb} Every point $(a,p)$ of $\fX_\C$ gives some simple object $(r_a,s_p)\in\Irr(\Sk_\C(S^1)^\cent)$, and associated simple objects of disconnected points are nonisomorphic. Hence
\end{lst}
\[
|\pi_0\fX_\C|\leq|\Irr(\Sk_\C(S^1)^\cent)|=\dim\undfldl[2]{\C}{S^1\times S^1}.
\]
\end{lemma}
\begin{pf}
\itemstep{\ref{XCsca}}
First observe that since $1_a=a\times I$, we have $\End_{\fX_\C}(1_a)=\fldl[2]{\C}{D^2}[r_a]$, so $a$ is simple if and only if $a$ is minimal. Now let $a$ and $b$ be minimal 0-fields in $\C$. As $[\xi\times D^1]=\id_\xi$ for 1-morphisms $\xi\in\fX_\C(a\to b)=\fld[1]{\C}{D^1}[\fcj{a}\amalg b]$, since a linear category is nonzero precisely when some object has nonzero identity we have
\begin{align*}
\text{$a$ and $b$ are connected in $\fX_\C$}
&
\iff
\fX_\C(a\to b)\neq0 
\\&
\iff 
[\xi\times D^1]\neq0 \text{ for some $\xi\in\fld[1]{\C}{D^1}[\fcj a\amalg b]$}
\\&
\overset{\text{def}}{\iff}
\text{some $\xi\in\fld[1]{\C}{D^1}[\fcj a\amalg b]$ is nondegenerate}
\\&
\underset{\eqref{lemcon}}{\iff}
a\sim b\text{ in $\C$}.
\end{align*}

\itemstep{\ref{XCscb}}
\cref{pointsconnect} gives the first claim. For the second claim, let $(a,p)$ and $(b,q)$ be disconnected points. Since $p\otimes\id_r\otimes q=0$ for all $r\in\fX_\C(a\to b)$, surjectivity of $\Gamma$ and \eqref{eq:annular-insertion} give $\Gamma(q)\circ g\circ\Gamma(p)=0$ for all $g\in\Sk_\C(S^1)(r_a\to r_b)$. As $s_p\leq\Gamma(p)$ and $s_q\leq\Gamma(q)$, it follows that $s_q\,\Sk_\C(S^1)(r_a\to r_b)\,s_p=0$, whence $(r_a,s_p)\not\cong(r_b,s_q)$ by Schur's Lemma. The assignment of the simple $(r_a,s_p)$ to one representative point per component gives an injection $\pi_0\fX_\C\hookrightarrow\Irr(\Sk_\C(S^1)^\cent)$. Finally, the equality $|\Irr(\Sk_\C(S^1)^\cent)|=\dim\undfldl[2]{\C}{S^1\times S^1}$ follows from \cref{lem:closing-up} with $Y\coloneq S^1$.
\end{pf}

\begin{proposition}
\label{unitary implies finite unitary cylinder 2-categories}
For a finite unitary disk-like 2-category $\C$, $\fX_\C$ is a finite pivotal $\Cstar$-2-category.
\end{proposition}

\begin{pf}
  Fix objects $a_1,\dots,a_\ell\in\fX_\C$ and let $\cL\coloneq\cL(a_1,\dots,a_\ell)$.

  \itemstep{$|\Irr(\cL^\cent)|<\infty$.}
  Since $\Irr(\cL^\cent)=\bigsqcup_{i,j=1}^{\ell}\Irr(\fX_\C(a_i\to a_j)^\cent)$ and each $\fX_\C(a_i\to a_j)$ is finite unitary by \cref{hom 1-categories are finite unitary}, $|\Irr(\cL^\cent)|<\infty$. Thus $\cL^\cent$ is multifusion.

\itemstep{$|\pi_0\fX_\C|<\infty$.}
  Indeed, $|\pi_0\fX_\C|\underset{\eqref{lem:XC-simple-connected}}{\leq}|\Irr(\Sk_\C(S^1)^\cent)|\underset{\eqref{lem:closing-up}}{=}\dim\undfldl[2]{\C}{S^1\times S^1}\underset{\ref{2psnF}}{<}\infty$.
\end{pf}

\subsubsection{\texorpdfstring{$\fX_\C$}{X(C)} as a proto-3-Hilbert space}

We can finally give the construction of a proto-3-Hilbert space from a unitary disk-like 2-category.

\begin{lemma}
\label{cor:cl-slide}
For all ${}_aX_b\in\fX_\C(a\to b)$ and $f\in\End_{\fX_\C}(X)$, $\cl_a(\tr_L^\vee(f))=\cl_b(\tr_R^\vee(f))$.
\end{lemma}

\begin{pf}
Let $f\in\End_{\fX_\C}({}_aX_b)$ and write $\tr_R^\vee(f)$ (resp. $\tr_L^\vee(f)$) for the string diagram on $D^2$ given by a single $X$-labeled loop enclosing an $a$-labeled (resp. $b$-labeled) region with an $f$-labeled vertex on the $X$-loop; this description determines a unique element of $\fldl[2]{\C}{D^2}[b\times \partial D^2]$ (resp. $\fldl[2]{\C}{D^2}[a\times \partial D^2]$) since such a diagram is unique up to isotopy-rel-boundary. Observe that by cutting out a $b$-labeled disk from the string diagram $\cl_b\tr_R^\vee(f)\in\undfldl[2]{\C}{S^2}$ and re-gluing, we obtain $\cl_a\tr_L^\vee(f)\in\undfldl[2]{\C}{S^2}$. Thus $\cl_b\tr_R^\vee(f)=\cl_a\tr_L^\vee(f)$ as 2-fields on $S^2$.
\end{pf}

\begin{construction}
\label{cstr:proto-3-Hilbert spaces from unitary disk-like 2-categories}
Let $\C$ be a unitary disk-like 2-category. We equip the strictly pivotal $\Cstar$-2-category $\fX_\C$ constructed above with the following spherical weight. For each $a\in\fX_\C$, define $\Psi_a^{\fX_\C}\colon\fX_\C(1_a\Rightarrow 1_a)\to\bC$ by 
\[
\Psi^{\fX_\C}_a(f)
\coloneq
\psn^\C(\cl_a(f))
=
\bkt{a\times D^2}{f}_{D^2,a\times S^1}
=\psn^\C\left(
\begin{tkz}[scale=0.355]
    \def\R{2.2}
    \def\E{0.55}
        \fill[gray!10] (0,\R) arc (90:270:\R) arc (-90:90:{\E} and \R);
        \foreach \w in {0.1, 0.3,...,0.9} {
            \draw[Green, thick, opacity=0.65] (0,\R) arc (90:270:{\w*\R} and \R);
        }
        \foreach \h in {-0.9, -0.8,..., 0.9} {
            \draw[Green, thick, opacity=0.65] ({-\R * sqrt(1 - (\h)*(\h))}, {\h*\R}) -- (0, {\h*\R});
        }
        \fill[gray!50] (0,0) ellipse ({\E} and \R);
        \draw[very thin] (0,\R) arc (90:270:\R); 
        \draw[thick, Green] (0,0) ellipse ({\E} and \R); 
\end{tkz} 
\blt
\begin{tkz}[scale=0.355]
\def\R{2.2}
\def\E{0.55}
\fill[\filf] (0,\R) arc (90:-90:\R) arc (270:90:{\E} and \R);
\draw[very thin] (0,\R) arc (90:-90:\R);
\draw[thick, Green] (0,\R) arc (90:270:{\E} and \R);
\draw[thick, Green, dashed] (0,\R) arc (90:-90:{\E} and \R);
\node[Blue!60!black] at (\R/2, 0) {\Large$f$};
\end{tkz}
\right)
=
\psn^\C\left(
\begin{tkz}[scale=0.355]
    \def\R{2.2}
    \def\E{0.55}
    \begin{scope}
        \fill[gray!10] (0,\R) arc (90:270:\R) arc (270:450:{\E} and \R);
        \foreach \w in {0.1, 0.3,..., 0.9} {
            \draw[Green, thick, opacity=0.65] (0,\R) arc (90:270:{\w*\R} and \R);
        }
        \foreach \h in {-0.9, -0.8,..., 0.9} {
            \draw[Green, thick, opacity=0.65] ({-\R * sqrt(1-(\h)*(\h))}, {\h*\R}) -- (0, {\h*\R});
        }
    \end{scope}
    \fill[\filf] (0,\R) arc (90:-90:\R) arc (270:90:{\E} and \R);
    \draw[very thin] (0,0) circle (\R); 
    \draw[thick, Green] (0,\R) arc (90:270:{\E} and \R);
    \draw[thick, Green, dashed] (0,\R) arc (90:-90:{\E} and \R);
    \node[Blue!60!black] at (\R/2, 0) {\Large$f$};
\end{tkz} 
\right).
\]
\end{construction}

\begin{lemma}
\cref{cstr:proto-3-Hilbert spaces from unitary disk-like 2-categories} indeed defines a proto-3-Hilbert space.
\end{lemma}

\begin{pf}
By \cref{lemXCunitary}, $\fX_\C$ is a pivotal $\Cstar$-2-category with its canonical UAF. To see $\Psi^{\fX_\C}$ is faithful, take an $f\in\fX_\C(1_a\Rightarrow 1_a)$ and isotope $f^\dag=\fcj{\iota_*f}$ to the opposite hemisphere to obtain $\Psi_a^{\fX_\C}(f^\dag\circ f)=\bkt{\id_{1_a}}{f^\dag\circ f}_{D^2,a\times S^1}\underset{\eqref{pairing-composition compatibility}}{=}\bkt{f}{f}_{D^2,a\times S^1}$, which by \ref{2psnP} is nonnegative with equality if and only if $f=0$.

To see $\Psi^{\fX_\C}$ is spherical, observe that for a 2-endomorphism $f\in\End_{\fX_\C}({}_aX_b)$, we have
\[
\Psi_a^{\fX_\C}(\tr_L^\vee(f))=\psn^\C(\cl_a(\tr_L^\vee(f)))\underset{\eqref{cor:cl-slide}}{=}\psn(\cl_b(\tr_R^\vee(f)))=\Psi_b^{\fX_\C}(\tr_R^\vee(f)).\qedhere
\]
\end{pf}

\subsection{From traditional to disk-like: the unitary disk-like 2-category \texorpdfstring{$\C^\fX$}{C(X)}}

\subsubsection{\texorpdfstring{$\fX$}{X}-string diagrams and the disk-like 2-category \texorpdfstring{$\C^\fX$}{C(X)}}

\begin{definition}
Let $\fX$ be a pivotal dagger 2-category.
\begin{lst} 
\item An \defn{$\fX$-string diagram} on a 1-ball $Y\in\Disk_1$ is a pair $\xi=(\Gamma,\lambda)$ where $\Gamma$ is a 1D string diagram stratification of $Y$ (see \cref{cstr:disk-like 1-categories from ordinary 1-categories}) and where $\lambda$ is a labeling of $\Gamma$ using the objects and $1$-morphisms of $\fX$.

\item An \defn{$\fX$-string diagram} on a 2-ball $X\in\Disk_2$ is a pair $\xi = (\Gamma,\lambda)$, where
\begin{itemize} 
\item $\Gamma$ is a \defn{2D string diagram stratification} of $X$, that is, a finite embedded graph in $X$ whose edges intersect $\partial X$ transversely and that decomposes $X$ into faces (also called 2-cells or regions) $\{f_j\}$, edges (also called 1-cells) $\{e_k\}$ each equipped with a transverse orientation $\varepsilon_k\in\{\pm1\}$, and interior vertices (also called 0-cells) $\{v_l\}$ each equipped with a \defn{link parametrization}, i.e., a homeomorphism $\upsilon_l \colon \mathrm{lk}(v_l)\to\partial D^2$ of the link of $v_l$; and 
\item $\lambda$ is an \defn{$\fX$-labeling} of $\Gamma$, which consists of an object label $\lambda_j \in \fX$ on each 2-cell, a 1-morphism label $\lambda_{e_k} \colon \lambda_{f_-} \to \lambda_{f_+}$ on each 1-cell where $\varepsilon_k$ identifies the adjacent source face $f_-$ and target face $f_+$, and a 2-morphism label $\lambda_{v_l} \colon \lambda_{e_1} \otimes \cdots \otimes \lambda_{e_m} \Rightarrow \lambda_{e'_1} \otimes \cdots \otimes \lambda_{e'_{m'}}$ at each 0-cell $v_l$, where the link parametrization $\upsilon_l\colon \mathrm{lk}(v_l)\to \partial D^2$ determines via the standard decomposition of $\partial D^2$ an ordering of the source strands $e_1<\cdots<e_m$ and the target strands $e'_1<\cdots<e'_{m'}$ ordered by $\partial_-\mathrm{lk}(v_l)$ and $\orev{\partial_+\mathrm{lk}(v_l)}$ respectively.
\end{itemize}
For a homeomorphism $\varphi\colon X\xrightarrow{\sim} X'$, we define $\varphi_*\xi\coloneq(\varphi_*\Gamma,\lambda\circ \varphi^{-1})$.
\end{lst}
We identify two $\fX$-string diagrams on a $k$-ball $W$ for $k\in\{1,2\}$ if one can be obtained from the other by erasing strands or vertices labeled by identity morphisms of $\fX$.
\end{definition}

\begin{construction}
\label{cstr:disk-like 2-categories from ordinary pivotal 2-categories}
A rigid $\Cstar$-2-category $\fX$ equipped with a choice of UAF $\vee$ gives a disk-like $2$-category $\C^\fX$ as follows.
\begin{lst}
    \item[\ref{DCk}] A field on a 0-ball $P$ is a labeling of $P$ by an object of $\fX$. A field on a 1-ball $Y$ is an $\fX$-string diagram $\xi=(\Gamma,\lambda)$ on $Y$. A pure $2$-field on a $2$-ball $X$ with boundary condition $\xi\in\undfld[1]{\C^\fX}{\partial X}$ is an $\fX$-string diagram $(\Gamma,\lambda)$ on $X$ whose restriction $(\Gamma|_{\partial X},\lambda|_{\partial X})$ to $\partial X$ equals $\xi$.
    \item[\ref{Dbdk}] For a string diagram $\xi=(\Gamma,\lambda)\in\fld{\C^\fX}{X}[c]$, define $\partial\xi\coloneq(\Gamma|_{\partial X},\lambda|_{\partial X})$.
    \item[\ref{DGk}] For compatible fields $\xi_1=(\Gamma_1,\lambda_1)$ and $\xi_2=(\Gamma_2,\lambda_2)$, define $\xi_1\blt\xi_2\coloneq(\Gamma_1\cup\Gamma_2,\lambda_1\cup\lambda_2)$.
    \item[\ref{Dpik}] For a pinched product map $\pi\colon X\to Y$ and a field $\xi\in\fld{\C^\fX}{Y}$, define $\pi^*\xi\coloneq(\pi^*\Gamma,\lambda\circ\pi)$. When $Y$ is a $0$-ball with $0$-field $a$, $\pi^*a\coloneq(\varnothing,a)$ is the string diagram on $X$ with a single 2-cell labeled by~$a$.
    \item[\ref{Drefl}] Reflection acts trivially on 0- and 1-fields. For a pure 2-field $\xi=(\Gamma,\lambda)$ on a 2-ball $X$, define $\fcj{\xi}\coloneq(\orev\Gamma,\lambda^\dag)$ where $\orev\Gamma$ is obtained from $\Gamma$ by postcomposing each link parametrization $\upsilon_v$ with the equatorial reflection $r\colon\partial D^2\to\partial D^2$, and where $\lambda^\dag$ agrees with $\lambda$ on 2- and 1-cell labels but replaces each 2-morphism label $\lambda_v$ by $\lambda_v^\dag$.
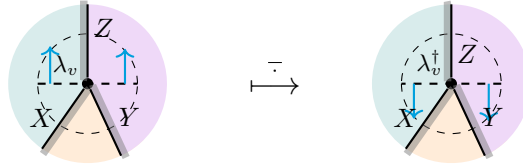
\[
    \begin{tkz}[scale=1.05]
        \fill[teal!15] (0,0) -- (90:1cm) arc (90:235:1cm) -- cycle;
        \fill[orange!15] (0,0) -- (235:1cm) arc (235:315:1cm) -- cycle;
        \fill[violet!15] (0,0) -- (-65:1cm) arc (-65:90:1cm) -- cycle;
        \node[circle, fill, inner sep=1.5pt, label={[label distance=-1.25mm]above left:{\scriptsize$\lambda_v$}}] (v) at (0,0) {};
        \pic[scale=2.2] at (0,0){link};
        \draw[thick, frameR] (235:1cm) -- (v) node[pos=.325, above left=-1.25mm, font=\footnotesize] {$X$};
        \draw[thick, frameR] (-65:1cm) -- (v) node[pos=.325, above right=-.5mm, font=\footnotesize,black] {$Y$};
        \draw[thick, frameR] (90:1cm) -- (v) node[pos=.325, right=-.75mm, font=\footnotesize] {$Z$};
    \end{tkz}
    \qquad
    \overset{\;\fcj{\,\cdot\,}\;}{\longmapsto}
    \qquad
    \begin{tkz}[scale=1.05]
        \fill[teal!15] (0,0) -- (90:1cm) arc (90:235:1cm) -- cycle;
        \fill[orange!15] (0,0) -- (235:1cm) arc (235:315:1cm) -- cycle;
        \fill[violet!15] (0,0) -- (-65:1cm) arc (-65:90:1cm) -- cycle;
        \node[circle, fill, inner sep=1.5pt, label={[label distance=-1.25mm]above left:{\scriptsize$\lambda_v^\dag$}}] (v) at (0,0) {};
        \pic[scale=2.2,rotate=180] at (0,0){link};
        \draw[thick, frameR] (235:1cm) -- (v) node[pos=.325, above left=-1.25mm, font=\footnotesize] {$X$};
        \draw[thick, frameR] (-65:1cm) -- (v) node[pos=.325, above right=-.5mm, font=\footnotesize,black] {$Y$};
        \draw[thick, frameL] (v) -- (90:1cm) node[pos=.325, right=-.75mm, font=\footnotesize] {$Z$};
    \end{tkz}
    \]
    \item[\ref{DCU}] The subspace $\fldlU{\C^\fX}{X}[c]$ of local relations is generated by formal differences $\xi-\xi'$ of string diagrams $\xi,\xi'$ that are related in any of the following ways.
    \begin{lst}
        \item (\emph{Isotopies rel boundary}) $\xi'=H(-,1)_*\xi$ for an isotopy-rel-boundary $H$ of $X$.
        \item (\emph{Adjacent vertex composition}) Two adjacent interior 0-cells $P_1,P_2$ connected by an interior 1-cell $E$ (with no other 0-cells on $E$) are replaced by a single 0-cell labeled by the 2-morphism obtained by composing $\lambda_{P_1}$ and $\lambda_{P_2}$ along $\lambda_E$, and $E$ is deleted.
        \item (\emph{Identity deletion}) An interior 0-cell labeled by an identity 2-morphism is deleted and the adjacent 1-cells are merged, or an interior 1-cell labeled by an identity 1-morphism is deleted, merging the adjacent 2-cells.
        \item (\emph{Pivotal reparameterization}) The link parametrization $\upsilon_v$ at a 0-cell $v$ is replaced by a different parametrization $\upsilon_v'$ transferring (rotating) 1-cells between the incoming and outgoing sides, with each transferred 1-cell having its label $X$ replaced by $X^\vee$ and the 2-morphism label $\lambda_v$ replaced by its image $\lambda_v'$ under the canonical isomorphism via the (co)evaluation maps. In short, $(\Gamma,\lambda)\sim(\Gamma',\lambda')$ where $\Gamma'$ differs from $\Gamma$ only in $\upsilon_v$. For example, in the following illustration $\lambda'_v=(\lambda_v\otimes\id_{Y^\vee})\circ(\id_X\otimes\coev_Y)$.
              \[
                  \begin{tkz}[scale=1.05]
                      \fill[teal!15] (0,0) -- (75:1cm) arc (75:235:1cm) -- cycle;
                      \fill[orange!15] (0,0) -- (225:1cm) arc (225:325:1cm) -- cycle;
                      \fill[violet!15] (0,0) -- (-35:1cm) arc (-35:75:1cm) -- cycle;
                      \draw[thick, frameR] (-135:1cm) -- (0,0) node[pos=.325, above left=-1.25mm, font=\footnotesize] {$X$};
                      \draw[thick, frameR] (-35:1cm) -- (0,0) node[pos=.275, above right=-.75mm, font=\footnotesize,black] {$Y$};
                      \draw[thick, frameR] (0,0) -- (75:1cm) node[pos=.675, right=-.75mm, font=\footnotesize] {$Z$};
                      \node[circle, fill, inner sep=1.5pt,label={[label distance=-1mm,xshift=-1.5mm]above:{\scriptsize$\lambda_v$}}] (v) at (0,0) {};
                      \draw (0,0) pic[scale=2]{link};
                  \end{tkz}
                  \qquad
                  \longmapsto
                  \qquad
                  \begin{tkz}[scale=1.05]
                      \fill[teal!15] (0,0) -- (75:1cm) arc (75:225:1cm) -- cycle;
                      \fill[orange!15] (0,0) -- (225:1cm) arc (225:325:1cm) -- cycle;
                      \fill[violet!15] (0,0) -- (-35:1cm) arc (-35:75:1cm) -- cycle;
                      \draw[thick, frameR] (-135:1cm) -- (0,0) node[pos=.325, above left=-1.25mm, font=\footnotesize] {$X$};
                      \draw[thick, frameR] (-35:1cm) -- (0,0) node[pos=.325, above right=-.5mm, font=\footnotesize,black] {$Y$};
                      \draw[thick, frameR] (0,0) -- (75:1cm) node[pos=.675, right=-.75mm, font=\footnotesize] {$Z$};
                      \node[circle, fill, inner sep=1.5pt,label={[label distance=-1mm,xshift=0mm,yshift=0.5mm]left:{\scriptsize$\lambda_v'$}}] (v) at (0,0) {};
                      \draw (0,0) pic[scale=2,rotate=-70]{link};
                  \end{tkz}
              \]
    \end{lst}
\end{lst}
\end{construction}

\begin{warning}
Often in the literature one comes across diagrams where cups and caps are equipped not with a normal framing (that is, a transverse orientation), but a tangential orientation, that is, 
\[
\ev_X=
\begin{tkz}[scale=.5]
   \fill[rounded corners,orange!15] (-0.5,0) rectangle (2.5, 1.75);
   \fill[teal!15] (0,0)--++(0,.25)arc(180:0:1)--++(0,-.25)--cycle;
   \draw[mid<=.125,mid<=.5,mid<=.95](0,0)--++(0,.25)arc(180:0:1)--++(0,-.25);
\end{tkz}
\qquad\text{and}\qquad
\coev_X=
\begin{tkz}[scale=.5,yscale=-1]
   \fill[rounded corners,teal!15] (-0.5,0) rectangle (2.5, 1.75);
   \fill[orange!15] (0,0)--++(0,.25)arc(180:0:1)--++(0,-.25)--cycle;
   \draw[mid<=.125,mid<=.5,mid<=.95](0,0)--++(0,.25)arc(180:0:1)--++(0,-.25);
\end{tkz}
\]
instead of 
\[
\ev_X=
\begin{tkz}[scale=.5]
   \fill[rounded corners,orange!15] (-0.5,0) rectangle (2.5, 1.75);
   \fill[teal!15] (0,0)--++(0,.25)to[out=90,in=90,looseness=1.75]++(2,0)--++(0,-.25)--cycle;
   \draw[frameL](0,0)--++(0,.25)to[out=90,in=90,looseness=1.75]++(2,0)--++(0,-.25);
\end{tkz}
\qquad\text{and}\qquad
\coev_X=
\begin{tkz}[scale=.5,yscale=-1]
   \fill[rounded corners,teal!15] (-0.5,0) rectangle (2.5, 1.75);
   \fill[orange!15] (0,0)--++(0,.25)to[out=90,in=90,looseness=1.75]++(2,0)--++(0,-.25)--cycle;
   \draw[frameL](0,0)--++(0,.25)to[out=90,in=90,looseness=1.75]++(2,0)--++(0,-.25);
\end{tkz}\,.
\]
However, in the presence of a reflection/dagger structure, this may be misleading when we think of the dagger as a reflection about some horizontal axis. For instance,
\[
\left(\begin{tkz}[scale=.5]
   \fill[rounded corners,orange!15] (-0.5,0) rectangle (2.5, 1.75);
   \fill[teal!15] (0,0)--++(0,.5)arc(180:0:1)--++(0,-.5)--cycle;
   \draw[mid>=.125,mid>=.5,mid>=.95](0,0)--++(0,.5)arc(180:0:1)--++(0,-.5);
\end{tkz}\right)^\dag
\quad=\quad
\begin{tkz}[scale=.5,yscale=-1]
   \fill[rounded corners,orange!15] (-0.5,0) rectangle (2.5, 1.75);
   \fill[teal!15] (0,0)--++(0,.5)arc(180:0:1)--++(0,-.5)--cycle;
   \draw[mid<=.125,mid<=.5,mid<=.95](0,0)--++(0,.5)arc(180:0:1)--++(0,-.5);
\end{tkz}
\;\quad\neq\quad\;
\begin{tkz}[scale=.5,yscale=-1]
   \fill[rounded corners,orange!15] (-0.5,0) rectangle (2.5, 1.75);
   \fill[teal!15] (0,0)--++(0,.5)arc(180:0:1)--++(0,-.5)--cycle;
   \draw[mid>=.125,mid>=.5,mid>=.95](0,0)--++(0,.5)arc(180:0:1)--++(0,-.5);
\end{tkz}\,.
\]
This is alleviated by the normal framings:
\[
\left(\begin{tkz}[scale=.5]
   \fill[rounded corners,orange!15] (-0.5,0) rectangle (2.5, 1.75);
   \fill[teal!15] (0,0)--++(0,.25)to[out=90,in=90,looseness=1.75]++(2,0)--++(0,-.25)--cycle;
   \draw[frameL](0,0)--++(0,.25)to[out=90,in=90,looseness=1.75]++(2,0)--++(0,-.25);
\end{tkz}\right)^\dag
\quad=\quad
\begin{tkz}[scale=.5,yscale=-1]
   \fill[rounded corners,orange!15] (-0.5,0) rectangle (2.5, 1.75);
   \fill[teal!15] (0,0)--++(0,.25)arc(180:0:1)--++(0,-.25)--cycle;
   \draw[frameR](0,0)--++(0,.25)arc(180:0:1)--++(0,-.25);
\end{tkz}
\;\quad\neq\quad\;
\begin{tkz}[scale=.5,yscale=-1]
   \fill[rounded corners,orange!15] (-0.5,0) rectangle (2.5, 1.75);
   \fill[teal!15] (0,0)--++(0,.25)to[out=90,in=90,looseness=1.75]++(2,0)--++(0,-.25)--cycle;
   \draw[frameL](0,0)--++(0,.25)to[out=90,in=90,looseness=1.75]++(2,0)--++(0,-.25);
\end{tkz}\,.
\]
This should not be too surprising, since it is the normal framings that indicate the Poincar\'e-dual pasting diagram, and the dagger reverses the arrows on the pasting diagram.
\end{warning}

\begin{remark}
We could instead choose to draw the normal framings of the 1-strata as transverse arrows, as illustrated below.
\[
\begin{tkz}[scale=1.05]
  \fill[teal!15] (0,0) -- (75:1cm) arc (75:225:1cm) -- cycle;
  \fill[orange!15] (0,0) -- (225:1cm) arc (225:325:1cm) -- cycle;
  \fill[violet!15] (0,0) -- (-35:1cm) arc (-35:75:1cm) -- cycle;
  \draw[thick,frameR] (-135:1cm) -- (0,0) node[pos=.325, above left=-1.25mm, font=\footnotesize] {$X$};
  \draw[thick,frameR] (-35:1cm) -- (0,0) node[pos=.325, above right=-.5mm, font=\footnotesize,black] {$Y$};
  \draw[thick,frameR] (0,0) -- (75:1cm) node[pos=.675, right=-.75mm, font=\footnotesize] {$Z$};
  \node[circle, fill, inner sep=1.5pt,label={[label distance=-1mm,xshift=0mm,yshift=0.5mm]left:{\scriptsize$\lambda_v'$}}] (v) at (0,0) {};
  \draw (0,0) pic[scale=1.5,rotate=-70]{link};
\end{tkz}
\qquad\qquad\rightsquigarrow\qquad\qquad
\begin{tkz}[scale=1.05]
  \fill[teal!15] (0,0) -- (75:1cm) arc (75:225:1cm) -- cycle;
  \fill[orange!15] (0,0) -- (225:1cm) arc (225:325:1cm) -- cycle;
  \fill[violet!15] (0,0) -- (-35:1cm) arc (-35:75:1cm) -- cycle;
  \draw[thick] (-135:1cm) --pic[pos=0.3,rotate=-135+90]{tro} (0,0) node[pos=.325, above left, font=\footnotesize] {$X$};
  \draw[thick] (-35:1cm) --pic[pos=0.3,rotate=-35+90]{tro} (0,0) node[pos=.325, above right=-.5mm, font=\footnotesize,black] {$Y$};
  \draw[thick] (0,0) --pic[pos=0.7,rotate=75-90]{tro} (75:1cm) node[pos=.675, right=-.75mm, font=\footnotesize] {$Z$};
  \node[circle, fill, inner sep=1.5pt,label={[label distance=-1mm,xshift=0mm,yshift=0.5mm]left:{\scriptsize$\lambda_v'$}}] (v) at (0,0) {};
  \draw (0,0) pic[scale=1.5,rotate=-70]{link};
\end{tkz}
\]
\end{remark}

\subsubsection{Evaluation maps for string diagrams}
\begin{facts} 
\label{eval-facts-2}
If $\fX$ is a pivotal $\dag$-2-category, then the following facts hold.
\begin{lst} 
\item For each 1-ball $Y$ and objects $a,b\in\fX$, there is a set map
$$\eval^\fX_{Y;\,a\to b}\colon\fld{\C^\fX}{Y}[\fcj a\amalg b]\to\mathrm{Obj}(\fX(a\to b))$$
sending an $\fX$-string diagram on $Y$ to the 1-morphism it represents in $\fX$:
$$\eval^\fX_{Y;\,a\to b}(\xi)\coloneq\lambda'_{v_1}\otimes\lambda'_{v_2}\otimes\cdots\otimes\lambda'_{v_m}$$
where $V_\xi=\{v_1<\cdots<v_m\}$ and
$$\lambda'_v\coloneq\begin{cases}\lambda_v&\text{if }\varepsilon_v=1,\\\lambda^\vee_v&\text{if }\varepsilon_v=-1.\end{cases}$$
If $\xi$ has no vertices, $\eval(\xi)\coloneq 1_a$. Since 1-fields are below the top dimension, $\eval$ is well-defined on pure 1-fields. Since $\fX$ is strict (\cref{strictnessconvention}), $\eval$ satisfies the strict equalities $\eval(\xi\otimes\eta)=\eval(\xi)\otimes\eval(\eta)$ and $\eval(1_a)=1_a$. Since there are no choices to be made and $\eval^\fX_{Y;\,a\to b}(\xi)=\eval^\fX_{Y';\,a\to b}(\xi)$ for any two 1-balls $Y$ and $Y'$, we may drop the above decorations and unambiguously write $\eval$ for the above map.

\item For all 2-balls $X$, $\fX$-string diagrams $c\in\undfld[1]{\C^\fX}{\partial X}$, and splittings $\partial X=\partial_-X\cup\partial_+X\eqcolon\partial_{\pm}X$ along which $c$ splits, there is a set map
\[
\eval^\fX_{X;\partial_{\pm}X}\colon\fld{\C^\fX}{X}[\fcj {c_-}\cup c_+]\to\fX(\eval(c_-)\Rightarrow \eval(c_+))
\]
sending a pure $\fX$-string diagram $\alpha$ on $X$ to the 2-morphism $\eval^\fX_{X;\partial_{\pm}X}(\alpha)$ it represents in $\fX$ with respect to this splitting. For example, see \cite{pivbicatSNpaper}. Alternatively, one could use Morse functions to decompose the string diagram into ``elementary pieces'' and to then compose these to obtain the morphism; one then shows that this assignment is independent of the choices made during this procedure.

\item\label{fact:eval-ker2}
 \emph{Local relations are $\ker(\eval)$.}
 The local relations $\fldlU{\C^\fX}{X}[c]$ of \cref{cstr:disk-like 2-categories from ordinary pivotal 2-categories} coincide with $\ker(\eval)$. Consequently, the linear extension of $\eval$ descends to a well-defined injective linear map $\bar\eval\colon\fldl{\C^\fX}{X}[c]\to\fX(\eval(c_-)\Rightarrow \eval(c_+))$. 
 By considering the cone string diagram over a single $\alpha$-labeled vertex at the center of $D^2$, we easily see $\overline{\eval}$ is also surjective.
 Also observe that by local relations $\eval(\beta\circ\alpha)=\eval(\beta)\circ\eval(\alpha)$ for vertically compatible 2-fields $\alpha$ and $\beta$ and that by strictness (see \cref{strictnessconvention}) $\eval(\alpha\otimes\beta)=\eval(\alpha)\otimes\eval(\beta)$ for horizontally compatible 2-fields $\alpha$ and $\beta$. 

\item\label{fact:eval-dag2}
\emph{The map $\eval$ preserves $\dag$.} 
By a standard Morse-theoretic argument, $\eval$ preserves the dagger on 2-fields.
\end{lst}
\end{facts}

\begin{facts}
\label{anglebracket}
For a 1-morphism $X\in\fX(a\to b)$, write $\langle X\rangle\in\fld{\C^\fX}{D^1}[\fcj{a}\amalg b]$ for the single-vertex string diagram on $D^1$ with positively oriented vertex $v$ at $0$ labeled by $X$. For a 2-morphism $f\in\fX(X\Rightarrow Y)$, write
\[
\langle f\rangle\coloneq\bar\eval^{-1}(f)\in\fX_{\C^\fX}(\langle X\rangle\Rightarrow\langle Y\rangle)
\]
for the (quotient-level) 2-field on $D^2$ represented by the cone diagram over the single vertex at the center $(0,0)$ of $D^2$ with incoming and outgoing boundary fields given by $\partial_-\langle f\rangle=\langle X\rangle$ and $\partial_+\langle f\rangle=\langle Y\rangle$ respectively with respect to the standard splitting $\partial_{\pm}D^2$ of $\partial D^2$. The link parameterization at $v$ agrees with the standard splitting $\partial_{\pm}D^2$ of $\partial D^2$ in that $\upsilon_v$ is the radial identification
$\mathrm{lk}(v)\xrightarrow{\sim}\partial D^2$, so that
$\eval^\fX_{D^2,\partial_{\pm}D^2}(\langle f\rangle)=f$. Observe that $\langle-\rangle$ has the following properties.
\begin{lst}
\item $\langle\alpha\circ\beta\rangle=\langle\alpha\rangle\circ\langle\beta\rangle$ for all 2-morphisms $\alpha$, $\beta$ in $\fX$.
\item $\langle\id_X\rangle=\id_{\langle X\rangle}$ for all 1-morphisms $X$ in $\fX$.
\item For all composable 1-morphisms ${}_aX_b$ and ${}_bY_c$ in $\fX$, the compositor $\langle-\rangle^2_{X,Y}\colon\langle{}_aX_b\rangle\otimes\langle{}_bY_c\rangle\overset\sim\to\langle{}_aX\otimes Y_c\rangle$ represented by the obvious trivalent vertex diagram is unitary.
\item $\langle\id_{1_a}\rangle=\id_{1_{a}}$ and $\langle 1_a\rangle=1_a$ for objects $a\in\fX$. (Indeed, recall that we identify string diagrams with those obtained from them by erasing identity morphisms.)
\item 
$\langle-\rangle$ is UAF-preserving.
\item 
$\langle-\rangle$ preserves $\dag$.
\item If $\fX$ has a spherical weight $\Psi$, then $\Psi^{\fX_{\C^\fX}}_a(\langle\alpha\rangle)=\Psi^\fX_a(\alpha)$ for all $a\in\fX$ and $\alpha\in\End_\fX(1_a)$.
\end{lst}
\end{facts}

The above facts immediately imply the following corollary.

\begin{corollary} 
The identity-on-objects functor $\langle-\rangle\colon\fX\to\fX_{\C^\fX}$ defined in \cref{anglebracket} is an equivalence of pivotal $\Cstar$-2-categories, and is isometric when $\fX$ is a proto-3-Hilbert space. Moreover, $\eval\circ\langle-\rangle=\id_\fX$, so $\eval$ is a (weak) inverse to $\langle-\rangle$.
\end{corollary}

\begin{lemma}
\label{lem-eval-equiv-2}
Let $\fX$ be a pivotal $\Cstar$-2-category. If 2-fields $\alpha,\beta\in\fld[2]{\C^\fX}{X}$ on a 2-ball $X$ with $\partial\alpha=\partial\beta$ satisfy $\eval(\alpha)=\eval(\beta)$, then $[\alpha]=[\beta]$. If 1-fields $\xi,\eta\in\fld[1]{\C^\fX}{W}[\fcj{a}\amalg b]$ on a 1-ball $W$ satisfy $\eval(\xi)\cong^\dag\eval(\eta)$ in $\fX$, then $\xi\cong^\star\eta$ in $\C^\fX$. In particular, $\xi\cong^\star\langle\eval(\xi)\rangle$ for 1-fields $\xi$ on $D^1$.
\end{lemma}

\begin{pf}
The first assertion is trivial since $\bar\eval$ is injective by \cref{eval-facts-2}\ref{fact:eval-ker2}. Now assume $\xi$ and $\eta$ are 1-fields in $\C^\fX$, which we may assume are $D^1$-shaped by choosing homeomorphisms, and let $\alpha\in\fX(\eval(\xi)\Rightarrow\eval(\eta))$ be a unitary isomorphism. Set $\langle\alpha\rangle\coloneq\bar\eval^{-1}(\alpha)\in\fld{\A_{\C^\fX}(W,\fcj{a}\amalg b)}{I}[\fcj{\xi}\amalg \eta]=\fld{\C^\fX}{W\times I}[\fcj{a\times I}\cup\fcj{\xi}\cup(b\times I)\cup\eta]$, so that $\eval(\langle\alpha\rangle^\star)=\alpha^\dag$ by \cref{eval-facts-2}\ref{fact:eval-dag2}. Writing $\langle\beta\rangle\circ\langle\alpha\rangle\coloneq\rho^2_*(\langle\alpha\rangle\blt_{D^1}\langle\beta\rangle)$ (2-composition in $\fX_\C$), we have
\[
\eval(\langle\alpha\rangle\bullet_{D^1}\langle\alpha\rangle^\star)=\alpha^\dag\circ\alpha=\id_{\eval(\xi)}=\eval(\xi\times I),
\]
so the first assertion gives $\langle\alpha\rangle\bullet_{D^1}\langle\alpha\rangle^\star=\xi\times I$, and similarly $\langle\alpha\rangle^\star\bullet_{D^1}\langle\alpha\rangle=\eta\times I$. Thus $\xi\cong^\star\eta$ via $\langle\alpha\rangle$. Now taking $\eta\coloneq\langle\eval(\xi)\rangle$ and $\alpha=\id_{\eval\xi}$ gives the last assertion.
\end{pf}

\subsubsection{The sphere trace on \texorpdfstring{$\C^\fX$}{C(X)}}

\

\begin{construction} 
\label{cstr:unitary disk-like 2-categories from proto-3-Hilbert spaces}
\label{def:Phi-eX-S2}
Let $(\fX,\vee,\Psi)$ be a proto-3-Hilbert space. We construct a unitary disk-like 2-category $(\C^{\fX},\psn^{\C^\fX})$ with underlying disk-like 2-category $\C^{\fX}$ from \cref{cstr:disk-like 2-categories from ordinary pivotal 2-categories}.
Consider a pure $2$-field (string diagram) $\alpha\in\undfld[2]{\C^\fX}{S^2}$. Choose a $2$-ball $X_0\hookrightarrow S^2\setminus(V_\alpha\cup E_\alpha)$ and set $X\coloneq S^2\setminus X_0$, and let $f_{X_0}\in F_\alpha$ be the face of $\alpha$ containing $X_0$. Let $a_X$ denote the object of $\fX$ such that $\alpha|_{\partial X}=a_X\times\partial X$. Now choose any two-ball splitting $\partial X=\partial_-X\cup\partial_+X$, and set
\begin{equation*}
  \tld\psn^{\C^\fX}(\alpha)\coloneq\Psi^{\fX}_{a_X}(\eval^\fX_{X;\partial_{\pm}X}(\alpha|_X)).
\end{equation*}
\[
\begin{tkz}[scale=1]
  \def\R{2.2}
  \coordinate (u1) at (-0.30,0.62);   \coordinate (u2) at (0.12,-0.68);
  \coordinate (u3) at (-1.13,-1.35);  \coordinate (u4) at (-1.52,0.83);
  \coordinate (w1) at (0.65,1.61);    \coordinate (w2) at (1.61,-0.65);
  \fill[Fd] (0,0) circle (\R);
  \fill[Fa] (u3) .. controls (-2.05,-1.00) and (-2.03,0.64) .. (u4)
                 .. controls (-1.35,1.98)  and (0.01,1.95)  .. (w1)
                 .. controls (2.04,1.52)   and (2.06,-0.56) .. (w2)
                 .. controls (1.33,-1.82)  and (-0.67,-1.75).. cycle;
  \fill[Fb] (u1) .. controls (0.30,0.24)  and (-0.30,-0.26).. (u2)
                 .. controls (0.65,-1.08) and (1.18,-1.04) .. (w2)
                 .. controls (2.06,-0.56) and (2.04,1.52)  .. (w1)
                 .. controls (0.26,1.30)  and (0.22,0.82)  .. cycle;
  \fill[Fc] (u2) .. controls (-0.26,-1.12) and (-0.70,-1.25).. (u3)
                 .. controls (-0.67,-1.75) and (1.33,-1.82) .. (w2)
                 .. controls (1.18,-1.04)  and (0.65,-1.08) .. cycle;
  \fill[Fe] (u4) .. controls (-1.20,1.12) and (-0.72,1.02).. (u1)
                 .. controls (0.22,0.82)  and (0.26,1.30) .. (w1)
                 .. controls (0.01,1.95)  and (-1.35,1.98).. cycle;
  \draw[eqt] (-\R,0) arc[start angle=180, end angle=360, x radius=\R, y radius=0.50];
  \draw[eqt, dash pattern=on 0.9pt off 1.5pt]
        (\R,0) arc[start angle=0, end angle=180, x radius=\R, y radius=0.50];
  \draw[strand,coL] (u3) .. controls (-2.05,-1.00) and (-2.03,0.64) .. (u4);
  \draw[strand,coL] (u4) .. controls (-1.35,1.98)  and (0.01,1.95)  .. (w1);
  \draw[strand,coR] (w1) .. controls (2.04,1.52)   and (2.06,-0.56) .. (w2);
  \draw[strand,coL] (w2) .. controls (1.33,-1.82)  and (-0.67,-1.75).. (u3);
  \draw[strand,coL] (u4) .. controls (-1.20,1.12)  and (-0.72,1.02) .. (u1);
  \draw[strand,coR] (u1) .. controls (0.22,0.82)   and (0.26,1.30)  .. (w1);
  \draw[strand,coL] (u2) .. controls (-0.26,-1.12) and (-0.70,-1.25).. (u3);
  \draw[strand,coR] (u2) .. controls (0.65,-1.08)  and (1.18,-1.04) .. (w2);
  \draw[strand,coL] (u1) .. controls (0.30,0.24)   and (-0.30,-0.26).. (u2);
  \lk[79]{(u1)}{0.20}   \lk[139]{(u2)}{0.20} \lk[176]{(u3)}{0.20}
  \lk[-118]{(u4)}{0.20} \lk[275]{(w1)}{0.20} \lk[-45]{(w2)}{0.20}
  \draw[bd] (0,0) circle (\R);
  \coordinate (X0) at (-0.92,-0.02);
  \fill[white] (X0) circle (0.32);
  \draw[dp,rotate=-100] ($(X0)+(30:0.32)$) arc (30:210:0.32);
  \draw[dm,rotate=-100] ($(X0)+(210:0.32)$) arc (210:390:0.32);
  \node[tl] at (X0) {$X_0$};
  \node[fl] at (-1.6,-0.6) {$a$};
  \node[fl] at (1.02,0.42)   {$b$};
  \node[fl] at (0.42,-1.25)  {$c$};
  \node[fl] at (0.36,-1.88)  {$d$};
  \node[fl] at (-0.52,1.36)  {$e$};
\end{tkz}
\qquad\qquad\rightsquigarrow\qquad\qquad
\begin{tkz}[scale=0.75]
  \coordinate (P) at (-3.9,0);        \coordinate (Q) at (3.9,0);
  \coordinate (v1) at (-1.755,1.08);  \coordinate (v2) at (1.495,1.00);
  \coordinate (v3) at (2.275,-0.60);  \coordinate (v4) at (-1.365,-1.16);
  \coordinate (z1) at (-0.715,-0.12); \coordinate (z2) at (0.845,0.20);
  \fill[Fa] (P) .. controls (-1.27,-3.13) and (1.27,-3.13) .. (Q)
                .. controls ( 1.27, 3.13) and (-1.27, 3.13).. cycle;
  \fill[Fb] (v1) .. controls (-1.105,1.64) and (0.715,1.60) .. (v2)
                 .. controls (1.755,0.68)  and (1.235,0.52) .. (z2)
                 .. controls (0.325,0.44)  and (-0.195,-0.44).. (z1)
                 .. controls (-1.235,0.28) and (-2.015,0.60).. cycle;
  \fill[Fd] (v3) .. controls (1.755,-1.40) and (0.195,-1.64).. (v4)
                 .. controls (-1.495,-0.76) and (-1.105,-0.52).. (z1)
                 .. controls (-0.195,-0.44) and (0.325,0.44) .. (z2)
                 .. controls (1.495,-0.12) and (2.015,-0.12).. cycle;
  \fill[Fc] (v2) .. controls (2.665,0.84)  and (3.185,0.12) .. (v3)
                 .. controls (2.015,-0.12) and (1.495,-0.12).. (z2)
                 .. controls (1.235,0.52)  and (1.755,0.68) .. cycle;
  \fill[Fe] (v4) .. controls (-2.795,-0.84) and (-2.925,0.44).. (v1)
                 .. controls (-2.015,0.60)  and (-1.235,0.28).. (z1)
                 .. controls (-1.105,-0.52) and (-1.495,-0.76).. cycle;
  \draw[strand,coR] (v1) .. controls (-1.105,1.64)  and (0.715,1.60)  .. (v2);
  \draw[strand,coR] (v2) .. controls (2.665,0.84)   and (3.185,0.12)  .. (v3);
  \draw[strand,coR] (v3) .. controls (1.755,-1.40)  and (0.195,-1.64) .. (v4);
  \draw[strand,coR] (v4) .. controls (-2.795,-0.84) and (-2.925,0.44) .. (v1);
  \draw[strand,coL] (v1) .. controls (-2.015,0.60)  and (-1.235,0.28) .. (z1);
  \draw[strand,coR] (v4) .. controls (-1.495,-0.76) and (-1.105,-0.52).. (z1);
  \draw[strand,coL] (v2) .. controls (1.755,0.68)   and (1.235,0.52)  .. (z2);
  \draw[strand,coL] (v3) .. controls (2.015,-0.12)  and (1.495,-0.12) .. (z2);
  \draw[strand,coL] (z1) .. controls (-0.195,-0.44) and (0.325,0.44)  .. (z2);
  \draw[dm] (P) .. controls (-1.27,-3.13) and (1.27,-3.13) .. (Q);
  \draw[dp] (Q) .. controls ( 1.27, 3.13) and (-1.27, 3.13).. (P);
  \lk[51]{(v1)}{0.20}  \lk[316]{(v2)}{0.20} \lk[48]{(v3)}{0.20}
  \lk[135]{(v4)}{0.20} \lk[274]{(z1)}{0.20} \lk[97]{(z2)}{0.20}
  \node[fl] at (-3.05,0.45)  {$a$};
  \node[fl] at (0.07,0.92)   {$b$};
  \node[fl] at (1.92,0.24)   {$c$};
  \node[fl] at (0.26,-0.84)  {$d$};
  \node[fl] at (-1.85,0.12)  {$e$};
  \node[fl,anchor=north,Orange!80!black] at (0,-2.44) {\scriptsize$\partial_-D^2$};
  \node[Green!45,anchor=south]        at (0, 2.44) {\scriptsize$\partial_+D^2$};
\end{tkz}
\]
\end{construction}

\begin{lemma}
    \label{lemw-dn=2}
    The linear functional $\tld\psn^{\C^\fX}$ constructed in \cref{cstr:unitary disk-like 2-categories from proto-3-Hilbert spaces}
    is well-defined and descends to a linear functional $\psn^{\C^\fX}\colon \undfldl[2]{\C^\fX}{S^2}\to\bC$ satisfying \ref{2psnP}.
\end{lemma}

\begin{pf}
\itemstep{$\tld\psn^{\C^\fX}(\alpha)$ is independent of the splitting of $\partial X$.}
Suppose we instead chose $\partial X=\partial_-'X\cup\partial_+'X$. Let $\varphi\colon X\to X$ be any homeomorphism of $X$ sending the splitting $\partial X=\partial_-'X\cup\partial_+'X$ to the splitting $\partial X=\partial_-X\cup\partial_+X$. By the Alexander Trick there is an isotopy $h$ with $h_0=\varphi$ and $h_1=\id_X$. Since $(h_t|_{\partial X})_*(\alpha|_{\partial X})=(h_t|_{\partial X})_*(a_X\times\partial X)=a_X\times\partial X$ for all $t$, by \ref{2CTi2} $\varphi_*\alpha|_X-\alpha|_X\in\fldlU{\C^\fX}{X}=\ker\eval_X$, so
$\eval_{X,\partial_\pm X}(\varphi_*\alpha|_X-\alpha|_X)=0$.
Thus
\[
  \Psi_{a_X}(\eval_{X,\partial_\pm X}(\varphi_*\alpha|_X))=\Psi_{a_X}(\eval_{X,\partial_\pm X}(\alpha|_X)).
\]
The left side equals $\Psi_{a_X}(\eval_{X,\partial_\pm'X}(\alpha|_X))$, since $\varphi$ sends $\partial_\pm'X$ to $\partial_\pm X$ and $\Psi_{a_X}$ is invariant under rotations of (evaluations of) string diagrams thanks to the fact that $\vee$ is a UAF for $\fX$. 

\itemstep{$\tld\psn^{\C^\fX}$ is independent of the choice of $X$ (equivalently, of $X_0$).}
Suppose we instead choose a $2$-ball $X_0'\subset\Int(f_{X_0'})$ in some $2$-cell $f_{X_0'}\in F_\alpha$, set $X'\coloneq S^2\setminus X_0'$, and let $a_{X'}\in\fX$ be such that $\alpha|_{\partial X'}=a_{X'}\times\partial X'$.

\emph{Case 1: $X_0,X_0'\subset\Int(f_{X_0})$.} Choose an isotopy of $S^2$ sending $S^2\setminus X_0'$ to $S^2\setminus X_0$. As isotopy is a local relation, this preserves $\tld\psn$.

\emph{Case 2: $f_{X_0}$ and $f_{X_0'}$ are adjacent.} Then there is $e\in E_\alpha$ with $e$ adjacent to both $f_{X_0}$ and $f_{X_0'}$. Set $a_{X_0}\coloneq\lambda_e\in\fX(a\to b)$; without loss of generality we may assume $a=a_X$ and $b=a_{X'}$. By independence of the boundary splitting, there are 2-balls $\tld X\subset X$, $\tld X'\subset X'$ and boundary splittings $\partial\tld X=\partial_-\tld X\cup\partial_+\tld X$ and $\partial\tld X'=\partial_-\tld X'\cup\partial_+\tld X'$ such that
\begin{lst}
  \item\label{evall1} $\eval_{X',\partial_\pm X'}(\alpha|_{X'})=\tr_R^\vee(\eval_{\tld X',\partial_\pm\tld X'}(\alpha|_{\tld X'}))$,
  \item\label{evall2} $\eval_{X,\partial_\pm X}(\alpha|_X)=\tr_L^\vee(\eval_{\tld X,\partial_\pm\tld X}(\alpha|_{\tld X}))$, and
  \item\label{evall3} $\eval_{\tld X,\partial_\pm\tld X}(\alpha|_{\tld X})=\eval_{\tld X',\partial_\pm\tld X'}(\alpha|_{\tld X'})$.
\end{lst}
Thus
\begin{multline*}
  \Psi_{a_{X'}}(\eval_{X',\partial_\pm X'}(\alpha|_{X'}))
  \underset{\text{\ref{evall1}}}{=}\Psi_{a_{X'}}(\tr_R^\vee(\eval_{\tld X',\partial_\pm\tld X'}(\alpha|_{\tld X'})))
  \underset{\text{\ref{evall3}}}{=}\Psi_{a_{X'}}(\tr_R^\vee(\eval_{\tld X,\partial_\pm\tld X}(\alpha|_{\tld X})))\\
  =\Psi_{a_X}(\tr_L^\vee(\eval_{\tld X,\partial_\pm\tld X}(\alpha|_{\tld X})))
  \underset{\text{\ref{evall2}}}{=}\Psi_{a_X}(\eval_{X,\partial_\pm X}(\alpha|_X))
\end{multline*}
where the third equality uses sphericality of $\Psi$.

\emph{Case 3: $f_{X_0}$ and $f_{X_0'}$ are not adjacent.} 
Choose a path of 2-cells, i.e., a sequence $(X_0=X_0^1,X_0^2,\dots,X_0^N=X_0')$ of $2$-balls in $S^2$ such that for all $1\leq j\leq N-1$, $f_{X_0^j}$ and $f_{X_0^{j+1}}$ are adjacent $2$-cells in $\alpha$. Then repeat the argument in the previous case with $X_0$ and $X_0'$ replaced first by $X_0^1$ and $X_0^2$ respectively, then replace these by $X_0^2$ and $X_0^3$ respectively, and so on, which gives the result.

\itemstep{$\tld\psn^{\C^\fX}(\lU{\C^\fX})=0$.}
As $\lU{\C^\fX}=\ker\eval$, each $u\in\fldlU{\C^\fX}{S^2}$ has $\tld\psn(u)=\Psi_{a_X}(\eval(u))=\Psi_{a_X}(0)=0$.

\itemstep{$\psn^{\C^\fX}$ satisfies \ref{2psnP}.}
Let $X$ be a 2-ball, let $c\in\undfld[1]{\C^\fX}{\partial X}$, let $f\in\fldl[2]{\C^\fX}{X}[c]$, and write $f=\sum_i\lambda_i[\alpha_i]$ for string diagrams $\alpha_i\in\fld[2]{\C^\fX}{X}[c]$ such that $\fcj f\blt_c f$ lives on $S^2=\orev X\cup_{\partial X}X$; by linearity we can treat $f$ as a string diagram itself in what follows. Since $c$ has finitely many vertices, we can choose $\partial_{\pm}X$ so that $\partial_-X$ is vertex-free, say with single face label $a\in\fX$. Now choose $X_0$ to be any 2-ball neighborhood of $\partial_-X$ in $S^2$ disjoint from the underlying string diagram strata $V_{\alpha_i}\cup E_{\alpha_i}$ for each $i$. 
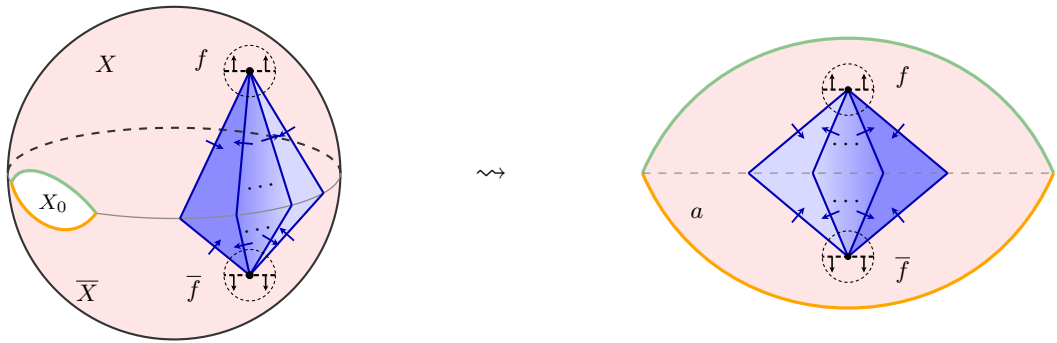
\[
\begin{tkz}
  \def\R{2.2}
  \coordinate (F)  at (1.00, 1.35);
  \coordinate (Fb) at (1.00,-1.35);
  \coordinate (w1) at (1.977,-0.263);
  \coordinate (w2) at (1.556,-0.424);
  \coordinate (w3) at (0.824,-0.556);
  \coordinate (w4) at (0.077,-0.600);
  \coordinate (p)  at (-2.152,-0.125);
  \coordinate (q)  at (-1.033,-0.530);
  \fill[Fa] (0,0) circle (\R);
  \fill[blue!15] (F) -- (w1) -- (Fb) -- (w2) -- cycle;
  \shade[left color=blue!45, right color=blue!15]
                 (F) -- (w2) -- (Fb) -- (w3) -- cycle;
  \fill[blue!45] (F) -- (w3) -- (Fb) -- (w4) -- cycle;
  \draw[bd] (0,0) circle (\R);
  \draw[bd,dashed] (\R,0) arc[start angle=0, end angle=180, x radius=\R, y radius=0.6];
  \draw[eqt] (-\R,0) arc[start angle=180, end angle=360, x radius=\R, y radius=0.6];
  \foreach \w in {w1,w3} {
    \draw[strand,coR] (F) -- (\w);
    \draw[strand,coR] (\w) -- (Fb);}
  \foreach \w in {w2,w4} {
    \draw[strand,coL] (F) -- (\w);
    \draw[strand,coL] (\w) -- (Fb);}
  \node[fl,rotate=10] at (1.17,-0.18) {$\cdots$};
  \node[fl,rotate=10] at (1.13,-0.75) {$\cdots$};
  \lk{(F)}{0.34}
  \lk[180]{(Fb)}{0.34}
  \node[fl] at (0.36,1.50) {$f$};
  \node[fl] at (0.26,-1.52) {$\fcj f$};
  \fill[white, opacity=1]
      (p) .. controls (-2.00,0.25) and (-1.36,-0.10) .. (q)
          .. controls (-1.40,-1.02) and (-2.04,-0.60) .. cycle;
  \draw[dp] (p) .. controls (-2.00,0.25) and (-1.36,-0.10) .. (q);
  \draw[dm] (q) .. controls (-1.40,-1.02) and (-2.04,-0.60) .. (p);
  \node[tl] at (-1.6,-0.4) {$X_0$};
  \node[fl] at (-0.90,1.45)  {$X$};
  \node[fl] at (-1.15,-1.55) {$\orev X$};
\end{tkz}
\qquad\qquad\rightsquigarrow\qquad\qquad
\begin{tkz}[scale=0.85]
  \coordinate (P)  at (-3.2,0);  \coordinate (Q)  at (3.2,0);
  \coordinate (c1) at (-1.55,0); \coordinate (c2) at (-0.55,0);
  \coordinate (c3) at ( 0.55,0); \coordinate (c4) at ( 1.55,0);
  \coordinate (F)  at (0,-1.30); \coordinate (Fb) at (0,1.30);
  \fill[Fa]
      (P) .. controls (-2.0,-2.8) and (2.0,-2.8) .. (Q)
          .. controls ( 2.0, 2.8) and (-2.0, 2.8) .. cycle;
  \fill[blue!15] (F) -- (c1) -- (c2) -- cycle;
  \fill[blue!15] (Fb) -- (c1) -- (c2) -- cycle;
  \shade[left color=blue!15, right color=blue!45] (F) -- (c2) -- (c3) -- cycle;
  \shade[left color=blue!15, right color=blue!45] (Fb) -- (c2) -- (c3) -- cycle;
  \fill[blue!45] (F) -- (c3) -- (c4) -- cycle;
  \fill[blue!45] (Fb) -- (c3) -- (c4) -- cycle;
  \draw[dm] (P) .. controls (-2.0,-2.8) and (2.0,-2.8) .. (Q);
  \draw[dp] (Q) .. controls ( 2.0, 2.8) and (-2.0, 2.8) .. (P);
  \draw[dashed,gray] (P) -- (Q);
  \foreach \c in {c1,c3} {
    \draw[strand,coR] (F) -- (\c);
    \draw[strand,coR] (\c) -- (Fb);}
  \foreach \c in {c2,c4} {
    \draw[strand,coL] (F) -- (\c);
    \draw[strand,coL] (\c) -- (Fb);}
  \lk[180]{(F)}{0.40}
  \lk{(Fb)}{0.40}
  \node[fl] at (0.85,-1.50) {$\fcj f$};
  \node[fl] at (0.85, 1.50) {$f$};
  \node[fl] at (-2.35,-0.62) {$a$};
  \node[fl] at (0,-0.45) {$\cdots$}; \node[fl] at (0,0.45) {$\cdots$};
\end{tkz}
\]
Splitting $\partial X_0$ as illustrated, we find
\[
    \bkt{f}{f}_{X,c}=\psn^{\C^\fX}(\fcj f\blt_{c}f)\defeq{}\Psi_a(\eval(f)^\dag\circ\eval(f)),
\]
which by faithfulness and positivity of $\Psi_a^\fX$ is nonnegative and equals zero if and only if $\eval(f)=0$, i.e., if and only if $f=0$. Thus $\bkt--_{X,c}$ is positive-definite.
\end{pf}

\begin{notation}
\label{not:cl-a}
For $a\in\fX$ and $f\in \End_\fX(1_a)$, write $\cl_a(f)\coloneq\cl_{a}(\langle f\rangle)\in\undfldl[2]{\C^\fX}{S^2}$ (\cref{def:closure}) where $\langle f\rangle\in\fldl[2]{\C^\fX}{D^2}[a\times S^1]$ is the cone $2$-field from \cref{anglebracket}.
\end{notation}

\begin{lemma} 
\label{lem:CX-finiteness}
For a finite proto-3-Hilbert space $\fX$, the unitary disk-like 2-category $\C^\fX$ is finite.
\end{lemma}

\begin{pf}
By \cref{cor:psnF-weaker}, it suffices to show this when $X$ is the 2-disk $D^2$, the annulus $S^1\times D^1$, and the 2-sphere $S^2$.

\itemstep{The 2-disk $D^2$.}
As $\fldl[2]{\C^\fX}{D^2}[c]$ is a vector space of 2-morphisms in $\fX$---choose a splitting of $c$ as $c=c_1\blt c_2$ (which exists by \ref{2CS2}) into source and target and a homeomorphism $D^2\to D^2$ sending $c_1$ and $c_2$ to the incoming and outgoing (lower and upper) boundaries of $D^2$---it is finite-dimensional.

\itemstep{The annulus $A\coloneq S^1\times D^1$.} 
Fix $c = \fcj{c_-}\amalg c_+\in\undfld[1]{\C^\fX}{\partial A}$. Choose a properly embedded arc $Y\hookrightarrow A$ transversely intersecting each of the two boundary circles of $A$ exactly once, say with endpoints landing in regions labeled by the objects $a,b\in\fX$, and consider the gluing data that forms $A$ by gluing the rectangle $X\coloneq A\setminus Y$ along $Y$ and $\orev{Y}$, so that $X_\glu=A$ and $X\cong D^2$. By \cref{eval-facts-2}, $\Sk_{\C^\fX}(Y,\fcj a\amalg b)\cong^\dag\fX(a\to b)$ is a finite pre-2-Hilbert space since $\fX$ is finite. Write $\Irr(\fX(a\to b)^\cent)=\{(r_i,p_i)\}_{i=1}^N$. Then \cref{glu-span} gives $\undfldl[2]{\C^\fX}{A}[c]=\sum_{i=1}^N\Gamma(p_iV_{r_i}p_i)$, so $\dim\undfldl[2]{\C^\fX}{A}[c]\leq\sum_{i=1}^N\dim (p_iV_{r_i}p_i)<\infty$ where each $V_{r_i}$ is finite-dimensional again by the $D^2$-case above, as $X\cong D^2$.

\itemstep{The 2-sphere $S^2$.}
It suffices to show $\dim\undfldl[2]{\C^\fX}{S^2}\leq|\pi_0\fX|$, as the right side is finite by finiteness of $\fX$. Since $\End_\fX(1_a)$ is commutative, $\End_\fX(1_a)\cong\bC^{m_a}$ for some $m_a\in\bZ_{\geq 0}$. Observe that $\End_\fX(1_a)$ is spanned by its set $P_a$ of minimal projections, which satisfy $\Psi_a(p)>0$ for all $p\in P_a$ by faithfulness and positivity of $\Psi$. Since $\undfldl[2]{\C^\fX}{S^2}$ is spanned by $\{\cl_a(p)\mid a\in\fX,p\in P_a\}$, it suffices to show for any finite set of objects $a_1,\dots,a_\ell$ that $\dim\bC\{\cl_{a_i}(p)\mid 1\leq i\leq\ell,\,p\in P_{a_i}\}\leq|\pi_0\fX|$.
We claim that for any $1\leq i,j\leq \ell$, $p\in P_{a_i}$, $q\in P_{a_j}$, and $X\in\fX(a_i\to a_j)$ with $p\otimes\id_X\otimes q\neq 0$, one has $\cl_{a_i}(p)\in\bC\,\cl_{a_j}(q)$. First note $\tr^\vee_L(p\otimes \id_X\otimes q)=p\circ\tr^\vee_L(\id_X\otimes q)\in p\End_\fX(1_{a_i})=\bC p$ by minimality (this is where we use commutativity), and similarly $\tr^\vee_R(p\otimes \id_X\otimes q)\in\bC q$. Thus $\tr^\vee_L(p\otimes \id_X\otimes q)=\mu p$ and $\tr^\vee_R(p\otimes \id_X\otimes q)=\nu q$ for some $\mu,\nu\in\bC$. As $p\otimes \id_X\otimes q$ is a nonzero projection, $\mu\Psi_{a_i}(p)=\Psi_{a_i}(\tr^\vee_L(p\otimes \id_X\otimes q))>0$ and $\nu\Psi_{a_{j}}(q)=\Psi_{a_j}(\tr^\vee_R(p\otimes \id_X\otimes q))>0$. As $p$ and $q$ are also nonzero projections, $\Psi_{a_i}(p)>0$ and $\Psi_{a_{j}}(q)>0$. Thus both $\mu$ and $\nu$ must be nonzero, so \cref{cor:cl-slide} gives $\cl_{a_i}(p)\in\bC\,\cl_{a_{j}}(q)$, i.e., the closures of all points in a fixed component of $\fX$ span a 1D vector space. It follows that $\dim\undfldl[2]{\C^\fX}{S^2}\leq |\pi_0\fX|$.
\end{pf}

\subsection{The round-trip equivalences}

\subsubsection{The first round-trip: \texorpdfstring{$\fX_{\C^\fX}\cong^\dag\fX$}{X(C(X)) = X}}

\begin{convention} 
\label{strictnessconvention}
By the strictification result \cref{strictification}, we may always assume our dagger 2-categories are strict and that each induces a strict pivotal structure; see \cref{strictification-sec} for more details. Thus we will henceforth assume such strictness when we say ``pivotal dagger 2-category.''
\end{convention}

\begin{proposition} 
\label{XCXcongX-n-2}
For a pivotal $\Cstar$-2-category $\fX$, the evaluation maps from \cref{eval-facts-2} assemble into a strictly UAF-preserving strict $\dag$-functor $\eval\colon\fX_{\C^\fX}\to\fX$. When $\fX$ is equipped with a spherical weight $\Psi$, $\eval$ is an isometric equivalence of proto-3-Hilbert spaces $\eval\colon \fX_{\C^\fX}\overset{\cong^\dag}\to\fX$.
\end{proposition} 

\begin{pf}
Define $\eval$ on objects $a\in\fX_{\C^\fX}$ by $\eval(a)\coloneq a$, on 1-morphisms $\xi\in\fX_{\C^\fX}(a\to b)$ by $\eval(\xi)\coloneq\eval^\fX_{I}(\xi)$, and on 2-morphisms $\alpha\in\fX_{\C^\fX}({}_a\xi_b\Rightarrow{}_a\eta_b)$ by
$\eval(\alpha)\coloneq\eval^\fX_{D^2,\partial_{\pm}D^2}(\alpha)$. This is a strict dagger 2-functor by \cref{eval-facts-2}. Moreover, $\eval$ is UAF-preserving: for ${}_a\xi_b\in\fX_{\C^\fX}(a\to b)$, we have
\[
\eval(\xi^\vee)=\eval(\fcj{\iota_*\xi})=(\lambda'_{v_m})^\vee\otimes(\lambda'_{v_{m-1}})^\vee\otimes\cdots\otimes(\lambda'_{v_1})^\vee=\eval(\xi)^\vee.
\]
To see $\eval$ is essentially surjective on 1-morphisms, observe that any single-vertex string diagram $\xi_X$ on $I$ with positively oriented vertex $v$ labeled by a 1-morphism $X$ has $\eval(\xi_X)=X$. Finally, \cref{eval-facts-2}\ref{fact:eval-ker2} shows $\eval$ is fully faithful. 

Now suppose $\fX$ has a spherical weight $\Psi$. Since $\eval$ is the identity on objects, it is trivially isometrically essentially surjective, so we need only show $\eval$ preserves the spherical weights. For $a\in\fX_{\C^{\fX}}$ and $\alpha\in\fX_{\C^{\fX}}(1_a\Rightarrow 1_a)$,
\[
\Psi^{\fX_{\C^\fX}}_a(\alpha)\underset{\eqref{cstr:proto-3-Hilbert spaces from unitary disk-like 2-categories}}{\defeq{}}\psn^{\C^\fX}(\cl_a(\alpha))\underset{\eqref{cstr:unitary disk-like 2-categories from proto-3-Hilbert spaces}}{\defeq{}}\Psi^\fX_a(\eval(\alpha)).\qedhere
\]
\end{pf}

\subsubsection{The second round-trip: \texorpdfstring{$\C^{\fX_\C}\cong^\dag\C$}{C(X(C)) = C}}

Unlike in the $n=1$ exposition, here we refer the reader to \cref{sec:unnrestricted} for the definition of an unrestricted disk-like functor and isometric weak equivalence.

Recall (see \cref{sec:liftable}) that $\tld\C^{\fX_\C}$ is the disk-like 2-category obtained from $\C^{\fX_\C}$ by only considering string diagrams on 2-balls whose vertices are labeled by $2$-morphisms of $\fX_\C$ that admit pure $2$-field representatives in $\C$. We call such $2$-morphisms \defn{liftable}. By \cref{lem:tilde-inclusion}, the inclusion $\iota\colon\widetilde{\C}^{\fX_\C}\to\C^{\fX_\C}$ is an isometric weak equivalence.

Fix a disk-like 2-category $\C$. For $k\in\{1,2\}$, denote by $\iota^{(k)}\colon D^k\to\orev{D^k}$ the orientation-preserving homeomorphism given by $(x_1,\ldots,x_k)\mapsto(x_1,\ldots,-x_k)$ and define the involutions $\sigma^{\dag_k}\coloneq\fcj{\iota^{(k)}_{*}\sigma}=\orev{\iota^{(k)}}_*\fcj\sigma\in\fld[k]{\C}{D^k}$, which satisfy $\iota^{(k)}_{*}\sigma^{\dag_k}=\fcj\sigma$ and $[\sigma^{\dag_k}]=[\sigma]^{\dag_k}$. 
For each liftable 2-morphism $f$ of $\fX_\C$, choose a pure 2-field $\sigma_f\in\fld[2]{\C}{D^2}$ with $[\sigma_f]=f$ and for each vertex $v$ of a 2D string diagram $\alpha$ set $\sigma_v\coloneq\sigma_{\eval_{N_v,\partial_{\pm}N_v}(\alpha|_{N_v})}$ where the splitting $\partial_{\pm}N_v$ is the image of the standard splitting $\partial_\pm D^2$ under the homeomorphism $\theta_v\colon D^2\to N_v$. By the same computation as in \cref{vertid} but with $\dag$ replaced with $\dag_2$, we have $(\orev{\theta_v}\circ\iota^{(2)})_*\sigma_v^{\dag_2}=\fcj{(\theta_v)_*\sigma_v}$. Also observe that for a 1-field $\sigma_e\in\fld[1]{\C}{D^1}$, a parameterization $\theta_e\colon D^1\times\hat e\to N_e$ of the truncated edge $\hat e$, and the canonical projection $\pi_e\colon D^1\times\hat e\to D^1$, we have
\begin{equation*} 
\begin{multlined}
\fcj{(\theta_e)_*\pi_{e}^*\sigma_e}
\underset{\ref{2Cfcj2}}{=}
\orev{\theta_e}_*\fcj{\pi_{e}^*\sigma_e}
\underset{\ref{2Cfcj2}}{=}
\orev{\theta_e}_*(\orev{\pi_{e}})^*\fcj{\sigma_e}
=
\orev{\theta_e}_*(\orev{\pi_{e}})^*\iota^{(1)}_*\sigma_e^{\dag_1}
\\
\qquad\qquad\underset{\ref{2Cpi2}}{=}
\orev{\theta_e}_*(\iota^{(1)}\times\id_{\hat e})_*\pi_{e}^*\sigma_e^{\dag_1}
\underset{\ref{2Cvarphi2}}{=}
(\orev{\theta_e}\circ(\iota^{(1)}\times\id_{\hat e}))_*\pi_{e}^*\sigma_e^{\dag_1}
\end{multlined}
\end{equation*}
where the fourth equality uses that $\orev{\pi_{e}}\circ(\iota^{(1)}\times\id_{\hat e})=\iota^{(1)}\circ\pi_{e}$ as maps $D^1\times\hat e\to\orev{D^1}$.

\begin{construction}
Define an unrestricted disk-like functor $\eL\colon\tld\C^{\fX_\C}\to\C$ on 0-fields by $\eL(a)\coloneq a$, on 1-fields exactly as in \cref{eLn1construction} but with $\dag$ replaced with $\dag_1$, and on 2-fields $\alpha\in\fld[2]{\tld\C^{\fX_\C}}{X}$ by the assignment $\mathrm{PS}(\alpha,\Gamma_\alpha)\to\fld[2]{\C}{X}$ given by
\[
\eL(\alpha,\Theta)
\coloneq
\glu
\{
\{(\theta_v)_*\sigma_v\}_{v\in V_\alpha}
\,\cup\,
\{(\theta_e)_*\pi_{e}^*\sigma_e\}_{e\in E_\alpha}\cup\{\pi_f^*\lambda_f\}_{f\in F_\alpha}
\}
\]
where
\begin{itemize}
    \item for $v\in V_\alpha$, $\sigma_v$ is defined above;
    \item for $e\in E_\alpha$, $\sigma_e\coloneq\sigma_{\lambda_e}^{(\varepsilon_e)}$ for $\sigma_{\lambda_e}\coloneq\lambda_e$ (a 1-morphism of $\fX_\C$ being a pure 1-field of $\C$ on $D^1$), $\sigma_{\lambda}^{(+1)}\coloneq \sigma_\lambda$ and $\sigma_{\lambda}^{(-1)}\coloneq \sigma_\lambda^{\dag_1}$, and where $\pi_e\colon D^1\times\hat{e}\to D^1$ for $e\in E_\alpha$ is the canonical projection; and
    \item for $f\in F_\alpha$, $\pi_f\colon\hat{f}\to\pt$ is the unique pinched product map from the truncated face $\hat f$.
\end{itemize}
For example,
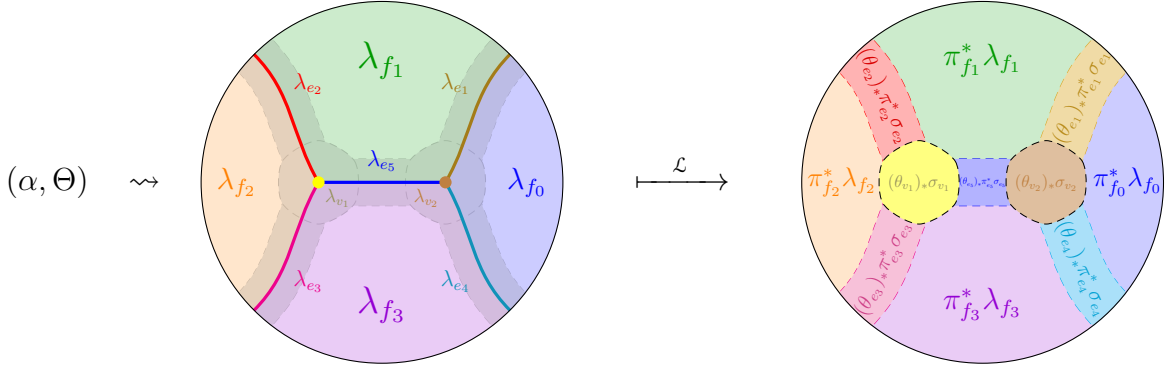
\[
(\alpha,\Theta)
\quad\rightsquigarrow\quad
\begin{tkz}[scale=1.2,transform shape]
    \clip circle (2);
    \fill[green!60!black!20] (0.7,0) to[out=60, in=225] (45:2) arc (45:135:2) to[out=315, in=120] (-0.7,0) to[out=0, in=180] (0.7,0) -- cycle;
    \fill[orange!20] (-0.7,0) to[out=120, in=315] (135:2) arc (135:225:2) to[out=45, in=240] (-0.7,0) -- cycle;
    \fill[violet!20] (-0.7,0) to[out=240, in=45] (225:2) arc (225:315:2) to[out=135, in=300] (0.7,0) to[out=180, in=0] (-0.7,0) -- cycle;
    \fill[blue!20] (0.7,0) to[out=300, in=135] (315:2) arc (315:405:2) to[out=225, in=60] (0.7,0) -- cycle;
    \begin{scope}[transparency group, opacity=0.2]
        \draw[line width=18pt, densely dashed, draw=black] (0.7,0) to[out=60, in=225] (45:2);
        \draw[line width=18pt, densely dashed, draw=black] (-0.7,0) to[out=120, in=315] (135:2);
        \draw[line width=18pt, densely dashed, draw=black] (-0.7,0) to[out=240, in=45] (225:2);
        \draw[line width=18pt, densely dashed, draw=black] (0.7,0) to[out=300, in=135] (315:2);
        \draw[line width=18pt, densely dashed, draw=black] (-0.7,0) to[out=0, in=180] (0.7,0);
        \draw[line width=17.6pt, solid, draw=gray] (0.7,0) to[out=60, in=225] (45:2);
        \draw[line width=17.6pt, solid, draw=gray] (-0.7,0) to[out=120, in=315] (135:2);
        \draw[line width=17.6pt, solid, draw=gray] (-0.7,0) to[out=240, in=45] (225:2);
        \draw[line width=17.6pt, solid, draw=gray] (0.7,0) to[out=300, in=135] (315:2);
        \draw[line width=17.6pt, solid, draw=gray] (-0.7,0) to[out=0, in=180] (0.7,0);
        \filldraw[fill=gray, draw=black, densely dashed, very thin] plot[smooth cycle] coordinates {(-0.7, 0.46) (-0.3, 0.23) (-0.3, -0.23) (-0.7, -0.46) (-1.1, -0.23) (-1.1, 0.23)};
        \filldraw[fill=gray, draw=black, densely dashed, very thin] plot[smooth cycle] coordinates {(0.7, 0.46) (1.1, 0.23) (1.1, -0.23) (0.7, -0.46) (0.3, -0.23) (0.3, 0.23)};
    \end{scope}
    \draw[very thick, Goldenrod!75!black] (0.7,0) to[out=60, in=225] node[pos=0.6, above left=-2pt] {\scalebox{0.6}{$\lambda_{e_1}$}} (45:2);
    \draw[very thick, red] (-0.7,0) to[out=120, in=315] node[pos=0.6, above right=-2pt] {\scalebox{0.6}{$\lambda_{e_2}$}} (135:2);
    \draw[very thick, magenta] (-0.7,0) to[out=240, in=45] node[pos=0.6, below right=-2pt] {\scalebox{0.6}{$\lambda_{e_3}$}} (225:2);
    \draw[very thick, cyan!70!black] (0.7,0) to[out=300, in=135] node[pos=0.6, below left=-2pt] {\scalebox{0.6}{$\lambda_{e_4}$}} (315:2);
    \draw[very thick, blue] (-0.7,0) to[out=0, in=180] node[pos=0.5, above=-1pt] {\scalebox{0.6}{$\lambda_{e_5}$}} (0.7,0);
    \draw[thick] (0,0) circle (2);
    \node[circle, fill=black, inner sep=1.35pt, fill=Yellow] at (-0.7,0) {};
    \node[circle, fill=black, inner sep=1.35pt, fill=brown] at (0.7,0) {};
    \node[text=green!60!black] at (0, 1.35) {$\lambda_{f_1}$};
    \node[text=blue] at (1.6, 0) {\scalebox{0.8}{$\lambda_{f_0}$}};
    \node[text=orange] at (-1.6, 0) {\scalebox{0.8}{$\lambda_{f_2}$}};
    \node[text=violet] at (0, -1.35) {$\lambda_{f_3}$};
    \node[text=yellow!50!black, below right=-2pt] at (-0.7, 0) {\scalebox{0.5}{$\lambda_{v_1}$}};
    \node[text=brown, below left=-2pt] at (0.7, 0) {\scalebox{0.5}{$\lambda_{v_2}$}};
\end{tkz}
\qquad\xmapsto{\quad\eL\quad}\qquad
\begin{tkz}[scale=1.2,transform shape]
    \clip circle (2);
    \fill[green!60!black!20] (0.7,0) to[out=60, in=225] (45:2) arc (45:135:2) to[out=315, in=120] (-0.7,0) to[out=0, in=180] (0.7,0) -- cycle;
    \fill[orange!20] (-0.7,0) to[out=120, in=315] (135:2) arc (135:225:2) to[out=45, in=240] (-0.7,0) -- cycle;
    \fill[violet!20] (-0.7,0) to[out=240, in=45] (225:2) arc (225:315:2) to[out=135, in=300] (0.7,0) to[out=180, in=0] (-0.7,0) -- cycle;
    \fill[blue!20] (0.7,0) to[out=300, in=135] (315:2) arc (315:405:2) to[out=225, in=60] (0.7,0) -- cycle;
    \draw[line width=18pt, densely dashed, draw=Goldenrod] (0.7,0) to[out=60, in=225] (45:2);
    \draw[line width=17.6pt, solid, draw=Goldenrod!40] (0.7,0) to[out=60, in=225] (45:2);
    \draw[line width=18pt, densely dashed, draw=red] (-0.7,0) to[out=120, in=315] (135:2);
    \draw[line width=17.6pt, solid, draw=red!30] (-0.7,0) to[out=120, in=315] (135:2);
    \draw[line width=18pt, densely dashed, draw=magenta] (-0.7,0) to[out=240, in=45] (225:2);
    \draw[line width=17.6pt, solid, draw=magenta!30] (-0.7,0) to[out=240, in=45] (225:2);
    \draw[line width=18pt, densely dashed, draw=cyan] (0.7,0) to[out=300, in=135] (315:2);
    \draw[line width=17.6pt, solid, draw=cyan!30] (0.7,0) to[out=300, in=135] (315:2);
    \draw[line width=18pt, densely dashed, draw=blue] (-0.7,0) to[out=0, in=180] (0.7,0);
    \draw[line width=17.6pt, solid, draw=blue!30] (-0.7,0) to[out=0, in=180] (0.7,0);
    \filldraw[fill=Yellow!50, draw=black, densely dashed, very thin] plot[smooth cycle] coordinates {(-0.7, 0.46) (-0.3, 0.23) (-0.3, -0.23) (-0.7, -0.46) (-1.1, -0.23) (-1.1, 0.23)};
    \filldraw[fill=brown!50, draw=black, densely dashed, very thin] plot[smooth cycle] coordinates {(0.7, 0.46) (1.1, 0.23) (1.1, -0.23) (0.7, -0.46) (0.3, -0.23) (0.3, 0.23)};
    \draw[thick] (0,0) circle (2);
    \node[text=green!60!black] at (0, 1.35) {\scalebox{0.8}{$\pi_{f_1}^*\lambda_{f_1}$}};
    \node[text=blue] at (1.6, 0) {\scalebox{0.7}{$\pi_{f_0}^*\lambda_{f_0}$}};
    \node[text=orange] at (-1.55, 0) {\scalebox{0.7}{$\pi_{f_2}^*\lambda_{f_2}$}};
    \node[text=violet] at (0, -1.35) {\scalebox{0.8}{$\pi_{f_3}^*\lambda_{f_3}$}};
    \path (0.7,0) to[out=60, in=225] node[pos=0.65, sloped, text=Goldenrod!85!black] {\scalebox{0.55}{$(\theta_{e_1})_*\pi_{e_1}^*\sigma_{e_1}$}} (45:2);
    \path (-0.7,0) to[out=120, in=315] node[pos=0.65, sloped, text=red!85!black] {\scalebox{0.55}{$(\theta_{e_2})_*\pi_{e_2}^*\sigma_{e_2}$}} (135:2);
    \path (-0.7,0) to[out=240, in=45] node[pos=0.65, sloped, text=magenta!85!black] {\scalebox{0.55}{$(\theta_{e_3})_*\pi_{e_3}^*\sigma_{e_3}$}} (225:2);
    \path (0.7,0) to[out=300, in=135] node[pos=0.65, sloped, text=cyan!85!black] {\scalebox{0.55}{$(\theta_{e_4})_*\pi_{e_4}^*\sigma_{e_4}$}} (315:2);
    \path (-0.65,0) to[out=0, in=180] node[pos=0.5, sloped, text=blue] {\scalebox{0.25}{$(\theta_{e_5})_*\pi_{e_5}^*\sigma_{e_5}$}} (0.65,0);
    \node[text=yellow!50!black] at (-0.7, 0) {\scalebox{0.45}{$(\theta_{v_1})_*\sigma_{v_1}$}};
    \node[text=brown] at (0.7, 0) {\scalebox{0.45}{$(\theta_{v_2})_*\sigma_{v_2}$}};
\end{tkz}
\]
where we have suppressed the normal framings (the transverse orientations on the edges and the link parameterizations of the vertices) of the string diagram strata. 
\end{construction}

Verification that $\eL$ is an unrestricted disk-like functor is a routine computation as in the $n=1$ case and is left to the reader.

\begin{proposition}
\label{CXCcongC-n2}
For a unitary disk-like 2-category $\C$, there is an isometric weak equivalence $\C\cong^\dag\C^{\fX_\C}$ witnessed by a zig-zag $\C\overset\eL\leftarrow\widetilde{\C}^{\fX_\C}\overset\iota\rightarrow\C^{\fX_\C}$.
\end{proposition}

\begin{pf}
First note $\eL$ is surjective on objects. To see $\eL$ is unitarily essentially surjective on 1-fields, fix $\Xi\in\eL_{\mathrm{PS}}(c)$ and a 1-field $\eta\in\fld[1]{\C}{Y}[\eL(c,\Xi)]$ and let $\xi_\eta \in \fld[1]{\widetilde{\C}^{\fX_\C}}{Y}[c]$ (assuming $Y=D^1$ by \cref{uresequvgen}) be the single-vertex 1-field whose vertex $v$ is labeled by the 1-morphism $\eta$ of $\fX_\C$. Then for any $\Theta \in \eL_{\mathrm{PS}}(\xi_\eta)$ with $\partial\Theta=\Xi$, the 1-field $\eL(\xi_\eta, \Theta) \in \fld[1]{\C}{Y}[ \eL(c,\Xi)]$ is obtained by gluing product fields to the boundary of the 1-field $\theta_{v*}\eta$, so in particular $\eta$ and $\eL(\xi_\eta, \Theta)$ are related by an extended isotopy (namely a collaring map), whence $\eL(\xi_\eta, \Theta) \cong^\star \eta$ in $\C$ by \ref{CTc2}. Finally, $\widehat{\eL}_{X,c} \colon \fldl[2]{\widetilde{\C}^{\fX_\C}}{X}[c] \to \fldl[2]{\C}{X}[\eL(c)]$ is an isomorphism: injectivity holds because $\ker\widehat{\eL}_{X,c} = \fldlU{\widetilde{\C}^{\fX_\C}}{X}[c] = \ker(\eval \circ \iota)$, and surjectivity holds because for any liftable $f \in \fldl[2]{\C}{X}[\eL(c)]$, any single-vertex cone diagram $\alpha_f$ has $[\eL(\alpha_f,\Theta)] = f$ by extended-isotopy invariance \ref{CTc2}, and liftable 2-fields span by \cref{facts:liftable-properties}\ref{facts:liftable-properties-1}.

To see $\eL$ is isometric, let $\alpha$ be an $\fX_\C$-string diagram on $S^2$. Choose a 2-ball $X_0 \hookrightarrow S^2 \setminus (V_\alpha \cup E_\alpha)$ in an $a$-labeled face $f_0$, set $f \coloneq \eval(\alpha|_{S^2 \setminus X_0}) \in \fX_\C(1_a \Rightarrow 1_a)$ with lift $\sigma_f$, and let $\alpha_P$ be the single-vertex string diagram with vertex label $f$, so that $[\alpha] = [\alpha_P]$. Then for any parameterized splitting $\Theta \in \eL_{\mathrm{PS}}(\alpha)$, the parameterized splitting $\Theta_P = \{\theta_v \colon D^2 \to N_v,\ \id_{S^2 \setminus N_v}\}$ for $\alpha_P$ satisfies
\[
\eL(\alpha,\Theta) 
= \eL(\alpha_P,\Theta_P)
= \glu\{(\theta_v)_*\sigma_f, \pi_{f_0}^* a\}
= \cl_a(\sigma_f)
\]
in $\undfldl[2]{\C}{S^2}$, so
applying $\psn^\C$ gives
$
\psn^\C([\eL(\alpha,\Theta)]) 
% \\& 
= \psn^\C(\cl_a(\sigma_f))
% \\& 
= \Psi^{\fX_\C}_a(f)
% \\& 
= \psn^{\C^{\fX_\C}}([\alpha])
$.
\end{pf}

\subsection{Representations and Hom-spaces}

\subsubsection{The spherical weight on functor 2-categories}

\begin{example}[Spherical weight on functor 2-categories {\cite[Defn. 5.1]{bases}}]
  \label{spherical weight on functor 2-categories}
  Suppose $\fX$ and $\fY$ are finite 3-Hilbert spaces. On the unitary 2-category $\Hom(\fX\to \fY)$ of UAF-preserving dagger functors, we define
  \begin{equation}
    \label{eq:3HilbOnFun}
    \Psi^{\Hom}_F(m\colon \id_F\Rrightarrow \id_F)
    \coloneq
    \sum_{b\in\pi_0\fX}
    \frac{d_b}{D_{\Omega_b}}
    \Psi^\fY_{F(b)}(m_b)
    \qquad\quad
    \forall\,F\in\Hom(\fX\to \fY)
  \end{equation}
where $d_b\coloneq d_{1_b} = \Psi^\fX_b(\id_{1_b})$, $D_{\Omega_b} \coloneq \sum_{X\in \Irr(\Omega_b)} d_X^2$, and $\pi_0\fX$ denotes a choice of one simple object from each component of $\fX$ (see \cref{def:connected}).
Equipped with the UAF from \cite[Defn. 5.1]{bases}, $\Psi^{\Hom}$ makes $\Hom(\fX\to\fY)$ a finite 3-Hilbert space \cite[Lem. 5.3 and Prop. 5.4]{bases}.
\end{example}

The following records some immediate consequences of our above work.

\begin{corollary}\label{lem:local-cauchy-rep}
  For a finite proto-3-Hilbert space $(\fX,\vee,\Psi)$ and a pre-3-Hilbert space $(\fY,\vee,\Psi^\fY)$, restriction along $\iota_{1\cent}\colon\fX\hookrightarrow\fX^{1\cent}$ induces an equivalence of unitary 2-categories
  \[
    \iota_{1\cent}^*\colon\Hom(\fX^{1\cent}\to\fY) \xrightarrow{\;\sim\;} \Hom(\fX\to\fY).
  \]
\end{corollary}

\begin{corollary}\label{cor:proto-hom-weight}
  For finite proto-3-Hilbert spaces $\fX$ and $\fY$, the inclusion $\Hom(\fX\to\fY)\hookrightarrow\Hom(\fX\to\fY^\cent)\cong^\dag\Hom(\fX^\cent\to\fY^\cent)$ restricts the structure of \cref{spherical weight on functor 2-categories} to a finite proto-3-Hilbert space structure on $\Hom(\fX\to\fY)$, given by $\Psi^{\Hom}_F(m)=\sum_{b\in\pi_0\fX^\cent}\frac{d_b}{D_{\Omega_b}}\Psi^{\fY^\cent}_{F^\cent(b)}(\tld m_b)$ where $F^\cent$ and $\tld m$ are the extensions of $F$ and $m$ respectively from the universal property of completion \cite[\S4]{3Hilb}.
\end{corollary}

\begin{corollary}
\label{cor:XCXcongX implies homs-2}
For finite proto-3-Hilbert spaces $\fX$ and $\fY$, we have isometric equivalences of finite proto-3-Hilbert spaces
\[
\Hom(\fX\to\fY)\cong^\dag\Hom(\fX_{\C^\fX}\to\fY)\cong^\dag\Hom(\fX\to\fX_{\C^\fY})\cong^\dag\Hom(\fX_{\C^\fX}\to\fX_{\C^\fY}).
\]
\end{corollary}

\subsubsection{Representations and the unitary Yoneda embedding}

\begin{corollary}\label{prop:rep-is-3hilb}
  If $(\fX,\vee,\Psi)$ is a finite proto-3-Hilbert space, then $\Rep(\fX)\coloneq\Hom(\fX\to 2\Hilb)$ is a finite 3-Hilbert space.
\end{corollary}

\begin{pf}
  By \cref{cstr:local-cauchy-3hilb}\ref{lcc2}, the local completion $\fX^{1\cent}$ is a finite pre-3-Hilbert space, so by \cref{lem:local-cauchy-rep} we have an equivalence of unitary 2-categories
  \[
    \iota_{1\cent}^*\colon\Hom(\fX^{1\cent}\to 2\Hilb) \xrightarrow{\;\sim\;} \Hom(\fX\to2\Hilb).
  \]
  Since $\fX^{1\cent}$ is a finite pre-3-Hilbert space, \cite[Props. 4.28 and 4.43]{3Hilb} give an equivalence of unitary 2-categories
  \[
    \iota_\cent^*\colon\Hom(\fX^\cent\to 2\Hilb) \xrightarrow{\;\sim\;} \Hom(\fX^{1\cent}\to 2\Hilb)
  \]
  where $\fX^\cent\coloneq\mathsf{H}^*\mathsf{Alg}(\Hilb_\boxplus(\fX^{1\cent}))$ is a finite 3-Hilbert space by \cite[Lem. 4.42 and Prop. 4.44]{3Hilb} and \cref{rem:points-completion}. By \cref{spherical weight on functor 2-categories}, $\Hom(\fX^\cent\to 2\Hilb)$ has a canonical 3-Hilbert space structure with spherical weight $\Psi^{\Hom^\cent}$. By construction, restriction along $\iota\coloneq\iota_\cent\circ\iota_{1\cent}$ is UAF-preserving, so by \cite[Lem. 4.23]{3Hilb} $\iota^*$ preserves Hilbert direct sums. Finally, since $\Hstar$-monads in $\Hom(\fX\to2\Hilb)$ split pointwise (in the 3-Hilbert space $2\Hilb$), each splitting of an $\Hstar$-monad in $\Hom(\fX^\cent\to 2\Hilb)$ is already in $\Hom(\fX\to2\Hilb)$, so $(\Hom(\fX\to2\Hilb),\vee^{\Hom^{\cent}},\Psi^{\Hom^{\cent}})$ is complete.
\end{pf}

\begin{corollary}
    \label{3hilb equiv to its reps}
    If $\fX$ is a finite 3-Hilbert space, then there is an isometric equivalence $\fX\cong^\dag\Rep(\fX^{1\op})$.
\end{corollary}

\begin{pf}
By \cite[Thm. D]{bases}, the isometric Yoneda embedding $\yo\colon\fX\xrightarrow{\sim}\Hom(\fX^{1\op}\to 2\Hilb)=\Rep(\fX^{1\op})$, $a\mapsto\fX(-\to a)$, is an isometric equivalence with inverse $\!\yo^{\!-1}$
from \cite[Cor. 4.16]{bases}.
\end{pf}

\subsection{Functoriality and the functor-category equivalence}

\subsubsection{Disk-like \texorpdfstring{(2,$\ell$)}{(2,l)}-transfors and their skeletonizations}
For the definition of disk-like $(2,\ell)$-transfors for $\ell\in\{0,1,2\}$ (where \defn{$(2,0)$-transfor} is another name for a disk-like functor), see \cref{def: DL k-transfor}. See also \cref{def: functor disk-like n-category} for the definition of the disk-like 2-category $\eHom(\C{\to}\D)$ of disk-like functors, $(2,1)$-transfors, and $(2,2)$-transfors.

\subsubsection{Functors and transfors from the disk-like side}

Fix disk-like 2-categories $\C$ and $\D$.

\begin{construction}[Functors from disk-like functors, \texorpdfstring{$n=2$}{n is 2}]
\label{dag-thing-ii-1}
If $\eF\colon\C\to\D$ is a disk-like functor, then $\eF$ induces a canonical strictly UAF-preserving strict $2$-functor $\fX_\eF\colon\fX_\C\to\fX_\D$ given by $\eF$ as follows.

Define $\fX_\eF$ on objects $a\in\Obj(\fX_\C)=\fld[0]{\C}{\pt}$ by $\fX_\eF(a)\coloneq\eF(a)$, on 1-morphisms $X\in\fX_\C(a\to b)=\fld[1]{\C}{I}[\fcj{a}\amalg b]$ by $\fX_\eF(X)\coloneq\eF(X)$, and on $2$-morphisms $f\in\fX_\C(X\Rightarrow Y)=\fldl{\C}{D^2}[\fcj{X}\cup Y]$ by $\fX_\eF(f)\coloneq\eF(f)$.
\end{construction}

\begin{lemma}
\label{dag-thing-ii-1-pf}
\cref{dag-thing-ii-1} indeed defines a (strict) functor that (strictly) preserves the UAF.
\end{lemma}

\begin{pf}
\itemstep{$\fX_\eF$ is a $2$-functor.} For $1$-morphisms $X,Y\in\fX_\C$,
  \[
      \fX_\eF(X\otimes Y)
      =\eF(\rho^1_*(X\blt Y))
      =\rho^1_*(\eF(X)\blt\eF(Y))
      =\fX_\eF(X)\otimes \fX_\eF(Y)
  \]
  and $\fX_\eF(1_a)=\eF(a\times I)=\eF(a)\times I=1_{\fX_\eF(a)}$ (so $\fX_\eF$ is strict). For composable 2-morphisms $f$ and $g$, 
  \[
      \fX_\eF(g\circ f)
      =\eF(\rho^2_*(f\blt g))
      =\rho^2_*(\eF(f)\blt\eF(g))
      =\fX_\eF(g)\circ \fX_\eF(f),
  \]
  and $\fX_\eF(\id_X)=\eF(X\times I)=\eF(X)\times I=\id_{\fX_\eF(X)}$. Computing similarly,
  \[
      \fX_\eF(f\otimes g)
      =\eF(\rho^1_*(f\blt g))
      =\rho^1_*(\eF(f)\blt\eF(g))
      =\fX_\eF(f)\otimes \fX_\eF(g).
  \]
  And $\fX_\eF$ is dagger, since
  $\fX_\eF(f^\dag)
      =\eF(\fcj{\iota^{(2)}_*f})
      =\fcj{\eF(\iota^{(2)}_*f)}
      =\fcj{\iota^{(2)}_*\eF(f)}
      =\fX_\eF(f)^\dag$.

\itemstep{$\fX_\eF$ is UAF-preserving.} 
Indeed, for 1-morphisms $X$ in $\fX_\C$, we have
\[
  \fX_\eF(X^\vee)
  \defeq{}\eF(\fcj{\iota^{(1)}_*X})
  \underset{\ref{axiom:eF-refl}}{=}\fcj{\eF(\iota^{(1)}_*X)}
  \underset{\ref{axiom:eF-nat}}{=}\fcj{\iota^{(1)}_*\eF(X)}
  \defeq{}\fX_\eF(X)^\vee.\qedhere
\]
\end{pf}

\begin{lemma}[Change of parameterized splitting]
\label{lem:PS-change}
Let $\eF\colon\C\to\D$ be an unrestricted disk-like functor of disk-like 2-categories, let ${}_aX_b\in\fX_\C(a\to b)$, and let $\Theta,\Theta',\Theta''\in\eF_{\mathrm{PS}}(X)$. Then the 2-morphisms
\[
u_{\Theta,\Theta'}\coloneq\widehat\eF^{\,\orev\Theta\blt\Theta'}_{D^2,\,\fcj X\cup X}(X\times I)\in\fX_\D(\eF(X,\Theta)\Rightarrow\eF(X,\Theta'))
\]
satisfy
\begin{lst}
\item\label{cl:PS-change-id} $u_{\Theta,\Theta}=\id_{\eF(X,\Theta)}$,
\item\label{cl:PS-change-trans} $u_{\Theta',\Theta''}\circ u_{\Theta,\Theta'}=u_{\Theta,\Theta''}$, and
\item\label{cl:PS-change-dag} $u_{\Theta,\Theta'}^\dag=u_{\Theta',\Theta}$.
\end{lst}
In particular, each $u_{\Theta,\Theta'}$ is unitary with inverse $u_{\Theta,\Theta'}^\dag=u_{\Theta',\Theta}$.
\end{lemma}

\begin{pf}
\itemstep{\ref{cl:PS-change-id}.} Choose a pinched product map $\pi\colon D^2\to D^1$ and $\Lambda\in\eF_{\mathrm{PS}}(\pi^*X)$ with $\partial\Lambda=\orev\Theta\blt\Theta$ (which exists by \ref{axiom:Gamma-products}). Then \ref{axiom:FTheta-products} gives $\eF(\pi^*X,\Lambda)=\pi^*\eF(X,\Theta)$, so $u_{\Theta,\Theta}=[\pi^*\eF(X,\Theta)]=\id_{\eF(X,\Theta)}$.

\itemstep{\ref{cl:PS-change-trans}.} Let $\sigma'$ and $\sigma''$ be pure representatives of $\id_X$. Then for any parameterized splittings $\Lambda'\in\eF_{\mathrm{PS}}(\sigma')$ and $\Lambda''\in\eF_{\mathrm{PS}}(\sigma'')$ with $\partial\Lambda'=\orev\Theta\blt\Theta'$ and $\partial\Lambda''=\orev{\Theta'}\blt\Theta''$, we have
\begin{align*}
u_{\Theta',\Theta''}\circ u_{\Theta,\Theta'}
&
\defeq
\rho^2_*(\eF(\sigma',\Lambda')\blt_{\eF(X,\Theta')}\eF(\sigma'',\Lambda''))
\\&
\underset{\ref{axiom:FTheta-nat}\,\&\,\ref{axiom:FTheta-gluing}}{=}
\eF(\rho^2_*(\sigma'\blt_X\sigma''),\rho^2_*(\Lambda'\blt_{\Theta'}\Lambda''))
=u_{\Theta,\Theta''}
\end{align*}
where the last equality holds by \cref{dldldl}, as $\rho^2_*(\sigma'\blt_X\sigma'')$ is a pure representative of $\id_X\circ\id_X=\id_X$, and $\partial(\rho^2_*(\Lambda'\blt_{\Theta'}\Lambda''))=\orev\Theta\blt\Theta''$.

\itemstep{\ref{cl:PS-change-dag}.} Indeed, $u_{\Theta,\Theta'}^\dag=[\fcj{\iota^{(2)}_*\eF(\sigma',\Lambda')}]\underset{\ref{axiom:FTheta-nat}\,\&\,\ref{axiom:FTheta-refl}}{=}[\eF(\fcj{\iota^{(2)}_*\sigma'},\fcj{\iota^{(2)}_*\Lambda'})]=u_{\Theta',\Theta}$ again by \cref{dldldl}: $\fcj{\iota^{(2)}_*\sigma'}$ is a pure representative of $X\times I=(X\times I)^{\dag_2}$, and $\partial(\fcj{\iota^{(2)}_*\Lambda'})=\orev{\Theta'}\blt\Theta$. 

Combining the above items, we obtain $u_{\Theta,\Theta'}^\dag\circ u_{\Theta,\Theta'}=u_{\Theta,\Theta}=\id$ and $u_{\Theta,\Theta'}\circ u_{\Theta,\Theta'}^\dag=u_{\Theta',\Theta'}=\id$, so $u_{\Theta,\Theta'}$ is unitary. 
\end{pf}

\begin{construction}[Functors from unrestricted disk-like functors, \texorpdfstring{$n=2$}{n is 2}]
\label{dag-thing-ii-1-ur}
If $\eF\colon\C\to\D$ is an unrestricted disk-like functor, then a choice of parameterized splitting $\Theta_X\in\eF_{\mathrm{PS}}(X)$ for each 1-morphism ${}_aX_b\in\fX_\C(a\to b)$ determines a canonical UAF-preserving dagger functor $\fX_\eF\colon\fX_\C\to\fX_\D$ as follows. 

\begin{itemize}
    \item On objects $a\in\fX_\C$, $\fX_\eF(a)\coloneq\eF(a)$.
    \item On 1-morphisms ${}_aX_b\in\fX_\C(a\to b)$, $\fX_\eF(X)\coloneq\eF(X,\Theta_X)$.
    \item On 2-morphisms $f\in\fX_\C({}_aX_b\Rightarrow{}_aY_b)$, $\fX_\eF(f)\coloneq\widehat\eF_{D^2,\fcj{X}\cup Y}(f)$.
    \item For 0-fields $a\in\fX_\C$, the unitors are the identity 2-morphisms $\fX_\eF(1_a)=1_{\eF(a)}$.
    \item For 1-morphisms ${}_aX_b$ and ${}_bY_c$ in $\fX_\C$, define the tensorator $(\fX_\eF)^2_{X,Y}\in\fX_\D(\fX_\eF(X)\otimes \fX_\eF(Y)\Rightarrow\fX_\eF(X\otimes Y))$ by the unitary
    \[
    (\fX_\eF)^2_{X,Y}\coloneq u_{\rho^1_*(\Theta_X\blt_b\Theta_Y),\Theta_{X\otimes Y}}
    \]
    defined in \cref{lem:PS-change}.
\end{itemize}
\end{construction} 

\begin{lemma}
\label{dag-thing-ii-1-ur-pf}
\cref{dag-thing-ii-1-ur} indeed defines a UAF-preserving $\dag$-functor.
\end{lemma}

\begin{pf}
That $\fX_\eF$ preserves vertical composition, identities, daggers, and horizontal composition of 2-morphisms follows from the same computations as the proof of \cref{dag-thing-ii-1-pf} but with $\eF(-)$ replaced by $\eF(-,\Theta)$ and the axioms \ref{axiom:eF-nat}, \ref{axiom:eF-gluing}, \ref{axiom:eF-products}, \ref{axiom:eF-refl} replaced by \ref{axiom:FTheta-nat}, \ref{axiom:FTheta-gluing}, \ref{axiom:FTheta-products}, \ref{axiom:FTheta-refl} respectively, while the tensorators $(\fX_\eF)^2_{X,Y}$ are unitary by \cref{lem:PS-change}. We leave the coherence equations to the reader.
To see $\fX_\eF$ is UAF-preserving, observe that the canonical isomorphism $\delta_X\colon\fX_\eF(X^\vee)\Rightarrow \fX_\eF(X)^\vee$ from \cref{UAF-preserving} equals the unitary $u_{\Theta_{X^\vee},\Theta_X^\vee}$ from \cref{lem:PS-change}, where $\Theta_X^\vee\coloneq\fcj{\iota^{(1)}_*\Theta_X}$. 
\end{pf}

The following lemma is the 2-categorical generalization (or horizontal categorification) of the standard fact that a monoidal natural transformation between monoidal functors out of a rigid monoidal category is automatically invertible (see \cite[\S2.10]{EGNO15}): when $\fX$ and $\fY$ are deloopings of monoidal categories, $F$ and $G$ are monoidal functors and $\eta$ is a monoidal natural transformation. 

\begin{lemma}
\label{autounitary}
For rigid 2-categories $\fX$ and $\fY$, suppose $\eta \colon F \Rightarrow G$ is known to satisfy all the conditions of a natural transformation between functors $F, G \colon \fX \to \fY$ except for the invertibility of the naturators. Then the naturators are automatically invertible. 
\end{lemma}

\begin{pf}
We utilize the overlay graphical calculus \cite{CP22} with
\[
a=\begin{tkz}[scale=1.25]
    \fill[rounded corners,lightgray!25](-.3,-.3) rectangle (.3,.3);
\end{tkz}\,, 
\qquad 
b=\begin{tkz}[scale=1.25]
    \fill[rounded corners,lightgray!75](-.3,-.3) rectangle (.3,.3);
\end{tkz}\,,
\qquad
X=\begin{tkz}[scale=1.25]
    \clip[rounded corners](-.3,-.3) rectangle (.3,.3);
    \fill[lightgray!25](-.3,-.3) rectangle (.3,.3);
    \fill[lightgray!75](0,-.3) rectangle (.3,.3);
    \draw(0,-0.3)--(0,0.3);
\end{tkz}\,,
\qquad
F=\begin{tkz}[scale=1.25] 
    \fill[rounded corners,primedregion=white,draw=black,dotted](-.3,-.3) rectangle (.3,.3);
\end{tkz}\,, 
\qquad
G=\begin{tkz}[scale=1.25] 
    \fill[rounded corners,boxregion=white,draw=black,dotted](-.3,-.3) rectangle (.3,.3);
\end{tkz}\,,
\qquad\text{and}\qquad
\eta=\begin{tkz}[scale=1.25]
    \clip[rounded corners](-.3,-.3) rectangle (.3,.3);
    \fill[primedregion=white](-.3,-.3) rectangle (.3,.3);
    \fill[boxregion=white](0,-.3) rectangle (.3,.3);
    \draw[red,very thick](0,-0.3)--(0,0.3);
    \draw[rounded corners,draw=black,dotted](-.3,-.3) rectangle (.3,.3);
\end{tkz}
\]
so that
\[
  \eta_X \;=
  \begin{tkz}[scale=0.55,xscale=1]
    \begin{scope}
      \clip[rounded corners] (0,0) rectangle (3,2);
      \fill[lightgray!25] (0,0) rectangle (3,2);
      \fill[lightgray!75] (0.8,0) to[out=90,in=270] (2.2,2) -- (3,2) -- (3,0) -- cycle;
      \fill[pattern=primeddots] (2.2,0) to[out=90,in=270] (0.8,2) -- (0,2) -- (0,0) -- cycle;
      \fill[pattern=primedbox] (2.2,0) to[out=90,in=270] (0.8,2) -- (3,2) -- (3,0) -- cycle;
    \end{scope}
    \draw (0.8,0) to[out=90,in=270] (2.2,2);
    \draw[red, thick] (2.2,0) to[out=90,in=270] (0.8,2);
    \node[below] at (0.8,0) {\scriptsize $F(X)$};
    \node[above] at (2.2,2) {\scriptsize $G(X)$};
    \node[above] at (0.8,2) {\scriptsize $\eta_a$};
    \node[below] at (2.2,0) {\scriptsize $\eta_b$};
  \end{tkz}
\quad\qquad\text{and}\quad\qquad
  \eta_{X^\vee} \;=
  \begin{tkz}[scale=0.55,xscale=1.05]
    \begin{scope}
      \clip[rounded corners] (0,0) rectangle (3,2);
      \fill[lightgray!75] (0,0) rectangle (3,2);
      \fill[lightgray!25] (0.8,0) to[out=90,in=270] (2.2,2) -- (3,2) -- (3,0) -- cycle;
      \fill[pattern=primeddots] (2.2,0) to[out=90,in=270] (0.8,2) -- (0,2) -- (0,0) -- cycle;
      \fill[pattern=primedbox] (2.2,0) to[out=90,in=270] (0.8,2) -- (3,2) -- (3,0) -- cycle;
    \end{scope}
    \draw (0.8,0) to[out=90,in=270] (2.2,2);
    \draw[red, thick] (2.2,0) to[out=90,in=270] (0.8,2);
    \node[below,xshift=-2mm] at (0.8,0) {\scriptsize $F(X^\vee)$};
    \node[above,xshift=2mm] at (2.2,2) {\scriptsize $G(X^\vee)$};
    \node[above,xshift=-2mm] at (0.8,2) {\scriptsize $\eta_b$};
    \node[below,xshift=2mm] at (2.2,0) {\scriptsize $\eta_a$};
  \end{tkz}
  \!\!\!\!.
\]
We claim $\eta_X$ has inverse given by
\[
\beta\coloneq
\begin{tkz}[scale=0.4,xscale=1.4]
    \begin{scope}
      \clip[rounded corners] (0,-4) rectangle (4.8,2.4);
      \fill[lightgray!25] (0,-4) rectangle (4.8,2.4);
      \fill[lightgray!75] (4.2,-4) -- (4.2,1.0) arc(0:180:0.6) -- (3.0,-0.4) to[out=270,in=90] (1.8,-2.0) -- (1.8,-2.6) arc(0:-180:0.6) -- (0.6,2.4) -- (4.8,2.4) -- (4.8,-4) -- cycle;
      \fill[pattern=primeddots] (3.0,-4) -- (3.0,-2.0) to[out=90,in=270] (1.8,-0.4) -- (1.8,2.4) -- (0,2.4) -- (0,-4) -- cycle;
      \fill[pattern=primedbox] (3.0,-4) -- (3.0,-2.0) to[out=90,in=270] (1.8,-0.4) -- (1.8,2.4) -- (4.8,2.4) -- (4.8,-4) -- cycle;
    \end{scope}
    \draw (4.2,-4) -- (4.2,1.0) arc(0:180:0.6) -- (3.0,-0.4) to[out=270,in=90] (1.8,-2.0) -- (1.8,-2.6) arc(0:-180:0.6) -- (0.6,2.4);
    \draw[red, thick] (3.0,-4) -- (3.0,-2.0) to[out=90,in=270] (1.8,-0.4) -- (1.8,2.4);
    \node[below] at (3.0,-4) {\scriptsize $\eta_a$};
    \node[below] at (4.2,-4) {\scriptsize $G(X)$};
    \node[above] at (0.6,2.4) {\scriptsize $F(X)$};
    \node[above] at (1.8,2.4) {\scriptsize $\eta_b$};
  \end{tkz}
  .
\]
Indeed,
\[
  \eta_X \circ\beta
  =
  \begin{tkz}[scale=0.3,xscale=1.8]
    \begin{scope}
      \clip[rounded corners] (0,-4) rectangle (4.8,5);
      \fill[lightgray!25] (0,-4) rectangle (4.8,5);
      \fill[lightgray!75] (4.2,-4) -- (4.2,1.0) arc(0:180:0.6) -- (3.0,-0.4) to[out=270,in=90] (1.8,-2.0) -- (1.8,-2.6) arc(0:-180:0.6) -- (0.6,2.2) to[out=90,in=270] (1.8,3.8) -- (1.8,5) -- (4.8,5) -- (4.8,-4) -- cycle;
      \fill[pattern=primeddots] (3.0,-4) -- (3.0,-2.0) to[out=90,in=270] (1.8,-0.4) -- (1.8,2.2) to[out=90,in=270] (0.6,3.8) -- (0.6,5) -- (0,5) -- (0,-4) -- cycle;
      \fill[pattern=primedbox] (3.0,-4) -- (3.0,-2.0) to[out=90,in=270] (1.8,-0.4) -- (1.8,2.2) to[out=90,in=270] (0.6,3.8) -- (0.6,5) -- (4.8,5) -- (4.8,-4) -- cycle;
    \end{scope}
    \draw (4.2,-4) -- (4.2,1.0) arc(0:180:0.6) -- (3.0,-0.4) to[out=270,in=90] (1.8,-2.0) -- (1.8,-2.6) arc(0:-180:0.6) -- (0.6,2.2) to[out=90,in=270] (1.8,3.8) -- (1.8,5);
    \draw[red, thick] (3.0,-4) -- (3.0,-2.0) to[out=90,in=270] (1.8,-0.4) -- (1.8,2.2) to[out=90,in=270] (0.6,3.8) -- (0.6,5);
    \node[below] at (3.0,-4) {\scriptsize $\eta_a$};
    \node[below] at (4.2,-4) {\scriptsize $G(X)$};
    \node[above] at (0.6,5) {\scriptsize $\eta_a$};
    \node[above] at (1.8,5) {\scriptsize $G(X)$};
  \end{tkz}
  \underset{\text{(naturality)}}{=}\;
  \begin{tkz}[scale=0.5,xscale=1.5]
    \begin{scope}
      \clip[rounded corners] (0,0) rectangle (4.8,4);
      \fill[lightgray!25] (0,0) rectangle (4.8,4);
      \fill[lightgray!75] (4.2,0) -- (4.2,2.6) arc(0:180:0.6) -- (3.0,1.4) arc(360:180:0.6) -- (1.8,4) -- (4.8,4) -- (4.8,0) -- cycle;
      \fill[pattern=primeddots] (0.6,0) -- (0.6,4) -- (0,4) -- (0,0) -- cycle;
      \fill[pattern=primedbox] (0.6,0) -- (0.6,4) -- (4.8,4) -- (4.8,0) -- cycle;
    \end{scope}
    \draw (4.2,0) -- (4.2,2.6) arc(0:180:0.6) -- (3.0,1.4) arc(360:180:0.6) -- (1.8,4);
    \draw[red, thick] (0.6,0) -- (0.6,4);
    \node[below] at (0.6,0) {\scriptsize $\eta_a$};
    \node[below] at (4.2,0) {\scriptsize $G(X)$};
    \node[above] at (0.6,4) {\scriptsize $\eta_a$};
    \node[above] at (1.8,4) {\scriptsize $G(X)$};
  \end{tkz}
  \underset{\text{(zig-zag)}}{=}\;
  \begin{tkz}[scale=0.5,xscale=1.5]
    \begin{scope}
      \clip[rounded corners] (0,0) rectangle (2.4,4);
      \fill[lightgray!25] (0,0) rectangle (2.4,4);
      \fill[lightgray!75] (1.8,0) -- (1.8,4) -- (2.4,4) -- (2.4,0) -- cycle;
      \fill[pattern=primeddots] (0.6,0) -- (0.6,4) -- (0,4) -- (0,0) -- cycle;
      \fill[pattern=primedbox] (0.6,0) -- (0.6,4) -- (2.4,4) -- (2.4,0) -- cycle;
    \end{scope}
    \draw (1.8,0) -- (1.8,4);
    \draw[red, thick] (0.6,0) -- (0.6,4);
    \node[below] at (0.6,0) {\scriptsize$\eta_a$};
    \node[below,xshift=2mm] at (1.8,0) {\scriptsize $G(X)$};
    \node[above] at (0.6,4) {\scriptsize $\eta_a$};
    \node[above,xshift=2mm] at (1.8,4) {\scriptsize$G(X)$};
  \end{tkz}
  \;=\;
  \id_{\eta_a\otimes G(X)}
\]
and
\[
  \beta\circ\eta_X
  =
  \begin{tkz}[scale=0.3,xscale=1.8]
    \begin{scope}
      \clip[rounded corners] (0,-6) rectangle (4.8,4);
      \fill[lightgray!25] (0,-6) rectangle (4.8,4);
      \fill[lightgray!75] (3.0,-6) -- (3.0,-5.2) to[out=90,in=270] (4.2,-3.6) -- (4.2,2.0) arc(0:180:0.6) -- (3.0,0.2) to[out=270,in=90] (1.8,-1.4) -- (1.8,-2.0) arc(0:-180:0.6) -- (0.6,4) -- (4.8,4) -- (4.8,-6) -- cycle;
      \fill[pattern=primeddots] (4.2,-6) -- (4.2,-5.2) to[out=90,in=270] (3.0,-3.6) -- (3.0,-1.4) to[out=90,in=270] (1.8,0.2) -- (1.8,4) -- (0,4) -- (0,-6) -- cycle;
      \fill[pattern=primedbox] (4.2,-6) -- (4.2,-5.2) to[out=90,in=270] (3.0,-3.6) -- (3.0,-1.4) to[out=90,in=270] (1.8,0.2) -- (1.8,4) -- (4.8,4) -- (4.8,-6) -- cycle;
    \end{scope}
    \draw (3.0,-6) -- (3.0,-5.2) to[out=90,in=270] (4.2,-3.6) -- (4.2,2.0) arc(0:180:0.6) -- (3.0,0.2) to[out=270,in=90] (1.8,-1.4) -- (1.8,-2.0) arc(0:-180:0.6) -- (0.6,4);
    \draw[red, thick] (4.2,-6) -- (4.2,-5.2) to[out=90,in=270] (3.0,-3.6) -- (3.0,-1.4) to[out=90,in=270] (1.8,0.2) -- (1.8,4);
    \node[below] at (3.0,-6) {\scriptsize$F(X)$};
    \node[below] at (4.2,-6) {\scriptsize $\eta_b$};
    \node[above] at (0.6,4) {\scriptsize$F(X)$};
    \node[above] at (1.8,4) {\scriptsize$\eta_b$};
  \end{tkz}
  \underset{\text{(naturality)}}{=}
  \begin{tkz}[scale=0.5,xscale=1.5]
    \begin{scope}
      \clip[rounded corners] (0,0) rectangle (4.8,4);
      \fill[lightgray!25] (0,0) rectangle (4.8,4);
      \fill[lightgray!75] (3.0,0) -- (3.0,2.6) arc(0:180:0.6) -- (1.8,1.4) arc(360:180:0.6) -- (0.6,4) -- (4.8,4) -- (4.8,0) -- cycle;
      \fill[pattern=primeddots] (4.2,0) -- (4.2,4) -- (0,4) -- (0,0) -- cycle;
      \fill[pattern=primedbox] (4.2,0) -- (4.2,4) -- (4.8,4) -- (4.8,0) -- cycle;
    \end{scope}
    \draw (3.0,0) -- (3.0,2.6) arc(0:180:0.6) -- (1.8,1.4) arc(360:180:0.6) -- (0.6,4);
    \draw[red, thick] (4.2,0) -- (4.2,4);
    \node[below] at (3.0,0) {\scriptsize $F(X)$};
    \node[below] at (4.2,0) {\scriptsize $\eta_b$};
    \node[above] at (0.6,4) {\scriptsize $F(X)$};
    \node[above] at (4.2,4) {\scriptsize $\eta_b$};
  \end{tkz}
\underset{\text{(zig-zag)}}{=}
  \begin{tkz}[scale=0.6,xscale=1.5]
    \begin{scope}
      \clip[rounded corners] (0,0) rectangle (2.4,4);
      \fill[lightgray!25] (0,0) rectangle (2.4,4);
      \fill[lightgray!75] (0.6,0) -- (0.6,4) -- (2.4,4) -- (2.4,0) -- cycle;
      \fill[pattern=primeddots] (1.8,0) -- (1.8,4) -- (0,4) -- (0,0) -- cycle;
      \fill[pattern=primedbox] (1.8,0) -- (1.8,4) -- (2.4,4) -- (2.4,0) -- cycle;
    \end{scope}
    \draw (0.6,0) -- (0.6,4);
    \draw[red, thick] (1.8,0) -- (1.8,4);
    \node[below] at (0.6,0) {\scriptsize $F(X)$};
    \node[below] at (1.8,0) {\scriptsize $\eta_b$};
    \node[above] at (0.6,4) {\scriptsize $F(X)$};
    \node[above] at (1.8,4) {\scriptsize $\eta_b$};
  \end{tkz}
  =
  \id_{F(X)\otimes \eta_b}.\qedhere
\]
\end{pf}

\begin{construction}[Natural transformations from (2,1)-transfors]
    \label{dag-thing-ii-2}
  A disk-like (2,1)-transfor $\eN\in\fld{\eFun(\C{\to}\D)}{I}[\fcj{\eF}\amalg\eG]$ for disk-like functors $\eF,\eG\colon\C\to\D$ induces a canonical natural transformation $\fX_\eN\colon\fX_\eF\Rightarrow\fX_\eG$ as follows.
  For each $a \in \Obj(\fX_\C) = \fld[0]{\C}{\pt}$, set
\[
    (\fX_\eN)_a \coloneq
    \begin{tkz}[scale=0.85]
    \draw[\colX,thick,mid<] (0,1) -- node[left]{\scriptsize$\eN(a)$} (0,0);
    \end{tkz}
    \;\;\in\;\; 
    \fld[1]{\D}{I}[ \fcj{\eF(a)} \amalg \eG(a)] = \fX_\D(\fX_\eF(a) \to \fX_\eG(a)).
\]
For each $X \in \fX_\C(a \to b) = \fld[1]{\C}{I}[ \fcj{a} \amalg b]$, define $\fX_\eN$ to have 2-cell components $(\fX_\eN)_X$ given by
\[
(\fX_\eN)_X
\coloneq
\begin{tkz}[scale=0.85]
    \fill[\colllGr] (0,0) rectangle (2,2);
    \draw[\colY,thick,mid>] (0,2) -- node[above=-1mm]{\scriptsize$\eF(X)$} (2,2);
    \draw[blue,thick,mid>] (2,2) -- node[right]{\scriptsize$\eN(b)$} (2,0);
    \draw[\colX,thick,mid>] (0,2) -- node[left]{\scriptsize$\eN(a)$} (0,0);
    \draw[red,thick,mid>] (0,0) -- node[below=-1mm]{\scriptsize$\eG(X)$} (2,0);
    \node at (1,1){\scriptsize$\eN(X)$};
\end{tkz}
\!\overset{\text{reparameterize}}{\rightsquigarrow}\!
\begin{tkz}[scale=1]
    \fill[\colllGr] (0,0) .. controls (0.5,\hgt) and (1.5,\hgt) .. (2,0) .. controls (1.5,-\hgt) and (0.5,-\hgt) .. (0,0);    
    \draw[\colY,thick,mid>] (0,0) .. controls (0.25,0.5*\hgt) and (0.625,0.75*\hgt) .. node[above left=-1mm]{\scriptsize$\eF(X)$}(1,0.75*\hgt);
    \draw[blue,thick,mid>] (1,0.75*\hgt) .. controls (1.375,0.75*\hgt) and (1.75,0.5*\hgt) .. node[above right=-1mm]{\scriptsize$\eN(b)$} (2,0);
    
    \draw[\colX,thick,mid>] (0,0) .. controls (0.25,-0.5*\hgt) and (0.625,-0.75*\hgt) .. node[below left=-1mm]{\scriptsize$\eN(a)$}(1,-0.75*\hgt);
    \draw[red,thick,mid>] (1,-0.75*\hgt) .. controls (1.375,-0.75*\hgt) and (1.75,-0.5*\hgt) .. node[below right=-1mm]{\scriptsize$\eG(X)$} (2,0);
    \node at (current bounding box.center){\scriptsize$\eN(X)$};
\end{tkz}
\in\fX_\D(\eN(a)\otimes\eG(X)\Rightarrow\eF(X)\otimes\eN(b))
\]
where we reparameterize to $D^2$ by some pre-fixed homeomorphism $\rho\colon I\times I\to D^2$ as above.
\end{construction}

\begin{lemma}
\label{dag-thing-ii-2-pf}
    \cref{dag-thing-ii-2} indeed defines a natural transformation.
\end{lemma}

\begin{pf}
The naturators $(\fX_\eN)_X$ are automatically invertible by \cref{autounitary}. To see $(\fX_\eN)_X$ is unitary, observe that $(\fX_\eN)_X^{-1}$ as defined in \cref{autounitary} and $(\fX_\eN)_X^\dag$ are related by a boundary-fixing orientation-preserving homeomorphism rotating the bulk by a one-click rotation and thus are representatives of the same 2-morphism in $\fX_\D(\eta_a\otimes G(X)\Rightarrow F(X)\otimes\eta_b)$ by \ref{CTi}. Next, the unit condition $(\fX_\eN)_{1_a}=\id_{\eN(a)}$ is automatic by \ref{axiom:eT-products}. Since $\fX_\eF$ and $\fX_\eG$ are strict, the 1-composition condition simplifies to $(\fX_\eN)_{X\otimes Y}=(\fX_\eN)_X\otimes(\fX_\eN)_Y$ for all compatible 1-morphisms $X,Y$ in $\fX_\C$, which is a straightforward 
check: 
\[
(\fX_\eN)_{X\otimes Y} = \eN(\rho^1_*(X \bullet_b Y)) = \rho^1_*(\eN(X \bullet_b Y)) = \rho^1_*(\eN(X) \bullet_{\eN(b)} \eN(Y)) = (\fX_\eN)_X \otimes (\fX_\eN)_Y.
\]

It remains to prove the naturality condition. Consider a pure 2-field $f\in\fld[2]{\C}{D^2}[\fcj{X}\cup Y]$. Gluing $\eF(f)$ on $D^2\times\{0\}$, $\eG(f)$ on $D^2\times\{1\}$, $(\fX_\eN)_X$ on $\partial_-D^2\times I$, and $(\fX_\eN)_Y$ on $\partial_+D^2\times I$ produces a pure 2-field on the 2-sphere $\partial(D^2\times I)$, to which we apply the naturality axiom \ref{axiom:eT-naturality} of the $(2,1)$-transfor $\eN$: the splittings $\partial D^2=\partial_-D^2\cup\partial_+D^2$ and $\partial I=\{0\}\amalg\{1\}$ cut this 2-sphere into the two 2-balls $Z_1\coloneq(D^2\times\{0\})\cup(\partial_-D^2\times I)$ and $Z_2\coloneq(\partial_+D^2\times I)\cup(D^2\times\{1\})$, on which the pure 2-fields $\xi$ and $\zeta$ live respectively, and \ref{axiom:eT-naturality} gives a homeomorphism $\varphi$ with $\varphi_*[\xi]=[\zeta]$. Thus by reparameterizing to $D^2$ we get the following equality of (quotient-level) 2-fields in $\fldl[2]{\D}{D^2}$, i.e., of 2-morphisms in $\fX_\D$.
\[
  (\eF(f)\otimes\id_{\eN(b)})\circ(\fX_\eN)_X
  =
  \begin{tkz}[scale=0.35, baseline=(current bounding box.center)]
    \fill[lightgray!30] (0,3) to[bend right=30] (4,3) -- (4,0) to[bend left=30] (0,0) -- cycle;
    \fill[orange!25] (0,3) to[bend left=30] (4,3) to[bend left=30] (0,3);
    \begin{scope}[mid>]
      \draw (0,3) to[bend left=30] node[above] {\tiny $\eF(Y)$} (4,3);
      \draw[thick,dotted] (0,3) to[bend right=30] 
      (4,3);
      \draw (0,3) -- node[left] {\tiny $\eN(a)$} (0,0);
      \draw (4,3) -- node[right] {\tiny $\eN(b)$} (4,0);
      \draw (0,0) to[bend right=30] node[below] {\tiny $\eG(X)$} (4,0);
    \end{scope}
    \node at (2,1) {\scriptsize $(\fX_\eN)_X$};
    \node at (2,3) {\scriptsize $\eF(f)$};
  \end{tkz}
  =
  \begin{tkz}[scale=0.35, baseline=(current bounding box.center)]
    \fill[lightgray!30] (0,3) to[bend left=30] (4,3) -- (4,0) to[bend right=30] (0,0) -- cycle;
    \fill[orange!25] (0,0) to[bend left=30] (4,0) to[bend left=30] (0,0);
    \begin{scope}[mid>]
      \draw (0,3) to[bend left=30] node[above] {\tiny $\eF(Y)$} (4,3);
      \draw[dotted,thick] (0,0) to[bend left=30] 
      (4,0);
      \draw (0,0) to[bend right=30] node[below] {\tiny $\eG(X)$} (4,0);
      \draw (0,3) -- node[left] {\tiny $\eN(a)$} (0,0);
      \draw (4,3) -- node[right] {\tiny $\eN(b)$} (4,0);
    \end{scope}
    \node at (2,2) {\scriptsize $(\fX_\eN)_Y$};
    \node at (2,0) {\scriptsize $\eG(f)$};
  \end{tkz}
  =
  (\fX_\eN)_Y\circ(\id_{\eN(a)}\otimes\eG(f)).\qedhere
\]
\end{pf}

\begin{construction}[Modifications from (2,2)-transfors]
\label{dag-thing-ii-3}
If $\eM\in\fldl{\eFun(\C{\to}\D)}{D^2}[\fcj{\eS}\cup\eT]$ is a $D^2$-shaped (2,2)-transfor with choice of boundary splitting $\fcj{\eS}\cup\eT$ for (2,1)-transfors $\eT,\eS\in\fld{\eFun(\C{\to}\D)}{D^1}[\fcj{\eF}\amalg \eG]$ of disk-like functors $\eF,\eG\colon\C\to\D$, then $\eM$ determines a canonical modification 
$\fX_\eM\colon\fX_\eS\Rrightarrow\fX_\eT$ whose component at an object 
$a\in\fX_\C$ is given by
\[
    (\fX_\eM)_a\coloneq \eM(a)=
    \begin{tkz}[scale=1.35]
    \def\hgt{0.6}
        \fill[\colllGr](0,0) .. controls (0.5,\hgt) and (1.5,\hgt) .. (2,0).. controls (1.5,-\hgt) and (0.5,-\hgt) .. (0,0);
        \draw[\colY,thick,mid>] (0,0) .. controls (0.5,\hgt) and (1.5,\hgt) .. node[above]{\scriptsize$\eT(a)$} (2,0); % node[midway, above] {\tiny $D^1_\out$};
        \draw[\colX,thick,mid<] (0,0) .. controls (0.5,-\hgt) and (1.5,-\hgt) ..  node[below]{\scriptsize${\eS(a)}$} (2,0); % node[midway, below] {\tiny $D^1_\inn$};
        \filldraw[\colG,thick] (0,0) circle (1pt) node[left]{\scriptsize${\eF(a)}$};
        \filldraw[\colG] (2,0) circle (1pt) node[right]{\scriptsize$\eG(a)$};
        \node at (current bounding box.center){\scriptsize$\eM(a)$};
    \end{tkz}
    \in\fX_\D(\eS(a)\Rightarrow\eT(a)).
\]
\end{construction}

As the proof of the following lemma is the same argument employed in the proofs of \cref{dag-thing-i-2-pf} and \cref{dag-thing-ii-2-pf}, we leave it to the reader.
\begin{lemma}
\label{dag-thing-ii-3-pf}
 \cref{dag-thing-ii-3} indeed defines a modification.
\end{lemma}

\begin{lemma}
\label{CcongD implies XCcongXD-n2}
If $\eF\colon\C\to\D$ is a (possibly unrestricted) isometric weak equivalence of unitary disk-like 2-categories, then $\fX_\eF$ is an isometric equivalence of proto-3-Hilbert spaces. 
\end{lemma}

\begin{pf}
By \cref{dag-thing-ii-1-ur-pf}, $\fX_\eF$ is a UAF-preserving dagger functor.

\itemstep{Isometrically essentially surjective on objects.} By \ref{WEurk} at $k=0$, for each $b\in\fX_\D$ there is $a\in\fX_\C$ with $\eF(a)\cong^\star b$ in $\D$, which by \cref{isomequiviff} is an isometric equivalence $\fX_\eF(a)\cong^\dag b$ in $\fX_\D$.

\itemstep{Essentially surjective on 1-morphisms.} Let $\eta\in\fX_\D(\eF(a)\to\eF(b))$. Since $\eF_{\mathrm{PS}}(\fcj a\amalg b)=\{\varnothing\}$, by \ref{WEurk} at $k=1$ there is a 1-morphism ${}_aX_b\in\fX_\C(a\to b)$ and a $\Theta\in\eF_{\mathrm{PS}}(X)$ with $\eF(X,\Theta)\cong^\star\eta$ in $\D$, which is just a unitary $2$-isomorphism $\eF(X,\Theta)\Rightarrow\eta$ in $\fX_\D$. By composing with $u_{\Theta_X,\Theta}$ from \cref{lem:PS-change}, we get a unitary 2-morphism $\fX_\eF(X)\Rightarrow\eta$.

\itemstep{Locally fully faithful.} Applying \ref{WEurn} to the boundary splitting $\Xi=\orev{\Theta_X}\blt\Theta_Y$ gives that $\fX_\eF\colon\fX_\C(X\Rightarrow Y)\to\fX_\D(\fX_\eF(X)\Rightarrow\fX_\eF(Y))$ is a linear isomorphism.

\itemstep{Isometric.} Let $a\in\fX_\C$, $f\in\End_{\fX_\C}(1_a)$, and choose a pure 2-field $\sigma_f$ with $[\sigma_f]=f$ and $\Theta\in\eF_{\mathrm{PS}}(\sigma_f)$. As \ref{axiom:Gamma-products} gives that $\eF_{\mathrm{PS}}(\fcj{a\times D^2})=\{\varnothing\}$, \ref{axiom:FTheta-gluing} and \ref{axiom:FTheta-products} imply $\eF(\cl_{a}\sigma_f,\varnothing\blt\Theta)=\cl_{\eF(a)}\eF(\sigma_f,\Theta)$. Hence
\[
\Psi^{\fX_\D}_{\eF(a)}(\fX_\eF(f))
\defeq{}\psn^\D([\cl_{\eF(a)}\eF(\sigma_f,\Theta)])
=\psn^\D([\eF(\cl_{a}\sigma_f,\varnothing\blt\Theta)])
=\psn^\C([\cl_{a}\sigma_f])
\defeq{}\Psi^{\fX_\C}_a(f),
\]
where the third equality uses that $\eF$ is isometric.
\end{pf}

The following lemma follows exactly as it did in the $n=1$ case (see \cref{lem:C-preserves-equivs}); the only new ingredient is showing isometric essential surjectivity and unitary essential surjectivity respectively on 0- and 1-fields, but this follows immediately from \cref{isomequiviff}.

\begin{lemma}
\label{2lem:C-preserves-equivs}
If $F\colon\fX\to\fY$ is an isometric equivalence of proto-3-Hilbert spaces, then $\C^F$ is an isometric weak equivalence of unitary disk-like 2-categories.
\end{lemma}

\subsubsection{Disk-like functors and transfors from the traditional side}

\

\begin{construction}[Disk-like functors from functors, \texorpdfstring{$n=2$}{n is 2}]
\label{DL 2-functor from 2-functor-1}
A UAF-preserving functor $F\colon \fX\to\fY$ induces a canonical disk-like functor $\C^F\colon\C^\fX\to\C^\fY$ sending $\fX$-string diagrams to the same string diagram but with $F$ applied to the labels. For example, 
\[
\zeta
=
\begin{tkz}[baseline=-.5ex,scale=0.275,font=\scriptsize]
    \def\tkzW{10.5}
    \def\sh{0.75}
    \def\XgI{3.5}
    \def\XgII{7.0}
    \pgfmathsetmacro{\midXOne}{\XgI/2}
    \pgfmathsetmacro{\midXTwo}{(\XgI+\XgII)/2}
    \pgfmathsetmacro{\midXThree}{(\XgII+\tkzW)/2}
    \begin{scope}
    \clip[rounded corners] (0,-\sh) rectangle (\tkzW,\sh);
    \fill[lightgray!25] (0,-\sh) rectangle (\XgI,\sh);
    \fill[lightgray!50] (\XgI,-\sh) rectangle (\XgII,\sh);
    \fill[lightgray!75] (\XgII,-\sh) rectangle (\tkzW,\sh);
    \end{scope}
    \draw[thick] (\XgI,-\sh) -- (\XgI,\sh);
    \draw[thick] (\XgII,-\sh) -- (\XgII,\sh);
        \pic at (\XgI,0){tro};
    \pic[rotate=180] at (\XgII,0){tro};
    \node[above] at (\XgI,\sh) {$X$};
    \node[above] at (\XgII,\sh) {$W$};
    \node at (\midXOne,0) {\tiny$a$};
    \node at (\midXTwo,0) {\tiny$b$};
    \node at (\midXThree,0) {\tiny$c$};
\end{tkz}
\in\fld[1]{\C^\fX}{Y}[\fcj{a}\amalg c]
\;\;\rightsquigarrow\;\;
\C^F(\zeta)
=
\begin{tkz}[scale=0.275,font=\scriptsize,baseline=-.5ex]
    \def\tkzW{10.5}
    \def\sh{0.75}
    \def\XgI{3.5}
    \def\XgII{7.0}
    \pgfmathsetmacro{\midXOne}{\XgI/2}
    \pgfmathsetmacro{\midXTwo}{(\XgI+\XgII)/2}
    \pgfmathsetmacro{\midXThree}{(\XgII+\tkzW)/2}
    \begin{scope}
    \clip[rounded corners] (0,-\sh) rectangle (\tkzW,\sh);
    \fill[primedregion=lightgray!25] (0,-\sh) rectangle (\XgI,\sh);
    \fill[primedregion=lightgray!50] (\XgI,-\sh) rectangle (\XgII,\sh);
    \fill[primedregion=lightgray!75] (\XgII,-\sh) rectangle (\tkzW,\sh);
    \end{scope}
    \draw[thick] (\XgI,-\sh) -- (\XgI,\sh);
    \draw[thick] (\XgII,-\sh) -- (\XgII,\sh);
    \pic at (\XgI,0){tro};
    \pic[rotate=180] at (\XgII,0){tro};
    \node[above] at (\XgI,\sh) {$F(X)$};
    \node[above] at (\XgII,\sh) {$F(W)$};
    \node at (\midXOne,0) {\tiny$F(a)$};
    \node at (\midXTwo,0) {\tiny$F(b)$};
    \node at (\midXThree,0) {\tiny$F(c)$};
\end{tkz}
\in\fld[1]{\C^\fY}{Y}[\fcj{F(a)}\amalg F(c)]
\]
(which are 1-fields, but we thicken them to illustrate their constituent sub-1-fields by shadings or patterns)
and
\[
\alpha
\;=\;
\begin{tkz}[scale=0.85]
    \def\lkr{0.15}
    \begin{scope}
    \clip (0,0) circle (1.3);
    \fill[lightgray!25] (-1.3,-1.3) rectangle (-0.5,1.3);
    \fill[lightgray!50] (-0.5,-1.3) rectangle (0.5,1.3);
    \fill[lightgray!75] (0.5,-1.3) rectangle (1.3,1.3);
    \end{scope}
    \draw[thick] (-0.5,-1.2) --pic{tro} (-0.5,-0.35);
    \draw[thick] (-0.5,-0.35) --pic{tro} (-0.5,1.2);
    \draw[thick] (0.5,-1.2) --pic{tro} (0.5,0.35);
    \draw[thick] (0.5,0.35) --pic[rotate=180]{tro} (0.5,1.2);
    \draw[thick](-.5,-.35)--(.5,.35);  
    \pic[rotate=-35] at (-0.5,-0.35) {link};
    \pic[rotate=-110] at (0.5,0.35) {link};
    \fill (-0.5,-0.35) circle (1.5pt);
    \fill (0.5,0.35) circle (1.5pt);
    \node[font=\scriptsize, below] at (-0.5,-1.3) {$X$};
    \node[font=\scriptsize, above] at (-0.5,1.3) {$Y$};
    \node[font=\scriptsize, below] at (0.5,-1.2) {$Z$};
    \node[font=\scriptsize,above] at (0.5,1.2) {$W$};
    \node[font=\scriptsize, left=1.5pt] at (-0.5,-0.35) {$f$};
    \node[font=\scriptsize, above left=0.75pt] at (0.5,0.35) {$g$};
    \node[gray,font=\scriptsize] at (-0.9,0) {$a$};
    \node[gray,font=\scriptsize] at (0,-0.75) {$b$};
    \node[gray,font=\scriptsize] at (0.9,0) {$c$};
\end{tkz}
\in\fld[2]{\C^\fX}{X}
\qquad\rightsquigarrow\qquad
\C^F(\alpha)
\;=\;
\begin{tkz}[scale=0.85]
    \def\lkr{0.15}
    \begin{scope}
    \clip (0,0) circle (1.3);
    \fill[primedregion=lightgray!25] (-1.3,-1.3) rectangle (-0.5,1.3);
    \fill[primedregion=lightgray!25] (-1.3,-1.3) rectangle (-0.5,1.3);
    \fill[primedregion=lightgray!50] (-0.5,-1.3) rectangle (0.5,1.3);
    \fill[primedregion=lightgray!75] (0.5,-1.3) rectangle (1.3,1.3);
    \end{scope}
    \draw[thick] (-0.5,-1.2) --pic{tro} (-0.5,-0.35);
    \draw[thick] (-0.5,-0.35) --pic{tro} (-0.5,1.2);
    \draw[thick] (0.5,-1.2) --pic{tro} (0.5,0.35);
    \draw[thick] (0.5,0.35) --pic[rotate=180]{tro} (0.5,1.2);
    \draw[thick](-.5,-.35)--(.5,.35);  
    \pic[rotate=-35] at (-0.5,-0.35) {link};
    \pic[rotate=-110] at (0.5,0.35) {link};
    \fill (-0.5,-0.35) circle (1.5pt);
    \fill (0.5,0.35) circle (1.5pt);
    \node[font=\scriptsize, below] at (-0.5,-1.3) {$F(X)$};
    \node[font=\scriptsize, above] at (-0.5,1.3) {$F(Y)$};
    \node[font=\scriptsize, below] at (0.5,-1.2) {$F(Z)$};
    \node[font=\scriptsize,above] at (0.5,1.2) {$F(W)$};
    \node[font=\scriptsize, left=0pt] at (-0.5,-0.35) {$F(f)$};
    \node[font=\scriptsize, above left=0pt] at (0.5,0.35) {$F(g)$};
    \node[gray,font=\scriptsize] at (-0.9,0) {$F(a)$};
    \node[gray,font=\scriptsize] at (0,-0.75) {$F(b)$};
    \node[gray,font=\scriptsize] at (0.9,0) {$F(c)$};
\end{tkz}
\in\fld[2]{\C^\fY}{X}.
\]
As 0-fields on a 0-ball $P$, we may illustrate the assignment $F\mapsto\C^F$ as
\[
F
\;=\;
\begin{tkz}[scale=1.25]
    \fill[rounded corners,lightgray!25](-.3,-.3) rectangle (.3,.3);
    \node at (0,0){\scriptsize$F$};
\end{tkz}
\;\in\;
\fld[0]{\C^{\Hom(\fX\to\fY)}}{P}
\qquad\rightsquigarrow\qquad
\C^F
\;=\;
\begin{tkz}[scale=1.25] 
    \fill[rounded corners,primedregion=white,draw=black,dotted](-.3,-.3) rectangle (.3,.3);
    \node at (0,0){\scriptsize$F$};
\end{tkz}
\;\in\;
\fld[0]{\eHom(\C^\fX{\to}\C^\fY)}{P}
\]
(where we are thickening these 0-fields to 2D to illustrate their field content). 
\end{construction}

\begin{construction}[(2,1)-transfors from \texorpdfstring{$\Hom(\fX\to\fY)$}{Hom}-string diagrams]
\label{DL 2-functor from 2-functor-2}
A $\Hom(\fX\to\fY)$-string diagram $\xi$ on a 1-ball $J$ defines a $J$-shaped disk-like $(2,1)$-transfor $\C^\xi \in \eFun(\C^\fX{\to}\C^\fY)$ as follows. 

For $0$-fields $a \in \fld[0]{\C^\fX}{P} = \Obj(\fX)$, define $\C^\xi(a)$ to be the $\fY$-string diagram $\xi(a)$ on $J$ obtained from $\xi$ by replacing edge labels (functors) $F_e$ with $F_e(a)$ and vertex labels (natural transformations) $\eta_v$ with their $a$-component $(\eta_v)_a$. 

For $\zeta \in \fld{\C^\fX}{Y}$, the $2$-field $\C^\xi(\zeta) \in \fld[2]{\C^\fY}{Y \times J}$ is given by the product string diagram $\zeta\times\xi$, which is illustrated by the following example. For
\[
\xi
\;=\;
\begin{tkz}[scale=0.38]
    \def\tkzH{6.5}
    \def\sw{0.4}
    \def\YoI{2.2}
    \def\YoII{4.3}
    \pgfmathsetmacro{\midYOne}{\YoI/2}
    \pgfmathsetmacro{\midYTwo}{(\YoI+\YoII)/2}
    \pgfmathsetmacro{\midYThree}{(\YoII+\tkzH)/2}
    \begin{scope}
    \clip[rounded corners] (-\sw,0) rectangle (\sw,\tkzH);
    \fill[primedregion=white,draw=black,dotted] (-\sw,0) rectangle (\sw,\YoI);
    \fill[boxregion=white] (-\sw,\YoI) rectangle (\sw,\YoII);
    \fill[white] (-\sw,\YoII) rectangle (\sw,\tkzH);
    \fill[starregion=white]
        (-\sw,\YoII) rectangle (\sw,\tkzH);
    \end{scope}
    \draw[dotted, rounded corners] (-\sw,0) rectangle (\sw,\tkzH);
    \draw[thick] (-\sw,\YoI) -- (\sw,\YoI);
    \draw[thick] (-\sw,\YoII) -- (\sw,\YoII);
    \pic[thick,scale=1.25,rotate=90] at (0,\YoI){tro};
    \pic[thick,scale=1.25,rotate=90] at (0,\YoII){tro};
    \node[right=3pt] at (\sw,\YoI) {\scriptsize$\eta_1$};
    \node[right=3pt] at (\sw,\YoII) {\scriptsize$\eta_2$};
    \node[right=3pt] at (\sw,\midYOne) {\tiny$F_1$};
    \node[right=3pt] at (\sw,\midYTwo) {\tiny$F_2$};
    \node[right=3pt] at (\sw,\midYThree) {\tiny$F_3$};
\end{tkz}
\;\in\fld[1]{\C^{\Hom(\fX\to\fY)}}{J}[\,\fcj{F_1}\amalg F_3]
\qquad\text{and}\qquad
\zeta
\;=\;
\begin{tkz}[baseline=-.5ex,scale=0.275,font=\scriptsize]
    \def\tkzW{10.5}
    \def\sh{0.6}
    \def\XgI{3.5}
    \def\XgII{7.0}
    \pgfmathsetmacro{\midXOne}{\XgI/2}
    \pgfmathsetmacro{\midXTwo}{(\XgI+\XgII)/2}
    \pgfmathsetmacro{\midXThree}{(\XgII+\tkzW)/2}
    \begin{scope}
    \clip[rounded corners] (0,-\sh) rectangle (\tkzW,\sh);
    \fill[lightgray!25] (0,-\sh) rectangle (\XgI,\sh);
    \fill[lightgray!50] (\XgI,-\sh) rectangle (\XgII,\sh);
    \fill[lightgray!75] (\XgII,-\sh) rectangle (\tkzW,\sh);
    \end{scope}
    \draw[thick] (\XgI,-\sh) -- (\XgI,\sh);
    \draw[thick] (\XgII,-\sh) -- (\XgII,\sh);
        \pic at (\XgI,0){tro};
    \pic[rotate=180] at (\XgII,0){tro};
    \node[above] at (\XgI,\sh) {$X_1$};
    \node[above] at (\XgII,\sh) {$X_2$};
    \node at (\midXOne,0) {\tiny$a$};
    \node at (\midXTwo,0) {\tiny$b$};
    \node at (\midXThree,0) {\tiny$c$};
\end{tkz}
\in\fld[1]{\C^\fX}{Y}[\fcj{a}\amalg c]
\]
(where again we are thickening these 1-fields to allow for shadings and patterns), the pure 2-field $\C^\xi(\zeta)\in\fld[2]{\C^\fY}{Y\times J}$ is the product string diagram
\[
\C^\xi(\zeta)
\;\;=\;\;\;
\begin{tkz}[font=\scriptsize,scale=0.85]
    \begin{scope}
    \clip[rounded corners] (0,0) rectangle (6.0,3.2);
    \fill[primedregion=lightgray!25] (0,0) rectangle (2.0,1.0);
    \fill[primedregion=lightgray!50] (2.0,0) rectangle (4.0,1.0);
    \fill[primedregion=lightgray!75] (4.0,0) rectangle (6.0,1.0);
    \fill[boxregion=lightgray!25] (0,1.0) rectangle (2.0,2.2);
    \fill[boxregion=lightgray!50] (2.0,1.0) rectangle (4.0,2.2);
    \fill[boxregion=lightgray!75] (4.0,1.0) rectangle (6.0,2.2);
    \fill[lightgray!25] (0,2.2) rectangle (2.0,3.2);
    \fill[starregion=lightgray!25] (0,2.2) rectangle (2.0,3.2);
    \fill[lightgray!50] (2.0,2.2) rectangle (4.0,3.2);
    \fill[starregion=lightgray!50] (2.0,2.2) rectangle (4.0,3.2);
    \fill[lightgray!75] (4.0,2.2) rectangle (6.0,3.2);
    \fill[starregion=lightgray!75](4.0,2.2) rectangle (6.0,3.2);
    \end{scope}
    \draw[thick] (0,1.0) --pic[rotate=90]{tro} (2.0,1.0) --pic[rotate=90]{tro} (4.0,1.0) --pic[rotate=90]{tro} (6.0,1.0);
    \draw[thick] (0,2.2) --pic[rotate=90]{tro} (2.0,2.2) --pic[rotate=90]{tro} (4.0,2.2) --pic[rotate=90]{tro} (6.0,2.2);
    \draw[thick] (2.0,0) --pic{tro} (2.0,1.0) --pic{tro} (2.0,2.2) --pic{tro} (2.0,3.2);
    \draw[thick] (4.0,0) --pic[rotate=180]{tro} (4.0,1.0) --pic[rotate=180]{tro} (4.0,2.2) --pic[rotate=180]{tro} (4.0,3.2);
    \pic[rotate=-45,scale=0.75] at (2.0,1.0) {link};
    \pic[rotate=45,scale=0.75]  at (4.0,1.0) {link};
    \pic[rotate=-45,scale=0.75] at (2.0,2.2) {link};
    \pic[rotate=45,scale=0.75]  at (4.0,2.2) {link};
    \filldraw[fill=black,thick] (2.0,1.0) circle (1.5pt);
    \filldraw[fill=black,thick] (4.0,1.0) circle (1.5pt);
    \filldraw[fill=black,thick] (2.0,2.2) circle (1.5pt);
    \filldraw[fill=black,thick] (4.0,2.2) circle (1.5pt);
    \node[left,font=\tiny] at (0,1.0) {$(\eta_1)_a$};
    \node[left,font=\tiny] at (0,2.2) {$(\eta_2)_a$};
    \node[right,font=\tiny] at (6.0,1.0) {$(\eta_1)_c$};
    \node[right,font=\tiny] at (6.0,2.2) {$(\eta_2)_c$};
    \node[below,font=\tiny] at (2.0,0) {$F_1(X_1)$};
    \node[below,font=\tiny] at (4.0,0) {$F_1(X_2)$};
    \node[above,font=\tiny] at (2.0,3.2) {$F_3(X_1)$};
    \node[above,font=\tiny] at (4.0,3.2) {$F_3(X_2)$};
\end{tkz}
\]
whose $(\eta_i,X_j)$-vertex is labeled by the naturator $(\eta_i)_{X_j}$ and has link parameterization determined by the transverse orientations of the corresponding $\eta_i$- and $X_j$-strands.

As 1-fields on a 1-ball $J$, we may illustrate the assignment $\xi\mapsto\C^\xi$ as
\[
\xi
\;=\;
\begin{tkz}[scale=0.4]
    \def\tkzH{6.5}
    \def\sw{0.5}
    \def\YoI{2.2}
    \def\YoII{4.3}
    \pgfmathsetmacro{\midYOne}{\YoI/2}
    \pgfmathsetmacro{\midYTwo}{(\YoI+\YoII)/2}
    \pgfmathsetmacro{\midYThree}{(\YoII+\tkzH)/2}
    \begin{scope}
    \clip[rounded corners] (-\sw,0) rectangle (\sw,\tkzH);
    \fill[lightgray!25] (-\sw,0) rectangle (\sw,\YoI);
    \fill[lightgray!50] (-\sw,\YoI) rectangle (\sw,\YoII);
    \fill[white] (-\sw,\YoII) rectangle (\sw,\tkzH);
    \fill[lightgray!75]
        (-\sw,\YoII) rectangle (\sw,\tkzH);
    \end{scope}
    \draw[thick] (-\sw,\YoI) -- (\sw,\YoI);
    \draw[thick] (-\sw,\YoII) -- (\sw,\YoII);
    \pic[thick,scale=1.25,rotate=90] at (0,\YoI){tro};
    \pic[thick,scale=1.25,rotate=90] at (0,\YoII){tro};
    \node[right=2pt] at (\sw,\YoI) {\scriptsize$\eta_1$};
    \node[right=2pt] at (\sw,\YoII) {\scriptsize$\eta_2$};
    \node[xshift=-2mm] at (\sw,\midYOne) {\tiny$F_1$};
    \node[xshift=-2mm] at (\sw,\midYTwo) {\tiny$F_2$};
    \node[xshift=-2mm] at (\sw,\midYThree) {\tiny$F_3$};
\end{tkz}
\;\in\;
\fld[1]{\C^{\Hom(\fX\to\fY)}}{J}
\qquad\rightsquigarrow\qquad
\C^\xi
\;=\;
\begin{tkz}[scale=0.4]
    \def\tkzH{6.5}
    \def\sw{0.5}
    \def\YoI{2.2}
    \def\YoII{4.3}
    \pgfmathsetmacro{\midYOne}{\YoI/2}
    \pgfmathsetmacro{\midYTwo}{(\YoI+\YoII)/2}
    \pgfmathsetmacro{\midYThree}{(\YoII+\tkzH)/2}
    \begin{scope}
    \clip[rounded corners] (-\sw,0) rectangle (\sw,\tkzH);
    \fill[primedregion=white,draw=black,dotted] (-\sw,0) rectangle (\sw,\YoI);
    \fill[boxregion=white] (-\sw,\YoI) rectangle (\sw,\YoII);
    \fill[white] (-\sw,\YoII) rectangle (\sw,\tkzH);
    \fill[starregion=white]
        (-\sw,\YoII) rectangle (\sw,\tkzH);
    \end{scope}
    \draw[dotted, rounded corners] (-\sw,0) rectangle (\sw,\tkzH);
    \draw[thick] (-\sw,\YoI) -- (\sw,\YoI);
    \draw[thick] (-\sw,\YoII) -- (\sw,\YoII);
    \pic[thick,scale=1.25,rotate=90] at (0,\YoI){tro};
    \pic[thick,scale=1.25,rotate=90] at (0,\YoII){tro};
    \node[right=2pt] at (\sw,\YoI) {\scriptsize$\eta_1$};
    \node[right=2pt] at (\sw,\YoII) {\scriptsize$\eta_2$};
    \node[right=2pt] at (\sw,\midYOne) {\tiny$F_1$};
    \node[right=2pt] at (\sw,\midYTwo) {\tiny$F_2$};
    \node[right=2pt] at (\sw,\midYThree) {\tiny$F_3$};
\end{tkz}
\;\in\;
\fld[1]{\eHom(\C^\fX\to\C^\fY)}{J}
\]
(where again we are thickening these 1-fields to 2D to illustrate their field content). Naturality \ref{axiom:eT-naturality} holds for a pure 2-field $f$ in $\C^\fX$ because the two hemisphere 2-fields have equal evaluations---for the standard splitting this is the naturality 2-cell of $\eval(\xi)$ at $f$, and for arbitrary splittings this follows by rotating both sides as in the proof of \cref{21lemma} below---so they are identified by \cref{lem-eval-equiv-2}.
\end{construction}

\begin{construction}[(2,2)-transfors from \texorpdfstring{$\Hom(\fX\to\fY)$}{Hom}-string diagrams]
\label{DL 2-functor from 2-functor-3}
For each $\Hom(\fX\to\fY)$-string diagram $\alpha$ on a 2-ball $X$ with region labels $F_r$ (functors), edge labels $\eta_e$ (natural transformations), and vertex labels $m_v$ (modifications), we define $\C^\alpha$ to be the $X$-shaped disk-like $(2,2)$-transfor whose component at $a \in \fX$ is the $\fY$-string diagram $\C^\alpha(a)$ on $X$ obtained from $\alpha$ by replacing each region label $F_r$ with $F_r(a)$, each edge label $\eta_e$ with $(\eta_e)_a$, and each vertex label $m_v$ with $(m_v)_a$. For example,
\[
\alpha
=
\begin{tkz}[scale=0.75]
    \begin{scope}
    \clip (0,0) circle (1.3);
    \fill[primedregion=white,draw=black,dotted] (-1.3,-1.3) rectangle (0,1.3);
    \fill[boxregion=white] (0,-1.3) rectangle (1.3,1.3);
    \fill[starregion=white] (0,0) -- (0,2)--(150:2)--cycle;
    \end{scope}
    \draw[dotted] (0,0) circle (1.3);
    \draw[thick] (0,-1.3) --pic{tro} (0,0);
    \draw[thick] (0,0) --pic{tro} (0,1.3);
    \draw[thick] (0,0) --pic[rotate=30]{tro} (150:1.3)node[above left=-1mm,overlay]{\scriptsize$\gamma$};
    \fill (0,1.3) circle (1.5pt);
    \fill (0,-1.3) circle (1.5pt);
    \fill (150:1.3) circle (1.5pt);
    \pic[rotate=30] at (0,0){link};
    \node[font=\scriptsize] at (0,1.5) {$\eta$};
    \node[font=\scriptsize] at (0,-1.5) {$\mu$};
    \fill (0,0) circle (1.5pt);
    \node[below left] at (0,0){$m$};
\end{tkz}
\in\fld[2]{\C^{\Hom(\fX\to\fY)}}{D^2}[c]
\;\rightsquigarrow\;
\C^\alpha(a)
=
\begin{tkz}[scale=0.75]
    \begin{scope}
    \clip (0,0) circle (1.3);
    \fill[primedregion=lightgray!25] (-1.3,-1.3) rectangle (0,1.3);
    \fill[boxregion=lightgray!25] (0,-1.3) rectangle (1.3,1.3);
    \fill[starregion=lightgray!25] (0,0) -- (0,2)--(150:2)--cycle;
    \end{scope}
    \draw[dotted] (0,0) circle (1.3);
    \draw[thick] (0,-1.3) --pic{tro} (0,0);
    \draw[thick] (0,0) --pic{tro} (0,1.3);
    \draw[thick] (0,0) --pic[rotate=30]{tro} (150:1.3)node[above left=-1mm,overlay]{\scriptsize$\gamma_a$};
    \fill (0,1.3) circle (1.5pt);
    \fill (0,-1.3) circle (1.5pt);
    \fill (150:1.3) circle (1.5pt);
    \pic[rotate=30] at (0,0){link};
    \node[font=\scriptsize] at (0,1.5) {$\eta_a$};
    \node[font=\scriptsize] at (0,-1.5) {$\mu_a$};
    \fill (0,0) circle (1.5pt);
    \node[below left] at (0,0){\scriptsize$m_a$};
\end{tkz}
\in\fld[2]{\C^\fY}{D^2}[c(a)].
\]
As 2-fields on a 2-ball $D^2$, we may illustrate this assignment $\alpha\mapsto\C^\alpha$ below.
\[
\alpha \;=\; \begin{tkz}[scale=0.75]     \begin{scope}     \clip (0,0) circle (1.3);     
\fill[lightgray!25] (-1.3,-1.3) rectangle (0,1.3);     \fill[lightgray!50] (0,-1.3) rectangle (1.3,1.3);     \fill[lightgray!75] (0,0) -- (0,2)--(150:2)--cycle;     \end{scope}     
\draw[thick] (0,-1.3) --pic{tro} (0,0);     \draw[thick] (0,0) --pic{tro} (0,1.3);     \draw[thick] (0,0) --pic[rotate=30]{tro} (150:1.3)node[above left=-1mm,overlay]{\scriptsize$\gamma$};     \fill (0,1.3) circle (1.5pt);     \fill (0,-1.3) circle (1.5pt);     \fill (150:1.3) circle (1.5pt);     \pic[rotate=30] at (0,0){link};     \node[font=\scriptsize] at (0,1.5) {$\eta$};     \node[font=\scriptsize] at (0,-1.5) {$\mu$};     \fill (0,0) circle (1.5pt);     \node[below left] at (0,0){\scriptsize$m$};     \node[font=\tiny] at (-0.6,-0.6) {$F_1$};     \node[font=\tiny] at (0.6,0) {$F_2$};     \node[font=\tiny] at (-0.4,0.8) {$F_3$}; \end{tkz} \in\fld[2]{\C^{\Hom(\fX\to\fY)}}{D^2} \quad\rightsquigarrow\quad \C^\alpha \;=\; 
\begin{tkz}[scale=0.75]  
\begin{scope}    
\clip (0,0) circle (1.3);    
\fill[primedregion=white,draw=black,dotted] (-1.3,-1.3) rectangle (0,1.3);     
\fill[boxregion=white] (0,-1.3) rectangle (1.3,1.3);     
\fill[starregion=white] (0,0) -- (0,2)--(150:2)--cycle;     \end{scope}     
\draw[dotted] (0,0) circle (1.3);     
\draw[thick] (0,-1.3) --pic{tro} (0,0);     
\draw[thick] (0,0) --pic{tro} (0,1.3);     
\draw[thick] (0,0) --pic[rotate=30]{tro} (150:1.3)node[above left=-1mm,overlay]{\scriptsize$\gamma$};     
\fill (0,1.3) circle (1.5pt);     
\fill (0,-1.3) circle (1.5pt);     
\fill (150:1.3) circle (1.5pt);     
\pic[rotate=30] at (0,0){link};     
\node[font=\scriptsize] at (0,1.5) {$\eta$};     
\node[font=\scriptsize] at (0,-1.5) {$\mu$};     
\fill (0,0) circle (1.5pt);     \node[below left] at (0,0){\scriptsize$m$};     
\node[font=\tiny] at (-0.6,-0.6) {$F_1$};     
\node[font=\tiny] at (0.6,0) {$F_2$};     
\node[font=\tiny] at (-0.4,0.8) {$F_3$}; 
\end{tkz} 
\in \fld[2]{\eHom(\C^\fX{\to}\C^\fY)}{D^2}
\]
Here \ref{axiom:eT-naturality} holds for a pure 1-field in $\C^\fX$ because the two hemisphere 2-fields have equal evaluations by the modification axiom of the $m_v$, exactly as in the proof of \cref{21lemma} below, and thus they are identified by \cref{lem-eval-equiv-2}.
\end{construction}

Observe that for pivotal $\Cstar$-2-categories $\fX$ and $\fY$, a disk-like functor $\eF\colon\C^\fX\to\C^\fY$, and an $\fX$-string diagram $\xi\in\fld[k]{\C^\fX}{W}$ on a $k$-ball $W$ for $k\in\{1,2\}$, by \ref{axiom:eF-products} the underlying string diagram stratification $\Gamma_{\eF(\xi)}$ of $W$ is a subset of the underlying string diagram stratification $\Gamma_\xi$ of $W$.

\begin{lemma}[$\C^{(-)}$ is unitarily essentially surjective on (2,0)-transfors]
\label{lem:ess-surj-functors}
For pivotal $\Cstar$-2-categories $\fX$ and $\fY$ and a disk-like functor $\eF\colon\C^\fX\to\C^\fY$, define $F^\eF\colon\fX\to\fY$ by
\begin{align*}
F^\eF(a)&\coloneq\eF(a)\in\fY,
\\
F^\eF({}_aX_b)&\coloneq\eval^\fY(\eF\langle X\rangle)\in\fY(F^\eF(a)\to F^\eF(b)),
\\
F^\eF({}_X f_Y)&\coloneq\eval^\fY_{D^2,\partial_\pm D^2}(\eF\langle f\rangle)\in\fY(F^\eF(X)\Rightarrow F^\eF(Y)).
\end{align*}
Then $F^\eF$ is a UAF-preserving $\dag$-functor and $\C^{F^\eF}\cong^\star\eF$ in $\eHom(\C^\fX{\to}\C^\fY)$.
\end{lemma}

\begin{pf}
As $F^\eF=\eval\circ\fX_\eF\circ\langle-\rangle$ is a composition of UAF-preserving dagger functors, $F^\eF\in\Hom(\fX\to\fY)$. For the second assertion, we will exhibit a unitary equivalence $\eN\in\fld[1]{\eHom(\C^\fX{\to}\C^\fY)}{I}[\fcj{\C^{F^\eF}}\amalg\eF]$. Since $\C^{F^\eF}(a)=F^\eF(a)=\eF(a)$ for all $a\in\fX$, we set $\eN(a)\coloneq \eF(a)\times I$. On a 1-field (1D $\fX$-string diagram) $\xi\in\fld[1]{\C^\fX}{Y}[\fcj{a}\amalg b]$, we define $\eN(\xi)\in\fld[2]{\C^\fY}{Y\times I}[\fcj{\C^{F^\eF}(\xi)}\amalg\eF(\xi)]$ by the $\fY$-string diagram on $Y\times I$ given by the double cone of $\eF(\xi)$ and $\C^{F^\eF}(\xi)$ onto their common evaluation $E\coloneq \eval(\C^{F^\eF}(\xi))=\eval(\eF(\xi))$; to see this equality, take any splitting of $\xi$ into single-vertex string diagrams $\langle X_i\rangle$, so that $\eval(\eF\langle X_i\rangle)=F^\eF(X_i)=\eval(\C^{F^\eF}\langle X_i\rangle)$ by definition of $F^\eF$ and strictness of $\eval$.
In more detail, the edges of $\eN(\xi)$ are the strands $\{v\}\times I$ for each vertex $v$ in $\xi$, and the vertices of $\eN(\xi)$ comprise the set $\{(v,j/3)\mid j\in\{1,2\},\text{ $v$ is a vertex of $\xi$}\}$, with each $(v,j/3)$ given the ``upward'' link parameterization (illustrated below) and the label $\id_{E_i}$ where $E_i\coloneq\eval(\eF\langle X_i\rangle)$ and $X_i$ is the morphism labeling $v$ in $\xi$.

For example, for
\[
\xi
\;=\;
\begin{tkz}[baseline=-.5ex,scale=0.425,font=\scriptsize]
    \def\tkzW{10.5}
    \def\sh{0.45}
    \def\XgI{3.5}
    \def\XgII{7.0}
    \pgfmathsetmacro{\midXOne}{\XgI/2}
    \pgfmathsetmacro{\midXTwo}{(\XgI+\XgII)/2}
    \pgfmathsetmacro{\midXThree}{(\XgII+\tkzW)/2}
    \begin{scope}
    \clip[rounded corners] (0,-\sh) rectangle (\tkzW,\sh);
    \fill[lightgray!25] (0,-\sh) rectangle (\XgI,\sh);
    \fill[lightgray!50] (\XgI,-\sh) rectangle (\XgII,\sh);
    \fill[lightgray!75] (\XgII,-\sh) rectangle (\tkzW,\sh);
    \end{scope}
    \draw[thick] (\XgI,-\sh) -- (\XgI,\sh);
    \draw[thick] (\XgII,-\sh) -- (\XgII,\sh);
    \pic at (\XgI,0){tro};
    \pic[rotate=180] at (\XgII,0){tro};
    \node[above] at (\XgI,\sh) {$X_1$};
    \node[above] at (\XgII,\sh) {$X_2$};
    \node at (\midXOne,0) {\scriptsize$a$};
    \node at (\midXTwo,0) {\scriptsize$b$};
    \node at (\midXThree,0) {\scriptsize$c$};
\end{tkz}
\in
\fld[1]{\C^\fX}{Y}[\fcj{a}\amalg c]
\quad\text{and}\quad
\C^{F^\eF}
=\begin{tkz}[scale=1.35]
    \fill[rounded corners,primedregion=white,draw=black,dotted](-.3,-.3) rectangle (.3,.3);
    \node at (0,0){\scriptsize$F^\eF$};
\end{tkz}
\in\fld[0]{\eHom(\C^\fX{\to}\C^\fY)}{P},
\]
we have
\[
\eN(\xi)
\;\coloneq\;
\begin{tkz}[scale=0.55, font=\sffamily\small]
    \def\tkzW{10.5}
    \def\tkzH{4.5}
    \def\vxI{3.5}
    \def\vxII{7.0}
    \def\vbot{1.5}
    \def\vtop{3.0}
    \begin{scope}
        \clip[rounded corners] (0.125*\tkzW,0) rectangle (.9*\tkzW,\tkzH);
        
        \fill[primedregion=lightgray!25] (0,0) rectangle (\vxI,\tkzH);
        \fill[primedregion=lightgray!50] (\vxI,0) rectangle (\vxII,\tkzH);
        \fill[primedregion=lightgray!75] (\vxII,0) rectangle (\tkzW,\tkzH);
    \end{scope}
    
    \draw[thick, RoyalBlue!80!black] (\vxI,0) --pic[rotate=180]{tro} (\vxI,\vbot);
    \draw[thick, RoyalBlue!80!black] (\vxI,\vbot) --pic{tro} (\vxI,\vtop);
    \draw[thick, RoyalBlue!80!black] (\vxI,\vtop) --pic[rotate=180]{tro} (\vxI,\tkzH);
    
    \draw[thick, RoyalBlue!80!black] (\vxII,0) --pic[rotate=180]{tro} (\vxII,\vbot);
    \draw[thick, RoyalBlue!80!black] (\vxII,\vbot) --pic{tro} (\vxII,\vtop);
    \draw[thick, RoyalBlue!80!black] (\vxII,\vtop) --pic{tro} (\vxII,\tkzH);

    \filldraw[fill=white, draw=RoyalBlue!80!black, thick] (\vxI,\vbot) circle (2pt) pic{link};
    \filldraw[fill=white, draw=RoyalBlue!80!black, thick] (\vxI,\vtop) circle (2pt) pic{link};
    
    \filldraw[fill=white, draw=RoyalBlue!80!black, thick] (\vxII,\vbot) circle (2pt) pic{link};
    \filldraw[fill=white, draw=RoyalBlue!80!black, thick] (\vxII,\vtop) circle (2pt) pic{link};
    
    \node[right=2pt] at (\vxI, {(\vbot+\vtop)/2}) {\scriptsize$E_1$};
    
    \node[left=2pt] at (\vxII, {(\vbot+\vtop)/2}) {\scriptsize$E_2$};
    
    \node[left=8pt] at (\vxI,\vbot) {\scriptsize$\eval$};
    \node[left=8pt] at (\vxI,\vtop) {\scriptsize$\langle-\rangle$};
    
    \node[right=8pt] at (\vxII,\vbot) {\scriptsize$\eval$};
    \node[right=8pt] at (\vxII,\vtop) {\scriptsize$\langle-\rangle$};
    
    \node[above] at (\vxI, \tkzH) {\scriptsize$\eF(\langle X_1 \rangle)$};
    \node[below] at (\vxI, 0) {\scriptsize$F^\eF(X_1)$};
    \node[above] at (\vxII, \tkzH) {\scriptsize$\eF(\langle X_2 \rangle)$};
    \node[below] at (\vxII, 0) {\scriptsize$F^\eF(X_2)$};
    
\end{tkz}
\;\in\;
\fld[2]{\C^\fY}{Y{\times} I}[(\fcj{\C^{F^\eF}(\xi)}\amalg\eF(\xi))\cup(\fcj{\eF(a){\times} I}\amalg \eF(c){\times} I)].
\]
To see $\eN$ is unitary, note that because $\eN(a)=\eF(a)\times I$ is a product field, $\beta_\eN(a)=(\eN(a)\times S^1)\cup(\eF(a)\times D^2)=\eF(a)\times D^2=(\C^{F^\eF}\times D^2)(a)$ and $\beta_{\eN^\star}(a)=\eF(a)\times D^2=(\eF\times D^2)(a)$ for all $a$. Thus the (co)isometry bubble conditions (\cref{2coisometryunitarydl}) hold pointwise, and hence in $\eFun(\C^\fX{\to}\C^\fY)$ by \ref{DCU}, so $\eN$ is unitary.
\end{pf}

\begin{lemma}[$\C^{(-)}$ is unitarily essentially surjective on (2,1)-transfors]
\label{21lemma}
Let $F,G\colon\fX\to\fY$ be UAF-preserving $\dag$-functors between pivotal $\Cstar$-2-categories $\fX$ and $\fY$. If $\eN\in\fld[1]{\eHom(\C^\fX{\to}\C^\fY)}{I}[\fcj{\C^F}\amalg \C^G]$, then
\begin{align*}
\eta_a^\eN&\coloneq\eval(\eN(a))\in\fY(F(a)\to G(a)),
\\
\eta_{{}_aX_b}^\eN&\coloneq\eval_{D^2,\partial_{\pm}D^2}(\eN\langle X\rangle)\in\fY(\eta_a^\eN\otimes G(X)\Rightarrow F(X)\otimes\eta_b^\eN)
\end{align*}
defines a natural transformation $\eta^\eN\colon F\Rightarrow G$ with $\C^{\langle\eta^\eN\rangle}\cong^\star\eN$ in $\eHom(\C^\fX{\to}\C^\fY)$.
\end{lemma}

\begin{pf}
    We first observe that $\eta^\eN$ is a natural transformation as a whiskering $\eta^\eN\coloneq\eval\circ\fX_\eN\circ\langle-\rangle$ of the natural transformation $\fX_\eN$. For the second assertion, define $\eS\in\fld{\eHom(\C^\fX{\to}\C^\fY)}{D^2}[\fcj{\C^{\langle\eta^\eN\rangle}}\cup\eN]$ on 0-fields $a$ in $\C^\fX$ by the cone string diagram on $D^2$ with incoming boundary field $\C^{\langle\eta^\eN\rangle}(a)$ and outgoing boundary field $\eN(a)$ with vertex label given by the identity 2-morphism $\id_{\eval(\eN(a))}$ in $\fY$. For example, 
    \[
\eN(a) =
\begin{tkz}[baseline=-.5ex,scale=0.465,xscale=0.65,font=\scriptsize]
    \def\tkzW{12}
    \def\sh{0.5}
    \def\XgI{3}
    \def\XgII{6}
    \def\XgIII{9}
    \begin{scope}
    \clip[rounded corners] (0,-\sh) rectangle (\tkzW,\sh);
    \fill[primedregion=lightgray!25] (0,-\sh) rectangle (\XgI,\sh);
    \fill[orange!10] (\XgI,-\sh) rectangle (\XgII,\sh);
    \fill[orange!30] (\XgII,-\sh) rectangle (\XgIII,\sh);
    \fill[boxregion=lightgray!25] (\XgIII,-\sh) rectangle (\tkzW,\sh);
    \end{scope}
    \draw[thick] (\XgI,-\sh) -- (\XgI,\sh);
    \draw[thick] (\XgII,-\sh) -- (\XgII,\sh);
    \draw[thick] (\XgIII,-\sh) -- (\XgIII,\sh);
    \pic[rotate=180] at (\XgI,0){tro};
    \pic at (\XgII,0){tro};
    \pic[rotate=180] at (\XgIII,0){tro};
    \node[above] at (\XgI,\sh) {$Y_1$};
    \node[above] at (\XgII,\sh) {$Y_2$};
    \node[above] at (\XgIII,\sh) {$Y_3$};
    \node at (1.5,0) {\scriptsize$F(a)$};
    \node at (4.5,0) {\scriptsize$y_1$};
    \node at (7.5,0) {\scriptsize$y_2$};
    \node at (10.5,0) {\scriptsize$G(a)$};
\end{tkz}
\in\fld[1]{\C^\fY}{I}
\,\rightsquigarrow\,
\eS(a) \coloneq
\begin{tkz}[baseline=-0.5ex, x=1cm,y=1cm, scale=0.9]
\coordinate (R_center) at (0, 0);
\coordinate (R_left)   at (-1.6, 0);
\coordinate (R_right)  at (1.6, 0);
\path (R_left) to[out=60, in=120] coordinate[pos=0.25] (Top_1) coordinate[pos=0.5] (Top_2) coordinate[pos=0.75] (Top_3) (R_right);
\path (R_left) to[out=-60, in=-120] coordinate[pos=0.5] (Bot_Mid) (R_right);
\begin{scope}
    \clip (R_left) to[out=60, in=120] (R_right) to[out=-120, in=-60] (R_left) -- cycle;
    \fill[primedregion=lightgray!25] (R_center) -- ($(R_center)!3!(Top_1)$) -- (-2, 3) -- (-2,-3) -- ($(R_center)!3!(Bot_Mid)$) -- cycle;
    \fill[orange!10] (R_center) -- ($(R_center)!3!(Top_1)$) -- ($(R_center)!3!(Top_2)$) -- cycle;
    \fill[orange!30] (R_center) -- ($(R_center)!3!(Top_2)$) -- ($(R_center)!3!(Top_3)$) -- cycle;
    \fill[boxregion=lightgray!25] (R_center) -- ($(R_center)!3!(Top_3)$) -- (2, 3) -- (2,-3) -- ($(R_center)!3!(Bot_Mid)$) -- cycle;
\end{scope}
\path (R_center) -- (Top_1) coordinate[pos=0.75] (T1t);
\path (R_center) -- (Top_2) coordinate[pos=0.75] (T2t);
\path (R_center) -- (Top_3) coordinate[pos=0.75] (T3t);
\path (Bot_Mid) -- (R_center) coordinate[pos=0.25] (Bt);
\draw[thick] (R_center) -- (Top_1);
\draw[thick] (R_center) -- (Top_2);
\draw[thick] (R_center) -- (Top_3);
\draw[thick] (Bot_Mid) -- (R_center);
\pic[rotate=235] at (T1t) {tro};
\pic at (T2t) {tro};
\pic[rotate=125] at (T3t) {tro};
\pic at (Bt) {tro};
\pic[scale=1.5] at (R_center) {link};
\fill (R_center) circle (1.5pt);
\node[font=\footnotesize,left=-1pt,yshift=-7pt] at (R_center) {\tiny$\id_{\eta^\eN_a}$};
\fill (Top_1) circle (1.5pt);
\fill (Top_2) circle (1.5pt);
\fill (Top_3) circle (1.5pt);
\fill (Bot_Mid) circle (1.5pt);
\node[above=1pt, font=\scriptsize] at (Top_1) {$Y_1$};
\node[above=1pt, font=\scriptsize] at (Top_2) {$Y_2$};
\node[above=1pt, font=\scriptsize] at (Top_3) {$Y_3$};
\node[below=1pt, font=\scriptsize] at (Bot_Mid) {$\eta^\eN_a$};
\end{tkz}
\in\fld[2]{\C^\fY}{D^2}[\fcj{\C^{\langle \eta^\eN \rangle}(a)}\cup\eN(a)].
\]
It is a straightforward check that $\eS$ satisfies the axioms of a (2,2)-transfor other than naturality, so we will only verify \ref{axiom:eT-naturality}. Fix a 1-field $\xi$ in $\C^\fX$. An arbitrary two-ball splitting of $\partial D^2$ differs from the standard one $\partial D^2=\partial_-D^2\cup\partial_+D^2$ by transferring pieces of $(\partial\eS)(\xi)$ across the equator of $S^1$, which changes the evaluations of both hemispheres by the same rotations of the pivotal graphical calculus (\cref{eval-facts-2}), so the general case follows from the standard one. For the standard splitting $\partial_\pm D^2$ of $\partial Q\cong\partial D^2$, the two hemisphere 2-fields in the statement of \ref{axiom:eT-naturality} are $\eN(\xi)\blt\eS(a)$ and $\eS(b)\blt\C^{\langle\eta^\eN\rangle}(\xi)$, so by \cref{lem-eval-equiv-2} it suffices to show these have equal evaluations.

By the gluing axiom for $\eS$ and strictness of $\eval$, it suffices to assume $\xi=\langle X\rangle$ for some 1-morphism ${}_aX_b$ in $\fX$.
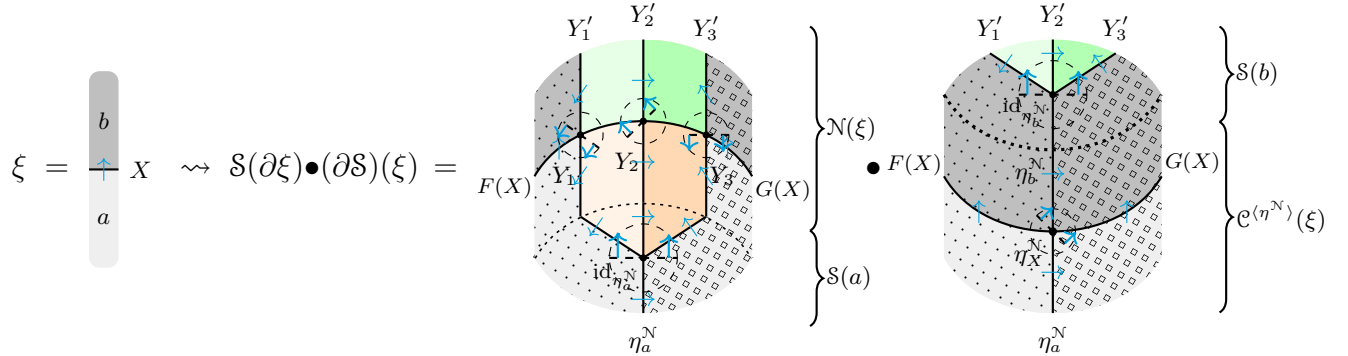
\[
\xi
\;=\;
\begin{tkz}[baseline=-.5ex,scale=.65,font=\scriptsize]
    \def\tkzH{2}
    \def\sw{0.3}
    \begin{scope}
    \clip[rounded corners] (-\sw,-\tkzH) rectangle (\sw,\tkzH);
    \fill[lightgray!25] (-\sw,-\tkzH) rectangle (\sw,0);
    \fill[lightgray] (-\sw,0) rectangle (\sw,\tkzH);
    \end{scope}
    \draw[thick] (-\sw,0) -- (\sw,0);
    \pic[rotate=90] at (0,0){tro};
    \node[right] at (\sw,0) {$X$};
    \node at (0,{-\tkzH/2}) {\scriptsize$a$};
    \node at (0,{\tkzH/2}) {\scriptsize$b$};
\end{tkz}
\;\rightsquigarrow\;
\eS(\partial\xi)\blt(\partial\eS)(\xi)
\,=
\begin{tkz}[scale=0.9, baseline=0]
    \begin{scope}[shift={(0,-1.2)}]
        \coordinate (R_center) at (0, 0);
        \coordinate (R_left)   at (-1.6, 0);
        \coordinate (R_right)  at (1.6, 0);
        \path (R_left) to[out=60, in=120] coordinate[pos=0.25] (Top_1) coordinate[pos=0.5] (Top_2) coordinate[pos=0.75] (Top_3) (R_right);
        \path (R_left) to[out=-60, in=-120] coordinate[pos=0.5] (Bot_Mid) (R_right);
        
        \begin{scope}
            \clip (R_left) to[out=60, in=120] (R_right) to[out=-120, in=-60] (R_left) -- cycle;
            \fill[primedregion=lightgray!25] (R_center) -- ($(R_center)!3!(Top_1)$) -- (-3, 3) -- (-3,-3) -- ($(R_center)!3!(Bot_Mid)$) -- cycle;
            \fill[orange!10] (R_center) -- ($(R_center)!3!(Top_1)$) -- ($(R_center)!3!(Top_2)$) -- cycle;
            \fill[orange!30] (R_center) -- ($(R_center)!3!(Top_2)$) -- ($(R_center)!3!(Top_3)$) -- cycle;
            \fill[boxregion=lightgray!25] (R_center) -- ($(R_center)!3!(Top_3)$) -- (3, 3) -- (3,-3) -- ($(R_center)!3!(Bot_Mid)$) -- cycle;
        \end{scope}
        
        \draw[very thick, dotted] (R_left) to[out=60, in=120] (R_right);
        
        \path (R_center) -- (Top_1) coordinate[pos=0.75] (T1t);
        \path (R_center) -- (Top_2) coordinate[pos=0.75] (T2t);
        \path (R_center) -- (Top_3) coordinate[pos=0.75] (T3t);
        \path (Bot_Mid) -- (R_center) coordinate[pos=0.25] (Bt);
        
        \draw[thick] (R_center) -- (Top_1);
        \draw[thick] (R_center) -- (Top_2);
        \draw[thick] (R_center) -- (Top_3);
        \draw[thick] (Bot_Mid) -- (R_center);
        
        \pic[rotate=235] at (T1t) {tro};
        \pic at (T2t) {tro};
        \pic[rotate=125] at (T3t) {tro};
        \pic at (Bt) {tro};
        
        \pic[scale=1.5] at (R_center) {link};
        \fill (R_center) circle (1.5pt);
        \node[font=\tiny, left=-3pt, yshift=-7pt] at (R_center) {$\id_{\eta^\eN_a}$};
        
        \node[below=1pt, font=\scriptsize] at (Bot_Mid) {$\eta^\eN_a$};
    \end{scope}
    \begin{scope}
        \path (-1.6,-1.2) to[out=60,in=120] coordinate[pos=0.25](B1) coordinate[pos=0.5](B2) coordinate[pos=0.75](B3) (1.6,-1.2);
        \path (-1.6,1.2)  to[out=60,in=120] coordinate[pos=0.25](T1) coordinate[pos=0.5](T2) coordinate[pos=0.75](T3) (1.6,1.2);
        \path (-1.6,0)    to[out=60,in=120] coordinate[pos=0.25](C1) coordinate[pos=0.5](C2) coordinate[pos=0.75](C3) coordinate[pos=0.01](FX) coordinate[pos=1](GX) coordinate[pos=0.15](TroL) coordinate[pos=0.85](TroR) (1.6,0);

        \begin{scope}
            \clip (-1.6,-1.2) to[out=60,in=120] (1.6,-1.2) -- (1.6,1.2) to[out=120,in=60] (-1.6,1.2) -- cycle;
            \fill[primedregion=lightgray!25] (-2,-2) rectangle (T1 |- 0,2);
            \fill[orange!10] (T1 |- 0,-2) rectangle (T2 |- 0,2);
            \fill[orange!30] (T2 |- 0,-2) rectangle (T3 |- 0,2);
            \fill[boxregion=lightgray!25] (T3 |- 0,-2) rectangle (2,2);
            \begin{scope}
                \clip (-1.6,0) to[out=60,in=120] (1.6,0) -- (1.6,2.5) -- (-1.6,2.5) -- cycle;
                \fill[primedregion=lightgray] (-2,-2) rectangle (T1 |- 0,2);
                \fill[green!10] (T1 |- 0,-2) rectangle (T2 |- 0,2);
                \fill[green!30] (T2 |- 0,-2) rectangle (T3 |- 0,2);
                \fill[boxregion=lightgray] (T3 |- 0,-2) rectangle (2,2);
            \end{scope}
        \end{scope}

        \draw[thick] (-1.6,0) to[out=60,in=120] (1.6,0);
        \pic[rotate=90] at (TroL) {tro};
        \pic[rotate=90] at (TroR) {tro};

        \draw[thick] (B1) -- (T1) coordinate[pos=0.25] (V1B) coordinate[pos=0.75] (V1T);
        \draw[thick] (B2) -- (T2) coordinate[pos=0.25] (V2B) coordinate[pos=0.75] (V2T);
        \draw[thick] (B3) -- (T3) coordinate[pos=0.25] (V3B) coordinate[pos=0.75] (V3T);

        \pic[rotate=235] at (V1B) {tro};
        \node[left=-1mm] at (V1B) {\scriptsize$Y_1$};
        \node[left=-1mm] at (V2B) {\scriptsize$Y_2$};
        \node[right=-1mm] at (V3B) {\scriptsize$Y_3$};
        
        \pic at (V2B) {tro};
        \pic[rotate=125] at (V3B) {tro};
        \pic[rotate=235] at (V1T) {tro};
        \pic at (V2T) {tro};
        \pic[rotate=125] at (V3T) {tro};

        \pic[rotate=145] at (C1) {link};
        \pic[rotate=45] at (C2) {link};
        \pic[rotate=180] at (C3) {link};

        \fill (C1) circle (1.5pt);
        \fill (C2) circle (1.5pt);
        \fill (C3) circle (1.5pt);
        
        \node[above=1pt, font=\scriptsize] at (T1) {$Y_1'$};
        \node[above=1pt, font=\scriptsize] at (T2) {$Y_2'$};
        \node[above=1pt, font=\scriptsize] at (T3) {$Y_3'$};

        \node[below left=-1mm, font=\tiny] at (FX) {$F(X)$};
        \node[below right=-1mm, font=\tiny] at (GX) {$G(X)$};
    \end{scope}

    \draw[brace] (2.45,2.2) -- node[right=2pt]{\scriptsize$\eN(\xi)$} (2.45,-0.8);
    \draw[brace] (2.45,-0.8) -- node[right=2pt]{\scriptsize$\eS(a)$} (2.45,-2.2);
\end{tkz}
\!\!\!\!\blt\!
\begin{tkz}[scale=.9, baseline=0]
    \path (-1.6,0) to[out=-60,in=-120] 
        coordinate[pos=0.01](FX) 
        coordinate[pos=0.2] (TroL)
        coordinate[pos=0.5] (C)
        coordinate[pos=0.8] (TroR) 
        coordinate[pos=1] (GX) 
        (1.6,0);
    \path (-1.6,-1.2) to[out=-60,in=-120] coordinate[pos=0.5] (Bot) (1.6,-1.2);
    \path (-1.6,1.2) to[out=-60,in=-120] coordinate[pos=0.5] (Top) (1.6,1.2);
    
    \begin{scope}
        \clip (-1.6,1.2) to[out=-60,in=-120] (1.6,1.2) -- (1.6,-1.2) to[out=-120,in=-60] (-1.6,-1.2) -- cycle;
        
        \fill[primedregion=lightgray!25] (-2,-2) rectangle (0,2);
        \fill[boxregion=lightgray!25] (0,-2) rectangle (2,2);
        
        \begin{scope}
            \clip (-1.6,0) to[out=-60,in=-120] (1.6,0) -- (1.6,2.5) -- (-1.6,2.5) -- cycle;
            \fill[primedregion=lightgray] (-2,-2) rectangle (0,2);
            \fill[boxregion=lightgray] (0,-2) rectangle (2,2);
        \end{scope}
    \end{scope}
    
    \draw[thick] (-1.6,0) to[out=-60,in=-120] (1.6,0);
    \pic[rotate=90] at (TroL) {tro};
    \pic[rotate=90] at (TroR) {tro};
    
    \draw[thick] (Bot) -- (Top);
    
    \path (Bot) -- (C) coordinate[pos=0.5] (VB);
    \path (C) -- (Top) coordinate[pos=0.5] (VT);
    
    \pic at (VB) {tro};
    \pic[rotate=-45] at (C) {link};
    
    \filldraw[fill=black] (C) circle (1.5pt);
    \node[below left=-1pt, font=\scriptsize] at (C) {$\eta^\eN_X$};
    \node[above left=-1mm,font=\tiny] at (FX) {$F(X)$};
    \node[above right=-1mm, font=\tiny] at (GX) {$G(X)$};
    \node[below=1pt, font=\scriptsize] at (Bot) {$\eta^\eN_a$};

    \begin{scope}[shift={(0,1.2)}]
        \coordinate (bCenter) at (0, 0);
        \path (-1.6,0) to[out=60,in=120] coordinate[pos=0.25] (bT1) coordinate[pos=0.5] (bT2) coordinate[pos=0.75] (bT3) (1.6,0);
        \path (-1.6,0) to[out=-60,in=-120] coordinate[pos=0.5] (bBot) (1.6,0);
        
        \begin{scope}
            \clip (-1.6,0) to[out=60,in=120] (1.6,0) to[out=-120,in=-60] (-1.6,0) -- cycle;
            
            \fill[primedregion=lightgray] (bCenter) -- ($(bCenter)!3!(bT1)$) -- (-3, 3) -- (-3,-3) -- ($(bCenter)!3!(bBot)$) -- cycle;
            \fill[green!10] (bCenter) -- ($(bCenter)!3!(bT1)$) -- ($(bCenter)!3!(bT2)$) -- cycle;
            \fill[green!30] (bCenter) -- ($(bCenter)!3!(bT2)$) -- ($(bCenter)!3!(bT3)$) -- cycle;
            \fill[boxregion=lightgray] (bCenter) -- ($(bCenter)!3!(bT3)$) -- (3, 3) -- (3,-3) -- ($(bCenter)!3!(bBot)$) -- cycle;
        \end{scope}
        
        \draw[very thick, dotted] (-1.6,0) to[out=-60,in=-120] (1.6,0);
        
        \path (bCenter) -- (bT1) coordinate[pos=0.75] (bT1t);
        \path (bCenter) -- (bT2) coordinate[pos=0.75] (bT2t);
        \path (bCenter) -- (bT3) coordinate[pos=0.75] (bT3t);
        \path (bBot) -- (bCenter) coordinate[pos=0.25] (bBt);
        
        \draw[thick] (bBot) -- (bCenter);
        \draw[thick] (bCenter) -- (bT1);
        \draw[thick] (bCenter) -- (bT2);
        \draw[thick] (bCenter) -- (bT3);
        
        \pic[rotate=235] at (bT1t) {tro};
        \pic at (bT2t) {tro};
        \pic[rotate=125] at (bT3t) {tro};
        \pic[yshift=-.5cm] at (bBt) {tro};
        
        \pic[scale=1.5] at (bCenter) {link};
        \fill (bCenter) circle (1.5pt);
        \node[font=\tiny, below left=-3pt] at (bCenter) {$\id_{\eta^\eN_b}$};
        \node[below left=-1pt, font=\scriptsize] at (bBot) {$\eta^\eN_b$};
        
        \node[above=1pt, font=\scriptsize] at (bT1) {$Y_1'$};
        \node[above=1pt, font=\scriptsize] at (bT2) {$Y_2'$};
        \node[above=1pt, font=\scriptsize] at (bT3) {$Y_3'$};
    \end{scope}
    
    \draw[brace] (2.45,2.2) -- node[right=2pt]{\scriptsize$\eS(b)$} (2.45,0.8);
    \draw[brace] (2.45,0.8) -- node[right=2pt]{\scriptsize$\C^{\langle \eta^\eN \rangle}(\xi)$} (2.45,-2.0);
\end{tkz}
\]
It remains to show that the above two diagrams evaluate to the same 2-morphism in $\fY$. And indeed, with respect to the evident boundary splitting of $X\times\{0\}\cup \partial_-X\times I\cong X\times\{1\}\cup \partial_+X\times I$ as illustrated in the running example above, $\eN(\xi)\blt\eS(a)$ (the left diagram) has evaluation $\eval(\eN(\xi))\circ\eval(\langle\id_{\eta_a^\eN}\rangle)=\eval(\eN(\xi))\circ\id_{\eta_a^\eN}=\eval(\eN(\xi))$ and $\eS(b)\blt\C^{\langle\eta^\eN\rangle}(\xi)$ (the right diagram) has evaluation $\eval(\langle\id_{\eta_b^\eN}\rangle)\circ\eval(\langle\eta_{X}^\eN\rangle)=\id_{\eta_b^\eN}\circ \eval(\langle\eta_{X}^\eN\rangle)=\eval(\eN(\xi))$, where the last equality is by \cref{lem-eval-equiv-2} and the definition of $\eta^\eN$. (In other words, if we evaluate the above two diagrams with respect to the ``read from the bottom to the top''-boundary splitting of $X\times\{0\}\cup \partial_-X\times I$ and $X\times\{1\}\cup \partial_+X\times I$ as illustrated above, then we get the same 2-morphism of $\fY$.) Therefore $\eS$ satisfies \ref{axiom:eT-naturality}, so $\eS$ is a valid (2,2)-transfor. To see $\eS$ is a unitary equivalence, we need to prove $\beta_\eS\defeq{}\eS\blt\eS^\star$ equals $\C^{\langle\eta^\eN\rangle}\times I$ and $\beta_{\eS^\star}\defeq{}\eS^\star\blt\eS$ equals $\eN\times I$ as 2-fields in $\eHom(\C^\fX{\to}\C^\fY)$. By definition of $\lU{\eFun(\C^\fX{\to}\C^\fY)}$, it suffices to check this pointwise on 0-fields. And indeed, for each $a\in\fX$, 
\(
\eval((\eS\blt\eS^\star)(a))
=\id_{\eta_a^\eN}^\dag\circ\id_{\eta_a^\eN}
=\id_{\eta_a^\eN}
=\eval((\C^{\langle\eta^\eN\rangle}\times I)(a))
\)
and $\eval((\eS^\star\blt\eS)(a))=\eval((\eN\times I)(a))$ similarly. Thus $\C^{\langle\eta^\eN\rangle}\cong^\star\eN$.
\end{pf}

\begin{lemma}[$\C^{(-)}$ is surjective on (2,2)-transfors]
\label{lem:surj-2-fields}
For $\Hom(\fX\to\fY)$-string diagrams $\xi_{\pm}$ on $\partial_\pm D^2$, $c\coloneq \xi_-\cup\xi_+\in\undfld[1]{\C^{\Hom(\fX\to\fY)}}{\partial D^2}$, and any pure (2,2)-transfor $\eP\in\fld[2]{\eFun(\C^\fX{\to}\C^\fY)}{D^2}[\C^{c}]$,
\[
m^\eP_a \coloneq\eval_{D^2,\partial_{\pm}D^2}^\fY(\eP(a)) \in\fY(\eval(\xi_-(a))\Rightarrow\eval(\xi_+(a)))
\]
defines a modification $m^\eP\colon\eval(\xi_-)\Rightarrow\eval(\xi_+)$, and $[\C^{\langle m^\eP\rangle}]=[\eP]$ in $\eFun(\C^\fX{\to}\C^\fY)$.
\end{lemma}

\begin{pf}
As $m^\eP\coloneq\eval\circ\fX_\eP\circ\langle-\rangle$ is a whiskering of the modification $\fX_\eP$, $m^\eP$ is a modification. To see $[\C^{\langle m^\eP\rangle}]=[\eP]$, first note that for each $0$-field $a$,
$$\mathrm{eval}(\C^{\langle m^\eP\rangle}(a))=\eval(\langle m^\eP_a\rangle)= m^\eP_a=\mathrm{eval}(\eP(a)),$$
where the second equality holds because the evaluation of a positively oriented single-vertex string diagram is its label (where by positively oriented we mean the link parameterization on the vertex agrees with the standard splitting $\partial_{\pm}D^2$). Thus by \cref{eval-facts-2}\ref{fact:eval-ker2}, $[\C^{\langle m^\eP\rangle}(a)]=[\eP(a)]$ for all $a$, so by definition of $\lU{\eFun(\C^\fX{\to}\C^\fY)}$, $[\C^{\langle m^\eP\rangle}]=[\eP]$ in $\eHom(\C^\fX{\to}\C^\fY)$. 
\end{pf}

\subsubsection{The disk-like \texorpdfstring{$\eH\mathrm{om}$}{H} and traditional \texorpdfstring{$\Hom$}{Hom} equivalence}

The following is the $n=2$ analog of \cref{tsf}, which allowed us to transfer sphere traces along weak equivalences of disk-like 1-categories. The proof is exactly the same as in the $n=1$ version, except for surjectivity we need to use essential surjectivity on both 0- and 1-fields.

\begin{proposition}
\label{tsf2}
If $\eF \colon \eC \to \eD$ is a weak equivalence of disk-like 2-categories and $\psn^\C\colon \undfldl{\C}{S^2} \to \bC$ is a linear functional satisfying \ref{2psnP}, then the following hold.
\begin{lst} 
\item\label{tsf2a} $\widehat\eF_{S^2}$ is an isomorphism.
\item\label{tsf2b} $\psn^{\eD} \coloneq \psn^{\eC} \circ \widehat{\eF}_{S^2}^{-1} \colon \undfldl{\D}{S^2} \to \bC$ satisfies \ref{2psnP}.
\item\label{tsf2c} $(\eD, \psn^{\eD})$ is a
unitary disk-like 2-category, and is finite if and only if $\eC$ is.
\end{lst}
\end{proposition}

\begin{pf}
\itemstep{Proof of \ref{tsf2a}.} We first show $\widehat\eF_{S^2}$ is surjective. Let $\{\alpha_i\}$ be a decomposition representing a 2-field $\alpha$ in the colimit $\undfldl[2]{\D}{S^2}$. Let $\Gamma=(V,E,F)$ be the (unframed) string diagram stratification of the splitting (so in particular $\Gamma$ is transverse to $\alpha$). Write $\ell_e\coloneq\alpha|_e$ for edges $e\in E$ and $d_v\coloneq \alpha|_v$ for vertices $v\in V$, and write $v_0(e),v_1(e)$ for the (possibly equal) endpoints of $e$. Throughout, we illustrate the argument by considering a local region of the splitting $\{X_i\}$ of $S^2$ and the corresponding region for $\alpha=[\{\alpha_i\}]$ as illustrated below.
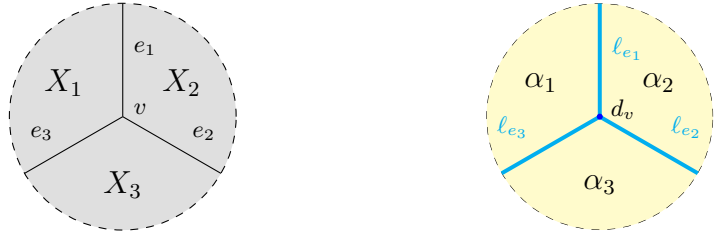
\def\ellblue{cyan}
\def\dvblue{blue}
\def\Fcvred{red}
\def\uvpurple{violet!50!white}
\def\mueorange{orange}
\def\webrown{brown}
\def\alphyellow{yellow!25!white}
\[
\begin{tkz}[scale=0.75]
    \fill[gray!25!white,draw=black,dashed] (0,0) circle (2);
    \draw(0,0)--node[right,pos=.6]{\scriptsize$e_1$}(90:2);
    \draw(0,0)--node[above right,pos=.6]{\scriptsize$e_2$}(-30:2);
    \draw(0,0)--node[above left,pos=.6]{\scriptsize$e_3$}(180+30:2);
    \node at (0,0)[right,yshift=1mm]{\scriptsize$v$};
    \node at (180-30:1.2) {$X_1$};
    \node at (30:1.2) {$X_2$};
    \node at (-90:1.2) {$X_3$};
\end{tkz}
\qquad\qquad\qquad\qquad
\begin{tkz}[scale=0.75]
    \clip (0,0)circle(2);
    \fill[\alphyellow,draw=black,dashed] (0,0) circle (2);
    \draw[\ellblue,line width=1.5pt](0,0)--node[right,pos=.6]{\scriptsize$\ell_{e_1}$}(90:2);
    \draw[\ellblue,line width=1.5pt](0,0)--node[above right,pos=.6]{\scriptsize$\ell_{e_2}$}(-30:2);
    \draw[\ellblue,line width=1.5pt](0,0)--node[above left,pos=.6]{\scriptsize$\ell_{e_3}$}(180+30:2);
    \node at (0,0)[right,yshift=1mm]{\scriptsize$d_v$};
    \fill[\dvblue] (0,0) circle(1.5pt);
    \node at (180-30:1.2) {$\alpha_1$};
    \node at (30:1.2) {$\alpha_2$};
    \node at (-90:1.2) {$\alpha_3$};
\end{tkz}
\]

We seek a 2-field $\gamma$ on $S^2$ in $\C$ such that $\widehat\eF_{S^2}(\gamma)=\alpha$. Since $\eF$ is unitarily essentially surjective on 0-fields, for each $v\in V$ there is a 0-field $c_v$ in $\C$ and a unitary equivalence $u_v\in\fldl[1]{\D}{I}[\fcj{\eF(c_v)}\amalg d_v]$. 
Thus
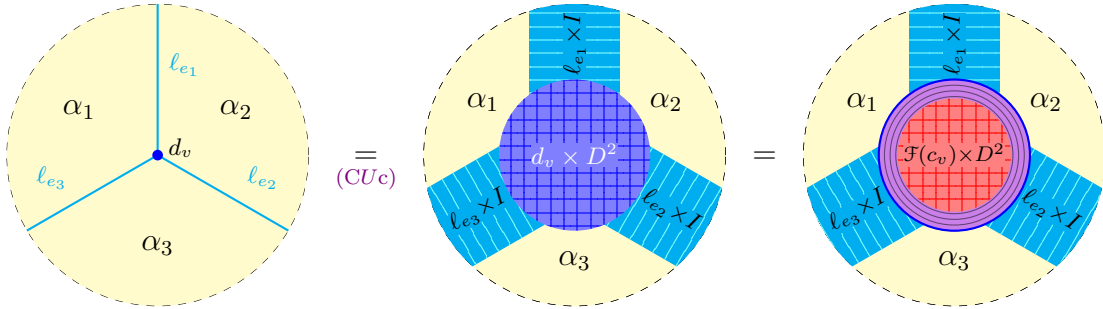
\[
\begin{tkz}
    \clip (0,0)circle(2);
    \fill[\alphyellow,draw=black,dashed] (0,0) circle (2);
    \draw[\ellblue,thick](0,0)--node[right,pos=.6]{\scriptsize$\ell_{e_1}$}(90:2);
    \draw[\ellblue,thick](0,0)--node[above right,pos=.6]{\scriptsize$\ell_{e_2}$}(-30:2);
    \draw[\ellblue,thick](0,0)--node[above left,pos=.6]{\scriptsize$\ell_{e_3}$}(180+30:2);
    \node at (0,0)[right,yshift=1mm]{\scriptsize$d_v$};
    \fill[\dvblue] (0,0) circle(2pt);
    \node at (180-30:1.2) {$\alpha_1$};
    \node at (30:1.2) {$\alpha_2$};
    \node at (-90:1.2) {$\alpha_3$};
\end{tkz}
\;\;\underset{\ref{CTc2}}{=}\;\;
\begin{tkz}
    \clip (0,0)circle(2);
    \fill[\alphyellow,draw=black,dashed] (0,0) circle (2);
    \draw[dashed] (0,0) circle (2);
    \fill[\ellblue, rotate=210, region={lines,color=\ellblue!50!white,angle=120}] (0,-0.6) rectangle (2,0.6);
    \fill[\ellblue, rotate=90, region={lines,color=\ellblue!50!white,angle=0}] (0,-0.6) rectangle (2,0.6);
    \fill[\ellblue, rotate=-30, region={lines,color=\ellblue!50!white,angle=60}] (0,-0.6) rectangle (2,0.6);
    \fill[\dvblue!50!white,region={grid,color=\dvblue,angle=0}](0,0) circle (1);
    \node at (180-30:1.4) {$\alpha_1$};
    \node at (30:1.4) {$\alpha_2$};
    \node at (-90:1.4) {$\alpha_3$};
    \path(0,0)--node[pos=.735,rotate=90,fill=\ellblue,inner sep=0pt]{\scriptsize$\ell_{e_1}{\times} I$}(90:2);
    \path(0,0)--node[pos=.735,rotate=-30,fill=\ellblue,inner sep=0pt]{\scriptsize$\ell_{e_2}{\times} I$}(-30:2);
    \path(0,0)--node[pos=.735,rotate=30,fill=\ellblue,inner sep=0pt]{\scriptsize$\ell_{e_3}{\times} I$}(180+30:2);
    \node[white,fill=\dvblue!50!white,inner sep=0pt] at (0,0){\scriptsize$d_v\times D^2$};
\end{tkz}
\;\;=\;\;
\begin{tkz}
    \clip (0,0)circle(2);
    \fill[\alphyellow,draw=black,dashed] (0,0) circle (2);
    \draw[dashed] (0,0) circle (2);
    \fill[\ellblue, rotate=210, region={lines,color=\ellblue!50!white,angle=120}] (0,-0.6) rectangle (2,0.6);
    \fill[\ellblue, rotate=90, region={lines,color=\ellblue!50!white,angle=0}] (0,-0.6) rectangle (2,0.6);
    \fill[\ellblue, rotate=-30, region={lines,color=\ellblue!50!white,angle=60}] (0,-0.6) rectangle (2,0.6);
    \fill[\uvpurple] (0,0) circle (1);
    \foreach \x in {1,0.925,...,0.7} {
        \draw[\uvpurple!50!black] (0,0) circle (\x);
    }
    \draw[\dvblue,thick] (0,0) circle (1);
    \fill[\Fcvred!50!white,region={grid,color=\Fcvred,angle=0}](0,0) circle (0.75);
    \node at (180-30:1.4) {$\alpha_1$};
    \node at (30:1.4) {$\alpha_2$};
    \node at (-90:1.4) {$\alpha_3$};
    \path(0,0)--node[pos=.735,rotate=90,fill=\ellblue,inner sep=0pt]{\scriptsize$\ell_{e_1}{\times} I$}(90:2);
    \path(0,0)--node[pos=.735,rotate=-30,fill=\ellblue,inner sep=0pt]{\scriptsize$\ell_{e_2}{\times} I$}(-30:2);
    \path(0,0)--node[pos=.735,rotate=30,fill=\ellblue,inner sep=0pt]{\scriptsize$\ell_{e_3}{\times} I$}(180+30:2);
    
    \node[fill=\Fcvred!50!white,inner sep=0pt] at (0,0) {\tiny$\eF(c_v){\times} D^2$};
\end{tkz}
\]
where the purple region is $u_v\times S^1$. Observe that for each $e\in E$, the 1-field $\ell'_e\coloneq u_{v_0(e)}\blt \ell_e\blt u^\star_{v_1(e)}$ has boundary in the image of $\eF$.
Thus, since $\eF$ is unitarily essentially surjective on 1-fields, there is a 1-field $\mu_e\in\fldl[1]{\C}{I}[\fcj{c_{v_0(e)}}\amalg c_{v_1(e)}]$ and a unitary equivalence $w_e\in\fldl[2]{\D}{I\times I}[\fcj{\eF(\mu_e)}\amalg\ell'_e]$, so the above decomposition is given by the following.
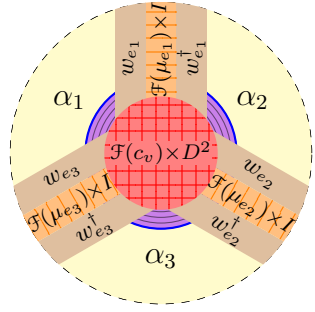
\[
\begin{tkz}
    \clip (0,0)circle(2);
    \fill[\alphyellow,draw=black,dashed] (0,0) circle (2);
    \draw[dashed] (0,0) circle (2);
    \fill[\uvpurple] (0,0) circle (1);
    \foreach \x in {1,0.925,...,0.7} {
        \draw[\uvpurple!50!black] (0,0) circle (\x);
    }
    \draw[\dvblue,thick] (0,0) circle (1);
    \fill[\webrown!50!white, rotate=210] (0,-0.6) rectangle (2,-0.2);
    \fill[\mueorange!50!white, rotate=210, region={lines,color=\mueorange,angle=120}] (0,-0.2) rectangle (2,0.2);
    \fill[\webrown!50!white, rotate=210] (0,0.2) rectangle (2,0.6);

    \fill[\webrown!50!white, rotate=90] (0,-0.6) rectangle (2,-0.2);
    \fill[\mueorange!50!white, rotate=90, region={lines,color=\mueorange,angle=0}] (0,-0.2) rectangle (2,0.2);
    \fill[\webrown!50!white, rotate=90] (0,0.2) rectangle (2,0.6);

    \fill[\webrown!50!white, rotate=-30] (0,-0.6) rectangle (2,-0.2);
    \fill[\mueorange!50!white, rotate=-30, region={lines,color=\mueorange,angle=60}] (0,-0.2) rectangle (2,0.2);
    \fill[\webrown!50!white, rotate=-30] (0,0.2) rectangle (2,0.6);
    
    \fill[\Fcvred!50!white,region={grid,color=\Fcvred,angle=0}](0,0) circle (0.75);

    \node at (180-30:1.4) {$\alpha_1$};
    \node at (30:1.4) {$\alpha_2$};
    \node at (-90:1.4) {$\alpha_3$};
    
    \path(0,0)--node[pos=.685,rotate=90,inner sep=0pt]{\tiny$\eF({\mu_{e_1}}){\times} I$}(90:2);
    \path(0,0)--node[pos=.685,rotate=-30,inner sep=0pt]{\tiny$\eF(\mu_{e_2}){\times} I$}(-30:2);
    \path(0,0)--node[pos=.685,rotate=30,inner sep=0pt]{\tiny$\eF(\mu_{e_3}){\times} I$}(180+30:2);
    
    \path[shift={(90-90:0.4cm)}](0,0)--node[pos=.635,rotate=90,fill=\webrown!50!white,inner sep=0pt]{\scriptsize$w_{e_1}^\star$}(90:2);
    \path[shift={(-30-90:0.4cm)}](0,0)--node[pos=.635,rotate=-30,fill=\webrown!50!white,inner sep=0pt]{\scriptsize$w_{e_2}^\star$}(-30:2);
    \path[shift={(180+30-90:0.4cm)}](0,0)--node[pos=.635,rotate=30,fill=\webrown!50!white,inner sep=0pt]{\scriptsize$w_{e_3}$}(180+30:2);
    
    \path[shift={(90+90:0.4cm)}](0,0)--node[pos=.635,rotate=90,fill=\webrown!50!white,inner sep=0pt]{\scriptsize$w_{e_1}$}(90:2);
    \path[shift={(-30+90:0.4cm)}](0,0)--node[pos=.635,rotate=-30,fill=\webrown!50!white,inner sep=0pt]{\scriptsize$w_{e_2}$}(-30:2);
    \path[shift={(180+30+90:0.4cm)}](0,0)--node[pos=.635,rotate=30,fill=\webrown!50!white,inner sep=0pt]{\scriptsize$w_{e_3}^\star$}(180+30:2);
    \node[fill=\Fcvred!50!white,inner sep=0pt] at (0,0) {\tiny$\eF(c_v){\times} D^2$};
\end{tkz}
\]
Now let $\tld\alpha_i$ denote $\alpha_i$ glued to all adjacent $w_e$, $w_e^{\star}$, and $u_v\times I$, with each remaining product field $\eF(\mu_e)\times I$ and each $\eF(c_v)\times D^2$ collapsed via a collar map onto its corresponding edge and vertex, and let $\tld X_i$ denote the ball on which $\tld\alpha_i$ lives. Thus by extended-isotopy invariance the above decomposition is the following.
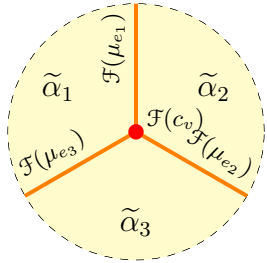
\[
\begin{tkz}[scale=0.85]
    \clip (0,0)circle(2);
    \fill[\alphyellow,draw=black,dashed] (0,0) circle (2);
    \clip (0,0) circle (2);
    \draw[dashed] (0,0) circle (2);
    \draw[\mueorange, line width=1.5pt] (0,0) -- (90:2);
    \draw[\mueorange, line width=1.5pt] (0,0) -- (-30:2);
    \draw[\mueorange, line width=1.5pt] (0,0) -- (180+30:2);
    \fill[\Fcvred](0,0) circle (0.12);
    \node at (180-30:1.4) {$\tld{\alpha}_1$};
    \node at (30:1.4) {$\tld{\alpha}_2$};
    \node at (-90:1.4) {$\tld{\alpha}_3$};
    \path(0,0)--node[pos=.65,rotate=90,above]{\scriptsize$\eF(\mu_{e_1})$}(90:2);
    \path(0,0)--node[pos=.65,rotate=-30,above]{\scriptsize$\eF(\mu_{e_2})$}(-30:2);
    \path(0,0)--node[pos=.65,rotate=30,above]{\scriptsize$\eF(\mu_{e_3})$}(180+30:2);
    \node[right,yshift=1.5mm] at (0,0) {\scriptsize$\eF(c_v)$};
\end{tkz}
\]
Now $\partial\tld\alpha_i$ is in the image of $\eF$, so since $\widehat\eF$ is a linear isomorphism on balls, there is a unique 2-field $\gamma_i$ on $\tld X_i$ in $\C$ such that $\widehat\eF(\gamma_i)=\tld\alpha_i$. The above argument/procedure therefore shows that $[\{\widehat\eF_{\tld X_i}(\gamma_i)\}]=[\{\tld\alpha_i\}]=\alpha$ where the second equality holds because the passage from $\{\alpha_i\}$ to $\{\tld\alpha_i\}$ used only \ref{CTc2}, the unitarity identities for the $u_v$ and $w_e$, and gluing, each of which preserves the image of the splitting in the colimit. Thus $\gamma\coloneq[\{\gamma_i\}]$ has $\widehat\eF_{S^2}(\gamma)=\alpha$, so $\widehat\eF_{S^2}$ is surjective.

Injectivity of $\widehat\eF_{S^2}$, as well as the proofs of \ref{tsf2b} and \ref{tsf2c}, now follow by the same arguments as they did for $n=1$.
\end{pf}

\begin{remark}
Notice that the above argument goes through essentially unchanged to show that $\widehat\eF_{X,c}$ is a linear isomorphism for all 2-manifolds $X$ and all $c\in\undfldl{\C}{\partial X}$.
\end{remark}

\begin{remark}
From the proof of \cref{tsf} and \cref{tsf2}, the general-$n$ argument that a disk-like weak equivalence $\eF\colon\C\to\D$ gives linear isomorphisms $\widehat\eF_{X,c}\colon\fldl{\C}{X}[c]\to\fldl{\D}{X}[\eF(c)]$ follows from the evident inductive procedure.
\end{remark}

\begin{theorem}
  \label{DL-thing-ii}
  For proto-3-Hilbert spaces $\fX$ and $\fY$, the constructions $F \smash{\overset{\eqref{DL 2-functor from 2-functor-1}}\longmapsto}\C^F$, $\xi\overset{\eqref{DL 2-functor from 2-functor-2}}\longmapsto\C^\xi$, and $\alpha\overset{\eqref{DL 2-functor from 2-functor-3}}\longmapsto\C^\alpha$ define a weak equivalence
  \[
    \C^{(-)}\colon \C^{{\Hom}(\fX\to\fY)} \overset\sim\longrightarrow\eFun(\C^\fX{\to}\C^\fY).
  \]
  Thus, when $\fX$ and $\fY$ are moreover finite, so that $\Hom(\fX\to\fY)$ carries the spherical weight of \cref{cor:proto-hom-weight}, $\eFun(\C^\fX{\to}\C^\fY)$ is unitary with its induced sphere trace via \cref{tsf2}\ref{tsf2c}.
\end{theorem}

\begin{pf}
Aside from showing \ref{axiom:eF-local-relations}, we leave the straightforward verification that $\C^{(-)}$ is a disk-like functor to the reader.

\itemstep{\ref{axiom:eF-local-relations}.} For $u\in\bC\{\fld[2]{\C^{{\Hom}(\fX\to\fY)}}{X}[c]\}$ with $2$-ball $X$ and $c\in\undfld[1]{\C^{{\Hom}(\fX\to\fY)}}{\partial X}$,
\begin{align*}
u\in\fldlU{\C^{{\Hom}(\fX\to\fY)}}{X}[c]
&\iff\eval^{\Hom(\fX\to\fY)}(u)=0\\
&\iff(\eval^{\Hom(\fX\to\fY)}(u))_a=0
&\forall a\in\fX&\\
&\iff\eval^\fY(u(a))=0
&\forall a\in\fX&\\  
&\iff u(a)\in\fldlU{\C^\fY}{X}[c(a)]
&\forall a\in\fX&\\
&\iff
\C^u\in\fldlU{\eFun(\C^\fX{\to}\C^\fY)}{X}[\C^{c}].
\end{align*}

\itemstep{$\C^{(-)}$ is unitarily essentially surjective on $0$-fields.}
For a disk-like functor $\eF\colon\C^\fX\to\C^\fY$, by \cref{lem:ess-surj-functors} the UAF-preserving dagger functor $F^\eF\in\Hom(\fX\to\fY)$ satisfies $\C^{F^\eF}\cong^\star\eF$.

\itemstep{$\C^{(-)}$ is unitarily essentially surjective on $1$-fields.} 
Let $F,G\in\Hom(\fX\to\fY)$ and let $\eN\in\fld[1]{\eFun(\C^\fX{\to}\C^\fY)}{Y}[\fcj{\C^F}\amalg\C^G]$ be a $(2,1)$-transfor on a 1-ball $Y$. Choose a homeomorphism $\varphi\colon D^1\to Y$, and let $\eta^{\varphi^{-1}_*\eN}\colon F\Rightarrow G$ be the natural transformation produced by \cref{21lemma} applied to the $D^1$-shaped $(2,1)$-transfor $\varphi^{-1}_*\eN$, so that $\C^{\langle\eta^{\varphi^{-1}_*\eN}\rangle}\cong^\star\varphi^{-1}_*\eN$.
Setting $\xi\coloneq\varphi_*\langle\eta^{\varphi^{-1}_*\eN}\rangle$, the homeomorphism $\varphi^{-1}\colon Y\to D^1$ sends $\xi\mapsto\langle\eta^{\varphi^{-1}_*\eN}\rangle$ and $\eN\mapsto\varphi^{-1}_*\eN$, so $\C^\xi\cong^\star\eN$. Thus $\C^{(-)}$ is unitarily essentially surjective on 1-fields.

\itemstep{$\C^{(-)}$ is a linear isomorphism on $2$-fields.} 
Faithfulness follows by following the chain of equivalences in the above proof of \ref{axiom:eF-local-relations} in reverse, which gives that $\C^u=0$ implies $u=0$. To see fullness, note that if $\eP\in\fldl[2]{\eHom(\C^\fX{\to}\C^\fY)}{X}[\fcj{\C^\xi}\amalg\C^\eta]$, then for any choice of homeomorphism $\varphi\colon D^2\to X$ we have $\varphi^{-1}_*[\eP]=[\C^{\langle m^{\varphi^{-1}_*\eP}\rangle}]$ by \cref{lem:surj-2-fields}. Applying $\varphi_*$, we conclude $[\eP]=[\C^{\varphi_*\langle m\rangle}]$.

\itemstep{$\C^{(-)}$ is isometric.} This holds by definition of the sphere trace on $\eFun(\C^\fX{\to}\C^\fY)$: $\psn^{\eFun}\coloneq\psn^{\C^{\Hom(\fX\to\fY)}}$ via \cref{tsf2}\ref{tsf2b}.
\end{pf}

\begin{corollary}
\label{dag-thing-ii-equiv}
\label{XFun-2}
For finite proto-3-Hilbert spaces $\fX$ and $\fY$, there is a canonical isometric equivalence of finite proto-3-Hilbert spaces
\[
    \fX_{\eFun(\C^\fX\to\C^\fY)} \overset\sim\longrightarrow \Hom(\fX\to\fY).
\]
\end{corollary}

\begin{pf}
    Composing the isometric equivalence $\fX_{\C^{(-)}}$ with the inverse of $\eval$ gives the claim:
    \[
    \fX_{\eFun(\C^\fX\to\C^\fY)}
    \xleftarrow[\substack{\fX_{\C^{(-)}}\\\eqref{DL-thing-ii}\,\&\,\eqref{CcongD implies XCcongXD-n2}}]{\cong^\dag}\fX_{\C^{\Hom(\fX\to\fY)}}
    \xrightarrow[\substack{\eval\\\eqref{XCXcongX-n-2}}]{\cong^\dag}
    \Hom(\fX\to\fY).\qedhere
    \]
\end{pf}

\subsection{The \texorpdfstring{$(2+1)$D}{(2+1)D} Unitary Cobordism Hypothesis}

Let $\C$ be a finite unitary disk-like 2-category, so that $\fX_\C$ is a finite proto-3-Hilbert space by \cref{cstr:proto-3-Hilbert spaces from unitary disk-like 2-categories} and \cref{unitary implies finite unitary cylinder 2-categories}. A \defn{unitary representation} of $\C$ is a disk-like functor $\eR\colon \C^{\fX_\C^{1\op}}\to 2\eHilb$ where $2\eHilb\coloneq\C^{2\Hilb}$. Unitary representations form a disk-like 2-category $\eRep(\C)\coloneq\eFun(\C^{\fX_\C^{1\op}}{\to}2\eHilb)$; see \cref{def: functor disk-like n-category}. By \cref{DL-thing-ii} and \cref{tsf2}, the induced spherical weight on $\Rep(\fX_\C^{1\op})=\Hom(\fX_\C^{1\op}\to2\Hilb)$
from \cref{prop:rep-is-3hilb}
equips $\eRep(\C)$ with a sphere trace making it finite unitary. 

\begin{remark}
\label{rem:eRep-strings-2}
As was the case for $n=1$ (see \cref{rem:eRep-strings}), for any finite proto-3-Hilbert space $\fX$ there are canonical isometric weak equivalences
\[
\eRep(\C^\fX)
\underset{\eqref{DL-thing-ii}}{\cong^\dag}
\C^{\Hom(\fX_{\C^\fX}^{1\op}\to2\Hilb)}
\underset{\eqref{2lem:C-preserves-equivs}\&\eqref{XCXcongX-n-2}}{\cong^\dag}
\C^{\Hom(\fX^{1\op}\to2\Hilb)}
\underset{\eqref{DL-thing-ii}}{\cong^\dag}
\eFun(\C^{\fX^{1\op}}{\to}2\eHilb).
\]
\end{remark}

\begin{definition}
\label{def:completeDL2cat}
    A finite unitary disk-like 2-category $\C$ is called \defn{complete} if there is an isometric weak equivalence $\eRep(\C)\cong^\dag\C$.
\end{definition}

\THREEDUCH*

\begin{pf}
The assertions that $(\fX_{\C^\fX},\vee,\Psi^{\fX_{\C^\fX}})\cong^\dag(\fX,\vee,\Psi)$ and $(\C^{\fX_\C},\psn^{\C^{\fX_\C}})\cong^\dag(\C,\psn)$ are \cref{XCXcongX-n-2,CXCcongC-n2} respectively, so it remains to show these equivalences preserve completeness in that (i) $\fX_{\eRep(\C)}$ is complete for all finite unitary disk-like 2-categories $\C$ and (ii) $\C^\fX$ is complete for all finite 3-Hilbert spaces $\fX$. 

For (i), observe that 
\[
\fX_{\eRep(\C)}\defeq{}\fX_{\eFun(\C^{\fX_\C^{1\op}}\to\C^{2\Hilb})} \;\underset{\mathclap{\eqref{XFun-2}}}{\cong^\dag}\; \Hom(\fX_\C^{1\op}\to2\Hilb)\defeq{}\Rep(\fX_\C^{1\op}).
\]
As $\fX_\C$ is a finite pivotal $\Cstar$-2-category by \cref{unitary implies finite unitary cylinder 2-categories}, it follows from \cref{prop:rep-is-3hilb} that $\Rep(\fX_\C^{1\op})$ is a 3-Hilbert space, whence $\fX_{\eRep(\C)}$ is a finite 3-Hilbert space, i.e., is complete.

For (ii), we have
\[
\eRep(\C^\fX)\underset{\eqref{rem:eRep-strings-2}}{\cong^\dag}\eFun(\C^{\fX^{1\op}}{\to}2\eHilb)\underset{\eqref{DL-thing-ii}}{\cong^\dag}\C^{\Rep(\fX^{1\op})}\underset{\eqref{3hilb equiv to its reps}\,\&\,\eqref{2lem:C-preserves-equivs}}{\cong^\dag}\C^\fX
\]
where the third equivalence applies $\C^{(-)}$ to the isometric equivalence $\!\yo^{\!-1}\colon\Rep(\fX^{1\op})\xrightarrow{\sim}\fX$ from \cref{3hilb equiv to its reps}. Thus $\C^\fX$ is complete.
\end{pf}

\subsection{Application: Classification of (oriented) (2+1)D unitary TQFTs as \mbox{(\texorpdfstring{$\Hstar$}{H*}-)Morita} classes of \texorpdfstring{$\Hstar$}{H*}-multifusion categories}
\label{classificationHstarmFC}

In this subsection, we categorify the results of \cref{classificationHstaralg} to classify fully extended oriented 3D unitary TQFTs by \defn{$\Hstar$-/isometric Morita equivalence classes} of $\Hstar$-multifusion categories, a notion defined below.

Recall from \cite[Defn. 3.1]{3Hilb} that an \defn{$\Hstar$-multifusion category} is a unitary multifusion category $\cA$ equipped with a choice of UAF $\vee$ on $\rB\cA$ and a spherical weight on $\rB\cA$. By \cite[Rmk. 3.2]{3Hilb}, every $\Hstar$-multifusion category $\cA$ has a canonical underlying 2-Hilbert space $L^2\cA$ with trace $\Tr^{L^2\cA} \coloneq \psi^\cA\circ\tr^\vee_L=\psi^\cA\circ\tr^\vee_R$ where $\tr_L^\vee$ and $\tr^\vee_R$ are the left and right traces $\bigoplus_{a\in \cA}\End(a)\to\End_\cA(1_\cA)$ respectively.

Let $\cB$ be an $\Hstar$-multifusion category and write $1_\cB = \bigoplus_{i=1}^r 1_i$ for the decomposition of $1_\cB$ into simple objects. A right $\cB$-module $\cM_\cB$ is \defn{faithful} if $\cM \triangleleft 1_i \neq 0$ for all $1 \le i \le r$. 

For a unitary monoidal category $\cA$, a (right) \defn{$\cA$-module category} is a finite semisimple unitary category $\cM$ equipped with a bilinear functor $\triangleleft\colon \cM\times\cA\to\cM$ and unitary natural coherence isomorphisms $(-\triangleleft -)\triangleleft - \cong -\triangleleft(-\otimes -)$ and $-\triangleleft 1_\cA \cong -$. If $\cA$ is an $\Hstar$-multifusion category, a (right) \defn{$\cA$-module} is a right $\cA$-module category further equipped with a module trace in the sense of \cite[Defn. 3.4]{3Hilb}.
By \cite[Cor. 3.13]{3Hilb}, for an $\Hstar$-multifusion category $\cA$ and a right $\cA$-module $\cM_\cA$, there always exists a faithful $\cA$-module trace $\Tr^{\cM_\cA}$, and if $\cA$ and $\cM$ are both indecomposable then $\Tr^{\cM_\cA}$ is unique up to a positive scalar. We denote by $\Mod(\cA)$ the 3-Hilbert space of unitary right $\cA$-modules equipped with $\cA$-module traces \cite[Defn. 3.6]{3Hilb}, and we equip $\Mod(\cA)$ with the UAF given by unitary adjunction and the spherical weight $\Psi^{\Mod(\cA)}$ of \cite[Eqn. (7)]{bases}. By \cite[Cor. 4.52]{3Hilb} every 3-Hilbert space is isometrically equivalent to $\Mod(\cA)$ for some $\Hstar$-multifusion category $\cA$.

\begin{definition}[Commutant of an $\Hstar$-multifusion category]
\label{def:commutant}
Let $\cB$ be an $\Hstar$-multifusion category and let $\cM_\cB$ be a (right) $\cB$-module. The \defn{commutant} of $\cB$ for $\cM_\cB$ is the $\Hstar$-multifusion category $\cB' \coloneq \End(\cM_\cB)$, with UDF given by unitary adjunction and spherical weight $\psi^{\cB'}$ given by $\Psi^{\Mod(\cB)}_{\cM_\cB}$, where $\Psi^{\Mod(\cB)}$ is the renormalized spherical weight \cite[(7)]{bases}.
\end{definition}

\begin{lemma}
\label{lem:daggerstar}
Suppose $\cA$ and $\cB$ are $\Hstar$-multifusion categories and $\alpha\colon \cA \to \cB$ is a UAF-preserving dagger monoidal equivalence with unitor $u_\alpha\colon 1_\cB \to \alpha(1_\cA)$. Then for all $a \in \cA$ and $f \in \End_\cA(a)$,
\begin{equation}
\label{eq:dagger-star}
\tr^{\vee}(\alpha(f))=u_\alpha^\dag \circ \alpha(\tr^{\vee}(f)) \circ u_\alpha.
\end{equation}
\end{lemma}
\begin{pf}
Since $\rB\alpha$ is UAF-preserving, the canonical isomorphisms $\delta_a\colon\alpha(a^\vee)\to\alpha(a)^\vee$ are unitary. We compute that
\[
\begin{aligned}
u_\alpha^\dag\alpha(\tr^\vee(f))u_\alpha&=\alpha(\ev_a) (\alpha^2_{a^\vee,a})^{-1} \alpha^2_{a^\vee,a}\alpha(\id_{a^\vee}\otimes f)(\alpha^2_{a^\vee,a})^{-1}\alpha^2_{a^\vee,a} (\alpha^2_{a^\vee,a})^{-1}\alpha(\ev_a^\dag)
\\&
=\ev_{\alpha(a)}(\delta_a\otimes\id_{\alpha(a)}) \alpha(\id_{a^\vee}\otimes f) (\delta_a^\dag\otimes\id_{\alpha(a)})\ev_{\alpha(a)}^\dag
\\&
=\ev_{\alpha(a)}(\delta_a\otimes\id_{\alpha(a)})(\alpha^2_{a^\vee,a})^{-1} \alpha^2_{a^\vee,a}(\id_{\alpha(a^\vee)}\otimes\alpha(f))(\delta_a^\dag\otimes\id_{\alpha(a)})\ev_{\alpha(a)}^\dag
\\&
=\ev_{\alpha(a)}(\delta_a\otimes\id_{\alpha(a)})(\id_{\alpha(a^\vee)}\otimes\alpha(f))(\delta_a^\dag\otimes\id_{\alpha(a)})\ev_{\alpha(a)}^\dag
\\&
=\ev_{\alpha(a)}((\delta_a\delta_a^\dag)\otimes\alpha(f))\ev_{\alpha(a)}^\dag
\\&
=\tr^\vee(\alpha(f))
\end{aligned}
\]
where the first and second equalities are by unitarity and naturality of the tensorators $\alpha^2_{a^\vee,a}$ respectively, and the third equality is by unitarity of the tensorators $\alpha^2_{a^\vee,a}$ and of $\delta_a$.
\end{pf}

Recall from \cite[Rmk. 3.30]{bases} that each $\Hstar$-multifusion category $\cA$ comes with a canonical \defn{$\Hilb$-valued trace} $\mathrm{TR}^{\cA}\colon \cA \to \Hilb$ given by the corepresentable functor $\cA(1_\cA \to -)$, which at objects $a\in\cA$ is given by the Hilbert space
\[
\mathrm{TR}^{\cA}(a) \coloneq \cA(1_\cA \to a)
\]
with the inner product
\begin{equation}
\label{eq:tr-inner-product}
\langle f |g \rangle_{\mathrm{TR}^{\cA}(a)} \coloneq \psi^{\cA}(\tr^\vee(f^\dag \circ g)) = \Tr^{L^2\cA}_{1_\cA}(f^\dag \circ g).
\end{equation}

\begin{lemma}
\label{lem:isometric-tfae}
For $\Hstar$-multifusion categories $\cA$ and $\cB$ and a dagger monoidal equivalence $\alpha\colon \cA \to \cB$ with (unitary) unitor $u_\alpha\colon 1_\cB \to \alpha(1_\cA)$, the following are equivalent.
\begin{lst}
\item \label{item:phi-unitary} The canonical natural isomorphism $\phi^\alpha\colon \mathrm{TR}^{\cA} \Rightarrow \mathrm{TR}^{\cB}\circ \alpha$, given at $a \in \cA$ by
\[
\phi_a^\alpha\colon \mathrm{TR}^{\cA}(a) \to \mathrm{TR}^{\cB}(\alpha(a)), \qquad f \mapsto \alpha(f)\circ u_\alpha,
\]
is unitary.
\item \label{item:psi-invariance} $\psi^{\cA}(\eta) = \psi^{\cB}(u_\alpha^\dag \circ \alpha(\eta) \circ u_\alpha)$ for all $\eta \in \End_{\cA}(1_\cA)$.
\item \label{item:trace-preserved} $\Tr^{L^2\cA} = \Tr^{L^2\cB}\circ \alpha$, i.e., for all $a \in \cA$ and $f \in \End_{\cA}(a)$, $\Tr^{L^2\cA}_a(f) = \Tr^{L^2\cB}_{\alpha(a)}(\alpha(f))$.
\end{lst}
\end{lemma}

\begin{definition} 
\label{def:HmFC-isometric-equivalence}
We will say that a dagger monoidal equivalence $\alpha\colon \cA\to\cB$ of $\Hstar$-multifusion categories $\cA$ and $\cB$ is \defn{isometric} if it satisfies any (hence all) of the conditions in \cref{lem:isometric-tfae}.
\end{definition}

\begin{pf}[Proof of {\cref{lem:isometric-tfae}}]~

\itemstep{\ref{item:phi-unitary}$\Rightarrow$\ref{item:psi-invariance}.}
Since $\phi_a^\alpha$ is unitary, for all $f,g,a$ we have $\langle \phi_a^\alpha(f) |\phi_a^\alpha(g)\rangle_{\mathrm{TR}^{\cB}(\alpha(a))} = \langle f |g \rangle_{\mathrm{TR}^{\cA}(a)}$, i.e., $\psi^{\cB}(u_\alpha^\dag \alpha(f)^\dag \alpha(g) u_\alpha) = \psi^{\cA}(f^\dag g)$, so \ref{item:psi-invariance} follows by taking $a = 1_\cA$, $f = \id_{1_\cA}$, and $g = \eta$.

\itemstep{\ref{item:trace-preserved}$\Rightarrow$\ref{item:phi-unitary}.}
Let $a \in \cA$ and $g,h \in \mathrm{TR}^{\cA}(a) = \cA(1_\cA \to a)$, so that $g^\dag h \in \End_\cA(1_\cA)$. Note that $\phi_a^\alpha$ is a linear isomorphism since $\alpha$ is an equivalence, and
\begin{multline*}
\langle g |h \rangle_{\mathrm{TR}^{\cA}(a)} \underset{\eqref{eq:tr-inner-product}}{=} \Tr^{L^2\cA}_{1_\cA}(g^\dag h) \underset{\ref{item:trace-preserved}}{=} \Tr^{L^2\cB}_{\alpha(1_\cA)}(\alpha(g)^\dag \circ \alpha(h)) 
\\
\underset{\ref{dagTr1}}{=} \Tr^{L^2\cB}_{1_\cB}(u_\alpha^\dag \circ \alpha(g)^\dag \circ \alpha(h) \circ u_\alpha) \underset{\eqref{eq:tr-inner-product}}{=} \langle \phi_a^\alpha(g) |\phi_a^\alpha(h) \rangle_{\mathrm{TR}^{\cB}(\alpha(a))},
\end{multline*}
so $\phi_a^\alpha$ is unitary. 

\itemstep{\ref{item:psi-invariance}$\Rightarrow$\ref{item:trace-preserved}.}
First note that by \ref{item:psi-invariance}, the delooping $\rB\alpha\colon \rB\cA \to \rB\cB$ is UAF-preserving by \cite[Prop. 4.9]{3Hilb}, and is thus an isometric equivalence of pre-3-Hilbert spaces. Thus, by applying \ref{item:psi-invariance} with $\eta = \tr^{\vee}(f)$ we obtain
\[
\Tr_a^{L^2\cA}(f) \defeq{} \psi^{\cA}(\tr^\vee(f)) \underset{\ref{item:psi-invariance}}{=} \psi^{\cB}(u_\alpha^\dag \alpha(\tr^\vee(f)) u_\alpha) \underset{\eqref{lem:daggerstar}}{=} \psi^{\cB}(\tr^\vee(\alpha(f))) \defeq{} \Tr_{\alpha(a)}^{L^2\cB}(\alpha(f)),
\]
which is \ref{item:trace-preserved}. 
\end{pf}

\begin{corollary}
\label{cor:isometric-uaf}
If $\alpha\colon \cA \to \cB$ is an isometric equivalence of $\Hstar$-multifusion categories, then it is automatically UAF-preserving, i.e., the delooping $\rB\alpha\colon \rB\cA \to \rB\cB$ is UAF-preserving and thus an isometric equivalence of pre-3-Hilbert spaces.
\end{corollary}

\begin{definition}
\label{def:hstar-morita}
We will say that $\Hstar$-multifusion categories $\cA$ and $\cB$ are \defn{$\Hstar$-Morita equivalent} if there is a faithful (right) $\cB$-module $\cM_\cB$ and an isometric equivalence of $\Hstar$-multifusion categories $\alpha\colon \cA \to \End(\cM_\cB)$ in the sense of \cref{def:HmFC-isometric-equivalence}.
\end{definition}

In what follows, we adopt the notation $\Omega_{\cM_\cB} \coloneq \End(\cM_\cB)$.

\begin{proposition}
\label{prop:mod-iso-iff-morita}
For $\Hstar$-multifusion categories $\cA$ and $\cB$, $\Mod(\cA) \cong^\dag \Mod(\cB)$ isometrically as 3-Hilbert spaces if and only if $\cA$ and $\cB$ are $\Hstar$-Morita equivalent.
\end{proposition}

\begin{pf}
\itemstep{($\Leftarrow$)} Let $\alpha\colon \cA \to \Omega_{\cM_\cB}$ be an isometric equivalence of $\Hstar$-multifusion categories. By \cref{cor:isometric-uaf}, $\rB\alpha\colon \rB\cA \to \rB\Omega_{\cM_\cB}$ is an isometric equivalence of pre-3-Hilbert spaces, so by the universal property of completion we get an isometric equivalence $(\rB\alpha)^\cent \colon (\rB\cA)^\cent \xrightarrow{\ \cong^\dag\ } (\rB\Omega_{\cM_\cB})^\cent$. Since $\Mod(\cB) = \bigboxplus_j \Mod(\cB_j)$ for the indecomposable summands $\cB_j$ of $\cB$, we may assume without loss of generality that $\cB$ is indecomposable. Since now $\fX \cong^\dag (\rB\Omega_x)^\cent$ for any object $x$ in a connected 3-Hilbert space $\fX$ by \cite[Exm. 3.29]{bases}, there is an isometric equivalence $(\rB\Omega_{\cM_\cB})^\cent \cong^\dag \Mod(\cB)$, and $(\rB\cA)^\cent \cong^\dag \Mod(\cA)$ by \cite[Exm. 3.28]{bases}. Composing with $(\rB\alpha)^\cent$ above, we obtain $\Mod(\cA) \cong^\dag \Mod(\cB)$.

\itemstep{($\Rightarrow$)} Let $F\colon \Mod(\cA) \to \Mod(\cB)$ be an isometric equivalence. Since $L^2\cA_\cA$ generates $\Mod(\cA)$, its image $\cM_\cB \coloneq F(L^2\cA_\cA)$ generates $\Mod(\cB)$ and is thus faithful. As $F$ is an isometric equivalence of 3-Hilbert spaces, it restricts to an isometric equivalence of 2-Hilbert spaces $\alpha_F'\colon \Omega_{L^2\cA_\cA} \to \Omega_{\cM_\cB}$, which is monoidal with unitary tensorators and unitors inherited from $F$. As $\Omega_{L^2\cA_\cA} = \End(L^2\cA_\cA)$ has UDF given by unitary adjunction and the commutant spherical weight $\psi^{\cA'}$, this is an isometric equivalence of $\Hstar$-multifusion categories.

Assuming without loss of generality that $\cA$ is indecomposable (as for $\cB$ above), the ``left multiplication'' $L\colon \cA \to \Omega_{L^2\cA_\cA}$ given by $a \mapsto a \otimes -$ is an isometric equivalence of 2-Hilbert spaces by \cite[Prop. 3.9]{3Hilb}, and thus is an isometric equivalence of $\Hstar$-multifusion categories by \cref{lem:isometric-tfae}\ref{item:trace-preserved}. Thus $\alpha_F \coloneq \alpha_F' \circ L\colon \cA \to \Omega_{\cM_\cB}$ is an isometric equivalence of $\Hstar$-multifusion categories.
\end{pf}

\begin{proposition}
\label{n2nonfailure}
For indecomposable $\Hstar$-multifusion categories $\cA$ and $\cB$, if $F\colon \Mod(\cA)\to\Mod(\cB)$ is an ordinary/algebraic equivalence in $\mathsf{H^*mFC}$, then there is some $\lambda\in\bR_{>0}$ for which the rescaled functor $F_\lambda\colon(\cM,\Tr^\cM)\mapsto(F(\cM),\lambda^{1/2}\Tr^{F(\cM)})$ is an isometric equivalence of 3-Hilbert spaces $\Mod(\cA)\to\Mod(\cB)$. Thus, applying this to each indecomposable summand separately (with one scalar for each), every ordinary/algebraic Morita equivalence can be augmented to an isometric/$\Hstar$-Morita equivalence.
\end{proposition}
\begin{pf}
Since $\Mod (\cA)$ and $\Mod(\cB)$ are connected 3-Hilbert spaces, by \cite[Rmk. 4.50]{3Hilb} their spherical weights are unique up to a positive scalar. Thus there is some $\lambda\in\bR_{>0}$ for which $\Psi^{\Mod(\cA)} = \lambda(\Psi^{\Mod(\cB)}\circ F)$. Define $F_\lambda\colon \Mod(\cA)\to\Mod(\cB)$ on objects by 
\[
F_\lambda(\cM,\Tr^\cM) \coloneq (F(\cM), \lambda^{1/2}\Tr^{F(\cM)})
\]
and on 1- and 2-morphisms by $F$. Then $F_\lambda$ is a dagger 2-equivalence since $F$ is, and $F_\lambda$ is isometric: indeed, $\Tr^{F_\lambda(\cM)}_m=\lambda^{1/2}\Tr^{F(\cM)}_m$ and
\[
d_m^{F_\lambda(\cM)}=\Tr^{F_\lambda(\cM)}_m(\id_m)
=\lambda^{1/2}\Tr^{F(\cM)}_m(\id_m)=\lambda^{1/2}d_m^{F(\cM)},
\]
so
\begin{align*}
\Psi^{\Mod(\cB)}_{F_\lambda(\cM)}(F\eta)
&
\underset{\text{\cite[(7)]{bases}}}{=}
\frac{\dim(\End_\cB(1_\cB))^2}{\FPdim(\cB)\psi^\cB(\id_{1_\cB})}
\sum_{m\in\pi_0 F_\lambda(\cM)}  d_m^{F_\lambda(\cM)}\Tr^{F_\lambda(\cM)}_m((F\eta)_m)
\\&
=
\frac{\lambda\dim(\End_\cB(1_\cB))^2}{\FPdim(\cB)\psi^\cB(\id_{1_\cB})}
\sum_{m\in\pi_0 F(\cM)}  
d_m^{F(\cM)}\Tr^{F(\cM)}_m((F\eta)_m)
\\&
\underset{\text{\cite[(7)]{bases}}}{=}
\lambda\Psi^{\Mod(\cB)}_{F(\cM)}(F\eta)=\Psi^{\Mod(\cA)}_\cM(\eta).\qedhere
\end{align*}
\end{pf}

\cref{n2nonfailure} suggests the expectation, which could serve as a desideratum for a future (algebraic) definition of finite $4$-Hilbert spaces: every ordinary/algebraic equivalence of objects in a finite $4$-Hilbert space should be augmentable to an isometric equivalence. 

We can think of the ``rescaling'' $F\mapsto F_\lambda$ in the above proof as a notion of \emph{polar decomposition} for morphisms between 3-Hilbert spaces. We expect that due to the Euler characteristic of spheres, such a polar decomposition is possible for morphisms between finite $n$-Hilbert spaces precisely when $n$ is odd.

\begin{remark}
Despite the above fact, isometric equivalence of \emph{$\Hstar$-multifusion categories} is nonetheless strictly stronger than unitary monoidal equivalence. Indeed, for $\lambda\in(0,1)$ let $\Hilb[2]_\lambda\coloneq(\Hilb^{\boxplus2},\psi_\lambda)$ where $\psi_\lambda\coloneq(\lambda,1-\lambda)$. Since every unitary monoidal autoequivalence of $\Hilb^{\boxplus2}$ permutes the two blocks, $\Hilb[2]_\lambda$ and $\Hilb[2]_\mu$ are isometrically equivalent if and only if $\{\lambda,1-\lambda\}=\{\mu,1-\mu\}$, but all of them share the underlying unitary multifusion category $\Hilb^{\boxplus2}$.
\end{remark}

\THREEDCLASS*

\begin{pf}
  By \cref{THREEDUCH}, it suffices to classify finite 3-Hilbert spaces up to isometric equivalence. Given a finite 3-Hilbert space $\fX$, choose an ONB $\pi_0\fX = \{b_i\}$ in the sense of \cite[Defn. 4.11]{bases}
  and set $b\coloneq\boxplus_i b_i$. Then $\cA\coloneq\Omega_b = \End_\fX(b)$ is an $\Hstar$-multifusion category by \cite[Exm. 4.5]{3Hilb} and $\fX \cong^\dag (\rB\Omega_b)^\cent\cong^\dag\Mod(\cA)$ by \cite[Rmk. 4.13 and Exm. 3.28]{bases}. Now apply \cref{prop:mod-iso-iff-morita}.
\end{pf}

\crefname{appendix}{Appendix}{Appendices}
\Crefname{appendix}{Appendix}{Appendices}
\appendix
\crefalias{section}{appendix}
\crefalias{subsection}{appendix}
\crefalias{subsubsection}{appendix}

\section{Gluing and the construction of the path integral}
\label{sec:gluing-appendix}

\subsection{Gluing lemmas}
\label{subsec:gluing-lemmas}

The results below and in the following subsection can be found in the literature in various places, for instance, in \cite{W06} and \cite{W21}, though with varying levels of detail. Despite this, we have found it useful to work out the details and present the results in the disk-like framework.

Throughout this appendix we freely use the language of \cref{sec:1+1D}: see \cref{skein 1-categories} for the ordinary/traditional 1-category $\Sk_\C(Y,\eta)$, and \cref{def:unitarycategory,def:finiteunitarycategory} for the definitions of a unitary category and a finite unitary category respectively.  

\begin{definition}[Gluing data]
\label{def:gluing-data}
Let $\C$ be any disk-like $n$-category. ($n$D) \defn{gluing data} $\cG=(X,Y,R,\xi,\eta)$ for $\C$ consists of an $n$-manifold $X$ with
\[
  \partial X = R\cup_{\partial\orev{Y}\amalg\partial Y}(\orev Y\amalg Y)
\]
for $(n-1)$-manifolds $Y$ and $R$ together with $\eta\in\undfld[n-2]{\C}{\partial Y}$ and $\xi\in\undfld[n-1]{\C}{R}[\fcj\eta\amalg\eta]$. 
We write $\Sk_\C(Y,\eta)\coloneq\cX_{\A_\C(Y,\eta)}$, we let $X_{\glu}$ (resp. $\xi_\glu$) denote $X$ (resp. $\xi$) glued along $Y$ and $\orev Y$, and for $r,r'\in\undfld[n-1]{\C}{Y}[\eta]$ and a projection $p\in\End_{\Sk_\C(Y,\eta)}(r)$ we define
\[
  V_{r,r'}\coloneq\undfldl[n]{\C}{X}[\xi\cup(\fcj r\amalg r')],\quad V_r\coloneq V_{r,r},
  \quad\text{and}\quad
  pV_rp\coloneq\{\fcj{p}\blt_{\orev Y}v\blt_Yp\mid v\in V_r\}\subset V_r.
\]
Finally, we let $\Gamma^\cG=\Gamma\colon\bigoplus_{r\in\undfld[n-1]{\C}{Y}[\eta]}V_r\to\undfldl[n]{\C}{X_{\glu}}[\xi_{\glu}]$ denote the gluing map.
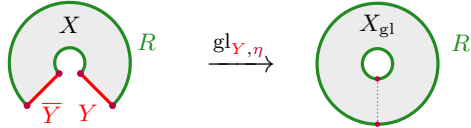
\[
\begin{tkz}[scale=2,cap=round,yscale=1]
  \draw\centerarc(0,0)(225:-45:0.4)coordinate(end1);
  \draw\centerarc(0,0)(225:-45:0.1)coordinate(end2);
  \draw(end1)--(end2);
  \filldraw[draw=black,fill=gray!15](225:.1)--(225:.4)arc(225:-45:.4)--(-45:.1)arc(-45:225:.1);
  \draw (90:.25) node{\scriptsize$X$};
  \draw[very thick,ForestGreen](225:.4)arc(225:-45:.4);
  \draw[very thick,ForestGreen](225:.1)arc(225:-45:.1);
  \draw[very thick,red](225:.1)--(225:.4) (-45:.1)--(-45:.4);
  \draw (250:.35)node[red]{\scriptsize$\orev{Y}$};
  \draw (-70:.35)node[red]{\scriptsize$Y$};
  \draw (225:.4) node[circle,fill=purple,inner sep=0pt,minimum width=.5ex]{};
  \draw (225:.1) node[circle,fill=purple,inner sep=0pt,minimum width=.5ex]{};
  \draw (-45:.4) node[circle,fill=purple,inner sep=0pt,minimum width=.5ex]{};
  \draw (-45:.1) node[circle,fill=purple,inner sep=0pt,minimum width=.5ex]{};
  \draw(15:.525)node[ForestGreen]{\scriptsize$R$};
  \path[very thick, dotted,orange](-90:.1)--(-90:.4);
\end{tkz}
\quad
\xrightarrow{\glu_{\color{red}Y\color{black},\color{purple}\eta}}\quad
\begin{tkz}[scale=2,yscale=1]
  \draw[fill=gray!15] (0,0) circle(0.4);
  \draw[fill=white] (0,0) circle(0.1);
  \draw[densely dotted,red](-90:.1)--(-90:.4);
  \draw (90:.25) node{\scriptsize$X_{\glu}$};
  \draw[very thick,ForestGreen](0,0)circle(.4);
  \draw[very thick,ForestGreen](0,0)circle(.1);
  \draw(14:.565)node[ForestGreen]{\scriptsize$R$};
  \fill[purple] (-90:.1) circle (0.5pt) (-90:.4) circle (0.5pt);
\end{tkz}
\]
\end{definition}

Below we will use the following diagrammatic shorthand for an $n$-field $v\in V_{r,r'}$.
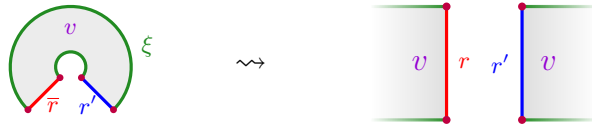
\[
\begin{tkz}[scale=2,cap=round,yscale=1]
  \draw\centerarc(0,0)(225:-45:0.4)coordinate(end1);
  \draw\centerarc(0,0)(225:-45:0.1)coordinate(end2);
  \draw(end1)--(end2);
  \fill[draw=black,fill=gray!15](225:.1)--(225:.4)arc(225:-45:.4)--(-45:.1)arc(-45:225:.1);
  \draw (90:.25) node[violet]{\scriptsize$v$};
  \draw[very thick,ForestGreen](225:.4)arc(225:-45:.4);
  \draw[very thick,ForestGreen](225:.1)arc(225:-45:.1);
  \draw[very thick,red](225:.1)--(225:.4);
  \draw[very thick,blue](-45:.1)--(-45:.4);
  \draw (245:.275)node[red]{\scriptsize$\fcj{r}$};
  \draw (-65:.275)node[blue]{\scriptsize$r'$};
  \draw (225:.4) node[circle,fill=purple,inner sep=0pt,minimum width=.5ex]{};
  \draw (225:.1) node[circle,fill=purple,inner sep=0pt,minimum width=.5ex]{};
  \draw (-45:.4) node[circle,fill=purple,inner sep=0pt,minimum width=.5ex]{};
  \draw (-45:.1) node[circle,fill=purple,inner sep=0pt,minimum width=.5ex]{};
  \draw(15:.525)node[ForestGreen]{\scriptsize$\xi$};
  \path[very thick, dotted,orange](-90:.1)--(-90:.4);
\end{tkz}
\qquad\rightsquigarrow\qquad\quad
\begin{tkz}
\def\tkzht{1.5}
\def\tkzwd{1}
\fill[path fading=west,gray!25](-\tkzwd,0)rectangle(0,\tkzht);
\draw[path fading=west,very thick,ForestGreen](0,0)--(-\tkzwd,0);
\draw[path fading=west,very thick,ForestGreen] (0,\tkzht)--(-\tkzwd,\tkzht);
\draw[very thick,red] (0,0)--node[right]{\scriptsize$r$}(0,\tkzht);
\fill[purple](0,0)circle(1.5pt);
\fill[purple](0,\tkzht)circle(1.5pt);
\node[violet] at (-0.35*\tkzwd,0.5*\tkzht){$v$};
\end{tkz}
\begin{tkz}[xscale=-1]
\def\tkzht{1.5}
\def\tkzwd{1}
\fill[path fading=east,gray!25](-\tkzwd,0)rectangle(0,\tkzht);
\draw[path fading=east,very thick,ForestGreen](0,0)--(-\tkzwd,0);
\draw[path fading=east,very thick,ForestGreen] (0,\tkzht)--(-\tkzwd,\tkzht);
\draw[very thick,blue] (0,0)--node[left]{\scriptsize$r'$}(0,\tkzht);
\fill[purple](0,0)circle(1.5pt);
\fill[purple](0,\tkzht)circle(1.5pt);
\node[violet] at (-0.35*\tkzwd,0.5*\tkzht){$v$};
\end{tkz}
\]

\begin{lemma}
\label{glu-slide}
Let $\C$ be a disk-like $n$-category and let $\cG=(X,Y,R,\xi,\eta)$ be gluing data for $\C$. Then for all $r,r'\in\undfld[n-1]{\C}{Y}[\eta]$, all $e\in\Sk_\C(Y,\eta)(r\to r')$ and all $v\in V_{r',r}$, 
\[
    \Gamma(v\blt_Y e)=\Gamma(\fcj{e^\dag}\blt_{\orev{Y}} v),
\]
that is,
\[
\Gamma\!\left(
\begin{tkz}
\def\tkzht{1.5}
\def\tkzwd{1}
\def\tkzblg{1.4}
\fill[path fading=west,gray!25](-\tkzwd,0)rectangle(0,\tkzht);
\draw[path fading=west,very thick,ForestGreen](0,0)--(-\tkzwd,0);
\draw[path fading=west,very thick,ForestGreen] (0,\tkzht)--(-\tkzwd,\tkzht);
\fill[purple](0,0)circle(1.5pt);
\fill[purple](0,\tkzht)circle(1.5pt);
\node[violet] at (-0.35*\tkzwd,0.5*\tkzht){$v$};
\draw[very thick,red] (0,0)--(0,\tkzht);
\end{tkz}
\quad
\begin{tkz}[xscale=-1]
\def\tkzht{1.5}
\def\tkzwd{1}
\def\tkzblg{1.4}
\fill[path fading=east,gray!25](-\tkzwd,0)rectangle(0,\tkzht);
\draw[path fading=east,very thick,ForestGreen](0,0)--(-\tkzwd,0);
\draw[path fading=east,very thick,ForestGreen] (0,\tkzht)--(-\tkzwd,\tkzht);
\fill[purple](0,0)circle(1.5pt);
\fill[purple](0,\tkzht)circle(1.5pt);
\node[violet] at (-0.35*\tkzwd,0.5*\tkzht){$v$};
\begin{scope}
\clip(0,0)to[bend right=80,distance=\tkzblg cm](0,\tkzht)--cycle;
\fill[brown!15](0,0) rectangle (5,\tkzht);
\node[purple] at ({0.3*\tkzblg}, {0.5*\tkzht}){\scriptsize$e$};
\end{scope}
\draw[very thick,blue] (0,0)--(0,\tkzht);
\draw[red,very thick](0,0)to[bend right=80,distance=\tkzblg cm](0,\tkzht);
\end{tkz}
\right)
\quad
=
\quad
\Gamma\!\left(
\begin{tkz}
\def\tkzht{1.5}
\def\tkzwd{1}
\def\tkzblg{1.4}
\fill[path fading=west,gray!25](-\tkzwd,0)rectangle(0,\tkzht);
\draw[path fading=west,very thick,ForestGreen](0,0)--(-\tkzwd,0);
\draw[path fading=west,very thick,ForestGreen] (0,\tkzht)--(-\tkzwd,\tkzht);
\fill[purple](0,0)circle(1.5pt);
\fill[purple](0,\tkzht)circle(1.5pt);
\node[violet] at (-0.35*\tkzwd,0.5*\tkzht){$v$};
\begin{scope}
\clip(0,0)to[bend right=80,distance=\tkzblg cm](0,\tkzht)--cycle;
\fill[brown!15](0,0) rectangle (5,\tkzht);
\node[purple] at ({0.3*\tkzblg}, {0.5*\tkzht}){\scriptsize$e^\dag$};
\end{scope}
\draw[very thick,red] (0,0)--(0,\tkzht);
\draw[blue,very thick](0,0)to[bend right=80,distance=\tkzblg cm](0,\tkzht);
\end{tkz}
\quad
\begin{tkz}[xscale=-1]
\def\tkzht{1.5}
\def\tkzwd{1}
\def\tkzblg{1.4}
\fill[path fading=east,gray!25](-\tkzwd,0)rectangle(0,\tkzht);
\draw[path fading=east,very thick,ForestGreen](0,0)--(-\tkzwd,0);
\draw[path fading=east,very thick,ForestGreen] (0,\tkzht)--(-\tkzwd,\tkzht);
\fill[purple](0,0)circle(1.5pt);
\fill[purple](0,\tkzht)circle(1.5pt);
\node[violet] at (-0.35*\tkzwd,0.5*\tkzht){$v$};
\draw[very thick,blue] (0,0)--(0,\tkzht);
\end{tkz}
\right)
\]
where $e^\dag\coloneq\fcj{(\id_Y\times\iota)_*e}\in\Sk_\C(Y,\eta)(r'\to r)$ is the dagger of $e$ in the skein 1-category $\Sk_\C(Y,\eta)$.
\end{lemma}

\begin{pf}
  Recall that $\undfldl[n]{\C}{X_\glu}[ \xi_\glu]$ is constructed as the colimit over permissible ball decompositions of $X_\glu$. Observe that $\Gamma(v\blt_Y e)$ and $\Gamma(\fcj{e^\dag}\blt_{\orev{Y}} v)$ represent the same colimit class of $\undfldl[n]{\C}{X_\glu}[ \xi_\glu]$: they are both coarsenings of the splitting of $X_\glu$ as $X$ glued to $Y\times D^1$ along the gluing locus $\orev{Y}\amalg Y$, so the result follows.
\end{pf}

\begin{lemma}
\label{glu-surj}
If $\C$ is any disk-like $n$-category and $\cG=(X,Y,R,\xi,\eta)$ is any gluing data for $\C$, then the gluing map $\Gamma^\cG$ is surjective.
\end{lemma}

\begin{pf}
Let $f \in \undfldl[n]{\C}{X_{\glu}}[ \xi_{\glu}]$, say $f = \sum_i \lambda_i [\alpha_i]$ for $\alpha_i \in \undfld[n]{\C}{X_{\glu}}[ \xi_{\glu}]$. For each $i$, choose $\widetilde\alpha_i \in \undfld[n]{\C}{X_{\glu}}[ \xi_{\glu}]$ isotopic-rel-boundary to $\alpha_i$ that is transverse to $Y$ (i.e., has a string locus transverse to $Y$), which exists by \ref{CS}, since isotopies act on string loci in the obvious way. Then $[\widetilde\alpha_i] = [\alpha_i]$ by \ref{CTi}, so $f = \sum_i \lambda_i [\widetilde\alpha_i]$. Cutting $\widetilde\alpha_i$ along $Y$ exhibits $[\widetilde\alpha_i]\in\Gamma(V_{r_i})$ for $r_i\coloneq\widetilde\alpha_i|_Y\in\undfld[n-1]{\C}{Y}[\eta]$. Thus $\Gamma$ is surjective.
\end{pf}

\begin{remark}
In the proof of \cref{glu-surj}, we appealed to the splitting axiom \ref{CS}. However, Kevin Walker has informed the author that a much weaker condition would suffice; see \cite[Appendix C]{RW}.
\end{remark}

\begin{lemma}[Compression]
\label{glu-compress}
Let $\C$ be a disk-like $n$-category, let $\cG=(X,Y,R,\xi,\eta)$ be gluing data for $\C$, let $r\in\undfld[n-1]{\C}{Y}[\eta]$, and let $\{q_j\}$ be a complete set of orthogonal projections in $\End_{\Sk_\C(Y,\eta)}(r)$ (note this requires $\End_{\Sk_\C(Y,\eta)}(r)$ to be finite-dimensional). Then
\begin{equation}
\label{eq:compress}
  \Gamma(v)=\sum_j\Gamma(\fcj{q_j}\blt_{\orev Y}v\blt_Y q_j)
  \qquad\forall v\in V_r.
\end{equation}
In particular, $\Gamma(V_r)\subset\sum_j\Gamma(q_jV_rq_j)$.
\end{lemma}
 
\begin{pf}
Fix $v\in V_r$. Then
\begin{align*}
\Gamma(v)
=
\Gamma\!\left(
\begin{tkz}
\def\tkzht{1.5}
\def\tkzwd{0.8}
\fill[path fading=west,gray!25](-\tkzwd,0)rectangle(0,\tkzht);
\draw[path fading=west,very thick,ForestGreen](0,0)--(-\tkzwd,0);
\draw[path fading=west,very thick,ForestGreen] (0,\tkzht)--(-\tkzwd,\tkzht);
\draw[very thick,red] (0,0)--(0,\tkzht);
\fill[purple](0,0)circle(1.5pt);
\fill[purple](0,\tkzht)circle(1.5pt);
\node[violet] at (-0.35*\tkzwd,0.5*\tkzht){$v$};
\end{tkz}
\right.
&\;\,
\left.
\begin{tkz}[xscale=-1]
\def\tkzht{1.5}
\def\tkzwd{0.8}
\fill[path fading=east,gray!25](-\tkzwd,0)rectangle(0,\tkzht);
\draw[path fading=east,very thick,ForestGreen](0,0)--(-\tkzwd,0);
\draw[path fading=east,very thick,ForestGreen] (0,\tkzht)--(-\tkzwd,\tkzht);
\draw[very thick,red] (0,0)--(0,\tkzht);
\fill[purple](0,0)circle(1.5pt);
\fill[purple](0,\tkzht)circle(1.5pt);
\node[violet] at (-0.35*\tkzwd,0.5*\tkzht){$v$};
\end{tkz}
\right)
\!=
\sum_{j,k}
\Gamma\!\left(
\begin{tkz}
\def\tkzht{1.5}
\def\tkzwd{0.8}
\def\tkzblg{1.2}
\fill[path fading=west,gray!25](-\tkzwd,0)rectangle(0,\tkzht);
\draw[path fading=west,very thick,ForestGreen](0,0)--(-\tkzwd,0);
\draw[path fading=west,very thick,ForestGreen] (0,\tkzht)--(-\tkzwd,\tkzht);
\fill[purple](0,0)circle(1.5pt);
\fill[purple](0,\tkzht)circle(1.5pt);
\node[violet] at (-0.35*\tkzwd,0.5*\tkzht){$v$};
\begin{scope}
\clip(0,0)to[bend right=80,distance=\tkzblg cm](0,\tkzht)--cycle;
\fill[brown!15](0,0) rectangle (5,\tkzht);
\node[brown] at ({0.3*\tkzblg}, {0.5*\tkzht}){\scriptsize$q_k$};
\end{scope}
\draw[very thick,red] (0,0)--(0,\tkzht);
\draw[red,very thick](0,0)to[bend right=80,distance=\tkzblg cm](0,\tkzht);
\end{tkz}
\!\!\!\!\!\!\!
\begin{tkz}[xscale=-1]
\def\tkzht{1.5}
\def\tkzwd{0.8}
\def\tkzblg{1.2}
\fill[path fading=east,gray!25](-\tkzwd,0)rectangle(0,\tkzht);
\draw[path fading=east,very thick,ForestGreen](0,0)--(-\tkzwd,0);
\draw[path fading=east,very thick,ForestGreen] (0,\tkzht)--(-\tkzwd,\tkzht);
\fill[purple](0,0)circle(1.5pt);
\fill[purple](0,\tkzht)circle(1.5pt);
\node[violet] at (-0.35*\tkzwd,0.5*\tkzht){$v$};
\begin{scope}
\clip(0,0)to[bend right=80,distance=\tkzblg cm](0,\tkzht)--cycle;
\fill[brown!15](0,0) rectangle (5,\tkzht);
\node[brown] at ({0.3*\tkzblg}, {0.5*\tkzht}){\scriptsize$q_j$};
\end{scope}
\draw[very thick,red] (0,0)--(0,\tkzht);
\draw[red,very thick](0,0)to[bend right=80,distance=\tkzblg cm](0,\tkzht);
\end{tkz}
\right)
\!\underset{\eqref{glu-slide}}{=}
\sum_{j,k} 
\Gamma\!\left(
\begin{tkz}
\def\tkzht{1.5}
\def\tkzwd{0.8}
\fill[path fading=west,gray!25](-\tkzwd,0)rectangle(0,\tkzht);
\draw[path fading=west,very thick,ForestGreen](0,0)--(-\tkzwd,0);
\draw[path fading=west,very thick,ForestGreen] (0,\tkzht)--(-\tkzwd,\tkzht);
\fill[purple](0,0)circle(1.5pt);
\fill[purple](0,\tkzht)circle(1.5pt);
\node[violet] at (-0.35*\tkzwd,0.5*\tkzht){$v$};
\draw[very thick,red] (0,0)--(0,\tkzht);
\end{tkz}
\!\!\!\!
\begin{tkz}[xscale=-1]
\def\tkzht{1.5}
\def\tkzwd{0.8}
\def\tkzblg{1.6}
\def\tkzmid{0.8}
\fill[path fading=east,gray!25](-\tkzwd,0)rectangle(0,\tkzht);
\draw[path fading=east,very thick,ForestGreen](0,0)--(-\tkzwd,0);
\draw[path fading=east,very thick,ForestGreen] (0,\tkzht)--(-\tkzwd,\tkzht);
\fill[purple](0,0)circle(1.5pt);
\fill[purple](0,\tkzht)circle(1.5pt);
\node[violet] at (-0.35*\tkzwd,0.5*\tkzht){$v$};
\begin{scope}
\clip(0,0)to[bend right=80,distance=\tkzblg cm](0,\tkzht)--cycle;
\fill[brown!15](0,0) rectangle (5,\tkzht);
\node[brown] at ({0.3*\tkzmid}, {0.5*\tkzht}){\scriptsize$q_j$};
\node[brown] at ({0.55*\tkzblg}, {0.5*\tkzht}){\scriptsize$q_k$};
\end{scope}
\draw[very thick,red] (0,0)--(0,\tkzht);
\draw[red,very thick](0,0)to[bend right=80,distance=\tkzmid cm](0,\tkzht);
\draw[red,very thick](0,0)to[bend right=80,distance=\tkzblg cm](0,\tkzht);
\end{tkz}
\right)
\\
&=
\sum_{j} 
\Gamma\!\left(
\begin{tkz}
\def\tkzht{1.5}
\def\tkzwd{0.8}
\fill[path fading=west,gray!25](-\tkzwd,0)rectangle(0,\tkzht);
\draw[path fading=west,very thick,ForestGreen](0,0)--(-\tkzwd,0);
\draw[path fading=west,very thick,ForestGreen] (0,\tkzht)--(-\tkzwd,\tkzht);
\fill[purple](0,0)circle(1.5pt);
\fill[purple](0,\tkzht)circle(1.5pt);
\node[violet] at (-0.35*\tkzwd,0.5*\tkzht){$v$};
\draw[very thick,red] (0,0)--(0,\tkzht);
\end{tkz}
\!\!
\begin{tkz}[xscale=-1]
\def\tkzht{1.5}
\def\tkzwd{0.8}
\def\tkzblg{1.6}
\def\tkzmid{0.8}
\fill[path fading=east,gray!25](-\tkzwd,0)rectangle(0,\tkzht);
\draw[path fading=east,very thick,ForestGreen](0,0)--(-\tkzwd,0);
\draw[path fading=east,very thick,ForestGreen] (0,\tkzht)--(-\tkzwd,\tkzht);
\fill[purple](0,0)circle(1.5pt);
\fill[purple](0,\tkzht)circle(1.5pt);
\node[violet] at (-0.35*\tkzwd,0.5*\tkzht){$v$};
\begin{scope}
\clip(0,0)to[bend right=80,distance=\tkzblg cm](0,\tkzht)--cycle;
\fill[brown!15](0,0) rectangle (5,\tkzht);
\node[brown] at ({0.3*\tkzmid}, {0.5*\tkzht}){\scriptsize$q_j$};
\node[brown] at ({0.55*\tkzblg}, {0.5*\tkzht}){\scriptsize$q_j$};
\end{scope}
\draw[very thick,red] (0,0)--(0,\tkzht);
\draw[red,very thick](0,0)to[bend right=80,distance=\tkzmid cm](0,\tkzht);
\draw[red,very thick](0,0)to[bend right=80,distance=\tkzblg cm](0,\tkzht);
\end{tkz}
\right)
= \sum_j 
\Gamma\!\left(
\begin{tkz}
\def\tkzht{1.5}
\def\tkzwd{0.8}
\def\tkzblg{1.4}
\fill[path fading=west,gray!25](-\tkzwd,0)rectangle(0,\tkzht);
\draw[path fading=west,very thick,ForestGreen](0,0)--(-\tkzwd,0);
\draw[path fading=west,very thick,ForestGreen] (0,\tkzht)--(-\tkzwd,\tkzht);
\fill[purple](0,0)circle(1.5pt);
\fill[purple](0,\tkzht)circle(1.5pt);
\node[violet] at (-0.35*\tkzwd,0.5*\tkzht){$v$};
\begin{scope}
\clip(0,0)to[bend right=80,distance=\tkzblg cm](0,\tkzht)--cycle;
\fill[brown!15](0,0) rectangle (5,\tkzht);
\node[brown] at ({0.3*\tkzblg}, {0.5*\tkzht}){\scriptsize$q_j$};
\end{scope}
\draw[very thick,red] (0,0)--(0,\tkzht);
\draw[red,very thick](0,0)to[bend right=80,distance=\tkzblg cm](0,\tkzht);
\end{tkz}
\!\!\!\!
\begin{tkz}[xscale=-1]
\def\tkzht{1.5}
\def\tkzwd{0.8}
\def\tkzblg{1.4}
\fill[path fading=east,gray!25](-\tkzwd,0)rectangle(0,\tkzht);
\draw[path fading=east,very thick,ForestGreen](0,0)--(-\tkzwd,0);
\draw[path fading=east,very thick,ForestGreen] (0,\tkzht)--(-\tkzwd,\tkzht);
\fill[purple](0,0)circle(1.5pt);
\fill[purple](0,\tkzht)circle(1.5pt);
\node[violet] at (-0.35*\tkzwd,0.5*\tkzht){$v$};
\begin{scope}
\clip(0,0)to[bend right=80,distance=\tkzblg cm](0,\tkzht)--cycle;
\fill[brown!15](0,0) rectangle (5,\tkzht);
\node[brown] at ({0.3*\tkzblg}, {0.5*\tkzht}){\scriptsize$q_j$};
\end{scope}
\draw[very thick,red] (0,0)--(0,\tkzht);
\draw[red,very thick](0,0)to[bend right=80,distance=\tkzblg cm](0,\tkzht);
\end{tkz}
\right).
\qedhere
\end{align*}
\end{pf}

\begin{lemma}[Corner invariance]
\label{glu-corners}
Let $\C$ be a disk-like $n$-category, let $\cG=(X,Y,R,\xi,\eta)$ be gluing data for $\C$, and suppose $\Sk_\C(Y,\eta)$ is unitary. Then for any $r,r'\in\Sk_\C(Y,\eta)$ and any projections $p\in\End_{\Sk_\C(Y,\eta)}(r)$ and $p'\in\End_{\Sk_\C(Y,\eta)}(r')$ with $(r,p)\cong(r',p')$ in
$\Sk_\C(Y,\eta)^\cent$,
\begin{equation}
\label{corners-agree}
  \Gamma(pV_rp)=\Gamma(p'V_{r'}p').
\end{equation}
\end{lemma}

\begin{pf}
Suppose $(r,p)$ and $(r',p')$ are isomorphic objects of $\Sk_\C(Y,\eta)^\cent$, so there is some
$u\in p'\Sk_\C(Y,\eta)(r\to r')p$ with $u^\dag u=p$ and $uu^\dag=p'$. Let $v\in pV_rp$, so that $v=\fcj{p}\blt_{\orev Y}v\blt_Yp$. Then
\[
\begin{aligned}
\Gamma(v)
&
=
\Gamma\!\left(
\begin{tkz}
\def\tkzht{1.5}
\def\tkzwd{1}
\def\tkzblg{1.4}
\fill[path fading=west,gray!25](-\tkzwd,0)rectangle(0,\tkzht);
\draw[path fading=west,very thick,ForestGreen](0,0)--(-\tkzwd,0);
\draw[path fading=west,very thick,ForestGreen] (0,\tkzht)--(-\tkzwd,\tkzht);
\fill[purple](0,0)circle(1.5pt);
\fill[purple](0,\tkzht)circle(1.5pt);
\node[violet] at (-0.35*\tkzwd,0.5*\tkzht){$v$};
\begin{scope}
\clip(0,0)to[bend right=80,distance=\tkzblg cm](0,\tkzht)--cycle;
\fill[brown!15](0,0) rectangle (5,\tkzht);
\node[brown] at ({0.3*\tkzblg}, {0.5*\tkzht}){\scriptsize$p$};
\end{scope}
\draw[very thick,red] (0,0)--(0,\tkzht);
\draw[red,very thick](0,0)to[bend right=80,distance=\tkzblg cm](0,\tkzht);
\end{tkz}
\begin{tkz}[xscale=-1]
\def\tkzht{1.5}
\def\tkzwd{1}
\def\tkzblg{1.4}
\fill[path fading=east,gray!25](-\tkzwd,0)rectangle(0,\tkzht);
\draw[path fading=east,very thick,ForestGreen](0,0)--(-\tkzwd,0);
\draw[path fading=east,very thick,ForestGreen] (0,\tkzht)--(-\tkzwd,\tkzht);
\fill[purple](0,0)circle(1.5pt);
\fill[purple](0,\tkzht)circle(1.5pt);
\node[violet] at (-0.35*\tkzwd,0.5*\tkzht){$v$};
\begin{scope}
\clip(0,0)to[bend right=80,distance=\tkzblg cm](0,\tkzht)--cycle;
\fill[brown!15](0,0) rectangle (5,\tkzht);
\node[brown] at ({0.3*\tkzblg}, {0.5*\tkzht}){\scriptsize$p$};
\end{scope}
\draw[very thick,red] (0,0)--(0,\tkzht);
\draw[red,very thick](0,0)to[bend right=80,distance=\tkzblg cm](0,\tkzht);
\end{tkz}
\right)
=
\Gamma\!\left(
\begin{tkz}
\def\tkzht{1.5}
\def\tkzwd{1}
\def\tkzblg{0.8}
\fill[path fading=west,gray!25](-\tkzwd,0)rectangle(0,\tkzht);
\draw[path fading=west,very thick,ForestGreen](0,0)--(-\tkzwd,0);
\draw[path fading=west,very thick,ForestGreen] (0,\tkzht)--(-\tkzwd,\tkzht);
\fill[purple](0,0)circle(1.5pt);
\fill[purple](0,\tkzht)circle(1.5pt);
\node[violet] at (-0.35*\tkzwd,0.5*\tkzht){$v$};
\begin{scope}
\clip(0,0)to[bend right=80,distance=\tkzblg cm](0,\tkzht)--cycle;
\fill[brown!15](0,0) rectangle (5,\tkzht);
\node[brown] at ({0.3*\tkzblg}, {0.5*\tkzht}){\scriptsize$p$};
\end{scope}
\draw[very thick,red] (0,0)--(0,\tkzht);
\draw[red,very thick](0,0)to[bend right=80,distance=\tkzblg cm](0,\tkzht);
\end{tkz}
\begin{tkz}[xscale=-1]
\def\tkzht{1.5}
\def\tkzwd{1}
\def\tkzblg{1.6}
\def\tkzmid{0.8}
\fill[path fading=east,gray!25](-\tkzwd,0)rectangle(0,\tkzht);
\draw[path fading=east,very thick,ForestGreen](0,0)--(-\tkzwd,0);
\draw[path fading=east,very thick,ForestGreen] (0,\tkzht)--(-\tkzwd,\tkzht);
\fill[purple](0,0)circle(1.5pt);
\fill[purple](0,\tkzht)circle(1.5pt);
\node[violet] at (-0.35*\tkzwd,0.5*\tkzht){$v$};
\begin{scope}
\clip(0,0)to[bend right=80,distance=\tkzblg cm](0,\tkzht)--cycle;
\fill[brown!15](0,0) rectangle (5,\tkzht);
\node[brown] at ({0.3*\tkzmid}, {0.5*\tkzht}){\scriptsize$p$};
\node[brown] at ({0.55*\tkzblg}, {0.5*\tkzht}){\scriptsize$p$};
\end{scope}
\draw[very thick,red] (0,0)--(0,\tkzht);
\draw[red,very thick](0,0)to[bend right=80,distance=\tkzmid cm](0,\tkzht);
\draw[red,very thick](0,0)to[bend right=80,distance=\tkzblg cm](0,\tkzht);
\end{tkz}
\right)
\\&
=
\Gamma\!\left(
\begin{tkz}
\def\tkzht{1.5}
\def\tkzwd{1}
\def\tkzblg{1.6}
\def\tkzmid{0.8}
\fill[path fading=west,gray!25](-\tkzwd,0)rectangle(0,\tkzht);
\draw[path fading=west,very thick,ForestGreen](0,0)--(-\tkzwd,0);
\draw[path fading=west,very thick,ForestGreen] (0,\tkzht)--(-\tkzwd,\tkzht);
\fill[purple](0,0)circle(1.5pt);
\fill[purple](0,\tkzht)circle(1.5pt);
\node[violet] at (-0.35*\tkzwd,0.5*\tkzht){$v$};
\begin{scope}
\clip(0,0)to[bend right=80,distance=\tkzblg cm](0,\tkzht)--cycle;
\fill[purple!25](0,0) rectangle (5,\tkzht);
\clip(0,0)to[bend right=80,distance=\tkzmid cm](0,\tkzht)--cycle;
\fill[purple!15](0,0) rectangle (5,\tkzht);
\end{scope}
\node[purple] at ({0.35*\tkzmid}, {0.5*\tkzht}){\scriptsize$u$};
\node[purple] at ({0.55*\tkzblg}, {0.5*\tkzht}){\scriptsize$u^\dag$};
\draw[very thick,red] (0,0)--(0,\tkzht);
\draw[blue,very thick](0,0)to[bend right=80,distance=\tkzmid cm](0,\tkzht);
\draw[red,very thick](0,0)to[bend right=80,distance=\tkzblg cm](0,\tkzht);
\end{tkz}
\begin{tkz}[xscale=-1]
\def\tkzht{1.5}
\def\tkzwd{1}
\def\tkzblg{2.4}
\fill[path fading=east,gray!25](-\tkzwd,0)rectangle(0,\tkzht);
\draw[path fading=east,very thick,ForestGreen](0,0)--(-\tkzwd,0);
\draw[path fading=east,very thick,ForestGreen] (0,\tkzht)--(-\tkzwd,\tkzht);
\fill[purple](0,0)circle(1.5pt);
\fill[purple](0,\tkzht)circle(1.5pt);
\node[violet] at (-0.35*\tkzwd,0.5*\tkzht){$v$};
\begin{scope}
\clip(0,0)to[bend right=85,distance=\tkzblg cm](0,\tkzht)--cycle;
\fill[purple!25](0,0) rectangle (5,\tkzht);
\clip(0,0)to[bend right=85,distance=0.75*\tkzblg cm](0,\tkzht)--cycle;
\fill[purple!15](0,0) rectangle (5,\tkzht);
\clip(0,0)to[bend right=85,distance=0.5*\tkzblg cm](0,\tkzht)--cycle;
\fill[purple!25](0,0) rectangle (5,\tkzht);
\clip(0,0)to[bend right=85,distance=0.25*\tkzblg cm](0,\tkzht)--cycle;
\fill[purple!15](0,0) rectangle (5,\tkzht);
\end{scope}
\node[purple] at ({0.085*\tkzblg}, {0.5*\tkzht}){\scriptsize$u$};
\node[purple] at ({0.275*\tkzblg}, {0.5*\tkzht}){\scriptsize$u^\dag$};
\node[purple] at ({0.45*\tkzblg}, {0.5*\tkzht}){\scriptsize$u$};
\node[purple] at ({0.65*\tkzblg}, {0.5*\tkzht}){\scriptsize$u^\dag$};
\draw[very thick,red] (0,0)--(0,\tkzht);
\draw[blue,very thick](0,0)to[bend right=85,distance=0.25*\tkzblg cm](0,\tkzht);
\draw[red,very thick](0,0)to[bend right=85,distance=0.5*\tkzblg cm](0,\tkzht);
\draw[blue,very thick](0,0)to[bend right=85,distance=0.75*\tkzblg cm](0,\tkzht);
\draw[red,very thick](0,0)to[bend right=85,distance=\tkzblg cm](0,\tkzht);
\end{tkz}
\right)
% \\&
=
\Gamma\!\left(
\begin{tkz}
\def\tkzht{1.5}
\def\tkzwd{1}
\def\tkzblg{1.6}
\def\tkztw{1.1}
\def\tkzon{0.55}
\fill[path fading=west,gray!25](-\tkzwd,0)rectangle(0,\tkzht);
\draw[path fading=west,very thick,ForestGreen](0,0)--(-\tkzwd,0);
\draw[path fading=west,very thick,ForestGreen] (0,\tkzht)--(-\tkzwd,\tkzht);
\fill[purple](0,0)circle(1.5pt);
\fill[purple](0,\tkzht)circle(1.5pt);
\node[violet] at (-0.35*\tkzwd,0.5*\tkzht){$v$};
\begin{scope}
\clip(0,0)to[bend right=80,distance=\tkzblg cm](0,\tkzht)--cycle;
\fill[purple!15](0,0) rectangle (5,\tkzht);
\clip(0,0)to[bend right=80,distance=\tkztw cm](0,\tkzht)--cycle;
\fill[purple!25](0,0) rectangle (5,\tkzht);
\clip(0,0)to[bend right=80,distance=\tkzon cm](0,\tkzht)--cycle;
\fill[purple!15](0,0) rectangle (5,\tkzht);
\end{scope}
\node[purple] at ({0.3*\tkzon}, {0.5*\tkzht}){\scriptsize$u$};
\node[purple] at ({\tkzon + 0.125*(\tkztw-\tkzon)}, {0.5*\tkzht}){\scriptsize$u^\dag$};
\node[purple] at ({0.9*\tkztw + 0*(\tkzblg-\tkztw)}, {0.5*\tkzht}){\scriptsize$u$};
\draw[very thick,red] (0,0)--(0,\tkzht);
\draw[blue,very thick](0,0)to[bend right=80,distance=\tkzon cm](0,\tkzht);
\draw[red,very thick](0,0)to[bend right=80,distance=\tkztw cm](0,\tkzht);
\draw[blue,very thick](0,0)to[bend right=80,distance=\tkzblg cm](0,\tkzht);
\end{tkz}
\begin{tkz}[xscale=-1]
\def\tkzht{1.5}
\def\tkzwd{1}
\def\tkzblg{1.6}
\def\tkztw{1.1}
\def\tkzon{0.55}
\fill[path fading=east,gray!25](-\tkzwd,0)rectangle(0,\tkzht);
\draw[path fading=east,very thick,ForestGreen](0,0)--(-\tkzwd,0);
\draw[path fading=east,very thick,ForestGreen] (0,\tkzht)--(-\tkzwd,\tkzht);
\fill[purple](0,0)circle(1.5pt);
\fill[purple](0,\tkzht)circle(1.5pt);
\node[violet] at (-0.35*\tkzwd,0.5*\tkzht){$v$};
\begin{scope}
\clip(0,0)to[bend right=80,distance=\tkzblg cm](0,\tkzht)--cycle;
\fill[purple!15](0,0) rectangle (5,\tkzht);
\clip(0,0)to[bend right=80,distance=\tkztw cm](0,\tkzht)--cycle;
\fill[purple!25](0,0) rectangle (5,\tkzht);
\clip(0,0)to[bend right=80,distance=\tkzon cm](0,\tkzht)--cycle;
\fill[purple!15](0,0) rectangle (5,\tkzht);
\end{scope}
\node[purple] at ({0.3*\tkzon}, {0.5*\tkzht}){\scriptsize$u$};
\node[purple] at ({\tkzon + 0.1*(\tkztw-\tkzon)}, {0.5*\tkzht}){\scriptsize$u^\dag$};
\node[purple] at ({0.9*\tkztw + 0*(\tkzblg-\tkztw)}, {0.5*\tkzht}){\scriptsize$u$};
\draw[very thick,red] (0,0)--(0,\tkzht);
\draw[blue,very thick](0,0)to[bend right=80,distance=\tkzon cm](0,\tkzht);
\draw[red,very thick](0,0)to[bend right=80,distance=\tkztw cm](0,\tkzht);
\draw[blue,very thick](0,0)to[bend right=80,distance=\tkzblg cm](0,\tkzht);
\end{tkz}
\right)
\\&
\underset{\eqref{glu-slide}}{=}
\Gamma\!\left(
\begin{tkz}
\def\tkzht{1.5}
\def\tkzwd{1}
\def\tkzblg{1.6}
\def\tkzmid{0.6}
\fill[path fading=west,gray!25](-\tkzwd,0)rectangle(0,\tkzht);
\draw[path fading=west,very thick,ForestGreen](0,0)--(-\tkzwd,0);
\draw[path fading=west,very thick,ForestGreen] (0,\tkzht)--(-\tkzwd,\tkzht);
\fill[purple](0,0)circle(1.5pt);
\fill[purple](0,\tkzht)circle(1.5pt);
\node[violet] at (-0.35*\tkzwd,0.5*\tkzht){$v$};
\begin{scope}
\clip(0,0)to[bend right=80,distance=\tkzblg cm](0,\tkzht)--cycle;
\fill[yellow!30](0,0) rectangle (5,\tkzht);
\clip(0,0)to[bend right=80,distance=\tkzmid cm](0,\tkzht)--cycle;
\fill[purple!15](0,0) rectangle (5,\tkzht);
\end{scope}
\node[purple] at ({0.3*\tkzmid}, {0.5*\tkzht}){\scriptsize$u$};
\node[brown] at ({\tkzmid + 0.15*(\tkzblg-\tkzmid)}, {0.5*\tkzht}){\scriptsize$p'$};
\draw[very thick,red] (0,0)--(0,\tkzht);
\draw[blue,very thick](0,0)to[bend right=80,distance=\tkzmid cm](0,\tkzht);
\draw[blue,very thick](0,0)to[bend right=80,distance=\tkzblg cm](0,\tkzht);
\end{tkz}
\begin{tkz}[xscale=-1]
\def\tkzht{1.5}
\def\tkzwd{1}
\def\tkzblg{1.6}
\def\tkzmid{0.6}
\fill[path fading=east,gray!25](-\tkzwd,0)rectangle(0,\tkzht);
\draw[path fading=east,very thick,ForestGreen](0,0)--(-\tkzwd,0);
\draw[path fading=east,very thick,ForestGreen] (0,\tkzht)--(-\tkzwd,\tkzht);
\fill[purple](0,0)circle(1.5pt);
\fill[purple](0,\tkzht)circle(1.5pt);
\node[violet] at (-0.35*\tkzwd,0.5*\tkzht){$v$};
\begin{scope}
\clip(0,0)to[bend right=80,distance=\tkzblg cm](0,\tkzht)--cycle;
\fill[yellow!30](0,0) rectangle (5,\tkzht);
\clip(0,0)to[bend right=80,distance=\tkzmid cm](0,\tkzht)--cycle;
\fill[purple!15](0,0) rectangle (5,\tkzht);
\end{scope}
\node[purple] at ({0.35*\tkzmid}, {0.5*\tkzht}){\scriptsize$u$};
\node[brown] at ({\tkzmid + 0.15*(\tkzblg-\tkzmid)}, {0.5*\tkzht}){\scriptsize$p'$};
\draw[very thick,red] (0,0)--(0,\tkzht);
\draw[blue,very thick](0,0)to[bend right=80,distance=\tkzmid cm](0,\tkzht);
\draw[blue,very thick](0,0)to[bend right=80,distance=\tkzblg cm](0,\tkzht);
\end{tkz}
\right)
\;\;=\;\;
\Gamma(\fcj{p'}\blt_{\orev Y}(\fcj{u}\blt_{\orev Y}v\blt_{Y}u)\blt_{Y}p')
\;\;\in\;\;
\Gamma(p'V_{r'}p').
\end{aligned}
\]
A similar argument shows $\Gamma(\fcj{p'}\blt_{\orev Y} v'\blt_Y p') \in \Gamma(pV_r p)$ for all $v' \in V_{r'}$, so \eqref{corners-agree} holds.
\end{pf}

\begin{corollary}
\label{glu-span}
Let $\C$ be any disk-like $n$-category and let $\cG=(X,Y,R,\xi,\eta)$ be gluing data for $\C$. If $\Sk_\C(Y,\eta)$ is unitary, say with $\Irr(\Sk_\C(Y,\eta)^\cent)=\{(r_i,p_i)\}_{i\in I}$ (possibly infinite), then
\begin{equation}
\label{eq:star-generalized}
  \undfldl[n]{\C}{X_{\glu}}[\xi_{\glu}]=\sum_{i\in I} \Gamma(p_iV_{r_i}p_i)
\tag{$\ast$}
\end{equation}
(finite linear combinations). Diagrammatically,
\[
\qquad\qquad
\begin{tkz}
\def\tkzht{1.5}
\def\tkzwd{2.5}
\fill[path fading=west,gray!25](-\tkzwd,0)rectangle(0,\tkzht);
\draw[path fading=west,very thick,ForestGreen](0,0)--(-\tkzwd,0);
\draw[path fading=west,very thick,ForestGreen] (0,\tkzht)--(-\tkzwd,\tkzht);
\fill[purple](0,0)circle(1.5pt);
\fill[purple](0,\tkzht)circle(1.5pt);
\end{tkz}
\mkern-5mu
\begin{tkz}[xscale=-1]
\def\tkzht{1.5}
\def\tkzwd{2.5}
\fill[path fading=east,gray!25](-\tkzwd,0)rectangle(0,\tkzht);
\draw[path fading=east,very thick,ForestGreen](0,0)--(-\tkzwd,0);
\draw[path fading=east,very thick,ForestGreen] (0,\tkzht)--(-\tkzwd,\tkzht);
\draw[dotted,red] (0,0)--(0,\tkzht);
\fill[purple](0,0)circle(1.5pt);
\fill[purple](0,\tkzht)circle(1.5pt);
\end{tkz}
\qquad\quad=\qquad
\sum_{i\in I}\;
\Gamma\!\left(\;
\begin{tkz}
\def\tkzht{1.5}
\def\tkzwd{1}
\def\tkzblg{1.4}
\fill[path fading=west,gray!25](-\tkzwd,0)rectangle(0,\tkzht);
\draw[path fading=west,very thick,ForestGreen](0,0)--(-\tkzwd,0);
\draw[path fading=west,very thick,ForestGreen] (0,\tkzht)--(-\tkzwd,\tkzht);
\fill[purple](0,0)circle(1.5pt);
\fill[purple](0,\tkzht)circle(1.5pt);
\begin{scope}
\clip(0,0)to[bend right=80,distance=\tkzblg cm](0,\tkzht)--cycle;
\fill[brown!15](0,0) rectangle (5,\tkzht);
\node[brown] at ({0.3*\tkzblg}, {0.5*\tkzht}){\scriptsize$p_i$};
\end{scope}
\draw[very thick,red] (0,0)--
(0,\tkzht);
\draw[red,very thick](0,0)to[bend right=80,distance=\tkzblg cm](0,\tkzht);
\end{tkz}
\!\!\!\!
\begin{tkz}[xscale=-1]
\def\tkzht{1.5}
\def\tkzwd{1}
\def\tkzblg{1.4}
\fill[path fading=east,gray!25](-\tkzwd,0)rectangle(0,\tkzht);
\draw[path fading=east,very thick,ForestGreen](0,0)--(-\tkzwd,0);
\draw[path fading=east,very thick,ForestGreen] (0,\tkzht)--(-\tkzwd,\tkzht);
\fill[purple](0,0)circle(1.5pt);
\fill[purple](0,\tkzht)circle(1.5pt);
\begin{scope}
\clip(0,0)to[bend right=80,distance=\tkzblg cm](0,\tkzht)--cycle;
\fill[brown!15](0,0) rectangle (5,\tkzht);
\node[brown] at ({0.3*\tkzblg}, {0.5*\tkzht}){\scriptsize$p_i$};
\end{scope}
\draw[very thick,red] (0,0)--
(0,\tkzht);
\draw[red,very thick](0,0)to[bend right=80,distance=\tkzblg cm](0,\tkzht);
\end{tkz}
\;\right).
\]
\end{corollary}

\begin{pf}
Fix $r\in \undfld[n-1]{\C}{Y}[\eta]$. The reverse inclusion in \cref{eq:star-generalized} is immediate, so it suffices to show $\Gamma(V_r)\subset\sum_{i\in  I}\Gamma(p_iV_{r_i}p_i)$. As $\Sk_\C(Y,\eta)$ is unitary, $\End_{\Sk_\C(Y,\eta)}(r)$ is a finite-dimensional $\Cstar$-algebra, so $\End_{\Sk_\C(Y,\eta)}(r)$ admits some complete set of minimal orthogonal projections $\{q_j\}_{j=1}^{m}$. Fix $1\leq j\leq m$. Since $q_j$ is minimal, $\End_{\Sk_\C(Y,\eta)^\cent}((r,q_j))\defeq{} q_j\End_{\Sk_\C(Y,\eta)}(r)q_j=\bC q_j$, so $(r,q_j)\in\Irr(\Sk_\C(Y,\eta)^\cent)$, i.e., $(r,q_j)\cong(r_{i_j},p_{i_j})$ in $\Sk_\C(Y,\eta)^\cent$ for some $i_j\in I$. Thus for every $v\in V_r$, 
\[
    \Gamma(v)
    \underset{\eqref{glu-compress}}{=}
    \sum_{j=1}^m\Gamma(\fcj{q_j}\blt_{\orev Y}v\blt_Yq_j)
    \;\in\;
    \sum_{j=1}^m\Gamma(q_jV_rq_j)
    \underset{\eqref{glu-corners}}{=}
    \sum_{j=1}^m\Gamma(p_{i_j}V_{r_{i_j}}p_{i_j})
    \subset
    \sum_{i\in I}\Gamma(p_iV_{r_i}p_i).\qedhere
\]
\end{pf}

\subsection{Finiteness and orthogonality from gluing}
\label{subsec:finiteness-self-gluing}

\

\begin{definition}[Path integral data]
\label{def:path-integral-data}
Let $\C$ be a disk-like $n$-category and let $\cG=(X,Y,R,\xi,\eta)$ be gluing data for $\C$. \defn{Path integral data} $\psi^\cG=(\psi_{Y\times I},\psi_X,\psi_{X_\glu})$ for $\cG$ consists of three linear functionals
\begin{align*}
\label{eq:pidata}
  \psi_{Y\times I}\colon&\undfldl[n]{\C}{\partial((Y\times I)\times I)}\longrightarrow\bC,
  \\
  \psi_X\colon&\undfldl[n]{\C}{\partial(X\times I)}\longrightarrow\bC,
  \\
  \psi_{X_\glu}\colon&\undfldl[n]{\C}{\partial(X_\glu\times I)}\longrightarrow\bC,
\end{align*}
(where here $-\times I$ means pinched so that $\partial(Y\times I)=\orev{Y}\cup_{\partial Y} Y$) whose induced pairings $\bkt fg_{Y\times I,\fcj r\amalg r'}\coloneq\psi_{Y\times I}(\fcj f\blt g)$, $\bkt vw_{X,\xi\cup(\fcj r\amalg r')}\coloneq\psi_X(\fcj v\blt w)$, and $\bkt ab_{X_\glu,\xi_\glu}\coloneq\psi_{X_{\glu}}(\fcj a\blt b)$ satisfy the following two conditions.
\begin{lst}[font=\upshape,leftmargin=0.6375in]
  \item[($\psi^\cG$Pos)]\label{ZP}
  For all $r,r'\in\undfld[n-1]{\C}{Y}[\eta]$, $\Sk_\C(Y,\eta)(r\to r')$ is finite-dimensional and $\bkt--_{Y\times I,\fcj r\amalg r'}$ is positive-definite on it.
  \item[($\psi^\cG$Glu)]\label{ZG}
  For all $r,r'\in\undfld[n-1]{\C}{Y}[\eta]$, $v\in V_r$, and $w\in V_{r'}$, \ref{Glu} holds in that
\end{lst}
\begin{equation}
\label{eq:pi-glu}
\bkt{\Gamma(v)}{\Gamma(w)}_{X_\glu,\xi_\glu}
=
\Omega_{\Sk_\C(Y,\eta)(r\to r')}\triangleright\psi_X\big((-)\blt(\fcj v\blt w)\big).
\end{equation}

We call $\psi^\cG$ (\defn{finite}) \defn{unitary} if the disk-like skein 1-category $\A_\C(Y,\eta)$ is (finite) unitary with sphere trace $\psi_{Y\times I}$ (see \cref{lem:reparam-gen}). 

For a linear functional $\psn\colon\undfldl[n]{\C}{S^n}\to\bC$, we will say that $\psi^\cG$ is $\psn$-\defn{normalized} if for any $(n,0)$-handlebody $X=\coprod_{i=1}^{m_0}X_i$ for $n$-balls $X_i$ and any $c\in\undfld[n-1]{\C}{\partial X}$, under the canonical isomorphism $\undfldl[n]{\C}{X}[c]\cong\bigotimes_{i=1}^{m_0}\undfldl[n]{\C}{X_i}[c|_{\partial X_i}]$ we have for simple tensors $\alpha=\alpha_1\otimes\cdots\otimes\alpha_{m_0}$ and $\beta=\beta_1\otimes\cdots\otimes\beta_{m_0}$ that $\bkt\alpha\beta_{X,c}=\prod_{i=1}^{m_0}\bkt{\alpha_i}{\beta_i}_{X_i,c|_{\partial X_i}}$, where the pairings $\bkt--_{X_i,c|_{\partial X_i}}$ are the $\psn$-pairings of \eqref{disk pairings}.
\end{definition}

The condition \ref{ZG} in \cref{def:path-integral-data} makes sense: each $\Sk_\C(Y,\eta)(r\to r')$ is finite-dimensional (by \ref{ZP}) and $\bkt--_{Y\times I,\fcj{r}\amalg r'}$ is a positive-definite sesquilinear form on it by \ref{ZP}, so the canonical element $\Omega_{\Sk_\C(Y,\eta)(r\to r')}$ exists.

\begin{remark}
\label{rmkm}
    By \cref{sec:1+1D}, path integral data $\psi^\cG$ is (finite) unitary if and only if the skein 1-category $\Sk_\C(Y,\eta)$ equipped with the trace $\Tr^{\Sk_\C(Y,\eta)}$ from \cref{cstr:pre-2-Hilbert spaces from unitary disk-like 1-categories} is a (finite) pre-2-Hilbert space (\cref{def:pre-2-Hilbert space}), and thus in particular is (finite) presemisimple.
\end{remark}

The following lemma is claimed in \cite[6.2.5]{W06} (for finite $I$). Together with \eqref{eq:star-generalized}, it will imply the dimension formula in \cref{glu-basis} below. 

\begin{lemma}[{\cite[6.2.5]{W06}}]
\label{glu-ip}
Let $\C$ be a disk-like $n$-category, let $\cG=(X,Y,R,\xi,\eta)$ be gluing data for $\C$, and let $\psi^\cG$ be unitary path integral data for $\cG$, with $\Irr(\Sk_\C(Y,\eta)^\cent)=\{(r_i,p_i)\}_{i\in I}$ (possibly infinite). Then for all $v_i\in p_iV_{r_i}p_i$ and $w_j\in p_jV_{r_j}p_j$, 
\begin{equation}
\label{ipeqn}
  \bkt{\Gamma(v_i)}{\Gamma(w_j)}_{X_{\glu},\xi_{\glu}}
  =
  \frac{\delta_{i=j}}{\bkt{p_i}{p_i}_{Y\times I,\fcj{r_i}\amalg r_i}}
  \bkt{v_i}{w_j}_{X,\xi\cup(\fcj{r_i}\amalg r_i)}.
\end{equation}
\end{lemma}

\begin{pf}
Fix $i,j\in I$. We first observe that for any $g\in\Sk_\C(Y,\eta)(r_i\to r_j)$,
\begin{equation}\label{(!)}
(\fcj{g}\amalg g) \bullet (\fcj{v_i}\bullet w_j)
=
(\fcj{g}\amalg g)
\bullet
(\fcj{p_i v_i p_i} \bullet p_j w_j p_j)
=
(\fcj{p_j g p_i}\amalg p_j g p_i) \bullet (\fcj{v_i}\bullet w_j).
\end{equation}
Indeed, both
\[
(\fcj{g}\amalg g)\blt(\fcj{p_i v_i p_i} \blt p_j w_j p_j)
\quad=\qquad
\makediagL[scale=1]{0}{1.5}
\quad
\makediagR[scale=1]{0}{1.5}
\]
and
\[
(\fcj{p_j g p_i}\amalg p_j g p_i)\blt(\fcj{v_i}\blt w_j)
\quad=\qquad
\makediagL[scale=1]{1.25}{0}
\quad
\makediagR[scale=1]{1.25}{0}
\]
are representatives of
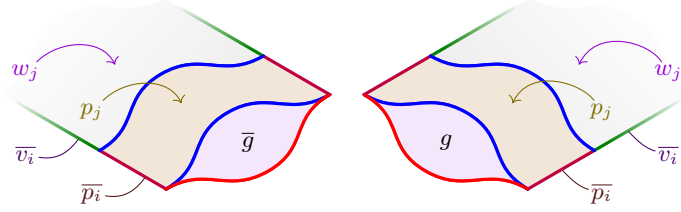
\[
\makediagL[scale=1.25]{0}{0}
\quad
\makediagR[scale=1.25]{0}{0}
\]
in the colimit. It follows that for any orthogonal basis $\{g_m\}$ of $\Sk_\C(Y,\eta)(r_i\to r_j)$,
\begin{equation}
\label{shared-start}
\begin{multlined}
\bkt{\Gamma(v_i)}{\Gamma(w_j)}_{X_{\glu}, \xi_{\glu}}
\;\defeq{}\;
\psi_{X_{\glu}}( \fcj{\Gamma(v_i)} \bullet \Gamma(w_j))
\\
\qquad\qquad\qquad
\underset{\ref{axiom:refl-gluing}\&\ref{CGa}}{=}
\psi_{X_{\glu}}( \Gamma(\fcj{v_i} \bullet w_j))
\underset{\ref{ZG}}{=}
\sum_m \frac{ \psi_{X}((\fcj{g_m}\amalg g_m) \bullet (\fcj{v_i}\bullet w_j )) }{ \bkt{g_m}{g_m}_{Y\times I,\fcj{r_i}\amalg r_j}}.
\end{multlined}
\end{equation}
Next observe that by semisimplicity of the completion $\Sk_\C(Y,\eta)^\cent$, Schur's Lemma \cite[Lem. 7.10.5]{UQSL} implies
\begin{equation}
\label{SchurApplication}
p_j\Sk_\C(Y,\eta)(r_i\to r_j)p_i
\defeq
\Sk_\C(Y,\eta)^\cent((r_i,p_i)\to (r_j,p_j))\cong\delta_{i=j}\bC.
\end{equation}
 
\itemstep{Case 1: $i=j$.} Since $p_i\neq0$ and the pairing on $\End_{\Sk_\C(Y,\eta)}(r_i)$ is positive-definite by unitarity, $\bkt{p_i}{p_i}_{Y\times I,\fcj{r_i}\amalg r_i}>0$, so we may choose an orthogonal basis $\{g_m\}$ of $\Sk_\C(Y,\eta)(r_i\to r_i)$ containing $p_i$. We claim
\begin{equation}
\label{corner}
p_i g_m p_i = \delta_{g_m=p_i}\, p_i\qquad\forall m.
\end{equation}
Indeed, since $p_i$ is minimal, $p_ig_mp_i=\lambda p_i$ for some $\lambda\in\bC$; if $g_m = p_i$, then $p_i p_i p_i = p_i=\delta_{g_m=p_i} p_i$, so $\lambda=1$, while if $g_m\neq p_i$, then $\bkt{g_m}{p_i}_{Y\times I,\fcj{r_i}\amalg r_i} = 0$ by orthogonality, so, using that $p_i$ is a projection and \cref{pairing-composition compatibility} (applied to $\A_\C(Y,\eta)$),
\[
\lambda\bkt{p_i}{p_i}=\bkt{p_i}{\lambda p_i}
=
\bkt{p_i}{p_i g_m p_i}
=
\bkt{p_i^\dag p_i p_i^\dag}{g_m}
=
\bkt{p_i p_i p_i}{g_m}
=
\bkt{p_i}{g_m} = 0,
\]
whence $\lambda=0$.
Next we observe that, using \eqref{(!)} with $g=g_m$ and $j=i$ (so $p_j=p_i$),
\begin{equation}
\label{ipcolimrep}
\begin{aligned}
(\fcj{g_m}\amalg g_m) \bullet (\fcj{v_i}\bullet w_j)
&
\underset{\eqref{(!)}}{=}
(\fcj{p_i g_m p_i}\amalg p_i g_m p_i) \bullet (\fcj{v_i}\bullet w_j )
\\&
\underset{\eqref{corner}}{=}
\delta_{g_m=p_i}\,
(\fcj{p_i}\amalg p_i) \bullet (\fcj{v_i}\bullet w_j )
\underset{\ref{CTc}}{=}\delta_{g_m=p_i}\,\fcj{v_i}\blt w_j,
\end{aligned}
\end{equation}
where the last equality uses $v_i=v_i\blt_Yp_i$ and $w_j=w_j\blt_Yp_i$. Thus
\[
\bkt{\Gamma(v_i)}{\Gamma(w_j)}_{X_{\glu}, \xi_{\glu}}
\underset{\eqref{shared-start}\&\eqref{ipcolimrep}}{=}
\frac{ \psi_{X}(\fcj{v_i}\bullet w_j) }{ \bkt{p_i}{p_i}_{Y\times I,\fcj{r_i}\amalg r_i}}
\;\;=\;\;
\frac{1}{\bkt{p_i}{p_i}_{Y\times I,\fcj{r_i}\amalg r_i}} \bkt{v_i}{w_j}_{X,\xi\cup(\fcj{r_i}\amalg r_i)}.
\]
 
\itemstep{Case 2: $i\neq j$.} Then by \eqref{SchurApplication} $p_j\Sk_\C(Y,\eta)(r_i\to r_j)p_i=0$, so $p_jgp_i=0$ for every $g\in\Sk_\C(Y,\eta)(r_i\to r_j)$. Thus for any orthogonal basis $\{g_m\}$ of $\Sk_\C(Y,\eta)(r_i\to r_j)$, we have $(\fcj{g_m}\amalg g_m) \bullet (\fcj{v_i}\bullet w_j)\underset{\eqref{(!)}}{=}(\fcj{p_j g_m p_i}\amalg p_j g_m p_i) \bullet (\fcj{v_i}\bullet w_j )=0$, so
\[
\bkt{\Gamma(v_i)}{\Gamma(w_j)}_{X_{\glu}, \xi_{\glu}}
\underset{\eqref{shared-start}}{=}
\sum_m \frac{ \psi_{X}((\fcj{g_m}\amalg g_m) \bullet (\fcj{v_i}\bullet w_j )) }{ \bkt{g_m}{g_m}_{Y\times I,\fcj{r_i}\amalg r_j}}
=0.
\]
This proves \eqref{ipeqn}.
\end{pf}

\begin{proposition}[Dimension formula]
\label{glu-basis}
Let $\C$ be a disk-like $n$-category, let $\cG=(X,Y,R,\xi,\eta)$ be gluing data for $\C$, and let $\psi^\cG$ be unitary path integral data for $\cG$, with $\Irr(\Sk_\C(Y,\eta)^\cent)=\{(r_i,p_i)\}_{i\in I}$ (possibly infinite). Further suppose that $\bkt--_{X,\xi\cup(\fcj{r_i}\amalg r_i)}$ is positive-definite on each $p_iV_{r_i}p_i$. Then
\begin{lst} 
\item $\bkt--_{X_{\glu},\xi_{\glu}}$ is positive-definite.
\item\label{stardirect} The sum \eqref{eq:star-generalized} is orthogonal-direct: $\undfldl[n]{\C}{X_\glu}[\xi_\glu]=\bigoplus_{i\in I}\Gamma(p_iV_{r_i}p_i)$.
\item
\label{dimformula}
We have $\dim\undfldl[n]{\C}{X_{\glu}}[\xi_{\glu}]=\sum_{i\in I}\dim(p_iV_{r_i}p_i)$.
\item\label{dfmbasis} If $\{f_{ij}\}_j$ is a $\bkt--_{X,\,\xi\cup(\fcj{r_i}\amalg r_i)}$-orthogonal basis of $p_iV_{r_i}p_i$ for each $i\in I$, then $\{\Gamma(f_{ij})\}_{i,j}$ is an orthogonal basis of $\undfldl[n]{\C}{X_{\glu}}[\xi_{\glu}]$.
\end{lst}
\end{proposition}

\begin{pf}
For each $i\in I$ let $\{f_{ij}\}_j$ be an orthogonal basis of the subspace $p_iV_{r_i}p_i\subset V_{r_i}$. By \eqref{ipeqn},
\[
\bkt{\Gamma(f_{ij})}{\Gamma(f_{kl})}_{X_\glu,\xi_\glu}
=
\frac{\delta_{i=k}\delta_{j=l}}{\bkt{p_i}{p_i}_{Y\times I,\fcj{r_i}\amalg r_i}}\bkt{f_{ij}}{f_{ij}}_{X,\xi\cup(\fcj{r_i}\amalg r_i)},
\]
which is nonzero when $(i,j)=(k,l)$ by definiteness of $\bkt--_{X,\xi\cup(\fcj{r_i}\amalg r_i)}$. Hence $\{\Gamma(f_{ij})\}_{ij}$ is an orthogonal set of nonzero vectors and so is linearly independent. Moreover, by \eqref{eq:star-generalized} $\{\Gamma(f_{ij})\}_{ij}$ spans $\undfldl[n]{\C}{X_{\glu}}[\xi_{\glu}]$ (by finite sums), giving \ref{dfmbasis}. The orthogonal-directness of the sum in \eqref{eq:star-generalized}, the positive-definiteness of $\bkt--_{X_\glu,\xi_\glu}$, and the dimension formula now immediately follow.
\end{pf}

The above tools allow for a quick proof of the following skein-theoretic result that Walker uses in \cite{W21}.

\begin{corollary}
\label{stadard-skthy}
Let $\C$ be a disk-like $n$-category and suppose we have $n$D gluing data $\cG=(X, Y, R, \xi, \eta)$ for $\C$ in which $X$ is a disjoint union $X = X_1 \amalg X_2$ for $n$-manifolds $X_1$ and $X_2$ and $X_\glu=X_1\cup_Y X_2$ such that $\partial X_1 = R_1 \cup_{\partial \orev Y} \orev{Y}$ and $\partial X_2 = R_2 \cup_{\partial Y} Y$, so that $R = R_1 \amalg R_2$ and $\xi = \xi_1 \amalg \xi_2$. 

Now suppose we have unitary path integral data $\psi^\cG$ for $\cG$ such that $\bkt--_{X,\xi\cup(\fcj{r_i}\amalg r_i)}$ is positive-definite on each $p_iV_{r_i}p_i$ (as in \cref{glu-basis}), so that $\Sk_\C(Y,\eta)$ is presemisimple with $\Irr(\Sk_\C(Y,\eta)^\cent)= \{(r_i,p_i)\}_{i\in I}$ (possibly infinite). Then the gluing map $\Gamma^\cG$ gives a canonical isometric isomorphism
\[
\undfldl[n]{\C}{X_{\glu}}[ \xi_{\glu}] \cong
\bigoplus_{i \in I}
(\undfldl[n]{\C}{X_1}[ \xi_1 \cup \fcj{r_i}]\blt_{\orev Y}\fcj{p_i})
\otimes_\bC
(p_i\blt_Y \undfldl[n]{\C}{X_2}[ r_i\cup \xi_2])
\]
where the right side is equipped with the orthogonal direct sum of the sesquilinear forms
\[
  \frac{1}{\bkt{p_i}{p_i}_{Y\times I,\,\fcj{r_i}\amalg r_i}}
  \bkt{-}{-}_{X_1,\,\xi_1\cup\fcj{r_i}}
  \cdot
  \bkt{-}{-}_{X_2,\,r_i\cup\xi_2}
\]
on the $i$-th summand.
\end{corollary}

\begin{pf}
Applying the $p_i$ gives $p_i V_{r_i} p_i \cong (\fcj{p_i} \blt_{\orev Y} \undfldl[n]{\C}{X_1}[ \xi_1 \cup \fcj{r_i}]) \otimes_\bC (\undfldl[n]{\C}{X_2}[r_i \cup \xi_2] \blt_Y p_i)$, which by \cref{glu-basis}\ref{stardirect} gives the result.
\end{pf}

\subsection{The Walker Extension Theorem preserves reflection positivity}

We first recall the Walker Extension Theorem from \cite{W06,W21}, which we state in the disk-like framework.

\begin{restatable}[Walker Extension Theorem {{\cite{W21}}}]{theorem}{WalkerExtensionTheorem}
    \label{WalkerTheorem}
    Suppose a disk-like $n$-category $\C$ with associated $(n+\e)$D TQFT $\C$ and a linear functional $\psi\colon
    \undfldl[n]{\C}{S^n}\to\bC$ satisfy the following conditions.
    \begin{lst}
        \item[(ND)]\label{ND} The functional $\psi$ induces nondegenerate pairings on $\fldl{\C}{D^n}[c]$ (via \eqref{disk pairings}) for all fields $c\in\undfld[n-1]{\C}{\partial D^n}$.
        \item[(0fSS)]\label{0fSS} The skein module $\undfldl{\C}{S^j\times D^{n-j}}[c]$ is finite-dimensional for all $0\leq j\leq n$ and all fields $c$ on $\bdy(S^j\times D^{n-j})$.
        \item[(1fSS)]\label{1fSS} The skein 1-category $\Sk_\C(S^j\times D^{n-j-1},d)\coloneq\cX_{\A_\C(S^{j}\times D^{n-j-1},d)}$ defined in \cref{skein 1-categories} is finite presemisimple for all $0\leq j\leq n-1$ and all $d\in\undfld[n-2]{\C}{\bdy(S^j\times D^{n-j-1})}$.
    \end{lst}
    Then there exists a unique path integral $\Z_\C$ for $\undC$ such that $\Z_\C(D^{n+1})=\psi$.
\end{restatable}

\begin{remark}\label{redundantIH3}
As claimed in \cite{W06}, if \ref{ND} is strengthened to require the pairings \eqref{disk pairings} to be positive-definite, then \ref{1fSS} is implied; this follows from the inductive proof below.
\end{remark}

An \defn{$(n+1,k)$-handle} is an $(n+1)$-ball $h$ together with a homeomorphism $\varphi_h\colon h\overset\sim\to D^{k}\times D^{n+1-k}$. We call $X_h\coloneq \varphi_h^{-1}(S^{k-1}\times D^{n+1-k})$ the \defn{attaching region} of $h$.

\begin{definition}
Let $M$ be an $(n+1)$-manifold. A \defn{handle decomposition} $\cH$ of $M$ is a filtration
\[
  \cH=\big(\varnothing = M_{-1} \subset M_0 \subset M_1 \subset \cdots \subset M_{n+1} = M\big)
\]
such that each $M_j$ is obtained from $M_{j-1}$ by attaching a \emph{finite} collection $\cH_j$ of $(n+1,j)$-handles. Each $h\in\cH_j$ with $h\cong D^j\times D^{n+1-j}$ is attached along an embedding of its attaching region $X_h=\varphi_h^{-1}(\partial D^j\times D^{n+1-j})\subset\partial h$ into $\partial M_{j-1}$, whose image is disjoint from those of the attaching regions of all other handles in $\cH_j$. (Thus $M_0$ is a disjoint union of $(n+1)$-balls.)
\end{definition}

Recall that an \defn{$(n+1,k)$-handlebody} is an $(n+1)$-manifold equipped with a handle decomposition where all handle indices are at most $k$. Since two manifolds are homeomorphic if and only if they admit handle decompositions that are related by a finite sequence of handle slides and handle cancellations, to prove \ref{Inv} it suffices to show the path integral is invariant under these operations. In \cite{W06,W21}, Walker constructs the path integral by handle attachments (\cref{cstr:path-integral-via-handle-decomp} below) and shows that this construction is well-defined, unique, and satisfies \ref{Inv} and \ref{Glu}. See also \cite{Hai}.

What remains for us is \ref{Rfl} and \ref{Pos}, which we verify simultaneously by going through Walker's construction and proving the following two statements by induction on the handle index $j$.
\begin{lst} 
\item[($\mathrm{Rfl}_j$)] $\Z_\C(\orev{M})(\fcj\psi)=\overline{\Z_\C(M)(\psi)}$ for every $(n+1,j)$-handlebody $M$ and all $\psi\in\undfldl[n]{\C}{\partial M}$.

\item[($\mathrm{Pos}_j$)] $\bkt{-}{-}_{Z,c}$ is positive-definite for every $(n,j)$-handlebody $Z$ and all $c\in\undfld[n-1]{\C}{\partial Z}$.
\end{lst}

\begin{construction}[Path integral via handle decomposition {\cite{W06,W21}}]
\label{cstr:path-integral-via-handle-decomp}
Let $\C$ be a finite unitary disk-like $n$-category and let $M$ be an $(n+1)$-manifold. Fix a handle decomposition
\(
    \cH=(\varnothing=M_{-1}\subset M_0\subset\cdots\subset M_{n+1}=M)
\)
and write $\cH_j=\{h^{(1)}_j,\dots,h^{(m_j)}_j\}$ for the collection of $(n+1,j)$-handles in $\cH$ and denote the attaching region of $h^{(k)}_j$ by $X^{(k)}_j\coloneq\varphi_{h^{(k)}_j}^{-1}(\partial D^j\times D^{n+1-j})\cong S^{j-1}\times D^{n+1-j}$. We construct the path integral on $M$ inductively over the handle index $j$ of $\cH$ as follows.

\itemstep{Base case: $0$-handles.} 
For the disjoint union $M_0=\amalg_{k=1}^{m_0}h_0^{(k)}$, under the canonical identifications
\[
\undfldl[n]{\C}{\partial({\amalg}_{k=1}^{m_0}h_0^{(k)})}
\underset{\eqref{lem:reparam-gen}}=
\undfldl[n]{\C}{{\amalg}_{i=1}^{m_0}S^n}
=
\bigotimes_{i=1}^{m_0}\undfldl[n]{\C}{S^n},
\]
we define $\Z_\C(M_0)\colon\undfldl[n]{\C}{\partial M_0}\to\bC$ on simple tensors $\alpha=\alpha_1\otimes\cdots\otimes\alpha_{m_0}$ by $\Z_\C(M_0)(\alpha)\coloneq\prod_{i=1}^{m_0}\psn^\C(\alpha_i)$. Now $(\mathrm{Rfl}_0)$ reduces to $\psn^\C(\fcj{\alpha_i})=\overline{\psn^\C(\alpha_i)}$ for each $i$, which follows from \ref{psnR}, while $(\mathrm{Pos}_0)$ is just \ref{psnP}.

\itemstep{Inductive step: $j$-handles for $1\leq j\leq n+1$.}
Suppose $\Z_\C(M)$ has been defined for all $(n+1,j-1)$-handlebodies and satisfies $(\mathrm{Pos}_{j-1})$ and $(\mathrm{Rfl}_{j-1})$. We construct $\Z_\C(M_j^{(k)})$ inductively over $1\leq k\leq m_j$ as follows. Define $M_j^{(0)}\coloneq M_{j-1}$ and set $M_j^{(k)}\coloneq M_j^{(k-1)}\cup_{X_j^{(k)}} h_j^{(k)}$ for each $h_j^{(k)}\in \cH_j$. 
\begin{itemize} 
\item Under the canonical identifications
\[
\undfldl[n]{\C}{\partial(M_j^{(k-1)}\amalg h_j^{(k)})}=\undfldl[n]{\C}{\partial M_j^{(k-1)}}\otimes\undfldl[n]{\C}{\partial h_j^{(k)}}
\underset{\eqref{lem:reparam-gen}}=
\undfldl[n]{\C}{\partial M_j^{(k-1)}}\otimes\undfldl[n]{\C}{S^n},
\]
we define $\Z_\C(M_j^{(k-1)}\amalg h_j^{(k)})$ on simple tensors $\alpha=\alpha_1\otimes\alpha_2$ by
\begin{equation}\label{istp1}
\Z_\C(M_j^{(k-1)}\amalg h_j^{(k)})(\alpha)\coloneq\Z_\C(M_j^{(k-1)})(\alpha_1)\cdot\psn^\C(\alpha_2).
\end{equation}
\item Let $R_1\coloneq\partial M_j^{(k-1)}\setminus X_j^{(k)}$ and $R_2\coloneq\partial h_j^{(k)}\setminus\orev{X_j^{(k)}}$.  Under the canonical identifications
\[
\undfldl[n]{\C}{\partial M_j^{(k)}}
=\undfldl[n]{\C}{R_1\cup_{\partial X_j^{(k)}} R_2}
=\bigoplus\nolimits_{c\in\undfld[n-1]{\C}{\partial X_j^{(k)}}} \undfldl[n]{\C}{R_1}[ c]\otimes\undfldl[n]{\C}{R_2}[ \fcj c],
\]
we have the following. Since $X_j^{(k)}\cong S^{j-1}\times D^{n+1-j}$ is an
$(n,j-1)$-handlebody, $(\mathrm{Pos}_{j-1})$ gives
positive-definiteness of $\bkt--_{X_j^{(k)},c}$, and moreover finiteness of $\C$ gives that $\undfldl[n]{\C}{X_j^{(k)}}[c]$ is finite-dimensional, so the canonical element $\Omega_{X^{(k)}_j,\,c}\coloneq\Omega_{\undfldl[n]{\C}{X_j^{(k)}}[c]}$ exists. We \emph{define} $\Z_\C(M_j^{(k)})\colon\undfldl[n]{\C}{\partial M_j^{(k)}}\to\bC$ directly by \eqref{eq:pi-glu}: for a fixed $c\in\undfld[n-1]{\C}{\partial X_j^{(k)}}$ and simple tensors $\alpha_1\otimes\alpha_2 \in \undfldl[n]{\C}{R_1}[ c]\otimes_\bC\undfldl[n]{\C}{R_2}[ \fcj c]$ with $\alpha_\glu =\alpha_1\blt_{\partial X_j^{(k)}}\alpha_2$, set
\begin{align*}
\Z_\C(M_j^{(k)})(\alpha_\glu)
&\coloneq
\Omega_{\undfldl[n]{\C}{X_j^{(k)}}[c]}\triangleright\Z_\C(M_j^{(k-1)}\amalg h_j^{(k)})((-)\blt(\alpha_1\amalg\alpha_2)) \\
&\underset{\eqref{istp1}}=
\Omega_{\undfldl[n]{\C}{X_j^{(k)}}[c]}\triangleright(\Z_\C(M_j^{(k-1)})((-)\blt\alpha_1)\cdot\psn^\C((-)\blt\alpha_2)).
\end{align*}
Then \ref{ZG} holds by construction.
\end{itemize}
Now set $\Z_\C(M_j)\coloneq \Z_\C(M_j^{(m_j)})$. It remains to verify $(\mathrm{Rfl}_j)$ and $(\mathrm{Pos}_j)$.

\itemstep{$(\mathrm{Rfl}_j)$.}
Note that there are two nested inductions in play: the outer induction on the handle index $j$, whose hypotheses are $(\mathrm{Rfl}_{j-1})$ and $(\mathrm{Pos}_{j-1})$, and an inner induction on the number $0\leq k\leq m_j$ of $j$-handles attached so far, whose hypothesis---denoted $(\mathrm{Rfl}_j^{k-1})$---is that $(\mathrm{Rfl}_j)$ holds for $M_j^{(k-1)}$.

We argue by induction on $0\leq k\leq m_j$. The base case $k=0$ follows from $(\mathrm{Rfl}_{j-1})$. Now let $1\leq k\leq m_j$ and assume $(\mathrm{Rfl}_j^{k-1})$. Let $\{e_\alpha\}$ be an orthogonal basis of $\undfldl[n]{\C}{X_j^{(k)}}[c]$. Since $X_j^{(k)}\times I$ is an $(n+1,j-1)$-handlebody and $\orev{X_j^{(k)}\times I}=\orev{X_j^{(k)}}\times I$, we have
\[
\bkt{\fcj{e_\alpha}}{\fcj{e_\beta}}_{\orev{X_j^{(k)}},\fcj c}
\underset{\ref{axiom:refl-gluing}\&(\mathrm{Rfl}_{j-1})}=
\overline{\bkt{e_\beta}{e_\alpha}_{X_j^{(k)},\,c}},
\]
so $\{\fcj{e_\alpha}\}$ is a $\bkt--_{\orev{X_j^{(k)}},\,\fcj c}$-orthogonal basis of $\undfldl[n]{\C}{\orev{X_j^{(k)}}}[\fcj c]$. Moreover,
\[
\|\fcj{e_\alpha}\|^2\underset{(\mathrm{Rfl}_{j-1})}{=}\overline{\|e_\alpha\|^2}\underset{(\mathrm{Pos}_{j-1})}{=}\|e_\alpha\|^2>0.
\]
Thus
\[
\begin{aligned}
\Z_\C(\orev{M_j^{(k)}})(\fcj\psi)
&=
\sum_\alpha
\frac{\Z_\C(\orev{M_j^{(k-1)}\amalg h_j^{(k)}})
(\fcj\psi\blt\fcj{e_\alpha}\blt e_\alpha)}
{\|e_\alpha\|^2}
\\
&=
\sum_\alpha
\frac{\overline{\Z_\C(M_j^{(k-1)}\amalg h_j^{(k)})
(\psi\blt e_\alpha\blt\fcj{e_\alpha})}}
{\|e_\alpha\|^2}
=
\overline{\Z_\C(M_j^{(k)})(\psi)}
\end{aligned}
\]
where the first equality is \eqref{eq:gluing-formula} applied to the orthogonal basis $\{\fcj{e_\alpha}\}$ of $\undfldl[n]{\C}{\orev{X_j^{(k)}}}[\fcj c]$, and the second equality uses the inner induction hypothesis $(\mathrm{Rfl}_j^{k-1})$, \ref{axiom:refl-gluing}, and
\ref{psnR}. 

\itemstep{$(\mathrm{Pos}_j)$.}
Let $Z$ be an $(n,j)$-handlebody. We induct on the number $m$ of $j$-handles of $Z$. If $m=0$, then $Z$ is an $(n,j-1)$-handlebody, so our claim follows from $(\mathrm{Pos}_{j-1})$. Now let $m\geq 1$, suppose the claim holds for $(n,j)$-handlebodies with $m-1$ $j$-handles, and write $Z=Z'\cup_Yh$ where $h$ is the final $j$-handle of $Z$ and $Z'$ is the $(n,j)$-handlebody obtained from $Z$ by omitting it. Since $Y\times D^1\cong S^{j-1}\times D^{n-j+1}$ is an $(n,j-1)$-handlebody, $(\mathrm{Pos}_{j-1})$ gives positive-definiteness of $\bkt--_{Y\times D^1,c}$. The pairing on $Z'\amalg h$ is positive-definite by the induction hypothesis and \ref{psnP}, so we obtain finite unitary path integral data for the resulting gluing data. Now \cref{glu-basis} gives
positive-definiteness of $\bkt--_{Z,c}$.

\itemstep{Final definition.}
Finally, we define $\Z_\C(M)\coloneq\Z_\C(M_{n+1})$.
\end{construction}

We have now proven the following corollary, which is the unitary version of the Walker Extension Theorem (\cref{WalkerTheorem}).

\UnitaryWalkerTheorem*

We further observe that the above proof of \cref{UnitaryWalkerTheorem} immediately implies the following.

\begin{corollary}
\label{cor:psnF-weaker}
If $\C$ is unitary, then $\C$ is finite if and only if $\undfldl[n]{\C}{S^j\times D^{n-j}}[\xi]$ is finite-dimensional for all $0\leq j\leq n$ and all boundary conditions $\xi\in\undfld[n-1]{\C}{\partial(S^j\times D^{n-j})}$.
\end{corollary}

The following is a well-known fact, which we recover from Walker's construction (\cref{cstr:path-integral-via-handle-decomp}).

\begin{proposition}
Let $\C$ be a finite unitary disk-like $n$-category and let $\lambda \in \bR_{>0}$. Let $\lambda\C$ denote the finite unitary disk-like $n$-category $\C$ equipped with the scaled sphere trace $\psn^{\lambda\C} \coloneq \lambda \psn^\C$. For any $(n+1)$-manifold $M$, the path integral satisfies
\[
\Z_{\lambda\C}(M) = \lambda^{\chi(M)} \Z_\C(M).
\]
\end{proposition}

\begin{pf}
\itemstep{Base case.} Following \cref{cstr:path-integral-via-handle-decomp},
\[
\Z_{\lambda\C}(M_0) = \prod_{i=1}^{m_0} \psn^{\lambda\C} = \lambda^{m_0} \prod_{j=1}^{|\mathcal{H}_0|} \psn^{\C} = \lambda^{|\mathcal{H}_0|} \Z_{\C}(M_0) = \lambda^{\chi(M_0)} \Z_{\C}(M_0).
\]

\itemstep{Induction step.} Suppose $\Z_{\lambda\C}(M_j^{(k-1)}) = \lambda^{\chi(M_j^{(k-1)})} \Z_{\C}(M_j^{(k-1)})$. Since
\[
\bkt--^{\lambda\C}_{X_h, c} = \lambda^{\chi(S^{j-1} \times D^{n+1-j})} \langle - | - \rangle^{\C}_{X_h, c} = \lambda^{\chi(S^{j-1})} \langle - | - \rangle^{\C}_{X_h, c} ,
\]
we have
\[
\Omega^{\lambda\C}_{X_h, c} = \sum_i \frac{\ket{e_i} \otimes \bra{e_i}}{\bkt{e_i}{e_i}^{\lambda\C}_{X_h, c}} = \lambda^{-\chi(S^{j-1})} \Omega^{\C}_{X_h, c},
\]
so attaching a $j$-handle $h \in \mathcal{H}_j$ along the attaching region $X_h \cong S^{j-1} \times D^{n+1-j}$ gives
\begin{align*}
\Z_{\lambda\C}(M_j^{(k-1)} \cup_{X_h} h) &= \Omega^{\lambda\C}_{X_h, c} \triangleright (\Z_{\lambda\C}(M_j^{(k-1)}) \cdot \psn^{\lambda\C}) \\
&= \lambda^{-\chi(S^{j-1})} \Omega^{\C}_{X_h, c} \triangleright (\lambda^{\chi(M_j^{(k-1)})} \Z_{\C}(M_j^{(k-1)}) \cdot \lambda \psn^{\C}) \\
&= \underbrace{\lambda^{\chi(M_j^{(k-1)}) - \chi(S^{j-1}) + 1}}_{\lambda^{\chi(M_j^{(k)})}} \underbrace{\Omega^{\C}_{X_h, c} \triangleright (\Z_{\C}(M_j^{(k-1)}) \cdot \psn^{\C})}_{\Z_{\C}(M_j^{(k-1)} \cup_{X_h} h)} ,
\end{align*}
as claimed.
\end{pf}

\section{Unrestricted disk-like functors}
\label{sec:unnrestricted}
\begin{definition}[Parameterized splitting]
\label{def:parameterized-splitting-filtered}
Let $\C$ be a disk-like $n$-category and let $W$ be a $k$-manifold.
A \defn{parameterized splitting} $\Theta$ of a pure $k$-field $\sigma\in\undfld[k]{\C}{W}$ with respect to an (unframed) string diagram stratification $\Gamma\hookrightarrow W$ is a collection of homeomorphisms
\[
  \Theta=\{\theta_s\colon D^{k-j}\times \check{s}\overset\sim\to N_s \mid s\in\mathrm{strata}(\Gamma),\ \dim(s)=j,\ 0\leq j\leq k\}
\]
defined inductively on $0\leq j\leq k$, where $\check{s}\coloneq s\setminus K_{j-1}$ for a $j$-stratum $s$, $K_{-1}\coloneq\varnothing$, and for $1\leq j\leq k$,
\[
K_{j-1}\coloneq
K_{j-2}\cup\bigcup_{\substack{(j-1)\text{-strata }s\in\Gamma}}\Int(N_s).
\]
The above data are subject to the following conditions.
\begin{lst}
  \item[($\Theta$1)] $N_s\subset W\setminus K_{j-1}$ for each $j$-stratum $s\in\Gamma$.
  \item[($\Theta$2)]\label{axiom:Theta-split} $\{N_s\}_{s\in\mathrm{strata}(\Gamma)}$ is a permissible ball decomposition of $W$ along which $\sigma$ splits.
  \item[($\Theta$3)]\label{axiom:Theta-core} $\theta_s|_{\{0\}\times\check{s}}=\id_{\check{s}}$; when $j=k$ this forces $\theta_s=\id_{\check{s}}$ and $N_s=\check{s}$.
  \item[($\Theta$4)]\label{axiom:Theta-bdy-compat} $\theta_s(D^{k-j}\times(\check{s}\cap\partial W))=N_s\cap\partial W$ for each $j$-stratum $s$ of $\Gamma$.
\end{lst}
\end{definition}

We write $\mathrm{PS}(\sigma,\Gamma)$ for the collection of all parameterized splittings $\Theta$ with respect to $\Gamma$ and $\mathrm{PS}(\sigma)$ for the collection of $\Theta$ that are parameterized splittings with respect to some unframed string diagram stratification. For $\Theta\in\mathrm{PS}(\sigma)$, we denote by $\Gamma_\Theta$ the corresponding unframed string diagram stratification.

\begin{notation}
Since $\theta_s=\id_s$ for every top-dimensional stratum $s\in\Gamma_\Theta$, a parameterized splitting is determined by its restriction to strata of codimension $\geq 1$. We suppress the top-dimensional data and write $\mathrm{PS}(\sigma,\Gamma_\sigma)=\{\varnothing\}$ when $\Gamma_\sigma$ has no strata of positive codimension.
\end{notation}

The data of $\C$ act on parameterized splittings as follows.
\begin{lst}  \item[($\Theta\varphi$)]\label{axiom:Theta-phi} For all homeomorphisms $\varphi\colon W\to W'$,
  \[
    \varphi_*\Theta\coloneq\{\varphi_*\theta_s\colon D^{k-j}\times\varphi(\check{s})\overset\sim\to\varphi(N_s)\mid\theta_s\in\Theta\}\in\mathrm{PS}(\varphi_*\sigma,\varphi(\Gamma_\Theta)),
  \]
  where $\varphi_*\theta_s\coloneq\varphi\circ\theta_s\circ(\id_{D^{k-j}}\times(\varphi|_{\check{s}})^{-1})$.
  \item[($\Theta\partial$)]\label{axiom:Theta-bdy} For all $\Theta \in \mathrm{PS}(\sigma, \Gamma)$, $\partial\Theta \coloneq \Theta \cap \partial W \in \mathrm{PS}(\partial\sigma, \Gamma \cap \partial W)$ consists of the restrictions $\theta_s|_{D^{k-j} \times \check{r}}$ for each boundary $(j-1)$-stratum component $r\hookrightarrow s\cap \partial W$ of a $j$-stratum $s \in \Gamma$. Observe that by \ref{axiom:Theta-bdy-compat}, the image of this restriction is $N_r$.
  \item[($\Theta$G)]\label{axiom:Theta-G} For all $\Theta_1\in\mathrm{PS}(\sigma_1,\Gamma_1)$ and $\Theta_2\in\mathrm{PS}(\sigma_2,\Gamma_2)$ such that each of (i) $\sigma_1$ and $\sigma_2$, (ii) $\Gamma_1$ and $\Gamma_2$, and (iii) $\Theta_1$ and $\Theta_2$ are compatible along a $(k-1)$-ball $E$ (the latter in that $(\partial\Theta_1)|_E=(\partial\Theta_2)|_E$), we define $\Gamma_1\blt_E\Gamma_2$ by merging corresponding strata intersecting $E$ and taking the disjoint union of the rest, and
  \[
    \Theta_1\blt_E\Theta_2\coloneq\{\theta_s\in\Theta_i\mid s\cap E=\varnothing,i=1,2\}\cup\{\theta_s\cup_E\theta_{s'}\mid s\cap E=s'\cap E\neq\varnothing\},
  \]
  where for corresponding $j$-strata $s$ and $s'$ that agree at $E$, the homeomorphism $\theta_s\cup_E\theta_{s'}\colon D^{k-j}\times(\check{s}\cup_E\check{s}')\overset\sim\to N_{s\cup_E s'}$ is defined by
  \[
    (\theta_s\cup_E\theta_{s'})(x,t)\coloneq
    \begin{cases}
      \theta_s(x,t) & \text{if }t\in\check{s},\\
      \theta_{s'}(x,t) & \text{if }t\in\check{s}'.
    \end{cases}
  \]
  \item[($\Theta\orev{\,\cdot\,}$)]\label{axiom:Theta-refl} For all $\Theta\in\mathrm{PS}(\sigma,\Gamma)$, $\orev\Theta\coloneq\{\orev{\theta_s}\coloneq\orev{\,\cdot\,}\circ\theta_s\circ(\iota^{(k-j)}\times\id_{\check{s}})\colon D^{k-j}\times\check{s}\to\orev{N_s}\mid\theta_s\in\Theta\}\in\mathrm{PS}(\fcj\sigma,\orev\Gamma)$, where $\iota^{(m)}\colon D^m\to\orev{D^m}$ negates the last coordinate.
\end{lst}

Finally, for a pinched product map $\pi\colon E\to W$, we will say a parameterized splitting $\widetilde\Theta\in\mathrm{PS}(\pi^*\sigma,\pi^{-1}(\Gamma_\Theta))$ is \defn{compatible with $\Theta$ over $\pi$} if $\pi\circ\widetilde\theta_{\pi^{-1}(s)}=\theta_s\circ(\id_{D^{k-j}}\times\pi)$ for each $j$-stratum $s$ of $\Gamma_\Theta$.

\begin{example}
When $\C$ is a disk-like $n$-category of string diagrams and $\xi\in\fld[k]{\C}{W}$, a parameterized splitting $\Theta$ along the underlying string diagram stratification $\Gamma_\xi$ is just a choice of parameterized regular neighborhoods of the strata of $\xi$.
\end{example}

\begin{definition}
A \defn{system of stratifications} for $\C$ is an assignment for each $0\leq k\leq n$ to each $k$-manifold $W$ and pure $k$-field $\sigma\in\undfld[k]{\C}{W}$ of an unframed string diagram stratification $\Gamma_\sigma\hookrightarrow W$ subject to the following conditions.
\begin{lst}
  \item[($\Gamma\varphi$)]\label{axiom:Gamma-nat} $\Gamma_{\varphi_*\sigma}=\varphi(\Gamma_\sigma)$ for all homeomorphisms $\varphi\colon W\to W'$.
  \item[($\Gamma\partial$)]\label{axiom:Gamma-bdy} $\Gamma_{\partial\sigma}=\Gamma_\sigma\cap\partial W$.
  \item[($\Gamma$G)]\label{axiom:Gamma-gluing} $\Gamma_{\sigma_1\blt_E\sigma_2}=\Gamma_{\sigma_1}\blt_E\Gamma_{\sigma_2}$ for all $\sigma_1$ and $\sigma_2$ compatible along $E$.
  \item[($\Gamma\orev{\,\cdot\,}$)]\label{ur-axiom:Gamma-refl} $\Gamma_{\fcj\sigma}=\orev{\Gamma_\sigma}$.
  \item[($\Gamma\pi$)]\label{axiom:Gamma-products} $\Gamma_{\pi^*\sigma}=\pi^{-1}(\Gamma_\sigma)$ for all pinched product maps $\pi$.
\end{lst}
\end{definition}

\begin{definition}[Unrestricted disk-like functor]
\label{def:unrestricted-functor}
An \defn{unrestricted disk-like functor} $\eF\colon\C\to\D$ consists of the following data.
\begin{lst}
  \item[(F$\Gamma$)]\label{axiom:FGamma-data} A system of stratifications $\sigma\mapsto\Gamma_\sigma$ on $\C$.
  \item[(F$\Theta$)]\label{axiom:FTheta-data} To each $k$-manifold $W$ and each pure $k$-field $\sigma$ on $W$, a set map $\eF(\sigma,-)\colon\eF_{\mathrm{PS}}(\sigma)\to\undfld[k]{\D}{W}$ valued in pure $k$-fields, where $\eF_{\mathrm{PS}}(\sigma)\coloneq\mathrm{PS}(\sigma,\Gamma_\sigma)$.
\end{lst}
By \ref{axiom:Gamma-nat}--\ref{axiom:Gamma-products}, the actions \ref{axiom:Theta-phi}--\ref{axiom:Theta-refl} define maps $\varphi_*\colon\eF_{\mathrm{PS}}(\sigma)\to\eF_{\mathrm{PS}}(\varphi_*\sigma)$, $\partial\colon\eF_{\mathrm{PS}}(\sigma)\to\eF_{\mathrm{PS}}(\partial\sigma)$, $-\blt_E-\colon\eF_{\mathrm{PS}}(\sigma_1)\times_{\eF_{\mathrm{PS}}(\partial\sigma_1|_E)}\eF_{\mathrm{PS}}(\sigma_2)\to\eF_{\mathrm{PS}}(\sigma_1\blt_E\sigma_2)$, and $\orev{\,\cdot\,}\colon\eF_{\mathrm{PS}}(\sigma)\to\eF_{\mathrm{PS}}(\fcj\sigma)$. The above data \ref{axiom:FGamma-data} and \ref{axiom:FTheta-data} are subject to the following conditions. 
\begin{lst}
  \item[(F$\Theta\varphi$)]\label{axiom:FTheta-nat} $\eF(\varphi_*\sigma,\varphi_*\Theta)=\varphi_*\eF(\sigma,\Theta)$.
  \item[(F$\Theta\partial$)]\label{axiom:FTheta-bdy} $\eF(\partial\sigma,\partial\Theta)=\partial\eF(\sigma,\Theta)$.
  \item[(F$\Theta$G)]\label{axiom:FTheta-gluing} $\eF(\sigma_1\blt_E\sigma_2,\Theta_1\blt_E\Theta_2)=\eF(\sigma_1,\Theta_1)\blt_E\eF(\sigma_2,\Theta_2)$ for gluings along any $(k-1)$-manifold $E$, whenever each of (i) $\sigma_1$ and $\sigma_2$, (ii) $\Theta_1$ and $\Theta_2$, and (iii) $\eF(\sigma_1,\Theta_1)$ and $\eF(\sigma_2,\Theta_2)$ are compatible along $E$.
  \item[(F$\Theta\orev{\,\cdot\,}$)]\label{axiom:FTheta-refl} $\eF(\fcj\sigma,\orev\Theta)=\fcj{\eF(\sigma,\Theta)}$.
  \item[(F$\Theta\pi$)]\label{axiom:FTheta-products} $\eF(\pi^*\sigma,\widetilde\Theta)=\pi^*\eF(\sigma,\Theta)$ for all pinched product maps $\pi\colon E\to W$ and all $\widetilde\Theta\in\eF_{\mathrm{PS}}(\pi^*\sigma)$ that are compatible with $\Theta\in\eF_{\mathrm{PS}}(\sigma)$ over $\pi$.
  
  \item[(F$\Theta$U)]\label{axiom:FTheta-local-relations} If $X$ is an $n$-ball, $\sigma\in\undfld[n-1]{\C}{\partial X}$, $\Xi\in\eF_{\mathrm{PS}}(\sigma)$, and $\sum_i\lambda_i\alpha_i\in\fldlU{\C}{X}[\sigma]$, then $\sum_i\lambda_i\eF(\alpha_i,\Theta_i)\in\fldlU{\D}{X}[\eF(\sigma,\Xi)]$ for any $\Theta_i\in\eF_{\mathrm{PS}}(\alpha_i)$ with $\partial\Theta_i=\Xi$.
\end{lst}
\end{definition}

\begin{remark}
If $W$ is a $0$-ball and $a\in\fld[0]{\C}{W}$, then $\eF_{\mathrm{PS}}(a)$ is a singleton. Thus we can unambiguously write $\eF(a)\in\fld[0]{\D}{W}$ for the value of $\eF$ on a $0$-field $a$.
\end{remark}

\begin{remark}
\label{dldldl}
By applying \ref{axiom:FTheta-local-relations} to the case $\sum_i\lambda_i\alpha_i=\alpha-\alpha=0$, we get that $[\eF(\alpha,\Theta)]=[\eF(\alpha,\Theta^\prime)]$ (that is, $\eF(\alpha,\Theta)-\eF(\alpha,\Theta^\prime)\in\fldlU{\D}{X}[\eF(\sigma,\Xi)]$) for all $\Theta,\Theta^\prime\in\eF_{\mathrm{PS}}(\alpha)$ with $\partial\Theta=\partial\Theta^\prime$. Thus $\eF$ descends to a linear map $\widehat{\eF}^{\,\Xi}_{X,\sigma}\colon\fldl[n]{\C}{X}[\sigma]\to\fldl[n]{\D}{X}[\eF(\sigma,\Xi)]$ for all $\Xi\in\eF_{\mathrm{PS}}(\sigma)$.
\end{remark}

\begin{remark} 
\label{dlurl}
Any disk-like functor can be viewed as an unrestricted disk-like functor by taking the empty system of stratifications, that is, $\Gamma_\sigma \coloneq \varnothing$ for all fields $\sigma$. Indeed, the unique element of $\mathrm{PS}(\sigma, \varnothing)$ is the empty collection of homeomorphisms.
\end{remark}

\subsection{Unrestricted equivalences of disk-like \texorpdfstring{$n$}{n}-categories}
\label{sec:unrestricted-equivalences}
Let $\C$ and $\D$ be disk-like $n$-categories.

\begin{definition}[Weak equivalence]
\label{def:unrestricted-weak-equivalence}
An unrestricted disk-like functor $\eF\colon\C\to\D$ is a \defn{weak equivalence} if it satisfies the following two conditions.
\begin{lst}
  \item[(UWE$k$)]\phantomsection\label{WEurk}
  ($\eF$ is \defn{unitarily essentially surjective on $k$-fields}.)
  For each $0\leq k<n$, each $k$-ball $W$, each $c\in\undfld[k-1]{\C}{\partial W}$, each $\Xi\in\eF_{\mathrm{PS}}(c)$, and each $\eta\in\fld[k]{\D}{W}[\eF(c,\Xi)]$, there is a $k$-field $\xi\in\fld[k]{\C}{W}[c]$ and a parameterized splitting $\Theta\in\eF_{\mathrm{PS}}(\xi)$ with $\partial\Theta=\Xi$ such that $\eF(\xi,\Theta)\cong^\star\eta$ in $\D$.
  \item[(UWE$n$)]\phantomsection\label{WEurn}
  For all $n$-balls $X$, $c\in\undfld[n-1]{\C}{\partial X}$, and $\Xi\in\eF_{\mathrm{PS}}(c)$, the linear map
  \[
    \begin{aligned} 
    \widehat{\eF}^{\,\Xi}_{X,c}\colon\fldl[n]{\C}{X}[c]&\longrightarrow\fldl[n]{\D}{X}[\eF(c,\Xi)],
    \\
    [\alpha]&\longmapsto[\eF(\alpha,\Theta)]
    \end{aligned}
  \]
  induced by \ref{axiom:FTheta-local-relations} for any $\Theta\in\eF_{\mathrm{PS}}(\alpha)$ with $\partial\Theta=\Xi$ (and well-defined by \cref{dldldl}) is a linear isomorphism.
\end{lst}
We will say $\C$ and $\D$ are \defn{weakly equivalent}, written $\C\cong\D$, if there is a zig-zag of weak equivalences between $\C$ and $\D$.

If moreover $(\C,\psn^\C)$ and $(\D,\psn^\D)$ are unitary disk-like $n$-categories, we call a weak equivalence $\eF\colon\C\to\D$ \defn{isometric} if $\psn^\C([\xi])=\psn^\D([\eF(\xi,\Theta)])$ for every pure $n$-field $\xi$ on $S^n$ and every $\Theta\in\eF_{\mathrm{PS}}(\xi)$.
We will say $\C$ and $\D$ are \defn{isometrically weakly equivalent}, written $\C\cong^\star\D$, if there is a zig-zag of isometric weak equivalences between $\C$ and $\D$.
\end{definition}

\begin{remark}
\label{uresequvgen}
As in the restricted case, by \ref{axiom:FTheta-nat} it is equivalent to require \ref{WEurk} and \ref{WEurn} to hold on the standard balls $D^k$ and $D^n$ respectively. When a disk-like functor $\eF$ is viewed as an unrestricted disk-like functor (see \cref{dlurl}), $\eF_{\mathrm{PS}}(\sigma)$ is a singleton for every $\sigma$, so \ref{WEurk}--\ref{WEurn} specialize to \ref{WEgenk}--\ref{WEgenn}, except that \ref{WEurk} moreover fixes the lift $c$ of $\partial\eta$. Thus a weak equivalence in the sense of \cref{def:unrestricted-weak-equivalence} is one in the sense of \cref{def:wkequiv}.
\end{remark}

\section{Disk-like \texorpdfstring{$n$}{n}-subcategories and liftable \texorpdfstring{$n$}{n}-fields}
\label{sec:liftable}
In this section, we define a notion of a disk-like $n$-subcategory and prove that we may always assume the vertices of string diagrams on $n$-balls are labeled by $n$-fields that admit pure $n$-field representatives. 

\begin{definition}
\label{def:disk-like subcategory}
Let $\D$ be a disk-like $n$-category. A \defn{disk-like $n$-subcategory} $\C$ of $\D$ is a choice of subsets of pure fields $\fld[k]{\C}{W} \subset \fld[k]{\D}{W}$ for each $k$-ball $W$, with local relations $\fldlU{\C}{X}[ c] \coloneq \fldlU{\D}{X}[ c] \cap \bC\{\fld[n]{\C}{X}[ c]\}$ and with the rest
of the disk-like $n$-category structure inherited from $\D$, such that for all $0 \leq k \leq n$, $\C$ is closed under
\begin{lst}
   \item boundary maps, i.e., $\partial\xi$ is in $\C$ whenever $\xi$ is;
   \item actions of homeomorphisms, i.e., $\varphi_*\xi$ is in $\C$ whenever $\xi$ is;
   \item reflection, i.e., $\fcj{\xi}$ is in $\C$ whenever $\xi$ is;
   \item pullbacks, i.e., $\pi^*\xi$ is in $\C$ whenever $\xi$ is;
   \item\label{dlsubcat3} gluing, i.e., if $\xi, \zeta$ are in $\C$ and are compatible along $Y$, then $\xi\blt_Y\zeta$ is in $\C$; and
   \item\label{dlsubcat4} splittings, i.e., if $\sigma$ is in $\C$ and is splittable along $Y$, say with $\sigma=\xi \blt_Y \zeta$, then $\xi, \zeta$ are in $\C$.
\end{lst}
\end{definition}

\begin{lemma}
\label{def:disk-like subcategory-proof}
A disk-like $n$-subcategory is a disk-like $n$-category.
\end{lemma}

\begin{pf}
Each axiom asserting an equality is obtained by including into $\D$ and using that $\D$ is disk-like. The $\lU{\C}$ axioms follow from the definition of $\lU{\C}$ combined with the closure of $\C$ under gluing, homeomorphisms, and pullbacks.

Gluing injectivity \ref{CGi} follows from inclusion into $\D$, since if $\xi_1 \blt_Y \eta_1 = \xi_2 \blt_Y \eta_2$ in $\fld[k]{\C}{W} \subset \fld[k]{\D}{W}$, injectivity in $\D$ forces $(\xi_1,\eta_1)=(\xi_2,\eta_2)$. The gluing map lands in $\C$ by \ref{dlsubcat3}. Finally, splittability \ref{CS} holds in $\C$: if $\sigma \in \C$ is splittable in $\D$ as $\xi \bullet_Y \zeta$, then \ref{dlsubcat4} ensures that $\xi,\zeta \in \C$, and hence that $\sigma$ is also splittable in $\C$.
\end{pf}

\begin{definition}
\label{def:liftable}
  We will call an $n$-field $f \in \undfldl[n]{\C}{X}[ c]$ \defn{liftable} if $f = [\tilde{f}]$ for some
  pure $n$-field $\tilde{f} \in \undfld[n]{\C}{X}[ c]$.
\end{definition}

The following facts are immediate consequences of the linearization/quotienting procedure that creates $n$-fields from pure $n$-fields.

\begin{facts}
\label{facts:liftable-properties}
  \begin{lst}
    \item[(i)]\label{facts:liftable-properties-1} Liftable $n$-fields span the vector spaces of $n$-fields $\fldl[n]{\C}{X}[ c]$.
    \item[(ii)]\label{facts:liftable-properties-2} Liftability of $n$-fields is preserved under gluing, homeomorphism actions, pullbacks (in that product $n$-fields are liftable), and reflection.
  \end{lst}
\end{facts}

\begin{construction}
Let $\eX$ be a pivotal $n$-category equipped with a bar structure $\fcj{\,\cdot\,}$ in the sense of \cite{W21}. For the purposes of this article, $\eX$ is a pre-2-Hilbert space $\cX$ (resp. a proto-3-Hilbert space $\fX$). 

Define the \defn{liftable disk-like $n$-subcategory} of $\C^{\eX_\C}$, denoted $\widetilde{\C}^{\eX_\C}$, as follows. For each $0\leq k< n$, set 
$\fld[k]{\widetilde{\C}^{\eX_\C}}{W}\coloneq\fld[k]{\C^{\eX_\C}}{W}$, and set
\[
\fld[n]{\widetilde{\C}^{\eX_\C}}{X}[c]\coloneq\{\xi\in\fld[n]{\C^{\eX_\C}}{X}[c]\mid \text{each vertex label of }\xi\text{ is liftable}\},
\]
with gluing, boundary, products, homeomorphism actions, and reflection structure inherited from $\C^{\eX_\C}$, and with local relations $\fldlU{\widetilde{\C}^{\eX_\C}}{X}[c]\coloneq\fldlU{\C^{\eX_\C}}{X}[c]\cap\bC\{\fld[n]{{\widetilde{\C}^{\eX_\C}}}{X}[c]\}$. If moreover $\C^{\eX_\C}$ has a sphere trace $\psn$ for which it is unitary, then define a sphere trace on $\fldl[n]{\widetilde{\C}^{\eX_\C}}{S^n}$ by $\psn\circ\iota$, where $\iota\colon\widetilde{\C}^{\eX_\C}\hookrightarrow\C^{\eX_\C}$ is the inclusion.
\end{construction}

\begin{lemma}
\label{lem:tilde-inclusion}
$\widetilde{\C}^{\eX_\C}$ is a disk-like $n$-category, and the inclusion $\iota\colon{\widetilde{\C}^{\eX_\C}}\hookrightarrow\C^{\eX_\C}$ is a weak equivalence.
If moreover $\C^{\eX_\C}$ has a sphere trace $\psn$ for which it is unitary, then $\iota$ is isometric.
\end{lemma}

\begin{pf}
By \cref{facts:liftable-properties}\ref{facts:liftable-properties-2}, the conditions of \cref{def:disk-like subcategory} follow, so $\widetilde{\C}^{\eX_\C}$ is a disk-like $n$-subcategory of $\C^{\eX_\C}$, and hence by \cref{def:disk-like subcategory-proof} a disk-like $n$-category. Note that the conditions of \cref{def:disk-like subcategory} also immediately give that $\iota$ is a disk-like functor.

It remains to show $\iota$ is an isometric weak equivalence. As $\iota=\id$ on $k$-fields for $0\leq k< n$, we need only show $\iota$ is a linear isomorphism on the vector spaces of $n$-fields $\fldl[n]{\widetilde{\C}^{\eX_\C}}{X}[c]$ for every $n$-ball $X$ and boundary condition $c\in\undfld[n-1]{\C^{\eX_\C}}{\partial X}$. Injectivity is immediate from the definition of ${\widetilde{\C}^{\eX_\C}}_U$. For surjectivity, first observe that by \cref{lem:eval-quotient} (resp. \cref{eval-facts-2}\ref{fact:eval-ker2}), the elements $[\xi]$ for $\xi\in\fld[n]{\C^{\eX_\C}}{X}[c]$ span $\fldl{\C^{\eX_\C}}{X}[c]$. By local relations, such a $[\xi]$ equals $[\xi_f]$ for some cone string diagram whose unique vertex is labeled by $f\in\eX_\C(a\to b)$. By \cref{facts:liftable-properties}\ref{facts:liftable-properties-1}, $f=\sum_k\alpha_k f_k$ for some $\alpha_k\in\bC$ and some liftable $n$-morphisms $f_k$ in $\eX_\C$. Thus $[\xi]=\sum_k\alpha_k[\xi_{f_k}]$ for some $\xi_{f_k}\in \fld[n]{\widetilde{\C}^{\eX_\C}}{X}[c]$. Note that injectivity of $\iota$ on all $n$-manifolds shows $\widetilde{\C}^{\eX_\C}$ satisfies \ref{psnF}, while \ref{psnP} is inherited from $\C^{\eX_\C}$ by construction. 
\end{pf}

\section{Strictification of pivotal dagger 2-categories}
\label{strictification-sec}
\tikzset{
  band/.style  ={draw=green!30, line width=16pt,  line join=round},
  wire/.style  ={draw=blue,     line width=1.2pt, line join=round},
  ribbon/.style={band, postaction={wire}},
  frame/.style ={draw=red,line width=1.5pt,line join=round},
  wire label/.style={font=\scriptsize, text=blue, inner sep=1pt},
  val label/.style ={font=\scriptsize, fill=white, inner sep=1pt},
}
\newcommand{\fr}[1]{\draw[frame] #1;}
\newcommand{\wl}[3]{\node[wire label, anchor=#1] at #2 {$#3$};}
\newcommand{\vl}[3]{\node[val label, anchor=#1] at #2 {$#3$};}
\NewDocumentCommand{\wbox}{o m m m}{
  \draw[cpn, fill=white] #2 rectangle #3;
  \node at \IfNoValueTF{#1}{($#2!0.5!#3$)}{#1} {#4};}

See \cref{sec:pivotal 2-category} for the definition of pivotal dagger 2-category. We first recall the statement of \cref{strictification}.

\strictification*

Our argument combines the approaches of \cite[Appendix B]{Fer24} and \cite[\S2]{NS07}. 

\begin{pf} 
Define $\widehat{\fX}$ as follows.
\begin{itemize} 
\item $\Obj(\widehat{\fX})\coloneq\Obj(\fX)$,
\item $\widehat{\fX}(a\to b)$ consists of ordered tuples $X=((X_1,\varepsilon^X_1),\dots,(X_{|X|},\varepsilon^X_{|X|}))$, possibly the empty tuple $\varnothing_{a}$ if $a=b$, where $|X|$ denotes the length of the tuple $X$, $\varepsilon^X_i\in\{\pm 1\}$, and
\[
\begin{cases}X_i\in\fX(C_{i-1}\to C_i)&\text{if }\varepsilon^X_i=+1,\\X_i\in\fX(C_i\to C_{i-1})&\text{if }\varepsilon^X_i=-1,\end{cases}
\]
for objects $a=C_0,C_1,\ldots,C_{|X|-1},C_{|X|}= b$. Define 1-composition by concatenation:
\[
({}_aX_b,\varepsilon^X)\otimes({}_bY_c,\varepsilon^Y)\coloneq((X_1,\varepsilon^X_1),\ldots,(X_{|X|},\varepsilon^X_{|X|}),(Y_1,\varepsilon^Y_1),\ldots,(Y_{|Y|},\varepsilon^Y_{|Y|})),
\]
with monoidal unit $1_a$ given by the empty tuple $\varnothing_a$.

For a 1-morphism $X\in\widehat{\fX}(a\to b)$, define its \defn{evaluation} in $\fX$ by $\mathbf{ev}(\varnothing_a)\coloneq 1_a$ for $a\in\widehat{\fX}$ and by
\[
\mathbf{ev}(X)\coloneq(\!\cdot\!\!\cdot\!\!\cdot((X^{\varepsilon^X_1}_1\otimes X^{\varepsilon^X_2}_2)\otimes X^{\varepsilon^X_3}_3)\cdots)\qquad\text{where}\qquad
X_i^{\varepsilon^X_i}
\coloneq 
\begin{cases}
X_i^{\vee}&\text{if }\varepsilon^X_i=-1,\\
X_i&\text{if }\varepsilon^X_i=+1.
\end{cases}
\]
For 1-morphisms ${}_aX_b$ and ${}_bY_c$ in $\widehat{\fX}$, define
\[
\mathbf{ev}^2_{X,Y}\in\fX(\mathbf{ev}(X)\otimes\mathbf{ev}(Y)\to\mathbf{ev}(X\otimes Y))
\]
by the unique reparameterization provided by the Mac Lane Coherence Theorem for 2-categories, which is unitary by unitarity of the coherence data in $\fX$.

\item Define $\widehat{\fX}(X\Rightarrow Y)\coloneq\fX(\mathbf{ev}(X)\to\mathbf{ev}(Y))$. For $f\in\widehat{\fX}(X\Rightarrow Y)$ and $g\in\widehat{\fX}(Y\Rightarrow Z)$, define $g\circ f$ by 2-composition in $\fX$. Define $f^\dag$ by the dagger in $\fX$.

For 2-morphisms $f\in\widehat{\fX}({}_aX_b\Rightarrow{}_aY_b)$ and $g\in\widehat{\fX}({}_bW_c\Rightarrow{}_bZ_c)$, define
\[
    f\hat\otimes g\coloneq\mathbf{ev}^2_{Y,Z}\circ(f\otimes g)\circ(\mathbf{ev}^2_{X,W})^{-1}
\]
where $\otimes$ denotes 1-composition in $\fX$.

\end{itemize}
By carrying the $\varepsilon$'s along through the proof of \cite[Prop. B.8]{Fer24}, we find that $\widehat{\fX}$ is a strict dagger 2-category and that $\mathbf{ev}\colon\widehat{\fX}\to\fX$ is an equivalence of dagger 2-categories; $\mathbf{ev}$ has unitary compositors and its unitors are identities.

For $X\in\widehat{\fX}(a\to b)$, we define $X^{\widehat\vee}\in\widehat{\fX}(b\to a)$ by
\[
X^{\widehat{\vee}}\coloneq((X_{|X|},-\varepsilon^X_{|X|}),\dots,(X_1,-\varepsilon^X_1)).
\]
Then the following facts are immediate:
\[
X^{\widehat{\vee}\widehat{\vee}}=X,
\qquad
(X\otimes Y)^{\widehat{\vee}}
=
Y^{\widehat{\vee}}\otimes X^{\widehat{\vee}},
\quad\text{and}\quad 
\varnothing_a^{\widehat{\vee}}=\varnothing_a.
\]

To define our (co)evaluation maps $\widehat{\mathrm{ev}}$ and $\widehat{\mathrm{coev}}$ in $\widehat{\fX}$, we first define for each 1-morphism $X\in\widehat{\fX}(a\to b)$ a unitary
\[
u_X\colon\mathbf{ev}(X^{\widehat{\vee}})\xrightarrow{\sim}\mathbf{ev}(X)^\vee
\]
by induction on the length $|X|$ as follows.
\begin{itemize} 
\item If $|X|=0$, so that $X=\varnothing_a$ for some $a\in\fX$, then set $u_X$ to be the canonical unitor $u_X\coloneq \vee^0_a\in\fX(1_a\to 1^\vee_a)$.
\item If $|X|=1$, set
\[
u_X\coloneq\begin{cases}\id_{X^\vee}&\text{if }\varepsilon^X_1=+1,\\\phi_X&\text{if }\varepsilon^X_1=-1.\end{cases}
\]
\item If $|X|\geq 2$, let $Y\coloneq((X_1,\varepsilon^X_1),\ldots,(X_{|X|-1},\varepsilon^X_{|X|-1}))$ and $Z\coloneq((X_{|X|},\varepsilon^X_{|X|}))$ and observe that $\mathbf{ev}(X)=\mathbf{ev}(Y\otimes Z)\cong\mathbf{ev}(Y)\otimes\mathbf{ev}(Z)$. Indeed, $|Z|=1$ and $\mathbf{ev}$ is defined using left-associated parenthesization. 
\end{itemize} 
Now define $u_X$ by the unitary composite
\begin{equation}\label{uX}
\begin{aligned}
u_X\colon\mathbf{ev}(X^{\hat\vee})=\mathbf{ev}(Z^{\hat\vee}\otimes Y^{\hat\vee})
\xrightarrow{\;\,(\mathbf{ev}^2_{Z^{\hat\vee},Y^{\hat\vee}})^{-1}\;\,}
&\,
\mathbf{ev}(Z^{\hat\vee})\otimes\mathbf{ev}(Y^{\hat\vee})\\
\xrightarrow{\quad\; u_Z\,\otimes\, u_Y\;\quad}
&\,\mathbf{ev}(Z)^\vee\otimes\mathbf{ev}(Y)^\vee
\\
\xrightarrow{\;\;\,\vee^2_{\mathbf{ev}(Y),\mathbf{ev}(Z)}\;\;\,}
&\,
(\mathbf{ev}(Y)\,\otimes\,\mathbf{ev}(Z))^\vee\cong\mathbf{ev}(Y\otimes Z)^\vee=\mathbf{ev}(X)^\vee.
\end{aligned}
\end{equation}
For $X\in\widehat{\fX}$, we will use the following graphical notation:
\[
\id_X=
\begin{tkz}[yscale=0.5]
\draw[ribbon] (0,0) -- (0,2);
\fr{(-0.28,0) -- (-0.28,2)}
\wl{west}{(0,1)}{X}
\end{tkz}
\qquad\text{and}\qquad
\id_{X^\vee}=
\begin{tkz}[yscale=0.5]
\draw[ribbon] (0,2) -- (0,0);
\fr{(0.28,2) -- (0.28,0)}
\wl{east}{(0,1)}{X^\vee}
\end{tkz}\;.
\]
Here the green shading with the red bar on the left (resp. right) indicates that $\mathbf{ev}$ (resp. $\vee\circ\mathbf{ev}$) has been applied to the underlying $\widehat{\fX}$-string diagram.

We will suppress tensorators of $\mathbf{ev}$ from the bottom and top of each diagram. We will also suppress the shadings for the objects in the regions/faces of the diagrams.
 
Now define
\[
\widehat{\mathrm{ev}}_X
\quad
\coloneq
\quad
\begin{tkz}[yscale=-1]
\draw[ribbon] (-1,2) -- (-1,0.5) arc(180:360:1) -- (1,2);
\fr{(-1.28,2) -- (-1.28,1)}
\fr{(-0.72,1) -- (-0.72,0.5) arc(180:360:0.72)}
\fr{(0.72,0.5) -- (0.72,2)}
\rbox{(-1,1)}{0.25}{0.2}{0.2}{$u_X$}
\wl{west}{(-1,1.55)}{X^{\hat\vee}}
\wl{west}{(1,1.55)}{X}
\vl{north}{(-1,2)}{\mathbf{ev}(X^{\widehat\vee})}
\wl{west}{(-1,1.55)}{X^{\hat\vee}}
\vl{north}{(1,2)}{\mathbf{ev}(X)}
\wl{west}{(-1.4,0.45)}{X}
\end{tkz}
\quad
=
\quad
\mathrm{ev}_{\mathbf{ev}(X)}\circ(u_X\otimes \id_{\mathbf{ev}(X)})\circ(\mathbf{ev}^2_{X^{\hat\vee},X})^{-1}
\]
and
\[
\widehat{\mathrm{coev}}_X
\quad
\coloneq
\quad
\begin{tkz}[yscale=-1]
\draw[ribbon] (1,-2) -- (1,-0.5) arc(0:180:1) -- (-1,-2);
\fr{(0.72,-2) -- (0.72,-1)}
\fr{(1.28,-1) -- (1.28,-0.5) arc(0:180:1.28)}
\fr{(-1.28,-0.5) -- (-1.28,-2)}
\rbox{(1,-1)}{0.25}{0.3}{0.3}{$u_X^{-1}$}
\wl{west}{(1,-1.55)}{X^{\hat\vee}}
\wl{west}{(-1,-1.55)}{X}
\wl{west}{(0.6,-0.5)}{X}
\vl{south}{(-1,-2)}{\mathbf{ev}(X)}
\vl{south}{(1,-2)}{\mathbf{ev}(X^{\hat\vee})}
\end{tkz}
\quad
=
\quad
\mathbf{ev}^2_{X,X^{\hat\vee}}\circ(\id_{\mathbf{ev}(X)}\otimes u^{-1}_X)\circ\mathrm{coev}_{\mathbf{ev}(X)}.
\]
These satisfy the zig-zag/snake equations. Indeed,
\[
\begin{tkz}
\draw[ribbon] (-1.5,2.5)node[above=-.8em]{$\mathbf{ev}(X)$} -- (-1.5,-1) arc(180:360:0.75) -- (0,1) arc(180:0:0.75) -- (1.5,-2.5)node[below=-.8em]{$\mathbf{ev}(X)$};
\fr{(-1.78,2.5) -- (-1.78,-1)}
\fr{(-1.78,-1) arc(180:360:1.03)}
\fr{(0.28,-1) -- (0.28,-0.75)}
\fr{(-0.28,-0.25) -- (-0.28,0.25)}
\fr{(0.28,0.75) -- (0.28,1)}
\fr{(0.28,1) arc(180:0:0.47)}
\fr{(1.22,1) -- (1.22,-2.5)}
\rbox{(0,0.5)}{0.3}{0.3}{0.3}{$u_X$}
\rbox{(0,-0.5)}{0.3}{0.35}{0.35}{$u_X^{-1}$}
\wl{west}{(-1.5,1)}{X}
\wl{east}{(0,1.05)}{X}
\wl{east}{(0,-1.05)}{X}
\wl{west}{(0,0)}{X^{\hat\vee}}
\wl{west}{(1.5,-1)}{X}
\end{tkz}
\quad
=
\quad
\begin{tkz}[yscale=1.25]
\draw[ribbon] (-1,2)node[above=-.8em]{$\mathbf{ev}(X)$} -- (-1,-0.5) arc(180:360:0.5) -- (0,0.5) arc(180:0:0.5) -- (1,-2)node[below=-.8em]{$\mathbf{ev}(X)$};
\fr{(-1.28,2) -- (-1.28,-0.5)}
\draw[frame,yscale=0.955](-1.28,-0.5) arc(180:360:0.78);
\fr{(0.28,-0.5) -- (0.28,0.535)}
\draw[frame,yscale=1.06](0.28,0.5) arc(180:0:0.22);
\fr{(0.72,0.535) -- (0.72,-2)}
\wl{west}{(-1,1)}{X}
\wl{east}{(0,0)}{X}
\wl{west}{(1,-1)}{X}
\end{tkz}
=
\quad
\id_X
\;\;\in\;\;
\widehat{\fX}(X\Rightarrow X)
\]
and
\[
\begin{tkz}[yscale=0.8]
\draw[ribbon] (-1.5,-2.5)node[below=-.8em]{$\mathbf{ev}(X^{\hat\vee})$} -- (-1.5,1) arc(180:0:0.75) -- (0,-1) arc(180:360:0.75) -- (1.5,2.5)node[above=-.8em]{$\mathbf{ev}(X^{\hat\vee})$};
\fr{(-1.78,-2.5) -- (-1.78,0)}
\fr{(-1.22,0) -- (-1.22,1)}
\draw[frame,yscale=0.975] (-1.22,1) arc(180:0:0.47);
\fr{(-0.28,1) -- (-0.28,-1.025)}
\draw[frame,yscale=1.025] (-0.28,-1) arc(180:360:1.03);
\fr{(1.78,-1.025) -- (1.78,0)}
\fr{(1.22,0) -- (1.22,2.5)}
\rbox{(-1.5,0)}{0.25}{0.3}{0.3}{$u_X$}
\rbox{(1.5,0)}{0.25}{0.3}{0.3}{$u_X^{-1}$}
\wl{west}{(-1.5,-0.6)}{X^{\hat\vee}}
\wl{east}{(-1.5,0.6)}{X}
\wl{west}{(0,0)}{X}
\wl{east}{(1.5,-0.6)}{X}
\wl{west}{(1.5,0.6)}{X^{\hat\vee}}
\end{tkz}
=
\begin{tkz}
\draw[ribbon] (0,-2)node[below=-.8em]{$\mathbf{ev}(X^{\hat\vee})$} -- (0,2)node[above=-.8em]{$\mathbf{ev}(X^{\hat\vee})$};
\fr{(-0.28,-2) -- (-0.28,-0.5)}
\fr{(0.28,-0.25) -- (0.28,0.25)}
\fr{(-0.28,0.5) -- (-0.28,2)}
\rbox{(0,-0.5)}{0.25}{0.3}{0.3}{$u_X$}
\rbox{(0,0.5)}{0.25}{0.3}{0.3}{$u_X^{-1}$}
\wl{west}{(0,-1)}{X^{\hat\vee}}
\wl{east}{(0,0)}{X}
\wl{west}{(0,1)}{X^{\hat\vee}}
\end{tkz}
=
\quad
\id_{X^{\hat\vee}}\in\widehat{\fX}(X^{\hat\vee}\Rightarrow X^{\hat\vee}).
\]
To see that $\dag\hat\vee=\hat\vee\dag$, observe that for any 2-morphism $f\in\widehat{\fX}({}_aX_b\Rightarrow{}_aY_b)$, the induced 2-morphism $f^{\hat\vee}$ given by
\[
f^{\widehat\vee}
\quad
=
\quad
\begin{tkz}
\draw[ribbon] (-1.5,-2.5)node[below=-0.8em]{$\mathbf{ev}(Y^{\hat\vee})$} -- (-1.5,1) arc(180:0:0.75) -- (0,-1) arc(180:360:0.75) -- (1.5,2.5)node[above=-0.8em]{$\mathbf{ev}(X^{\hat\vee})$};
\fr{(-1.78,-2.5) -- (-1.78,0.5)}
\fr{(-1.22,0.5) -- (-1.22,1)}
\fr{(-1.22,1) arc(180:0:0.47)}
\fr{(-0.28,1) -- (-0.28,-1)}
\fr{(-0.28,-1) arc(180:360:1.03)}
\fr{(1.78,-1) -- (1.78,-0.5)}
\fr{(1.22,-0.5) -- (1.22,2.5)}
\rbox{(-1.5,0.5)}{0.25}{0.3}{0.3}{$u_Y$}
\rbox{(0,0)}{0.25}{0.25}{0.25}{$f$}
\rbox{(1.5,-0.5)}{0.25}{0.3}{0.3}{$u_X^{-1}$}
\wl{west}{(-1.5,0)}{Y^{\hat\vee}}
\wl{east}{(-1.5,1)}{Y}
\wl{west}{(0,0.5)}{Y}
\wl{west}{(0,-0.5)}{X}
\wl{east}{(1.5,-1)}{X}
\wl{west}{(1.5,0.25)}{X^{\hat\vee}}
\end{tkz}
\]
satisfies
\[
f^{\widehat\vee\dag}
\quad=\quad
\begin{tkz}[yscale=-1]
\draw[ribbon] (-1.5,-2.5)node[above=-0.8em]{$\mathbf{ev}(Y^{\hat\vee})$} -- (-1.5,1) arc(180:0:0.75) -- (0,-1) arc(180:360:0.75) -- (1.5,2.5)node[below=-0.8em]{$\mathbf{ev}(X^{\hat\vee})$};
\fr{(-1.78,-2.5) -- (-1.78,0.5)}
\fr{(-1.22,0.5) -- (-1.22,1)}
\fr{(-1.22,1) arc(180:0:0.47)}
\fr{(-0.28,1) -- (-0.28,-1)}
\fr{(-0.28,-1) arc(180:360:1.03)}
\fr{(1.78,-1) -- (1.78,-0.5)}
\fr{(1.22,-0.5) -- (1.22,2.5)}
\rbox{(-1.5,0.5)}{0.25}{0.3}{0.3}{$u_Y^\dag$}
\rbox{(0,0)}{0.25}{0.25}{0.25}{$f^\dag$}
\rbox{(1.5,-0.5)}{0.25}{0.3}{0.3}{$(u_X^{-1})^\dag$}
\wl{west}{(-1.5,0)}{Y^{\hat\vee}}
\wl{east}{(-1.5,1)}{Y}
\wl{west}{(0,0.5)}{Y}
\wl{west}{(0,-0.5)}{X}
\wl{east}{(1.5,-1)}{X}
\wl{west}{(1.5,0.25)}{X^{\hat\vee}}
\end{tkz}
\quad
=
\quad
\begin{tkz}[yscale=-1]
\draw[ribbon] (-1.5,-2.5)node[above=-0.8em]{$\mathbf{ev}(Y^{\hat\vee})$} -- (-1.5,1) arc(180:0:0.75) -- (0,-1) arc(180:360:0.75) -- (1.5,2.5)node[below=-0.8em]{$\mathbf{ev}(X^{\hat\vee})$};
\fr{(-1.78,-2.5) -- (-1.78,0.5)}
\fr{(-1.22,0.5) -- (-1.22,1)}
\fr{(-1.22,1) arc(180:0:0.47)}
\fr{(-0.28,1) -- (-0.28,-1)}
\fr{(-0.28,-1) arc(180:360:1.03)}
\fr{(1.78,-1) -- (1.78,-0.5)}
\fr{(1.22,-0.5) -- (1.22,2.5)}
\rbox{(-1.5,0.5)}{0.25}{0.3}{0.3}{$u_Y^{-1}$}
\rbox{(0,0)}{0.25}{0.25}{0.25}{$f^\dag$}
\rbox{(1.5,-0.5)}{0.25}{0.3}{0.3}{$u_X$}
\wl{west}{(-1.5,0)}{Y^{\hat\vee}}
\wl{east}{(-1.5,1)}{Y}
\wl{west}{(0,0.5)}{Y}
\wl{west}{(0,-0.5)}{X}
\wl{east}{(1.5,-1)}{X}
\wl{west}{(1.5,0.25)}{X^{\hat\vee}}
\end{tkz}
\]
where we used that $u_Y$ and $u_X$ are unitary. By absorbing the $u^{-1}_Y$ and $u_X$ into the $f^\dag$-coupon by standard graphical calculus, then using that $\vee\dag=\dag \vee$ in $\fX$, and finally undoing the absorption, we get that the above is equal to
\[
\begin{tkz}
\draw[ribbon] (-1.5,-2.5)node[below=-0.8em]{$\mathbf{ev}(X^{\hat\vee})$} -- (-1.5,1) arc(180:0:0.75) -- (0,-1) arc(180:360:0.75) -- (1.5,2.5)node[above=-0.8em]{$\mathbf{ev}(Y^{\hat\vee})$};
\fr{(-1.78,-2.5) -- (-1.78,0.5)}
\fr{(-1.22,0.5) -- (-1.22,1)}
\fr{(-1.22,1) arc(180:0:0.47)}
\fr{(-0.28,1) -- (-0.28,-1)}
\fr{(-0.28,-1) arc(180:360:1.03)}
\fr{(1.78,-1) -- (1.78,-0.5)}
\fr{(1.22,-0.5) -- (1.22,2.5)}
\rbox{(-1.5,0.5)}{0.25}{0.3}{0.3}{$u_X$}
\rbox{(0,0)}{0.25}{0.3}{0.3}{$f^\dag$}
\rbox{(1.5,-0.5)}{0.25}{0.3}{0.3}{$u_Y^{-1}$}
\wl{west}{(-1.5,0)}{X^{\hat\vee}}
\wl{east}{(-1.5,1)}{X}
\wl{west}{(0,0.55)}{X}
\wl{west}{(0,-0.5)}{Y}
\wl{east}{(1.5,-1)}{Y}
\wl{west}{(1.5,0.15)}{Y^{\hat\vee}}
\end{tkz}
\quad
=
\quad
f^{\dag\widehat\vee},
\]
as desired.

The functor $\mathbf{ev}$ is UAF-preserving, since the canonical isomorphism in $\fX$ given by
\[
\begin{tkz}[yscale=0.8]
\draw[ribbon] (-1.5,-2.5)node[below=-0.8em]{$\mathbf{ev}(X^{\hat\vee})$} -- (-1.5,1) arc(180:0:0.75) -- (0,-1) arc(180:360:0.75) -- (1.5,2.5)node[above=-0.8em]{$\mathbf{ev}(X)^\vee$};
\fr{(-1.78,-2.5) -- (-1.78,0)}
\fr{(-1.22,0) -- (-1.22,1)}
\draw[frame,yscale=0.975] (-1.22,1) arc(180:0:0.47);
\fr{(-0.28,1) -- (-0.28,-1.025)}
\draw[frame,yscale=1.025] (-0.28,-1) arc(180:360:1.03);
\fr{(1.78,-1.025) -- (1.78,2.5)}
\rbox{(-1.5,0)}{0.3}{0.3}{0.3}{$u_X$}
\wl{west}{(-1.5,-0.7)}{X^{\hat\vee}}
\wl{east}{(-1.5,0.7)}{X}
\wl{west}{(0,0)}{X}
\wl{east}{(1.5,0)}{X}
\end{tkz}
\quad
=
\quad
\begin{tkz}
\draw[ribbon] (0,-2)node[below=-0.8em]{$\mathbf{ev}(X^{\hat\vee})$} -- (0,2)node[above=-0.8em]{$\mathbf{ev}(X)^\vee$};
\fr{(-0.28,-2) -- (-0.28,0)}
\fr{(0.28,0) -- (0.28,2)}
\rbox{(0,0)}{0.3}{0.3}{0.3}{$u_X$}
\wl{west}{(0,-0.5)}{X^{\hat\vee}}
\wl{east}{(0,0.5)}{X}
\end{tkz}
\quad
=
\quad
u_X
\]
is unitary. 

It only remains to show $\widehat{\fX}$ is strictly pivotal, i.e., that for all 1-morphisms ${}_aX_b\in\widehat{\fX}$, the 2-morphism $\hat\phi_X\in\widehat{\fX}(X^{\hat\vee\hat\vee}=X\Rightarrow X)$ given by
\[
\hat\phi_X
\quad
=
\!\!\!
\begin{tkz}
\draw[ribbon] (-1.5,-2.5)node[below=-0.8em]{$\mathbf{ev}(X)=\mathbf{ev}(X^{\hat\vee\hat\vee})$} -- (-1.5,1) arc(180:0:0.75) -- (0,-1) arc(180:360:0.75) -- (1.5,2.5)node[above=-0.8em]{$\mathclap{\mathbf{ev}(X)}$};
\fr{(-1.78,-2.5) -- (-1.78,1)}
\fr{(-1.78,1) arc(180:0:1.03)}
\fr{(0.28,1) -- (0.28,0)}
\fr{(-0.28,0) -- (-0.28,-1)}
\fr{(-0.28,-1) arc(180:360:1.03)}
\fr{(1.78,-1) -- (1.78,-0.5)}
\fr{(1.22,-0.5) -- (1.22,2.5)}
\rbox{(0,0)}{0.3}{0.3}{0.3}{$u_X$}
\rbox{(1.5,-0.5)}{0.3}{0.4}{0.4}{$u_{X^{\hat\vee}}^{-1}$}
\wl{west}{(-1.5,0.55)}{X}
\wl{east}{(0,0.55)}{X}
\wl{west}{(0,-0.5)}{X^{\hat\vee}}
\wl{east}{(1.5,-1)}{X^{\hat\vee}}
\wl{west}{(1.5,0.15)}{X=X^{\hat\vee\hat\vee}}
\end{tkz}
\]
equals $\id_X$.

\begin{claim}
\label{strictclaim}
Since $\hat\phi_{X\hat\otimes Y}=\hat\phi_X\hat\otimes\hat\phi_Y$, it suffices to check $\hat\phi_X=\id_X$ for $|X|\leq 1$.
\end{claim}

We assume \cref{strictclaim} and finish the proof before returning to prove the claim at the end.

When $|X|=0$, setting $u_{\varnothing_a}=\vee^0_a$ gives $\widehat{\mathrm{coev}}_{\varnothing_a}=(\vee^0_a)^{-1}\circ\vee^0_a=\id_{1_a}$, i.e., $\hat\vee^0_a=\id_{\varnothing_a}$, so by \ref{Cphiu} we have $\hat\phi_X=\id_{\varnothing_a}$ as claimed. When $|X|=1$ and $\varepsilon^X_1=+1$, so that $X=((X_1,+1))$ where $X_1$ is a 1-morphism in $\fX$, we have $u_X=\id_{X_1^\vee}$ and $u_{X^{\hat\vee}}=\phi_{X_1}$, so
\[
\hat\phi_X=
\begin{tkz}
\draw[band] (-1.5,-2.2)
node[below=-0.8em]{$\mathbf{ev}(X)$}
-- (-1.5,1) arc(180:0:0.75) -- (0,-1.4);
\draw[band] (1.5,-1.4) -- (1.5,2.5)node[above=-0.8em]{$\mathbf{ev}(X)$};
\fr{(-1.78,-2.2) -- (-1.78,1)}
\fr{(-1.78,1) arc(180:0:1.03)}
\fr{(0.28,1) -- (0.28,0)}
\fr{(-0.28,0) -- (-0.28,-1.4)}
\fr{(1.22,-1.4) -- (1.22,2.5)}
\draw[wire] (-1.5,-2.2) -- (-1.5,1) arc(180:0:0.75) -- (0,-1.4);
\draw[wire] (1.5,-1.5) --node[pos=0.1,right,wire label]{$X^{\hat\vee\hat\vee}{=}X$} (1.5,2.5);
\wbox{(-0.4,-2.0)}{(1.9,-1.4)}{$\widehat{\mathrm{coev}}_{X^{\hat\vee}}$}
\rbox{(0,0)}{0.3}{0.3}{0.3}{$\id_{X_1^\vee}$}
\rbox{(1.5,0.5)}{0.3}{0.3}{0.3}{$\phi_{X_1}^{-1}$}
\wl{west}{(-1.5,0.5)}{X}
\wl{east}{(0,0.5)}{X}
\wl{west}{(0,-0.5)}{X^{\hat\vee}}
\wl{west}{(1.5,1.5)}{X}
\end{tkz}
=
\begin{tkz}
\draw (-1.0,-2.2) -- (-1.0,1.0) arc(180:0:0.75) -- (0.5,-1.4);
\draw (2.0,-1.4) -- (2.0,2.25);
\wbox{(0.1,-2.0)}{(2.4,-1.4)}{$\mathrm{coev}_{X_1^\vee}$}
\rbox{(2.0,0.25)}{0.3}{0.4}{0.4}{$\phi_{X_1}^{-1}$}
\node[font=\scriptsize, anchor=east] at (-1.0,-1.75) {$X_1$};
\node[font=\scriptsize, anchor=east] at (0.5,-0.75) {$X_1^\vee$};
\node[font=\scriptsize, anchor=east] at (2.0,-0.75) {$X_1^{\vee\vee}$};
\node[font=\scriptsize, anchor=west] at (2.0,2) {$X_1$};
\end{tkz}
=\phi^{-1}_{X_1}\circ\phi_{X_1}=\id_{X_1}=\id_{\mathbf{ev}(X)}.
\]
When $|X|=1$ and $\varepsilon^X_1=-1$, so that $\mathbf{ev}(X)=X^\vee_1$, $\mathbf{ev}(X^{\hat\vee})=X_1$, $u_X=\phi_{X_1}$, and $u_{X^{\hat\vee}}=\id_{X^\vee_1}$, then
\[
\hat\phi_X\,=\!
\begin{tkz}
\draw[band] (-1.5,-2.2)
node[below=-0.8em]{$\mathbf{ev}(X)$}
-- (-1.5,1) arc(180:0:0.75) -- (0,-1.4);
\draw[band] (1.5,-1.4) -- (1.5,2.5)node[above=-0.8em]{$\mathbf{ev}(X)$};
\fr{(-1.78,-2.2) -- (-1.78,1)}
\fr{(-1.78,1) arc(180:0:1.03)}
\fr{(0.28,1) -- (0.28,0)}
\fr{(-0.28,0) -- (-0.28,-1.4)}
\fr{(1.22,-1.4) -- (1.22,2.5)}
\draw[wire] (-1.5,-2.2) -- (-1.5,1) arc(180:0:0.75) -- (0,-1.4);
\draw[wire] (1.5,-1.5) --node[pos=0.1,right,wire label]{$X^{\hat\vee\hat\vee}{=}X$} (1.5,2.5);
\wbox{(-0.4,-2.0)}{(1.9,-1.4)}{$\widehat{\mathrm{coev}}_{X^{\hat\vee}}$}
\rbox{(0,0)}{0.3}{0.3}{0.3}{$\phi_{X_1}$}
\rbox{(1.5,0.5)}{0.3}{0.3}{0.3}{$\id_{X_1^\vee}$}
\wl{west}{(-1.5,0.5)}{X}
\wl{east}{(0,0.5)}{X}
\wl{west}{(0,-0.5)}{X^{\hat\vee}}
\wl{west}{(1.5,1.5)}{X}
\end{tkz}
\!=\,\,
\begin{tkz}
\draw (-1.0,-2.2) -- (-1.0,1.0) arc(180:0:0.75) -- (0.5,-1.4);
\draw (2.0,-1.4) -- (2.0,2.25); 
\wbox{(0.1,-2.0)}{(2.4,-1.4)}{$\mathrm{coev}_{X_1}$}
\rbox{(0.5,0.25)}{0.3}{0.3}{0.3}{$\phi_{X_1}$}
\node[font=\scriptsize, anchor=west] at (0.5,1) {$X_1^{\vee\vee}$};
\node[font=\scriptsize, anchor=west] at (-1.0,-1.75) {$X_1^\vee$};
\node[font=\scriptsize, anchor=east] at (0.5,-0.75) {$X_1$};
\node[font=\scriptsize, anchor=east] at (2.0,1.75) {$X_1^\vee$};
\end{tkz}
=\phi^\vee_{X_1}\circ\phi_{X_1^\vee}\underset{\ref{CphiV}}{=}\id_{X_1^\vee}
=
\id_{\mathbf{ev}(X)}
\]
as desired. This completes the proof of \cref{strictification}.
\end{pf}

\begin{pf}[Proof of \cref{strictclaim}]
\label{strictclaim-pf}
We first extend the graphical calculus to allow several blue strands.

\itemstep{$u$ splits over concatenation.} For composable ${}_aX_b$ and ${}_bY_c$ in $\widehat{\fX}$, we observe that $u_{X\hat\otimes Y}=u_Y\otimes u_X$ as 2-morphisms $\mathbf{ev}((X\hat\otimes Y)^{\hat\vee})=\mathbf{ev}(Y^{\hat\vee}\hat\otimes X^{\hat\vee})\Rightarrow \mathbf{ev}(X\hat\otimes Y)^\vee$ in $\fX$.
\begin{equation}\label{usplit}
\begin{tkz}
\draw[band] (-0.28,-1.5) -- (-0.28,1.5);
\draw[band] (0.28,-1.5) -- (0.28,1.5);
\fr{(-0.56,-1.5) -- (-0.56,-0.28)}
\fr{(0.56,0.28) -- (0.56,1.5)}
\draw[wire] (-0.05,-1.5) -- (-0.05,1.5);
\draw[wire] (0.05,-1.5) -- (0.05,1.5);
\wbox{(-0.8,-0.28)}{(0.8,0.28)}{$u_{X\hat\otimes Y}$}
\node[wire label, left=-3mm] at (-0.35,-0.9) {$Y^{\hat\vee}$};
\node[wire label, right=-3mm] at (0.35,-0.9) {$X^{\hat\vee}$};
\node[wire label, left=-3mm] at (-0.35,0.9) {$X$};
\node[wire label, right=-3mm] at (0.35,0.9) {$Y$};
\end{tkz}
\qquad=\qquad
\begin{tkz}[yscale=0.75]
\draw[ribbon] (0,-2)
-- (0,2);
\fr{(-0.28,-2) -- (-0.28,0)}
\fr{(0.28,0) -- (0.28,2)}
\rbox{(0,0)}{0.3}{0.2}{0.2}{\normalsize$u_Y$}
\wl{west}{(0,-1.25)}{Y^{\hat\vee}}
\wl{east}{(0,1.25)}{Y}
\end{tkz}
\;
\begin{tkz}[yscale=0.75]
\draw[ribbon] (0,-2)
-- (0,2);
\fr{(-0.28,-2) -- (-0.28,0)}
\fr{(0.28,0) -- (0.28,2)}
\rbox{(0,0)}{0.3}{0.2}{0.2}{\normalsize$u_X$}
\wl{west}{(0,-1.25)}{X^{\hat\vee}}
\wl{east}{(0,1.25)}{X}
\end{tkz}
\end{equation}
Indeed, we proceed by induction on $|Y|$. If $|X|=0$ or $|Y|=0$, the corresponding wire is empty and the coupon sitting on it is $\vee^0$, which we have suppressed, so the two sides agree. If $|Y|=1$, then \eqref{usplit} is just the definition of $u_{X\hat\otimes Y}$ from \cref{uX}. If $|Y|\geq 2$, write $Y=Y'\hat\otimes Z$ with $|Z|=1$. Then 
\[
\begin{tkz}
    \def\pa{(-0.56,-1.6) -- (-0.56,1.6)}
    \def\pb{(0,-1.6) -- (0,1.6)}
    \def\pc{(0.56,-1.6) -- (0.56,1.6)}
    \draw[band] \pa;
    \draw[band] \pb;
    \draw[band] \pc;
    \fr{(-0.84,-1.6) -- (-0.84,-0.28)}
    \fr{(0.84,0.28) -- (0.84,1.6)}
    \draw[wire] \pa;
    \draw[wire] \pb;
    \draw[wire] \pc;
    \wbox{(-1.1,-0.28)}{(1.1,0.28)}{$u_{X\hat\otimes Y}$}
    \wl{west}{(-0.48,-1.05)}{Z^{\hat\vee}}
    \wl{east}{(-0.64,1.05)}{Z}
    \wl{west}{(0.08,-1.05)}{Y'^{\hat\vee}}
    \wl{east}{(-0.04,1.05)}{Y'}
    \wl{west}{(0.64,-1.05)}{X^{\hat\vee}}
    \wl{east}{(0.48,1.05)}{X}
\end{tkz}
\quad=\quad
\begin{tkz}
    \def\pa{(-1.2,-1.6) -- (-1.2,1.6)}
    \def\pb{(0,-1.6) -- (0,1.6)}
    \def\pc{(0.56,-1.6) -- (0.56,1.6)}
    \draw[band] \pa;
    \draw[band] \pb;
    \draw[band] \pc;
    \fr{(-1.48,-1.6) -- (-1.48,-0.28)}
    \fr{(-0.92,0.28) -- (-0.92,1.6)}
    \fr{(-0.28,-1.6) -- (-0.28,-0.28)}
    \fr{(0.84,0.28) -- (0.84,1.6)}
    \draw[wire] \pa;
    \draw[wire] \pb;
    \draw[wire] \pc;
    \wbox{(-1.75,-0.28)}{(-0.65,0.28)}{$u_Z$}
    \wbox{(-0.4,-0.28)}{(0.96,0.28)}{$u_{X\hat\otimes Y'}$}
    \wl{west}{(-1.12,-1.05)}{Z^{\hat\vee}}
    \wl{east}{(-1.28,1.05)}{Z}
    \wl{west}{(0.08,-1.05)}{Y'^{\hat\vee}}
    \wl{east}{(-0.08,1.05)}{Y'}
    \wl{west}{(0.64,-1.05)}{X^{\hat\vee}}
    \wl{east}{(0.48,1.05)}{X}
\end{tkz}
\quad=\quad
\begin{tkz}
    \def\pa{(-1.2,-1.6) -- (-1.2,1.6)}
    \def\pb{(0,-1.6) -- (0,1.6)}
    \def\pc{(1.2,-1.6) -- (1.2,1.6)}
    \draw[band] \pa;
    \draw[band] \pb;
    \draw[band] \pc;
    \fr{(-1.48,-1.6) -- (-1.48,-0.28)}
    \fr{(-0.92,0.28) -- (-0.92,1.6)}
    \fr{(-0.28,-1.6) -- (-0.28,-0.28)}
    \fr{(0.28,0.28) -- (0.28,1.6)}
    \fr{(0.92,-1.6) -- (0.92,-0.28)}
    \fr{(1.48,0.28) -- (1.48,1.6)}

    \draw[wire] \pa;
    \draw[wire] \pb;
    \draw[wire] \pc;
    \wbox{(-1.75,-0.28)}{(-0.65,0.28)}{$u_Z$}
    \wbox{(-0.55,-0.28)}{(0.55,0.28)}{$u_{Y'}$}
    \wbox{(0.65,-0.28)}{(1.75,0.28)}{$u_X$}
    \wl{west}{(-1.12,-1.05)}{Z^{\hat\vee}}
    \wl{east}{(-1.28,1.05)}{Z}
    \wl{west}{(0.08,-1.05)}{Y'^{\hat\vee}}
    \wl{east}{(-0.08,1.05)}{Y'}
    \wl{west}{(1.28,-1.05)}{X^{\hat\vee}}
    \wl{east}{(1.12,1.05)}{X}
\end{tkz}
\quad=\quad
\begin{tkz}
    \def\pa{(-0.8,-1.6) -- (-0.8,1.6)}
    \def\pb{(-0.24,-1.6) -- (-0.24,1.6)}
    \def\pc{(0.96,-1.6) -- (0.96,1.6)}
    \draw[band] \pa;
    \draw[band] \pb;
    \draw[band] \pc;
    \fr{(-1.08,-1.6) -- (-1.08,-0.28)}
    \fr{(0.04,0.28) -- (0.04,1.6)}
    \fr{(0.68,-1.6) -- (0.68,-0.28)}
    \fr{(1.24,0.28) -- (1.24,1.6)}
    \draw[wire] \pa;
    \draw[wire] \pb;
    \draw[wire] \pc;
    \wbox{(-1.2,-0.28)}{(0.2,0.28)}{$u_Y$}
    \wbox{(0.5,-0.28)}{(1.42,0.28)}{$u_X$}
    \wl{west}{(-0.72,-1.05)}{Z^{\hat\vee}}
    \wl{east}{(-0.88,1.05)}{Z}
    \wl{west}{(-0.16,-1.05)}{Y'^{\hat\vee}}
    \wl{east}{(-0.32,1.05)}{Y'}
    \wl{west}{(1.04,-1.05)}{X^{\hat\vee}}
    \wl{east}{(0.88,1.05)}{X}
\end{tkz}
\]
where the first and third equalities are applications of the definition \cref{uX} and the second equality is by the inductive hypothesis.

\itemstep{$\hat\phi$ is multiplicative.} We compute that
\def\tkzscl{0.65}
\[
\begin{aligned}
\hat\phi_{X\hat\otimes Y}
&=
\scalebox{\tkzscl}{
\begin{tkz}
\def\pa{(-1.53,-3.4) -- (-1.53,1) arc(180:0:0.93) -- (0.33,-1) arc(180:360:0.83) -- (1.99,3.4)}
\def\pb{(-1.43,-3.4) -- (-1.43,1) arc(180:0:0.83) -- (0.23,-1) arc(180:360:0.93) -- (2.09,3.4)}
\draw[band] \pa;
\draw[band] \pb;
\fr{(-1.81,-3.4) -- (-1.81,1) arc(180:0:1.21) -- (0.61,0.28)}
\fr{(-0.05,-0.28) -- (-0.05,-1) arc(180:360:1.21) -- (2.37,-0.68)}
\fr{(1.71,-0.12) -- (1.71,3.4)}
\draw[wire] \pa;
\draw[wire] \pb;
\wbox[(0.28,0)]{(-0.3,-0.28)}{(1.06,0.28)}{$u_{X\hat\otimes Y}$}
\wbox[(2.14,-0.4)]{(1.36,-0.78)}{(2.92,0.08)}{$u^{-1}_{(X\hat\otimes Y)^{\hat\vee}}$}
\wl{center}{(-1.69,-1.4)}{X}
\wl{center}{(-1.27,-1.4)}{Y}
\wl{center}{(0.07,0.62)}{Y}
\wl{center}{(0.5,0.62)}{X}
\wl{center}{(0,-0.85)}{Y^{\hat\vee}}
\wl{center}{(0.63,-0.85)}{X^{\hat\vee}}
\wl{center}{(1.83,1.2)}{X}
\wl{center}{(2.25,1.2)}{Y}
\end{tkz}}
\!\underset{\eqref{usplit}}{=}
\scalebox{\tkzscl}{
\begin{tkz}
\def\pa{(-2.75,-3.4) -- (-2.75,1.1) arc(180:0:1.65) -- (0.55,-1.1) arc(180:360:0.55) -- (1.65,3.4)}
\def\pb{(-1.65,-3.4) -- (-1.65,1.1) arc(180:0:0.55) -- (-0.55,-1.1) arc(180:360:1.65) -- (2.75,3.4)}
\draw[band] \pa;
\draw[band] \pb;
\fr{(-3.03,-3.4) -- (-3.03,1.1) arc(180:0:1.93) -- (0.83,-0.17)}
\fr{(0.27,-0.73) -- (0.27,-1.1) arc(180:360:0.83) -- (1.93,0.17)}
\fr{(1.37,0.73) -- (1.37,3.4)}
\fr{(-1.93,-3.4) -- (-1.93,1.1) arc(180:0:0.83) -- (-0.27,0.73)}
\fr{(-0.83,0.17) -- (-0.83,-1.1) arc(180:360:1.93) -- (3.03,1.42)}
\fr{(2.47,1.98) -- (2.47,3.4)}
\draw[wire] \pa;
\draw[wire] \pb;
\wbox{(-1.0,0.17)}{(-0.1,0.73)}{$u_Y$}
\wbox{(0.1,-0.73)}{(1.0,-0.17)}{$u_X$}
\wbox{(1.05,0.2)}{(2.25,0.93)}{$u^{-1}_{X^{\hat\vee}}$}
\wbox{(2.15,1.38)}{(3.35,2.1)}{$u^{-1}_{Y^{\hat\vee}}$}
\wl{west}{(-2.68,-2.3)}{X}
\wl{west}{(-1.58,-2.3)}{Y}
\wl{east}{(-0.62,0.95)}{Y}
\wl{west}{(-0.48,-0.05)}{Y^{\hat\vee}}
\wl{east}{(0.48,0.95)}{X}
\wl{west}{(0.62,-1.0)}{X^{\hat\vee}}
\wl{east}{(1.61,-0.05)}{X^{\hat\vee}}
\wl{west}{(2.88,2.5)}{Y}
\end{tkz}
}
\!=
\scalebox{\tkzscl}{
\begin{tkz}
\draw[ribbon] (-1.2,-2) -- (-1.2,1) arc(180:0:0.6) -- (0,-1) arc(180:360:0.6) -- (1.2,2);
\fr{(-1.48,-2) -- (-1.48,1) arc(180:0:0.88) -- (0.28,0.25)}
\fr{(-0.28,-0.25) -- (-0.28,-1) arc(180:360:0.88) -- (1.48,-0.65)}
\fr{(0.92,-0.15) -- (0.92,2)}
\wbox{(-0.45,-0.25)}{(0.45,0.25)}{$u_X$}
\wbox[(1.2,-0.4)]{(0.65,-0.75)}{(1.65,-0.05)}{$u^{-1}_{X^{\hat\vee}}$}
\wl{center}{(-1.56+0.2,-1.4)}{X}
\wl{center}{(-0.36+0.2,0.6)}{X}
\wl{center}{(0.26,-0.7)}{X^{\hat\vee}}
\begin{scope}[xshift=3.3cm]
\draw[ribbon] (-1.2,-2) -- (-1.2,1) arc(180:0:0.6) -- (0,-1) arc(180:360:0.6) -- (1.2,2);
\fr{(-1.48,-2) -- (-1.48,1) arc(180:0:0.88) -- (0.28,0.25)}
\fr{(-0.28,-0.25) -- (-0.28,-1) arc(180:360:0.88) -- (1.48,-0.65)}
\fr{(0.92,-0.15) -- (0.92,2)}
\wbox{(-0.45,-0.25)}{(0.45,0.25)}{$u_Y$}
\wbox[(1.2,-0.4)]{(0.65,-0.75)}{(1.65,-0.05)}{$u^{-1}_{Y^{\hat\vee}}$}
\wl{center}{(-1.56+0.2,-1.4)}{Y}
\wl{center}{(-0.36+0.2,0.6)}{Y}
\wl{center}{(0.26,-0.7)}{Y^{\hat\vee}}
\end{scope}
\end{tkz}
}
=
\hat\phi_X\hat\otimes\hat\phi_Y.
\end{aligned}
\]
Since every ${}_aX_b\in\widehat{\fX}$ is the concatenation $X^{(1)}\hat\otimes\cdots\hat\otimes X^{(|X|)}$ for some length-1 $X^{(i)}\coloneq((X_i,\varepsilon^X_i))$, it follows that $\hat\phi_X=\hat\phi_{X^{(1)}}\hat\otimes\cdots\hat\otimes\hat\phi_{X^{(|X|)}}$.
\end{pf}

\bibliographystyle{alphaurl}
\bibliography{bibliography}

@book{EGNO15,
  author    = {Etingof, Pavel and Gelaki, Shlomo and Nikshych, Dmitri and Ostrik, Victor},
  title     = {Tensor Categories},
  series    = {Mathematical Surveys and Monographs},
  volume    = {205},
  publisher = {American Mathematical Society},
  address   = {Providence, RI},
  year      = {2015},
  isbn      = {978-1-4704-2024-6},
  doi       = {10.1090/surv/205},
  mrnumber  = {3242743},
}

@article{MW12,
  author        = {Morrison, Scott and Walker, Kevin},
  title         = {Blob homology},
  journal       = {Geometry \& Topology},
  year          = {2012},
  volume        = {16},
  number        = {3},
  pages         = {1481--1607},
  doi           = {10.2140/gt.2012.16.1481},
  eprint        = {1009.5025},
  archivePrefix = {arXiv},
  primaryClass  = {math.AT},
  mrnumber      = {2978449},
}

@misc{W21,
  author        = {Walker, Kevin},
  title         = {A universal state sum},
  year          = {2021},
  eprint        = {2104.02101},
  archivePrefix = {arXiv},
  primaryClass  = {math.QA},
}

@misc{W06,
  author       = {Walker, Kevin},
  title        = {{TQFT}s [early incomplete draft]},
  year         = {2006},
  note         = {Unpublished notes, version 1h, May 11, 2006},
  howpublished = {\url{https://canyon23.net/math/tc.pdf}},
}

@article{3Hilb,
  author        = {Chen, Quan and Ferrer, Giovanni and Hungar, Brett and Penneys, David and Sanford, Sean},
  title         = {Manifestly unitary higher {H}ilbert spaces},
  journal       = {Journal of the London Mathematical Society},
  volume        = {113},
  number        = {5},
  pages         = {e70532},
  year          = {2026},
  doi           = {10.1112/jlms.70532},
  eprint        = {2410.05120},
  archivePrefix = {arXiv},
  primaryClass  = {math.QA},
}

@article{SS23,
  author        = {Stehouwer, Luuk and Steinebrunner, Jan},
  title         = {Dagger categories via anti-involutions and positivity},
  journal       = {Theory and Applications of Categories},
  year          = {2024},
  volume        = {41},
  number        = {56},
  pages         = {2013--2040},
  doi           = {10.70930/tac/b1qms3v0},
  eprint        = {2304.02928},
  archivePrefix = {arXiv},
  primaryClass  = {math.CT},
}

@article{CHPJP22,
  author        = {Chen, Quan and Hern{\'a}ndez Palomares, Roberto and Jones, Corey and Penneys, David},
  title         = {Q-system completion for {C$^*$} 2-categories},
  journal       = {Journal of Functional Analysis},
  volume        = {283},
  number        = {3},
  pages         = {109524},
  year          = {2022},
  doi           = {10.1016/j.jfa.2022.109524},
  eprint        = {2105.12010},
  archivePrefix = {arXiv},
  primaryClass  = {math.OA},
}

@misc{UQSL,
  author  = {Ferrer, Giovanni and Kawagoe, Kyle and Penneys, David},
  title   = {Unitary Quantum Symmetries Lite},
  year    = {2026},
  url     = {https://people.math.osu.edu/penneys.2/UQSL/UQSL.html},
  note    = {Draft book in progress},
}

@article{FH21,
  author        = {Freed, Daniel S. and Hopkins, Michael J.},
  title         = {{R}eflection positivity and invertible topological phases},
  year          = {2021},
  journal       = {Geometry \& Topology},
  volume        = {25},
  number        = {3},
  pages         = {1165--1330},
  doi           = {10.2140/gt.2021.25.1165},
  eprint        = {1604.06527},
  archivePrefix = {arXiv},
  primaryClass  = {hep-th},
}

@misc{bases,
      title={Orthonormal bases for higher {H}ilbert spaces}, 
      author={Ferrer, Giovanni and Hungar, Brett and Penneys, David and Wesley, Greyson},
      year={2026},
      eprint={2608.11358},
      archivePrefix={arXiv},
      primaryClass={math.QA},
}

@book{KS77,
  author    = {Kirby, Robion C. and Siebenmann, Laurence C.},
  title     = {Foundational Essays on Topological Manifolds, Smoothings, and Triangulations},
  series    = {Annals of Mathematics Studies},
  number    = {88},
  publisher = {Princeton University Press},
  address   = {Princeton, NJ},
  year      = {1977},
  isbn      = {978-1-4008-8150-5},
  doi       = {10.1515/9781400881505},
}

@article{Lurie09,
  author        = {Lurie, Jacob},
  title         = {On the classification of topological field theories},
  journal       = {Current Developments in Mathematics},
  volume        = {2008},
  number        = {1},
  pages         = {129--280},
  year          = {2009},
  publisher     = {International Press of Boston},
  address       = {Somerville, MA},
  doi           = {10.4310/CDM.2008.v2008.n1.a3},
  eprint        = {0905.0465},
  archivePrefix = {arXiv},
  primaryClass  = {math.CT},
  mrnumber      = {2555928},
}

@article{BD95,
  author        = {Baez, John C. and Dolan, James},
  title         = {Higher-dimensional algebra and topological quantum field theory},
  journal       = {Journal of Mathematical Physics},
  volume        = {36},
  number        = {11},
  pages         = {6073--6105},
  year          = {1995},
  doi           = {10.1063/1.531236},
  eprint        = {q-alg/9503002},
  archivePrefix = {arXiv},
  primaryClass  = {math.QA},
  mrnumber      = {1355899},
}

@article{Baez97,
  author        = {Baez, John C.},
  title         = {Higher-dimensional algebra {II}. 2-{H}ilbert spaces},
  journal       = {Advances in Mathematics},
  volume        = {127},
  number        = {2},
  pages         = {125--189},
  year          = {1997},
  doi           = {10.1006/aima.1997.1617},
  eprint        = {q-alg/9609018},
  archivePrefix = {arXiv},
  primaryClass  = {math.QA},
  mrnumber      = {1448713},
}

@misc{FHJF24,
  author        = {Ferrer, Giovanni and Hungar, Brett and Johnson-Freyd, Theo and Krulewski, Cameron and M{\"u}ller, Lukas and {Nivedita} and Penneys, David and Reutter, David and Scheimbauer, Claudia and Stehouwer, Luuk and Vuppulury, Chetan},
  title         = {Dagger $n$-categories},
  year          = {2024},
  eprint        = {2403.01651},
  archivePrefix = {arXiv},
  primaryClass  = {math.CT},
}

@incollection{NS07,
  author        = {Ng, Siu-Hung and Schauenburg, Peter},
  title         = {Higher {F}robenius-{S}chur indicators for pivotal categories},
  booktitle     = {Hopf Algebras and Generalizations},
  series        = {Contemporary Mathematics},
  volume        = {441},
  pages         = {63--90},
  publisher     = {American Mathematical Society},
  address       = {Providence, RI},
  year          = {2007},
  doi           = {10.1090/conm/441/08500},
  eprint        = {math/0503167},
  archivePrefix = {arXiv},
  primaryClass  = {math.QA},
  mrnumber      = {2381536},
}

@misc{Fer24,
  author        = {Ferrer, Giovanni},
  title         = {Foundations for operator algebraic tricategories},
  year          = {2024},
  eprint        = {2404.05193},
  archivePrefix = {arXiv},
  primaryClass  = {math.OA},
}

@phdthesis{SP09,
  author        = {Schommer-Pries, Christopher J.},
  title         = {The Classification of Two-Dimensional Extended Topological Field Theories},
  school        = {University of California, Berkeley},
  type          = {Ph.{D}. thesis},
  year          = {2009},
  eprint        = {1112.1000},
  archivePrefix = {arXiv},
  primaryClass  = {math.AT},
  note          = {Revised version: arXiv v2 (2014)},
}

@article{pivbicatSNpaper,
  author        = {Fuchs, J\"{u}rgen and Schweigert, Christoph and Yang, Yang},
  title         = {String-net models for pivotal bicategories},
  journal       = {Theory and Applications of Categories},
  year          = {2025},
  volume        = {44},
  number        = {17},
  pages         = {474--543},
  issn          = {1201-561X},
  doi           = {10.70930/tac/5jrwhjw1},
  eprint        = {2302.01468},
  archivePrefix = {arXiv},
  primaryClass  = {math.QA},
}

@article{morrison2011higher,
  author        = {Morrison, Scott and Walker, Kevin},
  title         = {Higher categories, colimits, and the blob complex},
  journal       = {Proc. Natl. Acad. Sci. USA},
  volume        = {108},
  number        = {20},
  pages         = {8139--8145},
  year          = {2011},
  doi           = {10.1073/pnas.1018168108},
  eprint        = {1108.5386},
  archivePrefix = {arXiv},
  primaryClass  = {math.CT},
  mrnumber      = {2806651},
}

@misc{walkeruscslides,
  author       = {Walker, Kevin},
  title        = {{TQFT}s and contact structures},
  year         = {2008},
  note         = {Talk slides, USC, February 2008},
  howpublished = {\sloppy\url{https://canyon23.net/math/talks/USC\%20contact.pdf}},
}

@misc{DR18,
  author        = {Douglas, Christopher L. and Reutter, David J.},
  title         = {Fusion 2-categories and a state-sum invariant for 4-manifolds},
  year          = {2018},
  eprint        = {1812.11933},
  archivePrefix = {arXiv},
  primaryClass  = {math.QA},
  note          = {To appear in Mem. Amer. Math. Soc.},
}

@misc{Hai,
  author        = {Ha{\"\i}oun, Benjamin},
  title         = {Defining extended {TQFT}s via handle attachments},
  year          = {2024},
  eprint        = {2412.14649},
  archivePrefix = {arXiv},
  primaryClass  = {math.GT},
  note          = {To appear in Algebr. Geom. Topol.},
}

@article{Pen20,
  author        = {Penneys, David},
  title         = {Unitary dual functors for unitary multitensor categories},
  journal       = {Higher Structures},
  year          = {2020},
  volume        = {4},
  number        = {2},
  pages         = {22--56},
  doi           = {10.21136/HS.2020.09},
  eprint        = {1808.00323},
  archivePrefix = {arXiv},
  primaryClass  = {math.QA},
}

@article{W17aasen,
  author        = {Aasen, David and Lake, Ethan and Walker, Kevin},
  title         = {Fermion condensation and super pivotal categories},
  journal       = {Journal of Mathematical Physics},
  year          = {2019},
  volume        = {60},
  number        = {12},
  pages         = {121901},
  doi           = {10.1063/1.5045669},
  eprint        = {1709.01941},
  archivePrefix = {arXiv},
  primaryClass  = {cond-mat.str-el},
}

@misc{HPT24,
  author        = {Henriques, Andr\'{e} and {Nivedita} and Penneys, David},
  title         = {Complete {${W}^*$}-categories},
  year          = {2024},
  eprint        = {2411.01678},
  archivePrefix = {arXiv},
  primaryClass  = {math.OA},
}

@misc{Z24,
  author        = {Liu, Zhengwei},
  title         = {{Functional Integral Construction of Topological Quantum Field Theory}},
  year          = {2024},
  eprint        = {2409.17103},
  archivePrefix = {arXiv},
  primaryClass  = {math-ph},
}

@article{CP22,
  author  = {Chen, Quan and Penneys, David},
  title   = {Q-system completion is a 3-functor},
  journal = {Theory and Applications of Categories},
  volume  = {38},
  number  = {4},
  pages   = {101--134},
  year    = {2022},
  issn    = {1201-561X},
  note    = {\href{https://arxiv.org/abs/2106.12437}{arXiv:2106.12437 [math.QA]}}
}

@misc{Xu26,
      title={Frobenius Algebras and Dual Bimodules in Monoidal 2-Categories}, 
      author={Xu, Hao},
      year={2026},
      eprint={2606.02046},
      archivePrefix={arXiv},
      primaryClass={math.QA},
}

@article{ST12,
  author       = {Stolz, Stephan and Teichner, Peter},
  title        = {Traces in Monoidal Categories},
  journal = {Transactions of the American Mathematical Society},
  year         = {2012},
  volume       = {364},
  number       = {8},
  pages        = {4425--4464},
  doi          = {10.1090/S0002-9947-2012-05615-7},
  issn         = {0002-9947}
}

@misc{M26,
      title={Topological defects in reflection positive topological field theories}, 
      author={M{\"u}ller, Lukas},
      year={2026},
      eprint={2608.07217},
      archivePrefix={arXiv},
      primaryClass={math-ph},
}

@unpublished{MS26,
  author = {M{\"u}ller, Lukas and Stehouwer, Luuk},
  title  = {Reflection positivity for once extended topological field theories},
  year   = {2026},
  note   = {in preparation}
}

@misc{W09fields,
  author  = {Walker, Kevin},
  title   = {Fields, blobs and {TQFT}s},
  howpublished = {Talk at UC Berkeley FRG meeting, January},
  note    = {Slides available at \url{https://canyon23.net/math/talks/ucb\%20frg\%20200901.pdf}},
  year    = {2009}
}

@misc{W14premodular,
  author  = {Walker, Kevin},
  title   = {Premodular {TQFT}s},
  howpublished = {Talk at ESI, Vienna, February},
  note    = {Slides available at \url{https://canyon23.net/math/talks/ESI\%20201402.pdf}},
  year    = {2014}
}

@misc{RW,
  author = {Reutter, David and Walker, Kevin},
  title  = {Nonsemisimple {TQFT}s and handle induction},
  note   = {In preparation}
}

@misc{R20,
  author       = {Reutter, David},
  title        = {From non-unital skein theory to modified traces and non-semisimple {TQFT}s},
  howpublished = {Seminar talk, Hausdorff Institute for Mathematics, Bonn},
  year         = {2020},
  month        = dec,
  note         = {Slides at \url{https://drive.google.com/file/d/1IslyZdNdYoAaqe_GRZ60uLEN9IOiAiEM/view}}
}

@article{Ati88,
  author  = {Atiyah, Michael},
  title   = {Topological quantum field theories},
  journal = {Inst. Hautes {\'E}tudes Sci. Publ. Math.},
  volume  = {68},
  pages   = {175--186},
  year    = {1988},
  doi  = {10.1007/BF02698547},
  url  = {http://www.numdam.org/item/PMIHES_1988__68__175_0/},
}

@incollection{Seg04,
  author    = {Segal, Graeme},
  title     = {The definition of conformal field theory},
  booktitle = {Topology, Geometry and Quantum Field Theory},
  series    = {London Math. Soc. Lecture Note Ser.},
  volume    = {308},
  pages     = {421--577},
  publisher = {Cambridge Univ. Press},
  year      = {2004},
  doi       = {10.1017/CBO9780511526398.019}
}

@misc{LZ26,
      title={Axiomatization of the {L}evin--{W}en wave function}, 
      author={Liu, Zhengwei and Zhao, Zishuo},
      year={2026},
      eprint={2608.21109},
      archivePrefix={arXiv},
      primaryClass={math-ph},
}

@article{DW90,
  author  = {Dijkgraaf, Robbert and Witten, Edward},
  title   = {Topological gauge theories and group cohomology},
  journal = {Communications in Mathematical Physics},
  volume  = {129},
  number  = {2},
  pages   = {393--429},
  year    = {1990},
  doi     = {10.1007/BF02096988},
}

@incollection{Q95,
  author       = {Quinn, Frank},
  title        = {Lectures on Axiomatic Topological Quantum Field Theory},
  booktitle    = {Geometry and Quantum Field Theory},
  series       = {IAS/Park City Mathematics Series},
  volume       = {1},
  pages        = {323--453},
  year         = {1995},
  publisher    = {American Mathematical Society},
  address      = {Providence, RI},
  doi          = {10.1090/pcms/001/05}
}

@article{LR97,
  author  = {Longo, Roberto and Roberts, John E.},
  title   = {A theory of dimension},
  journal = {$K$-Theory},
  volume  = {11},
  number  = {2},
  pages   = {103--159},
  year    = {1997},
}

@article{GL19,
  author        = {Giorgetti, Luca and Longo, Roberto},
  title         = {Minimal index and dimension for 2-{$C^*$}-categories with finite-dimensional centers},
  journal       = {Communications in Mathematical Physics},
  volume        = {370},
  number        = {2},
  pages         = {719--757},
  year          = {2019},
  eprint        = {1805.09234},
  archivePrefix = {arXiv},
}

@article{V20,
  author  = {Verdon, Dominic},
  title   = {Unitary pseudonatural transformations},
  journal = {High. Struct.},
  volume  = {9},
  number  = {1},
  year    = {2025},
  pages   = {1--35},
  eprint  = {2004.12760},
}

@misc{Sel11,
  author = {Selinger, Peter},
  title  = {A survey of graphical languages for monoidal categories},
  year = {2011},
  note   = {New Structures for Physics, Lecture Notes in Physics 813, Springer (2011), 289--355},
  eprint = {0908.3347},
  archivePrefix = {arXiv},
}

@book{categoricalquantum,
  title     = {Categories for Quantum Theory: An Introduction},
  author    = {Heunen, Chris and Vicary, Jamie},
  year      = {2019},
  publisher = {Oxford University Press},
  series    = {Oxford Graduate Texts in Mathematics},
  volume    = {28},
  isbn      = {9780198739616}
}
\end{document}